\documentclass[10pt,twoside,openright]{book}

\pdfoutput=1

\usepackage[palatino]{quotchap}
\usepackage{a4wide}
\usepackage{amsmath,amssymb,amsthm}
\usepackage[sc]{mathpazo}
\usepackage[T1]{fontenc}
\usepackage{tikz,pgfplots}
\usepackage[english]{babel}
\usepackage{graphicx}
\usepackage{subcaption}
\usepackage[numbers]{natbib}
\usepackage{tabularx}
\usepackage{enumitem}
\usepackage{color} 
\usepackage{appendix}
\usepackage{tikz}
\usepackage{geometry}
\usepackage{bigstrut}
\usepackage{multirow}
\usepackage[]{algorithm2e}
\usepackage{bbm}
\usepackage{tocvsec2}
\maxtocdepth{section}

\makeatletter
\renewcommand{\@makechapterhead}[1]{%
  \chapterheadstartvskip%
  \vspace{-3cm}
  {\size@chapter{\sectfont\raggedleft
    {\chapnumfont
      \ifnum \c@secnumdepth >\m@ne%
        \if@mainmatter\thechapter\else\phantom{\thechapter}%
      \fi\else\phantom{\thechapter}\fi
      \par\nobreak}%
    {\raggedleft\advance\leftmargin10em\interlinepenalty\@M #1\par}}
  \nobreak\chapterheadendvskip}}
\makeatother

\newtheorem{proposition}{Proposition}[chapter]
\newtheorem{theorem}{Theorem}[chapter]
\newtheorem{lemma}{Lemma}[chapter]
\newtheorem{corollary}{Corollary}[chapter]
\newtheorem{assumption}{Assumption}[chapter]
\newtheorem{conjecture}{Conjecture}[chapter]

\theoremstyle{remark}
\newtheorem{remark}{Remark}[chapter]

\def \E {\mathbb{E}}
\def \P {\mathbb{P}}
\def \Var {{\rm Var}\,}
\def \F {\Phi}
\def \f {\varphi}
\def \Pois {{\rm Pois}}
\def \l {\lambda}
\def \sl {{s_\lambda}}
\def \ee {{\rm e}}
\def \eps {\varepsilon}
\def \N {\mathcal{N}}
\def \b {\beta}
\def \d {\delta}
\def \l {\lambda}
\def \m {\mu}
\def \g {\gamma}
\def \a {\alpha}

\def \QY {\hat Q_1^b(t)}
\def \QH {\hat Q_1^h(t)}
\def \CY {\hat Q_2^b(t)}
\def \CH {\hat Q_2^h(t)}

\def \qy {q_1^b}
\def \qh {q_1^h}
\def \cy {q_2^b}
\def \ch {q_2^h}

\def \dd {{\rm d}}

\newcommand{\la} {\lambda}
\newcommand{\aaa} {\alpha}
\newcommand{\ka} {\kappa}
\newcommand{\thh} {\theta}
\newcommand{\s} {\sigma}
\newcommand{\de} {\delta}
\newcommand{\e} {\varepsilon}

\newcommand{\hx} {\hat{\gamma}}
\newcommand{\hy} {\hat{\beta}}
\newcommand{\mui} {\mu_\infty^\star}
\newcommand{\muT} {\mu_T^\star}
\newcommand{\Xlm} {X_{\mu}}
\newcommand{\Qlm} {Q_{\mu}}

\newcommand{\klm} {\kappa_{\mu}}
\newcommand{\iy} {\infty}
\newcommand{\Pih} {\hat{\Pi}_{T}}
\newcommand{\muh} {\hat{\mu}_T^\star}

\newcommand{\Rt}{R_{\rm tot}}

\newcommand{\Rightarrowd}{{\;\buildrel{d}\over\Rightarrow\;}}
\newcommand{\equalD}{{\;\buildrel{d}\over= \;}}

\def \less[#1] {\stackrel{(#1)}{\leq}}

\newenvironment{chapterstart}
{
\vspace{-1cm}
\noindent\rule[0.5ex]{\linewidth}{2pt}
\vspace{0.5cm}

\begin{minipage}{0.75\textwidth}

}{

\end{minipage}

\vspace{1 cm}
\noindent\rule[0.5ex]{\linewidth}{2pt}
}

\usepackage{fancyhdr}

\definecolor{col1}{rgb}{0.8, 0.01568627450980392, 0.}
\definecolor{col2}{rgb}{0.9372549019607843, 0.6274509803921569, 0.1}
\definecolor{col3}{rgb}{0.9686274509803922, 0.87, 0.}
\definecolor{col4}{rgb}{0.3176470588235294, 0.58, 0.0784313725490196}
\definecolor{col5}{rgb}{0.22745098039215686, 0.23921568627450981, 
  0.43}
\definecolor{col6}{rgb}{0.5019607843137255, 0.16, 0.4470588235294118}

\usepackage{pgfplots}
    \pgfplotsset{
        table/search path={Introduction/Version_4,Chapter_2,Chapter_3,Chapter_4,Chapter_5,Chapter_6,Chapter_7}
    }
\pgfplotsset{compat=1.14}

\begin{document}

\renewcommand{\headrulewidth}{0pt}
\fancyhead[LE,RO]{ }

\vspace*{\fill}

\begin{center}
\vspace{-7mm}
{\huge \bfseries Asymptotic dimensioning of \\
stochastic service systems
}
\end{center}

\vspace*{\fill}

\newpage

\noindent
This work is part of the Free Competition grant no.~613.001.213, which is financed by the Netherlands Organisation
for Scientific Research (NWO).\\ \\*


\vspace*{\fill}

\noindent
A catalogue record is available from the Eindhoven University of Technology\\* Library.\\
ISBN: 978-90-386-4255-0\\

\noindent
This thesis is number D 205 of the thesis series of the Beta Research School for\\* Operations Management and Logistics.\\

\noindent
Cover design by Kim Mathijsen.

\newpage

\vspace*{1.6cm}

\begin{center}
{
\fontsize{12}{1.05}\selectfont

{\Large \bfseries
Asymptotic dimensioning of \\
stochastic service systems
}

\vspace{1.6 cm}

PROEFSCHRIFT

\vspace{1.6 cm}

ter verkrijging van de graad van doctor aan de\\*
Technische Universiteit Eindhoven, op gezag van de\\*
rector magnificus prof.dr.ir. F.P.T. Baaijens, voor\\*
een commissie aangewezen door het College voor\\*
Promoties, in het openbaar te verdedigen op\\*
 donderdag 18 mei 2017 om 16.00 uur

\vspace{1.6 cm}

door

\vspace{1.6 cm}

Britt Walthera Johanna Mathijsen

\vspace{0.8 cm}

geboren te 's-Hertogenbosch

}
\end{center}

\vspace*{\fill}

\newpage


\noindent
Dit proefschrift is goedgekeurd door de promotoren en de samenstelling van de promotiecommissie is als volgt: \\*

\begin{tabular}{ll}
voorzitter: & prof.dr.ir. B. Koren\\
1e promotor: & prof.dr. J.S.H. van Leeuwaarden\\
2e promotor: & prof.dr. A.P. Zwart\\
leden:		& prof.dr. S. Bhulai (Vrije Universiteit Amsterdam)\\
			& dr. J.E. Reed (New York University Stern)\\
			& prof.dr. A.G. de Kok \\
			& prof.dr.ir. O.J. Boxma \\
adviseur(s):	& dr. G.B. Yom-Tov (Technion)\\
 &
\end{tabular}

\vspace*{\fill}

\noindent
Het onderzoek dat in dit proefschrift wordt beschreven is uitgevoerd in overeenstemming met de TU/e Gedragscode Wetenschapsbeoefening.\\*

\frontmatter

\fancyhf{}
\fancyhead[LE,RO]{\thepage}
\fancyhead[LO]{\itshape\nouppercase{\rightmark}}
\fancyhead[RE]{\itshape\nouppercase{\leftmark}}
\renewcommand{\headrulewidth}{0.4pt}

\addcontentsline{toc}{chapter}{Acknowledgments}

\chapter*{Acknowledgments}

\vspace{-1cm}
\noindent\rule[0.5ex]{\linewidth}{2pt}
\vspace{0.5cm}

\noindent
The completion of this thesis would not have been possible without the many people who have supported and encouraged me along the way. I take this opportunity to express my gratitude towards them.

First and foremost, I am greatly indebted to my supervisors Johan van Leeuwaarden and Bert Zwart. Johan, I think I can safely say that you have been my guiding light throughout this academic journey. 
Ever since my bachelor project, your endless optimism, patience and sense of perspective have been exactly what I needed to keep me motivated through the ups and downs that come with research. 
Thank you for teaching me to be persistent, critical and positive. 
Bert, thank you for the continuous flow of ingenious ideas that came during our many discussions. 
I truly admire your drive, enthusiasm and passion for mathematics. 

Much of work presented in this thesis is the result of fruitful collaborations with some great researchers, to whom I wish to express my appreciation as well. 
First of all, Onno Boxma, your inspiring lectures on stochastic processes may well have marked the starting point of the path that led to this PhD, and I am very happy that you have been willing to be a member of my defense committee.
It has been delightful to work with you, together with Shaul Bar-Lev and David Perry, on the subject of Chapter 7.
Guido Janssen, thank you for getting me acquainted with the mathematical theory behind asymptotics. 
Your thoroughness, preciseness and willingness to share your expertise are much appreciated.
Galit Yom-Tov, thank you for the nice and intensive collaboration that resulted in Chapter 5, and for serving on my committee. 
Also, I hold nice memories on our joint effort to organize the YEQT workshop together with Jan-Pieter.

I am also thankful to Sandjai Bhulai, Ton de Kok and Josh Reed for being part of my defense committee, for their time to read this thesis and for providing me with valuable feedback. 

During my visits to the Technion, I have been warmly welcomed by Avi Mandelbaum and the people at the SEELab. 
It has been fascinating to get a glimpse of the goldmine of service system data that they have collected and I am grateful to them for showing me the way. 
Avi, thank you for kindly introducing me to the Israeli culture, and for your many advices on research and academic life. 

I owe many thanks to my office mates, Fabio, Fiona and Thomas, for creating the most pleasant and comforting work environment. 
I am pleasantly surprised you have been able to tolerate my ever-fluctuating stress levels, certainly in the last couple of months.
The fact that I have remained (reasonably) sane throughout these years is in large part thanks to you.  

Outside my office door, I have been lucky to be part of an amazing research group. 
I want to thank my fellow PhDs of the Stochastics section, including the former generations, for creating and sustaining such an amiable atmosphere. 
Further, I am grateful to Remco van der Hofstad and Marko Boon for granting me the opportunity to develop my teaching and lecturing skills, which has brought refreshing variation into my job as a PhD student. 

Having been a member of the departmental PhD Council, I got the chance to meet many bright young minds, and it has been great fun to collaborate with them to strengthen the PhD community. 
I wish all of them the best for the future.

When you stick around at the same university for almost nine years, you run the risk of befriending some mathematicians. 
I am nevertheless very happy to have met Christine, Jorg, Jorn, Laura, Mark, Thomas and many others at this place. 
Thank you all for making my time at TU/e very enjoyable.

Rik, thank you for holding my hand throughout these years. You have been my rock. 

Finally, to my parents Jan and Willemien and my sister Kim, I owe my deepest gratitude.
Your unconditional support and comfort have kept me grounded and realistic. 
I could never have done this without you. \\

\begin{flushright}
\textit{
Britt Mathijsen\\
March 2017}
\end{flushright}

\tableofcontents

\mainmatter

\chapter{Introduction}

\begin{chapterstart}
Stochastic service systems describe situations in which customers compete for service from scarce resources. Think of check-in lines at airports, waiting rooms in hospitals or queues in supermarkets, where the scarce resource is human manpower.
Next to these traditional settings, resource sharing is also important in large-scale service systems such as the internet, wireless networks and cloud computing facilities. 
In these virtual environments, geographical location does not play a restricting role on the system size, paving the way for the emergence of large-scale resource sharing networks.
This thesis investigates how to design large-scale systems in order to achieve economies-of-scale, by which we mean that the system is highly occupied and hence utilizes efficiently the expensive resources, while at the same time, the offered service levels remain high. 
In this introductory chapter, we give an overview of the available machinery that supports such principles and explain how this thesis contributes to the existing study of large-scale service systems. 
A crucial concept behind most of the results discussed in the chapter is the Central Limit Theorem (CLT) -- arguably one of the most important theorems in mathematics and science.
\end{chapterstart}
 
\newpage

\section{Service systems \& queueing theory}

\subsection{Quality vs. Efficiency}

Large-scale service systems take many shapes and forms. 
Classical examples of large-scale service systems include call centers \cite{Erlang1917,Palm1957,Whitt1999,Gans2003,Borst2004,Brown2005,Zeltyn2005,Bassamboo2009,Khudyakov2006} and communication systems \cite{Kleinrock1976,Anick1982,Kelly1985,Kleinrock2007,johanthesis}. 
More recently, congestion-related issues in health care facilities and cloud-computing facilities have received much attention \cite{Armony2015,Green2007,YomTov2010,Gupta2007,Tan2012}. 
In all settings, one can think of service systems as being composed of \textit{customers} and \textit{servers}. 
In call centers, customers typically call to request help from one of the agents (servers). 
In communication networks, the data packets are the customers and the communication channels are the servers.
In health care facilities, patients are the customers, and nurses/physicians are the servers. 
The system scale may refer to the size of the client base it caters to, or the magnitude of its capacity, or both. 
Next to the central notions of customers and servers, we emphasize that service systems are inherently stochastic, that is, subject to uncertainty.
Although arrival volumes can be anticipated to some extent over a certain planning horizon, for instance through historical data and forecasting methods, one cannot predict with certainty future arrival patterns. 
Moreover, service requirements are typically random as well, adding more uncertainty.
This intrinsic stochastic variability is a predominant cause of delay experienced by customers in the system.

Due to the inherent randomness in both their arrival and service processes, stochastic models have proved instrumental in both quantifying and improving the operational performance of service systems. 
Queueing theory and stochastics provide the mathematical tools to describe and evaluate these service systems. 
Queueing models are often able to capture and explain fundamental phenomena that are common across applications.

A standard model for service systems is the $M/GI/s$ queue, which we will refer to as the \textit{many-server} queue. 
This model assumes that customers arrive to the queue according to a Poisson process with rate $\lambda$, and customer service times are mutually independent and identically distributed (i.i.d.) samples from the distribution of a non-negative random variable $B$.
The parameter $s$ denotes the number of servers in the system, and hence restricts the number of simultaneous services. 
The case $s=1$ corresponds to a single-server queue. 

First principles say that the queueing process is stable, that is, the number of customers does not explode as time evolves, if and only if the expected workload $R := \lambda\E[B]$ brought into the system per time unit is strictly less than the system capacity. 
In other words, the \textit{utilization} of the queue, defined as $\rho := \lambda\E[B] / s$ should remain strictly below one.  
Naturally, a system manager prefers to operate at a utilization level close to one, so that resources are used efficiently. 
However, it is known that pushing the occupation levels to 100\% leads to an explosive increase in congestion.
That is, the expected queue length and customer waiting time increase indefinitely, thereby reducing the quality-of-service (QoS) and also customer satisfaction. 
These seemingly conflicting objectives give rise to a classical trade-off between customer satisfaction and costs of resources.

\subsection{Economies-of-scale}

Under the assumption that service times are exponentially distributed with mean $1/\mu$, the many-server queue reduces to the well-studied $M/M/s$ queue. 
Despite its simplicity, the analysis of the $M/M/s$ queue explains mathematically the distinctive traits of queues in general, such as the non-linear effect of utilization on the queue size, and pooling effects.

Let $W^{(s)}$ denote the waiting time of a customer and $Q^{(s)}$ the queue length (including the customers in service) in the steady-state $M/M/s$ queue. Without loss of generality, we fix $\mu=1$, so that $\rho = \l/s$.
A straightforward balance argument gives the stationary distribution:

\begin{equation}
\label{eq:MMs_stationary_distribution}
\pi_k := \P( Q^{(s)} = k ) 
= \left\{
\begin{array}{ll}
\pi_0\frac{\l^k}{k!}, & \text{if } k < s, \\
\pi_0\frac{\l^s}{s!}\,\rho^{k-s} & \text{if } k \geq s,
\end{array}
\right.
\end{equation}
where
\begin{equation*}
\pi_0 := \Big( \sum_{k=0}^{s-1} \frac{\l^k}{k!} + \frac{1}{1-\rho} \frac{\l^s}{s!}\Big)^{-1}.
\end{equation*} 
Natural QoS indicators include the expected waiting time $\E[W^{(s)}]$ and the delay probability $\P(W^{(s)}>0)$. 
Invoking the PASTA (Poisson arrivals see time averages) property \cite{Wolff1982}, we know that the delay probability equals the probability of the queue length being greater or equal to the number of servers $s$. 
Thus,
\begin{equation}
\label{eq:MMs_wait}
\P(W^{(s)} > 0) = \P(Q^{(s)}\geq s) = \frac{\l^s}{s!} \Big( (1-\rho) \sum_{k=0}^{s-1} \frac{\l^k}{k!} + \frac{\l^s}{s!} \Big)^{-1}.
\end{equation}
By Little's law, which says that $\E[(Q^{(s)}-s)^+] =\l\E[W^{(s)}]$, we furthermore have
\begin{equation}
\E[W^{(s)}] = \P(W^{(s)} > 0)\,\frac{1/s}{1-\rho}.
\label{eq:MMs_wait2}
\end{equation}
From these formulae, it is readily seen that $\P(W^{(s)} > 0) \to 1$ and $\E[W^{(s)}] \to \infty$ as $\rho \uparrow 1$ . That is, increasing $\lambda$ to $s$, while keeping the latter fixed, leads to a system in which all customers are delayed before service, and the expected delay before reaching a server increases to infinity. 

\begin{figure}
\centering
\begin{subfigure}{0.48\textwidth}
\centering
\begin{tikzpicture}[scale=0.75]
\begin{axis}[
	xmin = 0, 
	xmax = 10.2,
	ymin = 0,
	ymax = 1.05,
	axis line style={->},
	axis lines = left,
	xlabel = {$\to s$},
	ylabel = {$\E[W^{(s)}] / \E[W^{(1)}]$},
	y label style = {at = {(axis cs: -0.9,0.5)}},
	xscale=1,
	yscale=0.8,
	legend style = {at = {(axis cs: 10,1)},anchor = north east}]
	
\addplot[thick] table{
1	1.
2	0.375
3	0.19708
4	0.119601
5	0.0787172
6	0.0546017
7	0.0392985
8	0.0290713
9	0.0219674
10	0.0168832
};
\addplot[thick, dashed] table {
1	1.
2	0.888889
3	0.808989
4	0.745541
5	0.692641
6	0.647215
7	0.607422
8	0.572056
9	0.540278
10	0.511475
};

\addplot[thick, dotted] table {
1	1.
2	0.487179
3	0.318249
4	0.234584
5	0.1848
6	0.151858
7	0.128489
8	0.111075
9	0.0976122
10	0.0869037
};

\legend{{\small $\lambda=0.60$},{\small $\lambda = 0.80$},{\small $\lambda = 0.95$}}
\end{axis}
\end{tikzpicture}
\caption{Expected waiting time}
\end{subfigure}
\hspace{2mm}
\begin{subfigure}{0.48\textwidth}
\begin{tikzpicture}[scale=0.75]
\begin{axis}[
	xmin = 0, 
	xmax = 10.2,
	ymin = 0,
	ymax = 1.05,
	axis line style={->},
	axis lines = left,
	xlabel = {\small $\to s$},
	ylabel = {\small $\P(W^{(s)}>0) / \P(W^{(1)}>0)$},
	y label style = {at = {(axis cs: -0.9,0.5)}},
	xscale=1,
	yscale=0.8,
	legend style={at={(0.03,0.05)},anchor=south west}]
	
\addplot[thick] table {
1	1.
2	0.75
3	0.591241
4	0.478405
5	0.393586
6	0.32761
7	0.275089
8	0.23257
9	0.197707
10	0.168832
};
\addplot[thick, dashed] table {
1	1.
2	0.888889
3	0.808989
4	0.745541
5	0.692641
6	0.647215
7	0.607422
8	0.572056
9	0.540278
10	0.511475
};
\addplot[thick, dotted] table {
1	1.
2	0.974359
3	0.954746
4	0.938336
5	0.923999
6	0.911146
7	0.899423
8	0.888599
9	0.878509
10	0.869037
};

\legend{{\small $\lambda=0.60$},{\small $\lambda = 0.80$},{\small $\lambda = 0.95$}}
\end{axis}

\end{tikzpicture}
\caption{Probability of delay}
\end{subfigure}
\caption{Effects of resource pooling in the $M/M/s$ queue.}
\label{fig:waiting_time_pooling}
\end{figure}
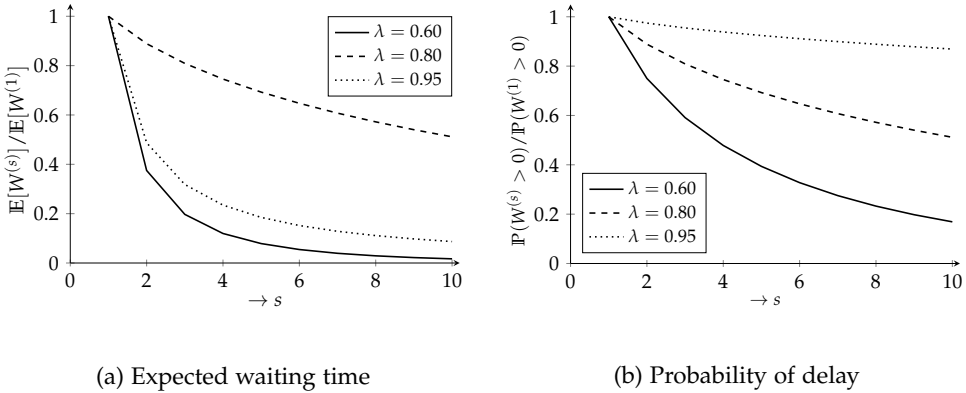

The $M/M/s$ queue also reveals the effect of \textit{resource pooling}.
To illustrate the operational benefits of sharing resources, we compare a system of $s$ separate $M/M/1$ queues, each serving a Poisson arrival stream with rate $\lambda<1$, against one $M/M/s$ queue facing arrival rate $\lambda s$. 
The two systems thus experience the same total workload and utilization, namely $\rho = \lambda$.
We fix the value of $\lambda$ and vary $s$. 
Obviously, the waiting time and queue length distribution in the first scenario are unaffected by the parameter $s$, since there is no interaction between the single-server queues. 
This lack of coordination tolerates a scenario of having an idle server, while the total number of customers in the system exceeds $s$, therefore wasting resource capacity. 
Such an event cannot happen in the many-server scenario, due to the central queue. 
This central coordination improves QoS. Indeed Figure \ref{fig:waiting_time_pooling} shows that the reduction in expected waiting time can be substantial.

\noindent
So pooling kills two birds with one stone: QoS for customers improves and the system efficiency increases. 

\subsection{Many-server scaling regimes}
\label{sec:intro_many_server_regimes}
Now that we know that economies-of-scale can be achieved, it is relevant to ask how to match capacity $s$ to a demand $\lambda$ in the setting where both $s$ and $\lambda$ become large.
The expressions in \eqref{eq:MMs_wait} and \eqref{eq:MMs_wait2} provide a starting point for finding such demand-capacity relations, particularly when we apply asymptotic analysis  for $s\to\infty$, \cite{Halfin1981,Borst2004,Reed2009}. 
Asymptotic theory of many-server systems relies on the prerequisite that the limiting behavior of the service system is determined by the way in which capacity $s$ is adjusted to demand, assuming demand grows large.
We illustrate this idea by investigating typical sample paths of the queue length process $Q = \{Q(t)\}_{t\geq 0}$ of an $M/M/s$ queue for increasing values of $\lambda$. 

\begin{figure}[b!]
\centering
\begin{tikzpicture}[scale = 0.8]
\begin{axis}[
	xmin = 0, 
	xmax = 8,
	ymin = 0,
	ymax = 10,
	axis line style={->},
	axis lines = left,
	xlabel = $\to t$,
	ylabel = {$Q(t)$},
	y label style = {at = {(axis cs: -0.5,6.6)}},
	xscale=1.5,
	yscale=0.8]
	
\addplot[draw = black, fill = col1, fill opacity = 0.2] table {
0	4
1.90023	4
1.90023	5
2.03028	5
2.03028	6
2.24549	6
2.24549	5
2.32361	5
2.32361	4
2.56087	4
2.56087	5
2.59708	5
2.59708	6
2.63833	6
2.63833	7
2.66002	7
2.66002	6
3.03161	6
3.03161	5
3.17381	5
3.17381	6
3.22186	6
3.22186	7
3.37289	7
3.37289	8
3.60466	8
3.60466	9
3.88682	9
3.88682	8
4.00394	8
4.00394	9
4.03197	9
4.03197	8
4.06907	8
4.06907	9
4.20129	9
4.20129	8
4.35099	8
4.35099	7
4.55316	7
4.55316	6
4.57055	6
4.57055	7
4.59485	7
4.59485	6
4.633	6
4.633	5
4.70679	5
4.70679	4
5.20235	4
5.20235	5
5.23937	5
5.23937	4
5.51929	4
5.51929	5
5.6359	5
5.6359	6
5.67363	6
5.67363	7
5.71895	7
5.71895	6
5.84485	6
5.84485	7
5.97328	7
5.97328	6
5.99596	6
5.99596	5
6.27451	5
6.27451	4
8.24592	4
};
\addplot[draw = black, fill = col4, fill opacity = 0.2] table {
0	4
0	3
0.169639	3
0.169639	2
0.570362	2
0.570362	3
0.756359	3
0.756359	2
1.13351	2
1.13351	1
1.60218	1
1.60218	2
1.63589	2
1.63589	3
1.73512	3
1.73512	4
4.78508	4
4.78508	3
4.87791	3
4.87791	4
5.26881	4
5.26881	3
5.30479	3
5.30479	4
6.37399	4
6.37399	3
6.39698	3
6.39698	2
6.58434	2
6.58434	3
6.66658	3
6.66658	2
6.67963	2
6.67963	3
6.85082	3
6.85082	2
6.88581	2
6.88581	3
7.17614	3
7.17614	4
7.3415	4
7.3415	3
7.62396	3
7.62396	2
7.96256	2
7.96256	1
8.24592	1
8.24592	2
8.24292	4

};
\addplot[dashed, very thick] coordinates { (0,4) (8.3,4) };
\end{axis}
\end{tikzpicture}
\caption{Sample path of the $M/M/s$ queue with $\lambda = 3$ and $s=4$.}
\label{fig:sample_path_small}
\end{figure}
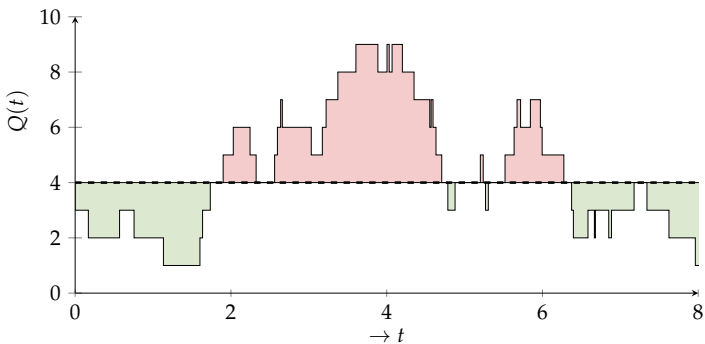

Figure \ref{fig:sample_path_small} depicts a sample path for $\lambda = 3$ and $s = 4$. 
The number of customers queueing at time $t$ is given by $(Q(t)-s)^+$ with $(\cdot)^+ := \max\{0,\cdot\}$.
The number of idle servers is given by $(s-Q(t))^+$. 
In Figure \ref{fig:sample_path_small}, the red and green area hence represent the cumulative queue length and cumulative number of idle servers, respectively, over the given time period.
Bearing in mind the dual goal of QoS and efficiency, we want to minimize both of these areas simultaneously.

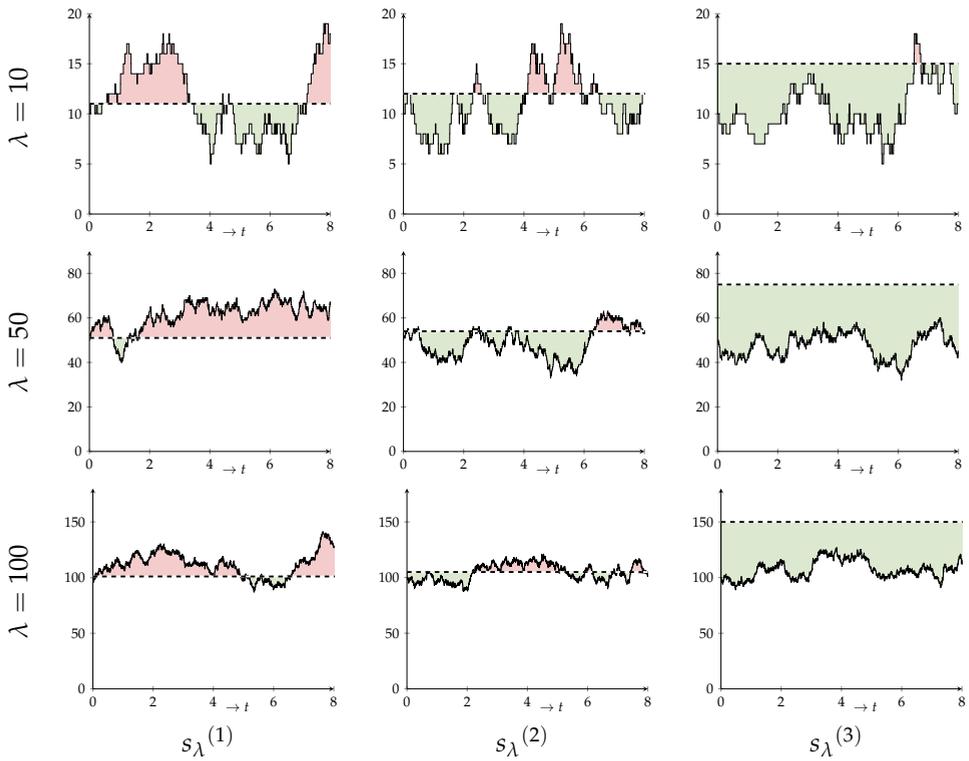
\begin{figure}
\centering
\begin{tabular}{cccc} 
\rotatebox{90}{\ \ \ \ \ \ \ \ \ \ \ \ $\lambda=10$}
&
\begin{tikzpicture}[scale = 0.58]
\begin{axis}[
	xmin = 0, 
	xmax = 8,
	ymin = 0,
	ymax = 20,
	axis line style={->},
	axis lines = left,
	x label style = {at = {(0.75,-0.05)}},
	xlabel = $\to t$,
	xscale=0.8,
	yscale=0.8]
\small	
\addplot[draw = black, fill = col1, fill opacity = 0.2] file {Introduction/Version_R1/tikz_tex/files/samplePath_10_ED_upper.txt};
\addplot[draw = black, fill = col4, fill opacity = 0.2] file {Introduction/Version_R1/tikz_tex/files/samplePath_10_ED_lower.txt};
\addplot[white, very thick] coordinates { (0,11) (8.3,11) };
\addplot[dashed, very thick] coordinates { (0,11) (8.3,11) };
\end{axis}
\end{tikzpicture}
&
\begin{tikzpicture}[scale = 0.58]
\begin{axis}[
	xmin = 0, 
	xmax = 8,
	ymin = 0,
	ymax = 20,
	axis line style={->},
	axis lines = left,
	x label style = {at = {(0.75,-0.05)}},
	xlabel = $\to t$,
	xscale=0.8,
	yscale=0.8]
\small
\addplot[draw = black, fill = col1, fill opacity = 0.2] file {Introduction/Version_R1/tikz_tex/files/samplePath_10_QED_upper.txt};
\addplot[draw = black, fill = col4, fill opacity = 0.2] file {Introduction/Version_R1/tikz_tex/files/samplePath_10_QED_lower.txt};
\addplot[white, very thick] coordinates { (0,12) (8.3,12) };
\addplot[dashed, very thick] coordinates { (0,12) (8.3,12) };
\end{axis}
\end{tikzpicture}
&
\begin{tikzpicture}[scale = 0.58]
\begin{axis}[
	xmin = 0, 
	xmax = 8,
	ymin = 0,
	ymax = 20,
	axis line style={->},
	axis lines = left,
	x label style = {at = {(0.75,-0.05)}},
	xlabel = $\to t$,
	xscale=0.8,
	yscale=0.8]
\small	
\addplot[draw = black, fill = col1, fill opacity = 0.2] file {Introduction/Version_R1/tikz_tex/files/samplePath_10_QD_upper.txt};
\addplot[draw = black, fill = col4, fill opacity = 0.2] file {Introduction/Version_R1/tikz_tex/files/samplePath_10_QD_lower.txt};
\addplot[white, very thick] coordinates { (0,15) (8.3,15) };
\addplot[dashed, very thick] coordinates { (0,15) (8.3,15) };
\end{axis}
\end{tikzpicture}
\\
\rotatebox{90}{\ \ \ \ \ \ \ \ \ \ \ $\lambda=50$}
&
\begin{tikzpicture}[scale = 0.58]
\begin{axis}[
	xmin = 0, 
	xmax = 8,
	ymin = 0,
	ymax = 90,
	axis line style={->},
	axis lines = left,
	x label style = {at = {(0.75,-0.05)}},
	xlabel = $\to t$,
	xscale=0.8,
	yscale=0.8]
\small	
\addplot[draw = black, fill = col1, fill opacity = 0.2] file {Introduction/Version_R1/tikz_tex/files/samplePath_50_ED_upper.txt};
\addplot[draw = black, fill = col4, fill opacity = 0.2] file {Introduction/Version_R1/tikz_tex/files/samplePath_50_ED_lower.txt};
\addplot[white, very thick] coordinates { (0,51) (8.3,51) };
\addplot[dashed, very thick] coordinates { (0,51) (8.3,51) };
\end{axis}
\end{tikzpicture}
&
\begin{tikzpicture}[scale = 0.58]
\begin{axis}[
	xmin = 0, 
	xmax = 8,
	ymin = 0,
	ymax = 90,
	axis line style={->},
	axis lines = left,
	x label style = {at = {(0.75,-0.05)}},
	xlabel = $\to t$,
	xscale=0.8,
	yscale=0.8]
\small	
\addplot[draw = black, fill = col1, fill opacity = 0.2] file {Introduction/Version_R1/tikz_tex/files/samplePath_50_QED_upper.txt};
\addplot[draw = black, fill = col4, fill opacity = 0.2] file {Introduction/Version_R1/tikz_tex/files/samplePath_50_QED_lower.txt};
\addplot[white, very thick] coordinates { (0,54) (8.3,54) };
\addplot[dashed, very thick] coordinates { (0,54) (8.3,54) };
\end{axis}
\end{tikzpicture}
&
\begin{tikzpicture}[scale = 0.58]
\begin{axis}[
	xmin = 0, 
	xmax = 8,
	ymin = 0,
	ymax = 90,
	axis line style={->},
	axis lines = left,
	x label style = {at = {(0.75,-0.05)}},
	xlabel = $\to t$,
	xscale=0.8,
	yscale=0.8]
\small	
\addplot[draw = black, fill = col4, fill opacity = 0.2] file {Introduction/Version_R1/tikz_tex/files/samplePath_50_QD_lower.txt};
\addplot[white, very thick] coordinates { (0,75) (8.3,75) };
\addplot[dashed, very thick] coordinates { (0,75) (8.3,75) };
\end{axis}
\end{tikzpicture}
\\
\rotatebox{90}{\ \ \ \ \ \ \ \ \ \ $\lambda=100$}
 &
\begin{tikzpicture}[scale = 0.58]
\begin{axis}[
	xmin = 0, 
	xmax = 8,
	ymin = 0,
	ymax = 180,
	axis line style={->},
	axis lines = left,
	x label style = {at = {(0.75,-0.05)}},
	xlabel = $\to t$,
	xscale=0.8,
	yscale=0.8]
\small	
\addplot[draw = black, fill = col1, fill opacity = 0.2] file {Introduction/Version_R1/tikz_tex/files/samplePath_100_ED_upper.txt};
\addplot[draw = black, fill = col4, fill opacity = 0.2] file {Introduction/Version_R1/tikz_tex/files/samplePath_100_ED_lower.txt};
\addplot[white, very thick] coordinates { (0,101) (8.3,101) };
\addplot[dashed, very thick] coordinates { (0,101) (8.3,101) };
\end{axis}
\end{tikzpicture}
&
\begin{tikzpicture}[scale = 0.58]
\begin{axis}[
	xmin = 0, 
	xmax = 8,
	ymin = 0,
	ymax = 180,
	axis line style={->},
	axis lines = left,
	x label style = {at = {(0.75,-0.05)}},
	xlabel = $\to t$,
	xscale=0.8,
	yscale=0.8]
\small	
\addplot[draw = black, fill = col1, fill opacity = 0.2] file {Introduction/Version_R1/tikz_tex/files/samplePath_100_QED_upper.txt};
\addplot[draw = black, fill = col4, fill opacity = 0.2] file {Introduction/Version_R1/tikz_tex/files/samplePath_100_QED_lower.txt};
\addplot[white, very thick] coordinates { (0,105) (8.3,105) };
\addplot[dashed, very thick] coordinates { (0,105) (8.3,105) };
\end{axis}
\end{tikzpicture}
&
\begin{tikzpicture}[scale = 0.58]
\begin{axis}[
	xmin = 0, 
	xmax = 8,
	ymin = 0,
	ymax = 180,
	axis line style={->},
	axis lines = left,
	x label style = {at = {(0.75,-0.05)}},
	xlabel = $\to t$,
	xscale=0.8,
	yscale=0.8]
\small
\addplot[draw = black, fill = col4, fill opacity = 0.2] file {Introduction/Version_R1/tikz_tex/files/samplePath_100_QD_lower.txt};
\addplot[white, very thick] coordinates { (0,150) (8.3,150) };
\addplot[dashed, very thick] coordinates { (0,150) (8.3,150) };
\end{axis}
\end{tikzpicture}\\
 &  \ \ $\sl^{(1)}$ & \ \ $\sl^{(2)}$ & \ \ $\sl^{(3)}$  
 \end{tabular}
\caption{Sample paths of the $M/M/s$ queue with $\lambda = 10,\,50$ and $100$ and $s$ set according to the three scaling rules in \eqref{eq:intro_three_scaling_rules}.}
\label{fig:sample_paths_lambda100}
\end{figure}

Next, we conduct a similar sample path experiment for increasing values of $\lambda$.
Since $s > \lambda$ is required for stability, the value of $s$ needs to be adjusted accordingly.
We propose three scaling rules:
\begin{equation}
\label{eq:intro_three_scaling_rules}
s^{(1)}_\lambda = \left[ \lambda + \beta \right ], \qquad
s^{(2)}_\lambda = \left[ \lambda + \beta\sqrt{\lambda} \right], \qquad
s^{(3)}_\lambda = \left[ \lambda + \beta\,\lambda \right],
\end{equation}
for some $\beta>0$, where $[\cdot]$ denotes the rounding operator. 
Note that these three rules differ in terms of overcapacity $s-\lambda$.
Figure \ref{fig:sample_paths_lambda100} depicts typical sample paths of the queue length process for increasing values of $\lambda$ for the three scaling rules with $\beta = 0.5$.

Observe that for all scaling rules, the stochastic fluctuations of the queue length processes relative to $\lambda$ decrease with the system size. 
Moreover, the paths in Figure \ref{fig:sample_paths_lambda100} appear to become smoother with increasing $\l$.
Of course, the actual sample path always consists of upward and downward jumps of size 1, but we will show how proper centering and scaling of the queue length process indeed gives rise to a \textit{diffusion process} in the limit as $\lambda\to\infty$.
Although the difference in performance of the three regimes is not yet evident for relatively small $\lambda$, clear distinctive behavior occurs for large $\lambda$. 

Under $\sl^{(1)}$, the majority of customers is delayed and server idle time is low, since $\rho = (1+\beta/\lambda)^{-1} \to 1$ as $\lambda \to \infty$.
Systems dimensioned according to this rule value server efficiency over customer satisfaction and  therefore this regime is in the literature also known as the \textit{efficiency-driven} (ED) regime \cite{Zeltyn2005}.

In contrast, the third scaling rule $s^{(3)}_\l$ yields a constant utilization level $\rho = 1/(1+\beta)$, which stays away from 1, even for large $\lambda$. 
Queues operating in this regime exhibit significant server idle times. 
Moreover, for the particular realization of the queueing processes for $\lambda = 50$ and $\lambda=100$ none of the customers waits.
This customer-centered regime is known as the \textit{quality-driven} (QD) regime \cite{Zeltyn2005}.

The scaling rule $s^{(2)}_\l$ is in some ways a combination of the other two regimes. 
First, we have $\rho = (1 +\beta/\sqrt{\lambda})^{-1} \to 1$ as $\lambda \to \infty$, which indicates efficient usage of resources as the system grows. 
The sample paths, however, indicate that only a fraction of the customers is delayed, and only small queues arise, which suggest good QoS. 
This regime is therefore called \textit{quality-and-efficiency driven} (QED) regime.
Since this scaling regime and the related \textit{square-root staffing rule}
\begin{equation}
\label{eq:square_root_staffing rule}
s_\l = \lambda + \beta\sqrt{\lambda}
\end{equation}
strikes the right balance between the two profound objectives of capacity allocation in service systems, we discuss in the next section the mathematical foundations of the QED regime and quantify the favorable properties revealed by Figure \ref{fig:sample_paths_lambda100}.

\section{The QED regime: two canonical examples} 
\label{sec:intro_QED_regime}

We saw in Figure \ref{fig:waiting_time_pooling} the advantageous effect of resource pooling and economies-of-scale in many-server systems.
In this section, we will explain how this is related to the Central Limit Theorem (CLT).

\begin{theorem}[Central Limit Theorem, e.g. {\cite[Thm.~27.1]{Billingsley1995}}]
Consider a sequence $X_1,X_2,\ldots,X_n$ of independent and identically distributed random variables having mean $\mu$ and positive variance $\sigma^2$. 
Then,
\[
\frac{\sum_{i=1}^n X_i - n\mu }{\sqrt{n}\sigma} \Rightarrowd \mathcal{N}(0,1), \qquad \text{for }n\to\infty.
\]
where $\Rightarrowd$ denotes convergence in distribution and $\mathcal{N}(0,1)$ is a random variable with standard normal distribution. 
\end{theorem}


We shall now apply the CLT to the delay probability in the $M/M/s$ queue. 
Striking the proper balance between queueing delay and server efficiency asymptotically, i.e.~balancing the green and red areas in Figure \ref{fig:sample_paths_lambda100}, in mathematical terms boils down to choosing a service level $s_\lambda$ such that both the delay probability $\P(Q^{(s_\lambda)} \geq s_\lambda)$ and $\P(Q^{(s_\lambda)} < s_\lambda)$ remain strictly smaller than 1 as $\lambda\to\infty$. 
In other words, one would like to see that $\P(Q^{(s_\lambda)} \geq s_\lambda)$ converges to a non-degenerate limit $\alpha \in (0,1)$ as $\lambda\to\infty$.

To get a feel for the natural scale of the queue, we first examine the situation with unlimited capacity. 
More precisely, let $Q^{(\infty)}$ be the number of customers in a steady-state $M/GI/\infty$ queue with mean service requirement $\E[B]=1$. 
Notice that in this infinite-server setting, $Q^{(\infty)}$ also represents the steady-state number of busy servers.
It is well known that $Q^{(\infty)}$ follows a Poisson distribution with mean equal to the expected workload, in our case $R = \l$. 
Moreover, if we assume that $\l$ is integer, then a Poisson random variable with rate $\l$ can be viewed as the sum of $\l$ i.i.d. Poisson random variables with rate 1. 
In other words, $Q^{(\infty)} = \sum_{i=1}^\l P_i$, where the $P_i$'s, $i=1,2,\ldots,n$, have Poisson distribution with unit mean and variance, and are mutually independent. 
\noindent
The CLT thus gives
\begin{equation}
\label{eq:infinite_server_tail}
\P(Q^{(\infty)} \geq x_\lambda ) 
= \P\left(\frac{Q^{(\infty)} -\lambda }{\sqrt{\lambda}} \geq \frac{ x_\lambda - \lambda}{\sqrt{\lambda}} \right)
\approx 1-\F\left( \frac{x_\lambda-\lambda}{\sqrt{\lambda}} \right),
\end{equation}
where $\Phi$ denotes the cumulative distribution function of the standard normal distribution for large $\l$.
Hence, the probability in \eqref{eq:infinite_server_tail} converges to a constant value away from both 0 and 1 if and only if $(x_\lambda - \lambda)/\sqrt{\lambda} \to x \in \mathbb{R}$, or equivalently $x_\lambda = \lambda + x \sqrt{\lambda} + o(\sqrt{\lambda})$, as $\lambda\to\infty$. 
Here, the relation $u(\l) = o(v(l))$ implies that $u(\l)/v(\l)\to 0$ as  $\lambda\to\infty$. 
Equation \eqref{eq:infinite_server_tail} also shows that the leading order of the random variable describing the queue length is $\l$, while the stochastic fluctuations are of order $\sqrt{\lambda}$.

If we now pretend, for a moment, that the infinite-server queue serves as a good approximation for the many-server queue with $s_\lambda$ servers, then \eqref{eq:infinite_server_tail} says that the steady-state probability of delay for $\sl = \l +\beta\sqrt\l$ obeys the Gaussian approximation
\begin{equation}
\label{eq:infinite_server_approx_delay}
\P(W^{(s_\lambda)}>0) = \P(Q^{(s_\lambda)} \geq s_\lambda ) \approx 1-\F(\beta),
\end{equation}
where $\F$ denotes the cumulative distribution function (cdf) of the standard normal distribution.
Of course, the infinite-server system ignores the one thing that makes a queueing system unique, namely that a queue is formed when all servers are busy. 
During these periods of congestion, customers will depart from a system with a finite number of servers at a slower pace than in its infinite-server counterpart.
So the approximation in \eqref{eq:infinite_server_approx_delay} is likely to overestimate the actual delay probability, and a more careful investigation of the queue length process in many-server settings is needed. Nevertheless, the infinite-server heuristic reveals that in a well-managed system, i.e.~queues are of acceptable length, the size at which the system operates is of the order $\l$, with fluctuations of order $\sqrt{\lambda}$. 
We shall now demonstrate through two canonical examples how these guessed natural scalings can be turned into mathematically rigorous statements.
Both examples which will play a key role in this thesis.

\subsection{The $M/M/s$ queue}
\label{sec:intro_MMsqueue}

\textbf{Converging delay probability.} 
Let $Q^{(s)}$ denote the steady-state number of customers in an $M/M/s$ queue with arrival rate $\lambda$ and mean service requirement 1, of which the probability distribution is given in \eqref{eq:MMs_stationary_distribution}. 
Halfin \& Whitt \cite{Halfin1981} showed that, just as the tail probability in the infinite-server setting, the delay probability in the $M/M/s$ queue converges under scaling \eqref{eq:square_root_staffing rule} to a value between 0 and 1.
Moreover, they showed that this is in fact the only scaling regime in which such a non-degenerate limit exists and identified its value. 
Let $\f$ denote the probability density function (pdf) of the standard normal distribution.

\begin{proposition}[{\cite[Prop.~2.1]{Halfin1981}}]
\label{prop:HalfinWhitt_delay_probability}
The probability of delay in the $M/M/s_\lambda$ queue has the non-degenerate limit
\begin{equation}
\lim_{\lambda\to\infty} \P( W^{(s_\lambda)} > 0 ) = \alpha \in (0,1)
\end{equation}
if and only if 
\begin{equation}
\label{eq:HalfinWhitt_scaling}
\lim_{\lambda\to\infty} (1-\rho_{s_\lambda}) \sqrt{s_\lambda} \to \beta, \quad \beta > 0,
\end{equation}
where $\alpha$ is given by
\begin{equation}\label{eq:HW_delay_prob}
\alpha = \left( 1+ \frac{\beta\,\F(\beta)}{\f(\beta)} \right)^{-1} =: g(\beta).
\end{equation}
\end{proposition} 
\begin{proof}
We first prove the sufficiency condition.
Rewrite \eqref{eq:MMs_wait} as
\begin{equation}
\label{eq:proof_HW_0}
\P( Q^{(s_\lambda)} \geq s_\lambda ) 
= \left( 1 + (1-\rho_{\sl})\frac{ \P(\Pois(\lambda) < \sl) }{\P(\Pois (\l) = \sl)}\right) ^{-1}.
\end{equation}
Similar to \eqref{eq:infinite_server_tail} we find
\begin{align}
\P(\Pois(\l) < \sl) 
&= \P\left(\frac{\Pois(\l)-\l}{\sqrt{\l}} < \frac{\sl-\l}{\sqrt{\l}}\right) \nonumber\\
&= \P\left(\frac{\Pois(\l)-\l}{\sqrt{\l}} < (1-\rho_{\sl})\,\frac{\sl}{\sqrt\l}\right)\nonumber\\
&= \P\left(\frac{\Pois(\l)-\l}{\sqrt{\l}} < (1-\rho_{\sl})\,\sqrt{\sl}\left(1+o(1)\right) \right) \to \F(\beta),
\label{eq:proof_HW_1}
\end{align}
for $\l\to\infty$.
Using Stirling's approximation, we get
\begin{align*}
\P(\Pois(\l)=s) &= {\rm e}^{-\l}\frac{\l^{\sl}}{\sl!}
\sim {\rm e}^{-\l} \l^{\sl}\cdot \frac{1}{\sqrt{2\pi\,\sl}} \left(\frac{\rm e}{\sl}\right)^{\sl} = \frac{1}{\sqrt{2\pi\sl}}\,\ee^{\sl-\l - \sl{\rm ln}(\rho_{\sl})}.
\end{align*}
Since ${\rm ln}(\rho_{\sl}) = -(1-\rho_{\sl}) - \tfrac{1}{2}(1-\rho_{\sl})^2 + o((1-\rho_{\sl})^2)$ we find that
\begin{equation}
\label{eq:proof_HW_2}
\frac{ \P(\Pois(\l) = s) }{ 1-\rho_{\sl} } 
= \frac{1}{(1-\rho_{\sl})\sqrt{\sl}} \, \frac{\ee^{ -\tfrac{1}{2}(1-\rho_{\sl})^2\sl + o\left((1-\rho_{\sl})^2\sl\right)}}{\sqrt{2\pi}}  \to \frac{1}{\beta}\, \frac{\ee^{{-}\tfrac{1}{2} \beta^2}}{\sqrt{2\pi}} = \frac{\f(\beta)}{\beta}. 
\end{equation}
Substituting \eqref{eq:proof_HW_1} and \eqref{eq:proof_HW_2} into \eqref{eq:proof_HW_0} gives \eqref{eq:HW_delay_prob}.
The necessary condition follows directly by the characterization of $\P( Q^{(s_\lambda)} \geq s_\lambda )$ as in \eqref{eq:proof_HW_0} by observing, through \eqref{eq:proof_HW_1} and \eqref{eq:proof_HW_2}, that the term
\begin{equation*}
(1-\rho_{\sl})\frac{ \P(\Pois(\lambda) < \sl) }{\P(\Pois (\l) = \sl)}
\end{equation*}
has a limiting value in $(0,\infty)$ only if $(1-\rho_{\sl})\sqrt{\sl} \to \beta$ for some $\beta>0$.  
\end{proof}

Observe that $g(\beta)$ is a strictly decreasing function on $(0,\infty)$ with $g(\beta) \to 1$ as $\beta\to 0$ and $g(\beta)\to 0$ for $\beta\to\infty$.
Thus all possible delay probabilities are achievable in the QED regime, which will prove useful for the dimensioning of systems (see Section \ref{sec:intro_dimensioning}). 

\begin{figure}
\centering
\begin{tikzpicture}[scale=0.8]
\begin{axis}[
	xmin = 0, 
	xmax = 40,
	ymin = 0,
	ymax = 1.05,
	grid = both, 
	axis line style={->},
	axis lines = left,
	xlabel = $\to \l$,
	ylabel = {$\P\left(Q^{(\sl)}\geq \sl\right)$},
	y label style = {at = {(-0.1,0.6)}},
	xscale=1,
	yscale=0.8]
\small
\addplot[very thick] table
{
0.904875	0.904875
1.86349	0.898824
2.83172	0.895881
3.80494	0.894042
4.78134	0.892747
5.76	0.89177
6.74038	0.890997
7.72211	0.890366
8.70496	0.889838
9.68873	0.889386
10.6733	0.888995
11.6586	0.888651
12.6444	0.888346
13.6308	0.888073
14.6177	0.887826
15.605	0.887602
16.5927	0.887397
17.5807	0.887209
18.5691	0.887035
19.5578	0.886874
20.5467	0.886724
21.5359	0.886584
22.5254	0.886453
23.5151	0.88633
24.505	0.886214
25.4951	0.886105
26.4854	0.886002
27.4758	0.885904
28.4665	0.885811
29.4573	0.885722
30.4482	0.885638
31.4393	0.885558
32.4305	0.885481
33.4219	0.885407
34.4134	0.885337
35.405	0.885269
36.3967	0.885205
37.3885	0.885142
38.3805	0.885082
39.3725	0.885024
40.3647	0.884969
41.3569	0.884915
42.3492	0.884863
43.3417	0.884813
44.3342	0.884764
45.3267	0.884717
46.3194	0.884671
47.3122	0.884627
48.305	0.884584
49.2979	0.884543
};
\addplot[dashed] coordinates{	(0,0.880287) (40,0.880287) } ;
\addplot[very thick] table {
0.609612	0.609612
1.40693	0.581007
2.25	0.567757
3.11722	0.55969
4.	0.554113
4.89389	0.549958
5.79623	0.546707
6.70527	0.544073
7.6198	0.541882
8.53893	0.540021
9.46198	0.538416
10.3884	0.537012
11.3179	0.53577
12.25	0.534662
13.1845	0.533664
14.1211	0.532761
15.0597	0.531937
16.	0.531181
16.942	0.530485
17.8854	0.529842
18.8303	0.529244
19.7765	0.528687
20.7238	0.528166
21.6723	0.527678
22.6219	0.527219
23.5724	0.526786
24.5239	0.526378
25.4763	0.525991
26.4295	0.525624
27.3835	0.525275
28.3383	0.524944
29.2938	0.524627
30.25	0.524325
31.2068	0.524037
32.1643	0.52376
33.1224	0.523496
34.0811	0.523242
35.0403	0.522997
36.	0.522763
36.9603	0.522537
37.921	0.522319
38.8822	0.522109
39.8439	0.521906
40.806	0.52171
41.7686	0.521521
42.7315	0.521338
43.6949	0.521161
44.6586	0.520989
45.6228	0.520822
46.5873	0.520661
};
\addplot[dashed] coordinates{ (0,0.504539) (40,0.504539) } ;
\addplot[very thick] table {
0.381966	0.381966
1.	0.333333
1.69722	0.312006
2.43845	0.29945
3.20871	0.290969
4.	0.284761
4.80742	0.27997
5.62772	0.276131
6.45862	0.272967
7.29844	0.270303
8.1459	0.268019
9.	0.266035
9.85995	0.264289
10.7251	0.262738
11.5949	0.261349
12.4689	0.260095
13.3467	0.258956
14.228	0.257915
15.1125	0.256959
16.	0.256078
16.8902	0.255261
17.783	0.254502
18.6782	0.253794
19.5756	0.253132
20.4751	0.25251
21.3765	0.251926
22.2798	0.251375
23.1849	0.250854
24.0917	0.250361
25.	0.249893
25.9098	0.249449
26.8211	0.249026
27.7337	0.248622
28.6477	0.248237
29.5628	0.247869
30.4792	0.247516
31.3967	0.247178
32.3153	0.246854
33.235	0.246543
34.1557	0.246243
35.0774	0.245955
36.	0.245677
36.9235	0.245409
37.8479	0.24515
38.7732	0.244901
39.6993	0.244659
40.6261	0.244426
41.5538	0.2442
42.4822	0.243981
43.4113	0.243768
};
\addplot[dashed] coordinates{ (0,0.223361) (40,0.223361) };

\end{axis}

\node[right] at (6.9,3.85) { \small $g(0.1)$ };
\node[right] at (6.9,2.2) { \small $g(0.5)$ };
\node[right] at (6.9,1) { \small $g(1)$ };

\end{tikzpicture}
\caption{The delay probability $\P(Q^{(\sl)} \geq \sl)$ with $\sl = [ \l + \beta \sqrt{\l} ]$ for $\beta = 0.1,\ 0.5,$ and 1 as a function of $\l$.}
\label{fig:delay_probs_HW_MMs}
\end{figure}

Although Proposition \ref{prop:HalfinWhitt_delay_probability} is an asymptotic result for $\l\to\infty$, Figure \ref{fig:delay_probs_HW_MMs} shows that $g(\beta)$ can serve as an accurate approximation for the delay probability for relatively small $\l$. 

From Proposition \ref{prop:HalfinWhitt_delay_probability}, it also follows that under \eqref{eq:HalfinWhitt_scaling}, 
\begin{equation}
\label{eq:halfinwhitt_wait}
\sqrt{\sl}\,\E[W^{(\sl)}] = \frac{\P(W^{(\sl)}>0)}{(1-\rho_\sl)\sqrt{\sl}} \to \frac{g(\beta)}{\beta} =: h(\beta), \qquad \text{ as }\l\to\infty,
\end{equation}
where we have used the characterization of $\E[W^{(\sl)}]$ in \eqref{eq:MMs_wait2}.
This implies that in the QED regime, the expected waiting time vanishes at rate $1/\sqrt{\sl}$ as $\l\to\infty$. 
By Little's law this implies that the expected queue length is $O(\sqrt{\sl})$. 
By the relation $u(\l) = O(v(\l))$ we mean that $\limsup_{\l\to\infty} u(\l)/v(\l)< \infty$. 
%
%
While these are all steady-state results, similar statements can be made for the entire queue-length process, as shown next.

The theoretical results of the QED regime we presented here are based on steady-state queueing analysis. 
But at the heart of the QED theory lies a much deeper result in which the entire queue-length process, over all points in time, converges to some other limiting process.
\\*

\noindent\textbf{Process-level convergence.}
Obtaining rigorous statements about stochastic-process limits poses considerable mathematical challenges. 
Rather than presenting the deep technical details of the convergence results, we give a heuristic explanation of how the limiting process arises and what it should look like.

The queue-length process $Q^{(\sl)}$ in Figure \ref{fig:sample_paths_lambda100} with scaling rule $\sl = [\l + \beta \sqrt{\l}]$ appears to concentrate around the level $\sl$. 
As argued before, the stochastic fluctuations are of order $\sqrt{\l}$, or equivalently $\sqrt{\sl}$.
For that reason, we consider the centered and scaled process
\begin{equation}
\label{eq:intro_scaled_queue_length_process}
X^{(\sl)}(t) := \frac{ Q^{(\sl)}(t) - \sl}{\sqrt{\sl}}, \qquad \text{ for\ all } t\geq 0,
\end{equation}
and ask what happens to this process as $\l\to\infty$. 
First, we consider the expected drift conditioned on $X^{(\sl)}(t) = x$. 
When $x> 0$, this corresponds to a state in which $Q^{(\sl)}>\sl$ and hence all servers are occupied. 
Therefore, the expected rate at which customers leave the system is $\sl$, while the arrival rate remains $\l$, so that the expected drift of $X^{(\sl)}(t)$ in $x>0$ satisfies
\[
\frac{\l - \sl}{\sqrt{\sl}} \to -\beta, \qquad \text{as }\l\to\infty,
\]
under scaling $\sqrt{\sl}(1-\rho_{\sl})\to \beta$ in \eqref{eq:HalfinWhitt_scaling}.
When $x\leq 0$, only $\sl + x\sqrt{\sl}$ servers are working, so that the net drift is
\[
\frac{\l - (\sl + x\sqrt{\sl} )}{\sqrt{\sl}} \to -\beta-x, \qquad \text{as }\l\to\infty.
\]
Now, imagine what happens to the sample paths of $\{X^{(\sl)}(t)\}_{t\geq 0}$ as we increase $\l$. 
Within a fixed time interval, larger $\l$ and $\sl$ will trigger more and more events, both arrivals and departures.
Also, the jump size at each event epoch decreases as $1/\sqrt{\sl}$ as a consequence of the scaling in \eqref{eq:intro_scaled_queue_length_process}.
Hence, there will be more events, each with a smaller impact, and in the limit as $\l\to\infty$, there will be infinitely many events of infinitesimally small impact. 
This heuristic explanation suggests that the process $X^{(\sl)}(t)$ converges to a stochastic-process limit, which is continuous, and has infinitesimal drift ${-}\beta$ above zero and ${-}\beta-x$ below zero. 
Figure \ref{fig:sample_paths_diffusion} visualizes the appearance of the suggested scaling limit as $\l$ and $\sl$ increase. 

\begin{figure}
\centering
\begin{subfigure}{0.48\textwidth}
\centering
\begin{tikzpicture}[scale = 0.7]
\begin{axis}[
	xmin = 0, 
	xmax = 5,
	ymin = -2.5,
	ymax = 2.5,
	axis line style={->},
	axis x line=middle,
	axis y line=left,
	ylabel near ticks,
	xlabel={},
 	ylabel={$X^{(\sl)}(t)$},
	xscale=1,
	yscale=1]	
\addplot[] file {Introduction/Version_R1/tikz_tex/files/sample_path_R5.txt};
\end{axis}
\end{tikzpicture}
\caption{$\l=5$}
\end{subfigure}
\begin{subfigure}{0.48\textwidth}
\centering
\begin{tikzpicture}[scale = 0.7]
\begin{axis}[
	xmin = 0, 
	xmax = 5,
	ymin = -2.5,
	ymax = 2.5,
	axis line style={->},
	axis x line=middle,
	axis y line=left,
	ylabel near ticks,
	xlabel={},
 	ylabel={$X^{(\sl)}(t)$},
	xscale=1,
	yscale=1]	
\addplot[] file {Introduction/Version_R1/tikz_tex/files/sample_path_R50.txt};
\end{axis}
\end{tikzpicture}
\caption{$\l = 50$}
\end{subfigure}
\begin{subfigure}{0.48\textwidth}
\centering
\begin{tikzpicture}[scale = 0.7]
\centering
\begin{axis}[
	xmin = 0, 
	xmax = 5,
	ymin = -2.5,
	ymax = 2.5,
	axis line style={->},
	axis x line=middle,
	axis y line=left,
	ylabel near ticks,
	xlabel={},
 	ylabel={$X^{(\sl)}(t)$},
 	xlabel style={right},
	xscale=1,
	yscale=1]	
\addplot[] file {Introduction/Version_R1/tikz_tex/files/sample_path_R100.txt};
\end{axis}
\end{tikzpicture}
\caption{$\l=100$}
\end{subfigure}
\begin{subfigure}{0.48\textwidth}
\centering
\begin{tikzpicture}[scale = 0.7]
\centering
\begin{axis}[
	xmin = 0, 
	xmax = 5,
	ymin = -2.5,
	ymax = 2.5,
	axis line style={->},
	axis x line=middle,
	axis y line=left,
	ylabel near ticks,
	xlabel={},
 	ylabel={$X^{(\sl)}(t)$},
 	xlabel style={right},
	xscale=1,
	yscale=1]	
\addplot[] file {Introduction/Version_R1/tikz_tex/files/sample_path_R500.txt};
\end{axis}
\end{tikzpicture}
\caption{$\l=500$}
\end{subfigure}
\caption{Sample paths of the normalized queue length process $X^{(\sl)}(t)$ with $\l = 5$, $\l=5$ and $\l=500$ and $\sl = [\l+0.5\sqrt{\l}]$.}
\label{fig:sample_paths_diffusion}
\end{figure}

The following theorem by Halfin \& Whitt \cite{Halfin1981} characterizes this scaling limit more formally.

\begin{theorem}
\label{thm:Halfin_Whitt_diffusion}
Let $X^{(\sl)}(0)\, \Rightarrowd X(0) \in \mathbb{R}$ and $\sqrt{\sl}(1-\rho_{\sl})\to\beta$. Then for all $t\geq 0$, 
\[
X^{(\sl)}(t) \Rightarrowd X(t),\qquad \text{ as }\l\to\infty,
\]
where $X(t)$ is the diffusion process with infinitesimal drift $m(x)$ given by
\[
m(x) = \left\{
\begin{array}{ll}
-\beta, & \text{if }x> 0,\\
-\beta-x, & \text{if } x \leq 0
\end{array}\right.
\]
and infinitesimal variance $\sigma^2(x) = 2$. 
\end{theorem}
The limiting diffusion process $\{X(t)\}_{t\geq 0}$ in Theorem \ref{thm:Halfin_Whitt_diffusion} is a combination of a negative-drift Brownian motion in the upper half plane and an Ornstein-Uhlenbeck (OU) process in the lower half plane. 
We refer to this hybrid diffusion process as the Halfin-Whitt diffusion.
Much is known for such diffusion processes with piecewise linear drift coefficient, see \cite{Leeuwaarden2012,Fralix2014}.
Its stationary distribution can for instance be derived, see e.g. \cite{Browne1995}.
\begin{proposition}
\label{thm:intro_HW_stationary_distribution}
Let $X(t) \Rightarrowd X(\infty)$ as $t\to\infty$ for a random variable $X(\infty)$ and $(1-\rho_{\sl})\sqrt{\sl}\to \beta$ for $\l\to\infty$. 
Then 
\begin{align}
\P(X(\infty) > 0 ) &= g(\beta),\\
\P(X(\infty) \geq x | X(\infty) > 0) &= {\rm e}^{-\beta x} ,\quad \text{for }x>0,\\
\P(X(\infty) \leq x | X(\infty) \leq 0 ) &= \frac{\F(\beta+x)}{\F(\beta)},\quad \text{for }x\leq 0.
\end{align}
\end{proposition}
\noindent
This result shows that as the system grows large, the $Q^{(\sl)}(t)$ concentrates around $\sl$, and the fluctuations are of order $\sqrt{\sl}$.
Moreover, Proposition \ref{thm:intro_HW_stationary_distribution} iterates the limiting values for the delay probability and scaled expected delay. Namely,
\[ \P\big(W^{({\sl})} > 0 \big) \rightarrow \P( X(\infty) > 0 ) = g(\beta)\]
and 
\[ \sqrt{\sl}\E[W^{(\sl)}] \approx \frac{\E[ Q^{(\sl)}]}{\sqrt{\sl}} \rightarrow \E[X(\infty)] = \int_0^\infty g(\beta){\rm e}^{-\beta x} \dd x = \frac{g(\beta)}{\beta},\]
For obvious reasons, the QED regime is also referred to as the Halfin-Whitt regime, and both these names are used interchangeably in this thesis.

\subsection{The $M/D/s$ queue}
\label{sec:intro_discrete_model}
We next consider a many-server queue with deterministic service requirements equal to one, a Poisson arrival process of rate $\l$ and $\sl$ servers.
We let $Q^{(\sl)}(t)$ be the process describing the number of customers in the system and only examine the process at discrete time epochs $t=0,1,2,\ldots$. 
In our analysis, we focus on the queue length process $Z^{(\sl)}(t) := (Q^{(\sl)}(t) - \sl)^+$.  

Since we discretize time, the number of new arrivals per time period is given by the sequence of i.i.d.~random variables $\{A_k\}_{k\geq 1}$, which has a Poisson distribution with mean $\l$.
At the start of the $k^{\rm th}$ period, $Z^{(\sl)}(k)$ customers are waiting. 
Because the service time of a customer is equal to the period length, all $\min\{Q^{(\sl)}(k),\sl\}$ customers in service at the beginning of the period will have left the system by time $t=k+1$. 
This implies that $\min\{Z^{(\sl)}(k),\sl\}$ of the waiting customers are taken into service during period $k$, but could not possibly have departed before its end, due to the deterministic service times. 
If $Z^{(\sl)}(k)<\sl$, then additionally $\min\{ A_k , \sl-Z^{(\sl)}(k) \}$ of the new arrivals are taken into service. 
This yields a total of $A_k$ arrivals, and $\min\{Z^{(\sl)}(k)+A_k,\sl\}$ departures from the queue during period $k$, which gives the Lindley type recursion \cite{Lindley1952}, with $Z^{(\sl)}(0) = 0$,
\begin{equation}
\label{eq:discrete_recursion}
Z^{(\sl)}(k+1) = Z^{(\sl)}(k) + A_k - \min\{Z^{(\sl)}(k)+A_k,\sl\} = \max\{ 0,Z^{(\sl)}(k) + A_k - \sl \}.
\end{equation}
The queue length process thus gives rise to a random walk with i.i.d.~steps of size
$(A^{(\sl)}-\sl)$, with a reflecting barrier at zero. We can iterate the recursion in \eqref{eq:discrete_recursion} to find
\begin{align}
Z^{(\sl)}(k+1) &= \max\left\{ 0 , Z^{(\sl)}(k) + A_k-\sl \right\} \nonumber\\
&= \max\left\{ 0 , \max\{ 0 , Z^{(\sl)}(k-1) + (A_{k-1}-\sl)\} + (A_k-\sl)\} \right\}\nonumber \\
&= \max\left\{ 0 , (A_k-\sl) , Z^{(\sl)}(k-1) + (A_k-\sl) + (A_{k-1}-\sl)\right\}\nonumber \\
&= \max_{0\leq j\leq k} \Big\{ \sum_{i=1}^j (A_{k-i}-\sl)\Big\} 
\equalD  \max_{0\leq j\leq k} \Big\{ \sum_{i=1}^j (A_i-\sl) \Big\},
\label{eq:max_randomwalk}
\end{align}
where the last equality in distribution holds due to the duality principle for random walks, see e.g.~\cite[Sec.~7.1]{Ross1996}.
For stability, the expected step size satisfies $\E[A_k - \sl] = \l-\sl < 0$.
We use the shorthand notation for the partial sum $S_k := \sum_{i=1}^k (A_i-\sl)$.
Let $Z^{(\sl)}(\infty):= \lim_{k\to\infty} Z^{(\sl)}(k)$ denote the stationary queue length in this $M/D/s$ queue, which can be shown to exist under our assumptions. 
The probability generating function (pgf) of $Z^{(\sl)}(\infty)$ can then be expressed in terms of the pgf of the positive parts of the partial sum:
\begin{equation}
\label{eq:Spitzers_identity}
\E[ w^{Z^{(\sl)}(\infty)} ] 
= \exp\Big\{  - \sum_{k=1}^\infty \frac{1}{k}\, (1- \E[w^{S_k^+}]) \Big\},\qquad |w|\leq 1.
\end{equation}
From \eqref{eq:Spitzers_identity}, which is a special case of Spitzer's identity~\cite{Spitzer1964}, we obtain for the mean queue length and empty-queue probability the expressions
\begin{align}
\E[Z^{(\sl)}(\infty)] &= \sum_{k=1}^\infty \frac{1}{k}\, \E[ S_k^+ ],\nonumber\\
\P(Z^{(\sl)}(\infty) = 0 ) &= \exp\Big\{ -\sum_{k=1}^\infty \frac{1}{k}\, \P( S_k^+ > 0 ) \Big\}.
\label{eq:spitzer_expressions}
\end{align}
Although explicit, the expressions in \eqref{eq:spitzer_expressions} reveal little of the structure of the queue length process. 
Hence, we again turn to asymptotics. \\

\noindent\textbf{Gaussian random walk.}
\label{sec:intro_gaussian_random_walk}
We take another look at the identity in \eqref{eq:max_randomwalk}, and ask ourselves what happens if $\l$ grows large. 
Since $\E[A_k-\sl] = \l-\sl = -\beta\sqrt{\l} + o(\sqrt{\l})$ under the QED scaling \eqref{eq:square_root_staffing rule}, it makes sense to consider the scaled queue length process $X^{(\sl)}(k) := Z^{(\sl)}(k)/\sqrt{\l}$ for all $k\geq 0$, with scaled steps $Y_k^{(\sl)} := (A_k-\sl)/\sqrt{\l}$. 
Dividing both sides of \eqref{eq:max_randomwalk} by $\sqrt{\l}$ then gives
\begin{equation}
X^{(\sl)}(k+1) = \max_{0\leq j\leq k} \Big\{ \sum_{i=1}^j Y^{(\sl)}_i \Big\}.
\end{equation}
Observe that $A_k \equalD \Pois(\l)$ with $\Pois(\l)$ a random variable with mean $\l$. 
Hence by the CLT 
\begin{equation*}
Y^{(\sl)}_k = \frac{ A_k - \sl }{\sqrt\l} = \frac{A_k-\l}{\sqrt\l} - \beta \ \Rightarrowd \ Y_k \equalD \N(-\beta,1),
\end{equation*}
for $\l\to\infty$, where $\N(-\beta,1)$ denotes a normally distributed random variable with mean $-\beta$ and standard deviation 1. 
So we expect the scaled queue length process to converge in distribution to a reflected random walk with normally distributed increments, i.e. a reflected \textit{Gaussian random walk}. 
Indeed, it is easily verified that \cite{Janssen2008a},
\begin{equation}
X^{(\sl)}(k)\ \Rightarrowd \  M_\beta(k) := \max_{0\leq j\leq k} \Big\{\sum_{i=1}^j Y_j \Big\}, \qquad \l\to\infty.
\end{equation}
Let $M_\beta:= \lim_{k\to\infty} M_\beta(k)$ denote the all-time maximum of a Gaussian random walk.
It can be shown that $M_\beta$ almost surely exists and that 
\[
X^{(\sl)}(\infty) := \lim_{k\to\infty} X^{(\sl)}(k) \Rightarrowd M_\beta,
\]
for instance by \cite[Prop.~19.2]{Spitzer1964} and \cite[Thm.~X6.1]{Asmussen2003}.
The following theorem can be proved using a similar approach as in \cite{Jelenkovic2004}.
(We prove this result in a more general setting in Chapter 3.)

\begin{theorem}
Let $X^{(\sl)}(\infty)$ be the scaled queue length in steady-state. If $(1-\rho_{\sl})\sqrt{\l}\to\beta$, then as $\l\to\infty$,
\begin{enumerate}
\item[\normalfont (i)] $X^{(\sl)}(\infty) \Rightarrowd M_\beta$,
\item[\normalfont (ii)] $\P(X^{(\sl)}(\infty) = 0) \to \P(M_\beta = 0)$,
\item[\normalfont (iii)] $\E[X^{(\sl)}(\infty)^k] \to \E[M_\beta^k]$, for any $k>0$.
\end{enumerate}
\end{theorem}

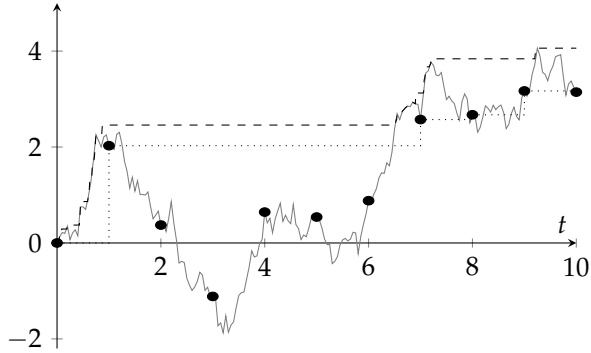
\begin{figure}
\centering
\begin{tikzpicture}[scale = 1]
\begin{axis}[
	xmin = 0, 
	xmax = 10,
	ymin = -2.2,
	ymax = 5,
	axis line style={->},
	axis x line=middle,
	axis y line=left,
	xlabel={$t$},
 	ylabel={},
	xscale=1,
	yscale=0.8]	
\addplot[gray] file {Introduction/Version_R1/tikz/Brownian_Motion_SamplePath/BM.txt};
\addplot[only marks,mark size = 2] file {Introduction/Version_R1/tikz/Brownian_Motion_SamplePath/GW.txt};
\addplot[dashed] file {Introduction/Version_R1/tikz/Brownian_Motion_SamplePath/maxBM.txt};
\addplot[dotted] file {Introduction/Version_R1/tikz/Brownian_Motion_SamplePath/maxGW.txt};
\end{axis}
\end{tikzpicture}
\caption{Brownian motion (gray) and embedded Gaussian random walk (marked) with their respective running maxima (dashed and dotted, respectively).}
\label{fig:BrownianMotion_vs_GaussianRW}
\end{figure}

The Gaussian random walk is well studied \cite{Siegmund1978,Chang1997,Janssen2006,Blanchet2006,Janssen2006} and there is an intimate connection with Brownian motion. 
The only difference, one could say, is that Brownian motion is a continuous-time process, whereas the Gaussian random walk only changes at discrete points in time.
If $\{B(t)\}_{t\geq 0}$ is a Brownian motion with drift $-\mu <0$ and infinitesimal variance $\sigma^2$ and $\{W(t)\}_{t \geq 0}$ is a random walk with $\N(-\mu,\sigma^2)$ steps and $B(0) = W(0)$, then $W$ can be regarded as the process $B$ embedded at equidistant time epochs. 
That is, $W(t) \equalD B(t)$ for all $t\in\mathbb{N}^+$. 
For the maximum of both processes this coupling implies
\begin{equation}
\max_{k\in \mathbb{N}^+} W(k) = \max_{k\in \mathbb{N}^+} B(k) \leq_{\rm st}
\max_{t\in \mathbb{R}^+} B(t),
\label{eq:max_inequality}
\end{equation}
where $\leq_{\rm st}$ denotes stochastic dominance. 
This difference in maximum is visualized in Figure \ref{fig:BrownianMotion_vs_GaussianRW}.
It is known that the all-time maximum of Brownian motion with negative drift $-\mu$ and infinitesimal variable $\sigma^2$ has an exponential distribution with mean $\sigma/2\mu$ \cite{Harrison1985}. 
Hence, \eqref{eq:max_inequality} implies that $M_\beta$ is stochastically upper bounded by an exponential random variable with mean $1/2\beta$.

Despite this easy bound, precise results for $M_\beta$ are more involved. Let $\zeta$ denote the Riemann zeta function, which is defined as, see e.g.~\cite[Eq.~25.2.1]{NIST},
\begin{equation}
\zeta(s) = \sum_{n=1}^\infty \frac{1}{n^s}.
\end{equation}

\begin{theorem}[{\cite[Thm.~1]{Chang1997} \& \cite[Thm.~2 \& 3]{Janssen2006}}]
For $0<\beta<2\sqrt{\pi}$,
\begin{align}
\P(M_\beta = 0) &= \sqrt{2}\beta\, \exp \left\{ \frac{\beta}{\sqrt{2\pi}} \sum_{l=0}^\infty 
\frac{\zeta(1/2-l)}{l!(2l+1)} \left(\frac{-\beta^2}{2}\right)^l \right\},\\
\E[M_\beta] &= \frac{1}{2\beta} + \frac{\zeta(1/2)}{\sqrt{2\pi}} + \frac{\beta}{4} 
+ \frac{\beta^2}{\sqrt{2\pi}} \sum_{l=0}^\infty 
\frac{\zeta(-1/2-l)}{l!(2l+1)(2l+2)} \left(\frac{-\beta^2}{2}\right)^l,\\
\Var M_\beta &= 
\frac{1}{4\beta^2} - \frac{1}{4} - \frac{2\,\zeta(-1/2)}{\sqrt{2\pi}}\beta - \frac{\beta^2}{24}\nonumber\\
&\qquad\qquad -
\frac{2\beta^3}{\sqrt{2\pi} } \sum_{l=0}^\infty 
\frac{\zeta(-3/2-l)}{l!(2l+1)(2l+2)(2l+3)} \Big(\frac{-\beta^2}{2}\Big)^l.
\end{align}
\end{theorem}

\subsection{Characteristics of the QED regime}
\label{sec:intro_characteristics}
Now that we have seen how the square-root staffing rule \eqref{eq:square_root_staffing rule} yields non-degenerate limiting behavior in two classical queueing models, we shall elaborate on how the QED regime gives rise to (at least) three desirable properties.
The first property relates to the efficient usage of resources, expressed as
\begin{equation}
\rho_{\sl} = \frac{\l}{\sl} = 1 - \frac{\beta}{\sqrt{\sl}} + O\big(1/\l\big), \tag{Efficiency}
\end{equation}
where we have used that $\sl = O(\l)$. 
The second distinctive property is the balance between QoS and efficiency:
\begin{equation}
\P(W^{(\sl)}>0) \to g(\beta), \qquad \text{and} \qquad \P(W^{(\sl)}>0) \to 1-\P(M_\beta=0), \tag{Balance} 
\end{equation}
as $\sl \to \infty$, in the $M/M/s$ queue and $M/D/s$ queue, respectively. 
The third property relates to good QoS:
\begin{equation}
\E[W^{(\sl)}] = \frac{h(\beta)}{\sqrt{\sl}} + o(1/\sqrt{\sl}) \qquad \text{and} \qquad \E[W^{(\sl)}] = \frac{\E[M_\beta]}{\sqrt{\sl}} + O(1/\sqrt{\sl}), \tag{QoS}
\end{equation}
in the $M/M/s$ queue and $M/D/s$ queue, respectively. 
Hence the expected waiting time vanishes at rate $1/\sqrt{\sl}$. 

Both limiting functions $g(\b)$ and $1-\P(M_\beta=0)$ can take all values in $(0,1)$ by tuning the parameter $\beta$.

Since the mathematical underpinning of these properties comes from the CLT, we can expect the properties to hold for a much larger class of models. 
These models should then be members of the same universality class (to which the CLT applies). 
Let us again show this by example. 

Consider a stochastic system in which demand per period is given by some random variable $A$, with mean $\mu_A$ and variance $\sigma_A^2<\infty$.
For systems facing large demand we propose to set the capacity according to the more general rule 
\[s = \mu_A + \beta\sigma_A,\]
which consists of a minimally required part $\mu_A$ and a variability hedge $\beta\sigma_A$. 
Assume that the workload brought into the system is generated by $n$ stochastically identical and independent sources.
Each source $i$ generates $A_{i,j}$ work in the $j^{th}$ period, with $\E[A_{i,j}] = \mu$ and $\Var\,A_{i,j} = \sigma^2$.
Then the total amount of work arriving to the system during one period is $A_j^{(n)} = \sum_{i=1}^n A_{i,j}$ with mean $n\mu$ and variance $n\sigma^2$. 
Assume that the system is able to process a deterministic amount of work $s_n$ per period and denote by $U^{(n)}(j)$ the amount of work left over at the end of period $j$. 
Then, 
\begin{equation}
 U^{(n)}(j+1) = \left( U^{(n)}(j) + A^{(n)}_j - s_n \right)^+.
 \end{equation} 
Given that $s_n >  \E[A^{(n)}_1] = n\mu$, the steady-state limit $U^{(n)} := \lim_{j\to\infty} U^{(n)}(j)$ exists and satisfies
\begin{equation}
U^{(n)} \equalD  \left( U^{(n)} + A^{(n)}_1 - s_n \right)^+. 
\label{eq:bulk_service_stationary_recursion}
\end{equation}
This framework is also known as the bulk service queue or the Anick-Mitra-Sondhi model \cite{Anick1982,Janssen2005,Janssen2008}.
In this scenario, increasing the system size is done by increasing $n$, the number of input flows. 
As we have seen before, it requires a rescaling of the process $U^{(n)}$ by an increasing function $c(n)$, in order to obtain a non-degenerate scaling limit $U := \lim_{n\to\infty} U^{(n)}/c(n)$. 
(We omit the technical details needed to justify the interchange of limits.)
From \eqref{eq:bulk_service_stationary_recursion} it becomes clear that the scaled increment
\begin{equation}
\frac{A^{(n)}_j - s_n}{c(n)} = \frac{\sum_{i=1}^n A_{i,j} - n\mu}{c(n)} + \frac{n\mu - s_n}{c(n)}
\end{equation}
only admits a proper limit if $c(n)$ is of the form $c(n) = O(\sqrt{n})$, by the virtue of the CLT, and $(s_n-n\mu)/c(n) \to \beta >0$ as $n\to\infty$. 
Especially for $c(n) = \sigma\sqrt{n}$, this reveals that $U$ has a non-degenerate limit, which is equal in distribution to the maximum of a Gaussian random walk with drift ${-}\beta$ and variance 1, if 
\[
s_n = n\mu+\beta \sqrt{n}\sigma + o(\sqrt{n}).
\]
Moreover, the results on the Gaussian random walk presented in Section \ref{sec:intro_gaussian_random_walk} are applicable to this model and the key features of the QED scaling carry over to this more general setting as well. 
In conclusion, the many-sources framework shows that the QED scaling finds much wider applications than queueing models with Poisson input only.

\subsection{Related literature}
We now provide a partial overview on the literature on heavy-traffic analysis in queueing theory and the QED regime in particular.\\
\\*

\noindent\textbf{Conventional heavy-traffic.}
Before the formal introduction of the Halfin-Whitt scaling regime in 1981, see \cite{Halfin1981}, the existing literature on the asymptotic analysis of queues mostly evolved around two types of scaling regimes: single-server and infinite-server regimes.

The idea of studying a sequence of queues in which the utilization approaches 100\%, i.e.~heavy-traffic, was first laid out by Kingman in the 1960s. 
In \cite{Kingman1961,Kingman1962} he showed how in the $GI/G/1$ queue, under mild conditions on the arrival and service processes, the scaled steady-state waiting time $(1-\rho)W^{(1)}$ converges to an exponentially distributed random variable. 
The notion that heavily loaded systems admit a scaling limit that is remarkably simple compared to the otherwise intractable pre-limit queueing systems triggered a surge of  research within the field of queueing theory in the 1960s and 1970s, see  \cite{Borovkov1965,Iglehart1970,Brumelle1971,Newell1973,Kollerstrom1974,Kollerstrom1979,Whitt1974} among others. 
These works conduct their asymptotic analysis in what we now call conventional heavy-traffic. 
That is, the service times and number of servers are held fixed, while the arrival rate approaches the critical value from below. 
A noteworthy result of these efforts is the extension of Kingman's findings to the $GI/G/s$, which finds that the scaled queue length $(1-\rho)Q^{(s)}$ converges in distribution to an exponential random variable with mean $(c_a^2+c_s^2)/2$, where $c_a$ and $c_s$ denote the coefficient of variation of the interarrival and service time distribution, respectively.
We remark that this limiting result is the key ingredient to the widely applied Kingman formula
\begin{equation}\label{eq:kingman}
\E[W^{(1)}] \approx \frac{\rho}{1-\rho} \cdot \frac{c_a^2+c_s^2}{2} \cdot \E[B],
\end{equation}
which serves as an approximation to the expected waiting time in the single-server queue.
The limit \eqref{eq:kingman} reveals that in the conventional heavy-traffic regime, the expected waiting time explodes as $\rho\to 1$.
Hence, efficient usage of resources is achieved, at the expense of poor QoS.

An alternative regime that received much attention, see e.g.~\cite{Iglehart1965,Borovkov1965,Iglehart1973,Iglehart1973a,Whitt1982}, fixes the service time distribution while increasing both the arrival rate $\l$ and the number of servers to infinity simultaneously, such that the ratio $\l/s$ remains constant. 
It has been shown that the sequence of queues under this scaling start resembling the behavior of infinite-server queues as $\l$ and $s$ grow. 
That is, the probability of a customer finding a queue on arrival is negligible. 
The sample paths in Figure \ref{fig:sample_paths_lambda100} are illustrative for this regime. 
Since the utilization level $\rho$ remains strictly away from one in the limit, this setting is in the literature typically not recognized as heavy-traffic. 

As Halfin \& Whitt indicate themselves, the QED regime in which service times are held fixed, and $\l$ and $s$ tend to infinite while satisfying $(1-\rho)\sqrt{s} \to \beta$, is a hybrid between the two aforementioned regimes. 
Namely, it considers the efficiency property of the conventional heavy-traffic scaling, and the good QoS levels from  infinite-server queues.\\
\\*
\noindent
\textbf{The $G/G/s$ queue in the QED regime.}
We have demonstrated in Section \ref{sec:intro_QED_regime} how to obtain QED limits for the $M/M/s$ queue and the $M/D/s$ queue. 
When one moves beyond the exponential and deterministic assumptions, establishing QED behavior becomes mathematically more challenging. 

The heavy-traffic analysis of the $G/G/s$ queue requires fundamentally different approaches than for Markovian queues. 
Most of the research conducted on the $G/G/s$ in the Halfin-Whitt regime evolves around the characterization of the stochastic process limit of the appropriately centered and scaled queueing process in terms of diffusion processes, under various assumptions on the model primitives. 
Puhalskii \& Reiman \cite{Puhalskii2000} analyzed the multi-class queue with phase-type service times in the Halfin-Whitt regime.
Heavy-traffic limits for queues in which service time distributions are lattice-based and/or have finite support are studied by Mandelbaum \& Mom\v{c}ilovi\'c \cite{Mandelbaum2008} and Gamarnik \& Mom\v{c}ilovi\'c \cite{Gamarnik2008}. 
Approaches through measure-valued processes are taken by Kang, Kaspi \& Ramanan \cite{Kaspi2011,Kang2012,Kaspi2013}.
The most general class of distributions is considered by Reed \cite{Reed2009} and Puhalskii \& Reed \cite{Puhalskii2010}, who impose no assumption on the service time distribution except for the existence of the first moment. 
For a survey on the techniques required for the analysis of process limits of $G/G/s$ queues, we refer the reader to \cite{Pang2007} and references therein.

Considerably less is known for the corresponding steady-state distribution of the $G/G/s$ queue in the QED regime. 
Namely, under the assumption of general service time distributions, truly infinite-dimensional limits arise, since the Markovian nature of the service time and `age' process can no longer be exploited. 
Works that have been able to characterize limiting behavior for the specific service time distribution classes include Jelenkovic et al.~\cite{Jelenkovic2004}, who assume deterministic service times, and Whitt \cite{Whitt2005}, who identifies the heavy-traffic limit in the case of hyperexponentially distributed service times. 
Progress in the understanding of steady-state behavior of $G/G/s$ queues in the Halfin-Whitt regime has been facilitated by Gamarnik \& Goldberg \cite{Goldberg,Gamarnik2013a}, who perform their analysis under the mild assumption that the service time distribution has finite $(2+\e)$ moment. 
A significant advance has been made by Aghanjani \& Ramanan \cite{Aghajani2016}, who identify the limit as the steady-state distribution of infinite-dimensional Markov process, given that the service time distribution has finite $(3+\eps)$ moment, while drawing upon previous results by Kang, Kaspi \& Ramanan \cite{Kaspi2011, Kang2012,Kaspi2013}.\\
\\*
\textbf{Model extensions.}
Many extensions to the standard many-server queue can be considered. 
A feature ubiquitous to service systems involving humans is customer abandonment \cite{Gans2003,Brown2005,Zeltyn2005,Mandelbaum2013}.
The $M/M/s+M$ queue introduced by Palm \cite{Palm1957}, also known as the Erlang-A model \cite{Garnett2002,Leeuwaarden2012}, acknowledges this feature by assigning every customer an exponentially distributed \textit{patience time} upon his arrival (denoted by $+M$ in the model definition). 
If a customer has not yet started receiving service by the expiration of his patience, he leaves the system.
Note that abandonments render queues stable under any load. 
Under QED scaling, the more general $G/G/s+G$ queue has received much attention under various modeling assumptions, see e.g.~\cite{Garnett2002,Gans2003,Whitt2006,Mandelbaum2009,Zeltyn2005,Mandelbaum2012a,Kang2012,Dai2010,Reed2012,Jennings2012,Zhang2013}. 
Noteworthy findings include the vanishing abandonment probability \cite{Garnett2002} and insensitivity of the patience time distribution as long as its density at 0, i.e.~the probability of abandoning immediately upon arrival, is fixed, as the system grows large under QED scaling.
Overviews of queues with abandonment and their asymptotic counterpart are given by Zeltyn \& Mandelbaum \cite{Zeltyn2005} and Dai \& He \cite{Dai2012} and Ward \cite{Ward2012}.

Other features that have been studied in the QED regime include multiple customer classes, see e.g.~\cite{Harrison2004,Atar2014,Gurvich2008,Gurvich2009,Tezcan2010}, or heterogeneous servers \cite{Armony2005,Armony2010,Mandelbaum2012b,Stolyar2010}.
These models are all interesting in their own respect and are fairly well-understood. 
Therefore, we choose to focus in this thesis on a different set of extensions, which will be discussed in Section \ref{sec:intro_beyond}.

\section{Dimensioning}
\label{sec:intro_dimensioning}
We adopt the term \textit{dimensioning} used by Borst, Mandelbaum \& Reiman~\cite{Borst2004} to say that the capacity of a service system is adapted to the load in order to reach certain performance levels.
In \cite{Borst2004} dimensioning refers to the staffing problem in a large-scale call center and key ingredients are the square-root staffing rule in \eqref{eq:square_root_staffing rule} and the QED regime. 
We now revisit the results in \cite{Borst2004} and its follow-up works to explain this connection to the QED regime. 

\subsection{Constraint satisfaction}
\label{sec:intro_constraint}
Consider the $M/M/s$ queue with arrival rate $\l$ and service rate $\mu$. 
A classical dimensioning problem is to determine the minimum number of servers $s$ necessary to achieve a certain target level of service, say in terms of waiting time.

\begin{figure}
\centering
\begin{subfigure}{0.48\textwidth}\centering
\begin{tikzpicture}[scale = 0.72]
\small
\begin{axis}[
	xmin = 0, 
	xmax = 1,
	ymin = 5,
	ymax = 12,
	axis line style={->},
	axis lines = left,
	legend cell align = left,
	xlabel = {\small $\to \eps$},
	ylabel = {},
	yscale = 0.8,
	legend style = {at = {(1,1.2)}, anchor = north east}]
	
\addplot[very thick] file {Introduction/Version_R1/tikz/Constraint_Satisfaction/lambda5_exact.txt};
\addplot[very thick, dashed, col1] file {Introduction/Version_R1/tikz/Constraint_Satisfaction/lambda5_asymptotic.txt};
\legend{{$s^*_\l(\eps)$},$s^{\rm srs}_\l(\eps)$}
\end{axis}
\end{tikzpicture}
\caption{$\l=5$}
\end{subfigure}
\begin{subfigure}{0.48\textwidth}\centering
\begin{tikzpicture}[scale = 0.72]
\small
\begin{axis}[
	xmin = 0, 
	xmax = 1,
	ymin = 10,
	ymax = 19,
	axis line style={->},
	axis lines = left,
	legend cell align = left,
	xlabel = {\small $\to \eps$},
	ylabel = {},
	yscale = 0.8,
	legend style = {at = {(1,1.2)}, anchor = north east}]
	
\addplot[very thick] file {Introduction/Version_R1/tikz/Constraint_Satisfaction/lambda10_exact.txt};
\addplot[very thick, dashed, col1] file {Introduction/Version_R1/tikz/Constraint_Satisfaction/lambda10_asymptotic.txt};
\legend{{$s^*_\l(\eps)$},$s^{\rm srs}_\l(\eps)$}
\end{axis}
\end{tikzpicture}
\caption{$\l=10$}
\end{subfigure}

\begin{subfigure}{0.48\textwidth}\centering
\begin{tikzpicture}[scale = 0.72]
\begin{axis}[
	xmin = 0, 
	xmax = 1,
	ymin = 100,
	ymax = 125,
	axis line style={->},
	axis lines = left,
	legend cell align = left,
	xlabel = {\small $\to \eps$},
	ylabel = {},
	yscale = 0.8,
	legend style = {at = {(1,1.2)}, anchor = north east}]
	
\addplot[very thick] file {Introduction/Version_R1/tikz/Constraint_Satisfaction/lambda100_exact.txt};
\addplot[very thick, dashed, col1] file {Introduction/Version_R1/tikz/Constraint_Satisfaction/lambda100_asymptotic.txt};
\legend{{$s^*_\l(\eps)$},$s^{\rm srs}_\l(\eps)$}
\end{axis}
\end{tikzpicture}
\caption{$\l=100$}
\end{subfigure}
\begin{subfigure}{0.48\textwidth}\centering
\begin{tikzpicture}[scale = 0.72]
\begin{axis}[
	xmin = 0, 
	xmax = 1,
	ymin = 500,
	ymax = 550,
	axis line style={->},
	axis lines = left,
	legend cell align = left,
	xlabel = {\small $\to \eps$},
	ylabel = {},
	yscale = 0.8,
	legend style = {at = {(1,1.2)}, anchor = north east}]
	
\addplot[very thick] file {Introduction/Version_R1/tikz/Constraint_Satisfaction/lambda500_exact.txt};
\addplot[very thick, dashed, col1] file {Introduction/Version_R1/tikz/Constraint_Satisfaction/lambda500_asymptotic.txt};
\legend{{$s^*_\l(\eps)$},$s^{\rm srs}_\l(\eps)$}
\end{axis}
\end{tikzpicture}
\caption{$\l=500$}
\end{subfigure}
\caption{Staffing levels as a function of the delay probability targets $\eps$.}
\label{fig:MMs_staffing_levels}
\end{figure}
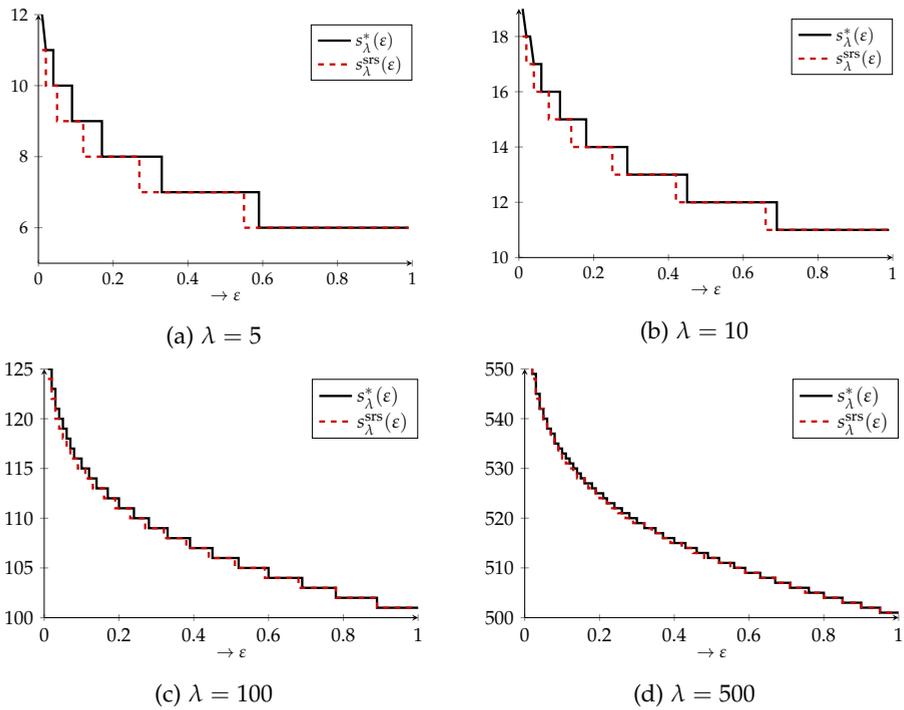

Suppose we want to determine the minimum number of servers such that the fraction of customers who are delayed in the queue is at most $\varepsilon\in(0,1)$. 
Hence we should find
\begin{equation}\label{eq:tagA}
s^{*}_\l(\eps) := \min \left\{s \geq \l\, |\, \P(W^{(s)}>0) \leq \eps \right\}. 
\end{equation}  
But alternatively, we can use the QED framework, which says that under \eqref{eq:HalfinWhitt_scaling},\ \  $\lim_{s\to\infty} \P(W^{(\sl)} > 0) = g(\beta)$ (see Proposition \ref{prop:HalfinWhitt_delay_probability}).
Then by \eqref{eq:tagA}, $s^*_\l(\eps)$ can be replaced by
\begin{equation}
s^{\rm srs}_\l(\eps) = \lceil \l + \beta^*(\eps) \sqrt{\l}\rceil,
\end{equation}
where $\beta^*(\eps)$ solves
\begin{equation}
g(\beta^*) = \eps.
\end{equation}
In Figure \ref{fig:MMs_staffing_levels} we plot the exact staffing level $s^*_\l(\eps)$ and the heuristically obtained staffing level $s^{\rm srs}_\l(\eps)$ as functions of $\eps$ for several loads $\l$. 

Observe that even for very small values of $\l$, the staffing function $s^{\rm srs}(\eps)$ coincides with the exact solution for almost all $\eps\in(0,1)$ and differs no more than by one server for all $\eps$.
Borst et al.~\cite{Borst2004} recognized this in their numerical experiments too, and Janssen, van Leeuwaarden \& Zwart \cite{Janssen2011} later confirmed this theoretically.
One can easily formulate other constraint satisfaction problems and reformulate them in the QED regime. 
For instance, constraints on the mean waiting time or the tail probability of the waiting time, e.g.~$\P(W^{(s)}>T)$, which are asymptotically approximated by $h(\beta)/\sqrt{\l}$ and $g(\beta)\ee^{-\beta \sqrt{\l} T}$, respectively.
See \cite{Borst2004} for more examples.

\subsection{Optimization}
\label{sec:intro_optimization}
\begin{figure}
\centering
\begin{subfigure}{0.48\textwidth}\centering
\begin{tikzpicture}[scale = 0.72]
\small
\begin{axis}[
	xmin = 0, 
	xmax = 5,
	axis y discontinuity = crunch,
	ymin = 5,
	ymax = 12,
	axis line style={->},
	axis lines = left,
	xlabel = {\small $\to r$},
	ylabel = {},
	yscale = 0.8,
	legend cell align = left,
	legend style = {at = {(1,1.2)}, anchor = north east}
	]
	
\addplot[very thick] file {Introduction/Version_R1/tikz/Optimization/lambda5_exact.txt};
\addplot[very thick, dashed, col1] file {Introduction/Version_R1/tikz/Optimization/lambda5_asymptotic.txt};
\legend{{$s^*_\l(r)$},$s^{\rm srs}_\l(r)$}
\end{axis}
\end{tikzpicture}
\caption{$\l=5$}
\end{subfigure}
\begin{subfigure}{0.48\textwidth}\centering
\begin{tikzpicture}[scale = 0.72]
\small
\begin{axis}[
	xmin = 0, 
	xmax = 5,
	ymin = 10,
	ymax = 19,
	axis line style={->},
	axis lines = left,
	xlabel = {\small $\to r$},
	ylabel = {},
	yscale = 0.8,
	legend cell align = left,
	axis y discontinuity=crunch,
	legend style = {at = {(1,1.2)}, anchor = north east}]
	
\addplot[very thick] file {Introduction/Version_R1/tikz/Optimization/lambda10_exact.txt};
\addplot[very thick, dashed, col1] file {Introduction/Version_R1/tikz/Optimization/lambda10_asymptotic.txt};
\legend{{$s^*_\l(r)$},$s^{\rm srs}_\l(r)$}
\end{axis}
\end{tikzpicture}
\caption{$\l=10$}
\end{subfigure}

\begin{subfigure}{0.48\textwidth}\centering
\begin{tikzpicture}[scale = 0.72]
\begin{axis}[
	xmin = 0, 
	xmax = 5,
	ymin = 100,
	ymax = 125,
	axis line style={->},
	axis lines = left,
	xlabel = {\small $\to r$},
	ylabel = {},
	yscale = 0.8,
	legend cell align = left,
	axis y discontinuity=crunch,
	legend style = {at = {(1,1.2)}, anchor = north east}]
	
\addplot[very thick] file {Introduction/Version_R1/tikz/Optimization/lambda100_exact.txt};
\addplot[very thick, dashed, col1] file {Introduction/Version_R1/tikz/Optimization/lambda100_asymptotic.txt};
\legend{{$s^*_\l(r)$},$s^{\rm srs}_\l(r)$}
\end{axis}
\end{tikzpicture}
\caption{$\l=100$}
\end{subfigure}
\begin{subfigure}{0.48\textwidth}\centering
\begin{tikzpicture}[scale = 0.72]
\begin{axis}[
	xmin = 0, 
	xmax = 5,
	ymin = 500,
	ymax = 550,
	axis line style={->},
	axis lines = left,
	xlabel = {\small $\to r$},
	ylabel = {},
	yscale = 0.8,
	legend cell align = left,
	axis y discontinuity=crunch,
	legend style = {at = {(1,1.2)}, anchor = north east}]
	
\addplot[very thick] file {Introduction/Version_R1/tikz/Optimization/lambda500_exact.txt};
\addplot[very thick, dashed, col1] file {Introduction/Version_R1/tikz/Optimization/lambda500_asymptotic.txt};
\legend{{$s^*_\l(r)$},$s^{\rm srs}_\l(r)$}
\end{axis}
\end{tikzpicture}
\caption{$\l=500$}
\end{subfigure}
\caption{Optimal staffing levels as a function of $r = a/q$.}
\label{fig:MMs_staffing_levels_optimization}
\end{figure}
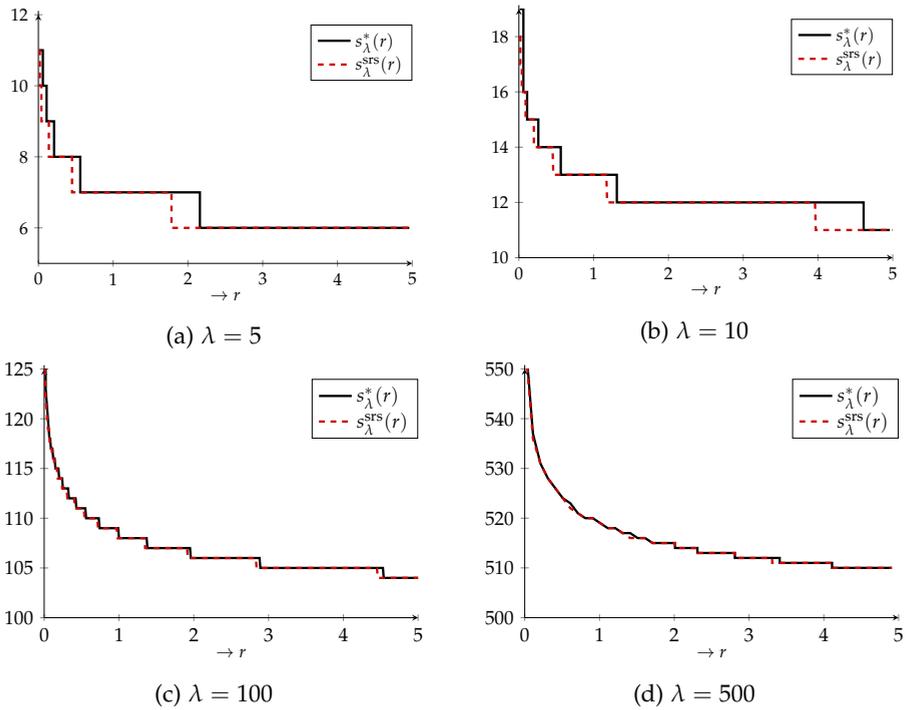

One can also consider optimization problems, for instance to strike the right balance between the costs for servers and costs incurred by customer dissatisfaction. 
More specifically, assume a salary cost of $a$ per server per unit time, and a penalty cost of $q$ per waiting customer per unit time, yielding the total cost function
\[
\bar{C}_\l(s) := a\,s + q\,\l\E[W^{(s)}]
\]
and then ask for the staffing level $s$ that minimizes $\bar{C}_\l(s)$.
Since $s>\l$, we have $\bar{C}_\l(s) > a\,\l$ for all feasible solutions $s$. 
Moreover, the minimizing value of $\bar{C}_\l$ is invariant with respect to scalar multiplication of the objective function. 
Hence we have to optimize
\begin{equation}
\label{eq:optimization_objective}
C_\l(s) = r\,(s-\l) + \l\E[W^{(s)}], \qquad r = a/q.
\end{equation}
Denote by $s^*_\l(r) := \arg\min_{s > \l} C_\l(s)$ the true optimal staffing level.
With $\sl = \l + \beta\sqrt{\l}$ and the QED limit in \eqref{eq:halfinwhitt_wait}, we can replace \eqref{eq:optimization_objective} by its asymptotic counterpart:
\begin{align*}
\frac{C_\l(\sl)}{\sqrt{\l}} = r\,\beta + \sqrt{\l} \E[W^{(s)}] \to r\,\beta + \frac{g(\beta)}{\beta} =: \hat{C}(\beta), \qquad \l\to\infty.
\end{align*} 
Once again we obtain a limiting objective function that is easier to work with than its exact pre-limit counterpart.
Hence, in the spirit of the asymptotic staffing procedure in the previous subsection, we propose the following method to determine the staffing level that minimizes overall costs.
First, (numerically) compute the value $\beta^*(r) = \arg\min_{\beta>0} \hat{C}(\beta)$, which is well-defined, because the function $\hat{C}(\beta)$ is strictly convex for $\beta>0$. 
Then, set $s^{\rm srs}_\l(r) = [ \l + \beta^*(r) \sqrt{\l} ]$. 
In Figure \ref{fig:MMs_staffing_levels_optimization} we compare the outcomes of this asymptotic staffing procedure against the true optima as a function of $r\in(0,\infty)$, for several values of $\l$. 
The staffing levels $s^{\rm srs}_\l(r)$ and $s^*_\l(r)$ are aligned for almost all $r$, and differ no more than one server for all instances. 

\subsection{Time-varying dimensioning}

So far we have only considered queues in which the model primitives are constant over time. 
In practice, though, the arrival rate can fluctuate and depends on the time of day, the day of the week, season or even larger time scales. 
It is therefore more realistic to describe these mostly predictable fluctuations through $\l(t)$, which represents the instantaneous arrival rate of the arrival process at time $t\in \mathbb{R}$. 
The existence of time-varying demand requires a re-evaluation of staffing levels throughout the planning horizon as well. 
That is, the number of servers $s(t)$ becomes a function of time, rather than a constant and this clearly asks for an adaptation of the dimensioning procedures in Sections \ref{sec:intro_constraint} and \ref{sec:intro_optimization}.

We explain the concept of time-varying staffing and the connection with the QED regime through the time-varying extension of the $M/M/s$ queue known as the $M_t/M/s_t$ queue, where the subscript $t$ refers to the time-varying nature of both the arrival process and the staffing level.
In this setting, customers arrive according to a non-homogeneous Poisson process with rate function $\l(t)$ and customers have exponentially distributed service times with mean $1/\mu$.
Under a constraint satisfaction strategy, we aim to find the staffing function $s(t)$ such that the delay probability is at most $\varepsilon\in(0,1)$ for all $t$. 
The analysis and optimization of time-varying many-server queueing systems is known to be intrinsically hard, but many approximation techniques and heuristic methods have been proposed throughout the years \cite{Green1991,Jennings1996}.

A natural but naive approach is the \textit{pointwise-stationary approximation} (PSA) \cite{Green1991}, which evaluates the system at time $t$ as if it were in steady-state with instantaneous parameters $\l=\l(t)$, $\mu$ and $s = s(t)$. 
Consequently, the analysis and optimization of queues is performed on steady-state performance metrics. 
Variants of the PSA method include the \textit{simple-stationary approximation} (SSA) \cite{Green2001}, which uses the long-term (moving) average arrival rate instead of the instantaneous arrival rate, and the \textit{stationary-independent-period-by-period approximation} (SIPP) \cite{Green2001}, which splits the time-horizon into multiple intervals and performs steady-state analysis with the averaged parameters in each of these intervals, among others.
PSA performs well in slowly varying environments with relatively short service times \cite{Green1991,Whitt1991}. 
However, when the model parameters fluctuate significantly, as is often the case in real-life systems, the accuracy of PSA can be poor, as we will see in the numerical experiment at the end of this section. 

The main reason why PSA, SSA and SIPP can fail is that these methods neglect that customers are actually residing in the system (being in service or waiting in the queue) for some time.
In contrast, staffing decisions should be based on the number of customers present in the system rather than the arrival rate at that particular time. 
Jennings et al.~\cite{Jennings1996} introduced a more sophisticated method that exploits the relation with infinite-server queues. 
We explain their idea in the context of the $M_t/M/s_t$ queue. 
By Eick et al. \cite{Eick1993}, the number of customers in the $M_t/M/\infty$ queue at time $t$ is Poisson distributed with mean
\begin{equation}
\label{eq:offered_load_eick}
R(t) = \E\left[ \l(t-B_e)\right] \E[B] = \int_0^\infty \l(t-u)\,\P(B>u)\, {\rm d}u = \int_0^\infty \l(t-u)\, \ee^{-\mu u} \,{\rm d}u.
\end{equation}
We remark that this result holds for more general service time distributions.
Now, recall that in large systems in the QED regime, the expected delay is negligible. 
Therefore, under these conditions, the many-server system may be approximated by the infinite-server approximation with offered load as in \eqref{eq:offered_load_eick}.
Accordingly, we can determine the staffing levels $s(t)$ for each $t$ based on steady-state $M/M/s$ measures with offered load $R=R(t)$.
Jennings et al. \cite{Jennings1996} proceed by exploiting the heavy-traffic results of Halfin-Whitt \eqref{eq:halfinwhitt_wait}.
In conjunction with the dimensioning scheme in Section \ref{sec:intro_constraint}, the authors propose to set
\begin{equation}
s(t) = \bigg\lceil R(t) + \beta^*(\eps) \sqrt{R(t)} \bigg\rceil,
\end{equation}
where $\beta^*(\eps)$ solves $g(\beta^*(\eps)) = \eps$.
Remark that the number of servers is rounded up to ensure that the achieved delay probability is indeed below $\eps$. 
This method was called in \cite{Jennings1996,Massey1994} the \textit{modified-offered-load} (MOL) approximation, and we adopt this term in this thesis. 

\begin{figure}
\centering
\begin{subfigure}{0.45\textwidth}
\centering
\begin{tikzpicture}[scale = 0.7]
\begin{axis}[
	xmin = 0, 
	xmax = 24,
	ymin = 0,
	ymax = 45,
	axis line style={->},
	axis lines = left,
	xlabel = {\small $\to t$},
	yscale = 0.8,
	legend style = {at = {(axis cs: 0.5,43)},anchor = north west}
	]
	
\addplot[very thick] file {Introduction/Version_R1/tikz/TimeVarying/arrival_rate.txt};
\addplot[very thick, col1] file {Introduction/Version_R1/tikz/TimeVarying/offered_load.txt};
\legend{{$\l(t)$},$R(t)$}
\end{axis}
\end{tikzpicture}
\caption{Arrival rate and offered load functions}
\end{subfigure}
\begin{subfigure}{0.45\textwidth}
\centering
\begin{tikzpicture}[scale = 0.7]
\begin{axis}[
	xmin = 0, 
	xmax = 24,
	ymin = 0,
	ymax = 60,
	axis line style={->},
	axis lines = left,
	xlabel = {\small $\to t$},
	yscale = 0.8,
	legend style = {at = {(axis cs: 0.5,59)},anchor = north west}
	]
	
\addplot[very thick] file {Introduction/Version_R1/tikz/TimeVarying/s_PSA.txt};
\addplot[very thick, col1] file {Introduction/Version_R1/tikz/TimeVarying/s_Jennings.txt};
\legend{PSA,MOL}
\end{axis}
\end{tikzpicture}
\caption{Staffing functions.}
\end{subfigure}
\caption{Time-varying parameters of a real-world emergency department.}
\label{fig:intro_example_arrival}
\end{figure}
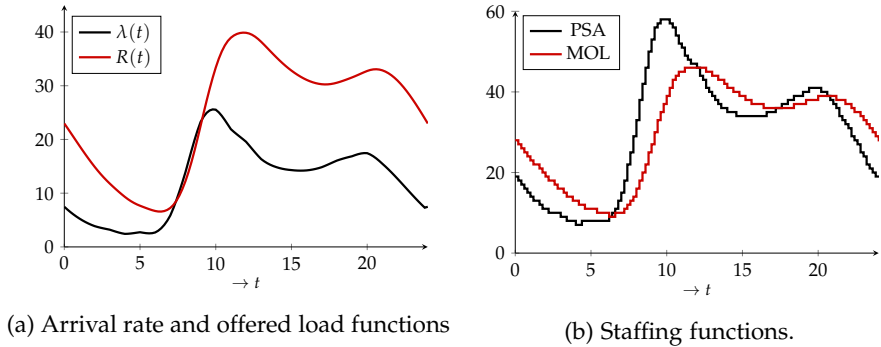

Let us demonstrate that this approximation scheme works.
Figure \ref{fig:intro_example_arrival}(a) shows an arrival rate pattern $\l(t)$ and corresponding offered load function $R(t)$ for $\mu=1/2$. 
This arrival rate stems from a real-world emergency department~\cite{Sinreich2005}. 
The resulting staffing level functions based on the PSA and MOL approximations with $\eps = 0.3$ are plotted in 
Figure \ref{fig:intro_example_arrival}(b).

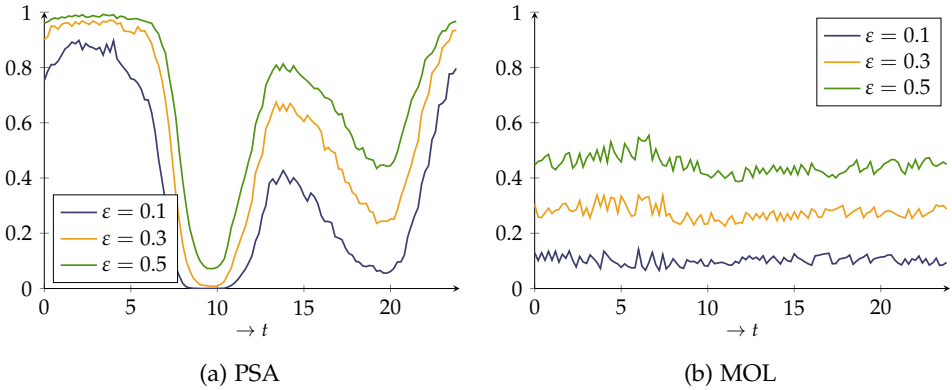
\begin{figure}
\centering
\begin{subfigure}{0.49\textwidth}
\begin{tikzpicture}[scale = 0.8]
\begin{axis}[
	xmin = 0, 
	xmax = 24,
	ymin = 0,
	ymax = 1,
	axis line style={->},
	axis lines = left,
	xlabel = {\small $\to t$},
	ylabel = {},
	yscale = 0.8,
	legend style = {at = {(axis cs: 0.5,0.02)},anchor = south west}]
	
\addplot[thick, col5] file {Introduction/Version_R1/tikz/TimeVarying/pdelay_e01_psa.txt};

\addplot[thick, col2] file {Introduction/Version_R1/tikz/TimeVarying/pdelay_e03_psa.txt};

\addplot[thick, col4] file {Introduction/Version_R1/tikz/TimeVarying/pdelay_e05_psa.txt};
\legend{{$\eps=0.1$},{$\eps=0.3$},{$\eps=0.5$}}
\end{axis}
\end{tikzpicture}
\caption{PSA}
\end{subfigure}
\begin{subfigure}{0.49\textwidth}
\begin{tikzpicture}[scale = 0.8]
\begin{axis}[
	xmin = 0, 
	xmax = 24,
	ymin = 0,
	ymax = 1,
	axis line style={->},
	axis lines = left,
	xlabel = {\small $\to t$},
	ylabel = {},
	yscale = 0.8,
	legend style = {at = {(axis cs: 23.5,0.98)},anchor = north east}]
	
\addplot[thick, col5] file {Introduction/Version_R1/tikz/TimeVarying/pdelay_e01_mol.txt};

\addplot[thick, col2] file {Introduction/Version_R1/tikz/TimeVarying/pdelay_e03_mol.txt};

\addplot[thick, col4] file {Introduction/Version_R1/tikz/TimeVarying/pdelay_e05_mol.txt};
\legend{{$\eps=0.1$},{$\eps=0.3$},{$\eps=0.5$}}
\end{axis}
\end{tikzpicture}
\caption{MOL}
\end{subfigure}
\caption{Probability of delay under staffing functions obtained through PSA and MOL approximations.}
\label{fig:intr_timevarying_simulation_results}
\end{figure}

Through simulation, we evaluate the delay probability as a function of time for $\eps = 0.1,\, 0.3$ and 0.5. 
In Figure \ref{fig:intr_timevarying_simulation_results} we see how the PSA approach fails to stabilize the performance of the queue, whereas the MOL method does stabilize around the target performance. 
The erratic nature of the delay probability as a function of time can be explained by rounding effects of the staffing level. 
Since this rather simple but elegant technique to address time-varying dimensioning is provably effective, we will adopt the underlying idea of the MOL method in various different settings in this thesis.

\section{Contributions}
\label{sec:intro_beyond}
We have explained how the QED regime can be used to dimension and staff large-scale service systems. 
The basic concepts, however, were explained for the relatively simple $M/M/s$ and $M_t/M/s_t$ queue.
Many real-world service systems have essential features that are not captured by these elementary models.
We will now discuss some of these features and address the need to consider more involved models and extend the existing QED theory. 

\subsection{Non-classical scaling regimes and pre-limit behavior} 
\label{sec:intro_novel_scalings}

The QED theory is centered around the scaling relation $\sqrt{\lambda}(1-\rho_\lambda) \to \beta$, or equivalently $s_\lambda = \lambda + \beta \sqrt{\lambda} + o(\sqrt{\lambda})$, for $\lambda\to\infty$.

It is worthwhile to study how pre-limit behavior of many-server queues is affected when one deviates from the square-root staffing rule.
Consider a novel family of heavy-traffic scaling regimes, described in terms of the parameter $\eta$ for which we assume that 
\begin{equation}
\label{eq:novel_scaling_rule}
 \lambda^\eta (1-\rho_\lambda) \to \beta, \qquad \text{as } \lambda\to\infty,\ \beta > 0.
 \end{equation} 
The parameter $\eta \geq 0$ defines a whole range of possible scaling regimes, including the classic case $\eta = 1/2$, as well as the cases $\eta=0$ and $\eta=1$ investigated in Subsection \ref{sec:intro_many_server_regimes}. 
In terms of a capacity sizing rule, the condition \eqref{eq:novel_scaling_rule} is tantamount to $s_\lambda = \lambda +\beta\,\lambda^{1-\eta}$.
This framework thus bridges the gap between the QD and QED regime if $\eta\in(0,1/2)$ and the QED and ED regime if $\eta\in(1/2,1)$, in the $M/M/s$ model.
Similar capacity sizing rules have been considered in \cite{Bassamboo2010,maman} for many-server systems with uncertain arrival rates. Hence, for $\eta\in(0,1/2)$ the variability hedge is relatively large, so that the regime parameterized by $\eta\in(0,1/2)$ can  be seen as \textit{moderate} heavy traffic: heavy-traffic conditions in which the full occupancy is reached more slowly, as a function of $\lambda$, than for classical heavy traffic. See \cite{Chang1996,Puhalskii1998,Puhalskii1999,Atar2012,Atar2014,Atar2015,Atar2016} for more details.
For opposite reasons the range $\eta\in(1/2,\infty)$ corresponds to \textit{extreme} heavy traffic due to a relatively small variability hedge. 

We use the insights of Section \ref{sec:intro_QED_regime} and the connection of the QED scaling to the CLT to argue intuitively that the following trichotomy in the qualitative system behavior as $\lambda\to\infty$ holds under scaling \eqref{eq:novel_scaling_rule}. 
For $\eta \in (0,1/2)$ the empty-system probability converges to $1$, because the order of the  variability hedge $\beta \l^{1-\eta}$ is greater than strictly necessary to accommodate the stochastic fluctuations in demand. 
Scalings in which $\eta\in(1/2,\infty)$, have adverse behavior, since stochastic fluctuations are not accounted for sufficiently, so that the probability of delay converges to 1.
The value $\eta=1/2$ is therefore the tipping point, at which the delay probability converges to a limit between 0 and 1.
Above and below this critical value, the asymptotic performance of the queue flips to either one of the extremes. 

In Chapter 2, we formalize this heuristic argument and conduct an asymptotic analysis to reveal the rate at which the limit of performance metrics is attained, depending on the parameters $\eta$ and $\beta$ and the system size $\l,\sl$.

\subsection{Overdispersed arrivals}
\label{sec:intro_overdispersion}
Until now we have considered queueing systems with perfect knowledge on the model primitives, including the mean demand per time period. For large-scale service systems, the dominant assumption in the literature is that demand arrives according to a non-homogeneous Poisson process, which in practice translates to the assumption that arrival rates are known for each basic time period (second, hour or day). 

Although natural and convenient from a mathematical viewpoint, the Poisson assumption often fails to be confirmed in practice. A deterministic arrival rate implies that the demand over any given period is a Poisson random variable, whose variance equals its expectation. A growing number of empirical studies of service systems shows that the variance of demand typically exceeds the mean significantly, see \cite{Avramidis:2004, Bassamboo2010, Bassamboo2009, Brown2005, Chen2001, Gans2003, Gurvich2010, koolejongbloed, kimwhitt,  maman, Mehrotra2010, Robbins2010, Steckley2009, Zan2012}. The feature that variability is higher than one expects from the Poisson assumption is referred to as \textit{overdispersion}.

Due to its inherent connection with the CLT, the dimensioning rule in \eqref{eq:square_root_staffing rule} relies heavily on the premise that the variance of the number of customers entering the system over a period of time is of the same order as the mean.
Subsequently, when stochastic models do not take into account overdispersion, resulting performance estimates are likely to be overoptimistic. The system then ends up being underprovisioned, which possibly causes severe performance problems, particularly in critical loading.

To deal with overdispersion, existing capacity sizing rules like the square-root staffing rule need to be modified in order to incorporate a correct hedge against (increased) variability.  
Following our findings in Section \ref{sec:intro_characteristics}, we propose a capacity allocation rule similar to \eqref{eq:square_root_staffing rule} in which the original variability hedge is replaced by an amount that is proportional to the square-root of the variance of the arrival process. 

In Chapter 3, we elaborate on this idea and show how to adapt the scaling of the queueing process appropriately to achieve QED-type behavior in the presence of overdispersion.

\subsection{Finite-size constraints}

The canonical examples in Section \ref{sec:intro_QED_regime} assume an infinite amount of waiting space.
Physical service systems, however, are sometimes limited in the number of customers that can be held in the system simultaneously.
For instance in a call center, the maximum number of clients in service or queueing is restricted  by the number of available trunk lines \cite{Khudyakov2006}, while in the emergency department of a hospital, the number of beds constrains the number of patients that can be admitted \cite{YomTov2010}. 
Depending on the practical setting and admission policy, if the maximum capacity, say $n$, is reached, newly arriving customers either leave the system immediately (blocking), reattempt  getting access later (retrials) or queue outside the facility (holding). 
In any case, expectations are that the queueing dynamics within the service facility are affected considerably in the presence of such additional capacity constraints.

We illustrate these implications through the $M/M/s/n$ queue, that is, the standard $M/M/s$ queue with additional property that a customer who finds upon arrival $n$ customers already present in the system, is deferred and considered lost. 
To avoid trivialities, let $n\geq s$. 
Since the expected workload reaching the servers is less than in the unconstrained scenario, one expects less congestion and resource utilization. 

Consider the $M/M/\sl/n_\l$ in the QED regime. 
So, let $\l$ increase while $\sl$ scales as $\sl=\l+\beta\sqrt{\l} + o(\sqrt{\l})$. 
We then ask how $n_\l$ should scale along with $\l$ and $\sl$ to maintain the non-degenerate behavior as seen in Section \ref{sec:intro_QED_regime}. 
We provide a heuristic answer. 
Let $Q^{(\sl,n_\l)}$ and $W^{(\sl,n_\l)}$ denote the number of customers in the system and the waiting time in the $M/M/\sl/n_\l$ queue in steady state.
Note through Proposition \ref{thm:intro_HW_stationary_distribution} that if there were no finite-size constraints, we would have, for $\l$ large,
\begin{align}
\P(Q^{(\sl)}\geq n_\l)  
&= \P\left(\frac{Q^{(\sl)}-\sl}{\sqrt{\sl}} \geq \frac{n_\l-\sl}{\sqrt{\sl}}\right) \nonumber \\
&\to
\left\{
\begin{array}{ll}
g(\beta), & \text{if }n_\l = \sl + o(\sl),\\
g(\beta)\,\ee^{-\beta \gamma}, & \text{if } n_\l = \sl+\gamma\sqrt{\sl} + o(\sqrt\sl),\\
0, & \text{if } n_\l = \sl+\Omega(\sqrt{\sl}),
\end{array}
\right.
\end{align}
as $\l\to\infty$ for some $\gamma>0$.
Here, the relation $u(\l) = \Omega(v(\l))$ implies $u(\l)/v(\l) >1$ for $\l\to\infty$.
Hence, asymptotically the finite-size effects only play a role if the extra variability hedge of $n_\l$ is of order $\sqrt{\sl}$ (or equivalently $o(\sqrt{\l})$). 
Furthermore, if the variability hedge is $o(\sqrt{\l})$, then we argue that asymptotically, all customers who do enter the system have probability of delay equal to zero. 
More formally, under the \textit{two-fold scaling rule} 
\begin{equation}
\label{eq:intro_twofold_scaling_rule}
\left\{
\begin{array}{ll}
\sl = \l + \beta\sqrt{\l} + o(\sqrt{\l}),\\
n_\l = \sl + \gamma \sqrt{\sl} + o(\sqrt{\l}),
\end{array}
\right.
\end{equation}
it is not difficult to deduce that, see e.g. \cite{masseywallace},
\begin{equation}
\P(W^{(\sl,n_\l)} > 0) \to \left( 1 + \frac{\beta\,\F(\beta)}{(1-\ee^{-\beta\gamma})\f(\beta)}\right)^{-1}, \quad \text{as } \l\to\infty,
\end{equation}
which is strictly smaller than $g(\beta)$ in \eqref{fig:delay_probs_HW_MMs}, but still bounded away from both 0 and 1.
Furthermore, the buffer size of the queue is $n_\l-\sl = \gamma\sqrt{\sl}$, so that by Little's law, the expected waiting time of an admitted customer is $O(1/\sqrt{\sl})$. 
Even though resource utilization in the $M/M/\sl/n_\l$ is less efficient than in the queue with unlimited waiting space, it can easily be shown that $\rho\to 1$ as $\l\to\infty$.
Hence, all three key characteristics of the QED regime are carried over to the finite-size setting if adhered to scaling \eqref{eq:intro_twofold_scaling_rule}.

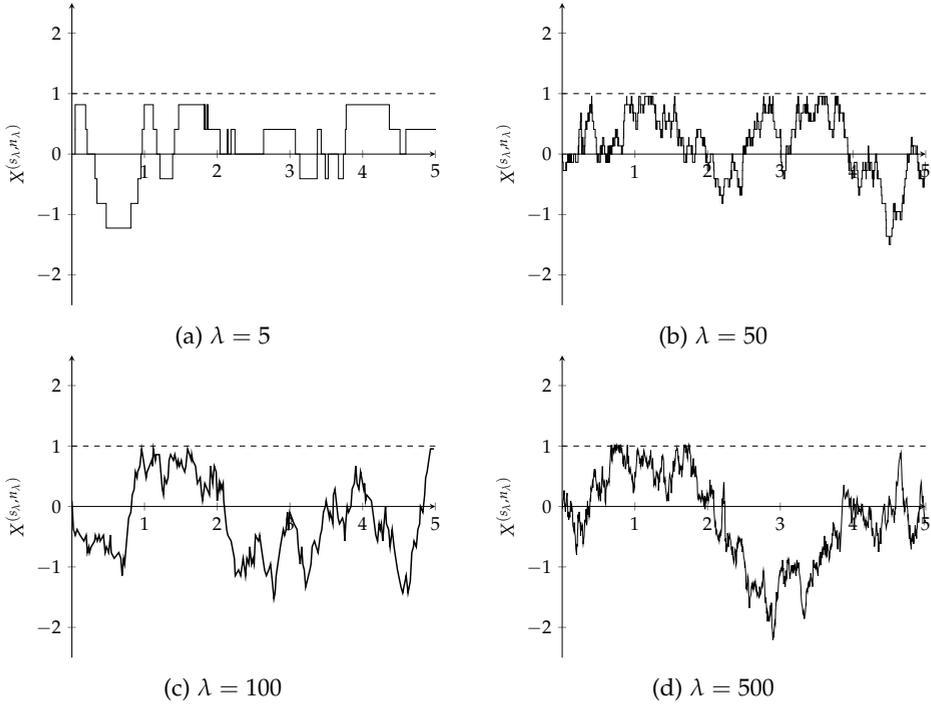
\begin{figure}
\centering
\begin{subfigure}{0.49\textwidth}
\centering
\begin{tikzpicture}[scale = 0.7]
\begin{axis}[
	xmin = 0, 
	xmax = 5,
	ymin = -2.5,
	ymax = 2.5,
	axis line style={->},
	axis x line=middle,
	axis y line=left,
	ylabel = {$X^{(\sl,n_\l)}$},
	xscale=1,
	yscale=1]	
\addplot[] file {Introduction/Version_R1/tikz/SamplePaths_MMsn/lambda5.txt};
\addplot[dashed] coordinates{ (0,1) (5.2,1) };
\end{axis}
\end{tikzpicture}
\caption{$\l=5$}
\end{subfigure}
\begin{subfigure}{0.49\textwidth}
\centering
\begin{tikzpicture}[scale = 0.7]
\begin{axis}[
	xmin = 0, 
	xmax = 5,
	ymin = -2.5,
	ymax = 2.5,
	axis line style={->},
	axis x line=middle,
	axis y line=left,
	ylabel = {$X^{(\sl,n_\l)}$},
	xscale=1,
	yscale=1]	
\addplot[] file {Introduction/Version_R1/tikz/SamplePaths_MMsn/lambda50.txt};
\addplot[dashed] coordinates{ (0,1) (5.2,1) };
\end{axis}
\end{tikzpicture}
\caption{$\l = 50$}
\end{subfigure}
\begin{subfigure}{0.49\textwidth}\centering
\begin{tikzpicture}[scale = 0.7]
\centering
\begin{axis}[
	xmin = 0, 
	xmax = 5,
	ymin = -2.5,
	ymax = 2.5,
	axis line style={->},
	axis x line=middle,
	axis y line=left,
 	xlabel style={right},
 	ylabel = {$X^{(\sl,n_\l)}$},
	xscale=1,
	yscale=1]	
\addplot[thick] file {Introduction/Version_R1/tikz/SamplePaths_MMsn/lambda100.txt};
\addplot[dashed] coordinates{ (0,1) (5.2,1) };
\end{axis}
\end{tikzpicture}
\caption{$\l=100$}
\end{subfigure}
\begin{subfigure}{0.49\textwidth}\centering
\begin{tikzpicture}[scale = 0.7]
\centering
\begin{axis}[
	xmin = 0, 
	xmax = 5,
	ymin = -2.5,
	ymax = 2.5,
	axis line style={->},
	axis x line=middle,
	axis y line=left,
 	xlabel style={right},
 	ylabel = {$X^{(\sl,n_\l)}$},
	xscale=1,
	yscale=1]	
\addplot[] file {Introduction/Version_R1/tikz/SamplePaths_MMsn/lambda500.txt};
\addplot[dashed] coordinates{ (0,1) (5.2,1) };
\end{axis}
\end{tikzpicture}
\caption{$\l=500$}
\end{subfigure}
\caption{Sample paths of the normalized queue length process $X^{(\sl,n_\l)}(t)$ with $\l = 5,\, 50,\, 100$ and $500$ under scaling \eqref{eq:intro_twofold_scaling_rule} with $\beta=0.5$ and $\gamma = 1$.}
\label{fig:sample_paths_MMsn}
\end{figure}

On a process level, adding a capacity constraint translates to adding a reflection barrier to the normalized queue length process $X^{(\sl,n_\l)} = (Q^{(\sl,n_\l)} -\sl ) /\sqrt{\sl}$, at $\gamma$, as is illustrated by the sample paths of $X^{(\sl,n_\l)}$ for three values of $\l$ in Figure \ref{fig:sample_paths_MMsn}.

\
It has been shown by \cite{masseywallace} that under \eqref{eq:intro_twofold_scaling_rule} 
\begin{equation}
\label{eq:asymptotic_blocking_prob}
\sqrt{\sl}\,\P({\rm block}) = \sqrt{\sl}\, \P(Q^{(\sl,n_\l)} = n_\l) \to f(\beta,\gamma), \quad \text{as } \l \to\infty,
\end{equation}
for a non-negative function $f$. 

The idea of the two-fold scaling in \eqref{eq:intro_twofold_scaling_rule} can be extended to settings in which the interior is in fact a network of queues, rather than the single-station setting discussed here, see \cite{Khudyakov2006,YomTov2010,Tan2012} for examples of such \textit{semi-open} queueing networks.

When customers retry getting access after being blocked initially, the QED analysis becomes much more difficult, and no explicit limiting results are known. 
Nevertheless, observe that the volume of blocked arrivals is by \eqref{eq:asymptotic_blocking_prob} of order $\sqrt{\l}$, the exact same magnitude as the variability hedge of both $\sl$ and $n_\l$. 
Therefore, retrials and holding customers have a non-negligible effect on the service levels within the facility in the QED regime. 
This will be the topic of Chapters 4 and 5.

\subsection{Pre-limit behavior}
The results on queues in the QED regime discussed in Section \ref{sec:intro_QED_regime} are in two ways of an asymptotic nature. 
First, the heavy-traffic limits prescribe the queueing dynamics for $\l,\sl\to\infty$. 
Real-world systems obviously do not experience infinite demand nor have infinite capacity, and hence the heavy-traffic limits only form an approximation for such finite-sized systems. 
Although these approximations are qualitatively insightful, the asymptotic analyses do not reveal much about their accuracy with respect to actual performance. 
For instance, we would like to know how fast the convergence takes place, and how inaccuracies in asymptotic approximations percolate into inaccuracies in pre-limit systems. 
To answer such questions, it would be helpful to have an asymptotic estimate for the difference between the (scaled) queueing process and its limiting counterpart, to be able to judge the error made by relying on asymptotic as opposed to actual performance evaluation. 
Characterization of the error term gives rise to so-called \textit{corrected diffusion approximations}~ \cite{Siegmund1978,Blanchet2006,Janssen2008}, which are refinements to heavy-traffic limits for finite systems, and are also related to \textit{universal approximations} \cite{Gurvich2014,Huang2016,Braverman2015,Braverman2015a}.
We will derive such corrected diffusion approximations in the context of the novel scaling regimes mentioned in Section \ref{sec:intro_novel_scalings} in Chapter 2.

Second, the bulk of queueing literature is concerned with the performance analysis and optimization of steady-state systems.
However, in practice, service systems certainty do not run infinitely long, which renders this assumption questionable.
Validation of the steady-state assumption is related to the \textit{relaxation time} of a queueing process \cite{Abate1987,Abate1988,relaxation,Leeuwaarden2011,Leeuwaarden2012,Gamarnik2013}, which prescribes the time it takes a system starting out of equilibrium to approximate its stationary distribution. 
In case the relaxation time is small, stationary performance evaluation is likely to be accurate. 
On the contrary, if the relaxation time is large, a time-dependent analysis of the queueing system is required in order to capture realistic behavior. 
Subsequently, we can investigate the implications of applying staffing principles that are based on steady-state performance metrics in settings which are inherently transient over the planning period.
We will touch upon this topic in Chapter 6.

\section{Outline of the thesis}

The remainder of this thesis builds upon the ideas behind the QED scaling regime exhibited in this introductory chapter, and is organized as follows. 

Chapter 2 is concerned with the analysis of the limiting behavior of queues in case one deviates from the square-root staffing principle as demand grows large. 
Using the bulk-service queue together with the many-sources paradigm as a vehicle, we derive corrected diffusion approximations for the performance metrics of pre-limit systems in these alternative scaling regimes. 
The work presented in Chapter 2 is based on \cite{Janssen2015}.

In Chapter 3, we also analyze the bulk-service queueing model, but with many correlated sources, so that demand becomes overdispersed. 
As we alluded to in Section \ref{sec:intro_overdispersion}, this requires an alternative scaling of the queue length process and associated staffing rule. 
This chapter exhibits the ideas of \cite{Mathijsen2016}.

In Chapter 4, we discuss how QED-type behavior prevails in simple settings in which the system size is finite, given appropriate capacity-sizing rules.
More specifically, we show how customer retrials can be incorporated heuristically into the performance analysis of finite-size systems in the QED regime. 
The content of this chapter is based on \cite{Leeuwaarden2015} and \cite{Leeuwaarden2016}. 

Building upon the insights gained in Chapter 4, we show in Chapter 5 how the approximation methods carry over to a more complex finite-size queueing system, inspired by delay analysis in a health care facility. 
We show how the QED scaling limits for this model offer surprisingly accurate approximations for realistic model parameters in systems of small to moderate size, and develop a staffing algorithm for dimensioning such systems. 
Chapter 5 is based on the ideas of \cite{Leeuwaarden2016a}.

Chapter 6 investigates the validity of a capacity allocation rule based on steady-state performance metrics in practical settings.
Namely, in realistic scenarios, the parameters of a queueing model are typically subject to change over the planning period. 
This asks for a more elaborate transient analysis of the queue dynamics, and an adaptation of the staffing level. 
In this chapter, we present how to do so appropriately in a single-server queueing model facing a L\'evy input process by prescribing a correction to the steady-state optimum, which has a square-root form. 
This chapter is based on \cite{Mathijsen2016a}.

Chapter 7 presents the analysis of an inventory model with backlogs, perishable goods and consumer impatience. 
This model resembles the inventory level of a blood bank, and can be regarded as a shot-noise model with both positive and negative jumps and exponential decay rates above and below zero.
Besides the derivation of the stationary distribution of the inventory level, we show how under appropriate scaling the process converges to an Ornstein-Uhlenbeck process. 
The latter allows for a more tractable approximate analysis of the model in case the number of blood deliveries and demand is large. 
Chapter 7 is based on \cite{Bar-Lev2015}.

\chapter{Novel heavy-traffic regimes}

\begin{chapterstart}
In this chapter, we introduce a family of heavy-traffic scalings for a large-scale service system meant to serve jobs coming from a large pool of customers. 
The scaling rules are inspired by the classical QED regime discussed in Chapter 1, but lead to a range of different system behaviors that includes the ED, QED and QD regime as special cases.
To determine the scaling limits, we describe the performance measures in terms of Pollaczek integrals and use asymptotic techniques to evaluate these integrals in the large-system limit. 
\end{chapterstart}

\begin{flushright}
Based on \\
\textbf{Novel heavy-traffic regimes for large-scale service systems}\\
\textit{Guido Janssen, Johan van Leeuwaarden \& Britt Mathijsen}\\
In \textit{SIAM Journal of Applied Mathematics, 75(2), 787-812 (2015)}
\end{flushright}

\newpage

\section{Introduction \& motivation}
We study the workload process of a system, experiencing stochastic demand and deterministic capacity $s_n$ per period, at equidistant time epochs. 
Demand is assumed to be generated by $n$ stochastically identical and independent sources. 
Let $A_{i,j}$ denote the workload brought into the system by source $i$ in period $j$, for which $\E[A_{i,j}] =\mu$ and $\Var A_{i,j} = \sigma^2$.
Then the total amount of demand arriving to the system in period $j$ is $A^{(n)}_j=\sum_{i=1}^n A_{i,j}$ with $\E[A^{(n)}_j] = n\mu$ and $\Var A^{(n)}_j = n\sigma^2$.
As explained in Chapter 1, a good capacity sizing rule for achieving economies-of-scale is $s_n = n\mu+\beta\sqrt{n}\sigma$ for some $\beta>0$.
If we denote the system utilization by $\rho_n := n\mu/s_n$, then this dimensioning rule in the bulk service queue with many sources is tantamount to the heavy-traffic scaling
\begin{equation}\label{bb1}
\sqrt{n}(1-\rho_n) \to \gamma = \frac{\beta\sigma}{\mu}, \qquad {\rm as }\ n\to\infty.
\end{equation}
Starting from this setting, we introduce a novel family described in terms of a parameter $\eta$ for which we assume that
\begin{equation}\label{bb}
n^{\eta}(1-\rho_n)\rightarrow \gamma, \quad {\rm as} \ n\to  \infty,  \ \gamma> 0.
\end{equation}
The parameter $\eta\geq 0$ defines a whole range of possible scaling regimes, including the classical case $\eta=1/2$. 
In terms of a capacity sizing rule for systems with many customers, the condition \eqref{bb} is tantamount to $s_n=n\mu+\beta \sigma n^{1-\eta}$. 
Similar capacity sizing rules have been considered in \cite{Bassamboo2010,maman} for many-server systems with uncertain arrival rates. Hence, for $\eta\in(0,1/2)$ the variability hedge is relatively large, so that the regime parameterized by $\eta\in(0,1/2)$ can  be seen as {\it moderate} heavy traffic: heavy-traffic conditions in which the full occupancy is reached more slowly, as a function of $n$, than for classical heavy traffic. For opposite reasons the range $\eta\in(1/2,\infty)$ corresponds to {\it extreme} heavy traffic due to a relatively small variability hedge. Note that the case $\eta=0$ does not lead to 100\% system utilization when $n\to\infty$.

In this chapter we show that economies-of-scale can be achieved for a large range of $\eta$, although the nature of the benefits obtained by operating on large scale depends on the precise capacity sizing rule (hence the parameter $\eta$). We quantify performance in terms of stationary  measures: The mean and variance of the congestion in the system, and the probability of an empty system. For these performance measures we derive heavy-traffic limits under the scalings \eqref{bb} that
 are relatively simple functions of only the first two moments of the demand per period. Such parsimonious expressions are useful for quantifying and improving system behavior. The heavy-traffic limits, however, provide also qualitative insight into the system behavior. Our asymptotic analysis shows that mean congestion is $O(n^\eta)$, which implies
that delays experienced by the customers are negligible for all values of $\eta\in [0,1)$, are roughly constant for $\eta=1$, and grow without bound for $\eta>1$.  We expect this qualitative behavior to be universal for a wide range of stochastic models to which the regime \eqref{bb} is applied.
We further show the existence of the following trichotomy as $n\to \infty$ under \eqref{bb}: For $\eta \in (0,1/2)$ the empty-system probability converges to $1$, for $\eta\in (1/2,1)$ it converges to $0$, while only for $\eta=1/2$ there is a limiting value in $(0,1)$. Hence, as expected, the system performance deteriorates with $\eta$, in a rather crude way for the empty-system probability, and in only a mild way for mean congestion levels. The regime \eqref{bb} thus presents a range of possible capacity sizing rules that all lead to economies-of-scale, and depending on what is the desired nature of performance for a particular service system, an appropriate $\eta$ can be selected. From the quantitative perspective, our detailed asymptotic analysis leads to more precise asymptotic estimates for the performance measures in heavy traffic, which reveal the exact manner in which congestion is influenced by $\eta$ and $\gamma$.\\
\\
\noindent\textbf{Motivating examples.}
The bulk service queueing model is one of the canonical models in queueing theory, having a wide range of applications in fields like digital communication, wireless networks, road traffic, reservation systems, health care and many more (see \cite{Bruneel1993} and \cite[Chap.~2]{johanthesis} for an overview). 
In road traffic, the basic model for congestion at an intersection, known as the fixed-cycle traffic-light queue \cite{Newell1960,Leeuwaarden2006}, is related to our discrete bulk service queue. 
Then $s_n$ represents the maximum number of delayed cars in front of a traffic light that can depart during one green period, while $A^{(n)}_j$ is the number of newly arriving cars during a consecutive green and red period. 

An example from health care is panel sizing  \cite{Zacharias2014}. 
Say a general practitioner has a pool of $n$ clients (typically in the order of 2,500~\cite{Green2008}), all of which are potential patients, and together require $A^{(n)}_j$ consults per day. 
Further assume that the practitioner can see a maximum number of $s_n$ patients per day. 
What is then an appropriate patient panel size $n$, which strikes a reasonable balance between accessing medical care in a timely manner and restricting the time that the practitioner sits idle? 
The panel size application is one of many examples of an appointment book, referring to some schedule of appointments for a fixed period, with capacity $s_n$ appointments per period and newly arriving number of appointments $A^{(n)}_j$ per period. 
See \cite{Dai2015} for another recent example of an appointment book in a health care setting, again in terms of our bulk service queue, with $A^{(n)}_j$ the new patients per day and $s_n$ the number of available beds. 

For all examples above, and many more, our new class of heavy-traffic scalings \eqref{bb} presents capacity sizing rules for which the expected performance can be quantified using the results in this chapter. This will be helpful in dimensioning the systems (How much capacity is needed to achieve a certain target performance?) while exploiting economies-of-scale. For appointment books, our model together with the capacity sizing rules \eqref{bb} is particularly relevant for {\it advanced access} \cite{Green2008}, a scheduling approach in health care designed to reduce delays by offering every patient a same-day appointment, regardless of the urgency of the problem. In that way, patients do not have to wait long for appointments, and practices do not waste capacity by holding appointments in anticipation of urgent cases.\\
\\*
\noindent\textbf{Pollaczek's formula.}
Next to the freedom to model different situations, another advantage of our model is that it is mathematically tractable, in the sense that it can be subjected to powerful mathematical methods from complex and asymptotic analysis. In order to establish the heavy-traffic limits we start from Pollaczek's formula for the transform of the stationary queue length distribution in terms of a contour integral. From this famous transform representation, contour integrals for the empty-system probability and the mean and variance of the congestion immediately follow. Contour integrals are often amenable to asymptotic evaluation (see e.g.~\cite{Cohen1982}), particularly for obtaining classical heavy-traffic asymptotics. 
We also subject the contour integral representations to asymptotic evaluation, but not under classical heavy-traffic scaling. 
This asymptotic analysis requires a {\it non-standard} saddle point method, tailored to the specific form of the integral expressions that arise under the capacity sizing rule \eqref{bb}. \\
\\*
\noindent
\textbf{Saddle point method.}
In complex analysis, the saddle point method in its standard form is a useful technique to estimate the asymptotic behavior of integrals of the form
\begin{equation}
\label{eq:sp_integral}
I(n) = \int_C h(z)\, \ee^{n f(z)}\, \dd z, 
\end{equation}
as $n\to\infty$, where $C$ is a contour in the complex plane, and $f(z)$ and $h(z)$ are functions that are analytic in some neighborhood of $C$. 
The main idea behind the saddle point method is that if the integrand in \eqref{eq:sp_integral} exhibits a sharp peak along the contour, then one may naturally expect that a small neighborhood around this peak provides the dominant contribution to the integral. 
More specifically, for large values of $n$, the function $f$ and its associated maximum $f(z^*)$ for $z^*\in C$ to a large extent determine the magnitude of the integrand (where $z^*$ is well-defined due to analyticity of $f$. 
In the setting of this chapter, $C$ is a closed curve, which implies that the value $z^*$ must be a \textit{saddle point} of $f$, i.e.~$f'(z^*) = 0$. 
Subsequently, one can replace $f(z)$ in \eqref{eq:sp_integral} by its Taylor expansion around $z^*$ and deduce through the Laplace method, see e.g.~\cite{debruijn}, that
\begin{equation*}
I(n) = \sqrt{2\pi}\,i\frac{h(z^*)\,\ee^{n f(z^*)}}{\sqrt{n\, |f''(z^*)|}}\Bigl( 1+ O(1/n)\Bigr),
\end{equation*}
as $n\to\infty$.
In Section \ref{spSec}, we show how the contour integrals describing stationary measures for the queue length, derived through Pollaczek's formula, can be reformulated into the shape of \eqref{eq:sp_integral}.
However, we will show that the saddle point method in its standard from cannot be applied to asymptotically characterize other stationary measures like the mean or mass at zero. 
Indeed, for our model the saddle point (the solution of \eqref{e21}) converges to one (as $n\to\infty$), which is a singular point of the integrand, and renders the standard saddle point method useless.
The non-standard saddle point method discussed in this chapter, originally proposed by \cite{debruijn}, is made specifically to overcome this complication. This leads to asymptotic expansions for the  performance measures, of which the limiting forms correspond to the heavy-traffic limits, and pre-limit forms present refined approximations for pre-limit systems ($n<\infty$) in heavy traffic. Such refinements to heavy-traffic limits are commonly referred to as {\em corrected diffusion approximations} \cite{Siegmund1978,Blanchet2006,Asmussen2003}.

\noindent{\bf Further connections to the literature.}
We now discuss two classes of stochastic systems for which the heavy-traffic regime \eqref{bb1} has been studied extensively, and for which our new family of regimes \eqref{bb} is largely unexplored. We discuss these classes because, despite the fact that the Pollaczek formula does not hold, we believe the qualitative results that we reveal for our particular model should to a large extent carry over to these settings as well, presenting some interesting avenues for further research (see Section \ref{subsec62}).

The first class concerns so-called  {\it nearly-deterministic} systems \cite{Sigman2011a,Sigman2011b}, denoted by $G_n/G_n/1$ system, where $G_n$ stands for {\it cyclic thinning} of order $n$, indicating that some  point process is thinned to contain only every $n$th point. As $n\to  \infty$, the $G_n/G_n/1$ systems approach the deterministic $D/D/1$ system. For $G_n/G_n/1$ systems, \cite{Sigman2011a} establishes stochastic-process limits, and \cite{Sigman2011b} derives heavy-traffic limits for stationary waiting times. In the framework of \cite{Sigman2011a,Sigman2011b}, our stochastic model corresponds to a $D/G_n/1$ queue, where the sequence of service times $\{A^{(n)}_j\}_{j\geq 1}$ follows from a cyclically thinned sequence of i.i.d.~random variables $A_{i,j}$. It follows from \cite[Theorem 3]{Sigman2011b} that the rescaled stationary waiting time process converges under \eqref{bb1} to a reflected Gaussian random walk. Hence, the performance measures of the nearly deterministic system, under \eqref{lind} and  \eqref{bb1}, should be well approximated by the performance measures of the reflected Gaussian random walk, giving rise to heavy-traffic approximations. This connection is discussed in detail in Section \ref{subsec3.2}. It seems likely that results similar to those presented  in this chapter can be obtained by applying the scaling \eqref{bb} to the nearly-deterministic systems in \cite{Sigman2011a,Sigman2011b}, and because Pollaczek's formula also applies to this setting, the non-standard saddle point method developed in this chapter can provide the appropriate methodology.

The second class concerns multi-server systems, and in particular the many-server regime. When we interpret $s_n$ as the number of servers, instead of capacity per time slot or order of thinning, the scaling \eqref{bb1} is similar to the QED or Halfin-Whitt regime for the $M/M/s_n$ system.
As we have reviewed in Chapter 1, the QED regime is characterized by a delay probability that converges to a non-degenerate limit away from both zero and one, and the mean delay is asymptotically negligible as the number of servers grows large. The QED regime \eqref{bb1} is naturally positioned in between the Quality-Driven (QD) regime and the Efficiency-Driven (ED) regime. In the QD regime, the load remains bounded away from 1, which corresponds to setting $\eta=0$ in \eqref{bb}. Hence, the range $\eta\in(0,1/2)$ bridges the gap between the QED regime and the QD regime. Likewise, the ED regime corresponds to setting $\eta=1$ in \eqref{bb}, so that the range $\eta\in(1/2,1]$ connects the QED regime and ED regime. For the birth-death process describing the $M/M/s_n$ system, Maman \cite{maman} introduced a scaling similar to \eqref{bb}, and called it the QED-$c$ regime, also bridging the ED and QD regimes. 
Theorem 4.1 of \cite{maman} says that the expected waiting time under the scaling $s_n = n\mu+\beta\sigma n^{1-\eta}$ is of order $s_n^{1-\eta}$, which is equivalent to the expected queue length being of order $n^\eta$ by Little's law. We should stress though that we expect the mathematical techniques that are needed to establish heavy-traffic results could be entirely different than in this chapter, because Pollaczek's formula does not apply to many-server settings. 

The specific model assumptions will determine to a large extent the appropriate methodology. Under Markovian assumptions leading to the $M/M/s_n$ system, simple exact solutions are available for the stationary distribution. This makes it possible to describe performance measures like the mean congestion directly in terms of real integrals. Where the saddle point method is used for integrals in the complex plane, the Laplace method (see e.g.~\cite{flajolet}) is used for real integrals. Hence, for the asymptotic evaluation of the $M/M/s_n$ system under the scaling \eqref{bb}, the Laplace method seems an appropriate methodology, although again one needs to deal with possible singularities in the integrand. For $G/D/s_n$ systems, which assume deterministic service times, it has been shown in \cite{Jelenkovic2004} that using a decomposition property the dynamics of this multi-server systems can be captured in terms of a single-server system. Hence, for these systems, Pollaczek's formula applies, and our saddle point method can most likely be applied to obtain heavy-traffic results in the regimes \eqref{bb}. Under more general conditions, for instance leading to a $G/G/s_n$ system, it is simply unclear at this stage how to obtain precise heavy-traffic approximations for \eqref{bb}, because a tractable description of the performance measures is not available; see Section 1.2.4 for details.\\
\\
\\*
\noindent{\bf Structure of the chapter.}
In Section \ref{sec1} we present in detail the model and the family of heavy-traffic scalings. In Section \ref{spSec} we introduce the saddle point method. In Section \ref{sec3} we apply the saddle point method for the mean congestion level. Theorem \ref{mainthm} gives for all heavy-traffic scalings the limiting behavior in terms of an integral expression. As a consequence, we show in
 Proposition \ref{prop1} that there are two types of heavy-traffic behavior, depending on whether $\eta\in(0,1/2)$ or $\eta\geq 1/2$.
In Section \ref{subsec3.2} we discuss for the case $\eta=1/2$ the connection with the Gaussian random walk and the Riemann zeta function.
In fact, we show that for all $\eta\geq 1/2$ there exists a connection between the integral expression in Theorem \ref{mainthm} and the Riemann zeta function.
In Section \ref{more} we apply the saddle point method to obtain several more heavy-traffic results, including refined heavy-traffic approximations for the mean congestion level, and the leading heavy-traffic behavior for the variance of the stationary congestion level and for the empty-system probability.
Finally, in Section \ref{numm} we confirm through numerical experiments the accuracy of our heavy-traffic approximations, and moreover show that under \eqref{bb}, various multi-server systems behave similar to our discrete bulk service queue.

\section{Model description \& heavy-traffic regimes}\label{sec1}
We thus consider a discrete stochastic model in which time is divided into periods of equal length. At the beginning of each period $j=1,2,3,...$ new demand $A^{(n)}_j$ arrives to the system. The demands per period $A^{(n)}_1,A^{(n)}_2,\ldots$ are assumed independent and equal in distribution to some non-negative integer-valued random variable $A^{(n)}$.
We will omit the superscript $(n)$ if no ambiguity is possible.
The system has a service capacity $s_n\in\mathbb{N}$ per period, so that the recursion
\begin{equation}
\label{lind}
Q(j+1) = \max\{Q(j) + A^{(n)}_j - s_n,0\},\qquad j=1,2,...,
\end{equation}
assuming $Q(0)=0$, gives rise to a Markov chain $\{Q(j)\}_{j\geq 1}$ that describes the congestion in the system over time. The probability generation function (pgf)
\begin{equation*} 
\tilde A(z)=\sum_{k=0}^{\infty} \mathbb{P}\big(A^{(n)}=k\big) z^k
\end{equation*}
of $A^{(n)}$ is assumed analytic in a disk $|z|<r$ with $r>1$, which implies that all moments of $A^{(n)}$ exist. We also assume that
\begin{equation} \label{e3}
\tilde A'(1)=\E[A^{(n)}_j]=\mu_A<s_n.
\end{equation}

Under the assumption \eqref{e3} the function $z^{s_n}-\tilde A(z)$ has exactly $s_n$ zeros in the closed unit disk, one of these being $z=1$ (see \cite{rouche}).
We further assume that $\mathbb{P}(A^{(n)}=j)>0$ for some $j>s_n$.
Under this assumption the function
$z^{s_n}-\tilde A(z)$ also has zeros outside $|z|\leq 1$, and we let $r_0$ be the minimum modulus of these zeros.
The number $r_0$ is the unique zero of $z^{s_n}-\tilde A(z)$ with real $z>1$; see e.g.~\cite{Janssen2005}.
Moreover, under assumption \eqref{e3} the stationary distribution $\lim_{j\to  \infty}\mathbb{P}\left(Q(j)=k\right)=\P(Q=k)$, $k=0,1,\ldots$ exists, with the random variable $Q$ defined as having this stationary distribution.

We let
\begin{equation*} 
\tilde Q(w)=\sum_{j=0}^{\infty}\mathbb{P}(Q=j)w^j
\end{equation*}
be the pgf of the stationary distribution. $\tilde Q(w)$ is analytic in $|w|<r_0$, and given by Pollaczek's formula (see e.g.~\cite{Abate1993, Cohen1982}).
In our discrete setting, we shall first derive a useful expression for $\tilde{Q}(w)$.

\begin{lemma}
For any $\eps>0$ with $1+\eps<r_0$,
\begin{equation} \label{e111}
\tilde Q(w)=\exp\Big(\frac{1}{2\pi i}\,\int_{|z|=1+\eps}\,{\rm ln}\Big(\frac{w-z}{1-z}\Big)\,\frac{(z^{s_n}-\tilde A(z))'}{z^{s_n}-\tilde A(z)}\,\dd z\Big)
\end{equation}
holds when $|w|<1+\eps$. 
\end{lemma}
\begin{proof}
We shall establish \eqref{e111} for any $w\in(1,1+\eps)$, and then the full result follows from analyticity of $\tilde{Q}(w)$ and of
\begin{equation*}
{\rm ln}\Big(\frac{w-z}{1-z}\Big)={\rm ln}\Big(\frac{1-w/z}{1-1/z}\Big)={-}\,\sum_{k=1}^{\infty}\,\frac1k\,\Big(\Big(\frac{w}{z}\Big)^k-\Big(\frac1z\Big)^k\Big)
\end{equation*}
in $w$, $|w|<1+\eps$ for any $z$ with $|z|=1+\eps$.

Our starting point is the formula, see \cite{Boudreau1962},
\begin{equation} \label{e113}
\tilde Q(w)=\frac{(s_n-\mu_A)(w-1)}{w^{s_n}-\tilde A(w)}\,\prod_{k=1}^{s_n-1}\,\frac{w-z_k}{1-z_k}
\end{equation}
that holds for all $w$, $|w|<r_0$, in which $z_1,\ldots,z_{s_n-1}$ are the $s_n-1$ zeros of $z^{s_n}-\tilde  A(z)$ in $|z|<1$. Fix $w\in(1,1+\eps)$. 
Then ${\rm ln}\,[(w-z)/(1-z)]$ is analytic in $z\in\mathbb{C}\backslash [1,w]$. 
It follows that
\begin{align}
 I_C &= \frac{1}{2\pi i}\,\int_{|z|=1+\eps}\,{\rm ln}\Big(\frac{w-z}{1-z}\Big)\,\frac{(z^{s_n}-\tilde A(z))'}{z^{s_n}-\tilde A(z)}\,\dd z \nonumber \\
&=~\sum_{k=1}^{s_n-1}\,{\rm ln}\Big(\frac{w-z_k}{1-z_k}\Big)+\frac{1}{2\pi i}\,\int_C\,{\rm ln}\Big(\frac{w-z}{1-z}\Big)\,\frac{(z^{s_n}-\tilde A(z))'}{z^{s_n}-\tilde A(z)}\,\dd z ,
\label{e114}
\end{align}
where $C$ is a contour encircling $[1,w]$ in the positive sense with none of the $z_k$'s in its interior. We let $\delta\in(0,\frac{w-1}{2})$ and we take $C$ the union of two line segments, from $1+\delta-i0$ to $w-\delta-i0$ and from $w-\delta+i0$ to $1+\delta-i0$, and two circles, of radius $\delta$ and encircling 1 and $w$ in positive sense.
A careful administration of the various contributions to the integral $I_C$ in \eqref{e114}, taking account of the branch cut $[1,w]$, yields
\begin{equation*}
I_C = {\rm ln }\left(\frac{(s_n-\mu_A)(w-1)}{w^s-\tilde A(w)}\right) + O(\delta\,{\rm ln}\, \delta ).
\end{equation*}
Using this in \eqref{e113} and letting $\delta \downarrow 0$, we get \eqref{e111} for $w\in(1,1+\varepsilon)$ and the proof is complete.
\end{proof} 

Using $\mathbb{P}(Q=0)=\tilde Q(0)$, $\mu_Q=\tilde Q'(1)$ and $\sigma_Q^2  =  \tilde Q''(1)+\tilde Q'(1)-(\tilde Q'(1))^2$,   it follows by straightforward manipulations that
\begin{align} \label{e6}
\mathbb{P}(Q=0)&=\exp\,\Big[\frac{1}{2\pi i}\,\int_{|z|=1+\eps}\,\ln\Big(\frac{z}{z-1}\Big)\frac{(z^{s_n}-\tilde A(z))'}{z^{s_n}-\tilde A(z)}\,\dd z\Big] , \\
\label{e7}
\mu_Q&=\frac{1}{2\pi i}\,\int_{|z|=1+\eps}\,\frac{1}{1-z}~\frac{(z^{s_n}-\tilde A(z))'}{z^{s_n}-\tilde A(z)}\,\dd z ,\\
\label{e8}
\sigma_Q^2 &= \frac{1}{2\pi i}\,\int_{|z|=1+\eps}\,\frac{-z}{(1-z)^2}~\frac{(z^{s_n}-\tilde A(z))'}{z^{s_n}-\tilde A(z)}\,\dd z .
\end{align}
Because $s_n$ appears directly in expressions \eqref{e6}-\eqref{e8}, we will be conducting our analysis with respect to $s_n$ rather than $n$. Note that this has no consequences for our results on the convergence speed of the performance metrics, since $s_n = O(n)$. Furthermore, we will omit the index $n$ when describing the capacity $s_n$ in the remainder of the chapter for brevity.

We next discuss in more detail the family of heavy-traffic scalings considered in this chapter, which combines two features. First, we have assumed that
 $A^{(n)}_j$ is in distribution equal to the sum of work generated by all sources, $A_{1,j}+...+A_{n,j}$, where the $A_{i,k}$ are for all $i$ and $k$ i.i.d.~copies of a random variable $X$, of which the pgf $\tilde X(z)=\sum_{k=0}^{\infty}\mathbb{P}(X=k)z^k$ has radius of convergence $r>1$, and
\begin{equation*}
0< \E[A^{(n)}] =n\mu = n \tilde X'(1)<s_n .
\end{equation*}
Hence
\begin{equation} \label{e10}
\vartheta:=\frac{n}{s_n}\in(0,1/\mu) .
\end{equation}
Second, we scale the system according to \eqref{bb}, for which we assume that
\begin{equation} \label{e11}
\rho_{s_n} =\vartheta\,\mu =1-\frac{\gamma}{s_n^\eta}
\end{equation}
in which $\gamma>0$ is bounded away from 0 and $\infty$ as $s_n\to  \infty$.
In the remainder of this chapter, we will omit the subscript in $s_n$. 
The condition that $\mathbb{P}(A^{(n)}=k)>0$ for some $k>s$ holds when the degree $d$ of $\tilde X(z)$ (with $d=\infty$ if $\tilde X(z)$ is not a polynomial) is such that $nd>s$.

To avoid certain complications when applying the saddle point method, we further assume that
\begin{equation} \label{e12}
|\tilde X(z)|<\tilde X(r_1) ,~~~~~~|z|=r_1\,,~~z\neq r_1 ,
\end{equation}
for any $r_1\in(0,r)$. This implies that $r_0$ is the unique zero of $z^s-\tilde A(z)$ on $|z|=r_0$.
This condition is related to Cram\'er's condition, see \cite[pp.~189 and 355]{Asmussen2003}, and it has also been used in \cite{relaxation}.
Condition \eqref{e12} holds when the set of all $j=0,1,\ldots$ such that $\mathbb{P}(X=k)>0$ is not contained in an  arithmetic progression with a ratio larger than one (see also \cite{rouche}).

\section{Non-standard saddle point method}\label{spSec}
\noindent
We illustrate our saddle point method for $\mu_Q$.
As a first step, we bring  \eqref{e7} in a form which is amenable to saddle point analysis.
\begin{lemma}
\begin{equation} \label{e18}
\mu_Q  =  \frac{s}{2\pi i}\,\int_{|z|=1+\eps}\,\frac{g'(z)}{z-1}~\frac{\exp(s\,g(z))}{1-\exp(s\,g(z))}\,\dd z
\end{equation}
with
\begin{equation} \label{e15}
g(z)={-}{\rm ln}\,z+\vartheta\,{\rm ln}(\tilde X(z)) .
\end{equation}
\end{lemma}
\begin{proof}
With $\tilde A(z)=\tilde X^n(z)$,
\begin{align} \label{e13}
\frac{(z^s-\tilde A(z))'}{z^s-\tilde A(z)} & =  \frac{s\,z^{s-1}-n\,\tilde X'(z)\,\tilde X^{n-1}(z)}{z^s-\tilde X^n(z)} \nonumber \\
& =  \frac{s}{z}-\frac{s}{z}\,\Big(\frac{n}{s}~\frac{z\,\tilde X'(z)}{\tilde X(z)}-1\Big)\,\frac{z^{-s}\,\tilde X^n(z)}{1-z^{-s}\,\tilde X^n(z)} .
\end{align}
Write
$
z^{-s}\,\tilde X^n(z)=\exp(s\,g(z))$.
Noting that
\begin{equation} \label{e16}
\frac{1}{2\pi i}\,\int_{|z|=1+\eps}\,\frac{s}{z}~\frac{1}{1-z}\,\dd z=0 ,
\end{equation}
and that
\begin{equation} \label{e17}
g'(z)=\frac1z\,\Big(\vartheta\,\frac{z\,\tilde X'(z)}{\tilde X(z)}-1\Big) ,
\end{equation}
gives \eqref{e18}. \end{proof}

Let us now explain how the standard saddle point method can be applied to \eqref{e18}.
Since
\begin{equation} \label{e19}
g(1)=g(r_0)=0~;~~~~~~g(z)<0\,,~~1<z<r_0 ,
\end{equation}
and by strict convexity of
\begin{equation*}
z^{-s}\,\tilde X^n(z)=z^{-s}\tilde A(z)=\sum_{k=0}^{\infty}\,a_k\,z^{k-s} ,~~~~~~z\in(0,r) ,
\end{equation*}
 $g(z)$ has a unique minimum on $[1,r_0]$. This minimum is found by solving $z\in[1,r_0]$ from $g'(z)=0$, and this yields the equation
\begin{equation} \label{e21}
\tilde X(z)=\vartheta\,z\,\tilde X'(z) .
\end{equation}
Denote the solution $z\in(1,r_0)$ of \eqref{e21} by $z_{\rm sp}$, and observe that $z_{\rm sp}$ is a saddle point of $g(z)$, explaining the notation. Thus, the saddle point method can be used for the integral in \eqref{e18} by taking $1+\eps=z_{\rm sp}$.

In the case that $\vartheta=n/s$ is bounded away from $1/\mu$ as $s\to  \infty$, we have that the minimum value of $g(z)$, $1\leq z\leq r_0$, is negative and bounded away from 0. Furthermore, $z_{\rm sp}$ is bounded away from 1, and the saddle point method can be applied in the classical way by replacing
\begin{equation*}
\frac{\exp(s\,g(z))}{1-\exp(s\,g(z))}~~~~{\rm by}~~~~\exp(s\,g(z)) ,
\end{equation*}
at the expense of an exponentially small relative error, and performing an expansion of $g'(z)/(z_{\rm sp}-1)=d_1(z-z_{\rm sp})+O((z-z_{\rm sp})^2)$ with $d_1=g''(z_{\rm sp})/(z_{\rm sp}-1)\neq 0$.
 Using that $g(z^{\ast})=(g(z))^{\ast}$, where the $^*$ denotes complex conjugation, it can be shown that
\begin{equation} \label{e23}
\mu_Q=\frac{\exp(s\,g(z_{\rm sp}))}{(z_{\rm sp}-1)^2\,\sqrt{2\pi s\,g''(z_{\rm sp})}}\,(1+O(s^{-1})) .
\end{equation}

We next explain why the standard saddle point method does not work for the heavy-traffic scaling considered in this chapter. 
Since we operate in \eqref{e11},
 $\vartheta\mu\to  1$ as $s\to  \infty$, and
\begin{align} \label{e24}
z_{\rm sp}-1&=\frac{\gamma}{a_2\,s^\eta}+O(s^{-2\eta}) ,\\
 \label{e25}
g(z_{\rm sp})&=\frac{-\gamma^2}{2a_2s^{2\eta}}+O(s^{-3\eta}) ,\\
 \label{e26}
g''(z_{\rm sp})&=a_2+O(s^{-\eta}) ,
\end{align}
where
\begin{equation} \label{e27}
a_2=\frac{\sigma^2}{\mu}-\frac{\gamma}{s^\eta}\,\Big(\frac{\sigma^2}{\mu}-1\Big) .
\end{equation}
Hence, $\exp(sg(z))$ near $z=z_{\rm sp}$ is (as $s\to  \infty$):
 vanishingly small when $\eta\in(0,1/2)$,
 bounded away from 1, but non-negligible when $\eta=1/2$,
and tending to 1 when $\eta\in(1/2,\infty)$.
Furthermore, $(z-1)^{-1}$ in \eqref{e18} is unbounded near $z=z_{\rm sp}$ as $s\to  \infty$. Therefore, an adaptation of the standard saddle point method is required, and the resulting asymptotic form of $\mu_Q$ will deviate significantly from the standard case \eqref{e23}. In particular, since $z_{\rm sp}\to  1$, this asymptotic form will contain information from $X(z)$ at $z=1$, rather than at a point away from 1 as is the case in \eqref{e23}.

The required adaptation of the saddle point method is modeled after a device developed in \cite[Sec.~5.12]{debruijn}. We use a substitution $z=z(v)$ in \eqref{e18} with real $v$ and $z(0)=z_{\rm sp}$ such that for sufficiently small $v$,
\begin{equation} \label{e29}
g(z(v))=g(z_{\rm sp})-\tfrac12\,v^2\,g''(z_{\rm sp}) .
\end{equation}
This is feasible, since
\begin{equation} \label{e30}
g(z)=g(z_{\rm sp})+\tfrac12\,g''(z_{\rm sp})(z-z_{\rm sp})^2\Big(1+\frac{g'''(z_{\rm sp})}{3g''(z_{\rm sp})}\,(z-z_{\rm sp})+...\Big)
\end{equation}
with $g''(z_{\rm sp})$ positive and bounded away from 0 as $s\to  \infty$. Hence, $z(v)$ can be found for small $v$ by inverting the equation
\begin{equation} \label{e31}
(z-z_{\rm sp})\Big(1+\frac{g'''(z_{\rm sp})}{3g''(z_{\rm sp})}\,(z-z_{\rm sp})+...\Big)^{1/2}=iv .
\end{equation}
By Lagrange's inversion theorem \cite{debruijn}, there is a $\delta>0$ (independent of $s$) such that
\begin{equation} \label{e32}
z(v)=z_{\rm sp}+iv+\sum_{k=2}^{\infty}\,c_k(iv)^k ,~~~~~~|v|<\delta ,
\end{equation}
with real coefficients $c_k$ (since $g(z)$ is real for real $z$) and
\begin{equation} \label{e33}
c_2={-}\,\frac{g'''(z_{\rm sp})}{6g''(z_{\rm sp})} .
\end{equation}
Thus
\begin{equation} \label{e34}
z(v)=z_{\rm sp}+iv-c_2\,v^2+O(v^3) ,~~~~~~|v|\leq\tfrac12\,\delta ,
\end{equation}
where the order term holds uniformly in $s$. The uniformity statement follows from an inspection of the usual argument
by which Lagrange's theorem is proved, noting that the inversion in \eqref{e29} with $g$ as in \eqref{e15} is considered for $\vartheta\to   1/\mu$, $z_{\rm sp}\to   1$ with radius
of convergence $r$ away from $1$.

By \eqref{e12} we can restrict the integration in \eqref{e18} to a fixed but arbitrarily small subset of $|z|=z_{\rm sp}$ near $z=z_{\rm sp}$, at the expense of an exponentially small error. Furthermore, by Cauchy's theorem and again at the expense of an exponentially small error, the integration path can be deformed in accordance with the transformation in \eqref{e29}--\eqref{e34}. Set
\begin{equation} \label{e35}
q(v)=g(z_{\rm sp})-\tfrac12\,v^2\,g''(z_{\rm sp})
\end{equation}
and note that from \eqref{e29},
\begin{equation*} 
g'(z(v))\,z'(v)={-}v\,g''(z_{\rm sp}) .
\end{equation*}
Then substituting $z=z(v)$ in \eqref{e18},  $\mu_Q$ is given with exponentially small error by
\begin{equation*}
\frac{s}{2\pi i}\,\int_{-\frac12\delta}^{\frac12\delta}\,\frac{g'(z(v))}{z(v)-1}~\frac{\exp(s\,g(z(v)))}{1-\exp(s\,g(z(v)))}z'(v)\,\dd v,
\end{equation*}
which gives the following result.
\begin{lemma} \label{lemma2} The mean stationary congestion level is given with exponentially small error by
\begin{equation} \label{e37}
\mu_Q =~\frac{-s}{2\pi i}\,g''(z_{\rm sp})\,\int_{-\frac12\delta}^{\frac12\delta}\,\frac{v}{z(v)-1}~\frac{\exp(s\,q(v))}{1-\exp(s\,q(v))}\,\dd v .
\end{equation}
\end{lemma}

In a similar fashion we get that $\mathbb{P}(Q=0)$ and $\sigma_Q^2$, see \eqref{e6} and \eqref{e8}, are given, both with exponentially small error, by
\begin{equation} \label{e39}
\frac{-s}{2\pi i}\,g''(z_{\rm sp})\,\int_{-\frac12\delta}^{\frac12\delta}\,v\,{\rm ln}\Big(\frac{z(v)}{z(v)-1}\Big)\frac{\exp(s\,q(v))}{1-\exp(s\,q(v))}\,\dd v
\end{equation}
and
\begin{equation} \label{e38}
\frac{-s}{2\pi i}\,g''(z_{\rm sp})\,\int_{-\frac12\delta}^{\frac12\delta}\,\frac{v\,z(v)}{(z(v)-1)^2}~\frac{\exp(s\,q(v))}{1-\exp(s\,q(v))}\,\dd v,
\end{equation}
respectively.

\section{Heavy-traffic limits for the mean congestion level} \label{sec3}
In this section we apply the non-standard saddle point method explained in Section \ref{spSec} to the Pollaczek integral representation for the mean stationary congestion level $\mu_Q$. In Section \ref{subsec3.1} we first derive an integral representation for the leading order behavior of $\mu_Q$  with a relative error of order $O(s^{-1})$, which serves as a heavy-traffic approximation in the regime $\rho_s=1-\gamma/s^\eta$ with $\eta>0$. We also consider separately the cases of moderate heavy traffic ($\eta\in(0,1/2)$) and extreme heavy traffic ($\eta\in(1/2,\infty)$), for which the integral representation leads to vastly different alternative expressions. We find that $\mu_Q\to   0$ more rapidly than any power of $1/s$ when $\eta\in(0,1/2)$. When $\eta\geq 1/2$ the saddle point method yields an integral representation with relative error $O(s^{-\min(1,\eta)})$.
In Section \ref{subsec3.2} we specialize this general result to the CLT case $\eta=1/2$, and make a connection with existing results.

\subsection{Leading order behavior in integral form} \label{subsec3.1}

 \begin{theorem}\label{mainthm}
The mean stationary congestion level is given by
\begin{equation} \label{e4.1}
 \mu_Q=\frac{2}{\pi}\,\sigma\,\sqrt{\dfrac{s}{2\mu}}\,\int_0^{\infty}\,\frac{t^2}{d^2(s)+t^2}~\frac{\exp({-}d^2(s)-t^2)}{1-\exp({-}d^2(s)-t^2)}\,\dd t\,\left(1+O({s^{{-}\min(1,\eta)}})\right)
\end{equation}
with $
d^2(s) = s^{1-2\eta}\gamma^2\mu/(2\sigma^2)$.
\end{theorem}

\begin{proof}
According to Lemma \ref{lemma2}, $\mu_Q$ is given with exponentially small error by \eqref{e37} with $q(v)$ given in \eqref{e35}. Since $z({-}v)=z^{\ast}(v)$ for real $v$, we have
\begin{align}
\frac{v}{z(v)-1}+\frac{-v}{z({-}v)-1} &= {-}2iv\,\frac{{\rm Im}(z(v))}{|z(v)-1|^2}\nonumber\\
&=\frac{-2iv^2+O(v^4)}{(z_{\rm sp}-1)^2+v^2-2c_2(z_{\rm sp}-1)\,v^2+O(v^4)} \nonumber \\
&=\frac{-2iv^2\left(1+O(v^2)\right)}{(z_{\rm sp}-1)^2+v^2-2c_2(z_{\rm sp}-1)\,v^2},
\label{e40}
\end{align}
for ${-}\tfrac{1}{2} \delta \leq v \leq \tfrac{1}{2} \delta$.
where \eqref{e34} and $c_k\in\mathbb{R}$ have been used. Using \eqref{e40} in \eqref{e37} and extending the integration range from $[{-}\tfrac12\delta,\tfrac12\,\delta]$ to $({-}\infty,\infty)$ while using symmetry of $q(v)$, we get that $\mu_Q$ is given with exponentially small error by
\begin{align} \label{e41}
\frac{s\,g''(z_{\rm sp})}{\pi}\,\int_0^{\infty}\,\frac{v^2\left(1+O(v^2)\right)}{(z_{\rm sp}-1)^2+v^2-2c_2(z_{\rm sp}-1)\,v^2}\frac{\exp(s\,q(v))}{1-\exp(s\,q(v))}\dd v .
\end{align}
With
\begin{equation} \label{e42}
B=\exp(s\,g(z_{\rm sp})) ,~~~~~~\alpha =g''(z_{\rm sp}),
\end{equation}
Equation \eqref{e41} takes the form
\begin{equation} \label{e43}
\frac{s\alpha }{\pi}\,\int_0^{\infty}\,\frac{v^2\left(1+O(v^2)\right)}{(z_{\rm sp}-1)^2+v^2-2c_2(z_{\rm sp}-1)\,v^2} \cdot \frac{B\,\exp({-}\tfrac12\,s\,\alpha \,v^2)}{1-B\,\exp({-}\tfrac12\,s\,\alpha \,v^2)}\dd v .
\end{equation}
Since $(z_{\rm sp} - 1)^2 = (\gamma/a_2)^2s^{{-}2\eta} + O(s^{-4\eta})$, see \eqref{e24}, the integrand in \eqref{e43} in leading order has the form
\begin{equation*}
\frac{B\,v^2\,\exp(-s\,D\,v^2)}{(v^2+C\,s^{-2\eta})(1-B\exp({-}s\,D\,v^2))},
\end{equation*}
and this is reminiscent of the integrand in \cite[Eq.~(5.12.3)]{debruijn} for the case  $\kappa=2\eta$. Proceeding as in \cite[Sec.~5.12]{debruijn}, the substitution $v=t\sqrt{{2}/(s\alpha )}$ brings \eqref{e43} into the form
\begin{equation} \label{e44}
\frac{2}{\pi}\sqrt{\tfrac12 s\alpha }\int_0^{\infty}\frac{t^2(1+O(t^2/s))}{\tfrac12 s\alpha (z_{\rm sp}-1)^2+t^2-2c_2(z_{\rm sp}-1)t^2} \,\frac{B\,\exp({-}t^2)}{1-B\,\exp({-}t^2)}\dd t .
\end{equation}
From \eqref{e24}--\eqref{e27} and \eqref{e42},
\begin{align}
\frac{2}{\pi}\,\sqrt{\frac{s\alpha }{2}} &= \frac{2}{\pi}\,\sigma_X\,\sqrt{\frac{s}{2\,\mu}}\,(1+O(s^{-\eta})),\label{y45}\\
\tfrac12\,s\,\alpha \,(z_{\rm sp}-1)^2 &= d^2(s) + O(s^{1-3\eta}),\label{y46}\\
2\,c_2(z_{\rm sp}-1) &= O(s^{-\eta}),\label{y47}\\
s\,g(z_{\rm sp}) &= -d^2(s) + O(s^{1-3\eta}),\label{y48}
\end{align}
where
\begin{equation} \label{y49}
d^2(s) = \frac{b_0^2}{s^{2\eta-1}},\quad b_0^2 := \frac{\gamma^2\mu}{2\,\sigma^2}.
\end{equation}
In the case that $2\eta-1<0$, we have that $\tfrac12\,s\,\alpha \,(z_{\rm sp}-1)^2 \to   \infty$ and that
\begin{equation} \label{y50}
B = \exp(s\,g(z_{\rm sp})) = O(\exp({-}b^2s^{1-2\eta}))
\end{equation}
for any $b\in(0,b_0)$. From \eqref{e44} it then follows that $\mu_Q = O(\exp({-}b^2\,s^{1-2\eta}))$ for any $b\in(0,b_0)$.
In the case that $2\,\eta-1\geq 0$, we have that $d^2(s)$ is bounded, and using that $1/s^{3\eta-1} = O(d^2(s)/s^\eta)$, we get
\begin{align*}
\tfrac12\,s\,\alpha \,(z_{\rm sp}-1)^2 + t^2-2\,c_2\,(z_{\rm sp}-1)\,t^2 
&= d^2(s) + t^2 + O\left(s^{-\eta}\,(d^2(s)+t^2)\right) \nonumber\\
&= \big(d^2(s)+t^2\big)\left(1+O(s^{-\eta})\right).
\end{align*}
Hence, in this case,
\begin{equation} 
\frac{t^2(1+O(t^2/s))}{\tfrac12 s\,\alpha (z_{\rm sp}-1)^2+t^2-2c_2(z_{\rm sp}-1)t^2} = \frac{t^2}{d^2(s)+t^2}\left(1+O(s^{-\eta})+O(t^2/s)\right).\label{y52}
\end{equation}
Furthermore,
\begin{align*}
1-B\,\exp(-t^2) &= 1-\exp({-}d^2(s)-t^2)\,\left(1+d^2(s)\,O(s^{-\eta})\right)\nonumber\\
&=(1-\exp({-}d^2(s)-t^2))\,\Big(1+\frac{d^2(s)}{\exp(d^2(s)+t^2)-1}O(s^{-\eta})\Big)\nonumber\\
&= (1-\exp({-}d^2(s)-t^2))\,(1+O(s^{-\eta})),
\end{align*}
It follows therefore that
\begin{equation} \label{y56}
\frac{B\,\exp({-}t^2)}{1-B\,\exp({-}t^2)} = \frac{\exp({-}d^2(s)-t^2)}{1-\exp({-}d^2(s)-t^2)}\,(1+O(s^{-\eta})).
\end{equation}
Combining the three items \eqref{y45}, \eqref{y52} and \eqref{y56}, we obtain for \eqref{e44} the result
\begin{equation*} 
\frac{2}{\pi}\,\sigma\,\sqrt{\frac{s}{2\,\mu}}  \int_0^{\infty}\frac{t^2}{d^2(s)+t^2} \cdot \frac{\exp({-}d^2(s)-t^2)}{1-\exp({-}d^2(s)-t^2)}\dd t
\left(1+O(s^{-\eta})+O(s^{-1})\right),
\end{equation*}
and this gives \eqref{e4.1}.
\end{proof}

Theorem \ref{mainthm} gives the leading-order behavior of $\mu_Q$ as $s\to  \infty$ with a relative error of $O(s^{{-}\min(1,\eta)})$. By considering in more detail the integral expressions, we obtain the following result, describing two different heavy-traffic behaviors.

\begin{proposition}\label{prop1}
If $\eta\in(0,1/2)$ the mean congestion level satisfies
\begin{equation*} 
\mu_Q=O\left(\exp(-b^2s^{1-2\eta})\right),
\end{equation*}
for any $b\in (0,b_0)$. If $\eta\in[1/2,\infty)$ the mean congestion level is given by
\begin{equation*} 
\mu_Q = s^\eta\,\frac{\sigma^2}{2\mu\gamma}\,\left(1+O(s^{\max(1/2-\eta,-1)})\right).
\end{equation*}
\end{proposition}
The first assertion in Proposition \ref{prop1} follows from the observation in \eqref{y50}, together with \eqref{e44}. The second assertion is based on a connection between the integral in Theorem \ref{mainthm} and the Riemann zeta function, which is explained in the next subsection.

\subsection{Classical heavy traffic and the Gaussian random walk}
\label{subsec3.2}
We now build on Theorem \ref{mainthm} to obtain further results for the  classical heavy traffic case $\eta=1/2$,
for which we know from \cite[Thm.~3]{Sigman2011b} that the rescaled congestion process converges under \eqref{bb1} to a reflected Gaussian random walk. The latter is defined as
 $(S_\beta(k))_{k\geq 0}$ with $S_\beta(0)=0$ and
\begin{equation*}
S_\beta(j)=Y_1+\ldots+Y_j
\end{equation*}
with $Y_1,Y_2,\ldots$ i.i.d.~copies of a normal random variable with mean $-\beta$ and variance 1.
Assume $\beta>0$ (negative drift), and denote the all-time maximum of this random walk by ${M}_\beta$.

 Denote by $Q^{(s)}_\infty$ the stationary congestion level for a fixed $s$ (that arises from taking
 $j\to   \infty$ in \eqref{lind}), and remember that we have assumed $\vartheta=n/s$ fixed.
Then, using $\rho_s=1-\gamma/\sqrt{s}$, with
 \begin{equation}\label{gammachoice}
 \gamma=\frac{\beta\sigma}{\mu\sqrt{\vartheta}},
\end{equation}
the spatially-scaled stationary congestion levels reach the limit
$Q^{(s)}_\infty/(\sigma\sqrt{n}) \Rightarrowd {M}_\beta$ as $s,n\to  \infty$ (see \cite{Jelenkovic2004,Sigman2011a,Sigman2011b}). From \cite[Thm.~4]{Sigman2011b} we then know that under the standard heavy-traffic scaling \eqref{bb1}
\begin{equation}
   \frac{\mathbb{E}[Q^{(s)}_\infty]}{\sigma\sqrt{n}}\to   \mathbb{E}[{M}_\beta], \quad {\rm as} \ s,n\to  \infty,
\end{equation}
   from which it follows that
   \begin{equation} \label{e48}
\mu_Q\approx \sigma\sqrt{n}\ \mathbb{E}[M_\beta].
\end{equation}
The random variable ${M}_\beta$ was studied in  \cite{Chang1997,Janssen2006}. In particular, \cite[Thm.~2]{Janssen2006} yields, for $\beta<2\sqrt{\pi}$,
\begin{equation*}
\mathbb{E}[{M}_\beta]= \frac{1}{2\beta}+\frac{\zeta(1/2)}{\sqrt{2\pi}}+\frac{\beta}{4}+\frac{\beta^2}{\sqrt{2\pi}}\sum_{r=0}^{\infty}\frac{\zeta(-1/2-r)}{r!(2r+1)(2r+2)}\left(\frac{-\beta^2}{2 }\right)^r,
\end{equation*}
where $\zeta$ denotes the Riemann zeta function, which is defined as, see (1.26).
Hence, for small values of $\beta$,
   \begin{equation} \label{estimate}
\mu_Q\approx \sigma\sqrt{n}\ \mathbb{E}[M_\beta] \approx \frac{\sigma\sqrt{n}}{2\beta} = \sqrt{s}\,\frac{\sigma^2}{2\mu\gamma}.
\end{equation}
We will now show how the approximation \eqref{estimate} follows from Theorem \ref{mainthm}, and also how similar steps give rise to Proposition \ref{prop1}.

Consider the integral
\begin{equation} \label{e49}
G_0(b)=G_1(b)-G_2(b)=\int_0^{\infty}\,\frac{t^2}{b^2+t^2}~\frac{\exp({-}b^2-t^2)}{1-\exp({-}b^2-t^2)}\,\dd t ,
\end{equation}
where $b>0$ and
\begin{equation} \label{e50}
G_1(b)=\int_0^{\infty}\,\frac{\exp({-}b^2-t^2)}{1-\exp({-}b^2-t^2)}\,\dd t\,,~~~~G_2(b)=\int_0^{\infty}\,\frac{b^2}{b^2+t^2}~\frac{\exp({-}b^2-t^2)}{1-\exp({-}b^2-t^2)}\dd t .
\end{equation}
We have, as in \cite[Sec.~2]{Janssen2006},
\begin{align} 
G_1(b) & =  \sum_{k=0}^{\infty}\:\int_0^{\infty}\,\exp({-}(k+1)(b^2+t^2))\,\dd t \nonumber \\
& =  \frac{\sqrt{\pi}}{2}\,\sum_{k=0}^{\infty}\,\frac{\ee^{-(k+1)b^2}}{\sqrt{k+1}} = \frac{\sqrt{\pi}}{2}\,\ee^{-b^2}\,\Phi(\ee^{-b^2},1/2,1) \nonumber \\
& =  \frac{\pi}{2b}+\frac{\sqrt{\pi}}{2}\,\sum_{r=0}^{\infty}\,\zeta(\tfrac12-r)\,\frac{({-}1)^r\,b^{2r}}{r!} , \label{e51}
\end{align}
where the last identity holds when $0<b<\sqrt{2\pi}$ and $\Phi(z,s,v)$ is Lerch's transcendent, which is defined as, see \cite[Eq.~25.14.1]{NIST},
\begin{equation*}
\Phi(z,s,v) = \sum_{n=0}^\infty \frac{z^n}{(v+n)^s}, \qquad \text{for }v\neq 0,{-}1,{-}2,\ldots, \ |z|<1; \, \Re s > 1,\ |z|=1.
\end{equation*}
As to $G_2(b)$, we make a connection with the complementary error function
\begin{equation*}
{\rm erfc}(z)=\frac{2}{\sqrt{\pi}}\,\int_z^{\infty}\,\ee^{-t^2}\,\dd t=\frac{2}{\pi}\,\ee^{-z^2}\,\int_0^{\infty}\,\frac{\ee^{-z^2t^2}}{1+t^2}\dd t ,
\end{equation*}
see \cite[Secs.~7.2 and 7.7.1]{NIST}. We thus compute
\begin{align} \label{e53}
G_2(b) & =  \sum_{k=0}^{\infty}\,\ee^{-(k+1)b^2}\,\int_0^{\infty}\,\frac{b^2}{b^2+t^2}\,\ee^{-(k+1)t^2}\dd t \nonumber \\
& =  \frac{\pi}{2}\,b\,\sum_{k=0}^{\infty}\,{\rm erfc}(b\,\sqrt{k+1}) .
\end{align}
From \cite[Eq.~(4.3) \& (4.23)]{Janssen2006},
\begin{equation} \label{e54}
\sum_{n=1}^{\infty}\,\frac{1}{\sqrt{2\pi}}\,\int_{\beta\sqrt{n}}^{\infty}\,\ee^{-x^2/2}\,dx= \frac{1}{2\beta^2}-\frac14-\frac{1}{\sqrt{2\pi}}\,\sum_{r=0}^{\infty} \frac{\zeta({-}1/2-r)({-}1/2)^r} {r!\,(2r+1)}\,\beta^{2r+1}
\end{equation}
in which $0<\beta<2\sqrt{\pi}$. Taking $\beta=b\,\sqrt{2}$ in \eqref{e54}, we get
\begin{equation} \label{e55}
G_2(b)=\frac{\pi}{4b}-\frac{\pi}{4}\,b-\sqrt{\pi}\,\sum_{r=0}^{\infty}\,\frac{\zeta({-}1/2-r)({-}1)^r\,b^{2r+2}}{r!\,(2r+1)}
\end{equation}
when $0<b<\sqrt{2\pi}$. The two results in \eqref{e51} and \eqref{e55} can be combined, as in \cite[Sec.~\ref{sec4}]{Janssen2006}, and this yields
\begin{equation} \label{e56}
G_0(b)=\frac{\pi}{4b}+\frac{\pi}{4}\,b+\frac{\sqrt{\pi}}{2}\,\zeta(1/2)+\sqrt{\pi}\,\sum_{r=0}^{\infty}\frac{\zeta({-}1/2-r)({-}1)^r\,b^{2r+2}}{r!\,(2r+1)(2r+2)}
\end{equation}
when $0<b<\sqrt{2\pi}$. 

Using \eqref{e56} in \eqref{e48}, we find that the leading order behavior of $\mu_Q$ is given as
\begin{equation} \label{e57}
\sigma_X\,\sqrt{\dfrac{s}{2\mu}}\,\left[\frac{1}{2b_0}+\frac{b_0}{2}+\frac{\zeta(1/2)}{\sqrt{\pi}}+\frac{2}{\sqrt{\pi}}\,\sum_{r=0}^{\infty}\,\frac{\zeta({-}1/2-r)({-}1)^r b_0^{2r+2}} {r!\,(2r+1)(2r+2)}\right]
\end{equation}
with relative error of $O(s^{-1/2})$ in which $b_0$ is given by \eqref{y49}. The expression \eqref{e57} is exactly equal to the right-hand side of \cite[Eq.~(4.25)]{Janssen2006} times $\sqrt{s}$ when we take there $\sigma=\mu=1$ and $\beta=b_0\,\sqrt{2}$.
Notice that, with $\gamma$ as in \eqref{gammachoice},
\begin{equation*}
\sigma\,\sqrt{\dfrac{s}{2\mu}}\frac{1}{2b_0}=\frac{\sigma\sqrt{n}}{2\beta},
\end{equation*}
which confirms the approximation \eqref{estimate}.

According to Theorem \ref{mainthm}, we have for $\eta\geq 1/2$,
\begin{equation*}
\mu_Q = \frac{2}{\pi}\,\sigma\,\sqrt{\frac{s}{2\,\mu}}G_0(d(s))\,\left(1+O(s^{{-}\min(1,\eta)})\right).
\end{equation*}
When $\eta=1/2$, so that $d(s) = b_0$ is independent of $s$, the series representation for $G_0$ in \eqref{e56} can be used, as long as $b_0\in(0,\sqrt{2\pi})$. When $\eta>1/2$, we have that $d(s) = b_0/s^{\eta-1/2}\to   0$ as $s\to  \infty$, and so this series representation can be used when $s$ is large enough. We then have from \eqref{e56} and $b_0^2 = \gamma^2\mu/2\,\sigma^2$, while replacing the whole series at the right-hand side by $O(b^2)$, for $\mu_Q$ the leading order behavior
\begin{equation} \label{y62}
s^\eta\left[\frac{\sigma^2}{2\,\gamma\,\mu}+\frac{\sigma\,\zeta(1/2)}{\sqrt{2\,\pi\,\mu}}\,\frac{1}{s^{\eta-1/2}}+\frac{1}{4}\,\gamma\,\frac{1}{s^{2\eta-1}}+O(s^{3/2-3\eta})\right]
\end{equation}
with relative error $O(s^{{-}\min(1,\eta)})$. Retaining the constant term $\sigma^2/(2\gamma\mu)$ and estimating the other terms between the brackets in \eqref{y62} as $O(s^{1/2-\eta})$, we get Proposition \ref{prop1}.
\section{More heavy-traffic results}\label{more}
In this section we apply the non-standard saddle point method to obtain several more heavy-traffic results. In Section \ref{subsec3.3} we derive refined heavy-traffic approximations for the mean congestion level by considering higher-order correction terms. In Section \ref{sec4} we derive the leading heavy-traffic behavior for the variance of the stationary congestion level, and in Section \ref{sec5} for the empty-system probability. To keep the developments tractable, we restrict Section \ref{subsec3.3} to $\eta=1/2$, and Section \ref{sec4} and Section \ref{sec5} to $\eta\in(0,1]$, although the same technique will work for all values $\eta>0$.

\subsection{Correction term for the mean congestion level for $\eta = 1/2$} \label{subsec3.3}
Our saddle point method not only establishes the leading-order heavy-traffic approximations, but also allows to derive refinements to these approximations. In this section we demonstrate how this works for the mean congestion level in the case $\eta=1/2$.

To obtain a refinement or correction term from \eqref{e44}, we must be more precise about the $O(s^{{-}\eta})$ terms that occur in the approximations in Section \ref{subsec3.1} for $\frac12\,s\,\alpha (z_{\rm sp}-1)^2$, $B$ and $\sqrt{s\,\alpha /2}$. When higher-order corrections are required, we should include higher-order terms in the approximations of these quantities, and be more specific about the $O(t^2/s)$ and $O(t^4/s)$  in the integrand in \eqref{e44}.

Let $g^{(i)}, \ i=1,2,...$~denote the $i^{\rm th}$ derivative of $g$ and define, see \eqref{e10} and \eqref{e15} with $\vartheta=(1-\gamma/s^\eta)\,\mu^{-1}$,
\begin{equation*} 
a_i=g^{(i)}(1);~~~~~~g(z)={-}{\rm ln}\,z+\vartheta\,{\rm ln}\,\tilde X(z) .
\end{equation*}
Dropping the $X$ from $\mu$ and $\sigma^2$ for brevity, we have
\begin{equation*} 
a_1={-}\,\frac{\gamma}{s^\eta} ,~~~~~~a_2=\frac{\sigma^2}{\mu}-\frac{\gamma}{s^\eta}\,\Big(\frac{\sigma^2}{\mu}-1\Big),
\end{equation*}
\begin{equation*} 
a_3={-}2+\Big(1-\frac{\gamma}{s^\eta}\Big)\Big(\frac{\tilde X'''(1)}{\tilde X'(1)}-3\tilde X''(1)+2(\tilde X'(1))^2\Big) .
\end{equation*}
For the purpose of finding a first-order correction term, we note that
\begin{align*} 
\alpha &=g''(z_{\rm sp})=a_2+(z_{\rm sp}-1)\,a_3+O(s^{-1}) ,\\
z_{\rm sp}-1&={-}\,\frac{a_1}{a_2}-\frac{a_3}{2a_2}\,\Big(\frac{a_1}{a_2}\Big)^2+O(s^{-3/2}) ,\\
c_2&={-}\,\frac{g'''(z_{\rm sp})}{6g''(z_{\rm sp})}={-}\,\frac{a_3}{6a_2}+O(s^{-1/2}) ,\\
g(z_{\rm sp})&={-}\,\frac{a_1^2}{2a_2}-\frac{a_3}{6a_2^3}\,a_1^3+O(s^{-2}) .
\end{align*}
This gives rise to
\begin{align} \label{e65}
\sqrt{\tfrac12\,s\,\alpha }&=\sigma\,\sqrt{\dfrac{s}{2\mu}}\,\Big(1+\frac{C_1}{\sqrt{s}}+O(s^{-1})\Big) ,\\
\tfrac12\,s\,\alpha (z_{\rm sp}-1)^2&=\frac{\gamma^2\,\mu}{2\sigma^2}+\frac{C_2}{\sqrt{s}}+O(s^{-1}) ,\\
\label{e67}
2c_2(z_{\rm sp}-1)&=\frac{C_3}{\sqrt{s}}+O(s^{-1}) ,\\
 \label{e68}
B=\exp(s\,g(z_{\rm sp}))&=\exp\Big({-}\,\frac{\gamma^2\,\mu}{2\sigma^2}\Big)\Big(1+\frac{C_4}{\sqrt{s}}+O(s^{-1})\Big) ,
\end{align}
with explicitly computable constants $C_1$, $C_2$, $C_3$, $C_4$. Remembering that $b_0^2=\gamma^2\mu/2\sigma^2$, see \eqref{y49}, we then get with errors of order $1/s$
\begin{align}
& \frac{t^2(1+O(t^2/s))}{\frac12\,s\,\alpha (z_{\rm sp}-1)^2+t^2-2c_2(z_{\rm sp}-1)\,t^2} \nonumber \\
& \qquad \qquad =~\frac{t^2}{b_0^2+t^2}-\frac{1}{\sqrt{s}}\,\Big((C_2+ b_0^2\,C_3)\,\frac{t^2}{(b_0^2+t^2)^2}-C_3\,\frac{t^2}{b_0^2+t^2}\Big) ,\label{e69}
\end{align}
and
\begin{equation} \label{e70}
\frac{B\,\exp({-}t^2)}{1-B\,\exp({-}t^2)}=\frac{\exp({-}b_0^2-t^2)}{1-\exp({-}b_0^2-t^2)}+\frac{C_4}{\sqrt{s}}~\frac{\exp({-}b_0^2-t^2)}{(1-\exp({-}b_0^2-t^2))^2} .
\end{equation}
Using \eqref{e65}, \eqref{e69} and \eqref{e70} in \eqref{e44} we get with an absolute error of order $1/\sqrt{s}$
\begin{align} \label{e71}
\mu_Q & =\frac{2}{\pi}\,\sigma\,\sqrt{\dfrac{s}{2\mu}}\,\Big(1+\frac{C_1}{\sqrt{s}}\Big)\nonumber \\
& \qquad\qquad \cdot \int_0^{\infty}\,\Big(\frac{t^2}{b_0^2+t^2}-\frac{1}{\sqrt{s}}\,\Big((C_2+b_0^2\,C_3)\,\frac{t^2}{(b_0^2+t^2)^2}-C_3\,\frac{t^2}{b_0^2+t^2}\Big)\Big) \nonumber \\
& \qquad\qquad\qquad \cdot~\Big(\frac{\exp({-}b_0^2-t^2)}{1-\exp({-}b_0^2-t^2)}+\frac{C_4}{\sqrt{s}}~\frac{\exp({-}b_0^2-t^2)}{(1-\exp({-}b_0^2-t^2))^2}\Big)\dd t \nonumber \\
& =\frac{2\sigma}{\pi}\,\sqrt{\dfrac{s}{2\mu}}\,G_0(b_0)\nonumber\\
& \qquad\qquad + ~\frac{2\sigma}{\pi}\,\sqrt{\dfrac{1}{2\mu}}\,\big((C_1+C_3)\,G_0(b_0)-(C_2+b_0^2\,C_3)\,G_3(b_0)+C_4\,G_4(b_0)\big) ,
\end{align}
where $G_0$ is as in \eqref{e49}, and
\begin{align} \label{e72}
G_3(b_0)&=\int_0^{\infty}\,\frac{t^2}{(b_0^2+t^2)^2}~\frac{\exp({-}b_0^2-t^2)}{1-\exp({-}b_0^2-t^2)}\,\dd t ,\\
\label{e73}
G_4(b_0)&=\int_0^{\infty}\,\frac{t^2}{b_0^2+t^2}~\frac{\exp({-}b_0^2-t^2)}{(1-\exp({-}b_0^2-t^2))^2}\,\dd t .
\end{align}
We shall express the integrals in \eqref{e72} and \eqref{e73} in terms of $\zeta$-functions. By partial integration
\begin{align} \label{e74}
G_3(b) & =  \frac12\,\int_0^{\infty}\,\frac{1}{b^2+t^2}~\frac{\exp({-}b^2-t^2)}{1-\exp({-}b_0^2-t^2)}\,\,\dd t \nonumber \\ 
&\qquad\qquad -~\int_0^{\infty}\,\frac{t^2}{b^2+t^2}~\frac{\exp({-}b^2-t^2)}{(1-\exp({-}b^2-t^2))^2}\,\dd t \nonumber \\
& =  \frac{1}{2b^2}\,G_2(b)-G_4(b) ,
\end{align}
see \eqref{e49} and \eqref{e73}. Since $G_2(b)$ is expressed in terms of $\zeta$-functions in \eqref{e55}, it is sufficient to consider $G_4(b)$.

As to $G_4(b)$,
\begin{equation*} 
G_4(b)=G_5(b)-G_6(b) ,
\end{equation*}
where
\begin{align*} 
G_5(b)&=\int_0^{\infty}\,\frac{\exp({-}b^2-t^2)}{(1-\exp({-}b^2-t^2))^2}\,\dd t ,\\
G_6(b)&=\int_0^{\infty}\,\frac{b^2}{b^2+t^2}~\frac{\exp({-}b^2-t^2)}{(1-\exp({-}b^2-t^2))^2}\,\dd t .
\end{align*}
We have, compare \eqref{e51},
\begin{align} \label{e78}
G_5(b) & =  \sum_{k=0}^{\infty}\,(k+1)\,\int_0^{\infty}\,\ee^{-(k+1)(b^2+t^2)}\,\dd t \nonumber \\[3.5mm]
& =  \frac{\sqrt{\pi}}{2}\,\ee^{-b^2}\,\Phi(\ee^{-b^2},{-}\tfrac12,1) =  \frac{\pi}{4b^3}+\frac{\sqrt{\pi}}{2}\,\sum_{r=0}^{\infty}\,\zeta({-}\tfrac12-r)\,\frac{({-}1)^r\,b^{2r}}{r!} ,
\end{align}
the last identity being valid when $0<b<\sqrt{2\pi}$. Next we have, compare \eqref{e53},
\begin{align*} 
G_6(b) & =  \sum_{k=0}^{\infty}\,(k+1)\,b^2\,\int_0^{\infty}\,\frac{\exp({-}(k+1)(b^2+t^2))}{b^2+t^2}\dd t \nonumber \\
& =  \frac{\pi}{2}\,b\,\sum_{k=0}^{\infty}\,(k+1)\,{\rm erfc}(b\,\sqrt{k+1}) .
\end{align*}
From \cite[Eq.~(5.4) \& (5.21)]{Janssen2006} we have
\begin{equation} \label{e80}
\sum_{n=1}^{\infty}\frac{n}{\sqrt{2\pi}}\int_{\beta\sqrt{n}}^{\infty}\ee^{-x^2/2}\,dx = \frac{3}{4\beta^4}-\frac{1}{24}-\frac{1}{\sqrt{2\pi}}\sum_{r=0}^{\infty}\frac{\zeta({-}3/2-r)({-}1/2)^r}{r!\,(2r+1)}\,\beta^{2r+1}
\end{equation}
when $0<\beta<2\sqrt{\pi}$. Taking $\beta=b\,\sqrt{2}$ in \eqref{e80}, we get
\begin{equation} \label{e81}
G_6(b)=\frac{3\pi}{16b^2}-\frac{\pi b}{24}-\sqrt{\pi}\,\sum_{r=0}^{\infty}\frac{\zeta({-}3/2-r)({-}1)^r}{r!\,(2r+1)}\,b^{2r+2}
\end{equation}
when $0<b<\sqrt{2\pi}$. The two results \eqref{e78} and \eqref{e81} can be combined, as in \cite[Sec.~5]{Janssen2006} and this yields
\begin{equation} \label{e82}
G_4(b)=\frac{\pi}{16b^3}+\frac{\pi b}{24}+\tfrac12\,\zeta({-}1/2)\,\sqrt{\pi}+\sqrt{\pi}\,\sum_{r=0}^{\infty}\frac{\zeta({-}3/2-r)({-}1)^r\,b^{2r+2}}{r!\,(2r+1)(2r+2)}
\end{equation}
when $0<b<\sqrt{2\pi}$.
Finally, we can rewrite
\begin{align}
\frac{1}{2b^2}G_2(b) &= \frac{\pi}{8b^3} - \frac{\pi}{8b} - \frac{\sqrt{\pi}}{2} \sum_{r=0}^\infty \frac{\zeta({-}1/2-r)(1-)^rb^{2r}}{r!(2r+1)} \nonumber\\
&= \frac{\pi}{8b^3} - \frac{\pi}{8b} - \frac{\sqrt{\pi}}{2} \sum_{r=-1}^\infty \frac{\zeta({-}3/2-r)(-1)^{r+1}b^{2r+2}}{(r+1)!(2r+3)} \nonumber \\
&= \frac{\pi}{8b^3} - \frac{\pi}{8b} - \tfrac{1}{2}\zeta(-1/2)\sqrt{\pi} 
+ \sqrt{\pi} \sum_{r=0}^\infty \frac{\zeta(-3/2-r)(-1)^r b^{2r+2}}{r!\,(2r+2)(2r+3)}
\label{e82a}
\end{align}
and use \eqref{e82} and \eqref{e82a} in \eqref{e74}, by which we obtain for $0<b<\sqrt{2\pi} $,
\begin{align} 
G_3(b) &=
\frac{\pi}{16b^3}-\frac{\pi}{8b}-\frac{\pi b}{24}-\zeta({-}1/2)\,\sqrt{\pi}  \nonumber \\
&\qquad +\sqrt{\pi} \sum_{r=0}^\infty \frac{ \zeta(-3/2-r)(-1)^rb^{2r+2}}{r!\,(2r+2)} \Big[ \frac{1}{2r+3}-\frac{1}{2r+1}\Big]\nonumber\\
 &=  \frac{\pi}{16b^3}-\frac{\pi}{8b}-\frac{\pi b}{24}-\zeta({-}1/2)\,\sqrt{\pi}  -~2\sqrt{\pi}\,\sum_{r=0}^{\infty}\,\frac{\zeta({-}3/2-r)({-}1)^r\,b^{2r+2}}{r!\,(2r+1)(2r+2)(2r+3)}.\label{e83}
\end{align}  
The right-hand side of \eqref{e83} equals the right-hand side of \cite[Eq.~(2.3)]{Janssen2006} multiplied by ${\pi}/{(2b)}$ with $\beta=b\,\sqrt{2}$. 
\subsection{Variance of the congestion level}\label{sec4}
We have from \eqref{e38} in Section \ref{sec1}, using the same approach and notation as in Section \ref{subsec3.1} for $\mu_Q$, that $\sigma_Q^2$ is given with exponentially small error by
\begin{equation} \label{e84}
\frac{-s\,\alpha }{2\pi i}\,\int_{-\frac12\delta}^{\frac12\delta}\,\frac{v\,z(v)}{(z(v)-1)^2}~\frac{B\,\exp({-}\frac12\,s\,\alpha \,v^2)}{1-B\,\exp({-}\frac12\,s\,\alpha \,v^2)}\dd v,
\end{equation}
with $B$ and $\alpha $ given in \eqref{e42}. From $z({-}v)=z^{\ast}(v)$ for real $v$ we now compute
\begin{equation*} 
\frac{z(v)}{(z(v)-1)^2}-\frac{z({-}v)}{(z({-}v)-1)^2}={-}2i\,\frac{|z(v)|^2-1}{|z(v)-1|^4}\,{\rm Im}(z(v)) ,
\end{equation*}
and so \eqref{e84} becomes
\begin{equation} \label{e86}
\frac{s\alpha }{\pi}\,\int_0^{\frac12\delta}\,\frac{|z(v)|^2-1}{|z(v)-1|^4}\,v\,{\rm Im}(z(v))\,\frac{B\,\exp({-}\frac12\,s\,\alpha \,v^2)}{1-B\,\exp({-}\frac12\,s\,\alpha \,v^2)}\dd v .
\end{equation}
From
\begin{equation*}
{\rm Im}(z(v))=v+O(v^3) ,\qquad |z(v)|^2-1=z_{\rm sp}^2-1+O(v^2) ,
\end{equation*}
we get for the expression in \eqref{e86}
\begin{equation} \label{y70}
\frac{s\alpha }{\pi}\,\int_0^{\frac{1}{2}\delta}\,\frac{v^2\,(z_{\rm sp}^2-1+O(v^2))(1+O(v^2))}{((z_{\rm sp}-1)^2+v^2 + O((z_{\rm sp}-1)\,v^2)+O(v^4))^2}
\frac{B\,\exp({-}\frac12\,s\,\alpha \,v^2)}{1-B\,\exp({-}\frac12\,s\,\alpha \,v^2)}\dd v.
\end{equation}
 When $2\eta-1<0$, we have as for the case of $\mu_Q$ in Section \ref{subsec3.1} that the whole expression in \eqref{y70} is $O(\exp({-}b^2\,s^{1-2\eta}))$ for any $b\in(0,b_0)$, as $s\to  \infty$. When $2\eta-1\geq 0$, we get as in the case of $\mu_Q$ after substitution $v = t\sqrt{{2}/{(s\,\alpha })}$ for the expression in \eqref{y70}
\begin{equation*} 
\frac{2}{\pi}\,\Big(\frac{s\,\alpha }{2}\Big)^{3/2}\,\int_0^\infty\frac{t^2\,(z_{\rm sp}^2-1+O(t^2/s))(1+O(t^2/s))}{(d^2(s)+t^2)^2\,(1+O(1/s^{\eta})+O(t^2/s))}~\frac{B\,\ee^{{-}t^2}}{1-B\,\ee^{{-}t^2}}\dd t.
\end{equation*}
When $2\eta-1\geq 0$, the leading order behavior of $\sigma_Q^2$ depends crucially on the factor $z_{\rm sp}^2-1+O(t^2/s)$, where
\begin{equation*} 
z_{\rm sp}^2-1 = \frac{2\,\gamma\,\mu}{\sigma^2\,s^\eta}\,\left(1+O(s^{-\eta})\right)
\end{equation*}
is dominant when $\eta<1$, while the $O(t^2/s)$ is dominant when $\eta>1$. In the case that $\eta\in(1/2,1)$, we get for the leading order behavior of $\sigma_Q^2$
\begin{align*}
\frac{2}{\pi}\,\Big(\frac{s\,\alpha }{2}\Big)^{3/2}& \,\frac{2\,\gamma\,\mu}{\sigma^2\,s^\eta}\,\int_0^\infty\frac{t^2}{(d^2(s)+t^2)^2}\cdot~\frac{\ee^{{-}d^2(s)-t^2}}{1-\ee^{{-}d^2(s)-t^2}}\dd t\,\left(1+O(s^{\eta-1})\right)\nonumber\\
&= \frac{\gamma\,\sigma}{\pi}\,\Big(\frac{2}{\mu}\Big)^{1/2}\,s^{3/2-\eta}\,G_3(d(s))\,\left(1+O(s^{\eta-1})\right),
\end{align*}
where \eqref{e26}, \eqref{e27} and \eqref{e42} have been used for $\alpha  = g''(z_{\rm sp})$ and where $G_3$ is given in \eqref{e72}.

When we insert the expansion \eqref{e83} for $G_3(b)$, with the whole series on the second line being $O(b^2)$, we get the leading order behavior of $\sigma_Q^2$ as
\begin{align}
s^{2\eta}\,\Big( \frac{\sigma^4}{4\,\gamma^2\mu^2}- \frac{\sigma^2}{4\,\mu}&\,\frac{1}{s^{2\eta-1}} - \Big(\frac{2\,\sigma^2}{\pi\,\mu}\Big)^{1/2}\,\frac{\gamma\,\zeta(-1/2)}{s^{3\eta-3/2}}\nonumber\\
& - \frac{\gamma^2}{24\,s^{5\eta-5/2}}+O(s^{1-4\eta})\Big)\,\left(1+O(s^{\eta-1})\right)\nonumber \\
&\ = s^{2\eta}\,\frac{\sigma^4}{4\,\gamma^2\,\mu^2}\,\Big(1+O(s^{\max(1-2\eta,\eta-1)})\Big)\label{y74}
\end{align}
when $\eta\in(1/2,1)$. For the case $\eta=1/2$, we get the leading order behavior, assuming $0<b_0<\sqrt{2\pi}$,
\begin{align}
\frac{\sigma^2 s}{\mu}\left[ \frac{1}{8\,b_0^2} - \frac{1}{4}-\frac{1}{12}\,b_0^2 - \frac{2\,\zeta(-1/2)}{\sqrt{\pi}}\,b_0- \frac{4}{\sqrt{\pi}}\,\sum_{r=0}^\infty \frac{\zeta(-3/2-r)\,(-1)^r\,b_0^{2r+3}}{r!\,(2r+1)\,(2r+2)\,(2r+3)} \right]\label{y75}
\end{align}
with relative error $O(s^{-1/2})$. The expression between brackets in \eqref{y75} coincides with the right-hand side of \cite{Janssen2006}, (2.3) with $\beta = b_0\,\sqrt{2}$.

This leads to the following two results.
\begin{theorem} \label{varthm}
For $\eta\in[1/2,1)$,
\begin{equation*} 
\sigma_Q^2 = \frac{\gamma\,\sigma_X}{\pi}\,\sqrt{\frac{2}{\mu}}\,s^{3/2-\eta}\,G_3(d(s))\,\left(1+O(s^{\eta-1})\right)
\end{equation*}
with $G_3$ given in \eqref{e72}.
\end{theorem}

\begin{proposition}\label{varprop}
For $\eta\in(0,1/2)$, and for all $b<b_0$,
\begin{equation*}
\sigma_Q^2 = O(\exp({-}b^2\,s^{1-2\eta})).
 \end{equation*}
 For $\eta = 1/2$, $\sigma_Q^2$ equals expression \eqref{y75} with relative error $O(s^{-1/2})$. For $\eta\in(1/2,1)$ and $b_0\in(0,\sqrt{2\pi})$, $\sigma_Q^2$ has the form in \eqref{y74}.
\end{proposition}

As in Section \ref{subsec3.3} for the mean congestion level
with $\eta=1/2$, it is possible to give a correction term which involves now integrals and series with $\zeta$-functions as considered in \cite[Secs.~4-5]{Janssen2007}.

\subsection{The empty-system probability} \label{sec5}

We have from \eqref{e6} by proceeding as in \eqref{e13}--\eqref{e17} that
\begin{align} \label{e100}
{\rm ln}\,[\P(Q=0)] & =  \frac{s}{2\pi i}\,\int_{|z|=1+\eps}\,{\rm ln}\Big(\frac{z}{z-1}\Big)\,\frac{g'(z)\,\exp(s\,g(z))}{1-\exp(s\,g(z))}\,\dd z \nonumber \\[3.5mm]
& =  \frac{1}{2\pi i}\,\int_{|z|=1+\eps}\,\frac{1}{z(z-1)}\,{\rm ln}\left(1-\exp(s\,g(z))\right)\,\dd z ,
\end{align}
where in the last step we used partial integration (noting that ${\rm Re}\,[g(z)]<0$ on $|z|=1+\eps$). Then, as in Section \ref{sec1} for $\mu_Q$, the last integral in \eqref{e100} is, with exponentially small error, given by
\begin{equation} \label{e101}
\frac{1}{2\pi i}\,\int_{-\frac12\delta}^{\frac12\delta}\,\frac{z'(v)}{z(v)(z(v)-1)}\,{\rm ln}\left(1-B\,\exp\big(-\tfrac{1}{2} s\,\alpha  v^2\big)\right)\,\dd v .
\end{equation}
Now for $v\geq0$ from $z({-}v)=z^{\ast}(v)$, $z'({-}v)={-}(z'(v))^{\ast}$,
\begin{align*} 
& \frac{z'(v)}{z(v)(z(v)-1)}+\frac{z'({-}v)}{z({-}v)(z({-}v)-1)}=2i\,{\rm Im}\,\Big[\frac{z'(v)}{z(v)(z(v)-1)}\Big] \nonumber \\
& \qquad \qquad =~2i\,{\rm Im}\,\Big[\frac{z'(v)\,z^{\ast}(v)(z^{\ast}(v)-1)}{|z(v)|^2\,|z(v)-1|^2}\Big] \nonumber \\
& \qquad \qquad =~2i\,\frac{z_{\rm sp}-1+O(v^2)}{(z_{\rm sp}+O(v^2))((z_{\rm sp}-1)^2+v^2-2c_2(z_{\rm sp}-1)\,v^2+O(v^4))}\,,
\end{align*}
where we used \eqref{e32} and the fact that $z_{\rm sp}$ and $c_2$ are real with $z_{\rm sp}>1$. Therefore, we get for the expression in \eqref{e101}
\begin{align} 
\frac{1}{\pi} &\int_0^{\frac{1}{2}\delta}\frac{1}{z_{\rm sp}{\rm +}O(v^2)}\frac{z_{\rm sp}-1+O(v^2)}{(z_{\rm sp}-1)^2+v^2+O((z_{\rm sp}-1)v^2)+O(v^4)} \nonumber \\
& \qquad \qquad \qquad \qquad \cdot {\rm ln}\left(1-B\exp(-\tfrac12 s\,\alpha  v^2)\right)\dd v.
\label{y77}
\end{align}
In the case that $2\eta-1<0$, we have as earlier that the whole expression in \eqref{y77} is $O(\exp({-}b^2\,s^{1-2\eta}))$ for any $b\in(0,b_0)$, as $s\to  \infty$. In the case that $2\eta-1\geq 0$, we substitute $v=t\sqrt{{s}/{(2\,\alpha )}}$, and we get as earlier for the expression \eqref{y77}, assuming also that $\eta<1$,
\begin{align*}
\frac{1}{\pi}&\,\sqrt{s\,\alpha /2}\,\int_0^{\infty}\frac{z_{\rm sp}-1+O(t^2/s)}{(d^2(s)+t^2)\,(1+O(s^{-\eta})+O(t^2/s))}\,{\rm ln}(1-B\,\ee^{-t^2}) \,\dd t\nonumber\\
&= \frac{1}{\pi}\,\int_0^{\infty}\frac{\sqrt{s\,\alpha /2} \ (z_{\rm sp}-1)}{d^2(s)+t^2}{\rm ln}(1-B\,\ee^{-t^2})\,\dd t\,\left(1+O(s^{\eta-1})\right)\nonumber\\
&= \frac{1}{\pi}\,\int_0^{\infty}\frac{d(s)}{d^2(s)+t^2}{\rm ln}(1-\ee^{{-}d^2(s)-t^2})\,\dd t\,\left(1+O(s^{\eta-1})\right).
\end{align*}
Here we also used \eqref{y46} and that $1/s^{3\eta-1} = O(d^2(s)/s^\eta)$, so that
\begin{equation*} 
(\tfrac12\,s\,\alpha )^{1/2}\,(z_{\rm sp}-1) = d(s)\,\left(1+O(s^{-\eta})\right) = d(s)\left(1+O(s^{\eta-1})\right),
\end{equation*}
since $\eta\geq 1/2$.

We have for $b>0$
\begin{align}
\frac{1}{\pi}&\int_0^\infty \frac{b}{b^2+t^2}\,{\rm ln}(1-\exp({-}b^2-t^2))\dd t =-\frac12\sum_{k=0}^{\infty}\,\frac{1}{k+1}\,{\rm erfc}(b\,\sqrt{k+1}) = -F(b\,\sqrt{2}),\label{y80}
\end{align}
where according to \cite[Eq.~(3.3) \& (3.12)]{Janssen2006} for $\beta>0$
\begin{align}
F(\beta) &= \sum_{n=1}^\infty\,\frac{1}{n}\,\frac{1}{\sqrt{2\pi}}\,\int_{\beta\sqrt{n}}^\infty \ee^{-x^2/2}dx\nonumber\\
&= -{\rm ln}\,\beta - \frac12\,{\rm ln}\,2 - \frac{1}{\sqrt{2\pi}}\,\sum_{r=0}^\infty \frac{\zeta(1/2-r)\,(-1/2)^r\,\beta^{2r+1}}{r!\,(2r+1)},\label{y81}
\end{align}
the last identity being valid for $0<\beta<2\sqrt{\pi}$.

Using \eqref{y81} with $\beta^2 = d^2(s)= b_0^2/s^{2\eta-1}$, with the entire series on the second line being $O(\beta)$, we get the leading order behavior of ${\rm ln}[\P(Q=0)]$ as
\begin{equation} \label{y82}
\Big({-}(\eta-1/2)\,{\rm ln}\,s+{\rm ln}(2\,b_0)+O(s^{1/2-\eta})\Big)\left(1+O(s^{\eta-1})\right)
\end{equation}
when $\eta\in(1/2,1)$. For $\eta = 1/2$, we get the leading order behavior, assuming $0<b_0<\sqrt{2\pi}$,
\begin{equation} \label{y83}
{\rm ln}(2\,b_0) + \frac{1}{\sqrt{\pi}}\,\sum_{r=0}^\infty \,\frac{\zeta(1/2-r)\,(-1)^r}{r!\,(2r+1)}\,b_0^{2r+1}
\end{equation}
with relative error $O(s^{-1/2})$. The expression \eqref{y83} coincides with ${\rm ln}[\mathbb{P}(M=0)]$ as given by \cite[Eq.~(2.1)]{Janssen2006} with $\beta = b_0\,\sqrt{2}$. The next two results summarize the above.

\begin{theorem} \label{emptythm}
For $\eta\in(1/2,1)$,
\begin{equation*} 
{\rm ln}[\mathbb{P}(Q=0)] = - F\big(d(s)\,\sqrt{2}\big)\left(1+O(s^{\eta-1})\right)
\end{equation*}
with $F$ given by \eqref{y81}.
\end{theorem}
\begin{proposition} \label{emptyprop}
For $\eta\in (0,1/2)$, and for all $b<b_0$,
\begin{equation*}
{\rm ln}[\mathbb{P}(Q=0)] = O(\exp({-}b^2\,s^{1-2\eta})).
 \end{equation*}
 For $\eta=1/2$, ${\rm ln}[\mathbb{P}(Q=0)]$ equals $-F(b_0\,\sqrt{2})$ with a relative error $O(1/\sqrt{s})$. For $\eta\in (1/2,1)$ and $0<b_0<\sqrt{2\pi}$, ${\rm ln}[\mathbb{P}(Q=0)]$ has leading order behavior as in \eqref{y82}.
\end{proposition}

As in Section \ref{subsec3.3} for the mean congestion level case
with $\eta=1/2$, it is possible to give a correction term which involves now the integrals in \eqref{y80} and \eqref{e51}.

\section{Numerical examples}\label{numm}
\subsection{Accuracy of the approximations}
In this subsection we present a numerical example that serves to illustrate the accuracy of the derived heavy-traffic approximations. Consider the Poisson case
\begin{equation*}
\tilde X(z)=\ee^{z-1},\quad \mu = \sigma^2 = 1.
\end{equation*}
We fix $\mu$ and vary $n$ with the value of $s$, according to
\begin{equation*}
\vartheta = \frac{n}{s} = 1-\frac{\gamma}{s^\eta}
\end{equation*}
for some $\gamma>0$ and $\eta\geq 1/2$. To calculate the exact value of the mean congestion level we use the expression, see \cite{Boudreau1962},
\begin{equation*}
\mu_Q=\frac{\sigma_A^2}{2(s-\mu_A)}-\frac{s-1+\mu_A}{2}+\sum_{k=1}^{s-1}\frac{1}{1-z_k}.
\end{equation*}
Here $z_1,\ldots,z_{s-1}$ are the zeros of $z^s-A(z)$ in $|z|<1$. We apply the method of successive substitution described in \cite{Janssen2005} to obtain accurate numerical approximations for $z_1,...,z_{s-1}$ and consequently $\mu_Q$.

From Theorem \ref{mainthm}, we find that the leading order behavior of $\mu_Q$ is given by
\begin{equation} \label{x18}
\frac{\sqrt{2s}}{\pi}\,G_0\Big(\frac{\gamma}{\sqrt{2}\,s^{\eta-\frac{1}{2}}}\Big).
\end{equation}
In order to find the correction terms, we proceed by setting $\eta = 1/2$. Deriving constants $C_1,C_2,C_3,$ and $C_4$ for our setting and substituting these into \eqref{e71},  we get for $\mu_Q$, with an absolute error of $O(s^{-1/2})$, the approximation
\begin{equation*}
\frac{\sqrt{2\,s}}{\pi}\Big(\Big(1-\frac{\gamma}{3\,\sqrt{s}}\Big)\,G_0(b_0)-\frac{\gamma^3}{3\,\sqrt{s}}\,(\,G_3(b_0)+G_4(b_0))\Big),
\end{equation*}
which by \eqref{e49} and \eqref{e74} reduces to
\begin{equation}\label{x20}
\frac{\sqrt{2\,s}}{\pi}\,G_0(b_0)-\frac{\sqrt{2}\,\gamma}{3\,\pi}\,G_1(b_0).
\end{equation}

\begin{table}
\centering
\begin{tabular}{|r|rrrr|}
\hline
$s$ & $\rho$ & $\mu_Q$ & \eqref{x18} & \eqref{x20} \bigstrut \\
\hline
 10 & 0.683 & 0.244 & 0.399 & 0.247 \bigstrut[t] \\
 20 & 0.776 & 0.410 & 0.565 & 0.412 \\
 50 & 0.858 & 0.739 & 0.893 & 0.741 \\
 100 & 0.900 & 1.110 & 1.263 & 1.111 \\
 200 & 0.929 & 1.633 & 1.787 & 1.634 \\
 500 & 0.955 & 2.672 & 2.825 & 2.673 \\
 1000 & 0.968 & 3.843 & 3.996 & 3.843 \bigstrut[b]\\
 \hline
\end{tabular}
\caption{Numerical results for $\gamma = 1$.}\label{tab:poisson1}
\end{table}
\begin{table}
\centering
\begin{tabular}{|r|rrrr|}
\hline
$s$ & $\rho$ & $\mu_Q$ & \eqref{x18} & \eqref{x20} \bigstrut \\
\hline
 10 & 0.968 & 13.707 & 14.046 &13.732 \bigstrut[t] \\
 20 & 0.977 & 19.533 & 19.865 &19.551\\
 50 & 0.985 & 31.084 & 31.409 &31.095\\
 100 & 0.990 & 44.097 & 44.419 &44.106\\
 200 & 0.992 & 62.499 & 62.819 &62.505\\
 500 & 0.995 & 99.008 & 99.325 &99.011\\
 1000 & 0.996 & 140.152 & 140.468 & 140.154 \bigstrut[b]\\
 \hline
\end{tabular}
\caption{Numerical results for $\gamma = 0.1$.}\label{tab:poisson2}
\end{table}
\begin{table}
\centering
\begin{tabular}{|r|rr|rr|rr|}
\cline{2-7}
 \multicolumn{1}{c|}{} & \multicolumn{2}{c|}{$\eta=0.6$} & \multicolumn{2}{c|}{$\eta=0.75$} & \multicolumn{2}{c|}{$\eta=0.9$} \bigstrut\\
 \hline
$s$ &  $\mu_Q$ & \eqref{x18} &  $\mu_Q$ & \eqref{x18} &  $\mu_Q$ & \eqref{x18} \bigstrut\\
\hline
 10 & 17.781 & 18.125 & 25.970 & 26.318 & 37.553 & 37.905 \bigstrut[t]  \\
 20 & 27.309 & 27.647 & 44.391 & 44.734 & 71.195 & 71.541 \\
 50 & 47.948 & 48.281 & 89.623 & 89.961 & 164.637 & 164.978 \\
 100 & 73.245 & 73.574 & 152.031 & 152.367 & 309.353 & 309.692 \\
 200 & 111.752 & 112.079 & 257.435 & 257.769 & 580.170 & 580.507 \\
 500 & 195.082 & 195.409 & 515.443 & 515.776 & 1329.581 & 1329.917 \\
 1000 & 297.122 & 297.448 & 870.524 & 870.857 & 2487.227 & 2487.562 \bigstrut[b]\\
 \hline
\end{tabular}
\caption{Numerical results for $\gamma=0.1$ and several values of $\eta$.}\label{tab:poisson3}
\end{table}
\noindent Numerical results for $\eta=1/2$ and various values of $s$ are given in Table \ref{tab:poisson1} and \ref{tab:poisson2}, for $ \gamma = 1$ and  $\gamma = 0.1$, respectively.
We note that for small $s$ the leading order approximation is still off by a significant amount, while the refinement only shows an error in the second decimal for $\gamma = 0.1$. This seems to justify the use of the correction term.
In Table \ref{tab:poisson3} we compare the approximation \eqref{x18} against the exact value of $\mu_Q$ for three values of $\eta\geq 1/2$ to assess the influence of $\eta$. Clearly, the leading order approximation is relatively accurate for all three scenarios. As expected, the mean congestion increases along with $\eta$, since utilization approaches 1 more rapidly in this case.
\subsection{Connection to other queueing models}\label{subsec62}
As argued in the introduction, we believe that the heavy-traffic behavior for the discrete model in this chapter will up to leading order be universal for a wide range of other models (when subjected to the same heavy-traffic regime \eqref{bb}). We shall now substantiate this for many-server systems, for which under \eqref{bb}, it turns out that the mean congestion is $O(s^\eta)$.  We compare the mean congestion level in our discrete queue with that in the multi-server systems $M/M/s$, $M/D/s$ and Gamma/Gamma/$s$, all with unit mean service time and occupation rate $1-\gamma/s^\eta$.

\begin{figure}
\centering
\begin{tikzpicture}
\begin{axis}[
	xmin = 2.5,
	xmax = 6.5,
	ymin = 0,
	ymax = 7.2,
	xlabel = {$\log(s)$},
	ylabel = {$\log(\mu_Q)$},
	x label style={at={(0.6,-0.1)}},
	y label style={at={(-0.06,0.7)}},
	axis line style={->},
	axis lines = left,
	legend cell align=left,
	legend style = {at = {(axis cs: 6.4,0.2)},anchor = south east},
	yscale = 0.8,
	xscale = 1
]

\addplot[thick,mark = o,mark options={scale=1.25}] table[x=log_s,y=mms] {Chapter_2/tikz/novel_figure1.txt};
\addplot[thick, mark=triangle,dashed,mark options={scale=1.25,solid}] table[x=log_s,y=mds] {Chapter_2/tikz/novel_figure1.txt};
\addplot[thick,mark=square,dotted,mark options={scale=1.25,solid}] table[x=log_s,y=ggs] {Chapter_2/tikz/novel_figure1.txt};

\legend{$M/M/s$,$M/D/s$,Gamma/Gamma/$s$}
\end{axis}

\end{tikzpicture}
\caption{$\mu_Q$ plotted against $s$ on log scale for 3 queues for $\eta=0.75$.}
  \label{fig1}
\end{figure}
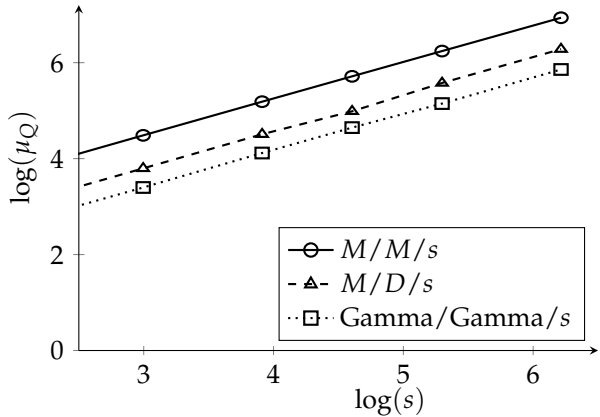
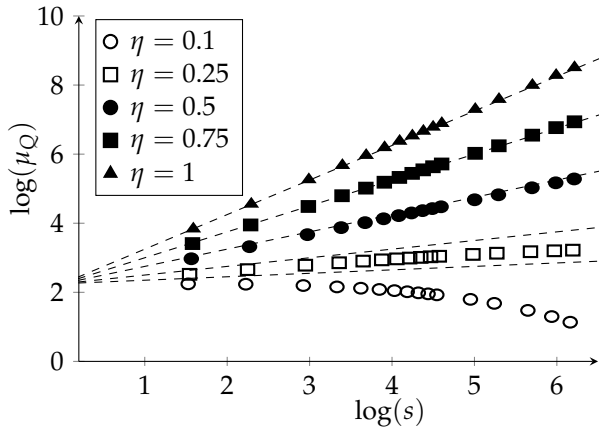
\begin{figure}
\centering
\begin{tikzpicture}
\begin{axis}[
	xmin = 0.2,
	xmax = 6.5,
	ymin = 0,
	ymax = 10,
	xlabel = {$\log(s)$},
	ylabel = {$\log(\mu_Q)$},
	x label style={at={(0.6,-0.1)}},
	y label style={at={(-0.06,0.7)}},
	axis line style={->},
	axis lines = left,
	legend cell align=left,
	legend style = {at = {(axis cs: 0.4,7.35)},anchor = west},
	yscale = 0.8,
	xscale = 1
]

\addplot[thick,only marks,mark = o,mark options={scale=1.25}] table[x=n01,y=m01] {Chapter_2/tikz/novel_figure2.txt};
\addplot[thick,only marks,mark = square,mark options={scale=1.25}] table[x=n025,y=m025] {Chapter_2/tikz/novel_figure2.txt};
\addplot[thick,only marks,mark = *,mark options={scale=1.25}] table[x=n05,y=m05] {Chapter_2/tikz/novel_figure2.txt};
\addplot[thick,only marks,mark = square*,mark options={scale=1.25}] table[x=n075,y=m075] {Chapter_2/tikz/novel_figure2.txt};
\addplot[thick,only marks,mark = triangle*,mark options={scale=1.25}] table[x=n1,y=m1] {Chapter_2/tikz/novel_figure2.txt};

\addplot[dashed] coordinates{ (0,2.25) (7,9.25) };
\addplot[dashed] coordinates{ (0,2.25) (7,7.5) };
\addplot[dashed] coordinates{ (0,2.25) (7,5.75) };
\addplot[dashed] coordinates{ (0,2.25) (7,4) };
\addplot[dashed] coordinates{ (0,2.25) (7,2.95) };
\legend{\ $\eta=0.1$,\ $\eta=0.25$,\ $\eta=0.5$,\ $\eta=0.75$,\ $\eta=1$}
\end{axis}

\end{tikzpicture}
  \caption{$\mu_Q$ of $M/M/s$ plotted against $s$ on log scale for different values of $\eta$.}
  \label{fig2}
\end{figure}
Figure \ref{fig1} shows on logarithmic scale the mean congestion levels for $\gamma=0.1$ and  $\eta=0.75$  under the specified scaling for three systems. We also display three lines with slope 0.75 for comparison, which confirms that mean congestion levels are of the order $s^\eta$, also in these multi-server system. Formally establishing this heavy-traffic behavior for these multi-server system is an important open problem and requires other mathematical approaches than the ones taken in this chapter (see the introduction for more details).

Figure \ref{fig2} shows the mean queue length in the $M/M/s$ system for several values of $\eta$, again on logarithmic scale, together with lines with slope $\eta$. For $\eta\geq 1/2$, we see the same $O(s^\eta)$ behavior, similar as for $\mu_Q$ in our discrete model. For $\eta<1/2$ the mean queue length decays, again in agreement with our results for $\mu_Q$. We note that this qualitative behavior of the $M/M/s$ system was also observed by \cite[Thm.~4.1]{maman}, by proving that the mean waiting time in the $M/M/s$ queue under \eqref{bb} is of order $1/s^{1-\eta}$, which by Little's law implies the mean queue length to be of order $s^\eta$.

\chapter{Overdispersion}

\begin{chapterstart}
Arrival processes to service systems often display fluctuations that are larger than anticipated under the Poisson assumption, a phenomenon that is referred to as \textit{overdispersion}.
Motivated by this, we analyze a class of discrete stochastic models for which we derive heavy-traffic approximations that are scalable in the system size.
Subsequently, we show how this leads to novel capacity sizing rules that acknowledge the presence of overdispersion. This, in turn, leads to robust approximations for performance characteristics of systems that are of moderate size and/or may not operate in heavy traffic. 
\end{chapterstart}

\begin{flushright}
Based on\\
\textbf{Robust heavy-traffic approximations for\\ service systems facing overdispersed demand}\\
\textit{Britt Mathijsen, Guido Janssen, Johan van Leeuwaarden \& Bert Zwart}\\
arxiv.org/abs/1512.05581
\end{flushright}
\newpage 

\section{Introduction}\label{intro}

In the previous chapter, we analyzed the scaling limit of a queueing model in which demand exhibits stochastic fluctuations that are asymptotically proportional to the square-root of the nominal load, while we deliberately chose to deviate from the square-root staffing principle by allocating a variability hedge that does not match the order of these fluctuations.
This chapter in some ways does the opposite. 
We assume the demand faced by the queueing system is more volatile than anticipated by the independent many-sources paradigm that leads to Poisson traffic models. 
As will become clear in this chapter, this in fact \emph{requires} an adaptation of the square-root staffing principle in order to maintain the desirable properties of the QED regime. 
We start by motivating our research through empirical evidence of the presence of so-called \emph{overdispersion} in arrival processes faced by service systems reported by recent literature. \\

\noindent
\textbf{Motivation.}
The bulk of the queueing literature assumes perfect knowledge about the model primitives, including the mean demand per time period. For large-scale service systems, like health care facilities, communication systems or call centers, the dominant assumption is that demand arrives according to a (non)homogeneous Poisson process, which in practice translates into the assumption that arrival rates are known for each basic time period (second, hour or day).
Although natural and convenient from a mathematical viewpoint, the Poisson assumption often fails to be confirmed in practice. A deterministic arrival rate implies that the demand over any given period is a Poisson random variable, whose variance equals its expectation. A growing number of empirical studies shows that the variance of demand typically deviates from the mean significantly. Recent work \cite{Kim2015b,Kim2015a} reports variance being strictly less than the mean in health care settings employing appointment booking systems. This reduced variability, known as underdispersion, can be accredited to the goal of the booking system to create a more predictable arrival pattern. 
On the other hand, in other scenarios with no control over the arrivals, the variance typically dominates the mean, see \cite{Avramidis:2004, Bassamboo2010, Bassamboo2009, Brown2005, Chen2001, Gans2015, Gurvich2010, koolejongbloed, kimwhitt,  maman, Mehrotra2010, Robbins2010, Steckley2009, Zan2012}. 
The feature that variability is higher than one expects from the Poisson assumption is referred to as overdispersion. The latter concept will be the center of our attention in this chapter. 

Stochastic models with the Poisson assumption have been widely applied to optimize capacity levels in service systems. The goal is to minimize operating costs while providing sufficiently high QoS in terms of performance measures such as mean delay or excess delay. When stochastic models, however, do not take into account overdispersion, resulting performance estimates are likely to be overoptimistic. The system then ends up being underprovisioned, which possibly causes severe performance problems, particularly under critical loading.
\\
\\*
\noindent
\textbf{Causes of overdispersion.} 
The literature discussed above proves that the presence of overdispersion is widespread across applications.
It however does not specify what causes the increased variability in the arrival process. 
We name two possible explanations. 

First, we revisit the many-sources characterization of demand inflow discussed in Chapter 2. 
Recall that in this setting, demand is generated by $n$ stochastically identical and independent sources, with $n$ large, so that workload arriving to the system in period $j$ is given by $A^{(n)}_j = \sum_{i=1}^n A_{i,j}$, where $A_{i,j}$, $i=1,2,\ldots,n$ are i.i.d.~random variables. 
This resulted in nominal workload $\mu_n = n\mu$ and $\sigma_n^2 = n\sigma^2$, thus both of order $n$. 
If we now relax the assumption on the (pairwise) independence of the sources, but rather consider the scenario in which these are positively correlated, then the nominal load remains to be equal to $n \mu$, while the variance of demand becomes
\begin{equation*}
\sigma_n^2 = \Var A_j^{(n)} = n\, \Var A_{1,j} + n(n-1)\,{\rm Cov}(A_{1,j},A_{2,j}),
\end{equation*}
which is of higher order than $n$ if $n\,{\rm Cov}(A_{1,j},A_{2,j}) \to \infty$ as $n\to\infty$.

A second interpretation of overdispersion in arrival processes relates to \emph{arrival rate uncertainty}. 
The canonical process for modeling the arrival process of a service system is the Poisson process with a given arrival rate $\lambda$. 
Since model primitives, in particular the arrival rate, are typically estimated through historical data, these are prone to be subject to forecasting errors. 
In the realm of Poisson processes, this inherent uncertainty can be acknowledged by viewing the arrival rate $\Lambda_n$ itself as being stochastic. The resulting doubly stochastic Poisson process, also known as Cox process (first presented in \cite{Cox1955}), implies that demand in a given interval $A_j$ follows a mixed Poisson distribution. 
In this case, the expected demand per period equals $\mu_n = \E[\Lambda_n]$, while the variance is $\sigma_n^2 = \E[\Lambda_n]+\Var\Lambda_n$.
By selecting the distribution of the mixing factor $\Lambda_n$, the magnitude of overdispersion can be made arbitrarily large, and only a deterministic $\Lambda_n$ leads to a true Poisson process.

The mixed Poisson model presents a useful way to fit both the mean and variance to real data, particularly in case of overdispersion. 
The mixing distribution can be estimated parametrically or non-parametrically, see \cite{koolejongbloed,maman}. 
A popular parametric family is the Gamma distribution, which gives rise to an effective data fitting procedure that uses the fact that a Gamma mixed Poisson random variable follows a negative binomial distribution. 
We will in this chapter adopt the assumption of a Gamma-Poisson mixture as the demand process.\\
\\*
\textbf{Adapted QED scaling.} To deal with overdispersion 
new models are needed, scaling rules must be adapted, and existing capacity sizing rules need to be modified in order to incorporate a correct hedge against (increased) variability. 
In this chapter, we consider an extension of the discrete queueing model of Chapter 2 that has a doubly stochastic Poisson process as input, $A_j\sim\,{\rm Pois}(\Lambda_n)$ and we identify the heavy-traffic regime in which it displays QED behavior. 
That is, it fits the three asymptotic characteristics in Section 1.2.3 of this thesis.
As we argued in that particular section, a sensible candidate capacity allocation rule is $s_n = \mu_n + \beta \sigma_n$ for some $\beta>0$, which is equivalent to the scaling 
\begin{equation*}
\frac{\mu_n}{\sigma_n}\,(1-\rho_n) \to \beta, \qquad \text{as }n\to\infty.
\end{equation*}
We will verify mathematically that this is asymptotically the appropriate choice.
Studies that have adressed similar capacity allocation problems with stochastic arrival rates include \cite{Kocaga2015, maman, Whitt1999, Whitt2006}.
Of the aforementioned papers, our work best relates to \cite{maman}, in the sense that we also assess the asymptotic performance of a queueing system having a stochastic arrival rate in heavy traffic.
We therefore expand the paradigm of the QED regime, in order to have it accommodate for overdispersed demand that follows from a doubly stochastic Poisson process.
\\
\\*
\textbf{Structure of the chapter}. The remainder of this chapter is structured as follows. Our model is introduced in Section \ref{modelSection} together with some preliminary results. 
In Section 3.3 we derive the classical heavy-traffic scaling limits for the queue length process in the presence of overdispersed arrivals both for the moments and the distribution itself.
Section 3.4 presents our main theoretic result, which provides a robust refinement to the heavy-traffic characterization of the queue length measures in pre-limit systems.
In Section 3.5, we describe the numerical results and demonstrate the heavy-traffic approximation for a real data set coming from a health care setting. Section 3.6 provides some concluding remarks.

\section{Model description}\label{modelSection}
We consider the same mathematical model as in Section 2.2, in which time is divided into periods of equal length. At the beginning of each period $j=1,2,3,...$ new demand $A^{(n)}_j$ arrives to the system. The demands per period $A^{(n)}_1,A^{(n)}_2,...$ are assumed independent and equal in distribution to some non-negative integer-valued random variable $A^{(n)}$.
The system has a service capacity $s_n\in\mathbb{N}$ per period, the steady-state queue length can be characterized as, see (1.27), 
\begin{equation}
   \label{mm3}
   Q^{(n)} \equalD \max_{k\geq 0}\Bigl\{\sum_{i=1}^k (A^{(n)}_i-s_n)\Bigr\}.
   \end{equation}
For brevity, we define $\mu_n:= \E [A^{(n)}_1]$ and $\sigma_n^2 = \Var A^{(n)}_1$.
The behavior of $Q^{(n)}$ predominantly depends on the characteristics of $A^{(n)}$ and $s_n$. As noted before, $\mu_n<s_n$ is a necessary condition for the maximum in \eqref{mm3} to be finite and consequently for the queue to be stable. Before continuing the analysis of $Q^{(n)}$, we impose a set of conditions on the asymptotic properties of $s_n,\mu_n$ and $\sigma_n$.

\begin{assumption}
\label{as1}
\ \\*
\vspace{-6mm}
\begin{enumerate}
\item[{\normalfont (a)}] {\rm (Asymptotic growth)}
\begin{equation*}
\mu_n,\sigma_n \to \infty, \quad \text{\rm for } n\to\infty.
\end{equation*} 
\item[{\normalfont (b)}] {\rm (Persistence of overdispersion)}
\begin{equation*}
\sigma_n^2/\mu_n \to \infty \quad \text{\rm for } n\to\infty.
\end{equation*}
\item[{\normalfont (c)}] {\rm (Heavy-traffic condition)} 
The utilization $\rho_n := \mu_n/s_n \to 1$ as $n\to\infty$, while
\begin{equation}\label{mm5}
s_n = \mu_n + \beta\, \sigma_n,
\end{equation}
for some $\beta > 0$. This is equivalent to requiring
\begin{equation}\label{mm4}
(1-\rho_n)\frac{\mu_n}{\sigma_n} \to \beta, \qquad \text{\rm for }n\to\infty.
\end{equation}
\end{enumerate}
\end{assumption}
\noindent
Assumption \ref{as1} is assumed to hold throughout the remainder of this chapter.
Since we are mainly interested in the system behavior in heavy traffic, it is appropriate to study the queue length process in a scaled form. Substituting $s_n$ as in Assumption \ref{as1}(c), and dividing both sides of \eqref{mm3} by $\sigma_n$, gives
\begin{equation}
\label{mm6}
\frac{Q^{(n)}}{\sigma_n} = \max_{k\geq 0} \Bigl\{{\sum_{i=1}^k} \Bigl(\frac{A^{(n)}_i-\mu_n}{\sigma_n} - \beta\Bigr)\Bigr\}.
\end{equation}
By defining $\hat{Q}^{(n)} := Q^{(n)}/\sigma_n$ and $\hat{A}^{(n)}_i := (A^{(n)}_i-\mu_n)/\sigma_n$, we see that the scaled queue length process is in distribution equal to the maximum of a random walk with i.i.d. increments distributed as $\hat{A}^{(n)}-\beta$. Besides $\E[\hat{A}^{(n)}] = 0$ and $\Var \hat{A}^{(n)}=1$, the scaled and centered arrival count $\hat{A}^{(n)}$ has a few other nice properties which we turn to later in this section.

The model in \eqref{mm3} is valid for any distribution of $A^{(n)}$, also for the original case where the number of arrivals follows a Poisson distribution with fixed parameter $\lambda_n$, but in that case Assumption \ref{as1}(b) does not hold. Instead, we assume $A^{(n)}$ to be Poisson distributed with uncertain arrival rate rendered by the  non-negative random variable $\Lambda_n$. This $\Lambda_n$ is commonly referred to as the \emph{prior} distribution, while $A^{(n)}$ is given the name of a Poisson mixture, see \cite{Grandell1997}. Given that the moment generation function of $\Lambda_n$, denoted by $M^\Lambda_n(\cdot)$, exists, we are able to express the probability generating function (pgf) of $A^{(n)}$ through the former. Namely,
\begin{equation}
\label{mm7}
\tilde{A}^{(n)}(z) = \E[\E[ z^{A^{(n)}} | \Lambda_n ] ] = \E[ \exp(\Lambda_n(z-1))] = M^\Lambda_n(z-1).
\end{equation}
From \eqref{mm7}, we get
\begin{equation}
\label{mm8}
\mu_n = \E[A^{(n)}] =  \E[\Lambda_n],\qquad
\sigma_n^2 = \Var A^{(n)} = \Var \Lambda_n + \E[\Lambda_n],
\end{equation}
so that $\mu_n<\sigma_n^2$ if $\Lambda_n$ is non-deterministic. Assumption \ref{as1}(b) hence translates to 
\[\Var \Lambda_n/\E[\Lambda_n]\rightarrow \infty, \qquad n\rightarrow\infty.\]
The next result relates the converging behavior of the centered and scaled $\Lambda_n$ to that of $\hat{A}^{(n)}$.
\begin{lemma}\label{gaussStep}
Let $\mu_n,\sigma_n^2\rightarrow\infty$ and $\sigma_n^2/\mu_n\rightarrow\infty$. If
\begin{equation*}
\hat{\Lambda}_n := \frac{\Lambda_n-\mu_n}{\sigma_n}\Rightarrowd \mathcal{N}(0,1), \qquad \text{\normalfont for } n\rightarrow\infty,
\end{equation*}
then $\hat{A}^{(n)}$ converges weakly to a standard normal variable as $n\rightarrow\infty$.
\end{lemma}
\noindent
The proof can be found in Appendix \ref{formalSec}.
The prevalent choice for  $\Lambda_n$ is the Gamma distribution. The Gamma-Poisson mixture turns out to provide a very good fit to arrival counts experienced by service systems, as was observed by \cite{koolejongbloed}. Assuming $\Lambda_n$ to be of Gamma type with scale and rate parameters $a_n$ and $1/b_n$, respectively, we get for the pgf of $A^{(n)}$:
\begin{equation}
\label{r0}
\tilde{A}^{(n)}(z) = \Bigl(\frac{1}{1+b_n(1-z)}\Bigr)^{a_n},
\end{equation}
in which we recognize the pgf of a negative binomial distribution with parameters $a_n$ and $1/(b_n+1)$, so that
\begin{equation*}
 \label{t21}
 \mu_n = a_nb_n,\qquad \sigma_n^2 = a_nb_n(b_n+1).
 \end{equation*}
Note that in the context of a Gamma prior, the restrictions in Assumption \ref{as1} reduce to only two rules. For completeness, we include the revised list below. 
\begin{assumption}\label{as2}
\ \\*
\vspace{-6mm}
\begin{enumerate}
\item {\rm (Asymptotic regime and persistence of overdispersion)}
\begin{equation*}
a_n, b_n \to \infty, \quad \text{\rm for } n\to\infty.
\end{equation*}
\item {\rm (Heavy-traffic condition)}
Let 
\begin{equation*}
s_n = a_n b_n + \beta \sqrt{a_n b_n(b_n+1)},
\end{equation*}
for some $\beta>0$, or equivalently
\begin{equation*}
(1-\rho_n)\sqrt{a_n} \to \beta, \quad \text{\rm for } n\to\infty.
\end{equation*}
\end{enumerate}
\end{assumption}
The next result follows from the fact that $\Lambda_n$ is a Gamma random variable:
\begin{corollary}\label{scaledLambdaLemma}
Let $\Lambda_n\sim\text{\normalfont Gamma}(a_n,1/b_n)$, $A^{(n)}\sim{\rm Pois }(\Lambda_n)$ and $a_n,b_n\rightarrow \infty$. Then $\hat{A}^{(n)}$ converges weakly to a standard normal random variable as $n\rightarrow \infty$.
\end{corollary}
\begin{proof}
By Lemma \ref{gaussStep}, it is sufficient to prove that $\hat{\Lambda}_n\Rightarrowd\mathcal{N}(0,1)$ for this particular choice of $\Lambda_n$.  
We do this by proving the pointwise convergence of the characteristic function (cf) of $\hat{\Lambda}_n$ to $\exp({-} t^2/2)$, the cf of the standard normal distribution. 
Let $\phi_{G}(\cdot)$ denote the characteristic function of a random variable $G$. By basic properties of the cf,
\begin{align*}
\phi_{\hat{\Lambda}_n}(t) &= \ee^{-i\mu_nt/\sigma_n}\,\phi_{\Lambda_n}(t/\sigma_n)
= \ee^{-i\mu_nt/\sigma_n} \Bigl(1-\frac{i b_nt}{\sigma_n}\Bigr)^{-a_n}\nonumber\\
&= \exp\Bigl[ -\frac{i\mu_nt}{\sigma_n}\, - a_n\,{\rm ln}\Bigl(1-\frac{i b_nt}{\sigma_n}\Bigr)\Bigr]\nonumber\\
\label{g13d}
&= \exp\Bigl[ -\frac{i\mu_nt}{\sigma_n} -a_n\Bigl( {-}\frac{i\,b_nt}{\sigma_n} + \frac{b_n^2t^2}{2\sigma_n^2} + O( b_n^3/\sigma_n^3)\Bigr)\Bigr]  \nonumber\\
&= \exp\Bigl[ -\frac{b_n\,t^2}{2(b_n+1)} + O\left(1/\sqrt{a_n}\right)\Bigr] \rightarrow \exp\big({-} t^2/2\big),
\end{align*}
for $n\rightarrow\infty$. By L\'evy's continuity theorem this implies $\hat{\Lambda}_n$ is indeed asymptotically standard normal. 
\end{proof}
 The characterization of the arrival process as a Gamma-Poisson mixture is of vital importance in later sections.\\
\\*
\noindent
\textbf{Expressions for the stationary distribution.} \label{expressionsSubsec}
Our main focus is on the stationary queue length distribution, denoted by 
\[\mathbb{P}(Q^{(n)}=i) =\lim_{k\rightarrow\infty} \mathbb{P}(Q^{(n)}(k)=i).\]
Denote the pgf of $Q^{(n)}$ by
\begin{equation*}
\label{t1}
\tilde{Q}^{(n)}(w) := \sum_{i=0}^\infty \P(Q^{(n)}=i) w^i.
\end{equation*}
Furthermore, let $\mu_Q := \E[Q^{(n)}]$ and $\sigma_Q^{2} := \Var Q^{(n)}$ denote the stationary mean and variance of the queue length, respectively. 
To avoid notational complexity, we omit the superscript $(n)$ in these definitions. 
To continue our analysis of $Q^{(n)}$, we need one more condition on $A^{(n)}$. 
\begin{assumption}\label{as3}
The pgf of $A^{(n)}$, denoted by $\tilde{A}^{(n)}(w)$, exists for $|z|<r_0$, for some $r_0>1$, so that all moments of $A^{(n)}$ are finite.
\end{assumption}

We next recall two characterizations of $\tilde{Q}^{(n)}(w)$ that play prominent roles in the remainder of our analysis. 
The first characterization of $\tilde{Q}^{(n)}(w)$ originates from a random walk perspective. As we see from \eqref{mm3}, the (scaled) stationary queue length is equal in distribution to the all-time maximum of a random walk with i.i.d. increments distributed as $A^{(n)}-\beta$ (or $\hat{A}^{(n)}-\beta$ in the scaled setting). Spitzer's identity, see e.g. \cite[Theorem VIII4.2]{Asmussen2003} and Section 1.2.2 of this thesis, then gives
\begin{equation*}
\label{t3}
\tilde{Q}^{(n)}(w) = \exp\left\{\sum_{k=1}^\infty \frac{1}{k}\,\Big(\E\Big[w^{\left(\sum_{i=1}^k \{A^{(n)}_i-s_n\}\right)^+}\Big]-1\Big)\right\},
\end{equation*}
where $(x)^+ = \max\{x,0\}$. Hence,
\begin{equation*}
\label{t4}
\mu_Q = \E[Q^{(n)}] = \tilde{Q}^{(n)\prime}(1) = \sum_{k=1}^\infty \frac{1}{k}\,\E\Bigl[ {\sum_{i=1}^k} (A^{(n)}_i - s_n) \Bigr]^+,
\end{equation*}
\begin{equation*}
\label{t4a}
\sigma^{2}_Q = \Var Q^{(n)} = \tilde{Q}^{(n)\prime\prime}(1)+Q^{(n)\prime}(1)-\left(\tilde{Q}^{(n)\prime}(1)\right)^2 = \sum_{k=1}^\infty \frac{1}{k}\,\E\Bigl[ \Big(\sum_{i=1}^k (A^{(n)}_i - s_n) \Big)^+\Bigr]^2,
\end{equation*}
\begin{align*}
\label{t5}
\P(Q^{(n)}=0) = \tilde{Q}_n(0) &= \exp\Bigl\{{-}{\sum_{k=1}^\infty}\frac{1}{k}\,\P\Bigl({\textstyle\sum_{i=1}^k} (A^{(n)}_i-s_n) > 0\Bigr) \Bigr\}.
\end{align*}
A second characterization follows from Pollaczek's formula, see \cite{Abate1993} and Section 2.2.2 of this thesis:
\begin{equation}
\label{t6}
\tilde{Q}^{(n)}(w) = \exp\Bigl\{ \frac{1}{2\pi i}\int_{|z|=1+\eps} {\rm ln}\Bigl(\frac{w-z}{1-z}\Bigr) \,\frac{(z^{s_n}-\tilde{A}^{(n)}(z))'}{z^{s_n}-\tilde{A}^{(n)}(z)}\dd z\Bigr\},
\end{equation}
which is analytic for $|w|<r_0$, for some $r_0>1$. Therefore, $\eps>0$ has to be chosen such that $|w|<1+\eps<r_0$. This gives
\begin{align}
\label{t7}
\mu_Q  &= \frac{1}{2\pi i} \int_{|z|=1+\eps} \frac{1}{1-z}\,\frac{(z^{s_n}-\tilde{A}^{(n)}(z))'}{z^{s_n}-\tilde{A}^{(n)}(z)} \dd z,\\
\label{t7a}
\sigma_Q^{2} &= \frac{1}{2\pi i} \int_{|z|=1+\eps} \frac{{-}z}{(1-z)^2}\,\frac{(z^{s_n}-\tilde{A}^{(n)}(z))'}{z^{s_n}-\tilde{A}^{(n)}(z)}\dd z,\\
\label{t8}
\P(Q^{(n)}=0) &= \exp\Bigl\{ \frac{1}{2\pi i}\int_{|z|=1+\eps} {\rm ln}\Bigl(\frac{z}{z-1}\Bigr) \,\frac{(z^{s_n}-\tilde{A}^{(n)}(z))'}{z^{s_n}-\tilde{A}^{(n)}(z)}\dd z\Bigr\}.
\end{align}

Pollaczek-type integrals like \eqref{t6}-\eqref{t8} first occurred in the work of Pollaczek on the classical single-server queue (see \cite{Abate1993,Cohen1982,Janssen2008} for historical accounts). These integrals are fairly straightforward to evaluate numerically and hence give rise to efficient algorithms for performance evaluation \cite{Abate1993,boon2017pollaczek}. The integrals also proved useful in establishing heavy-traffic results by asymptotic evaluation of the integrals in various heavy-traffic regimes \cite{Kingman1962,Cohen1982,Janssen2015,boon2017pollaczek2}, and in this paper we follow that approach for a heavy-traffic regime that is suitable for overdispersion.

\section{Heavy-traffic limits}

In this section we present the result on the convergence of the discrete process $\hat{Q}^{(n)}$ to a non-degenerate limiting process and of the associated stationary moments. The latter requires an interchange of limits. Using this asymptotic result, we derive two sets of approximations for $\mu_Q$, $\sigma^2_Q$ and $\P(Q^{(n)}=0)$, that capture the limiting behavior of $Q^{(n)}$. The first set provides a rather crude estimation for the first cumulants of the queue length process for any arrival process $A^{(n)}$ satisfying Assumption \ref{as1}. The second set, which is the subject of the next section, is derived for the specific case of a Gamma prior and is therefore expected to provide more accurate, robust approximations for the performance metrics.

We start by indicating how the asymptotic properties of the scaled arrival process give rise to a proper limiting random variable describing the stationary queue length. The asymptotic normality of $\hat{A}^{(n)}$ provides a link with the Gaussian random walk and nearly deterministic queues \cite{Sigman2011a,Sigman2011b}.
The main results in \cite{Sigman2011a,Sigman2011b} were obtained under the assumption that $\rho_n\sim 1-\beta/\sqrt{n}$, in which case it follows from \cite[Thm.~3]{Sigman2011b} that the rescaled stationary waiting time process converges to a reflected Gaussian random walk.

We shall also identify the Gaussian random walk as the appropriate scaling limit for our stationary system. However, since the normalized natural fluctuations of our system are given by $\mu_n/\sigma_n$ instead of $\sqrt{n}$, we assume that the load grows like $\rho_n \sim 1 - \frac{\beta}{\mu_n/\sigma_n}$. Hence, in contrast to \cite{Sigman2011a,Sigman2011b}, our systems' characteristics display larger natural fluctuations, due to the mixing factor that renders the arrivals. Yet, by matching this overdispersed demand with the appropriate hedge against variability, we again obtain Gaussian limiting behavior. This is not surprising, since we saw in Lemma \ref{gaussStep} that the increments start resembling Gaussian behavior for $n\rightarrow\infty$. The following result summarizes this.

\begin{theorem}
\label{gaussianThm}
Let $\Lambda_n$ be a non-negative random variable such that $(\Lambda_n-\mu_n)/\sigma_n$ is asymptotically standard normal, with $\mu_n$ and $\sigma_n$ as defined in \eqref{mm8}, and $\E[\Lambda_n^3]<\infty$ for all $n\in\mathbb{N}$. Then under Assumption \ref{as1}, for $n\rightarrow \infty$,
\begin{enumerate}
\item[{\rm (i)}] $\hat{Q}^{(n)} \Rightarrowd M_\beta$,
\item[{\rm (ii)}] $\mathbb{P}(Q^{(n)} = 0) \rightarrow \mathbb{P}(M_\beta=0)$,
\item[{\rm (iii)}] $\E[\hat{Q}^{(n)}] \rightarrow \E [M_\beta]$,
\item[{\rm (iv)}] $\Var \hat{Q}_n \rightarrow \Var\, M_\beta$,
\end{enumerate}
where $M_\beta$ is the all-time maximum of a random walk with i.i.d. normal increments with mean $-\beta$ and unit variance.
\end{theorem}
The proof of Theorem \ref{gaussianThm} is given in Appendix \ref{formalSec}. The following result shows that Theorem \ref{gaussianThm} also applies to Gamma mixtures, which is a direct consequence of Corollary \ref{scaledLambdaLemma}.
\begin{corollary}
Let $\Lambda_n\sim$ \normalfont{Gamma}$(a_n,b_n)$. Then under Assumption \ref{as2} the four convergence results of Theorem \ref{gaussianThm} hold true.
\end{corollary}

It follows from Theorem \ref{gaussianThm} that the scaled stationary queueing process converges under \eqref{mm4} to a reflected Gaussian random walk. Hence, the performance measures of the original system should be well approximated by the performance measures of the reflected Gaussian random walk, yielding heavy-traffic approximations.

Like our original system, the Gaussian random walk falls in the classical setting of the reflected one-dimensional random walk, whose behavior is characterized by both Spitzer's identity and Pollaczek's formula. In particular, Pollaczek's formula gives rise to contour integral expressions for performance measures that are easy to evaluate numerically, also in heavy-traffic conditions. The numerical evaluation of such integrals is considered in \cite{Abate1993}. For $\E [M_\beta]$ such an integral is as follows
\begin{equation}
\label{g13e}
\E [M_\beta] = {-}\frac{1}{\pi}\int_0^\infty {\rm Re}\Bigl[\frac{1-\phi(-z)}{z^2}\Bigr]\dd y,
\end{equation}
where $z=x+iy$ with an appropriately chosen real part $x$, with $\phi(z) = \exp(-\beta\,z+\tfrac12\,z^2)$, the Laplace transform of a normal random variable with mean $-\beta$ and unit variance. 
Note that this integral involves complex-valued functions with complex arguments. Similar Pollaczek-type integrals exist for $\mathbb{P}(M_\beta=0)$ and $\Var M_\beta$, see \cite{Abate1993}. The following result simply rewrites these integrals in terms of a real integral and uses the fact that the scaled queue length process mimics the maximum of the Gaussian random walk for large $n$.

\begin{corollary}\label{abateThm}
Under Assumption \ref{as1}, the leading order behavior of $\mathbb{P}(Q^{(n)}=0)$, $\mu_Q$ and $\sigma^2_Q$ as $n\to\infty$ are given by, respectively,
\begin{equation}
\label{h1a}
\exp\Bigl[\frac{1}{\pi} \int_0^\infty \frac{\beta/\sqrt{2}}{\tfrac12\beta^2+t^2}\,{\rm ln}\Bigl(1-e^{-\tfrac12\beta^2-t^2}\Bigr)\dd t\Bigr],
\end{equation}
\begin{equation}
\label{h1}
\frac{\sqrt{2}\sigma_n}{\pi}\int_0^\infty \frac{t^2}{\tfrac12\beta^2+t^2}\, \frac{\exp(-\tfrac12\beta^2- t^2)}{1-\exp(-\tfrac12 \beta^2 - t^2)} \dd t,
\end{equation}
\begin{equation}
\label{h1b}
\frac{\sqrt{2}\beta\sigma_n^2}{\pi}\,\int_0^\infty \frac{t^2}{(\tfrac12 \beta^2+t^2)^2}\frac{\exp(-\tfrac12\beta^2- t^2)}{1-\exp(-\tfrac12 \beta^2 - t^2)} \dd t.
\end{equation}

\end{corollary}

\begin{proof}
According to \cite[Eq.~(15)]{Abate1993},
\begin{equation*}
\label{z1}
{-}\,{\rm ln}\,[\mathbb{P}(M_\beta=0)] = c_0,\quad \E[M_\beta]\ = c_1, \quad \Var\, M_\beta = c_2,
\end{equation*}
where
\begin{equation*}
\label{z2}
c_n = \frac{(-1)^nn!}{\pi} \,{\rm Re}\Bigl[\int_0^\infty \frac{{\rm ln}\,(1-\exp(\beta\,z+\tfrac12 z^2))}{z^{n+1}} \dd y\Bigr],
\end{equation*}
in which $z={-}x+i\,y$, $y\geq 0$, and $x$ is any fixed number between 0 and $2\beta$. 
Take $x=\beta$, so that
\begin{equation*}
 \label{z3}
 \beta z+\tfrac12 z^2 = {-}\tfrac12\beta^2 - \tfrac12 y^2\leq 0,\quad y\geq 0.
 \end{equation*}
For $n=0$, this gives
\begin{align*}
c_0 &= \frac{1}{\pi}\,{\rm Re}\Bigl[\int_0^\infty \frac{{\rm ln}\,(1-\exp({-}\tfrac12 \beta^2-\tfrac12 y^2))}{{-}\beta+i\,y} \dd y\Bigr] \nonumber\\
&= {-}\frac{1}{\pi}\,\int_0^\infty \frac{\beta}{\beta^2+y^2}\,{\rm ln}\,(1-\exp({-}\tfrac12 \beta^2- \tfrac12 y^2)) \dd y\nonumber\\
\label{z4}
&= {-}\frac{1}{\pi}\,\int_0^\infty \frac{\beta/\sqrt{2}}{\tfrac12\beta^2+t^2}\,{\rm ln}\,(1-\exp({-}\tfrac12 \beta^2-t^2)) \dd t,
\end{align*}
where we used that
\begin{equation*}
\label{z5}
{\rm Re }\Bigl[\frac{1}{{-}\beta+i\, y}\Bigr] = \frac{{-}\beta}{\beta^2+y^2},
\end{equation*}
together with the substitution $y=t\sqrt{2}$. For $n=1,2,\ldots,$ partial integration gives
\begin{align*}
c_n &= \frac{(-1)^n n!}{\pi} \, {\rm Re}\Bigl[\int_0^\infty \frac{{\rm ln}(1-\exp(-\tfrac12\beta^2-\tfrac12 y^2))}{({-}\beta+i\,y)^{n+1}} \dd y\nonumber\\
&= \frac{(-1)^{n-1}(n-1)!}{\pi}\,{\rm Im}\Bigl[\int_0^\infty {\rm ln}(1-\exp(-\tfrac12\beta^2-\tfrac12 y^2))\dd \Bigl(\frac{1}{(-\beta+i\,y)^n}\Bigr)\Bigr]\nonumber\\
\label{z6}
&= {-}\frac{(-1)^{n-1}(n-1)!}{\pi} {\rm Im}\Bigl[ \int_0^\infty \frac{y}{(-\beta+i\,y)^n}\,\frac{\exp(-\tfrac12\beta^2-\tfrac12 y^2)}{1-\exp(-\tfrac12\beta^2-\tfrac12 y^2)}\dd y\Bigr],
\end{align*}
where we have used that
\begin{equation*}
\label{z7}
{\rm Im}\Bigl[\frac{{\rm ln}(1-\exp(-\tfrac12\beta^2-\tfrac12 y^2))}{(-\beta+i\,y)^n}\Bigr]\Bigl|_0^\infty\Bigr. = 0.
\end{equation*}
Using
\begin{equation*}
\label{z8}
\frac{1}{(-\beta+i\,y)^n} = (-1)^n\,\frac{(\beta+i\,y)^n}{(\beta^2+y^2)^n},
\end{equation*}
we then get
\begin{equation*}
\label{z9}
c_n = \frac{(n-1)!}{\pi}\,{\rm Im}\,\Bigl[\int_0^\infty \frac{y(\beta+i\,y)^n}{(\beta^2+y^2)^n}\,\frac{\exp(-\tfrac12\beta^2-\tfrac12 y^2)}{1-\exp(-\tfrac12\beta^2-\tfrac12 y^2)}\dd y\Bigr],
\end{equation*}
which after the substitution of $y=t\sqrt{2}$ gives
\begin{align}
c_1&=\frac{1}{\pi}\,\int_0^\infty \frac{y^2}{\beta^2+y^2}\,\frac{\exp(-\tfrac12\beta^2-\tfrac12 y^2)}{1-\exp(-\tfrac12\beta^2-\tfrac12 y^2)} \dd y \nonumber\\
\label{z10}
&= \frac{\sqrt{2}}{\pi}\,\int_0^\infty \frac{t^2}{\tfrac12 \beta^2+t^2}\,\frac{\exp(-\tfrac12\beta^2-t^2)}{1-\exp(-\tfrac12\beta^2-t^2)}\dd t,
\end{align}
\begin{align*}
c_2&=\frac{2\beta}{\pi}\,\int_0^\infty \frac{y^2}{(\beta^2+y^2)^2}\,\frac{\exp(-\tfrac12\beta^2-\tfrac12 y^2)}{1-\exp(-\tfrac12\beta^2-\tfrac12 y^2)} \dd y\nonumber\\
\label{z11}
&= \frac{\beta\sqrt{2}}{\pi}\,\int_0^\infty \frac{t^2}{(\tfrac12 \beta^2+t^2)^2}\,\frac{\exp(-\tfrac12\beta^2-t^2)}{1-\exp(-\tfrac12\beta^2-t^2)} \dd t.
\end{align*}
\end{proof}

\section{Robust heavy-traffic approximations}
We shall now establish robust heavy-traffic approximations for the canonical case of Gamma-POisson mixutres; see \eqref{r0}.
 As noted earlier, Gamma mixing yields an arrival process that has a negative binomial distribution, which allows us to establish the detailed asymptotic results in the next theorem.

\begin{theorem}\label{saddlepointThm}
Let $a_n,b_n$ and $s_n$ be as in Assumption \ref{as2}. Then the leading order behavior of $\mu_Q$ is given by
\begin{equation}
\label{r1}
\frac{\sqrt{2}\,\beta_n}{\pi}\Bigl(\frac{b_n+\rho_n}{1-\rho_n}\Bigr)\,\int_{0}^\infty \frac{t^2}{\tfrac12\beta^2_n+t^2}\,\frac{\exp({-}\tfrac12\beta^2_n-t^2)}{1-\exp({-}\tfrac12\beta^2_n-t^2)} \dd t\,(1+o(1)),
\end{equation}
where
\begin{equation}
\label{r2}
\beta_n^2 = s_n\Bigl(\frac{1-\rho_n}{b_n+1}\Bigr)^2\Bigl(1+\frac{b_n}{\rho_n}\Bigr).
\end{equation}
Furthermore, the leading order behavior of $\mathbb{P}(Q^{(n)}=0)$ and $\sigma^2_Q$ is given by
\begin{equation*}
\label{r3}
\exp\Bigl[\frac{1}{\pi}\,\frac{b_n+\rho_n}{b_n+1}\,\int_0^\infty \frac{\beta_n/\sqrt{2}}{\tfrac12\beta^2_n+t^2}\,{\rm ln}\,\Bigl(1-\ee^{{-}\tfrac12\beta^2_n-t^2}\Bigr)\dd t\Bigr],
\end{equation*}
and
\begin{equation}
\label{r4}
\frac{\beta_n^3/\sqrt{2}}{\pi}\Bigl(\frac{b_n+\rho_n}{1-\rho_n}\Bigr)^2\Bigl(\frac{b_n+1}{b_n+\rho_n}+1\Bigr)\int_0^\infty \frac{t^2}{(\tfrac12 \beta_n+t^2)^2}\, \frac{\exp({-}\tfrac12\beta_n-t^2)}{1-\exp({-}\tfrac12\beta_n^2-t^2)}\dd t,
\end{equation}
respectively.
\end{theorem}

The proof of Theorem \ref{saddlepointThm} requires asymptotic evaluation of the Pollaczek-type integrals \eqref{t6}-\eqref{t8}, for which shall use the \textit{non-standard} saddle point method---originally proposed by \cite{debruijn} and also applied in Chapter 2 of this thesis---to turn these contour integrals into practical approximations. 

In contrast to the setting of Chapter 2, both the relevant saddle point and the analyticity radius tend to one as $n\to\infty$, which is a singular point of the integrand, in the setting with overdispersion.
For the proof of Theorem \ref{saddlepointThm}, we therefore modify the special saddle point method developed in Chapter 2 to account for this circumstance.

\begin{proof} 
Our starting point is the probability generating function of the number of arrivals per time slot, given in \eqref{r0}, which is analytic for $|z|<1+1/b_n=:r$. Under Assumption \ref{as2}, we consider $\mu_Q$ as given in  \eqref{t7}. We set
\begin{equation}
\label{a7}
g(z) = -{\rm ln }\,z+\frac{1}{s_n}\,{\rm ln }\Bigl[\tilde{A}^{(n)}(z)\Bigr] = -{\rm ln }\,z - \frac{a_n}{s_n}\,{\rm ln }\left(1+(1-z)b_n\right),
\end{equation}
to be considered in the entire complex plane with branch cuts $(-\infty,0]$ and $[r,\infty)$. The relevant saddle point $z_{\rm sp}$ is the unique zero $z$ of $g'(z)$ with $z\in(1,r_0)$. Since
\begin{equation}
\label{a8}
g'(z) = -\frac{1}{z} + \frac{\rho_n}{1+(1-z)b_n},
\end{equation}
this yields,
\begin{equation}
\label{a9}
1+(1-z_{\rm sp})b_n = \rho_n z_{\rm sp},\quad {\rm i.e., } \quad z_{\rm sp} = 1+\frac{1-\rho_n}{\rho_n+b_n}.
\end{equation}
We then find 
\begin{equation}
\label{a10}
\mu_Q = \frac{s_n}{2\pi i} \int_{|z| = 1+\eps} \frac{g'(z)}{z-1}\,\frac{\exp(s_n\,g(z))}{1-\exp(s_n\,g(z))}\dd z,
\end{equation}
and we take here $1+\eps = z_{\rm sp}$. There are no problems with the branch cuts since we consider $\exp(s_ng(z))$ with integer $s_n$. \\

We continue as in Chapter 2, Section 3 and thus we intend to substitute $z=z(v)$ in the integral in \eqref{a10}, where $z(v)$ satisfies
\begin{equation*}
\label{k1}
g(z(v)) = g(z_{\rm sp})-\tfrac12\,v^2\,g''(z_{\rm sp}) =: q(v)
\end{equation*}
on a range ${-}\tfrac12\delta_n \leq v\leq \tfrac12 \delta_n$ with $\delta_n \to 0$ as $n\to\infty$. 
Note that, this range depends on $n$, whereas these bounds $\pm \tfrac{1}{2} \delta_n$ remained bounded away from zero in \cite{Janssen2015}. 
This severely complicates the present analysis. 
We consider the approximate representation
\begin{equation}
\label{k2}
\frac{-s_n\,g''(z_{\rm sp})}{2\pi i}\int_{-\tfrac12 \delta_n}^{\tfrac12 \delta_n}\frac{v}{z(v)-1}\,\frac{\exp(s_n\,q(v))}{1-\exp(s_n\, q(v))} \dd v
\end{equation}
of $\mu_Q$. We have to operate here with additional care, since both the analyticity radius $r=1+1/b_n$ and the saddle point $z_{\rm sp}$ outside zero $r_0$ tend to 1 as $n\rightarrow\infty$. Specifically, proceeding under the assumptions that $(1-\rho_n)^2a_n$ is bounded while $a_n\rightarrow\infty$ and $b_n\geq 1$, see Assumption \ref{as2}, we have from \eqref{a9} that
\begin{equation}\label{a19}
z_{\rm sp}-1=\frac{1-\rho_n}{b_n+\rho_n} = \frac{1-\rho_n}{b_n} + O\Bigl(\frac{1-\rho_n}{b^2_n}\Bigr),
\end{equation}
where the $O$-term is small compared to $(1-\rho_n)/b_n$ when $b_n\rightarrow\infty$. Next, we approximate $r_0$, using that $r_0>1$ satisfies
\begin{equation*}
\label{a20}
{-}{\rm ln}\, r_0 - \frac{\rho_n}{b_n}\, {\rm ln}\,(1+(1-r_0)b_n) = 0.
\end{equation*}
Write $r_0 = 1+u/b_n$, so that we get the equation
\begin{align*}
0 &= {-}{\rm ln}\,\left(1+\frac{u}{b_n}\right) - \frac{\rho_n}{b_n}\,{\rm ln }(1-u)\nonumber \\
\label{a21}
&= {-}\frac{u}{b_n}\Bigl(1-\rho_n-\tfrac12\Bigl(\frac{1}{b_n}+\rho_n\Bigr)u-\tfrac{1}{3}\Bigl(\frac{-1}{b^2_n}+\rho_n\Bigr)u^2+\cdots\Bigr),
\end{align*}
where we have used the Taylor expansion of ${\rm ln}(1+x)$ at $x=0$. Thus we find
\begin{equation*}
\label{a22}
u=\frac{2(1-\rho_n)}{\rho_n+1/b_n}+O(u^2) = 2(1-\rho_n)+O((1-\rho_n)^2)+O\Bigl(\frac{1-\rho_n}{b_n}\Bigr),
\end{equation*}
and so,
\begin{equation*}
\label{a23}
r_0 = 1+2\,\frac{1-\rho_n}{b_n}+O\Bigl(\frac{(1-\rho_n)^2}{b_n}\Bigr) + O\Bigl(\frac{1-\rho_n}{b^2_n}\Bigr).
\end{equation*}
In \eqref{k2} we choose $\delta_n$ so large that the integral has converged within exponentially small error using $\pm\delta_n$ as integration limits, and, at the same time, so small that there is a convergent power series 
\begin{equation}
\label{a26}
z(v) = z_{\rm sp}+iv+ \sum_{k=2}^\infty c_k(iv)^k, \qquad \text{for } |v| \leq \tfrac12 \delta_n.
\end{equation}
To achieve these goals, we supplement the information on $g(z)$, as given by $\eqref{a7}-\eqref{a9}$, by
\begin{equation}
\label{a27}
g''(z)=\frac{1}{z^2}+\frac{\rho_nb_n}{(1+(1-z)b_n)^2},\quad g''(1) = 1+\rho_nb_n,\quad g''(z_{\rm sp}) =\frac{1}{z_{\rm sp}^2}\Bigl(1+\frac{b_n}{\rho_n}\Bigr).
\end{equation}
Now
\begin{equation*}
\label{a36}
\exp(s_n\,q(v)) = \exp(s_n\,g(z_{\rm sp}))\exp(-\tfrac12\,s_n\,g''(z_{\rm sp})\,v^2),
\end{equation*}
and
\begin{equation*}
\label{a37} s_n\, g''(z_{\rm sp})v^2 = s_n\,b_nv^2(1+o(1)) = a_n(b_n\,v)^2(1+o(1)).
\end{equation*}
Therefore, \eqref{k2} approximates $\mu_Q$ with exponentially small error when we take $\tfrac12 \delta_n$ of the order $1/b_n$.

We next aim at showing that we have a power series for $z(v)$ as in \eqref{a26} that converges for $|v|\leq\tfrac12\delta_n$ with $\tfrac12\delta_n$ of the order $1/b_n$.

\begin{lemma}
Let
\begin{equation*}
\label{a38}
r_n:=\frac{1}{2\,b_n}-(z_{\rm sp} -1 ),\quad m_n:= \tfrac{2}{3}\rho_nr_n\sqrt{\frac{b_n+\rho_n^{-1}}{b_n+\rho_n}},
\end{equation*}
where we assume $r_n>0$. Then \eqref{a26} holds with real coefficients $c_k$ satisfying
\begin{equation}
\label{a39}
|c_k|\leq\frac{r_n}{m_n^k},\quad k=2,3,\ldots.
\end{equation}
\end{lemma}
\begin{proof}
We let
\begin{equation}
\label{a40}
G(z):=\frac{2(g(z)-g(z_{\rm sp}))}{g''(z_{\rm sp})(z-z_{\rm sp})^2}.
\end{equation}
Then $G(z_{\rm sp})=1$ and so we can write \eqref{k1} as
\begin{equation}
\label{a41}
F(z):=(z-z_{\rm sp})\sqrt{G(z)} = i v
\end{equation}
when $|z-z_{\rm sp}|$ is sufficiently small. Since $F(z_{\rm sp})=0$, $F'(z_{\rm sp})=1$, the B\"urmann-Lagrange inversion theorem implies validity of a power series as in \eqref{a26}, with real $c_k$ since $G(z)$ is positive and real for real $z$ close to $z_{\rm sp}$. We therefore just need to estimate the convergence radius of this series from below.

To this end, we start by showing that
\begin{equation}
\label{a42}
{\rm Re}[g''(z)] > \frac{4}{9}\,\rho_n^2\frac{b_n+\rho_n^{-1}}{b_n+\rho_n},\quad |z-z_{\rm sp}|\leq r_n.
\end{equation}
For this, we consider the representation
\begin{equation}
\label{a43}
G(z) = 2\int_{0}^1\int_0^1 \frac{g''(z_{\rm sp}+s\,t(z-z_{\rm sp}))}{g''(z_{\rm sp})} \,t\dd s\dd t.
\end{equation}
We have for $\zeta\in\mathbb{C}$ and $|\zeta-1|\leq 1/2b_{n}\leq 1/2$ from \eqref{a27} that
\begin{equation}
\label{a44}
{\rm Re}[g''(\zeta)] = {\rm Re}(1/\zeta^2) + \rho_nb_n\,{\rm Re}\Bigl[\Bigl(\frac{1}{1+(1-\zeta)b_n}\Bigr)^2\Bigr]\geq \tfrac{4}{9}(1+\rho_nb_n).
\end{equation}
To show the inequality in \eqref{a44}, it suffices to show that
\begin{equation}
\label{a45}
\min_{|\xi-1|\leq 1/2} {\rm Re}\Bigl(\frac{1}{\xi^2}\Bigr) = \frac{4}{9}.
\end{equation}
The minimum in \eqref{a45} is assumed at the boundary $|\xi-1|=1/2$, and for a boundary point $\xi$, we write
\begin{equation*}
\label{a46}
\xi= 1+\tfrac12\cos\theta+\tfrac12 i \sin\theta, \quad 0\leq \theta\leq 2\pi,
\end{equation*}
so that
\begin{equation*}
\label{a47}
{\rm Re}\Bigl(\frac{1}{\xi^2}\Bigr) = \frac{1+\cos\theta+\tfrac{1}{4}\cos 2\theta}{(\tfrac{5}{4}+\cos\theta)^2}.
\end{equation*}
Now
\begin{equation*}
\label{a48}
\frac{\dd}{d\theta} \Bigl[\frac{1+\cos\theta+\tfrac{1}{4}\cos2\theta}{(\tfrac{5}{4}+\cos\theta)^2}\Bigr] = \frac{\sin \theta\,(1-\cos \theta)}{4(\tfrac{5}{4}+\cos\theta)^3}
\end{equation*}
vanishes for $\theta=0,\pi,2\pi$, where ${\rm Re}(1/\xi^2)$ assumes the values $4/9$, 4, 4/9, respectively. This shows \eqref{a45}.

We use \eqref{a45} with $\xi=\zeta$ and with $\xi=1+(1-\zeta)b_n$, with
\begin{equation}
 \label{a49}
 \zeta = \zeta(s,t) = z_{\rm sp} + s t\,(z-z_{\rm sp}),\quad 0\leq s,\, t\leq 1,
 \end{equation}
 where we take $\zeta$ such that $|\zeta-1|\leq 1/2b_n$. It is easy to see from
 $1<z_{\rm sp}<1+1/2b_n$ that $|\zeta-1|\leq 1/2b_n$ holds when $|z-z_{\rm sp}|\leq r_n=1/2b_n-(z_{\rm sp}-1)$. We have, furthermore, from \eqref{a9} that $0<g''(z_{\rm sp})\leq 1+b_n/\rho_n$. Using this, together with \eqref{a44} where $\zeta$ is as in \eqref{a49}, yields
 \begin{equation*}
 \label{a50}
 {\rm Re}[G(z)] \leq \frac{4}{9}\,\frac{1+\rho_nb_n}{1+b_n/\rho_n}\,2\,\int_0^1\int_0^1 t\,\dd s\dd t = \tfrac{4}{9}\,\rho_n^2\,\frac{b_n+\rho_n^{-1}}{b_n+\rho_n}
 \end{equation*}
when $|z-z_{\rm sp}|\leq r_n$, and this is \eqref{a42}.
We therefore have from \eqref{a41} that
\begin{equation*}
\label{a51}
|F(z)|>r_n\cdot\frac{2}{3}\rho_n\sqrt{\frac{b_n+\rho_n^{-1}}{b_n+\rho_n}} = m_n,\quad |z-z_{\rm sp}|=r_n.
\end{equation*}
Hence, for any $v$ with $|v|\leq m_n$, there is exactly one solution $z=z(v)$ of the equation $F(z)-iv=0$ in $|z-z_{\rm sp}|\leq r_n$ by Rouch\'e's theorem. This $z(v)$ is given by
\begin{equation*}
\label{a52}
z(v) = \frac{1}{2\pi i}\,\int_{|z-z_{\rm sp}|=r_n} \frac{F'(z)\,z}{F(z)-iv}\dd z,
\end{equation*}
and depends analytically on $v$, $|v|\leq m_n$. From $|z(v)-z_{\rm sp}|\leq r_n$, we can finally bound the power series coefficients $c_k$ according to
\begin{equation*}
\label{a53}
|c_k| = \Bigl|\frac{1}{2\pi i}\int_{|iv|=m_n} \frac{z(v)-z_{\rm sp}}{(iv)^{k+1}}\dd(iv)\Bigr| \leq \frac{r_n}{m_n^k},
\end{equation*}
and this completes the proof of the lemma.
\end{proof}

\begin{remark}
We have $z_{\rm sp}-1=o(1/b_n)$, see \eqref{a19},  and so
\begin{equation*}
\label{a54}
r_n = \frac{1}{2b_n}(1+o(1)),\quad m_n = \frac{1}{3b_n}(1+o(1)),
\end{equation*}
implying that the radius of convergence of the series in \eqref{a26} is indeed of order $1/b_n$ (since we have assumed $b_n\geq 1$).
\end{remark}

We let $\delta_n=m_n$, and we write for $0\leq v\leq \tfrac12\delta_n$
\begin{equation*}
\label{a55}
\frac{v}{z(v)-1}+\frac{{-}v}{z({-}v)-1} = \frac{-2iv\,{\rm Im}(z(v))}{|z(v)-1|^2},
\end{equation*}
where we have used that all $c_k$ are real, so that $z(-v)=z(v)^*$, where $ ^*$ denotes the complex conjugate. Now from \eqref{a39} and realness of the $c_k$, we have
\begin{equation}
\label{a56}
{\rm Im}(z(v)) = v+\sum_{l=1}^\infty c_{2l+1}(-1)^l\,v^{2l+1} = v+O(v^3),
\end{equation}
and in similar fashion
\begin{equation}
\label{a57}
|z(v)-1|^2 = (z_{\rm sp}-1)^2+v^2+O((z_{\rm sp}-1)^2v^2) + O(v^4)
\end{equation}
when $0\leq v\leq \tfrac12\delta_n$. The order terms in \eqref{a56}-\eqref{a57} are negligible in leading order, and so we get for $\mu_{Q^{(n)}}$ via \eqref{k2} the leading order expression
\begin{equation*}
\label{a58}
\frac{{-}s_n\,g''(z_{\rm sp})}{2\pi i}\,\int_0^{\tfrac12\delta_n}\frac{{-}2iv^2}{(z_{\rm sp}-1)^2+v^2}\,\frac{\exp(s_n\,q(v))}{1-\exp(s_n\, q(v))}\dd v.
\end{equation*}
We finally approximate $q(v) = g(z_{\rm sp})-\tfrac12 g''(z_{\rm sp})v^2$.
There is a $z_1$, $1\leq z_1\leq z_{\rm sp}$ such that
\begin{equation*}
\label{a59}
g(z_{\rm sp}) = {-}\tfrac12(z_{\rm sp}-1)^2\,g''(z_1),
\end{equation*}
and, see \eqref{a19} and \eqref{a27},
\begin{equation*}
\label{a60}
g''(z_1) = g''(z_{\rm sp}) + O((1-\rho_n)b_n).
\end{equation*}
Hence
\begin{align}
s_n\,q(v) &= {-}\tfrac12 s_n\,g''(z_{\rm sp})\,[(z_{\rm sp}-1)^2+v^2]+O((1-\rho_n)b_ns_n(z_{\rm sp}-1)^2)\nonumber\\
&= {-}\tfrac12 s_n\,g''(z_{\rm sp})[(z_{\rm sp}-1)^2+v^2]+O((1-\rho_n)^2a_n),\label{a61}
\end{align}
where \eqref{a19} has been used and $a_nb_n = s_n(1+o(1))$ Therefore, the $O$-term in \eqref{a61} tends to 0 by our assumption that $(1-\rho_n)^2a_n$ is bounded. Thus, we get for $\mu_{Q^{(n)}}$ in leading order
\begin{equation}\label{a62}
\frac{s_n g''(z_{\rm sp})}{\pi} \int_{0}^{\tfrac12\delta_n}\frac{v^2}{(z_{\rm sp}-1)^2+v^2}\,
\frac{\exp(-\tfrac12 g''(z_{\rm sp})s_n((z_{\rm sp}-1)^2+v^2))}{1-\exp(-\tfrac12 g''(z_{\rm sp})s_n((z_{\rm sp}-1)^2+v^2))} \dd v,
\end{equation}
When we substitute $t=v\sqrt{s_n\,g''(z_{\rm sp})/2}$ and extend the integration in \eqref{a62} to all $t\geq 0$ (at the expense of an exponentially small error), we get for $\mu_{Q^{(n)}}$ in leading order
\begin{equation*}
\label{a63}
\frac{1}{\pi}\,\sqrt{2\,s_n\,g''(z_{\rm sp})}\,\int_{0}^\infty \frac{t^2}{\tfrac12\beta_n^2}\,\frac{\exp({-}\tfrac12\beta^2_n-t^2)}{1-\exp({-}\tfrac12\beta^2_n-t^2)}\dd t,
\end{equation*}
where
\begin{equation*}
\label{a64}
\beta^2_n = s_n\,g''(z_{\rm sp})(z_{\rm sp}-1)^2.
\end{equation*}
Now using \eqref{a9} and \eqref{a27}, we get the result of Theorem \ref{saddlepointThm}. A separate analysis of $\beta_n$ is provided in Subsection \ref{convRobust}.

\section{Numerical \& empirical studies}

A similar analysis, modeled after the one given in Chapter 2 gives under Assumption \ref{as1} the leading-order expression
\begin{equation}
\label{a65}
\frac{1}{z_{\rm sp} \pi}\int_0^\infty \frac{\beta_n/\sqrt{2}}{\tfrac12\beta_n^2+t^2}\,{\rm ln}(1-e^{-\tfrac12\beta_n^2-t^2})\dd t
\end{equation}
for ${\rm ln}\,\mathbb{P}(Q^{(n)}=0)$. Observe that the quantity in \eqref{a65} is negative, but bounded away from ${-}\infty$ when $\beta_n$ is bounded away from 0.
Furthermore, we find for the variance of $Q^{(n)}$ the approximation
\begin{equation*}
\label{a66}
\frac{\beta_n^3/\sqrt{2}}{\pi}\frac{z_{\rm sp}+1}{(z_{\rm sp}-1)^2}\int_0^\infty \frac{t^2}{(\tfrac12 \beta_n+t^2)^2}\, \frac{\exp({-}\tfrac12\beta_n-t^2)}{1-\exp({-}\tfrac12\beta_n^2-t^2)}\dd t.
\end{equation*}
\end{proof}
\noindent
Note that we can write \eqref{r1} as
\begin{equation*}
\label{ra1}
\mu_Q \approx \tilde{\sigma}_n\,\E[ M_{\beta_n}]\quad \text{and}\quad \sigma^2_Q \approx \tilde{\sigma}^2_n\, \Var M_{\beta_n}
\end{equation*}
with
\begin{equation}
\label{ra5}
\tilde{\sigma}_n = \beta_n \Bigl(\frac{b_n+\rho_n}{1-\rho_n}\Bigr).
\end{equation}

This robust approximation for $\mu_Q$ is suggestive of the following two properties that extend beyond the mean system behavior, and hold at the level of approximating the queue by $\sigma_n$ times the Gaussian random walk:

\begin{itemize}
\item[\rm (i)] At the process level, the space should be normalized with $\sigma_n$, as in \eqref{mm7}. The approximation \eqref{r1} suggests that it is better to normalize with $\tilde{\sigma}_n$. Although $\tilde \sigma_n\to\sigma_n$ for $n\to\infty$, the $\tilde \sigma_n$ is expected to lead to sharper approximations for finite $n$.
\item[\rm (ii)] Again at the process level, it seems better to replace the original hedge $\beta$ by the robust hedge $\beta_n$. This thus means that the original system for finite $n$ is approximated by a Gaussian random walk with drift $-\beta_n$. Apart from this approximation being asymptotically correct for $n\to \infty$, it is also expected to approximate the behavior better for finite $n$.
\end{itemize}

\subsection{Convergence of the robust hedge\label{convRobust}}
We next examine the accuracy of the heavy-traffic approximations for $\mu_Q$ and $\sigma^2_Q$, following Corollary \ref{abateThm} and Theorem \ref{saddlepointThm}. We expect the robust approximation to be considerably better than the classical approximation when $\beta_n$ and $\tilde{\sigma}_n$ differ substantially from their limiting counterparts. Before substantiating this claim numerically, we present a result on the convergence rates of $\beta_n$ to $\beta$ and $\tilde{\sigma}_n$ to $\sigma_n$.

\begin{proposition}\label{gammanProp}
Let $a_n,b_n$ and $s_n$ as in Assumption \ref{as2}. Then
\begin{equation}
\label{r3a}
\beta_n^2 = \beta^2\Bigl(1 - \frac{1}{1+b_n+\sigma_n/\beta}\Bigr).
\end{equation}
\end{proposition}
\begin{proof}
From \eqref{r2}, we have
\begin{align*}
\beta_n^2 &= s_n\Bigl(\frac{1-\rho_n}{b_n+1}\Bigr)^2\Bigl(1+\frac{b_n}{\rho_n}\Bigr)= \frac{1}{s_n}\Bigl(\frac{s_n-a_nb_n}{b_n+1}\Bigr)^2\Bigl(1+\frac{s_n}{a_n}\Bigr)\nonumber\\
\label{x1}
&= \frac{1}{s_n}\frac{\beta^2\,a_nb_n(b_n+1)}{(b_n+1)^2}\Bigl(1+\frac{s_n}{a_n}\Bigr) = \beta^2\,\frac{b_n}{b_n+1}\,\Bigl(1+\frac{a_n}{s_n}\Bigr) =:\beta^2\,\bar{F}_n.
\end{align*}
Now,
\begin{align*}
\bar{F_n} &= \frac{b_n}{b_n+1}\,\Bigl(1+\frac{a_n}{s_n}\Bigr) = \frac{b_n}{b_n+1}+\frac{1}{b_n+1}\,\frac{a_nb_n}{s_n}\nonumber\\
&= 1-\frac{1}{b_n+1}\,\Bigl(1-\frac{a_nb_n}{s_n}\Bigr) = 1-\frac{1}{b_n+1}\,\frac{\beta\,\sigma_n}{s_n}\nonumber\\
&= 1-\frac{1}{b_n+1}\,\frac{1}{1+\frac{\mu_n}{\beta\sigma_n}} = 1-\frac{1}{b_n+1+\frac{1}{\beta}\sqrt{a_nb_n(b_n+1)}},
\end{align*}
which together with $\sigma_n^2=a_nb_n(b_n+1)$ proves the proposition.
\end{proof}
Note that $\beta_n$ always approaches $\beta$ from below. Also, \eqref{r3a} shows that $b_n$ is the dominant factor in determining the rate of convergence of $\beta_n$.

\begin{proposition}\label{sigmanProp}
Let $\tilde{\sigma}_n$ as in \eqref{ra5}. Then
\begin{equation*}
\tilde{\sigma}_n = \sigma_n + b_n\beta_n + O(1).
\end{equation*}
\end{proposition}
\begin{proof}
Straightforward calculations give
\begin{align*}
\tilde{\sigma}_n &= \beta_n\,\Bigl(\frac{s_nb_n+a_nb_n}{s_n-a_nb_n}\Bigr) \nonumber\\
&= \frac{\beta_n}{\beta}\,\frac{b_n}{\sigma_n}\,(s_n+a_n)
= \frac{\beta_n}{\beta}\,\sqrt{\frac{b_n}{a_n(b_n+1)}}\left(a_n(b_n+1)+\beta\sqrt{a_nb_n(b_n+1)}\right)\nonumber\\
&= \frac{\beta_n}{\beta}\left(\sqrt{a_nb_n(b_n+1)}+\beta b_n\right) = \frac{\beta_n}{\beta}\,\sigma_n + \beta_n b_n.
\end{align*}
Applying Proposition \ref{gammanProp} together with the observation
\begin{equation*}
\sigma_n \sqrt{1 - \frac{1}{1+b_n+\sigma_n/\beta}} = \sigma_n(1 + O(1/\sqrt{a_n}b_n)) = \sigma_n + O(1)
\end{equation*}
yields the result.
\end{proof}

In Figure \ref{fig:convHedge}, we visualize the convergence speed of both parameters in case $\mu_n=n$, $\sigma_n = n^\delta$ with $\delta=0.7$ and $\beta=1$. This implies $a_n = n/(n^{2\delta}-1)$ and $b_n = n^{2\delta}-1$.

\begin{figure}
\centering
\begin{subfigure}{0.49\textwidth}
\centering
\begin{tikzpicture}[scale = 0.75]
\begin{axis}[
	xmin = 0,
	xmax = 200,
	ymin = 0,
	ymax = 1.05,
	xlabel = {$x$},
	ylabel = {$\tilde{\beta}_n/\beta_n$},
	y label style={at={(-0.09,0.75)}},
	axis line style={->},
	axis lines = left,
	legend cell align=left,
	legend style = {at = {(axis cs: 195,0.052)},anchor = south east},
	yscale = 0.8,
	xscale = 1
]

\addplot[thick,col1] table[x=n,y=d06] {Chapter_3/tikz/gamman.txt};
\addplot[thick,col4] table[x=n,y=d075] {Chapter_3/tikz/gamman.txt};
\addplot[thick,col5] table[x=n,y=d09] {Chapter_3/tikz/gamman.txt};

\addplot[dashed] coordinates { (0,1) (200,1) };
\legend{$\delta = 0.6$, $\delta=0.75$, $\delta=0.9$};
\end{axis}
\end{tikzpicture}

\caption{Convergence of $\beta_n$.}
\end{subfigure}
\begin{subfigure}{0.49\textwidth}
\centering
\begin{tikzpicture}[scale = 0.75]
\begin{axis}[
	xmin = 0,
	xmax = 200,
	ymin = 0,
	ymax = 1.7,
	xlabel = {$x$},
	ylabel = {$\tilde{\sigma}_n/\sigma_n$},
	y label style={at={(-0.09,0.75)}},
	axis line style={->},
	axis lines = left,
	legend cell align=left,
	legend style = {at = {(axis cs: 195,0.1)},anchor = south east},
	yscale = 0.8,
	xscale = 1
]

\addplot[thick,col1] table[x=n,y=d06] {Chapter_3/tikz/sigman.txt};
\addplot[thick,col4] table[x=n,y=d075] {Chapter_3/tikz/sigman.txt};
\addplot[thick,col5] table[x=n,y=d09] {Chapter_3/tikz/sigman.txt};

\addplot[dashed] coordinates { (0,1) (200,1) };
\legend{$\delta = 0.6$, $\delta=0.75$, $\delta=0.9$};
\end{axis}
\end{tikzpicture}

\caption{Convergence of $\tilde{\sigma}_n$.}
\end{subfigure}
\caption{}
\label{fig:convHedge}
\end{figure}

We observe that $\beta_n$ starts resembling $\beta$ fairly quickly, as predicted by Proposition \ref{gammanProp}; $\tilde{\sigma}_n$, on the other hand, converges extremely slowly to its limiting counterpart. Since $\mu_Q$ and $\sigma^2_Q$ are approximated by $\tilde{\beta}_n$ and $\tilde{\sigma}_n^2$, multiplied by a term that remains almost constant as $n$ grows, the substitution of $\sigma_n$ by $\tilde{\sigma}_n$, is essential for obtaining accurate approximations, as we illustrate further in the next subsection.

\subsection{Comparison between heavy-traffic approximations}
We set $\mu_n=n$ and $\sigma^2_n=n^{2\delta}$ with $\delta>\tfrac{1}{2}$, so that  $s_n = n+\beta n^{\delta}$, and $a_n =n/(n^{2\delta-1}-1)$ and $b_n = n^{2\delta-1}-1$.
\begin{table}
\centering
\begin{tabular}{|r|r|rrr|rrr|}
\hline
$s_n$     & $\rho_n$ & $\mu_Q$ & \eqref{h1} & \eqref{r1} & $\sigma_Q$ & \eqref{h1b} & \eqref{r4} \bigstrut \\
\hline
5     & 0.609 & 0.343 & 0.246 & 0.363 & 1.002 & 0.835 & 0.978 \bigstrut[t] \\
10    & 0.683 & 0.535 & 0.400 & 0.551 & 1.239 & 1.063 & 1.216 \\
50    & 0.815 & 1.405 & 1.168 & 1.405 & 1.995 & 1.817 & 1.971 \\
100   & 0.855 & 2.113 & 1.824 & 2.105 & 2.445 & 2.270 & 2.420 \\
500   & 0.920 & 5.446 & 5.006 & 5.412 & 3.923 & 3.762 & 3.899 \bigstrut[b] \\
\hline
\end{tabular}
\caption{Numerical results for the Gamma-Poisson case with $\beta=1$ and $\delta=0.6$.}
\label{gammaPoisson1}
\end{table}

\begin{table}
\centering
\begin{tabular}{|r|r|rrr|rrr|}
\hline
$s_n$     & $\rho_n$ & $\mu_Q$ & \eqref{h1} & \eqref{r1} & $\sigma_Q$ & \eqref{h1b} & \eqref{r4}  \bigstrut \\
\hline
5     & 0.550 & 0.462 & 0.284 & 0.479 & 1.162 & 0.896 & 1.130 \bigstrut[t]\\
10    & 0.587 & 0.852 & 0.521 & 0.855 & 1.570 & 1.213 & 1.528 \\
50    & 0.668 & 3.197 & 2.093 & 3.106 & 3.025 & 2.433 & 2.947 \\
100   & 0.700 & 5.561 & 3.784& 5.377 & 3.983 & 3.270 & 3.887\\
500   & 0.766 & 19.887 & 14.741 & 19.202 & 7.514 & 6.455 & 7.361 \bigstrut[b]\\
\hline
\end{tabular}
\caption{Numerical results for the Gamma-Poisson case with $\beta=1$ and $\delta=0.8$.}
\label{gammaPoisson2}
\end{table}

\begin{table}
\centering
\begin{tabular}{|r|r|rrr|rrr|}
\hline
$s_n$     & $\rho_n$ & $\mu_Q$ & \eqref{h1} & \eqref{r1} & $\sigma_Q$ & \eqref{h1b} & \eqref{r4} \bigstrut  \\
\hline
5     & 0.949 & 11.532 & 11.306 & 11.495 & 3.634 & 3.559 & 3.602 \bigstrut[t] \\
10    & 0.961 & 17.565 & 17.268 & 17.548 & 4.474& 4.398 & 4.444 \\
50    & 0.979 & 46.368 & 45.869 & 46.418 & 7.241 & 7.168 & 7.218 \\
100   & 0.984 & 70.340 & 69.735 & 70.430 & 8.910 & 8.839 & 8.888 \\
500   & 0.991 & 184.900 & 183.989 & 185.108 & 14.422 & 14.357 & 14.404 \bigstrut[b]\\
\hline
\end{tabular}
\caption{Numerical results for the Gamma-Poisson case with $\beta=0.1$ and $\delta=0.6$.}
\label{gammaPoisson3}
\end{table}

\begin{table}
\centering
\begin{tabular}{|r|r|rrr|rrr|}
\hline
$s_n$     & $\rho_n$ & $\mu_Q$ & \eqref{h1} & \eqref{r1} & $\sigma_Q$ & \eqref{h1b} & \eqref{r4}  \bigstrut \\
\hline
5     & 0.931 & 15.730 & 15.209 & 15.909 & 4.276 & 4.127 & 4.233 \bigstrut[t]\\
10    & 0.939 & 27.561 & 26.672 & 27.958 & 5.652 & 5.466 & 5.605 \\
50    & 0.955 & 100.660 & 97.967 & 102.070 & 10.760 & 10.476 & 10.698 \\
100   & 0.961 & 175.591 & 171.360 & 177.818 & 14.189 & 13.855 & 14.117 \\
500   & 0.971 & 638.097 & 626.346 & 644.105 & 26.963 & 26.490 & 26.864 \bigstrut[b]\\
\hline
\end{tabular}
\caption{Numerical results for the Gamma-Poisson case with $\beta=0.1$ and $\delta=0.8$.}
\label{gammaPoisson4}
\end{table}
Tables \ref{gammaPoisson1}-\ref{gammaPoisson4} present numerical results for various parameter values.
In these tables, we fixed $s_n$ to integer values, and use the associated value of $n$ in our calculations.
 The exact values of the performance measures are calculated using the method in Appendix \ref{numprocs}. 
Several conclusions are drawn from these tables. Observe that the heavy-traffic approximations based on the Gaussian random walk, \eqref{h1} and \eqref{h1b}, capture the right order of magnitude for both $\mu_Q$ and $\sigma_Q$. However, the values are off, in particular for small $s_n$ and relatively low $\rho_n := \E[A^{(n)}] / s_n$. The inaccuracy also increases with the level of overdispersion. In contrast, the approximations that follow from Theorem \ref{saddlepointThm}, \eqref{r1} and \eqref{r4} are remarkably accurate. Even for small systems with $s_n = 5$ or 10, the approximations for $\mu_Q$ are within 6$\%$ of the exact value for small $\rho_n$ and within $2\%$ for $\rho_n$ close to 1. For $\sigma_Q^2$, these percentages even reduce to $3\%$ and $1\%$, respectively. For larger values of $s_n$ these relative errors naturally reduce further. Overall, we observe that the approximations improve for heavily loaded systems, and the corrected approximations are particularly useful for systems with increased overdispersion.

\subsection{Capacity allocation in health care}
We next apply our model and robust approximations to real-life patient arrivals. We consider  emergency patients who require diagnostic tests at the radiology department of a hospital. Green \cite{Green2004} points out that patients at the radiology department can be roughly categorized into three groups: Inpatients, outpatients and emergency patients. Inpatient and outpatient arrivals are relatively predictable as these are usually scheduled by appointment. Emergency patients, on the other hand, are inherently unpredictable: They typically require urgent care and therefore timely access to testing facilities, as well as a quick assessment of the test results. This leads to prioritization of emergency patients over the other two groups in such settings, so that they do not experience any delay caused by the groups of lower priority. However, patients from the same top-priority group can still cause considerable congestion. A careful evaluation of capacity allocation is hence required, bearing in mind that additional sophisticated pieces of medical equipment are very costly.

In the setting we study, capacity is defined by the number of time slots available to perform radiology tests on emergency patients in a given time period, which we  set at 24 hours. As radiology tests are commonly performed in appointment slots of fixed length, the number of slots available per day is also indirectly fixed. In terms of our model parameters, see Section \ref{modelSection}, we have $s$ as the number of slots per day allocated to emergency patients, and $A(k)$ the number of test requests received by the department on day $k$. We omit the subscript $n$ in this section due to the absence of limits. Then $\E[Q]$ can be interpreted as the expected number of patients whose test is carried over to the next day. A more natural performance measure in this setting is the expected waiting time, namely the time between the physician's request and the actual start of the test. However, Little's law implies that there is a linear relation between the two, hence we choose to only evaluate $\E[Q]$.

The data set on which our empirical study is based originates from the emergency department of SKHospital, monitored over a period of 76 days from September to November 2012. We extracted information of ED patients referred to the radiology department by the ED physicians, which includes the time the test request was made and the exact test type performed. The two test types, X-ray and CT scans, are performed on different equipment and hence it makes sense to analyze their queueing processes separately.

\begin{figure}
\begin{center}
\begin{subfigure}{0.49\textwidth}
\centering
\begin{tikzpicture}[scale = 0.74]
\begin{axis}[
	xmin = 40,
	xmax = 110,
	ymin = 0,
	ymax = 1,
	xlabel = {$x$},
	ylabel = {$\P(A\leq x)$},
	y label style={at={(-0.1,0.6)}},
	axis line style={->},
	axis lines = left,
	legend cell align=left,
	legend style = {at = {(axis cs: 109,0.05)},anchor = south east},
	yscale = 0.8,
	xscale = 1
]

\addplot[thick] table[x=x,y=emp] {Chapter_3/tikz/xray.txt};
\addplot[thick,col1] table[x=x,y=poisson] {Chapter_3/tikz/xray.txt};
\addplot[thick,col4] table[x=x,y=fitted] {Chapter_3/tikz/xray.txt};
\addplot[thick] table[x=x,y=emp] {Chapter_3/tikz/xray.txt};

\legend{Empirical,Poisson,Gamma-Poisson};
\end{axis}
\end{tikzpicture}
\caption{X-ray}
\end{subfigure}
\begin{subfigure}{0.49\textwidth}
\centering
\begin{tikzpicture}[scale = 0.74]
\begin{axis}[
	xmin = 5,
	xmax = 32,
	ymin = 0,
	ymax = 1,
	xlabel = {$x$},
	ylabel = {$\P(A\leq x)$},
	y label style={at={(-0.1,0.6)}},
	axis line style={->},
	axis lines = left,
	legend cell align=left,
	legend style = {at = {(axis cs: 32,0.05)},anchor = south east},
	yscale = 0.8,
	xscale = 1
]

\addplot[thick] table[x=x,y=emp] {Chapter_3/tikz/ct.txt};
\addplot[thick,col1] table[x=x,y=poisson] {Chapter_3/tikz/ct.txt};
\addplot[thick,col4] table[x=x,y=fitted] {Chapter_3/tikz/ct.txt};
\addplot[thick] table[x=x,y=emp] {Chapter_3/tikz/xray.txt};
\addplot[thick] table[x=x,y=emp] {Chapter_3/tikz/ct.txt};

\legend{Empirical,Poisson,Gamma-Poisson};
\end{axis}
\end{tikzpicture}
\caption{CT scan}
\end{subfigure}
\end{center}
\caption{Empirical, fitted Poisson and fitted Gamma-Poisson cumulative distribution functions of the number of arrivals.}
\label{fig:fittedHospital}
\end{figure}
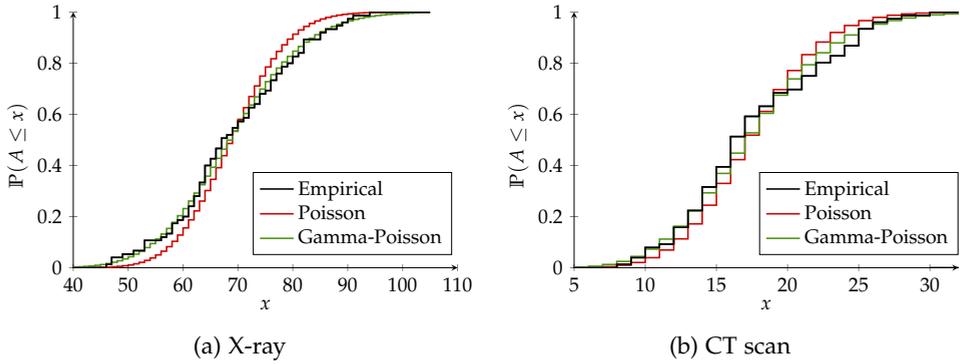

We refer to test requests as arrivals. The empirical cumulative distribution functions of the number of arrivals per day, for each type, are depicted by the black lines in Figure \ref{fig:fittedHospital}. The sample means equal 69.81 and 17.47, for the X-ray and CT scans respectively, whereas the sample variances are 121.8 and 26.12. These values suggest that fitting a Poisson distribution is inappropriate, which is visually backed up by the fitted Poisson cdf, depicted in Figure \ref{fig:fittedHospital} by the red line. To strengthen this claim, we tested both samples for the Poisson assumption using the \emph{dispersion test}, see Appendix \ref{statproc}, and obtained $p$-values equal $7.01\cdot 10^{-3}$ and $3.57\cdot 10^{-3}$ respectively, which allow us to safely reject the Poisson hypothesis in both cases.

In search for a better distributional fit with the arrivals count, we resort to Gamma-Poisson mixtures. Here we employ the procedure in \cite{koolejongbloed}, which is basically a maximum log-likelihood method, to obtain estimates for the parameters $a$ and $b$. This yields
\begin{equation*}
\label{parameterEstimators}
\hat{a}_{\rm X-ray} = 95.68,\quad \hat{b}_{\rm X-ray} = 0.7297,\quad \hat{a}_{\rm CT} = 37.19,\quad \hat{b}_{\rm CT} = 0.4698.
\end{equation*}
Applying the bootstrapping method to the data and the fitted model, also described in the appendix of \cite{koolejongbloed}, returns p-values that equal 0.7354 and 0.2120 for X-ray and CT scans, respectively. Therefore, the null hypothesis, stating that the data originated from a Gamma-Poisson mixture, cannot be rejected. The cdfs of the fitted Gamma-Poisson distributions, plotted in Figure \ref{fig:fittedHospital}, give visual confirmation of this claim as well. 
Naturally, we also compared the estimated densities to the empirical pdf of the data. However, these fail to give a convincing visual fit due to the relatively small sample size and are therefore omitted here.

We now have clear evidence that both the X-ray and CT scan facilities face an overdispersed arrival stream. In our final step of the empirical study we now evaluate the accuracy of our performance measure of interest $\E[Q]$, and the consequences of assessing system performance while ignoring the presence of overdispersion. We take the following approach: Trivially, we need to choose $s> \E[A]$ in order for the system to be stable. Hence, we vary $s$ from 70 to 80 for X-rays and from 18 to 24 for CT scans and simulate the queue length process by sampling the number of requests per day from the actual arrival counts. The number of iterations performed is $10^8$ for each configuration, excluding a warm-up interval of length $10^7$ (days). Around the mean of $Q$ obtained from this simulation, we create a 95\% confidence interval. Next, we approximate the expected stationary queue length under two scaling rules. Assuming Poisson arrivals, the appropriate capacity allocation rule would be $s=\hat{\mu}+\beta\sqrt{\hat{\mu}}$, for some $\beta>0$. Our novel capacity sizing rule prescribes $s = \hat{\mu} + \beta\hat{\sigma} = \hat{a}\hat{b}+\beta\sqrt{\hat{a}\hat{b}(\hat{b}+1)}$. We compute the first approximation based on square-root safety capacity sizing by deducing $\beta$ for each $s$, which we substitute in $\E[Q^{\rm srs}] \approx \sqrt{\hat{\mu}}\,\E[M_{\beta}]$. Similarly, we obtain $\beta$ from the new rule, and plug this value, together with the fitted parameters $\hat{a}$ and $\hat{b}$, into \eqref{r1}. The results are given in Tables \ref{tab:simXRay} and \ref{tab:simCT}. The last column shows the 95\% relative error bound of the second approximation.

\begin{table}[h]
\centering
\begin{tabular}{|r|r|rrrr|r|}
\hline
$s$     & $\rho$   & $\E[Q] \ (\pm $ conf. iv.) & $\E[Q^{\rm srs}]$ & $\eqref{h1}$ & $\eqref{r1}$    & rel. error  \bigstrut \\
\hline
70    & 0.997 & $328.313 \pm\ 6.6\cdot 10^{-2}$ & 186.664 & 324.231 & 325.221 & $9.6\cdot 10^{-3}$  \bigstrut[t]\\
71    & 0.983 & $45.073 \pm\ 1.0\cdot 10^{-2}$ & 24.946 & 45.290 & 45.308 & $5.4\cdot 10^{-3}$ \\
72    & 0.970 & $21.988 \pm\ 5.4\cdot 10^{-3}$ & 11.650 & 21.982 & 22.129 & $6.6\cdot 10^{-3}$ \\
73    & 0.956 & $13.546 \pm\ 3.6\cdot 10^{-3}$ & 6.847 & 13.455 & 13.649 & $7.8\cdot 10^{-3}$ \\
74    & 0.943 & $9.230 \pm\ 2.7\cdot 10^{-3}$ & 4.438 & 9.106 & 9.319 & $1.0\cdot 10^{-2}$ \\
75    & 0.931 & $6.655 \pm\ 2.1\cdot 10^{-3}$ & 3.031 & 6.513 & 6.731 & $1.2\cdot 10^{-2}$ \\
76    & 0.919 & $4.949 \pm\ 1.7\cdot 10^{-3}$ & 2.136 & 4.821 & 5.037 & $1.8\cdot 10^{-2}$ \\
77    & 0.907 & $3.755 \pm\ 1.4\cdot 10^{-3}$ & 1.534 & 3.650 & 3.861 & $2.8\cdot 10^{-2}$ \\
78    & 0.895 & $2.884 \pm\ 1.1\cdot 10^{-3}$ & 1.115 & 2.807 & 3.009 & $4.4\cdot 10^{-2}$ \\
79    & 0.884 & $2.230 \pm\ 1.0\cdot 10^{-3}$ & 0.816 & 2.183 & 2.374 & $6.5\cdot 10^{-2}$ \\
80    & 0.873  & $1.734 \pm\ 8.5\cdot 10^{-4}$ & 0.600 & 1.710 & 1.890 & $9.1\cdot 10^{-2}$ \bigstrut[b]\\
\hline
\end{tabular}%
\caption{Computational results for X-ray.}
\label{tab:simXRay}
\end{table}

\begin{table}
\centering
\begin{tabular}{|r|r|rrrr|r|}
\hline
$s$     & $\rho$   & $\E[Q]\ (\pm $ conf.iv.) & $\E[Q^{\rm srs}]$ & $\eqref{h1}$ & $\eqref{r1}$    & rel. error \bigstrut \\
\hline
18    & 0.970 & 22.116 $\pm\ 4.9\cdot 10^{-3}$ & 14.235 & 21.965 & 21.724 & $1.8\cdot 10^{-2}$  \bigstrut[t] \\
19    & 0.919 & 6.289 $\pm\ 1.7\cdot 10^{-3}$ & 3.640 & 5.941 & 6.040 & 4.0$\cdot 10^{-2}$ \\
20    & 0.873 & 3.101 $\pm\ 1.0\cdot 10^{-3}$ & 1.589 & 2.772 & 2.917 & 6.0$\cdot 10^{-2}$ \\
21    & 0.832 & 1.767  $\pm\ 6.6\cdot 10^{-4}$ & 0.800 & 1.507 & 1.658 & 6.1$\cdot 10^{-2}$ \\
22    & 0.794 & 1.066  $\pm\ 4.6\cdot 10^{-4}$ & 0.425 & 0.874 & 1.016 & 4.7$\cdot 10^{-2}$ \\
23    & 0.760 & 0.653  $\pm\ 3.3\cdot 10^{-4}$ & 0.230 & 0.522 & 0.649 & 7.1$\cdot 10^{-3}$\\
24    & 0.728 & 0.377  $\pm\ 2.3\cdot 10^{-4}$ & 0.124 & 0.315 & 0.424 & 1.2$\cdot 10^{-1}$ \bigstrut[b]\\
\hline
\end{tabular}%
\caption{Computational results for CT scan.}
\label{tab:simCT}
\end{table}

Based on these figures, we make several remarks. First, assuming the conventional regime (neglecting overdispersion) the approximation severely overestimates system performance in both queues. Because the square-root safety margin underestimates the stochastic fluctuations in the arrival process, the safety parameter $\beta$ is overestimated, which leads to a smaller magnitude of the approximated queue length process. This clearly illustrates the distorted view estimated performance characteristics can give under the wrong scaling.
Secondly, it is worth noticing the very good proximity of $\eqref{r1}$ to the values obtained through simulation. As we expected, the quality of the approximation deteriorates with increasing values of $s$. This makes sense, because we assumed the system to be in heavy traffic in the derivation of the formulas. What is surprising, on the other hand, is the fact that the approximation performs very well, even though the parameter $b$ is very small for these particular data sets, while the analysis of Theorem \ref{saddlepointThm} assumes that $a$ and $b$ are large. This demonstrates that the approximation scheme is remarkably robust and is able to capture the pre-limit behavior of these types of queues very well.

\section{Conclusion \& future research}
In this chapter, we proposed an adaptation to the square-root staffing rule for service systems facing overdispersed demand, using the bulk service queue as a vehicle for our analysis. 
Subsequently, we derive two sets of asymptotic approximations for the scaled steady-state queue length moments for large arrival volumes. 
The first set being based on the resemblance with the maximum of a Gaussian random walk, the second set being derived through a non-standard saddle point method, assuming arrivals follow a Gamma-Poisson mixture. 
Numerical experiments indicated that our approximations capture the pre-limit behavior well under different order of overdispersion, and are robust against any parameter estimation errors. 

Although our method provides a robust way to approximate and dimension queues with overdispersed arrival processes, we see some interesting directions for future research. 

First, we accentuate that our model is a discrete time queueing model in which a deterministic amount of workload can we handled within each period. 
This approach allowed us to use Pollaczek's formula as a starting point to obtain more refined asymptotic approximations for the performance indicators of the system.
In case we consider queueing models of birth-death-type, such as the $M/M/s$ queue, in the presence of overdispersion demand, different techniques are required to provide scaling limits and corresponding capacity allocation rules, see e.g.~\cite{maman}.  
Although we expect that, just as in the novel scaling regimes of Chapter 2, the asymptotic behavior of the bulk service queue and the multi-server queueing models to be similar, this needs proper analysis and understanding. 

Second, empirical work, see e.g. \cite{Avramidis:2004}, shows that in real-life settings, demand in consecutive time periods is often positively correlated, rather than independently distributed as assumed in this chapter. 
This correlation structure obviously alters the queue's dynamics and presumably requires an adaptation of the square-root staffing rule as well, making it a worthwhile direction for further analysis.

Last, we have only considered the analysis of the queueing model in steady state. 
Typical service systems however do not face a constant expected arrival rate, nor do they run infinitely long. 
Henceforth, it would be interesting to study the influence of overdispersion on the transient dynamics of the queue and to investigate the capacity allocation problem in scenarios with time-varying demand. 
The theory developed in this chapter may serve as a building block to tackle these more profound questions.

\section*{Appendix}

\addcontentsline{toc}{section}{\hspace{7.1mm} Appendix}

\begin{subappendices}

\settocdepth{chapter}

\section{Proofs of convergence results}
 \label{formalSec}
 
This section presents the details of the proof of Lemma \ref{gaussStep} and Theorem \ref{gaussianThm}, using the random walk perspective of the process $\{Q^{(n)}(k)\}_{k=0}^\infty$. This section is structured as follows. The next two lemmata are necessary for proving the first assertion of Theorem \ref{gaussianThm}, concerning the weak convergence of the scaled process to the maximum of the Gaussian random walk, which is summarized in Proposition \ref{prop6}. The two remaining propositions of this section show convergence of $\hat{Q}^{(n)}$ at the process level as well as in terms of the three characteristics.

Let us first fix some notation:
\begin{equation}
\label{b1}
Y^{(n)}_k := \hat{A}^{(n)}_k-\beta,\quad
S^{(n)}_k = \sum_{i=1}^k Y^{(n)}_i,
\end{equation}
with $S_0^{(n)} = 0$ and $k=1,2,...$. Then \eqref{mm6} can be rewritten as
\begin{equation}
\label{g5a}
 \hat{Q}^{(n)} \equalD \max_{k\geq 0} \Bigl\{{\textstyle \sum}_{i=1}^k Y^{(n)}_i\Bigr\} =: M_\beta^{(n)},
\end{equation}
Last, we introduce the sequence of independent normal random variables $Z_1,Z_2,\ldots$ with mean $\-\beta$ and unit variance 1, and
\begin{equation*}
M_\beta \equalD \max_{k\geq 0} \{{\textstyle \sum}_{i=1}^k Z_i\}.
\end{equation*}

\subsection{Proof of Lemma \ref{gaussStep}}
\begin{proof}
We show weak convergence of the random variable $\hat{A}^{(n)}$, as defined in Section \ref{modelSection}, to a standard normal random variable. Since $\hat{\Lambda}_n$ is asymptotically standard normal, its characteristic function converges pointwise to the corresponding limiting characteristic function, i.e.
\begin{equation}
 \label{g8}
 \lim_{n\rightarrow\infty} \phi_{\hat{\Lambda}_n}(t) = \lim_{n\rightarrow \infty} \ee^{-i\mu_nt/\sigma_n}\,\phi_{\Lambda_n}(t/\sigma_n) = \ee^{{-}t^2/2},\qquad \forall t\in \mathbb{R}.
 \end{equation}
Furthermore, by definition of $A^{(n)}$,
\begin{equation*}
\label{g9}
\phi_{A^{(n)}}(t) = \E\left[ \exp(\Lambda_n(\ee^{it}-1))\right] = \phi_{\Lambda_n}\left(-i(\ee^{it}-1)\right),
\end{equation*}
so that
\begin{equation}
\label{g10}
\phi_{\hat{A}_k^{(n)}}(t) = \ee^{-i\mu_nt/\sigma_n}\,\phi_{A_k^{(n)}}(t/\sigma_n) = \ee^{-i\mu_nt/\sigma_n}\phi_{\Lambda_n}\left(-i(\ee^{it/\sigma_n}-1)\right).
\end{equation}
Now fix $t\in\mathbb{R}$. By using
\begin{equation*}
\label{g11}
-i(\ee^{it/\sigma_n} - 1) = \frac{t}{\sigma_n} -\frac{it^2}{2\sigma_n^2} + O\left(t^3/\sigma_n^3\right),
\end{equation*}
we expand the last term in \eqref{g10},
\begin{equation*}
\label{g12}
\phi_{\Lambda_n}(t/\sigma_n) + \Bigl(-\frac{i\,t^2}{2\sigma_n^2}+O\left(t^3/\sigma_n^3\right)\Bigr)
\phi_{\Lambda_n}'(t/\sigma_n) + O\Bigl(\Bigl(-\frac{i\,t^2}{2\sigma_n^2}+O\left(\frac{t^3}{\sigma_n^3}\right)\Bigr)^2\phi_{\Lambda_n}''\Big(\frac{t}{\sigma_n}\Big)\Bigr)
\end{equation*}
\begin{equation*}
\label{g13}
= \phi_{\Lambda_n}(t/\sigma_n) - \Bigl(\frac{i\,t^2}{2\sigma_n^2}+O\left(t^3/\sigma_n^3\right)\Bigr)
\phi_{\Lambda_n}'(\zeta)
\end{equation*}
for some $\zeta$ such that $|\zeta - t/\sigma_n| < |i(1-\ee^{it/\sigma_n})-t/\sigma_n|$. Also,
\begin{align}
|\phi_{\Lambda_n}'(u)| &= \left|\frac{\d}{\dd u}\int_{-\infty}^\infty \ee^{iux}\dd F_{\Lambda_n}(x)\right| = \left|\int_{0}^{\infty} ix\,\ee^{iux}\dd F_{\Lambda_n}(x)\right| \nonumber\\
\label{g13a}
&\leq \int_{-\infty}^\infty |ix\,\ee^{iux}|\,\dd F_{\Lambda_n}(x) = \int_0^\infty x\dd F_{\Lambda_n}(x) = \mu_n
\end{align}
for all $u\in\mathbb{R}$. Hence, by substituting \eqref{g10},
\begin{align}
\left| \phi_{\hat{A}_k^{(n)}}(t)-\ee^{-i\mu_nt/\sigma_n}\phi_{\Lambda_n}(t/\sigma_n)\right| &= \left|\ee^{-i\mu_nt/\sigma_n}\,\left(\frac{i\,t^2}{2\sigma_n^2}+O(t^3/\sigma_n^3)\right)\,\phi_{\Lambda_n}'(\zeta)\right|\nonumber\\
& \leq \left(\frac{t^2}{2\sigma_n^2}+O(t^3/\sigma_n^3)\right) |\phi_{\Lambda_n}'(\zeta)|\nonumber\\
& = \frac{\mu_n t^2}{\sigma_n^2} + O\left(\frac{\mu_nt^3}{\sigma_n^3}\right),
\label{g13b}
\end{align}
which tends to zero as $n\rightarrow \infty$ by our assumption that $\mu_n/\sigma_n^2\rightarrow 0$.
Finally,
\begin{equation*}
\label{g13c}
\left| \phi_{\hat{A}_k^{(n)}}(t)-\ee^{-\tfrac12 t^2}\right| \leq \left| \phi_{\hat{A}_k^{(n)}}(t)-\ee^{-i\mu_nt/\sigma_n}\phi_{\Lambda_n}(t/\sigma_n)\right| +
\left| \ee^{-i\mu_nt/\sigma_n}\phi_{\Lambda_n}(t/\sigma_n) - \ee^{-\tfrac12 t^2}\right|,
\end{equation*}
in which both terms go to zero for $n\rightarrow \infty$, by \eqref{g8} and \eqref{g13b}. Hence $\phi_{\hat{A}^{(n)}_k}(t)$ converges to $\ee^{{-}t^2/2}$ for all $t\in\mathbb{R}$, so that we can conclude by L\'evy's continuity theorem that $\hat{A}_k^{(n)} \Rightarrowd \mathcal{N}(0,1)$.
\end{proof}
\subsection{Proof of Theorem \ref{gaussianThm}}
To secure convergence in distribution of $\hat{Q}^{(n)}$ to $M_\beta$, i.e. the maximum of a Gaussian random walk with negative drift, the first assertion of Theorem \ref{gaussianThm},
 the following property of the sequence $\{Y_k^{(n)}\}_{n\in\mathbb{N}}$ needs to hold.
\begin{lemma}\label{uilemma}
Let $Y^{(n)}_k$ be defined as in \eqref{b1} with $\mu_n,\sigma_n^2 < \infty$ for all $n\in\mathbb{N}$. Then the sequence $\{(Y_k^{(n)})^+\}_{n\in\mathbb{N}}$ is uniform integrable, i.e.
\begin{equation*}
\label{g14}
\lim_{K\rightarrow\infty}\sup_n \E\Big[Y^{(n)\,+}_k|\mathbbm{1}_{\{|Y^{(n)\,+}_k|\geq K\}}\Big] = 0.
\end{equation*}
\end{lemma}
\begin{proof}
Because the sequence $\{Y^{(n)}_k\}_{k\in\mathbb{N}}$ is i.i.d. for all $n$, we  omit the index $k$ in this proof. First, fix $K>0$ and note that
\begin{equation*}
\label{g15}
\E[|Y^{(n)+}|\mathbbm{1}{\{|Y^{(n)\,+}|\geq K\}}] = \E[Y^{(n)+}\mathbbm{1}{\{Y^{(n)+}\geq K\}}] = \E[Y^{(n)}\mathbbm{1}_{\{Y^{(n)}\geq K\}}].
\end{equation*}
This last expression can be bounded from above using the Cauchy-Schwarz inequality, so that
\begin{equation*}
\label{g16}
\E[Y^{(n)}\mathbbm{1}_{\{Y^{(n)}\geq K\}}] \leq \E[ Y^{(n)\,2}]^{1/2}\,\mathbb{P}(Y^{(n)}\geq K)^{1/2}.
\end{equation*}
By the definition of $Y^{(n)}$, we know $\E [Y^{(n)}] = -\beta$ and $\Var Y^{(n)} = \Var A^{(n)} / \sigma_n^2 = 1$. Using this information, we find
\begin{equation*}
   \label{g17}
   \E[Y^{(n)\,2}] = \Var Y^{(n)} + (\E[Y^{(n)}])^2 = 1+\beta^2
   \end{equation*}
   and
   \begin{align*}
   \mathbb{P}(Y^{(n)}\geq K )&=\mathbb{P}(Y^{(n)}+\beta\geq K+\beta) \leq \mathbb{P}(|Y^{(n)}+\beta|\geq K+\beta)\nonumber\\
   &\leq \frac{\Var Y^{(n)}}{(K+\beta)^2} = \frac{1}{(K+\beta)^2},
   \end{align*}
   where we used Chebyshev's inequality for the last upper bound. Therefore,
 \begin{align*}
\lim_{K\rightarrow \infty} \sup_n \E[|Y^{(n)\,+}|\mathbbm{1}_{\{|Y^{(n)\,+}|\geq K\}}] &=
\lim_{K\rightarrow \infty} \sup_n \E[Y^{(n)}\mathbbm{1}_{\{Y^{(n)}\geq K\}}]\nonumber\\
&\leq \lim_{K\rightarrow \infty} \sup_n \E[Y^{(n)\,2}]^{1/2}\,\mathbb{P}(Y^{(n)}\geq K )^{1/2}\nonumber\\
&\leq \lim_{K\rightarrow \infty} \frac{\sqrt{1+\beta^2}}{K+\beta} = 0.
\end{align*}

\end{proof}
By combining the properties proved in Lemma \ref{gaussStep} and \ref{uilemma} with Assumption \ref{as2}, the next result follows directly by \cite[Thm.~X6.1]{Asmussen2003}.
\begin{proposition}\label{maxRWprop}
Let $\hat{Q}^{(n)}$ as in \eqref{g5a}. Then
\begin{equation*}
 \hat{Q}^{(n)}\Rightarrowd M_\beta,\qquad {\rm as}\ n\rightarrow\infty.
\end{equation*}
\end{proposition}

Although Proposition \ref{maxRWprop} tells us that the properly scaled $Q^{(n)}$ converges to a non-degenerate limiting random variable, it does not cover the convergence of its mean, variance and the empty-queue probability. In order to secure convergence of these performance measures as well, we follow the approach similar to \cite{Sigman2011b}, using Assumptions \ref{as2} and \ref{as3}.

\begin{proposition}\label{prop6}
Let $\hat{Q}^{(n)}$ as in \eqref{g5a}, $\mu_n,\sigma_n^2 \rightarrow \infty$ such that both $\sigma_n^2/\mu_n\rightarrow \infty$ and $\E[\hat{A}^{(n)3}]$ $<\infty$. Then
\begin{align*}
\label{b16}
\mathbb{P}(\hat{Q}^{(n)}= 0)&\rightarrow \mathbb{P}(M_\beta = 0),\\
\E [\hat{Q}^{(n)}]&\rightarrow \E [M_\beta],\\
\Var \hat{Q}^{(n)}&\rightarrow \Var M_\beta,
\end{align*}
as $n\rightarrow\infty$.
\end{proposition}
\proof
First, we recall that $\hat{Q}^{(n)}\equalD M_\beta^{(n)}$ for all $n\in\mathbb{N}$, so that $\mathbb{P}(\hat{Q}^{(n)} = 0) = \mathbb{P}(M_\beta^{(n)}=0)$, $\E[\hat{Q}^{(n)}]=\E[M_\beta^{(n)}]$ and $\Var\,\hat{Q}^{(n)}=\Var\,M_\beta^{(n)}$ as defined in \eqref{b1}. Our starting point is Spitzer's identity, see \cite[p.~230]{Asmussen2003},
\begin{equation}
\label{b17}
\E[\ee^{it M_\beta^{(n)}}] = \exp\Bigl( \sum_{k=1}^\infty \frac{1}{k} (\E[\ee^{it(S^{(n)}_k)^+}]-1)\Bigr),
\end{equation}
with $S^{(n)}_k$ as in \eqref{b1}, and $M_\beta^{(n)}$ the all-time maximum of the associated random walk. Simple manipulations of \eqref{b17} give
\begin{align}
\label{y1}
{\rm ln}\,\mathbb{P}(M_\beta^{(n)} = 0) &=  -\sum_{k=1}^\infty \frac{1}{k}\,\mathbb{P}(S^{(n)}_k > 0),\\
\label{y2}
\E[M_\beta^{(n)}] &= \sum_{k=1}^\infty \frac{1}{k} \E[S^{(n)\,+}_k] = \sum_{k=1}^\infty \frac{1}{k}\int_0^\infty \mathbb{P}(S^{(n)}_k > x) \dd x,\\
\label{y3}
\Var M_\beta^{(n)} &= \sum_{k=1}^\infty \frac{1}{k} \E[(S^{(n)\,+}_k)^2] =\sum_{k=1}^\infty \frac{1}{k}\int_0^\infty \mathbb{P}(S^{(n)}_k > \sqrt{x}) \dd x.
\end{align}
By Lemma \ref{gaussStep}, we know
\begin{equation*}
\label{y4}
\mathbb{P}(S^{(n)}_k > y) = \mathbb{P}\left( {\sum_{i=1}^k} Y^{(n)}_i > y \right) \rightarrow \mathbb{P}\left({\textstyle\sum_{i=1}^k} Z_i > y\right),
\end{equation*}
for $n\rightarrow \infty$, where the $Z_i$'s are independent and identically normally distributed with mean $-\beta$ and variance 1.
Because equivalent expressions to \eqref{y1}-\eqref{y3} apply to the limiting Gaussian random walk, it is sufficient to show that the sums converge uniformly in $n$, so that we can apply dominated convergence to prove the result.

We start with the empty-queue probability. To justify interchangeability of the infinite sum and limit, note
\begin{equation*}
\label{y5}
 \mathbb{P}(S^{(n)}_k > 0) \leq  \mathbb{P}(|S^{(n)}_k+k\beta| > k\beta )\leq \frac{k}{\beta^2k^2} = \frac{1}{\beta^2k},
 \end{equation*}
where we used that $\E[ S^{(n)}_k] = k\E [Y^{(n)}_1] = -k\beta$ and $\Var S^{(n)}_k = k$ and apply Chebychev's inequality, so that
\begin{equation*}
\label{y6}
\sum_{k=1}^\infty \frac{1}{k}\mathbb{P}(S^{(n)}_k > 0) \leq \sum_{k=1}^\infty \frac{1}{\beta^2 k^2} < \infty, \qquad \forall n\in\mathbb{N}.
\end{equation*}
Hence,
\begin{align*}
\lim_{n\rightarrow\infty} {\rm ln}\,\mathbb{P}(\hat{Q}^{(n)}= 0) &= \lim_{n\rightarrow\infty}  - \sum_{k=1}^\infty \frac{1}{k}\mathbb{P}(S^{(n)}_k > 0) = -\sum_{k=1}^\infty \frac{1}{k} \lim_{n\rightarrow\infty}\mathbb{P}(S^{(n)}_k > 0)\nonumber\\
&= -\sum_{k=1}^\infty \frac{1}{k} \mathbb{P}({\textstyle\sum_{i=1}^k} Z_i > 0) = {\rm ln}\,\mathbb{P}(M_\beta = 0).
\end{align*}
Finding a suitable upper bound on $\frac{1}{k}\int_0^\infty \mathbb{P}(\hat{Q}^{(n)}>x) \dd x$ and $\frac{1}{k}\int_0^\infty \mathbb{P}(\hat{Q}^{(n)}>\sqrt{x}) \dd x$ requires a bit more work. We initially focus on the former, the latter follows easily. The following inequality from \cite{Nagaev1979} proves to be very useful:
\begin{equation}
\label{y8}
\mathbb{P}(\bar{S}(k)>y) \leq C_r\,\Bigl(\frac{k\,\sigma^2}{y^2}\Bigr)^2 + k\,\mathbb{P}(X>y/r),
\end{equation}
where $\bar{S}(k)$ is the sum of $k$ i.i.d. random variables distributed as $X$, with $\E[X] = 0$ and $\Var\, X=\sigma^2$, $y > 0$, $r>0$ and $C_r$ a constant only depending on $r$. We take $r=3$ for brevity in the remainder of the proof, although any $r>2$ will suffice. We  analyze the integral in two parts, one for the interval $(0,k)$ and one for $[k,\infty)$. For the first part, we have
\begin{align}
\label{y9}
\int_0^k\mathbb{P}(S^{(n)}_k>x) \dd x &=\int_0^k \mathbb{P}({\textstyle \sum_{i=1}^\infty}\hat{A}^{(n)}_i > x+k\beta)\dd x\, \leq\, \int_0^k \mathbb{P}({\textstyle \sum_{i=1}^\infty}\hat{A}^{(n)}_i > k\beta)\dd x  \nonumber\\
&= k\,\mathbb{P}({\textstyle \sum_{i=1}^k }\hat{A}^{(n)}_i >  k\beta) \,\leq\, \frac{C_3}{k^2\beta^6} + k^2\mathbb{P}(\hat{A}^{(n)}> \tfrac{1}{3}k),
\end{align}
where we used \eqref{y8} in the last inequality. 
Hence,
\begin{align}
\label{y10}
\sum_{k=1}^\infty\frac{1}{k}\, \int_0^k \mathbb{P}(S^{(n)}_k>x)\dd x &\leq \, \frac{C_3}{\beta^6}\sum_{k=1}^\infty k^{-3} +\sum_{k=1}^\infty k\,\mathbb{P}(\hat{A}^{(n)}>\tfrac{1}{3}k) \nonumber \\
&\leq C_1^*+\sum_{k=1}^\infty k\,\mathbb{P}(\hat{A}^{(n)}>\tfrac{1}{3}k).
\end{align}
With the help of the inequality (see \cite{Sigman2011b}),
\begin{equation}
\label{y11}
(b-a)a\,\mathbb{P}(X>b) \leq \int_a^b x\,\mathbb{P}(X>x) \dd x \qquad \text{\rm for } 0<a<b,
\end{equation}
we get by taking $a=(k-1)/3$ and $b=k/3$,
\begin{align}
\label{y12}
k\,\mathbb{P}(\hat{A}^{(n)}>\tfrac{1}{3}k) &\leq \frac{9\,k}{k-1}\int_{(k-1)/3}^{k/3} x\,\mathbb{P}(\hat{A}^{(n)}>x) \dd x \nonumber \\
&\leq 18\int_{(k-1)/3}^{k/3} x\,\mathbb{P}(\hat{A}^{(n)}>x) \dd x,
\end{align}
for $k\geq 2$. Since the tail probability for $k=1$ is obviously bounded by 1, this yields
\begin{align}
\label{y13}
\sum_{k=1}^\infty k\,\mathbb{P}(\hat{A}^{(n)}>\tfrac{1}{3}k) &\leq 1+18\sum_{k=2}^\infty\int_{(k-1)/3}^{k/3} x\,\mathbb{P}(\hat{A}^{(n)}>x) \dd x\nonumber\\
&\leq 1+ 18\int_{0}^{\infty} x\,\mathbb{P}(\hat{A}^{(n)}>x)\dd x  \leq 1+18\,\E[\hat{A}^{(n)2}] < \infty,
\end{align}
since $\hat{A}^{(n)}$ has finite variance by assumption. This completes the integral over the first interval. For the second part, we use \eqref{y8} again to find
\begin{align}
\label{y14}
\int_k^\infty \mathbb{P}(S^{(n)}_k>x)dx &=\int_k^\infty \mathbb{P}({\textstyle \sum_{i=1}^\infty}\hat{A}^{(n)} > x+k\beta)dx \leq \int_k^\infty \mathbb{P}({\textstyle \sum_{i=1}^\infty}\hat{A}^{(n)} > x)\dd x\nonumber \\
&\leq C_3\int_k^\infty \frac{k^2}{x^6} \dd x + k\int_k^\infty \mathbb{P}(\hat{A}^{(n)} >\tfrac{1}{3}x)\dd x\nonumber \\
&=  \frac{5 C_3}{k^3}+ k\int_k^\infty \mathbb{P}(\hat{A}^{(n)}>\tfrac{1}{3}x) \dd x.
\end{align}
So,
\begin{equation}
\label{y15}
\sum_{k=1}^\infty \frac{1}{k} \int_k^\infty \mathbb{P}(S^{(n)}_k>x)\dd x \leq C_2^* + \sum_{k=1}^\infty \int_k^\infty \mathbb{P}(\hat{A}^{(n)}_i>\tfrac{1}{3}x) \dd x,
\end{equation}
for some constant $C_2^*$. Last, we are able to upper bound the second term in \eqref{y15} by Tonelli's theorem:
\begin{align}
\sum_{k=1}^\infty \int_k^\infty \mathbb{P}(\hat{A}^{(n)}_i>\tfrac{1}{3}x) dx &\leq \int_1^\infty x\mathbb{P}(\hat{A}^{(n)}>\tfrac{1}{3}x) \dd x \nonumber\\
\label{y16}
&\leq 9\int_0^\infty y\mathbb{P}(\hat{A}^{(n)}>y) dy = 9\E[\hat{A}^{(n)2}] < \infty.
\end{align}
Combining the results in \eqref{y10}, \eqref{y13}, \eqref{y15} and \eqref{y16}, we find
\begin{equation*}
\label{y17}
\sum_{k=1}^\infty \frac{1}{k} \int_0^\infty \mathbb{P}(S^{(n)}_k>x)\dd x < \infty,
\end{equation*}
and thus
\begin{align*}
\lim_{n\rightarrow\infty} \E[\hat{Q}^{(n)}] &= \lim_{n\rightarrow\infty}  \sum_{k=1}^\infty \frac{1}{k}\int_0^\infty \mathbb{P}(S^{(n)}_k > x)\dd x \nonumber\\
&= \sum_{k=1}^\infty \frac{1}{k} \int_0^\infty\mathbb{P}({\textstyle\sum_{i=1}^k} Z_i > x)\dd x = \E [M_\beta].
\end{align*}
Finally, we show how the proof changes for the convergence of $\Var \hat{Q}^{(n)}$. The expressions for $\E [\hat{Q}^{(n)}]$ and $\Var \hat{Q}^{(n)}$ in \eqref{y1} and \eqref{y2} only differ in the term $\sqrt{x}$. Hence only minor modifications are needed to also prove convergence of the variance. Note that boundedness of the integral over the interval $(0,k)$ in \eqref{y9}-\eqref{y13} remains to hold when substituting $\sqrt{x}$ for $x$. \eqref{y14} changes into
\begin{align*}
\label{y18}
\int_k^\infty \mathbb{P}(S^{(n)}_k>\sqrt{x})dx &=\int_k^\infty \mathbb{P}({\textstyle \sum_{i=1}^\infty}\hat{A}^{(n)}_i > \sqrt{x}+k\beta)\dd x \nonumber \\
&\leq C_3\int_k^\infty \frac{k^2}{(\sqrt{x}+k\beta)^6} dx + k\,\int_k^\infty \mathbb{P}(\hat{A}^{(n)}>\tfrac{1}{3}\sqrt{x}) \dd x \nonumber\\
&\leq \frac{C_4^*}{k}+ k\,\int_k^\infty \mathbb{P}(\hat{A}^{(n)}>\tfrac{1}{3}\sqrt{x}) \dd x,
\end{align*}
for some constant $C_4^*$, so that
\begin{equation*}
\sum_{k=1}^\infty \frac{1}{k} \int_k^\infty \mathbb{P}(S^{(n)}_k>\sqrt{x})\dd x \leq C_4^* + \sum_{k=1}^\infty \int_k^\infty \mathbb{P}(\hat{A}^{(n)}>\tfrac{1}{3}\sqrt{x}) \dd x.
\end{equation*}
Lastly, we have
\begin{align*}
\sum_{k=1}^\infty \int_k^\infty \mathbb{P}(\hat{A}^{(n)}>\tfrac{1}{3}\sqrt{x}) \dd x &\leq \int_1^\infty x\mathbb{P}(\hat{A}^{(n)}>\tfrac{1}{3}\sqrt{x}) \dd x \nonumber\\
\label{y17a}
&\leq 18\int_0^\infty y^2\mathbb{P}(\hat{A}^{(n)}>y) \dd y = 18\,\E[\hat{A}^{(n)3}] < \infty.
\end{align*}
Therefore the sum describing the variance is also uniformly convergent in $n$, so that interchanging of infinite sum and limit is permitted and
\begin{align*}
\lim_{n\rightarrow\infty} \Var\,\hat{Q}^{(n)} &= \lim_{n\rightarrow\infty}  \sum_{k=1}^\infty \frac{1}{k}\int_0^\infty \mathbb{P}(S^{(n)}_k > \sqrt{x})\dd x \nonumber \\
&= \sum_{k=1}^\infty \frac{1}{k} \int_0^\infty\mathbb{P}({\textstyle\sum_{i=1}^k} Z_i > \sqrt{x})\dd x = \Var M_\beta.
\end{align*}

\section{Numerical procedures}\label{numprocs}
An alternative characterization of the stationary distribution is based on the analysis in \cite{Boudreau1962} and considers a factorization in terms of (complex) roots:
\begin{equation*}
\label{t9}
Q^{(n)}(w) = \frac{(s_n-\E [A^{(n)}])(w-1)}{w^{s_n}-\tilde{A}^{(n)}(w)}\,\prod_{k=1}^{s_n-1} \frac{w-z^n_k}{1-z^n_k},
\end{equation*}
where $z_1^n,z_2^n...,z_{s_n-1}^n$ are the $s_n-1$ zeros of $z^{s_n}-\tilde{A}^{(n)}(z)$, in $|z|<1$, yielding
\begin{equation*}
\label{c2}
\mu_Q  = \frac{\sigma_n^2}{2(s_n-\mu_n)}-\frac{s_n-1+\mu_n}{2} + \sum_{k=1}^{s_n-1} \frac{1}{1-z^n_k},
\end{equation*}
\begin{equation*}
\label{c3}
\mathbb{P}(Q^{(n)}=0) = \frac{s_n-\mu_A}{\tilde{A}^{(n)}(0)}\prod_{k=1}^{s_n-1}\frac{z^n_k}{z^n_k-1},
\end{equation*}
which for our choice of $\tilde{A}^{(n)}(z)$ becomes
\begin{equation*}
\label{c4}
\mu_Q = \frac{a_nb_n(b_n+1)}{2\beta\sqrt{a_n}b_n}-\frac{2a_nb_n+\beta\sqrt{a_nb_n(b_n+1)}-1}{2}+\sum_{k=1}^{s_n-1} \frac{1}{1-z^n_k},
\end{equation*}
\begin{equation*}
\label{c5}
\mathbb{P}(Q^{(n)}=0) = \beta \sqrt{a_nb_n(b_n+1)}(1+b_n)^{a_n}\prod_{k=1}^{s_n-1} \frac{z^n_k}{z^n_k-1},
\end{equation*}
where $z_1,...,z_{s_n-1}$ denote the zeros of $z^{s_n} - \tilde{A}^{(n)}(z)$ in $|z|<1$. A robust numerical procedure to obtain these zeros is essential for a base of comparison. We  discuss two methods that fit these requirements. The first follows directly from \cite{Janssen2005}. \\
\begin{lemma}\label{fixedIterLemma}
Define the iteration scheme
\begin{equation}
\label{c6}
z_k^{n,l+1} = w^n_k [\tilde{A}^{(n)}(z_k^{n,l})]^{1/s_n},
\end{equation}
with $w^n_k = \ee^{2\pi ik/s_n}$ and $z_k^{n,0}=0$ for $k=0,1,\ldots,s_{n-1}$. Then $z_k^{n,l} \rightarrow z_k^n$ for all $k=0,1,...,s_n-1$ for $l\rightarrow \infty$.
\end{lemma}

\begin{proof}
The successive substitution scheme given in \eqref{c6} is the fixed point iteration scheme described in \cite{Janssen2005}, applied to the pgf of our arrival process. The authors show that, under the assumption of $\tilde{A}^{(n)}(z)$ being zero-free within $|z|\leq 1$, the zeros can be approximated arbitrarily closely, given that the function $[\tilde{A}^{(n)}(z)]^{1/s_n}$ is a contraction for $|z|\leq 1$, i.e.
\begin{equation*}
\label{c7}
\Bigl|\frac{\dd}{\dd z}[\tilde{A}^{(n)}(z)]^{1/s_n}\Bigr| < 1.
\end{equation*}
In our case,
\begin{align}
\label{c8}
\Bigl|\frac{\dd}{\dd z}[\tilde{A}^{(n)}(z)]^{1/s_n}\Bigr| = \Bigl|\frac{\dd}{\dd z}\left(1+(1-z)b_n\right)^{-a_n/s_n}\Bigr| = \frac{a_nb_n}{s_n}\Bigl|1+(1-z)b_n\Bigr|^{-a_n/s_n-1},
\end{align}
where $a_nb_n/s_n = \rho_n$ is close to, but less than 1 and
\begin{align*}
\label{c9}
|1+(1-z)b_n| \geq |1+b_n|-|z|b_n = 1+(1-|z|)b_n \geq 1,
\end{align*}
when $|z|\leq 1$. Hence the expression in \eqref{c8} is less than 1 for all $|z|\leq 1$. Evidently, $\tilde{A}^{(n)}(z)$ is also zero-free in $|z|\leq 1$. Thus \cite[Lemma~3.8]{Janssen2005} shows that $z_k^{n,l}$ as in \eqref{c6} converges to the desired roots $z^n_k$ for all $k$ as $l$ tends to infinity.
\end{proof}

\begin{remark}
The asymptotic convergence rate of the iteration in \eqref{c6} equals \\
\noindent $\frac{\dd}{\dd z}[\tilde{A}^{(n)}(z)]^{1/s_n}$ evaluated at $z=z_k^n$. Hence, convergence is slow for zeros near 1 and fast for zeros away from 1.
\end{remark}

A different approach is based on the B\"urmann-Lagrange inversion formula.
\begin{lemma}\label{BLLemma}
Let $w^n_k = e^{2\pi ik/s_n}$ and $\alpha_n = a_n/s_n$. Then the zeros of $z^{s_n}-\tilde{A}^{(n)}(z)$ are given by
\begin{equation*}
z_k^n = \sum_{l=1}^\infty \frac{1}{l!}\,\frac{\beta[l\alpha_n+l-1)}{\beta(l\alpha_n)}\,\frac{b_n+1}{b_n}\Bigl(\frac{b_n}{(b_n+1)^{\alpha_n+1}}\Bigr)^l (w_k^n)^l,
\end{equation*}
for $k=0,1,...,s_n-1$.
\end{lemma}

\begin{proof}
Note that we are looking for $z$'s that solve
\begin{equation*}
\label{c10}
z\,[\tilde{A}^{(n)}(z)]^{-1/s_n} = z\left(1+(1-z)b_n\right)^{a_n/s_n} = w,
\end{equation*}
where $w = w_k = \ee^{2\pi i k/s_n}$. The B\"urmann-Lagrange formula for $z=z(w)$, as can be found in \cite[Sec.~2.2]{debruijn} for $z=z(w)$ is given by
\begin{align*}
z(w) &= \sum_{l=1}^\infty \frac{1}{l!}\,\left(\frac{\dd}{\dd z}\right)^{l-1}\left[\left(\frac{z}{z(1+(1-z)b_n)^{a_n/s_n}}\right)^l\right]_{z=0}\,w^l\nonumber\\
\label{c11}
&= \sum_{l=1}^\infty \frac{1}{l!}\,\left(\frac{\dd}{\dd z}\right)^{l-1}\left[\left(1+(1-z)b_n)^{-l\,a_n/s_n}\right)\right]_{z=0}\,w^l.
\end{align*}
Set $\alpha_n = a_n/s_n$. We compute
\begin{equation*}
\label{c1}
\left(\frac{\dd}{\dd z}\right)^{l-1}\left[ (1+(1-z)b_n)^{-l\alpha_n}\right]_{z=0} = \frac{\beta(l\alpha_n+l-1)}{\beta(l\alpha_n)}\,\frac{1+b_n}{b_n}\,\left(\frac{b_n}{(1+b_n)^{\alpha_n+1}}\right)^l.
\end{equation*}
With $c_n = b_n/(1+b_n)^{\alpha_n+1}$ and $d_n = (1+b_n)/b_n$, we thus have
\begin{equation*}
\label{c13}
z(w) = d_n\,\sum_{l=1}^\infty \frac{\beta(l\alpha_n+l-1)}{\beta(l+1)\beta(l\alpha_n)} c_n^l\,w^l.
\end{equation*}
By Stirling's formula
\begin{equation*}\label{c14}
\frac{\beta(l\alpha_n+l-1)}{\beta(l+1)\beta(l\alpha_n)}  = \frac{D}{l\sqrt{l}}\left(\frac{(\alpha_n+1)^{\alpha_n+1}}{\alpha_n^{\alpha_n}}\right)^l,
\end{equation*}
where $D=\alpha_n^{1/2}(\alpha_n+1)^{-3/2}(2\pi)^{-1/2}$. Now,
\begin{equation*}
\label{c15}
\frac{(\alpha_n+1)^{\alpha_n+1}}{\alpha_n^{\alpha_n}}\, c_n = \frac{(\alpha_n+1)^{\alpha_n+1}}{\alpha_n^{\alpha_n}}\cdot \frac{b_n}{(1+b_n)^{\alpha_n+1}} = \left(\frac{b_n+\rho_n}{b_n+1}\right)^{\rho_n/b_n + 1}\left(\frac{1}{\rho_n}\right)^{\rho_n/b_n}.
\end{equation*}
This determines the radius of convergence $r_{\rm BL}$ of the above series for $z(w)$:
\begin{equation}
\label{c16}
\frac{1}{r_{\rm BL}} := \left(\frac{b_n+\rho_n}{b_n+1}\right)^{\rho_n/b_n + 1}\left(\frac{1}{\rho_n}\right)^{\rho_n/b_n}.
\end{equation}
The derivative with respect to $\rho_n$ of the quantity
\begin{equation}
\label{c17}
\left(1+\frac{\rho_n}{b_n}\right) {\rm ln }\left(\frac{b_n+\rho_n}{b_n+1}\right)+\frac{\rho_n}{b_n}\,{\rm ln}\left(\frac{1}{\rho_n}\right)
\end{equation}
is given by
\begin{equation*}
\label{c18}
\frac{1}{b_n}{\rm ln }\Bigl(\frac{b_n+\rho_n}{b_n\rho_n+\rho_n}\Bigr) > 0
\end{equation*}
for $0<\rho_n<1$ and $b_n>0$. Furthermore, the quantity in \eqref{c17} vanishes at $\rho_n=1$ and is therefore negative for $0<\rho_n<1$ and $b_n>0$.
\begin{remark}
The formula for the radius of convergence in \eqref{c16} clearly shows the decremental effect of both having a large $b_n$ and of having $\rho_n$ close to 1. The quantities $\beta(l\alpha+l-1)/(\beta(l+1)\beta(l\alpha))$ in the power series for $z(w)$ are not very convenient for recursive computation, although normally $\alpha_n = a_n/s_n$ is a rational number.\end{remark}
\end{proof}

\section{Statistical procedures}\label{statproc}
To calibrate our model to real data, we now discuss some statistical procedures to show the presence of overdispersion and to estimate the parameters of the mixed Gamma-Poisson distribution. Let $x_1,...,x_n$ denote the observed arrival counts in consecutive time slots. These observations can be interpreted as realizations of the random variables $A_1,...,A_N$, and
 \begin{equation*}
\bar{a}_N=\frac{1}{N}\sum_{i=1}^N x_i, \qquad \bar{s}_N^2 = \frac{1}{N-1}\sum_{i=1}(x_i-\bar{x}_i)^2,
\end{equation*}
the sample mean and variance with equivalent random variables $\bar{A}_N$ and $S_N^2$, respectively. Several tests with null hypothesis that $x_1,...,x_N$ originate from a (constant rate) Poisson distribution were discussed by \cite{Brown2002}. We  mention two of them. The first is frequently referred to as the \emph{dispersion test}, and is based on the test statistic
\begin{equation*}
\label{dispTest}
D_N = \frac{(N-1)S_N^2}{\bar{A}_N},
\end{equation*}
which is approximately chi-squared distributed with $N-1$ degrees of freedom. When using a significance level $\alpha$, the critical value is equal to the $(1-\alpha)$-th quantile of chi-squared distribution $\chi^2_{N-1,1-\alpha}$. The second test, which is also used in \cite{koolejongbloed}, involves the test statistic
\begin{equation*}
\label{NStest}
T_N = \sqrt{N/2}\,\Bigl(\frac{S_N^2}{\bar{A}_N}-1\Bigr),
\end{equation*}
which is known as the Neyman-Scott test statistic. Under the null hypothesis $T_N$ tends to a standard normal random variable for large $N$. Hence both test statistics rely on the ratio of the sample variance and sample mean, which should be 1 if $A_1,...,A_N$ are indeed i.i.d. Poisson distributed. Excessive values of $D_N$ and $T_N$ therefore raise the suspicion of overdispersed arrivals.

Once either (or both) of the Poisson tests rejects the hypothesis of arrivals originating from a unicomponent Poisson process, we fit the data to the Gamma-Poisson mixture. Note that if we assume $A_i$ to be distributed as a Poisson random variable with random rate $\Lambda_i$, which is in turn Gamma distributed with parameters $a$ and $1/b$, then $A_i$ is in fact a negative binomial random variable with parameters $r = a$ and $p=b/(b+1)$. Finding estimators $\hat{a}$ and $\hat{b}$ therefore is equivalent to fitting a negative binomial distribution to the data to obtain $\hat{r}$ and $\hat{p}$, followed by retrieving $\hat{a} = \hat{r}$ and $\hat{b} = \hat{p}/(1-\hat{p})$. We proceed by applying the maximum likelihood estimation method described in \cite{koolejongbloed} to find $\hat{r}$ and $\hat{p}$. This method prescribes to set $\hat{r}$ to be the value of $r$ for which the \emph{profile loglikelihood function} defined by
\begin{equation*}
L(r) = \frac{1}{N}\,\sum_{i=1}^N\sum_{j=1}^{a_i} {\rm ln}(r+j+1)+r\,{\rm ln}\,r -(r+\bar{a}_N)\,{\rm ln}(r+\bar{a}_N),
\end{equation*}
is attained. Subsequently, $\hat{p} = \hat{r}/(\hat{r}+\bar{a}_N)$, so that $\hat{a} = \hat{r}$ and $\hat{b} = \hat{r}/\bar{a}_N$.

Finally, given the estimators $\hat{a}$ and $\hat{b}$, we need statistical evidence that the obtained Poisson mixture indeed fits the data reasonably well. Here we again cite \cite{koolejongbloed}, who give a method to retrieve the $p$-value for the goodness-of-fit based on bootstrap and Monte-Carlo simulation. In our experiments, we work with $10^6$ replications of the Monte-Carlo simulation to obtain the approximated $p$-value. We refer to the appendix of \cite{koolejongbloed} for further details on this method.

\resettocdepth 

\end{subappendices}

\chapter{Retrial queues in the QED regime}

\begin{chapterstart}
Large-scale queueing systems with retrying customers are intrinsically hard to evaluate analytically. 
We in this chapter explore and extend the asymptotic approximation technique proposed by Avram et al.~\cite{Avram2013}, that is able to characterize the impact of slow retrials in the QED regime, in three queueing models. 
The technique evolves around a fixed-point equation that quantifies the increased inflow due to retrials implicitly. 
We translate this fixed-point method into a powerful and elegant dimensioning procedure that is able to deal with both stationary and time-varying demand.
\end{chapterstart}
\begin{flushright}
Based on \\
\textbf{Delayed workload shifting in many-server systems}\\
\textit{Johan van Leeuwaarden, Britt Mathijsen \& Fiona Sloothaak}\\
In \textit{SIGMETRICS Performance Evaluation Review, 43(2), 10--12 (2015)}\\
and \\
\textbf{Cloud provisioning in the QED regime}\\
\textit{Johan van Leeuwaarden, Britt Mathijsen \& Fiona Sloothaak}\\
In \textit{Proceedings of the 9th EAI International Conference on Performance Evaluation Methodologies and Tools, 180--187 (2016)}
\end{flushright}

\newpage

\section{Introduction}
%

\textbf{Retrial queues.}
In the previous chapters, we analyzed queueing systems in which all arriving customers join the queue and stay until eventually completing service with one of the servers. 
From a practical perspective though, these assumptions are questionable. 
For instance, in call centers, customer impatience is known to play a crucial role in the queueing dynamics, see e.g.~\cite{Gans2003,Brown2005,Zeltyn2005}. 
Similar features are also seen in health care \cite{Dai2012,Armony2015}.
However, impatience may not be the only cause of customers leaving the system without being seen by a server. 
Physical constraints may force system managers to apply some sort of admission policy. 
The simplest example of such admission control is the busy-signal in call centers, in which arriving customers finding all servers busy are simply discarded. 
But more elaborate strategies can be considered. 
A straightforward relaxation of the busy-signal policy is to allow a finite amount of waiting space, and block customers who find a full waiting room upon arrival. 
Many other options, such as probabilistic and dynamic admission control policies may be considered, see e.g.~\cite{Janssen2013,Armony2004} and references therein.

Since customers arrive to the system for the purpose of getting assistance from one of the servers, it is reasonable to assume that these refused customers retry getting access tot the system at a later point in time. 
In fact, retrials are widely observed in telecommunication systems, see e.g.~\cite{Cohen1957,Falin1997,Mandelbaum2000,Aguir2008}, and customers typically repeat their attempt until successful. 
Naturally, retrials have a detrimental effect on the performance of the queueing system in terms of QoS, compared to the setting in which blocked customers do not return. 
Hence, one needs to account for their impact in both performance analysis and the staffing decisions. 

Unfortunately, the modeling of retrials is analytically challenging \cite{Cohen1957,Falin1997}, and numerical approaches become computationally infeasible as the number of servers increases, which is precisely the regime we are interested in. 
We therefore aim to tackle the performance analysis of such retrial systems in an asymptotic manner.
We do so through a clever technique that was recently documented by Avram et al.~\cite{Avram2013}.\\
\\*
\textbf{Fixed-point equation.}
In this chapter, we will show how the asymptotic approximation technique of \cite{Avram2013} can be extended to more complex large-scale retrial queueing systems in the QED regime.
In \cite{Avram2013}, the authors study the $M/M/s/s$ queue with slow retrials. 
That is to say, customers retry only after a (stochastic) delay period that is relatively long compared to the service time. 
Under this assumption, the authors combine QED limits with a fixed-point equation, which characterizes the impact of retrials implicitly. 
We summarize and reformulate their main ideas here for completeness. 

Consider the standard $M/M/s/s$ queue with arrival rate $\l$ and service rate $\mu=1$, so that the offered load is $R=\l$.
Customers finding upon arrival all $s$ servers busy retry after a stochastic delay with mean $1/\delta$. 
The first ingredient of the method is Cohen's equation. 
This result, first reported by Cohen \cite{Cohen1957}, says that the stationary distribution of a $M/M/s/s$ queue with retrials converges as $1/\delta \to\infty$ to that of a $M/M/s/s$ queue with increased arrival rate $R+\Omega$, where $\Omega$ is the unique positive solution to 
\begin{equation}
\Omega = \left( R+\Omega \right)\, B(R+\Omega,s),
\label{eq:cohens_equation}
\end{equation}
where $B(R,s)$ denotes the blocking probability in the $M/M/s/s$ queue with offered load $R$, i.e.~the Erlang-B formula:
\begin{equation}
B(R,s) := \frac{ R^s / s! }{\sum_{k=0}^s R^k/k!} = \frac{\P({\rm Pois}(R) = s)}{\P({\rm Pois}(R) \leq s)}.
\end{equation}
Equation \eqref{eq:cohens_equation} essentially equates the arrival volume generated by retrials, given by definition on the left-hand side, to the right-hand side which quantifies this volume as a fraction of customers blocked times the increased arrival volume. 
Indeed, the retrial stream can for long retrial times be considered as independent from the primary arrival stream, yielding a thinned Poisson process \cite{Cohen1957}. 

The second crucial observation is that under QED scaling, i.e. $s = R+\beta\sqrt{R} + o(\sqrt{R})$, we have
\begin{equation}
\sqrt{R}\cdot B( R, R+\beta\sqrt{R} ) \to \frac{\f(\b)}{\F(\b)} =: f_0(\beta),
\label{eq:erlangb_limit}
\end{equation} 
as $R\to\infty$ for all $\b\in\mathbb{R}$, see e.g.~Lemma 1 of \cite{Avram2013} or the proof of Proposition 1.1 in this thesis.
Hence, we heuristically deduce that $\Omega = \alpha\sqrt{R}$ for some $\alpha > 0$ and denote $\Rt = R + \Omega$ -- a rigorous argument can be found in \cite{Avram2013}. 
Rewrite $R = \Rt-\alpha \sqrt{\Rt} + o(\sqrt{\Rt})$ and note that $\Rt = O(R)$. 
Then \eqref{eq:cohens_equation} becomes
\begin{equation}
\alpha \sqrt{\Rt} = \Rt \cdot B \big( \Rt, \Rt+(\b-\alpha) \sqrt{\Rt} \big) + o(\sqrt{\Rt}).
\label{eq:cohen_equation_filled}
\end{equation}
Dividing both sides of \eqref{eq:cohen_equation_filled} by $\sqrt{\Rt}$ and letting $\Rt \to \infty$ then together with \eqref{eq:erlangb_limit} yields the fixed-point equation 
\begin{equation}
\alpha = f_0(\beta -\alpha).
\end{equation}
It can be shown that this fixed-point equation has a unique positive root for all $\beta > 0$, which can be computed numerically. 
Recalling Cohen's retrial queue characterization, Avram et al.~\cite{Avram2013} conclude that the loss system with slow retrials can in the QED regime be characterized in terms of the original loss system without retrials, but with a corrected QoS-parameter $\b-\alpha$. \\
\\*
\textbf{Structure of the chapter.}
The fixed-point method of Avram et al. provides a quick and elegant way to approximate the behavior of large-scale loss systems that experience retrials. 
In the remainder of this chapter we will explore if and how this technique extends to three more complex queueing settings.
These three models have in common that (i) they exhibit QED limiting behavior, which can be quantified explicitly, and (ii) the blocking probability is $O(1/\sqrt{R})$, so that the retrial volume is $O(\sqrt{R})$. 
Note that these were the two essential features for the fixed-point method to work.

In Section \ref{sec:basic_model} we describe a direct extension of the $M/M/s/s$ queue, in which some amount of waiting room is present. 
That is, we analyze the $M/M/s/n$ queue with retrials, where $n>s$. 
Naturally, this requires a scaling for both $s$ and $n$ as $R\to\infty$, which will become clear in this section. 
Motivated by a process related to cloud computing, we in Section \ref{sec:cloud_model} study a tandem queueing network, in which total number of concurrent admissions is limited.
Section \ref{sec:abandonments_retrials} analyzes a queueing model in which all customers are admitted upon arrival, but make the deliberate decision to abandon the queue and retry later in case their patience runs out.
In Section \ref{sec:retrial_dimensioning} we show how the fixed-point method together with QED scaling can be used for dimensioning purposes in both stationary and time-varying environments. 
We end the chapter in Section \ref{sec:retrial_conclusion} with some final remarks and suggestions for future research.

\section{The $M/M/s/n$ queue}
\label{sec:basic_model}
In this section, we discuss a simple extension of the loss model of \cite{Avram2013}, namely the $M/M/s/n$ queue with retrials with $n>s$, to expose typical behavior of retrial queues and the influence of the retrial rate $\delta$. 
Second, we illustrate the fixed-point method for this model and perform numerical experiments to verify its accuracy.

\subsection{Markov process}
\begin{figure}[h]
\centering
\begin{tikzpicture}[scale=1]

\draw[dashed] (0,0.25) rectangle (4.75,2.75) node[right,above] {$n$};
\draw[thick] (0.75,1) -- (2.5,1) -- (2.5,2) -- (0.75,2);
\draw[thick] (1.5,1.1) -- (1.5,1.9);
\draw[thick] (1.75,1.1) -- (1.75,1.9);
\draw[thick] (2,1.1) -- (2,1.9);
\draw[thick] (2.25,1.1) -- (2.25,1.9);
\draw[thick] (3.5,1.5) circle (0.5) node {$s$};

\draw[thick,->] (-2.25,1.5)  -- (1,1.5);
\draw[thick,->] (4.25,1.5) -- (5.5,1.5);

\node at (-2.25,1.8) {Pois$(\lambda)$};
\node at (3.5,2.35) {$\exp(\mu)$};

\draw[thick] (-0.5,1.5) -- (-0.5,1);
\draw[thick] (-0.5,1) to [in=0,out=270] (-1,0.5);
\draw[thick] (-1,0.5) to [in=270,out=180] (-1.5,1);
\draw[thick,->] (-1.5,1) to [in=180,out=90] (-1,1.45);

\node at (-2,0.4) {$\exp(\delta)$};

\end{tikzpicture}
\caption{An $M/M/s$ queue with space constraints and retrials.}
\label{fig:BasicModel}
\end{figure}
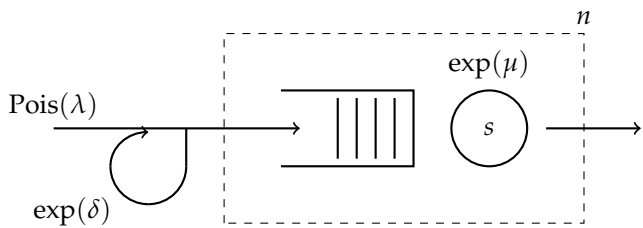

We consider the standard $M/M/s/n$ queue with arrival rate $\lambda$ and service rate $\mu$.
Without loss of generality, we set $\mu=1$ throughout this chapter, so that offered load $R$ equals $\lambda$. 
A customer that finds upon arrival a free server occupies this server immediately, while customers that meet more than $s$ but fewer than $n>s$ customers in the system are admitted and wait in a queue for a server to become available.
Customers who meet upon arrival $n$ customers  are not admitted directly, but will retry after an exponentially distributed time with mean $1/\delta$. 
Each initially blocked customer performs retrials until admitted eventually; see Figure \ref{fig:BasicModel}. 

\begin{figure}
\centering
\begin{tikzpicture}[scale = 0.75]

\draw[->] (-0.5,0) -- (11,0); 
	\node[right] at (11,0)  {$N(t)$};
\draw[->] (0,-0.5) -- (0,8.5);
	\node[above] at (0,8.5) {$Q(t)$};
\draw[dashed] (0,8) -- (11,8);
\draw[dashed] (0,5) -- (11,5);

\foreach \x in {0,1,...,10}
	\foreach \y in {0,1,...,8}
		\draw[color=gray,fill=gray] (\x,\y) circle (0.05cm);

\draw[fill] (0,0) circle (0.1cm);
\draw[fill] (0,5) circle (0.1cm);
\draw[fill] (0,8) circle (0.1cm);

\node[right] at (12,8) {$n$};
\node[right] at (12,5) {$s$};
\node[left] at (0,3) {$j$};
\node[below] at (7,0) {$i$};

\draw[thick,->] (0,0) to[bend left] coordinate (arc1) (0,1);
	\node[left] at (arc1) {$\lambda$};
\draw[thick,->] (0,5) to[bend left] coordinate (arc2) (0,6);
	\node[left] at (arc2) {$\lambda$};
\draw[thick,->] (0,6) to[bend left] coordinate (arc2) (0,5);
	\node[right] at (arc2) {$s$};
\draw[thick,->] (0,5) to[bend left] coordinate (arc2) (0,4);
	\node[right] at (arc2) {$s$};
\draw[thick,->] (0,7) to[bend left] coordinate (arc3) (0,8) ;
	\node[left] at (arc3) {$\lambda$};
\draw[thick,->] (0,8) to[bend left] coordinate (arc4) (1,8) ;
	\node[above] at (arc4) {$\lambda$};

\draw[fill] (7,3) circle (0.1cm);

\draw[thick, ->] (7,3) to[bend right] coordinate (arc5) (7,4);
	\node[left] at (arc5) {$\lambda$};
\draw[thick, ->] (7,3) to[bend left] coordinate (arc5) (7,2);
	\node[right] at (arc5) {$j$};	
\draw[thick, ->] (7,3) to[bend left] coordinate (arc5) (6,4);
	\node[left] at (arc5) {$i\delta$};
		
\draw[fill] (7,6) circle (0.1cm);
	
\draw[thick, ->] (7,6) to[bend right] coordinate (arc5) (7,7);
	\node[left] at (arc5) {$\lambda$};
\draw[thick, ->] (7,6) to[bend left] coordinate (arc5) (7,5);
	\node[right] at (arc5) {$s$};	
\draw[thick, ->] (7,6) to[bend left] coordinate (arc5) (6,7);
	\node[left] at (arc5) {$i\delta$};	

\draw[fill] (3,8) circle (0.1cm);

\draw[thick, ->] (3,8) to[bend left] coordinate (arc5) (4,8);
	\node[above] at (arc5) {$\lambda$};
\draw[thick, ->] (3,8) to[bend left] coordinate (arc5) (3,7);
	\node[right] at (arc5) {$s$};	

\end{tikzpicture}
\caption{Transition diagram of the Markov process $(Q(t),N(t))$.}
\label{fig:transition_diagram}
\end{figure}
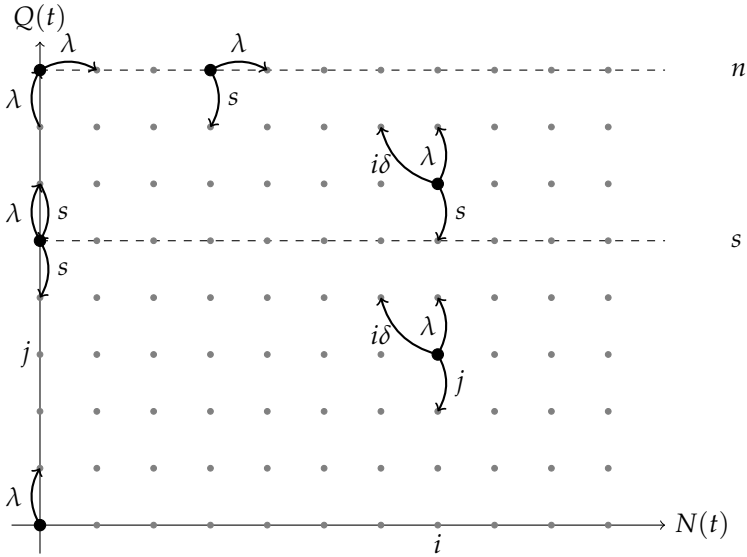

\noindent
\textbf{Quasi-birth-death process.} The system state can be described by a two-dimensional process $\{ (Q(t),N(t))\}_{t\geq 0}$ with $Q(t)$ the number of customers inside the system (either being served or waiting), and $N(t)$ the number of customers in the retrial orbit. 
Under the above assumptions, this process is a continuous-time Markov chain on the semi-infinite strip $\{0,1,\ldots,n\}\times \{0,1,\ldots\}$. 
Its transition diagram is presented in Figure \ref{fig:transition_diagram}.
From this diagram it is evident that the process is a quasi-birth-death (QBD) process.
Under stability condition $R>s$, the QBD structure of the process allows for numerical computation of the stationary distribution $\pi(i,j)$, where
\begin{equation*}
\pi(i,j) = \lim_{t\to\infty} \P\left( Q(t) = i, N(t) = j\right).
\end{equation*}
The stationary probability that an arriving customer has to wait or is blocked is given by, respectively,
\begin{align}
\P_r({\rm delay}) = \sum_{i=s}^n \sum_{j=0}^\infty \pi(i,j),\qquad 
\P_r({\rm block}) = \sum_{j=0}^\infty \pi(n,j).
\label{eq:performance_measures}
\end{align} 
Here, the subscript $r$ is meant to indicate that we consider the system with retrials.\\

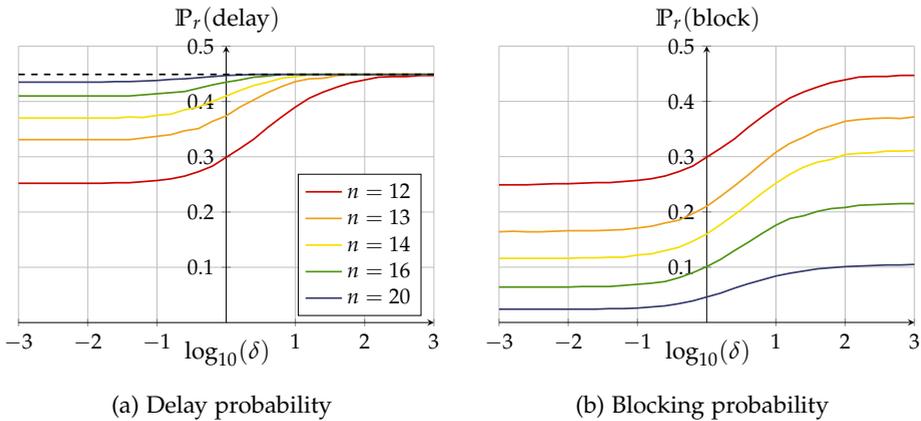
\begin{figure}
\centering
\begin{subfigure}{0.48\textwidth}
\centering
\begin{tikzpicture}[scale=0.8]
\begin{axis}[
	xmin = -3, 
	xmax = 3,
	ymin = 0,
	ymax = 0.5,
	ytick = {0,0.1,0.2,0.3,0.4,0.5},
	grid = both,
	axis line style={->},
	axis lines = middle,
	xscale=1,
	yscale=0.8,
	every axis x label/.style={at={(current axis.right of origin)},anchor=west},
	every axis y label/.style={at={(current axis.north west)},above=2mm},
	legend pos = south east]

\addplot[col1,thick] table[x=log_delta,y=n12] {Chapter_4/tikz/pdelay_small.txt};
\addplot[col2,thick] table[x=log_delta,y=n13] {Chapter_4/tikz/pdelay_small.txt};
\addplot[col3,thick] table[x=log_delta,y=n14] {Chapter_4/tikz/pdelay_small.txt};
\addplot[col4,thick] table[x=log_delta,y=n16] {Chapter_4/tikz/pdelay_small.txt};
\addplot[col5,thick] table[x=log_delta,y=n20] {Chapter_4/tikz/pdelay_small.txt};
\addplot[dashed,thick] coordinates { (-3,0.449) (3,0.449) };

\legend{{$n=12$},{$n=13$},{$n=14$},{$n=16$},{$n=20$}};

\end{axis}

\node at (3.45,-0.5) {\small $\log_{10}(\delta)$};
\node at (3.45,5) {\small $\P_r({\rm delay})$};

\end{tikzpicture}
\caption{Delay probability}
\end{subfigure}
\begin{subfigure}{0.48\textwidth}
\centering
\begin{tikzpicture}[scale=0.8]
\begin{axis}[
	xmin = -3, 
	xmax = 3,
	ymin = 0,
	ymax = 0.5,
	ytick = {0,0.1,0.2,0.3,0.4,0.5},
	grid = both,
	axis line style={->},
	axis lines = middle,
	xscale=1,
	yscale=0.8,
	every axis x label/.style={at={(current axis.north east)},anchor=west},
	every axis y label/.style={at={(current axis.north west)},above=2mm},
	legend pos = south east]

\addplot[col1,thick] table[x=log_delta,y=n12] {Chapter_4/tikz/pblock_small.txt};
\addplot[col2,thick] table[x=log_delta,y=n13] {Chapter_4/tikz/pblock_small.txt};
\addplot[col3,thick] table[x=log_delta,y=n14] {Chapter_4/tikz/pblock_small.txt};
\addplot[col4,thick] table[x=log_delta,y=n16] {Chapter_4/tikz/pblock_small.txt};
\addplot[col5,thick] table[x=log_delta,y=n20] {Chapter_4/tikz/pblock_small.txt};


\end{axis}

\node at (3.45,-0.5) {\small $\log_{10}(\delta)$};
\node at (3.45,5) {\small $\P_r({\rm block})$};
\end{tikzpicture}
\caption{Blocking probability}
\end{subfigure}
\caption{Performance metrics of the basic model with $R=10$ and $s=12$ as a function of $\log_{10}(\delta)$ for several $n$.}
\label{fig:influence_of_delta}
\end{figure}

\noindent
\textbf{Influence of retrial rate $\delta$.}
We first compute the stationary distribution of the Markov process numerically, in order to understand the influence of retrials on the queue performance. 
In particular, we investigate the effect of varying $\delta$ on the delay and blocking probability as defined in \eqref{eq:performance_measures}. 
In Figure \ref{fig:influence_of_delta} we fix $R=10$ and $s=12$, and plot the delay and blocking probability as a function of $\log_{10}(\delta)$ for several values of $n$. 
We see that the value of $\delta$ indeed does influence the performance of the queue, and its effect is particularly pronounced in systems with $n$ close to $s$. 
Both the delay probability and the blocking probability increase with $\delta$.
This can be explained as follows.
If a customer finds $n$ customers on arrival (or retrial) and hence gets blocked, she is more likely to find a less congested system in case she retries after a relatively long amount of time than a short retrial time, because the system might not yet have had enough time to recover from the congested period. 
Slow retrials hence create an opportunity to smooth out workload over time, resulting in better quality-of-service. 
Figure \ref{fig:influence_of_delta} also suggests that performance no longer changes if $\delta$ is decreased below $10^{-1}$ or increased beyond $10^2$.

Also, we note that the delay probability increases with $n$, and the blocking probability decreases with $n$, regardless of the value of $\delta$. 
Fewer customers get blocked if the waiting room $(n-s)$ increases. On the other hand, this allows more customers to enter the system, creating higher congestion levels.

Finally, notice that the delay probability approaches a constant as $\delta\to\infty$. 
In fact, this constant equals the delay probability in the standard $M/M/s$ queue, see Equation (1.2), which under these parameter settings equals 0.449 and is represented by the dashed horizontal line.
Indeed, when $\delta\to\infty$ blocked customers retry getting access to the system instantaneously and effectively create a queue (in random order) outside the system, which immediately fills up vacant spaces after service completions. 
Therefore, the $M/M/s/n$ queue with instant retrials essentially resembles the behavior of the $M/M/s$ queue. 
By similar reasoning, the blocking probability in the $M/M/s/n$ queue approaches as $\delta\to\infty$ the probability that the number of customers in the $M/M/s$ queue exceeds $n$. 

Figure \ref{fig:influence_of_delta} shows that the influence of retrials on congestion can be significant. 
For fast retrials, we are able to characterize the performance metrics through the standard multi-server queue. 
However, for slow retrials, say $\delta < 10^{-1}$, the system behavior is not comparable to that of the open $M/M/s$ queue. 

\subsection{QED regime}
\begin{figure}
\centering
\begin{subfigure}{0.48\textwidth}
\centering
\begin{tikzpicture}[scale=0.72]
\begin{axis}[
	xmin = 0, 
	xmax = 50,
	ymin = 0,
	ymax = 18,
	axis line style={->},
	axis lines = left,
	xscale= 1,
	yscale=0.8,
	xlabel = {$t$},
	legend pos = south east]

\addplot[col5,thick] table[x=t,y=Q] {Chapter_4/tikz/sample_path_lambda10.txt};
\addplot[col1,thick] table[x=t,y=N] {Chapter_4/tikz/sample_path_lambda10.txt};
\addplot[dashed,thick] coordinates { (0,12) (50,12) };
\addplot[dashed,thick] coordinates { (0,16) (50,16) };
\end{axis}

\node[anchor = west] at (6.8,3) {\small $s$};
\node[anchor = west] at (6.8,4) {\small $n$};
\end{tikzpicture}
\caption{$R=10$}
\end{subfigure}
\begin{subfigure}{0.48\textwidth}
\centering
\begin{tikzpicture}[scale=0.72]
\begin{axis}[
	xmin = 0, 
	xmax = 50,
	ymin = 0,
	ymax = 70,
	axis line style={->},
	axis lines = left,
	xscale= 1,
	yscale=0.8,
	xlabel = {$t$},
	legend pos = south east]

\addplot[col5,thick] table[x=t,y=Q] {Chapter_4/tikz/sample_path_lambda50.txt};
\addplot[col1,thick] table[x=t,y=N] {Chapter_4/tikz/sample_path_lambda50.txt};
\addplot[dashed,thick] coordinates { (0,62) (50,62) };
\addplot[dashed,thick] coordinates { (0,54) (50,54) };
\end{axis}

\node[anchor = west] at (6.8,3.5) {\small $s$};
\node[anchor = west] at (6.8,4) {\small $n$};
\end{tikzpicture}
\caption{$R=50$}
\end{subfigure}

\begin{subfigure}{0.48\textwidth}
\centering
\vspace{-2 mm}

\begin{tikzpicture}[scale=0.72]
\begin{axis}[
	xmin = 0, 
	xmax = 50,
	ymin = 0,
	ymax = 125,
	axis line style={->},
	axis lines = left,
	xscale= 1,
	yscale=0.8,
	xlabel = {$t$},
	legend pos = south east]

\addplot[col5,thick] table[x=t,y=Q] {Chapter_4/tikz/sample_path_lambda100.txt};
\addplot[col1,thick] table[x=t,y=N] {Chapter_4/tikz/sample_path_lambda100.txt};
\addplot[dashed,thick] coordinates { (0,105) (50,105) };
\addplot[dashed,thick] coordinates { (0,115) (50,115) };
\end{axis}

\node[anchor = west] at (6.8,3.85) {\small $s$};
\node[anchor = west] at (6.8,4.2) {\small $n$};
\end{tikzpicture}
\caption{$R=100$}
\end{subfigure}
\begin{subfigure}{0.48\textwidth}
\centering
\vspace{-2.5 mm}

\begin{tikzpicture}[scale=0.72]
\begin{axis}[
	xmin = 0, 
	xmax = 50,
	ymin = 0,
	ymax = 1200,
	axis line style={->},
	axis lines = left,
	xscale= 1,
	yscale=0.8,
	xlabel = {$t$},
	legend pos = south east]

\addplot[col5,thick] table[x=t,y=Q] {Chapter_4/tikz/sample_path_lambda1000.txt};
\addplot[col1,thick] table[x=t,y=N] {Chapter_4/tikz/sample_path_lambda1000.txt};
\addplot[dashed,thick] coordinates { (0,1017) (50,1017) };
\addplot[dashed,thick] coordinates { (0,1050) (50,1050) };

\end{axis}

\node[anchor = west] at (6.8,3.85) {\small $s$};
\node[anchor = west] at (6.8,4) {\small $n$};
\end{tikzpicture}
\caption{$R=1000$}
\end{subfigure}
\caption{Sample paths of $Q(t)$ (blue) and $N(t)$ (red) for increasing $R$ while $s$ and $n$ are scaled as in \eqref{eq:twofold_scaling_basic_model} with $\beta=0.5$ and $\gamma=1$ and retrial rate $\delta=0.1$.}
\label{fig:sample_paths_retrial}
\end{figure}
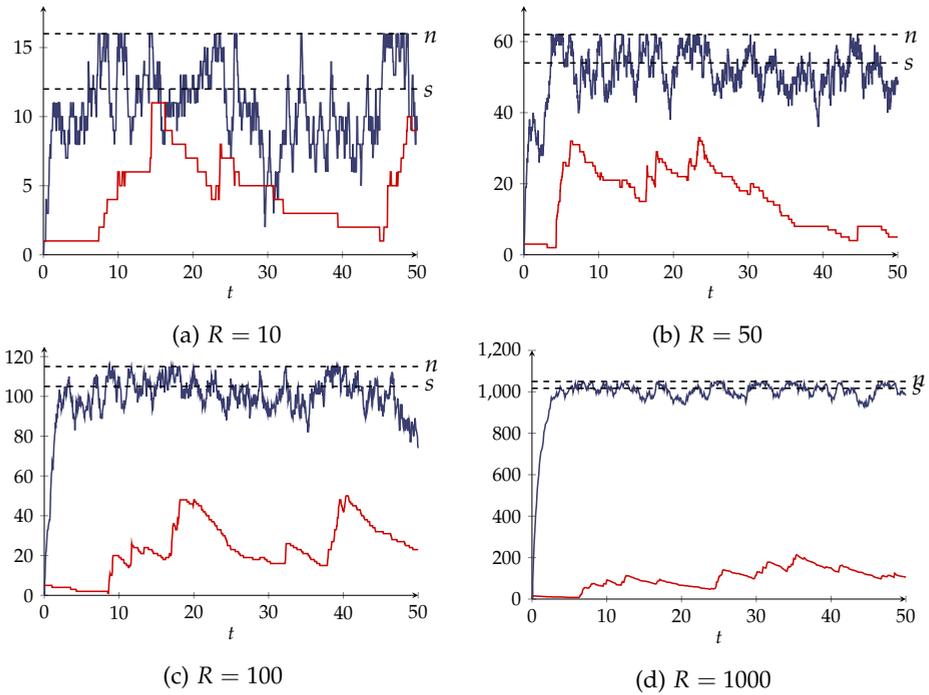

Following the approach of Avram et al.~\cite{Avram2013}, we choose to take a step back and consider the model in Figure \ref{fig:BasicModel} without the retrials first. 
When blocked customers are simply discarded, the process $\{ (Q(t),N(t))\}_{t\geq 0}$ reduces to that of the $M/M/s/n$ queue. 
In this case $N(t)=0$ and $Q(t)$ is a birth-death process with stationary distribution
\[
\pi(i) = \lim_{t\to\infty} \P(Q(t)=i) =
\left\{
\begin{array}{ll}
\pi(0)\, \frac{R^i}{i!}, & \text{if }i < s,\\
\pi(0)\, \frac{R^i}{s!s^{i-s}}, & \text{if }s\leq i \leq n,\\
\end{array}
\right.
\]
where 
\[\pi(0) = \Big( \sum_{i=0}^{s-1} \frac{R^i}{i!} + \sum_{i=s}^n \frac{R^i}{s!s^{i-s}} \Big)^{-1}.\]
Hence, 
\begin{equation}
\P({\rm delay}) = \pi(0) \sum_{i=s}^n \frac{R^i}{s!s^{i-s}}, 
\qquad 
\P({\rm block}) = \pi(0) \frac{R^n}{s!s^{n-s}}. 
\end{equation}
The $M/M/s/n$ queue is well understood.
In particular, Massey \& Wallace \cite{masseywallace} identified the asymptotic scaling regime for $s$ and $n$ under which QED-type behavior prevails.
Namely, under the two-fold scaling rule 
\begin{align}
s &= R + \beta\sqrt{R} + o(\sqrt{R}), \nonumber\\
n &= s + \gamma\sqrt{R} + o(\sqrt{R}), 
\label{eq:twofold_scaling_basic_model}
\end{align}
for $\beta\in\mathbb{R}$ and $\gamma>0$, they show that the delay probability converges to a value strictly between 0 and 1, while the blocking probability vanishes as $R\to\infty$.
Note that this is in line with our reasoning in Section 1.4.
In the next proposition, we cite the asymptotic results of \cite{masseywallace} for completeness. 
\begin{proposition}[\cite{masseywallace}]
If $s$ and $n$ scale according to \eqref{eq:twofold_scaling_basic_model}, then in the $M/M/s/n$ queue,
\begin{align}
\P({\rm delay}) &\to \frac{ 1-\ee^{-\beta\gamma}}{1-\ee^{-\beta\gamma}+\b\,\F(\b)/\f(\b)} =: g(\b,\g),\label{eq:limit_delay}
\\
\sqrt{R} \, \P({\rm block}) &\to \frac{\beta\,\ee^{-\beta\gamma}}{1-\ee^{-\b\g} + \b\,\F(\b)/\f(\b)} =: f(\b,\g),
\label{eq:limit_block}
\end{align}
as $R\to\infty$. 
\end{proposition}
Before turning to the asymptotic analysis of the model with retrials, we check empirically whether the scaling in \eqref{eq:twofold_scaling_basic_model} also achieves the desirable limiting behavior in case blocked customers are not discarded.  
In Figure \ref{fig:sample_paths_retrial} we plot sample paths $Q(t)$ and $N(t)$ in the system with retrials with $\beta = 0.5$ and $\gamma=1$ and slow retrials ($\delta = 0.1$) for increasing values of $R$. 

From these sample paths, we observe that indeed the server utilization approaches unity as $R$ tends to infinity, indicating efficient usage of resources. 
This should not be surprising, since although retrials occur, all customers eventually receive service, so that the server utilization equals $R/s = R/(R+\beta\sqrt{R}) \to 1$ as $R\to\infty$. 
Furthermore, we see that the number of customers in the system concentrates around the level $s$, implying a delay probability away from both 0 and 1.  
Observe that the order of magnitude of $N(t)$, the number of customers in the retrial orbit, is smaller than $Q(t)$ or $R$. 
This implies that as $R$ grows large, only a small fraction of customers ends up retrying. 
Naturally, the order of $N(t)$ also depends on the mean retrial time $1/\delta$.
It can be numerically verified that the expected retrial population grows linearly in $1/\delta$.
Last, observe that $N(t)$ is increasing only if $Q(t)=n$, which is visible through the surges in the sample paths of $N(t)$ in Figure \ref{fig:sample_paths_retrial}. 
This is illustrative for the dependency between the two coordinates of the process $\{(Q(t),N(t))\}_{t\geq 0}$ and therefore, we cannot expect to find a simple decoupling in the limit either. 
Instead, we propose to evaluate the model with retrials through a heuristic approach which builds upon the asymptotic behavior of the model without retrials.

\subsection{Fixed-point method}
\label{sec:fixed_point}


We continue to translate the ideas behind the fixed-point method by noting that due to \eqref{eq:limit_block}, the fraction of blocked customers is of order $1/\sqrt{R}$, which implies that the mean additional load due to retrials must be of order $\sqrt{R}$. 
We can thus assume that the total arrival rate $R_{\rm tot}$ takes the form $R_{\rm tot} = R+\alpha\sqrt{R}$ for some $\alpha>0$. 
Then, using that $R = O(R_{\rm tot})$, the first scaling rule in \eqref{eq:twofold_scaling_basic_model} is asymptotically equivalent with
\begin{equation}
 s = R_{\rm tot} + (\b-\alpha) \sqrt{R_{\rm tot}} + o(\sqrt{R_{\rm tot}}),
 \end{equation} 
while the scaling for $n$ remains unchanged. 
We thus argue that the retrial system in the QED regime mimics an $M/M/s/n$ queue with parameters $\beta_\alpha = \beta-\alpha$ and $\gamma$. 
Note that the volume of blocked users in this setting is $f(\beta-\alpha,\gamma)\sqrt{R_{\rm tot}}$. 
This quantity must equal the mean additional load $\alpha\sqrt{R} \sim \alpha\sqrt{R_{\rm tot}}$ and therefore we obtain the \textit{fixed-point equation}
\begin{equation}
\alpha = f(\beta-\alpha,\gamma).
\label{eq:fixed_point_basic}
\end{equation}
Numerically determining $\alpha$ is straightforward, particularly because it is uniquely defined.

\begin{lemma}
Equation \eqref{eq:fixed_point_basic} has a unique solution for all $\beta,\gamma>0$.
\end{lemma}

\begin{proof}
Let $h(\beta):=\f(\beta)/\F(\beta)$ and $w(\beta):=(1-\ee^{-\beta\gamma})/\beta$. Write
 \begin{align}\label{}
f(\beta):=f(\beta,\gamma)=\frac{(1-\beta w(\beta))h(\beta)}{1+w(\beta)h(\beta)},
\end{align}
so that
 \begin{align}\label{}
\beta+f(\beta)=\frac{\beta+h(\beta)}{1+w(\beta)h(\beta)}
\end{align}
For $h(\beta)$ it is known that, see \cite{Sampford1953}, for $\beta\in\mathbb{R}$,
 \begin{align}\label{}
h(\beta)>-\beta, \quad -1<h'(\beta)<0, \quad h''(\beta)>0,
\end{align}
so that $h(\beta)$ is non-negative and non-increasing in $\beta\in\mathbb{R}$, while $\beta+h(\beta)$ is positive and strictly increasing in $\beta\in\mathbb{R}$.
 Because $\ee^x\geq 1+x$,
  \begin{align}\label{}
w'(\beta)=\frac{\ee^{-\beta\gamma}}{\beta^2}(1+\beta\gamma-\ee^{\beta\gamma})\leq 0
\end{align}
 so $w(\beta)$ is also non-negative and non-increasing in $\beta\in \mathbb{R}$. It thus follows that $\beta+f(\beta)$ is strictly increasing in $\beta\in\mathbb{R}$. 
 Moreover, $\beta+f(\beta)\to 0$ as $\beta\to-\infty$ and $\beta+f(\beta)\to \infty$ as $\beta\to \infty$.  
 Let $\Delta=\beta-\alpha$, and rewrite \eqref{eq:fixed_point_basic} as
 \begin{align}\label{}
\beta=\Delta+f(\Delta).
\end{align}
Hence, for each fixed $\beta>0$ there is a unique solution $\Delta \in \mathbb{R}$ from which $\alpha=\beta-\Delta$ follows. 
\end{proof}

\noindent As a result, the delay probability $\P_r({\rm delay})$ and the blocking probability $P_r({\rm block})$ in the model with retrials can be approximated in the QED regime by
\begin{equation}
\P_r({\rm delay}) \approx g(\beta-\alpha,\gamma),
\qquad  
\P_r({\rm block}) \approx \alpha/\sqrt{R},
\label{eq:delay_approx}
\end{equation}
which should become more accurate as $R$ grows large. 

\begin{figure}
\centering
\begin{subfigure}{0.48\textwidth}
\centering
\begin{tikzpicture}[scale=0.75]
\begin{axis}[
	xmin = 101, 
	xmax = 120,
	ymin = 0,
	ymax = 1,
	axis line style={->},
	axis lines = left,
	xscale= 1,
	yscale= 0.8,
	xlabel = {$s$},
	legend cell align=left,
	legend style ={ at={(0.99,1.2)}, 
        anchor=north east,
        align = left}
        ]

\addplot[mark = *] table[x=s,y=sim] {Chapter_4/tikz/basic_model_accuracy_gamma_05.txt};
\addplot[thick,dashed] table[x=s,y=approx] {Chapter_4/tikz/basic_model_accuracy_gamma_05.txt};
\addplot[mark = x] table[x=s,y=block] {Chapter_4/tikz/basic_model_accuracy_gamma_05.txt};

\legend{{$\P_r({\rm delay})$},{$g(\beta-\alpha,\gamma)$}, {$g(\beta,\gamma)$}};
\end{axis}
\end{tikzpicture}
\caption{$\gamma = 0.5$ ($n-s=5$)}
\end{subfigure}
\begin{subfigure}{0.48\textwidth}
\centering
\begin{tikzpicture}[scale =0.75]
\begin{axis}[
	xmin = 101, 
	xmax = 120,
	ymin = 0,
	ymax = 0.25,
	axis line style={->},
	axis lines = left,
	xscale= 1,
	yscale= 0.8,
	xlabel = {$s$},
	scaled ticks=false, tick label style={/pgf/number format/fixed} ,
	legend cell align=left,
	legend style ={ at={(0.99,1.2)}, 
        anchor=north east,
        align = left}
        ]

\addplot[mark = *] table[x=s,y=sim] {Chapter_4/tikz/basic_model_pblock_gamma_05.txt};
\addplot[thick,dashed] table[x=s,y=approx] {Chapter_4/tikz/basic_model_pblock_gamma_05.txt};
\addplot[mark = x] table[x=s,y=block] {Chapter_4/tikz/basic_model_pblock_gamma_05.txt};

\legend{{$\sqrt{R}\,\P_r({\rm block})$},{$\alpha/\sqrt{R}$}, {$f(\beta,\gamma)$}};
\end{axis}
\end{tikzpicture}
\caption{$\gamma = 0.5$ ($n-s=5$)}
\end{subfigure}
\begin{subfigure}{0.48\textwidth}
\centering
\begin{tikzpicture}[scale =0.75]
\begin{axis}[
	xmin = 101, 
	xmax = 120,
	ymin = 0,
	ymax = 1,
	axis line style={->},
	axis lines = left,
	xscale= 1,
	yscale= 0.8,
	xlabel = {$s$},
	legend cell align=left,
	legend style ={ at={(0.99,1.2)}, 
        anchor=north east,
        align = left}
        ]

\addplot[mark = *] table[x=s,y=sim] {Chapter_4/tikz/basic_model_accuracy_gamma_1.txt};
\addplot[thick,dashed] table[x=s,y=approx] {Chapter_4/tikz/basic_model_accuracy_gamma_1.txt};
\addplot[mark = x] table[x=s,y=block] {Chapter_4/tikz/basic_model_accuracy_gamma_1.txt};

\legend{{$\P_r({\rm delay})$},{$g(\beta-\alpha,\gamma)$}, {$g(\beta,\gamma)$}};
\end{axis}
\end{tikzpicture}
\caption{$\gamma = 1$ ($n-s=10$)}
\end{subfigure}
\begin{subfigure}{0.48\textwidth}
\centering
\begin{tikzpicture}[scale =0.75]
\begin{axis}[
	xmin = 101, 
	xmax = 120,
	ymin = 0,
	ymax = 0.25,
	axis line style={->},
	axis lines = left,
	xscale= 1,
	yscale= 0.8,
	xlabel = {$s$},
	scaled ticks=false, tick label style={/pgf/number format/fixed} ,
	legend cell align=left,
	legend style ={ at={(0.99,1.2)}, 
        anchor=north east,
        align = left}
        ]

\addplot[mark = *] table[x=s,y=sim] {Chapter_4/tikz/basic_model_pblock_gamma_1.txt};
\addplot[thick,dashed] table[x=s,y=approx] {Chapter_4/tikz/basic_model_pblock_gamma_1.txt};
\addplot[mark = x] table[x=s,y=block] {Chapter_4/tikz/basic_model_pblock_gamma_1.txt};

\legend{{$\sqrt{R}\,\P_r({\rm block})$},{$\alpha/\sqrt{R}$}, {$f(\beta,\gamma)$}};
\end{axis}
\end{tikzpicture}
\caption{$\gamma = 1$ ($n-s=10$)}
\end{subfigure}
\begin{subfigure}{0.48\textwidth}
\centering
\begin{tikzpicture}[scale =0.75]
\begin{axis}[
	xmin = 101, 
	xmax = 120,
	ymin = 0,
	ymax = 1,
	axis line style={->},
	axis lines = left,
	xscale= 1,
	yscale= 0.8,
	xlabel = {$s$},
	legend cell align=left,
	legend style ={ at={(0.99,1.2)}, 
        anchor=north east,
        align = left}
        ]

\addplot[mark = *] table[x=s,y=sim] {Chapter_4/tikz/basic_model_accuracy_gamma_2.txt};
\addplot[thick,dashed] table[x=s,y=approx] {Chapter_4/tikz/basic_model_accuracy_gamma_2.txt};
\addplot[mark = x] table[x=s,y=block] {Chapter_4/tikz/basic_model_accuracy_gamma_2.txt};

\legend{{$\P_r({\rm delay})$},{$g(\beta-\alpha,\gamma)$}, {$g(\beta,\gamma)$}};
\end{axis}
\end{tikzpicture}
\caption{$\gamma = 2$ ($n-s=20$)}
\end{subfigure}
\begin{subfigure}{0.48\textwidth}
\centering
\begin{tikzpicture}[scale =0.75]
\begin{axis}[
	xmin = 101, 
	xmax = 120,
	ymin = 0,
	ymax = 0.25,
	axis line style={->},
	axis lines = left,
	xscale= 1,
	yscale= 0.8,
	xlabel = {$s$},
	scaled ticks=false, tick label style={/pgf/number format/fixed},
	legend cell align=left,
	legend style ={ at={(0.99,1.2)}, 
        anchor=north east,
        align = left}
        ]

\addplot[mark = *] table[x=s,y=sim] {Chapter_4/tikz/basic_model_pblock_gamma_2.txt};
\addplot[thick,dashed] table[x=s,y=approx] {Chapter_4/tikz/basic_model_pblock_gamma_2.txt};
\addplot[mark = x] table[x=s,y=block] {Chapter_4/tikz/basic_model_pblock_gamma_2.txt};

\legend{{$\sqrt{R}\,\P_r({\rm block})$},{$\alpha/\sqrt{R}$}, {$f(\beta,\gamma)$}};
\end{axis}
\end{tikzpicture}
\caption{$\gamma = 2$ ($n-s=20$)}
\end{subfigure}
\caption{Accuracy of the delay probability approximation in basic model with $R=100$ and $\delta = 0.01$.}
\label{fig:basic_model_accuracy}
\end{figure}
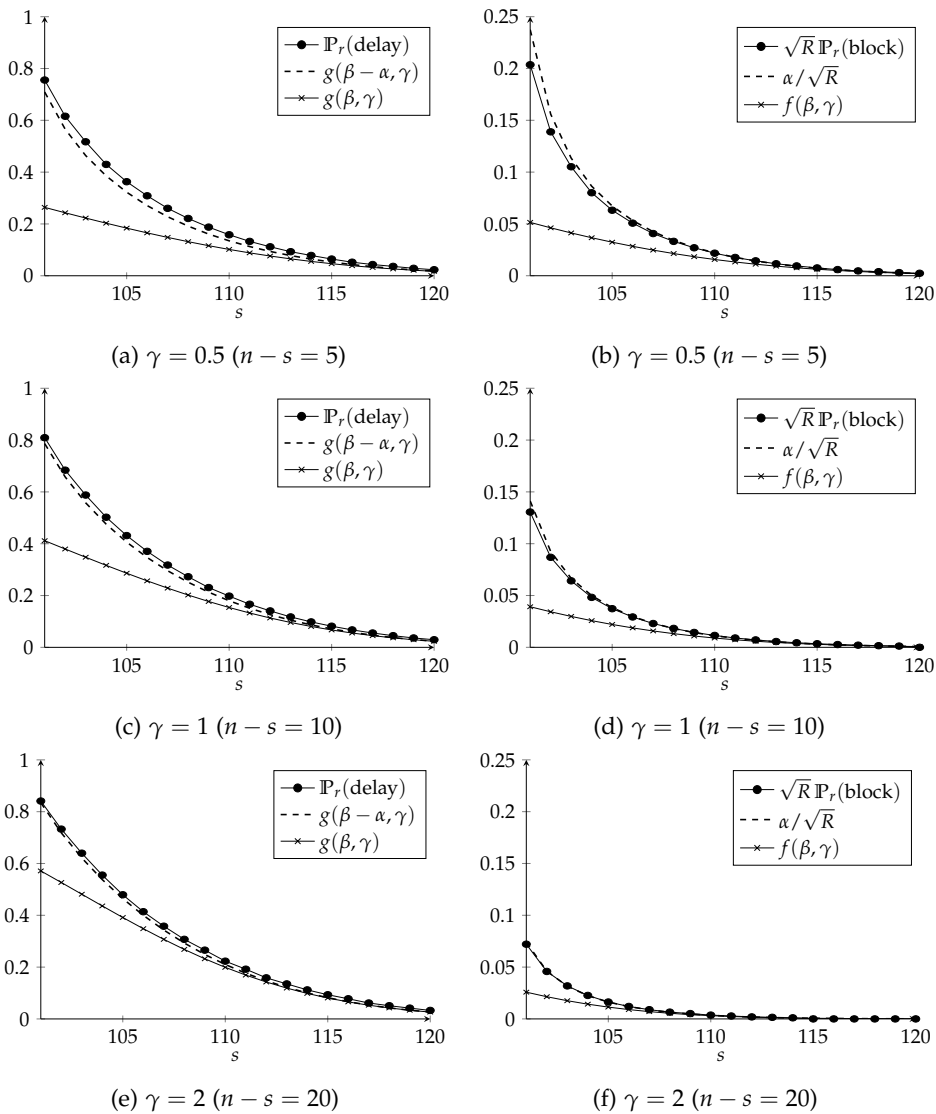

We next test the accuracy of the approximated delay probability in \eqref{eq:delay_approx} in the basic model with slow retrials against the true values obtained through simulation.
Given $R,s$, and $n$, we compute $\beta = (s-R)/\sqrt{R}$ and $\gamma = (n-s)/\sqrt{R}$ in order to approximate the delay and blocking probability as in \eqref{eq:delay_approx} with $\alpha$ as in \eqref{eq:fixed_point_basic}. 
First, we assess the quality of the fixed-point approximation for a large but finite system with $R=100$. 
In Figure \ref{fig:basic_model_accuracy}, we plot the simulated delay probability against the approximation as a function of $s$ (or equivalently $\beta$). 
We consider different values of $\gamma$, namely $\gamma=0.5$, $1$ and $2$, which corresponds to waiting room size $\gamma\sqrt{100} = 5,\, 10$ and 20, respectively.
For comparison, we also include $g(\beta,\gamma)$, the asymptotic delay probability in the system with no retrials, in these plots.

We observe that the heuristic is remarkably accurate in describing both the delay and blocking probability over all values of $\beta,\gamma>0$ considered here. 
The approximation improves as $\gamma$ increases.
Figure \ref{fig:basic_model_accuracy} also clearly illustrates the impact of retrials on the performance measures, which decreases with both $\beta$ and $\gamma$.

\begin{table}
\centering
\small
\begin{tabular}{|r|rr|rr||rr|rr|}
\cline{2-9}\multicolumn{1}{r|}{} & \multicolumn{4}{c||}{$(\b,\g) = (0.5,0.5)$} & \multicolumn{4}{c|}{$(\b,\g) = (1,0.5)$} \bigstrut\\
\hline
$R$    & $s$     & $n$     & $\P_r({\rm delay})$ & $\sqrt{R}\P_r({\rm bl.})$ & $s$     & $n$     & $\P_r({\rm delay})$ & $\sqrt{R}\P_r({\rm bl.})$  \bigstrut\\
\hline
5     & 6     & 7     & 0.5019 & 0.5982 & 7     & 8     & 0.2607 & 0.2610 \bigstrut[t]\\
10    & 12    & 14    & 0.3697 & 0.3679 & 13    & 15    & 0.2298 & 0.1966 \\
50    & 54    & 58    & 0.3509 & 0.4931 & 57    & 61    & 0.1765 & 0.2019 \\
100   & 105   & 110   & 0.3640 & 0.6336 & 110   & 115   & 0.1579 & 0.2178 \\
500   & 511   & 522   & 0.3460 & 0.6780 & 522   & 533   & 0.1482 & 0.2297 \\
1000  & 1016  & 1032  & 0.3333 & 0.6481 & 1032  & 1048  & 0.1412 & 0.2141 \bigstrut[b]\\
\hline
\multicolumn{1}{r|}{} & \multicolumn{2}{r|}{Approx} & 0.3225 & 0.6734 & \multicolumn{2}{r|}{Approx} & 0.1349 & 0.2206 \bigstrut\\
\cline{2-9}\end{tabular}%

\vspace{5 mm} 
\begin{tabular}{|r|rr|rr||rr|rr|}
\cline{2-9}\multicolumn{1}{r|}{} & \multicolumn{4}{c||}{$(\b,\g) = (0.5,1)$} & \multicolumn{4}{c|}{$(\b,\g) = (1,1)$} \bigstrut\\
\hline
$R$    & $s$     & $n$     & $\P_r({\rm delay})$ & $\sqrt{R}\P_r({\rm bl.})$ & $s$     & $n$     & $\P_r({\rm delay})$ & $\sqrt{R}\P_r({\rm bl.})$  \bigstrut\\
\hline
5     & 6     & 8     & 0.5337 & 0.4065 & 7     & 9     & 0.2866 & 0.1612  \bigstrut[t]\\
10    & 12    & 15    & 0.3932 & 0.2701 & 13    & 16    & 0.2472 & 0.1374 \\
50    & 54    & 61    & 0.3993 & 0.3171 & 57    & 64    & 0.2063 & 0.1183 \\
100   & 105   & 115   & 0.4333 & 0.3754 & 110   & 120   & 0.1971 & 0.1143 \\
500   & 511   & 533   & 0.4247 & 0.3986 & 522   & 544   & 0.1928 & 0.1202 \\
1000  & 1016  & 1048  & 0.4115 & 0.3689 & 1032  & 1064  & 0.1831 & 0.1088  \bigstrut[b]\\
\hline
\multicolumn{1}{r|}{} & \multicolumn{2}{r|}{Approx} & 0.4062 & 0.3828 & \multicolumn{2}{r|}{Approx} & 0.1798 & 0.1106 \bigstrut\\
\cline{2-9}\end{tabular}%

\caption{\normalsize Numerical results of the fixed-point method for the basic model as $R\to\infty$.}
\label{tab:basic_model_accuracy}
\end{table}

Table \ref{tab:basic_model_accuracy} furthermore shows how the accuracy of the approximations increases as $R$ increases. 
In this table, we used the simulated delay and blocking probability for systems of increasing size while adhering to the two-fold scaling rule of \eqref{eq:twofold_scaling_basic_model}.
The values of $s$ and $n$ are rounded to the nearest integer.

\section{Cloud model}
\label{sec:cloud_model}
The second model we consider in this chapter is inspired by cloud computing services. 
We shall see how our fixed-point heuristic helps cloud providers in their provisioning process.\\
\\*
\subsection{Practical context}
Cloud computing enables network access to a shared pool of configurable computing resources, allowing users (e.g.~companies, service providers) to store and process their data in third-party data centers, without investing in the operating equipment themselves.
At the foundation of cloud computing lies the idea of sharing resources to achieve economies-of-scale in terms of maximizing computing power usage and
reducing the overall cost of resources such as energy and
infrastructure. 
Cloud providers, such as Amazon EC2, Windows Azure and Rackspace \cite{Armbrust2010}, offer virtual machine (VM) provisioning, which allows users to request VM instances configured to their preference. 
In a service system context, the provider thus serves users by supplying them with a VM that matches their requirements, running on one of the cloud's physical machines.

\begin{figure}
\centering
\includegraphics[scale=0.33]{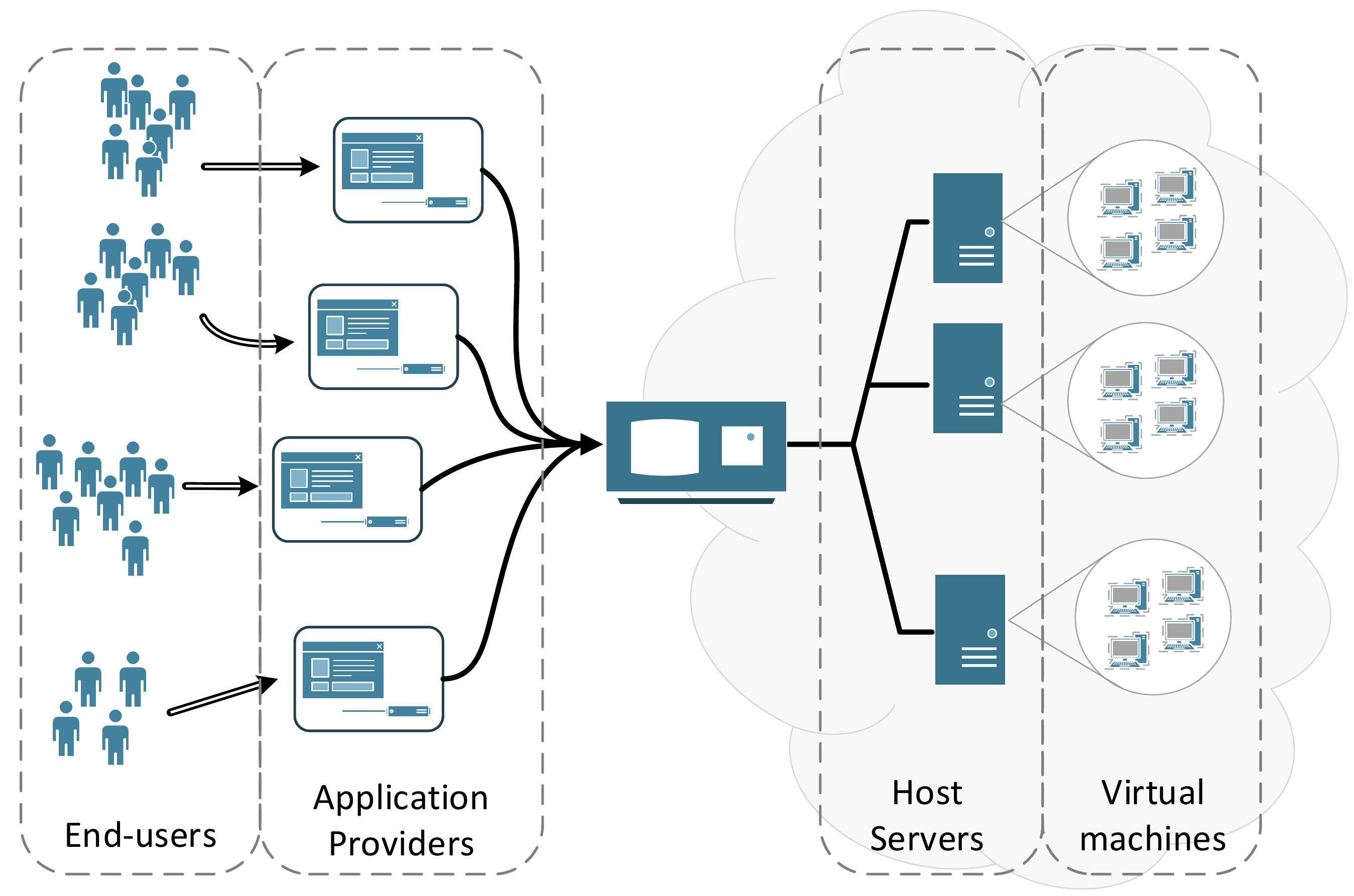}
\caption{Cloud provisioning process}
\label{fig:CloudScheme}
\end{figure}

Let us describe the cloud provisioning process in more detail; see Figure \ref{fig:CloudScheme}.
At the highest granularity level there are the \emph{end-users}, devices typically directly operated by humans, using an \emph{application provider} (AP), usually a company that provides software usage over the internet (e.g.~SaaS \cite{Mell2011}). To some extent, the AP will rely on a static set of computing resources, but certainly in case of sudden surges in workload, these might not be sufficient. When the AP recognizes the need for additional capacity, for instance by \emph{auto-scaling} procedures \cite{Amazon}, a VM request is submitted to the cloud provider. The request is handled by a \emph{host server} that starts the set-up of the VM with requested specifications. This includes elementary operations such as copying the VM image and assigning an IP address. Each server is able to host multiple VM instances in parallel, although the VMs in set-up need their dedicated attention, due to concurrency level constraints incurred by large I/O activities. Once the set-up is completed, the VM is ready for use, and the AP may start using the additional computing resources. 

Our focus lies on the capacity allocation within the cloud environment, so the right-hand side of Figure \ref{fig:CloudScheme}. 
Successful management of cloud systems requires the right scaling of both the number of host servers (denoted by $s$) at the first I/O queue and the maximum number of VMs (denoted by $n$) that can be hosted simultaneously.  Moreover, this needs to be done in a dynamic way in order to respond effectively to the time-varying demand. The capacity $n$ defines a hard constraint on whether a new VM request will be accepted immediately or not. Therefore, new requests will be delayed or even dropped if the available host capacity is insufficient, which is more likely to occur during periods in which the $s$ host servers are overloaded.\\
\\*
\subsection{Queueing model}
\begin{figure}
\centering
\begin{tikzpicture}[scale=1]

\draw[dashed] (0,0.25) rectangle (6.5,2.75) node[right,above] {$n$};
\draw[thick] (0.75,1) -- (2.5,1) -- (2.5,2) -- (0.75,2);
\draw[thick] (1.5,1.1) -- (1.5,1.9);
\draw[thick] (1.75,1.1) -- (1.75,1.9);
\draw[thick] (2,1.1) -- (2,1.9);
\draw[thick] (2.25,1.1) -- (2.25,1.9);
\draw[thick] (3.5,1.5) circle (0.5) node {$s$};
\draw[thick] (5.5,1.5) circle (0.5) node {$\infty$};

\draw[thick,->] (-2.25,1.5)  -- (1,1.5);
\draw[thick,->] (4.25,1.5) -- (4.75,1.5);
\draw[thick,->] (6.25,1.5) -- (7.55,1.5);

\node at (-2.25,1.8) {Pois$(\lambda)$};
\node at (3.5,2.35) {$\exp(\mu)$};
\node at (5.5,2.35) {$\exp(\kappa)$};

\draw[thick] (-0.5,1.5) -- (-0.5,1);
\draw[thick] (-0.5,1) to [in=0,out=270] (-1,0.5);
\draw[thick] (-1,0.5) to [in=270,out=180] (-1.5,1);
\draw[thick,->] (-1.5,1) to [in=180,out=90] (-1,1.45);

\node at (-2,0.4) {$\exp(\delta)$};

\end{tikzpicture}
\caption{Abstracted model of VM provisioning process.}
\label{fig:CloudModel}
\end{figure}
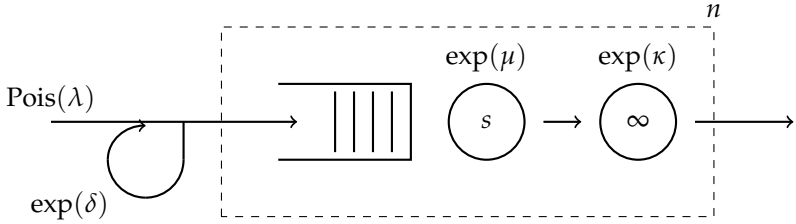

To describe the cloud system in mathematical terms, we extend the model proposed by Tan et al.~\cite{Tan2012}. 
Each host server may host a number of VM instances at the same time, yielding a total number of $n$ parallel VM instances. Requests, arriving to the system according to a Poisson process with rate $\l$, are granted only if one of these $n$ positions is available.
Otherwise, the user retries getting access after an exponentially distributed time with mean $1/\delta$. 
If granted, the request is assigned to a host server not busy initializing another VM instance, if available, or waits for one to become available. This start-up time is assumed to be exponentially distributed with mean $1/\mu$. On completion of the initialization phase, VM usage is initiated by the client. The VM continues to be occupied for a random amount of time, with mean $1/\kappa$, until release by the user. We note that the model of Tan et al. \cite{Tan2012} has three queues in tandem, one $M/M/s$ queue, followed by two $M/M/\infty$ queues that separately model a second initialization phase and the actual VM usage by the cloud user. We thus replace the two $M/M/\infty$ queues by one $M/G/\infty$ queue with an aggregated service time, which does not alter the system performance analysis. 
This yields the queueing model in Figure~\ref{fig:CloudModel}. 

\begin{remark}
We mention that a queueing model similar to the one in Figure \ref{fig:CloudModel} without retrials is analyzed by Khudyakov et al.~\cite{Khudyakov2006} in a telecommunication environment. 
In their work, $s$ and $n$ represent the number of agents and trunk lines in a call center. 
Although the order of the two queues is switched, the stationary analysis of their model and the cloud model is the same, due to the product-form structure of the stationary distribution.
In fact, Tan et al.~\cite{Tan2012} use the results of \cite{Khudyakov2006} in their asymptotic analysis.
\end{remark}

An exact analysis of the cloud model is again obstructed by the absence of a product-form solution in case of retrials. 
We therefore turn to the QED paradigm to approximate the system behavior as $R\to\infty$. 

Following the approach in \cite{Khudyakov2006,Tan2012}, we argue that the appropriate QED scaling for $s$ and $n$ should be
\begin{align}
s &= R + \beta\sqrt{R} + o(\sqrt{R}),  &\beta>0, \nonumber\\
n &= s + R/\kappa + \gamma\sqrt{R/\kappa} + o(\sqrt{R}), &\gamma>0,
\label{eq:twofold_scaling_cloud_model}
\end{align}
where $R = \lambda / \mu = \lambda$. 
To understand why this is indeed the correct scaling regime to obtain non-degenerate limiting behavior, we recall the arguments we presented in Section 1.3.
Namely, in order to achieve QED performance, one allocates the nominal workload brought towards the queue plus a variability hedge that is proportional to the square-root of this amount. 
For $s$, this results in the standard square-root staffing rule. 
For $n$, this is the sum of the capacity needed at the multi-server queue, i.e. $s$, and the consecutive infinite-server queue. 
Since the expected workload at the second queue equals $R/\kappa$, the capacity required at this stage equals $R/\kappa+\gamma\sqrt{R/\kappa}$ for some $\gamma$. 
In total, this yields the scaling for the number of VMs $n$ as in \eqref{eq:twofold_scaling_cloud_model}.
The limiting behavior of this queueing model without retrials is documented in \cite{Khudyakov2006,Tan2012}.

\begin{proposition}\label{prop:cloud_model_limits}
Let $s$ and $n$ in the cloud model of Figure \ref{fig:CloudModel} without retrials scale as in \eqref{eq:twofold_scaling_cloud_model}. Then, as $R\to\infty$,
\begin{align}\label{pdelay}
\P^c({\rm delay})&\to \frac{\xi_1-\xi_2}{\eta+\xi_1-\xi_2} =:g_c(\beta,\gamma),\\
\label{pblock}
\sqrt{R}\cdot \P^c({\rm block}) &\to \frac{\nu}{\eta+\xi_1-\xi_2} =:f_c(\beta,\gamma),
\end{align}
where 
\[
\eta = \int_{-\infty}^\b \F\Big(\gamma+(\beta-t)\sqrt{\kappa}\Big)\,\f(t)\, \dd t, \qquad
\xi_1 = \frac{ \f(\b)\F(\gamma) }{\b},
\]
\[
\xi_2 = \frac{1}{\b}\,\f\left(\sqrt{\beta^2+\gamma^2}\right)\ee^{\tfrac{1}{2}(\gamma-\beta/\sqrt{\kappa})^2} \F(\gamma-\beta/\sqrt{\kappa}),
\]
\[
\nu = \sqrt{\frac{\kappa}{1+\kappa}}\,\f\Big(\frac{\gamma+\beta\sqrt{\kappa}}{\sqrt{1+\kappa}}\Big) \F\Big(\frac{\beta-\gamma\sqrt{\kappa}}{\sqrt{1+\kappa}}\Big) + \beta\,\xi_2.
\]
\end{proposition}

\subsection{Fixed-point method}
Proposition \ref{prop:cloud_model_limits} shows that also in this model, the blocking probability vanishes at rate $1/\sqrt{R}$, making it amenable to our fixed-point method for retrials.
Let $\alpha\sqrt{R}$ the volume of retrials, so that a total arrival rate is $R_{\rm tot} = R + \alpha\sqrt{R}$, or equivalently $R = R_{\rm tot} - \alpha \sqrt{R_{\rm tot}} + o(\sqrt{R_{\rm tot}})$. 
Substituting this into the two-fold scaling rule in \eqref{eq:twofold_scaling_cloud_model} gives
\begin{align*}
s &= R_{\rm tot} + (\beta - \alpha ) \sqrt{R_{\rm tot}} + o(\sqrt{R_{\rm tot}}),\\
n &= s + \frac{R_{\rm tot}}{\kappa} + \Big( \gamma - \frac{\alpha}{\sqrt{\kappa}} \Big) \sqrt{\frac{R_{\rm tot}}{\kappa}} + o(\sqrt{R_{\rm tot}}).
\end{align*}
Accordingly, the constant $\alpha$ is defined as the solution of the fixed-point equation
\begin{equation}
\label{eq:fixed_point_cloud}
f_c\left(\beta - \alpha, \gamma - \alpha/\sqrt{\kappa}\right) = \alpha.
\end{equation}
Approximations for the delay and blocking probability in the cloud model with retrials are hence given by
\begin{equation}
\P_r^c({\rm delay}) \approx g_c\left(\beta-\alpha,\gamma-\alpha/\sqrt{\kappa}\right), \qquad \P_r^c({\rm block}) \approx \alpha/\sqrt{R}.
\label{eq:cloud_approximations}
\end{equation}
Note that in contrast to the fixed-point equation \eqref{eq:fixed_point_basic} for the basic model, the second argument $\gamma-\alpha/\sqrt{R}$ is also corrected.
\\
\\*
We next test the accuracy of the fixed-point equation for several instances. 
In Table \ref{tab:cloud_model_accuracy}, we present the simulation results for $\kappa=0.02$, 0.2 and 1, and two pairs of $(\beta,\gamma)$ for increasing $R$.
\begin{table}
\centering
\small
\begin{subtable}{0.99\textwidth}\centering
\begin{tabular}{|r|rr|rr||rr|rr|}
\cline{2-9}\multicolumn{1}{r|}{} & \multicolumn{4}{c||}{$(\beta,\gamma) = (0.5,1)$} & \multicolumn{4}{c|}{$(\beta,\gamma) = (1,1)$} \bigstrut\\
\hline
$R$     & $s$     & $n$     & $\P^c_r({\rm del})$ & $\sqrt{R}\,\P_r^c({\rm bl})$ & $s$     & $n$     & $\P^c_r({\rm del})$ & $\sqrt{R}\,\P_r^c({\rm bl})$ \bigstrut\\
\hline
5     & 6     & 13    & 0.5309 & 0.5540 & 7     & 14    & 0.2800 & 0.2426 \bigstrut[t]\\
10    & 12    & 25    & 0.3864 & 0.3810 & 13    & 26    & 0.2393 & 0.2164 \\
50    & 54    & 111   & 0.3904 & 0.4525 & 57    & 114   & 0.1965 & 0.2010 \\
100   & 105   & 215   & 0.4300 & 0.5474 & 110   & 220   & 0.1859 & 0.1952 \\
500   & 511   & 1033  & 0.4139 & 0.5586 & 522   & 1044  & 0.1787 & 0.2052 \\
1000  & 1016  & 2048  & 0.4003 & 0.5479 & 1032  & 2064  & 0.1660 & 0.1803 \bigstrut[b]\\
\hline
\multicolumn{1}{r|}{} & \multicolumn{2}{r|}{Approx} & 0.4029 & 0.5638 & \multicolumn{2}{r|}{Approx} & 0.1709 & 0.1992 \bigstrut\\
\cline{2-9}\end{tabular}%
\caption{$\kappa = 1$}
\end{subtable}

\vspace{5mm}

\begin{subtable}{0.99\textwidth}\centering
\begin{tabular}{|r|rr|rr||rr|rr|}
\cline{2-9}\multicolumn{1}{r|}{} & \multicolumn{4}{c||}{$(\beta,\gamma) = (0.5,1)$} & \multicolumn{4}{c|}{$(\beta,\gamma) = (1,1)$} \bigstrut\\
\hline
$R$     & $s$     & $n$     & $\P^c_r({\rm del})$ & $\sqrt{R}\,\P_r^c({\rm bl})$ & $s$     & $n$     & $\P^c_r({\rm del})$ & $\sqrt{R}\,\P_r^c({\rm bl})$ \bigstrut\\
\hline
5     & 6     & 36    & 0.5664 & 0.3049 & 7     & 37    & 0.3079 & 0.1457 \bigstrut[t]\\
10    & 12    & 69    & 0.4263 & 0.2227 & 13    & 70    & 0.2683 & 0.1410 \\
50    & 54    & 320   & 0.4444 & 0.2555 & 57    & 323   & 0.2293 & 0.1334 \\
100   & 105   & 627   & 0.4826 & 0.3085 & 110   & 632   & 0.2187 & 0.1379 \\
500   & 511   & 3061  & 0.4842 & 0.3235 & 522   & 3072  & 0.2182 & 0.1358 \\
1000  & 1016  & 6087  & 0.4630 & 0.2906 & 1032  & 6103  & 0.2039 & 0.1332 \bigstrut[b]\\
\hline
\multicolumn{1}{r|}{} & \multicolumn{2}{r|}{Approx} & 0.4687 & 0.3029 & \multicolumn{2}{r|}{Approx} & 0.2042 & 0.1326 \bigstrut\\
\cline{2-9}\end{tabular}%
\caption{$\kappa=0.2$}
\end{subtable}

\vspace{5mm}

\begin{subtable}{0.99\textwidth}\centering
\begin{tabular}{|r|rr|rr||rr|rr|}
\cline{2-9}\multicolumn{1}{r|}{} & \multicolumn{4}{c||}{$(\beta,\gamma) = (0.5,1)$} & \multicolumn{4}{c|}{$(\beta,\gamma) = (1,1)$} \bigstrut\\
\hline
$R$     & $s$     & $n$     & $\P^c_r({\rm del})$ & $\sqrt{R}\,\P_r^c({\rm bl})$ & $s$     & $n$     & $\P^c_r({\rm del})$ & $\sqrt{R}\,\P_r^c({\rm bl})$ \bigstrut\\
\hline
5     & 6     & 272   & 0.5836 & 0.1085 & 7     & 273   & 0.3217 & 0.0738 \bigstrut[t]\\
10    & 12    & 534   & 0.4456 & 0.0945 & 13    & 535   & 0.2822 & 0.0781 \\
50    & 54    & 2604  & 0.4683 & 0.1031 & 57    & 2607  & 0.2429 & 0.0764 \\
100   & 105   & 5176  & 0.5106 & 0.1064 & 110   & 5181  & 0.2345 & 0.0759 \\
500   & 511   & 25669 & 0.5130 & 0.1158 & 522   & 25680 & 0.2353 & 0.0795 \\
1000  & 1016  & 51240 & 0.4946 & 0.0975 & 1032  & 51256 & 0.2223 & 0.0744 \bigstrut[b]\\
\hline
\multicolumn{1}{r|}{} & \multicolumn{2}{r|}{Approx} & 0.4999 & 0.0862 & \multicolumn{2}{r|}{Approx} & 0.2207 & 0.0595 \bigstrut\\
\cline{2-9}\end{tabular}%
\caption{$\kappa=0.02$}
\end{subtable}
\caption{\normalsize Numerical results of the fixed-point method for the cloud model with slow retrials as $R\to\infty$.}
\label{tab:cloud_model_accuracy}
\end{table}

First, observe from Table \ref{tab:cloud_model_accuracy} that $n$ now lives on a different scale than $s$. 
This is required to facilitate the long sojourn time of customers in the second stage, which is proportional to $1/\kappa$, creating the need for larger system size. 
Besides that, the numerical results show that the fixed-point approximation is again remarkably accurate over a wide range of parameter settings.
Even for cloud systems as small as 50 servers, the fixed-point method gives accurate approximations.

\section{Abandonments}
\label{sec:abandonments_retrials}

Whereas in the basic model of Section \ref{sec:basic_model}, retrials were governed by the system architecture (arriving customers are requested to reattempt if $n$ customers are present in the system), we now consider a setting in which departures from the queue are customer-initiated. 
That is, customers deliberately decide to leave the queue to return for service at a later time.
Hence we consider a queueing system with abandonments and retrials. 

\subsection{The Erlang-A model}
The canonical model for abandonments is the $M/M/s+M$ or Erlang-A model \cite{Palm1957,Garnett2002}. 
The queueing dynamics of the Erlang-A model are similar to those in the $M/M/s$ queue, with the additional feature that each customer is assigned an i.i.d.~patience time, which is exponentially distributed with mean $1/\theta$.
If a customer's patience time expires before reaching an available server, she leaves (abandons) the system. 
As the number of customers in the Erlang-A queue is a birth-death process, its stationary distribution and associated performance measures are fairly well-understood, also in the QED regime \cite{Garnett2002,Zeltyn2005,Zhang2012}.
Most importantly to us, Garnett et al. \cite{Garnett2002} and Zeltyn \& Mandelbaum \cite{Zeltyn2005} identified the asymptotic delay and abandonment probability in the Erlang-A model under QED scaling. 

\begin{proposition}{ \cite[Thm.~4.1]{Zeltyn2005} }
\label{prop:abandonment_prop}
Let $s = R+\beta\sqrt{R} + o(\sqrt{R})$ for some $\b\in\mathbb{R}$. Then in the $M/M/s+M$ queue
\begin{align}
\P^a({\rm delay}) &\to \left( 1 + \sqrt{\theta}\,\frac{h(\beta/\sqrt{\theta})}{h({-}\beta)}\right)^{-1} =: g_a(\beta) \\
\sqrt{R}\,\P^a({\rm abandon}) &\to  
\frac{ \sqrt{\theta}\,h(\beta/\sqrt{\theta})- \beta }
{ 1 + \sqrt{\theta}\,h(\beta/\sqrt{\theta})/h({-}\beta)} =: f_a(\beta),
\end{align}
as $R\to\infty$ where $h(\beta) = \f(\beta)/\F({-}\beta)$.
\end{proposition}

We remark that in \cite{Zeltyn2005}, the QED limits for generally distributed patience time were derived. 
Although our heuristic also works for this more general setting, we focus on the exponentially distributed patience here to convey our main ideas. \\

\begin{remark}
Large-scale Markovian multi-server queues with abandonments and retrials have been thoroughly studied in a series of papers by Mandelbaum et al.~\cite{Mandelbaum1999,Mandelbaum1999a,Mandelbaum2002}. 
In these works, the authors consider a system with time-varying arrivals and a retrial rate that remains bounded away from zero, for which they deduce fluid and diffusion limits as the system grows large.
These limits provide approximations for the time-dependent queue length and virtual waiting time processes, including their means and variances. 
We in this section take a different approach by assuming $\delta\to 0$, which enables us to characterize the steady-state behavior of queues with abandonments and retrials. 
\end{remark}

\subsection{Fixed-point method}
Next, we include (slow) retrials. 
More specifically, we assume that customers who abandon the queue rejoin the queue after an exponentially distributed time with mean $1/\delta \gg 1$. 
Just as in the $M/M/s/n$ queue, the $M/M/s+M$ queue with retrials is analytically intractable, and therefore we apply our fixed-point method to approximate its performance in the QED regime. 

Observe through Proposition \ref{prop:abandonment_prop} that the fraction of customers leaving before receiving service is roughly $\alpha/\sqrt{R}$. 
Following the reasoning of Section \ref{sec:fixed_point}, the total arrival volume, consisting of new arrivals and reattempting customers, is $R_{\rm tot} = R+\alpha \sqrt{R}$, with
\begin{equation}
 \label{eq:abandon_fixed_point}
 \alpha = f_a(\beta - \alpha).
 \end{equation} 
Accordingly, this yields the following approximations for the delay and abandonment probability in the system with retrials
\begin{equation}
\label{eq:abandon_approximation}
\P_r^a({\rm delay}) \approx g_a(\beta-\alpha),\qquad \P_r^a({\rm abandon}) \approx \alpha / \sqrt{R}.
\end{equation}

We test our heuristic in the model with abandonment with parameters $R=100$, $\delta = 0.01$.
For $\theta = 0.2$, customers are quite patient, as they are willing to wait on average 5 times their expected service time. 
Customer abandonment becomes more dominant for the cases in which customers are reasonably patient $\theta=1$ and very impatient $\theta = 10$. 

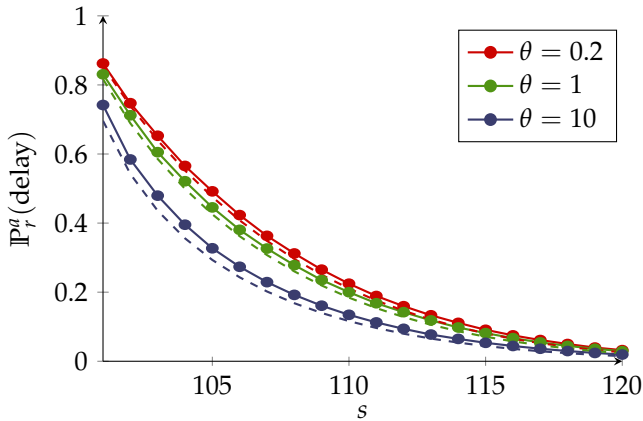
\begin{figure}
\centering
\begin{tikzpicture}
\begin{axis}[
	xmin = 101, 
	xmax = 120,
	ymin = 0,
	ymax = 1,
	axis line style={->},
	axis lines = left,
	xscale= 1,
	yscale= 0.8,
	xlabel = {$s$},
	ylabel = {$\P^a_r({\rm delay})$},
	legend cell align=left,
	legend style ={ at={(0.99,1.2)}, 
        anchor=north east,
        align = left}
        ]

\addplot[col1,thick,mark = *] table[x=s,y=sim_02] {Chapter_4/tikz/abandonments_delay.txt};
\addplot[col4,thick,mark = *] table[x=s,y=sim_1] {Chapter_4/tikz/abandonments_delay.txt};
\addplot[col5,thick,mark = *] table[x=s,y=sim_10] {Chapter_4/tikz/abandonments_delay.txt};

\addplot[col1, thick,dashed] table[x=s,y=approx_02] {Chapter_4/tikz/abandonments_delay.txt};
\addplot[col4, thick,dashed] table[x=s,y=approx_1] {Chapter_4/tikz/abandonments_delay.txt};
\addplot[col5, thick,dashed] table[x=s,y=approx_10] {Chapter_4/tikz/abandonments_delay.txt};

\legend{{$\theta = 0.2$},{$\theta = 1$}, {$\theta = 10$}};
\end{axis}
\end{tikzpicture}
\caption{Simulated (solid) and approximated (dashed) delay probability in the $M/M/s+M$ queue with retrials and $R=100$, $\delta=0.01$.}
\label{fig:abandonment_accuracy_delay}
\end{figure}

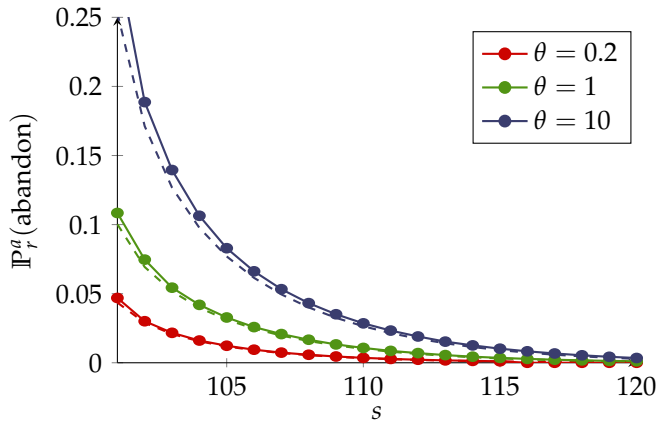
\begin{figure}
\centering
\begin{tikzpicture}
\begin{axis}[
	xmin = 101, 
	xmax = 120,
	ymin = 0,
	ymax = 0.25,
	axis line style={->},
	axis lines = left,
	xscale= 1,
	yscale= 0.8,
	xlabel = {$s$},
	legend cell align=left,
	ylabel = {$\P^a_r({\rm abandon})$},
	scaled ticks=false, tick label style={/pgf/number format/fixed},
	legend style ={ at={(0.99,1.2)}, 
        anchor=north east,
        align = left}
        ]

\addplot[col1,thick,mark = *] table[x=s,y=sim_02] {Chapter_4/tikz/abandonments_ab.txt};
\addplot[col4,thick,mark = *] table[x=s,y=sim_1] {Chapter_4/tikz/abandonments_ab.txt};
\addplot[col5,thick,mark = *] table[x=s,y=sim_10] {Chapter_4/tikz/abandonments_ab.txt};

\addplot[col1, thick,dashed] table[x=s,y=approx_02] {Chapter_4/tikz/abandonments_ab.txt};
\addplot[col4, thick,dashed] table[x=s,y=approx_1] {Chapter_4/tikz/abandonments_ab.txt};
\addplot[col5, thick,dashed] table[x=s,y=approx_10] {Chapter_4/tikz/abandonments_ab.txt};

\legend{{$\theta = 0.2$},{$\theta = 1$}, {$\theta = 10$}};
\end{axis}
\end{tikzpicture}
\caption{Simulated (solid) and approximated (dashed) abandonment probability in the $M/M/s+M$ queue with retrials and $R=100$, $\delta=0.01$.}
\label{fig:abandonment_accuracy_ab}
\end{figure}

In Figures \ref{fig:abandonment_accuracy_delay} and \ref{fig:abandonment_accuracy_ab}, we plot the simulated delay and (scaled) abandonment probability against approximations obtained through the fixed-point method. \\
\noindent 
Again, we see a very good match between approximated and actual values. 
As $\theta$ decreases, that is, customers become more patient, accuracy of the approximations improves. 
This makes sense, since the volume of retrials decreases, and the system behaves more and more like a standard $M/M/s$ queue.
 
In Table \ref{tab:numerical_accuracy_abandonments} we also check the asymptotic accuracy of the model with abandonments and retrials and see that the approximation indeed improves as $R$ increases. 

\begin{table}
\centering
\small
\begin{subtable}{0.99\textwidth}
\centering
\begin{tabular}{|r|r|rr||rr||rr|}
\cline{3-8}
\multicolumn{2}{c|}{} & \multicolumn{2}{c||}{$\theta = 0.2$} & \multicolumn{2}{c||}{$\theta = 1$} & \multicolumn{2}{c|}{$\theta = 10$} \bigstrut\\
\hline
$R$   & $s$   & $\P^a_r({\rm del})$ & $\sqrt{R}\,\P_r^a({\rm ab})$ & $\P^a_r({\rm del})$ & $\sqrt{R}\,\P_r^a({\rm ab})$ & $\P^a_r({\rm del})$ & $\sqrt{R}\,\P_r^a({\rm ab})$ \bigstrut\\
\hline
5     & 6     & 0.5703 & 0.1522 & 0.5423 & 0.4158 & 0.4766 & 1.0980 \bigstrut[t]\\
10    & 12    & 0.6612 & 0.2152 & 0.6277 & 0.5619 & 0.5429 & 1.5099 \\
50    & 54    & 0.5521 & 0.1517 & 0.5089 & 0.4009 & 0.3920 & 1.0251 \\
100   & 105   & 0.4896 & 0.1218 & 0.4456 & 0.3276 & 0.3282 & 0.8321 \\
500   & 511   & 0.4877 & 0.1228 & 0.4442 & 0.3302 & 0.3132 & 0.8135 \\
1000  & 1016  & 0.4992 & 0.1274 & 0.4472 & 0.3359 & 0.3148 & 0.8244 \bigstrut[b]\\
\hline
\multicolumn{2}{|r|}{Approx} & 0.4757 & 0.1182 & 0.4254 & 0.3120 & 0.2933 & 0.7695 \bigstrut\\
\hline
\end{tabular}
\caption{$\beta = 0.5$}
\end{subtable}

\vspace{5mm}

\begin{subtable}{0.99\textwidth}
\centering
\begin{tabular}{|r|r|rr||rr||rr|}
\cline{3-8}
\multicolumn{2}{c|}{} & \multicolumn{2}{c||}{$\theta = 0.2$} & \multicolumn{2}{c||}{$\theta = 1$} & \multicolumn{2}{c|}{$\theta = 10$} \bigstrut\\
\hline
$R$   & $s$   & $\P^a_r({\rm del})$ & $\sqrt{R}\,\P_r^a({\rm ab})$ & $\P^a_r({\rm del})$ & $\sqrt{R}\,\P_r^a({\rm ab})$ & $\P^a_r({\rm del})$ & $\sqrt{R}\,\P_r^a({\rm ab})$ \bigstrut\\
\hline
5     & 7     & 0.3130 & 0.0527 & 0.2908 & 0.1563 & 0.2382 & 0.4033 \bigstrut[t]\\
10    & 13    & 0.2732 & 0.0444 & 0.2503 & 0.1331 & 0.1917 & 0.3504 \\
50    & 57    & 0.2344 & 0.0373 & 0.2090 & 0.1120 & 0.1442 & 0.2974 \\
100   & 110   & 0.2244 & 0.0357 & 0.1999 & 0.1077 & 0.1355 & 0.2877 \\
500   & 522   & 0.2232 & 0.0355 & 0.1973 & 0.1079 & 0.1282 & 0.2842 \\
1000  & 1032  & 0.2210 & 0.0359 & 0.1979 & 0.1092 & 0.1287 & 0.2890 \bigstrut[b]\\
\hline
\multicolumn{2}{|r|}{Approx} & 0.2105 & 0.0335 & 0.1842 & 0.1005 & 0.1162 & 0.2642 \bigstrut\\
\hline
\end{tabular}%
\caption{$\beta=1$}
\end{subtable}
\caption{\normalsize Numerical results of the fixed-point method for the $M/M/s+M$ queue with slow retrials.}
\label{tab:numerical_accuracy_abandonments}
\end{table}

\begin{remark}
Note that even though the fixed-point approximations come close to the simulated values as $R$ increases, a small gap remains, especially notable in $\beta=0.5$. 
This can be attributed to both rounding errors and the heuristic assumption that the retrial stream is independent Poisson.
The latter is obviously false, as the retrial process naturally depends on the history of the external arrival and service processes.
Therefore, the fixed-point heuristic slightly underestimates congestion levels in the actual system. 
However, this error is relatively small and moreover in small to moderate-size systems negligible compared to the effects of rounding. 
\end{remark}

\begin{remark}
Our fixed-point heuristic easily extends to the case in which only a fraction of abandoning customers returns to the system later. 
If each customer who abandons decides (independent from others and his own retrial history) to return with probability $q\in[0,1]$, then the arrival stream due to retrial becomes $q\cdot\alpha\sqrt{R}$, so that the fixed-point becomes $f_a(\beta-q\alpha) = \alpha$. 
Approximations of the performance measures follow accordingly.
\end{remark}

\section{Dimensioning}
\label{sec:retrial_dimensioning}
The asymptotic QED expressions for the systems we considered in Sections \ref{sec:basic_model}-\ref{sec:abandonments_retrials} without retrials together with the corrections obtained through the fixed-point equation provide a method for dimensioning large-scale systems with retrials. 
For sufficiently large arrival volumes, we can tune the QoS-levels offered by the systems through the QoS-parameters $\b$ and $\g$. 
In this section we demonstrate how to do so in the cloud model, using the delay probability as a vehicle.
First we explore the procedure under stationary conditions, then in a time-varying environment.
The methods we propose easily translate to the two other model settings considered in this chapter, and the blocking probability.

\subsection{Stationary dimensioning}
We consider the dimensioning problem in the cloud model from a constraint satisfaction perspective. 
That is, given the offered load $R$, we search for the pair $(s,n)$ that realizes a target delay probability $\e\in(0,1)$.
Relying on the two-fold scaling in \eqref{eq:twofold_scaling_cloud_model}, this under large offered loads $R$ is tantamount to finding the pair $(\b,\g)$ that achieves asymptotic delay probability $\e$. 
In a system without retrials, attaining this target performance boils down to finding a pair $(\beta_\e,\g_\e)$ such that $g_c(\b_\e,\g_\e) = \e$. 
The fixed-point heuristic however tells us that the model with retrials performs slightly worse, namely as if the QoS-parameters were $(\beta_\e-\alpha,\g-\alpha/\sqrt{\kappa})$ for $\alpha>0$. 
Henceforth, to attain the target delay probability $\e$ in the limit with retrials, larger QoS parameters are required. 
To be precise, $\b^*_\e = \b_\e+\alpha$ and $\g^*_\e = \g_\e+\alpha/\sqrt{\kappa}$, with $\alpha$ satisfying $\alpha = f_c(\beta^*-\alpha, \gamma^*_\e -\alpha/\sqrt{\kappa}) = f_c(\beta_\e,\gamma_\e)$. 
Finally, we substitute $\beta^*_\e$ and $\gamma^*_\e$ in the scaling \eqref{eq:twofold_scaling_cloud_model} to obtain capacity levels $s$ and $n$. 
Altogether, this yields the QED dimensioning procedure in Algorithm \ref{alg:cloud_stationary}, in which $[\cdot]$ denotes the integer rounding operator.

\begin{algorithm}
\hspace{1cm}\rule{10cm}{1pt}\\
\hspace{1.1cm}\KwIn{Offered load $R$\\
\hspace{1.1cm}Expected time spent in seconds stage $1/\kappa$\\
\hspace{1.1cm}Target delay probability  $\e\in(0,1)$}
\hspace{1.1cm}\KwOut{Capacity levels $s$ and $n$.}
\vspace{-1mm}
\hspace{1cm}\rule{10cm}{0.5pt}\\
\vspace{-1mm}
\begin{enumerate}
\item[] \hspace{0.5cm} 1. Compute $(\beta_\e,\g_\e)$ such that $g_c(\beta_\e,\g_\e)=\e$.
\item[] \hspace{0.5cm} 2. Set $\b _\e^* = \b_{\e} + f_c(\beta_{\e},\g_\e)$ and $\g_\e^* = \g_\e + f_c(\beta_{\e},\g_\e)/\sqrt{\kappa}$.
\item[] \hspace{0.5cm} 3. Return $s = \lceil R + \beta^*_\e\sqrt{R} \rceil$ and $n = [s +R/\kappa + \gamma^*_\e \sqrt{R/\kappa}]$.
\end{enumerate}
\vspace{-3 mm}
\hspace{1cm}\rule{10cm}{1pt}\\
\vspace{2 mm}
\caption{Stationary dimensioning for cloud model with retrials.}
\label{alg:cloud_stationary}
\end{algorithm}

\begin{table}
\centering
\small
\begin{subtable}{0.99\textwidth}
\centering
\begin{tabular}{|r|rr|r||rr|r||rr|r|}
\cline{2-10}\multicolumn{1}{r|}{} & \multicolumn{3}{c||}{$\e = 0.10$} & \multicolumn{3}{c||}{$\e = 0.25$} & \multicolumn{3}{c|}{$\e = 0.40$} \bigstrut[t]\\
\cline{2-10}\multicolumn{1}{r|}{} & \multicolumn{3}{c||}{$(\b^*_\e,\g^*_\e) = (1.17,0.35)$} & \multicolumn{3}{c||}{$(\b^*_\e,\g^*_\e) = (0.73,0.57)$} & \multicolumn{3}{c|}{$(\b^*_\e,\g^*_\e) = (0.48,0.78)$} \bigstrut\\
\hline
$R$     & $s$     & $n$     & $\P^c_r({\rm del})$ & $s$     & $n$     & $\P^c_r({\rm del})$ & $s$     & $n$     & $\P^c_r({\rm del})$ \bigstrut\\
\hline
10    & 14    & 25    & 0.1231 & 13    & 25    & 0.2268 & 12    & 24    & 0.3717 \bigstrut[t]\\
50    & 59    & 111   & 0.0969 & 56    & 110   & 0.2289 & 54    & 110   & 0.3844 \\
100   & 112   & 215   & 0.1069 & 108   & 214   & 0.2426 & 105   & 213   & 0.4148 \\
500   & 527   & 1035  & 0.0994 & 517   & 1030  & 0.2486 & 511   & 1028  & 0.3996 \\
1000  & 1038  & 2049  & 0.0977 & 1024  & 2042  & 0.2442 & 1016  & 2041  & 0.3925 \bigstrut[b]\\
\hline
\end{tabular}%
\caption{$\kappa=1$}
\end{subtable}

\vspace{ 5mm }

\begin{subtable}{0.99\textwidth}
\centering
\begin{tabular}{|r|rr|r||rr|r||rr|r|}
\cline{2-10}\multicolumn{1}{r|}{} & \multicolumn{3}{c||}{$\e = 0.10$} & \multicolumn{3}{c||}{$\e = 0.25$} & \multicolumn{3}{c|}{$\e = 0.40$} \bigstrut[t]\\
\cline{2-10}
\multicolumn{1}{r|}{} & \multicolumn{3}{c||}{$(\b^*_\e,\g^*_\e) = (1.34,0.50)$} & \multicolumn{3}{c||}{$(\b^*_\e,\g^*_\e) = (0.87,0.65)$} & \multicolumn{3}{c|}{$(\b^*_\e,\g^*_\e) = (0.59,0.79)$} \bigstrut\\
\hline
$R$     & $s$     & $n$     & $\P^c_r({\rm del})$ & $s$     & $n$     & $\P^c_r({\rm del})$ & $s$     & $n$     & $\P^c_r({\rm del})$ \bigstrut\\
\hline
10    & 15    & 69    & 0.0893 & 13    & 68    & 0.2628 & 12    & 68    & 0.4238 \bigstrut[t]\\
50    & 60    & 318   & 0.1007 & 57    & 317   & 0.2210 & 55    & 318   & 0.3560 \\
100   & 114   & 625   & 0.0989 & 109   & 624   & 0.2504 & 106   & 624   & 0.4099 \\
500   & 530   & 3055  & 0.1051 & 520   & 3053  & 0.2449 & 514   & 3054  & 0.3857 \\
1000  & 1043  & 6078  & 0.0986 & 1028  & 6074  & 0.2459 & 1019  & 6075  & 0.3981 \bigstrut[b]\\
\hline
\end{tabular}%
\caption{$\kappa=0.2$}
\end{subtable}

\vspace{ 5mm }

\begin{subtable}{0.99\textwidth}
\centering
\begin{tabular}{|r|rr|r||rr|r||rr|r|}
\cline{2-10}\multicolumn{1}{r|}{} & \multicolumn{3}{c||}{$\e = 0.10$} & \multicolumn{3}{c||}{$\e = 0.25$} & \multicolumn{3}{c|}{$\e = 0.40$} \bigstrut[t]\\
\cline{2-10}
\multicolumn{1}{r|}{} & \multicolumn{3}{c||}{$(\b^*_\e,\g^*_\e) = (1.41,0.68)$} & \multicolumn{3}{c||}{$(\b^*_\e,\g^*_\e) = (0.93,0.74)$} & \multicolumn{3}{c|}{$(\b^*_\e,\g^*_\e) = (0.59,0.79)$} \bigstrut\\
\hline
$R$     & $s$     & $n$     & $\P^c_r({\rm del})$ & $s$     & $n$     & $\P^c_r({\rm del})$ & $s$     & $n$     & $\P^c_r({\rm del})$ \bigstrut\\
\hline
10    & 15    & 530   & 0.0996 & 13    & 530   & 0.2814 & 13    & 531   & 0.2816 \bigstrut[t]\\
50    & 60    & 2594  & 0.1146 & 57    & 2594  & 0.2416 & 55    & 2595  & 0.3805 \\
100   & 115   & 5163  & 0.0933 & 110   & 5163  & 0.2323 & 107   & 5163  & 0.3794 \\
500   & 532   & 25640 & 0.1001 & 521   & 25638 & 0.2516 & 515   & 25641 & 0.3873 \\
1000  & 1045  & 51198 & 0.1018 & 1030  & 51196 & 0.2437 & 1021  & 51199 & 0.3900 \bigstrut[b]\\
\hline
\end{tabular}%
\caption{$\kappa=0.02$}
\end{subtable}
\caption{\normalfont Results of the stationary dimensioning algorithm for $\varepsilon = 0.1,\, 0.25$ and $0.4$.}
\label{tab:cloud_dimensioning}
\end{table}

In Table \ref{tab:cloud_dimensioning} we performed this stationary dimensioning procedure for $\kappa = 1$, 0.2 and 0.02, and $\varepsilon = 0.1$, 0.25 and 0.4, and increasing offered loads $R$, and used simulation to obtain the actual delay probabilities. 
We immediately see that the procedure yields remarkably good results, that are very close to the target delay probabilities. 
\subsection{Time-varying dimensioning}

We next discuss how the parameters $s$ and $n$ can be adjusted in time-varying environments where the offered load $R(t)$ is a function of time. 
For this we use the mean-offered-load (MOL) method, which was developed in \cite{Jennings1996} to approximate and dimension the $M_t/G/s$ system by establishing a relation with the analytically tractable $M_t/G/\infty$ system. 
An underlying assumption of the MOL method is that a well-capacitated multi-server queue delays only a small portion of users and only for short periods. 
Therefore, the system can be approximated by an infinite-server system. 
The MOL approximation \cite{Jennings1996} combines the desirable QoS properties rendered by the QED regime with the analytic tractability of the $M/G/\infty$ queue, see \cite{Eick1993}, to establish a dynamic algorithm for choosing $s(t)$ that stabilizes the system behavior at some QoS-target.

To understand why the MOL approximation is likely to be accurate for the systems in this chapter, observe that under the QED scalings, the blocking probability vanishes asymptotically and hence the main assertion on which the MOL approximation is built continues to hold. 
Following the line of thought in \cite{Jennings1996}, we consider the number of users in a system with $s = n = \infty$ to obtain $R_1(t) = \E[R(t-S)]\E [S]$, where $S$ is the service requirement per customer taken to be unit exponentially distributed. Then,
\begin{equation}
R_1(t) = \int_0^\infty e^{-u}R(t-u) \,\dd u.
\label{eq:R1}
\end{equation}
Note that this transformation typically shifts and levels peaks in workload ahead in time with respect to those in $R(t)$. As the time-varying counterparts of $s$ in \eqref{eq:twofold_scaling_cloud_model}, we then get 
\begin{align}
\label{eq:sTimeVarying}
s(t) &= R_1(t) + \b\sqrt{R_1(t)}.
\end{align}
Secondly, the number of customers present in the system strongly depends on the number of customers in the second phase of the system, especially since $\kappa \ll \mu$. 
Therefore, we moreover need an approximation for the workload offered to the second queue as a function of $t$, which we denote by $R_2(t)$. 
Continuing the reasoning of MOL, we argue that $R_2(t)$ is equal to the output process of the first queue $R_1(t)$. 
Then,
\begin{align}
R_2(t) = \E[R_1(t-S_2)]\E[S_2] = \int_0^\infty\int_0^\infty  e^{-u-\kappa v} R(t-u-v) \dd u\, \dd v \label{eq:R2}
\end{align}
and the natural shape of $n(t)$ becomes
\begin{equation}
n(t) = s(t) + R_2(t) + \gamma\sqrt{R_2(t)}.
\end{equation}
Combining the above ingredients then leads to Algorithm \ref{alg:cloud_timevarying}.
\begin{algorithm}
\hspace{1cm}\rule{10cm}{1pt}\\
\hspace{1.1cm}\KwIn{Offered load function $R(t)$\\
\hspace{1.1cm}Expected time spent in second stage $1/\kappa$\\
\hspace{1.1cm}Target delay probability  $\e\in(0,1)$}
\hspace{1cm}\KwOut{Capacity levels $s(t)$ and $n(t)$ achieving $P_r({\rm delay}) = \e$.}
\vspace{-1mm}
\hspace{1cm}\rule{10cm}{0.5pt}\\
\vspace{-1mm}
\begin{enumerate}[noitemsep]
\item[] \hspace{0.5cm} 1. Compute $\beta^*_\e$ and $\g_\e^*$ according to Algorithm \ref{alg:cloud_stationary}.
\item[] \hspace{0.5cm} 2. Compute $R_1(t)$ and $R_2(t)$ as in \eqref{eq:R1} and \eqref{eq:R2}.
\item[] \hspace{0.5cm} 3. Return
\begin{align*}
\hspace{-0.5cm} s(t) &= \Big\lceil R_1(t)+\b_\e^*\sqrt{R_1(t)}\Big\rceil , \nonumber \\
\hspace{-0.5cm} n(t) &= \Big[s(t) + R_2(t) + \g_\e^*\sqrt{R_2(t)}\Big].
\end{align*}
\end{enumerate}
\vspace{-3 mm}
\hspace{1cm}\rule{10cm}{1pt}\\
\vspace{2 mm}
\caption{Time-varying dimensioning algorithm for cloud model with retrials.}
\label{alg:cloud_timevarying}
\end{algorithm}
Observe that if the service times at the infinite-server queue are relatively short compared with the rate of change of the load function, we have $R_2(t)\approx R_1(t)/\kappa$, so that $n(t)$ as in Algorithm \ref{alg:cloud_timevarying} shows resemblance with \ref{eq:twofold_scaling_cloud_model}.

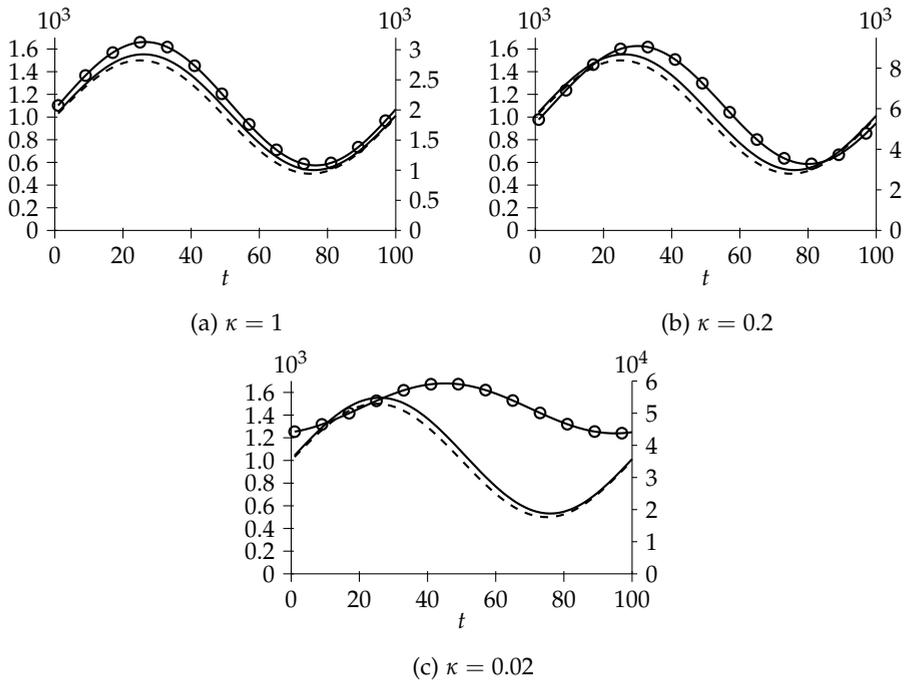
\begin{figure}
\centering
\begin{subfigure}{0.48\textwidth}
\begin{tikzpicture}[y=0.0015cm, x=0.045cm]
 	\small
	\draw (0,0) -- coordinate (x axis mid) (100,0);
    \draw (0,0) -- coordinate (y axis mid) (0,1700) node[above] {$10^3$};
	\draw (100,0) -- coordinate (y2 axis mid) (100,1700) node[above] {$10^3$};
    	\foreach \x in {0,20,...,100}
     		\draw (\x,1pt) -- (\x,-3pt)
			node[anchor=north] {\x};
    	\foreach \y in {0,0.2,0.4,0.6,0.8,1.0,1.2,1.4,1.6}
     		\draw (1pt,{1000*\y}) -- (-3pt,{1000*\y})
     			node[anchor=east] {\y};
     			
     	\foreach \z in {0,0.5,...,3.2}
     		\draw (100,{\z*531.25}) -- (101.5,{\z*531.25})
     			node[anchor=west] {{\z}};

	\node[below=0.4cm] at (x axis mid) {$t$};
	
	\draw[thick,dashed] plot
			file {Chapter_4/tikz/R_timevarying.dat};
	\draw[thick] plot
			file {Chapter_4/tikz/s_timevarying.dat};
	\draw[thick] plot[mark=o,mark repeat=8]
			file {Chapter_4/tikz/n_th1_timevarying.dat};
\end{tikzpicture}
\caption{$\kappa=1$}
\end{subfigure}
\begin{subfigure}{0.48\textwidth}
\begin{tikzpicture}[y=0.0015cm, x=0.045cm]
 	\small
	\draw (0,0) -- coordinate (x axis mid) (100,0);
    \draw (0,0) -- coordinate (y axis mid) (0,1700) node[above] {$10^3$};
	\draw (100,0) -- coordinate (y2 axis mid) (100,1700) node[above] {$10^3$};
    	\foreach \x in {0,20,...,100}
     		\draw (\x,1pt) -- (\x,-3pt)
			node[anchor=north] {\x};
    	\foreach \y in {0,0.2,0.4,0.6,0.8,1.0,1.2,1.4,1.6}
     		\draw (1pt,{1000*\y}) -- (-3pt,{1000*\y})
     			node[anchor=east] {\y};
     			
     	\foreach \z in {0,2,...,9}
     		\draw (100,{\z*178.947}) -- (101.5,{\z*178.947})
     			node[anchor=west] {{\z}};

	\node[below=0.4cm] at (x axis mid) {$t$};
	
	\draw[thick,dashed] plot
			file {Chapter_4/tikz/R_timevarying.dat};
	\draw[thick] plot
			file {Chapter_4/tikz/s_timevarying.dat};
	\draw[thick] plot[mark=o,mark repeat=8]
			file {Chapter_4/tikz/n_th02_timevarying.dat};
\end{tikzpicture}
\caption{$\kappa=0.2$}
\end{subfigure}
\begin{subfigure}{0.48\textwidth}
\begin{tikzpicture}[y=0.0015cm, x=0.045cm]
 	\small
	\draw (0,0) -- coordinate (x axis mid) (100,0);
    \draw (0,0) -- coordinate (y axis mid) (0,1700) node[above] {$10^3$};
	\draw (100,0) -- coordinate (y2 axis mid) (100,1700) node[above] {$10^4$};
    	\foreach \x in {0,20,...,100}
     		\draw (\x,1pt) -- (\x,-3pt)
			node[anchor=north] {\x};
    	\foreach \y in {0,0.2,0.4,0.6,0.8,1.0,1.2,1.4,1.6}
     		\draw (1pt,{1000*\y}) -- (-3pt,{1000*\y})
     			node[anchor=east] {\y};
     			
     	\foreach \z in {0,1,...,6}
     		\draw (100,{\z*283.333}) -- (101.5,{\z*283.333})
     			node[anchor=west] {{\z}};

	\node[below=0.4cm] at (x axis mid) {$t$};
	
	\draw[thick,dashed] plot
			file {Chapter_4/tikz/R_timevarying.dat};
	\draw[thick] plot
			file {Chapter_4/tikz/s_timevarying.dat};
	\draw[thick] plot[mark=o,mark repeat=8]
			file {Chapter_4/tikz/n_th002_timevarying.dat};

\end{tikzpicture}
\caption{$\kappa=0.02$}
\end{subfigure}
\caption{Arrival rate function $R(t)$ (dashed) and staffing functions $s(t)$ (solid) and $n(t)$ (o) for different values of $\kappa$. The left vertical axis refers to $R(t)$ and $s(t)$, where the right axis refers to $n(t)$.}
\label{fig:staffingCurves}
\end{figure}
\begin{figure}
\centering
\begin{subfigure}{0.48\textwidth}
\begin{tikzpicture}[y=6cm, x=0.0445cm]
\small
 	
	\draw (0,0) -- coordinate (x axis mid) (100,0);
    \draw (0,0) -- coordinate (y axis mid) (0,0.52) node[above]{$P_r({\rm delay})$};
    	\foreach \x in {0,20,...,100}
     		\draw (\x,1pt) -- (\x,-3pt)
			node[anchor=north] {\x};
    	\foreach \y in {0,0.1,0.2,0.3,0.4,0.5}
     		\draw (1pt,\y) -- (-3pt,\y)
     			node[anchor=east] {\y};
	\node[below=0.4cm] at (x axis mid) {$t$};
	
	\draw[dashed] (0,0.1) -- (100,0.1);
	\draw[dashed] (0,0.25) -- (100,0.25);
	\draw[dashed] (0,0.4) -- (100,0.4);	
	
	\draw[thick,color=col1] plot
			file {Chapter_4/tikz/time_th1_e01.dat};
	\draw[thick,color=col4] plot
			file {Chapter_4/tikz/time_th1_e025.dat};
	\draw[thick,color=col5] plot
			file {Chapter_4/tikz/time_th1_e04.dat};
\end{tikzpicture}
\caption{$\kappa = 1$}
\end{subfigure}
\begin{subfigure}{0.48\textwidth}
\begin{tikzpicture}[y=6cm, x=0.0445cm]
\small

	\draw (0,0) -- coordinate (x axis mid) (100,0);
    \draw (0,0) -- coordinate (y axis mid) (0,0.52) node[above]{$P_r({\rm delay})$};
    	\foreach \x in {0,20,...,100}
     		\draw (\x,1pt) -- (\x,-3pt)
			node[anchor=north] {\x};
    	\foreach \y in {0,0.1,0.2,0.3,0.4,0.5}
     		\draw (1pt,\y) -- (-3pt,\y)
     			node[anchor=east] {\y};
	\node[below=0.4cm] at (x axis mid) {$t$};
	
	\draw[dashed] (0,0.1) -- (100,0.1);
	\draw[dashed] (0,0.25) -- (100,0.25);
	\draw[dashed] (0,0.4) -- (100,0.4);	
	
	\draw[thick,color=col1] plot
			file {Chapter_4/tikz/time_th02_e01.dat};
	\draw[thick,color=col4] plot
			file {Chapter_4/tikz/time_th02_e025.dat};
	\draw[thick,color=col5] plot
			file {Chapter_4/tikz/time_th02_e04.dat};
\end{tikzpicture}
\caption{$\kappa =0.2$}
\end{subfigure}
\begin{subfigure}{0.48\textwidth}
\begin{tikzpicture}[y=6cm, x=0.0445cm]
\small

	\draw (0,0) -- coordinate (x axis mid) (100,0);
    \draw (0,0) -- coordinate (y axis mid) (0,0.52) node[above]{$P_r({\rm delay})$};
    	\foreach \x in {0,20,...,100}
     		\draw (\x,1pt) -- (\x,-3pt)
			node[anchor=north] {\x};
    	\foreach \y in {0,0.1,0.2,0.3,0.4,0.5}
     		\draw (1pt,\y) -- (-3pt,\y)
     			node[anchor=east] {\y};
	\node[below=0.4cm] at (x axis mid) {$t$};

	\draw[dashed] (0,0.1) -- (100,0.1);
	\draw[dashed] (0,0.25) -- (100,0.25);
	\draw[dashed] (0,0.4) -- (100,0.4);	
	
	\draw[thick,color=col1] plot
			file {Chapter_4/tikz/time_th002_e01.dat};
	\draw[thick,color=col4] plot
			file {Chapter_4/tikz/time_th002_e025.dat};
	\draw[thick,color=col5] plot
			file {Chapter_4/tikz/time_th002_e04.dat};
\end{tikzpicture}
\caption{$\kappa = 0.02$}
\end{subfigure}
\caption{Simulated time-dependent delay probabilities in the cloud model with $\d = 10^{-2}$, targets $\e = 0.1, \e=0.25$ and $\e=0.4$, and  capacity levels determined by Algorithm \ref{alg:cloud_timevarying}.}
\label{fig:timevarying}
\end{figure}
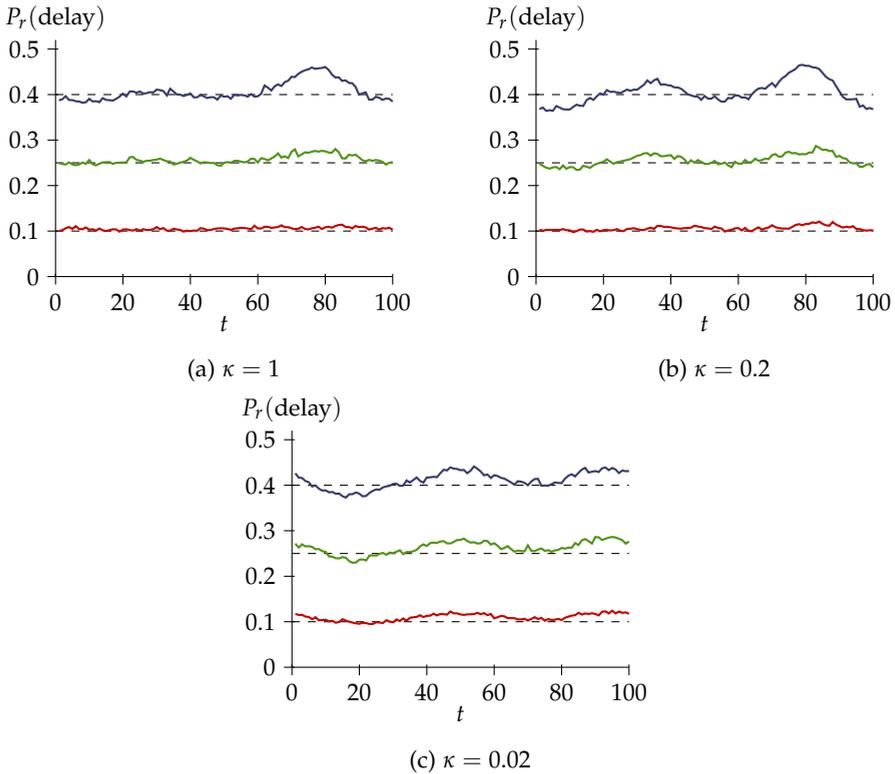

To illustrate Algorithm \ref{alg:cloud_timevarying} we consider the time-varying load
\begin{equation}
R(t) = a + b \,\sin\left(2\pi t / T\right),
\end{equation}
where we set $a=1000$ as the mode, $b=500$ as the amplitude and $T=100$ as the cycle length. 
This system experiences large fluctuations in load volume over the course of one cycle. 
Since $\mu=1$, this implies that one cycle on average consists of 100 service times at the host server queue. 
Due to relatively short service times with respect to the cycle length, the MOL approximation for the number of customers at the first queue is roughly equal to the original load, i.e.~$R_1(t) \approx R(t)$. 
These short services at the first queue compared to the cycle length are typical for cloud systems, in which case the cycle is usually one day.

First, we examine the functions $s(t)$ and $n(t)$ as prescribed by Algorithm \ref{alg:cloud_timevarying} for $\kappa=1,0.2,0.02$ and $\e = 0.25$. 
The resulting values are depicted in Figure~\ref{fig:staffingCurves} together with the arrival rate function. 
Note that $n(t)$ lives on a different scale than $s(t)$, and has its own vertical axis at the right side of the plots.
For small and hence realistic values of $\kappa$, the function $n(t)$ displays a shifted phase compared to the real-time offered load, due to the relatively long service time at the second station. 
The lag can be observed in~\eqref{eq:R2}. 
Hence, while the number of servers $s(t)$ allocated at time $t$ is almost in phase with the arrival rate $R(t)$, $n(t)$ undergoes a shift of its peak capacity somewhat ahead in time.
Observe also that $n(t)$ shows milder fluctuations when $\kappa$ decreases. This can be attributed to the added hedge which is of order $\sqrt{R/\kappa}$. 
Remark that the overcapacity is relatively small.  
This illustrates the economies-of-scale that can be achieved in these large-scale systems.
Next, we simulate the time-dependent process, given the staffing functions depicted in Figure~\ref{fig:staffingCurves}, as well as the staffing functions designed for the target delay probabilities $\e =0.1$ and $\e =0.4$ for the three values of $\kappa$. 
The results of the simulations are depicted in Figure \ref{fig:timevarying}. In all cases, the time-dependent delay probability only mildly fluctuates around the target.
As we increase $\e$, the stabilizing effect of the method weakens somewhat, which for other systems was also observed in \cite{Jennings1996}.

\section{Conclusion}
\label{sec:retrial_conclusion}
In this chapter, we studied the impact of retrying customers in large-scale systems in the QED regime. 
The presence of retrials has a detrimental effect on congestion-related performance, compared to systems in which customers are simply discarded upon blockage/abandonment. 
On the other hand, compared to similar systems without physical size restrictions or customer impatience, the performance gain can be substantial, if retrial times are relatively long compared to the service times. 
Namely, retrials prompt temporary release of pressure from the system by shifting workload ahead in time. 

Through our analysis, we have shown how the performance of large-scale queueing systems facing slow retrials can be approximated by appropriately combining a fixed-point technique with QED scaling.
We showed the remarkable accuracy of this approximation scheme in various retrial settings, that are otherwise intractable to analyze. 
As we discussed in Section \ref{sec:retrial_dimensioning}, our novel asymptotic analysis technique is furthermore a powerful and elegant tool for dimensioning large-scale systems with slow retrials, which is moreover amenable to deal with time-varying demand. \\
\\*

We illustrate a few directions for future research. 
As we explained before, the fixed-point method relies heavily on the premise that the blocking (or abandonment) probability vanishes at rate $1/\sqrt{R}$ in the QED regime, and on the availability of expressions for its limiting behavior.
Since this description likely fits a wide range of queueing models, we henceforth believe that our fixed-point method and the related dimensioning scheme find application beyond the three models we discussed here.

Secondly, in the dimensioning procedure of Section \ref{sec:retrial_dimensioning} we took a constraint satisfaction perspective in which we aimed to achieve a preset target QoS-level. 
As an alternative approach, one could define a cost function to quantify the trade-off between capacity costs and customer dissatisfaction.
Specifically, suppose a cost $c_1$ is associated with each server per unit of time, cost $c_2$ is charged for every waiting customer per time unit, and cost $c_3$ is the penalty for each blocked customer. 
Then in the $M/M/s/n$ queue with retrials in the QED regime, we use that $s = R+\beta\sqrt{R}$, the blocking probability is roughly $f(\beta-\alpha,\gamma)/\sqrt{R}$ and the expected waiting time is approximately $h(\beta-\alpha,\gamma)/\sqrt{R}$ for some function $h$, see \cite{masseywallace}, yielding total operational cost
\[
c_1\left( R + \beta\sqrt{R} \right) 
+ c_2\, R\, \frac{f(\beta-\alpha,\gamma)}{\sqrt{R}} 
+ c_3\, R\, \frac{h(\beta-\alpha,\gamma)}{\sqrt{R}} ,
\]
where $\alpha$ satisfies the fixed-point equation.
Hence, asymptotic dimensioning of the system boils down to finding the parameters $\beta^*$ and $\gamma^*$ that minimize
\[ c_1\beta^* + c_2\,f(\beta^*-\alpha^*,\gamma^*) + c_2\,h(\beta^*-\alpha^*,\gamma^*),\]
with corresponding fixed point $\alpha^*$. 
Solving this optimization problem is not straightforward and a detailed study of this and related asymptotic dimensioning problems is an interesting avenue for future research. 

Last, we remark that even though the fixed-point method works very well for systems with slow retrials, i.e.~$\delta\to 0$, it may also serve as an approximation to systems with short to moderate retrial times. 
In these scenarios, the method is likely to underestimate congestion levels as it ignores dependencies between the primary and retrial stream of arrivals.
In the extreme case that $\delta\to\infty$, that is, blocked customers retry immediately, the customers in the retrial orbit basically form a (random order) queue outside the service facility. 
When the inside of the facility consists of more than one queue, our fixed-point may be used as an heuristic approach to account for the increased workload that builds up outside the facility. 
We explore this heuristic idea in a health care context in the next chapter.

\chapter{Finite-size effects in critically dimensioned emergency departments}

\begin{chapterstart}
Motivated by health care systems with repeated services that have both personnel (nurse/physician) and space (beds) constraints, we study a restricted version of the Erlang-R model. The space restriction policies we account for are blocking or holding in a pre-entrant queue. We develop many-server approximations for the system performance measures when either policy applies, and explore the connection between them.
We show that capacity allocation of both resources should be determined simultaneously, and derive the methodology to determine it explicitly.
 We show that the system dynamics is captured by the fraction of needy time in the network, and that returning patients should be accounted for both in steady-state and time-varying conditions.
 We demonstrate the application of our policies in two case-studies of resource allocation in hospitals.
\end{chapterstart}

\begin{flushright}
Based on\\
\textbf{Finite-size effects in critically dimensioned emergency departments}\\
\textit{Johan van Leeuwaarden, Britt Mathijsen, Fiona Sloothaak \& Galit Yom-Tov}\\
Submitted to \textit{Operations Research}
\end{flushright}
\newpage

\section{Introduction}

In recent years, operations research techniques have received increased interest from the health care community, as they are able to design and improve workflow processes in health care facilities~\cite{Armony2015,Green2008,Bekker2009b,Hall2006,Hall2012}. 
Because these processes are typically stochastic in nature, it is common practice to use queueing theory for performance analyses and workforce planning. 
As a first step towards understanding the processes going on in health care environments, systems are commonly modeled after a single station queue, such as the $M/M/s$ (Erlang-C), $M/M/s/s$ (Erlang-B) or $M/M/s+M$ (Erlang-A) models, and fluid and diffusion approximations are used to provide insights into the process dynamics. 
However, simple single station models often fail to capture the more intricate dynamics of the settings specific to health care contexts. 
Prime examples include the flows of patients in a hospital from one medical ward to another \citep{Armony2015}, within the Emergency Department (ED) between different stages of treatment \citep{Junfei2015}, or between medical facilities \citep{zychlinski2016bedblocking}.  
Queueing networks can capture the dependency between several service stages and several types of resources. 
More specifically, we are interested in the ubiquitous feature, particularly present in health care environments, that patients during their stay in the system might require a specific resource multiple times, e.g.~physicians and nurses who treat patients several times during their stay in the medical wards \citep{Jennings2011} or the ED \citep{YomTov2014}, while multiple resources types are limited (e.g.\ medical staff and beds).
In this chapter, we concentrate on the dynamics within EDs. 

\begin{figure} 
\centering
\begin{subfigure}{0.48\textwidth}
\centering
\begin{tikzpicture}[scale = 0.55]
\footnotesize
\draw [thick, ->] (-2.75,4.5) -- (-1.25,4.5);
\draw [thick] (-1.5,5) -- (0,5) -- (0,4) -- (-1.5,4);
\draw [thick] (-0.25,4) -- (-0.25,5);
\draw [thick] (-0.5,4) -- (-0.5,5);
\draw [thick] (-0.75,4) -- (-0.75,5);
\draw [thick, ->] (0,4.5) -- (2.5,4.5);
\draw [thick] (2.25,5) -- (3.75,5) -- (3.75,4) -- (2.25,4);
\draw [thick] (3,4) -- (3,5);
\draw [thick] (3.25,4) -- (3.25,5);
\draw [thick] (3.5,4) -- (3.5,5); 
\draw [thick] (4.25,4.5) circle [radius=0.5] node {$s$} node[below= 0.3 cm] {needy} node[above=0.25cm] {$\exp(\mu)$};
\draw [thick,->] (4.75,4.5) -- (7.75,4.5);
\draw [thick,->] (5.75,4.5) -- (5.75,2) -- (4,2);
\draw [thick] (3.5,2) circle [radius=0.5] node {$\infty$} node[below=0.3cm] {content} node[above=0.25cm] {$\exp(\delta)$};
\draw [thick,->] (3,2) -- (1.5,2) -- (1.5,4.5);
\node [above] at (-1.5,2.9) {\footnotesize Pois($\lambda$)};
\node [below] at (5.5,2) {\footnotesize $p$};
\node [above] at (6.25,4.5) {\footnotesize $1-p$};
\draw [thick, dashed] (0.5,0.75) rectangle (7,5.9) node[right] {\footnotesize $n$}; 
\node [above] at (-0.75,5) {\footnotesize holding};
\end{tikzpicture}
\caption{Erlang-R model with holding.}
\label{fig:Erlang_R_holding}
\end{subfigure}
\begin{subfigure}{0.48\textwidth}
\begin{tikzpicture}[scale = 0.55]
\footnotesize
\draw [thick, ->] (-1.5,4.5) -- (2.5,4.5);
\draw [thick, ->] (0,4.5) -- (0,2.5) node[below left] {blocked};
\draw [thick] (2.25,5) -- (3.75,5) -- (3.75,4) -- (2.25,4);
\draw [thick] (3,4) -- (3,5);
\draw [thick] (3.25,4) -- (3.25,5);
\draw [thick] (3.5,4) -- (3.5,5); 
\draw [thick] (4.25,4.5) circle [radius=0.5] node {$s$} node[below= 0.3 cm] {needy} node[above=0.25cm] {$\exp(\mu)$};
\draw [thick, ->] (4.75,4.5) -- (7.75,4.5);
\draw [thick,->] (5.75,4.5) -- (5.75,2) -- (4,2);
\draw [thick] (3.5,2) circle [radius=0.5] node {$\infty$} node[below=0.3cm] {content} node[above=0.25cm] {$\exp(\delta)$};
\draw [thick,->] (3,2) -- (1.5,2) -- (1.5,4.5);
\node [above] at (-1.5,3.4) {\footnotesize Pois($\lambda$)};
\node [below] at (5.5,2) {\footnotesize $p$};
\node [above] at (6.25,4.5) {\footnotesize $1-p$};
\draw [thick, dashed] (0.5,0.75) rectangle (7,5.9) node[right] {\footnotesize $n$};
\end{tikzpicture}
\caption{Erlang-R model with blocking.}
\label{fig:Erlang_R_blocking}
\end{subfigure}
\caption{Restricted Erlang-R models with maximally $n$ customers in system.}
\label{fig:Erlang_R_model}
\end{figure}
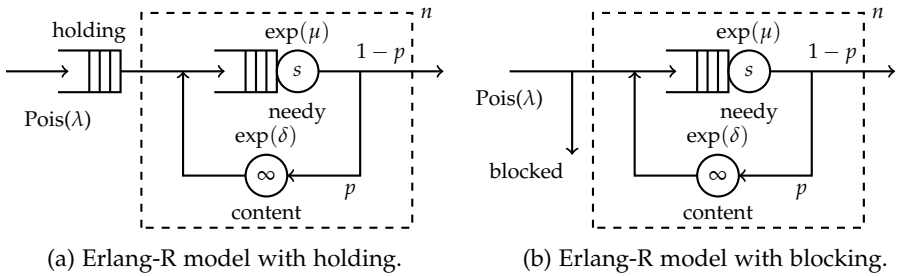

An often ignored yet essential feature of medical facilities concerns the restriction of the number of patients that can reside in the facility simultaneously. 
In Chapter 4, we already observed that finite-size restrictions can have a significant effect on the performance of queueing systems. 
In this chapter, we investigate the influence of such multiple restrictions on the network dynamics and the required staffing policies in the context of an ED. \\
\\*
\noindent
\textbf{The restricted Erlang-R model.}
The canonical model for service networks with returns is the Erlang-R model, introduced by Yom-Tov \& Mandelbaum~\citep{YomTov2014}.
In this open two-station model, customers arrive according to a Poisson process to an $M/M/s$. 
After service completion, the customer with probability $1-p$ leaves the system and with probability $p$ returns to the queue after a random delay. 
This delay is modeled as an infinite-server queue. 
A schematic visualization of the Erlang-R model is depicted in Figure \ref{fig:ErlangR}

in which customers, during their stay in the system receive a random number of services from the same pool of servers. 
Yom-Tov \& Mandelbaum \cite{YomTov2014} showed that such a simple network model can be used to determine staffing in an ED both in stable and time-varying conditions. 
Nevertheless, empirical studies report that some countries, such as the US, use a different operational mode that applies strict restrictions on entering the ED \citep{EDexperiment}. 
In typical US EDs, a patient will not enter the ED until both a bed and a physician are available to treat her. 
Those restrictions can be either physical (beds) restrictions or managerial ones --- for instance by imposing a patient-to-physician ratio.
In this work, we extend the Erlang-R model by enforcing a constraint on the maximum number of available places inside the facility. 
Our model hence incorporates two kinds of resource constraints: servers that provide the actual service and the maximum available places inside the service system. 
Both affect the system in a highly interdependent way. 
The model, presented in Figure \ref{fig:Erlang_R_model}, assumes $s$ servers and a maximum capacity of $n$ concurrent places. 
We assume that patients arrive according to a Poisson process with rate $\lambda$. 
In case a new arrival finds $n$ or more patients already present, we consider two options: either she waits outside the service facility in a holding queue until a vacant space becomes available (Figure \ref{fig:Erlang_R_holding}) or she is blocked (Figure \ref{fig:Erlang_R_blocking}), such as is the case when patients are sent to an alternative facility. 
Once a patient is admitted, she requires assistance from one of the $s$ servers for an exponentially distributed duration with mean $1/\mu$. 
Then, with probability $1-p$, the patient leaves the system or, with probability $p$, returns to service again after an exponentially distributed time with mean $1/\delta$. 
Following Jennings \& de V\'ericourt \cite{Jennings2011} and Yom-Tov  \& Mandelbaum \cite{YomTov2014}, we call patients {\it needy} when they require attention from one of the servers and {\it content} when they are in the delayed return phase. 
In addition, we call patients {\it holding} when they are waiting outside the facility for an available space. We assume that the arrival process, the needy times and content times are mutually independent. 
In the holding queue and the needy queue, we apply the First-Come-First-Served (FCFS) discipline.

As mentioned, we consider two versions of the finite-capacity constraint. 
The first version is called \emph{Erlang-R with holding}, in which patients wait for an available space in the system. 
The second version is called \emph{Erlang-R with blocking}, in which patients meeting a full system are blocked.
Naturally, intermediate scenarios can be constructed in which a proportion of the total arrival volume of patients indeed leaves upon finding a full system, while the rest joins the holding queue. 
While this chapter focuses on the two extreme cases, straightforward adaptations can fit these intermediate scenarios. \\
\\*
\noindent
\textbf{Examples of restricted Erlang-R.}
As noted before, an ED operated in the US can be modeled using a restricted Erlang-R model. Another health care example is medical units (MUs) in a hospital.
Such units specialize in specific types of illnesses (cardiology, oncology, etc.) and have limited resources such as nurses and beds. If the unit is full, new patients are either allocated to an alternative medical unit, i.e.\ blocked, or wait for an available bed. 
Both policies are problematic in terms of quality-of-care, because the personnel in the alternative unit (or the ED) may be less knowledgeable about the patient's medical condition and waiting in the ED was shown to increase mortality.
Moreover, ED waiting may reduce available capacity for treating ED patients \citep{Carmen2016,israelit}, hence endangering both the delayed patient as well as others. Both the number of personnel (nurses and physicians) and the number of beds impact service dynamics and quality-of-care. Research so far looked at the capacity allocation of those resources separately. Green \& Yankovic \cite{GY2011} and Jennings \& de V\'ericourt \cite{Jennings2008} looked at nurse staffing in medical units, while de Bruin et al.~\cite{Bekker2009b} looked into bed allocation. The unified model we suggest enables us to capture the dependency between those two decisions, and its impact on other medical units in the hospital.
At the same time, we capture the two most commonly used modes of operation --- blocking and holding of new patients.  \\
\\*
\noindent
\textbf{Two-fold square-root staffing rule.}
Our main goal is to provide staffing policies for the ED that  high resource utilization, while at the same time maintain good quality-of-care. 
This goal relates to the philosophy of the Quality-and-Efficiency-Driven (QED) regime that is the recurring theme of this thesis. 
In this chapter, we obtain asymptotic results for the Erlang-R model with blocking in the QED regime (Section \ref{sec:QED_limit_block}). 
Following \cite{Jennings2008}, we employ a two-fold QED staffing policy: $s=R_1 +\beta \sqrt{R_1}$  for the number of nurses and $n=R_1/r+\gamma \sqrt{R_1/r}$ for the number of patients in the system (beds), where $\beta$ and $\gamma$ are constants, $R_1$ is the offered load of the servers (nurses) and $r$ is the fraction of time a patient spends in the needy state. 
We establish limiting expressions for performance measures, such as the probability of delay and blocking, in the form of explicit functions that depend solely on $\beta$ and $\gamma$. 
In deriving these limit results, we use the available product-form solution for the stationary distribution.

Likewise, we pursue QED performance for the Erlang-R model with holding.
However, a direct analytic approach is obstructed by the absence of product-form solutions.
We provide two solutions for establishing QED behavior. 
First, we provide stochastic performance bounds that stay meaningful in the QED regime, which demonstrate the non-degenerate behavior of the two-fold scaling in the large-system limit. 
Second, we develop a heuristic method that quantifies the difference between the holding model and the blocking model. 
This method is to a large extent related to the asymptotic approximation method for retrial queues discussed in Chapter 4, in the sense that we approximate the model with holding through the model with blocking, yet with an increased arrival rate. 
The increase in arrival rate turns out to be the solution of a fixed-point equation.
Using our results on the asymptotic behavior of the model with blocking in the QED regime, we then obtain approximative QED performance measures for the model with holding. 
These theoretical findings ultimately yield algorithms for dimensioning and time-varying staffing. \\
\\*
\textbf{Structure of the chapter.}
We first review related literature on the subject of staffing in health care environments in Section \ref{sec:ed_literature}. 
In Section \ref{sec:modeldescription}, we introduce the mathematical models more formally, and deduce preliminary results on their stability conditions and relative performance. 
Section \ref{sec:QED_scaling} describes the scaling regime we use for our asymptotic study of the restricted Erlang-R models, and Sections \ref{sec:QED_limit_block} and \ref{sec:QED_limit_holding} present our main theoretical findings. 
We turn to dimensioning problems in Section \ref{sec:dimensioning}, and show how our asymptotic QED results can be used to make resource allocation decisions in realistic settings. 
Section \ref{sec:analysis_chapter5} is devoted to the numerical and comparative analysis of the restricted Erlang-R models, and also shows how our method can be applied in time-varying environments through a case study.
We summarize our findings and give directions for future research in Section \ref{sec:conclusion}.

\section{Literature review}
\label{sec:ed_literature}
Due to increasing demand and tightening budgets in health care, there is a growing need for efficient workforce management \citep{Green2008}. Personnel (nurse and physician) expenditure is one of the biggest factors in hospital costs \citep{Kazahaya2005}, and inadequate nursing levels have been mentioned as a significant factor in medical errors and ED overcrowding. In order to establish appropriate nursing levels, a staffing policy requires assessment of a wide range of variables, such as differing nurse expertise and patient acuity during the day. Current methods, such as the minimum nurse-to-patient ratios, are often too inflexible to capture those varying conditions. The American Hospital Association (AHA) and others call for dynamic staffing policies that can deal with the complex and evolving nature of health care \citep{AHA2007}.
Workforce management in health care systems has been studied extensively; see \cite{Denton2013,Hall2006,Hall2012} for overviews.
In recent years it has become apparent that queueing models can be helpful in developing staffing and routing recommendations, not just for large-scale service systems, but also for the small and complicated health care systems.

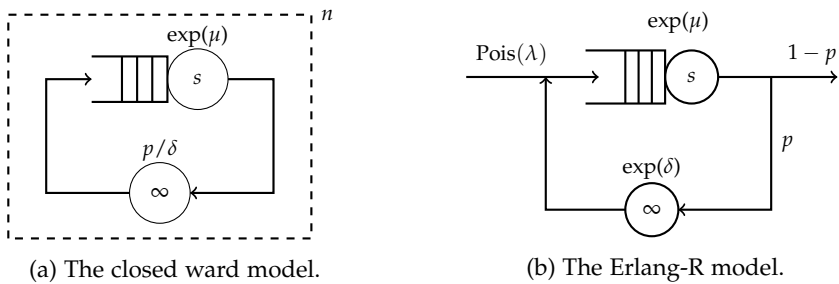
\begin{figure}
\begin{subfigure}{0.48\textwidth}
\centering
\begin{tikzpicture}[scale=1]
\draw [dashed, thick] (-0.5,-0.1) rectangle (3.5,2.85) node[right] {\footnotesize $n$};
\draw [thick,->] (1.1,0.5) -- (0,0.5) -- (0,2) -- (0.6,2);
\draw [thick,->] (2.4,2) -- (3,2) -- (3,0.5) -- (1.9,0.5);
\draw [thick] (0.6,1.7) -- (1.6,1.7) -- (1.6,2.3) -- (0.6,2.3);
\draw [thick] (1,1.7) -- (1,2.3);
\draw [thick] (1.2,1.7) -- (1.2,2.3);
\draw [thick] (1.4,1.7) -- (1.4,2.3);
\draw (2.4,2) -- (3,2) -- (3,0.5) -- (1.9,0.5);
\draw (1.5,0.5) circle [radius=0.4] node[above=0.3cm] {\footnotesize $p/\d$} ;
\draw (2,2) circle [radius=0.4] node[above=0.3cm] {\footnotesize exp($\mu$)};
\node at (1.5,0.5) {\footnotesize $\infty$};
\node at (2,2) {\footnotesize $s$};
\end{tikzpicture}
\caption{The closed ward model.}
\label{fig:Jennings}
\end{subfigure}
\begin{subfigure}{0.48\textwidth}
\centering
\begin{tikzpicture}[scale=0.7]
\draw [thick, ->] (0,4.5) node[above=0.3cm,right] {\footnotesize Pois$(\l)$} -- (2.5,4.5);
\draw [thick] (2.25,5) -- (3.75,5) -- (3.75,4) -- (2.25,4);
\draw [thick] (3,4) -- (3,5);
\draw [thick] (3.25,4) -- (3.25,5);
\draw [thick] (3.5,4) -- (3.5,5);
\draw [thick] (4.25,4.5) circle [radius=0.5];
\draw [thick, ->] (4.75,4.5) -- node[above=0.3cm,right] {\footnotesize $1-p$} (7,4.5);
\draw [thick,->] (5.75,4.5) -- node[right] {\footnotesize $p$} (5.75,2) -- (4,2);
\draw [thick] (3.5,2) circle [radius=0.5];
\draw [thick,->] (3,2) -- (1.5,2) -- (1.5,4.5);
\node at (4.25,4.5) {\footnotesize $s$};
\node at (3.5,2) {\footnotesize $\infty$};
\node [above] at (4,5.15) {\footnotesize exp($\mu$)};
\node [above] at (3.5,2.4) {\footnotesize exp($\delta$)};
\end{tikzpicture}
\caption{The Erlang-R model.}
\label{fig:ErlangR}
\end{subfigure}
\caption{Related queueing models.}
\end{figure}

The first to try such an approach through queueing models were Green et al.~\cite{Green2006,Green2008}, who used the single station stationary Erlang-C model to set staffing levels in EDs and panel sizes for clinics. Using a similar approach, Bekker \& de Bruin~\cite{Bekker2009a} used the Erlang-B model to determine bed allocation for medical wards. 
The first to observe the significant impact of interrupted services in a health care setting were Jennings \& de V\'ericourt \cite{Jennings2008,Jennings2011}. Motivated by the need to set nurse-to-patient ratios for internal wards, they considered a closed queueing system with $s$ nurses and $n$ beds. This is essentially the Erlang-C model with the additional restriction that a finite population of the $n$ patients requires care. In their model, all beds are always occupied, and patients alternate between two phases: the needy phase where patients require service of a nurse and the content phase where they do not; see Figure \ref{fig:Jennings}. The system dynamics of the restricted Erlang-R model are equivalent to those of the closed ward model of \cite{Jennings2008} if the holding queue would never be empty.

Campello et al.~\cite{Campello2016} analyzed a similar operational decision, referred to as ED case management, which determines the maximal number of patients a physician should handle in parallel. They also used queueing networks and analyzed the stationary distribution. Note that in practice such a decision is not only affected by operational measurements such as waiting times, but also by psychological constraints that limit physician capability to manage multiple tasks (patients) in parallel.
KC \cite{diwas} provided empirical evidence that physicians should not treat more than 6-7 patients at the same time. Therefore, many hospitals in the US restrict entrance to EDs even if beds are available if physicians are overloaded. 
We too consider such constraints, and analyze their impact on performance. We take a different approach than \cite{Campello2016}; instead of analyzing numerically steady-state distributions, we develop many-server approximations that can produce insight into the system dynamics, and can be incorporated into time-varying staffing procedures; see Section  \ref{sec:case_study}.  


The model in~\cite{Jennings2008,Jennings2011} was developed for modeling internal dynamics within an internal ward. However, in the ED, beds are not constantly occupied and the utilization level depends on the flow of patients that arrive from outside the system.
Yom-Tov \& Mandelbaum \cite{YomTov2014} highlight the interrupted services while accounting for the transient nature of patient's arrival process, and introduced the Erlang-R model as a model for an ED. The Erlang-R model is an open two-station queueing network that has the same layout as the restricted Erlang-R model, except that all patients find a bed available upon arrival, see Figure \ref{fig:ErlangR}. In both models patients experience the interrupted services, but the Erlang-R model has no further restrictions on the bed capacity, hence neglecting the finite-size effects. Yom-Tov \& Mandelbaum \cite{YomTov2014} showed, using a simulator tailored to an Israeli ED, that the complicated small ED dynamics can be captured using the relatively simple Erlang-R model, and hence, its recommendations can be implemented in ED workforce management.
Although the feature of interrupted services is present in many systems, it is particularly important for modeling EDs, because the duration of the interruption is typically much longer than the time patients require care from a nurse. This explains why the Erlang-R model is considered to be the canonical model for EDs. The restricted Erlang-R model with holding/blocking thus extends the Erlang-R model with finite-size constraints which, like interrupted services, are expected to have a decisive impact on performance.

\section{Models and performance measures}
\label{sec:modeldescription}

\subsection{Three-dimensional Markov process}
\label{sec:Markov_process}

Since in the restricted Erlang-R model described above the arrival process is taken Poisson, and all service and content times are assumed independent and exponential, the system can be characterized in terms of a Markov process.
Let $Q(t) = (H(t),Q_1(t),Q_2(t))$ represent the number of patients in the \emph{holding}, \emph{needy} and \emph{content} state at time $t$, respectively. 
In both variants, $n$ is the maximum number of patients admitted to system, we have $Q_1(t)+ Q_2(t)\leq n$ for all $t\geq 0$. 
Due to the absence of holding patients in the Erlang-R model with blocking,  $H(t)=0$ is enforced in this case, whereas $H(t)$ has unbounded support in the model with holding. 
This distinction requires us to explore the stationary distribution of the two variants separately. 
Before doing so, we introduce some additional notation. 
We define 
\begin{equation}
 R_1 := \frac{\l}{(1-p)\mu}, \qquad R_2 := \frac{p\l}{(1-p)\d}, 
 \label{eq:R1_R2}
 \end{equation}
where $R_1$ and $R_2$ can be interpreted as the offered workload brought towards the needy queue and the content (infinite-server) queue, respectively. 
Furthermore, we define 
\begin{equation}
r:= \frac{\d}{\d+p\m},
\label{eq:delta}
\end{equation} 
which is the fraction of time a patient spends in the needy state (in case she experienced no wait during her sojourn). \\
\\*

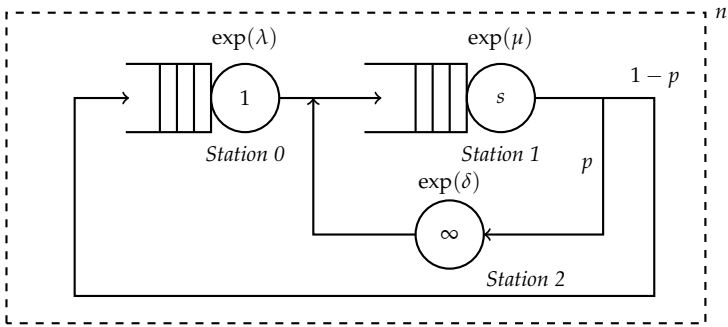
\begin{figure}
\centering
\begin{tikzpicture}[scale = 0.9]
\draw [thick] (-1.25,5) -- (0,5) -- (0,4) -- (-1.25,4);
\draw [thick] (0.5,4.5) circle [radius = 0.5] node {\footnotesize 1} node[above=0.5cm] {\footnotesize exp$(\l)$} 
	node[below =0.5cm] {\footnotesize \textit{Station 0}} ;
\draw [thick] (-0.25,4) -- (-0.25,5);
\draw [thick] (-0.5,4) -- (-0.5,5);
\draw [thick] (-0.75,4) -- (-0.75,5);
\draw [thick, ->] (1,4.5) -- (2.5,4.5);
\draw [thick] (2.25,5) -- (3.75,5) -- (3.75,4) -- (2.25,4);
\draw [thick] (3,4) -- (3,5);
\draw [thick] (3.25,4) -- (3.25,5);
\draw [thick] (3.5,4) -- (3.5,5);
\draw [thick] (4.25,4.5) circle [radius=0.5] node {\footnotesize $s$} node[above=0.5cm] {\footnotesize exp$(\mu)$}
node[below = 0.5cm] {\footnotesize \textit{Station 1}} ;
\draw [thick, ->] (4.75,4.5) -- node[right=0.8cm,above] {\footnotesize $1-p$} (6.5,4.5) -- (6.5,1.6) -- (-2,1.6) -- (-2,4.5) -- (-1.2,4.5);
\draw [thick,->] (5.75,4.5) -- node[left] {\footnotesize $p$} (5.75,2.5) -- (4,2.5);
\draw [thick] (3.5,2.5) circle [radius=0.5] node {\footnotesize $\infty$} node[above=0.4cm] {\footnotesize exp$(\d)$}
node[below right = 0.35cm] {\footnotesize \textit{Station 2}} ;;
\draw [thick,->] (3,2.5) -- (1.5,2.5) -- (1.5,4.5);

\draw [thick, dashed] (-3,1.2) rectangle (7.25,5.75) node[right] {\footnotesize $n$};
\end{tikzpicture}
\caption{The Erlang-R model with blocking viewed as a closed Jackson network.}
\label{fig:ErlangR_blocking}
\end{figure}

\noindent
\textbf{Erlang-R model with blocking.}
In case of the blocking model, $Q(t)$ reduces to a finite-state Markov process $Q(t) = (Q_1(t),Q_2(t))$, where $Q_1(t)+Q_2(t)\leq n$ for all $t\geq 0$. 
In fact, this is equivalent to the closed Jackson network depicted in Figure \ref{fig:ErlangR_blocking} with finite population $n$. 
Station 1 in Figure \ref{fig:ErlangR_blocking} is an $M/M/s$ queue with service rate $\mu$, modeling the number of needy patients $Q_1(t)$. 
Station 2 models the number of content patients $Q_2(t)$, and can therefore be represented as an infinite-server queue with service rate $\d$. 
A patient can enter the unit only if $Q_1(t)+Q_2(t)<n$. 
Station 0---a single-server queue---moderates this as it only produces output at rate $\l$ in case its queue length is positive, i.e.\ if $n-Q_1(t)-Q_2(t)>0$.

Observe that because patients finding a full network are blocked, the number of patients in the system cannot grow beyond $n$. 
Hence, the system is stable for all parameter settings, and hence a steady-state distribution exists. Moreover, the simplification of the model with blocking allows us to express the steady-state distribution of the system in explicit product-form. 
Let $\pi_b(j,k)$ denote the steady-state probabilities of having $j$ needy and $k$ content patients in the system. Then, 
\begin{equation}\label{eq:pih(i,j)}
\pi_b(j,k) = \left\{
\begin{array}{ll}
\pi_0\,\frac{1}{\kappa(j)}\,\frac{1}{k!}\cdot R_1^j\cdot R_2^k, & ~~~\text{if }j+k \leq n,\\
0, & ~~~\text{else,}
\end{array}\right.
\end{equation}
where 
\begin{equation*}
\kappa(j) := \left\{
\begin{array}{ll}
j! , & ~~\text{if }j \leq s,\\
s!\, s^{j-s}, &~~ \text{else,}
\end{array}\right.
\end{equation*}
and $
\pi_0^{-1} = \sum_{j+k\leq n} \frac{1}{\kappa(j)}\,\frac{1}{k!}\cdot R_1^j\cdot R_2^k$.\\
\\*
\textbf{Erlang-R model with holding.}
\label{ref:modelsoft}
The Erlang-R model with holding does not lead to a Jackson network with an elegant product-form solution for the steady-state distribution, because the holding queue cannot be modeled as a station that is independent from the other queues in the system.
However, we are able to describe the system as a two-dimensional Markov process without loss of information. 
To see this, define $N:= \{N(t)\}_{t\geq 0}$ with $N(t) := H(t)+Q_1(t) + Q_2(t)$, the total number of patients in the system (including the holding queue).  
Using the restriction $Q_1(t)+Q_2(t) \leq n$ together with the fact that no bed is left vacant if a patient is waiting in the holding queue, this yields
\begin{equation*}
H(t) = \left(N(t) - n\right)^+, \quad t\geq 0,
\end{equation*}
where $(\cdot)^+ := \max\{0,\cdot\}$.
For the same reason, $Q_2(t) = N(t) - Q_1(t)$ if $H(t)=0$, and $Q_2(t) = n-Q_1(t)$ otherwise. 
In other words,
\begin{equation*}
Q_2(t) = \min\{N(t),n\} - Q_1(t), \quad t \geq 0.
\end{equation*}
Therefore, we can express the state of all three queues in the Erlang-R model with holding using a two-dimensional Markov process $X:= \{X(t)\}_{t\geq 0}$, where 
\begin{equation*}
X(t) :=\left( N(t), Q_1(t) \right).
\end{equation*}
The process $X$ lives on the semi-infinite strip
 \begin{equation*}
X(t) \in \left\{\,(i,j)\, |\, j \leq \min\{i,n\}, i\in \mathbb{N}_0, j \in \{0,1,\ldots,n\}\, \right\},
\end{equation*}
and belongs to the class of Quasi-Birth-Death (QBD) processes.
The reader is referred to Appendix~\ref{app:QBDdescription} for a detailed description of this process, in terms of its transition diagram and generator matrix.

Contrary to the model with blocking, the system with holding \emph{can} become unstable in case capacity is insufficient to satisfy patient demand. 

\begin{proposition}\label{prop:StabilityCondition}
The Erlang-R model with holding is stable if and only if 
\begin{equation}
\frac{\lambda}{(1-p)\mu s} < \frac{ \sum_{i=0}^s \frac{i}{s}\binom{n}{i} \left(\frac{\d}{p\mu}\right)^i + \sum_{i=s+1}^n \binom{n}{i} \frac{i!}{s!} s^{s-i} \left(\frac{\d}{p\m}\right)^i}
{ \sum_{i=0}^s \binom{n}{i} \left(\frac{\d}{p\mu}\right)^i + \sum_{i=s+1}^n \binom{n}{i} \frac{i!}{s!} s^{s-i} \left(\frac{\d}{p\m}\right)^i}
=: \rho_{\max}(s,n).
\label{eq:StabilityCondition} 
\end{equation}
\end{proposition}
The proof is given in Appendix~\ref{app:stability} and follows from the general theory for QBD processes. 

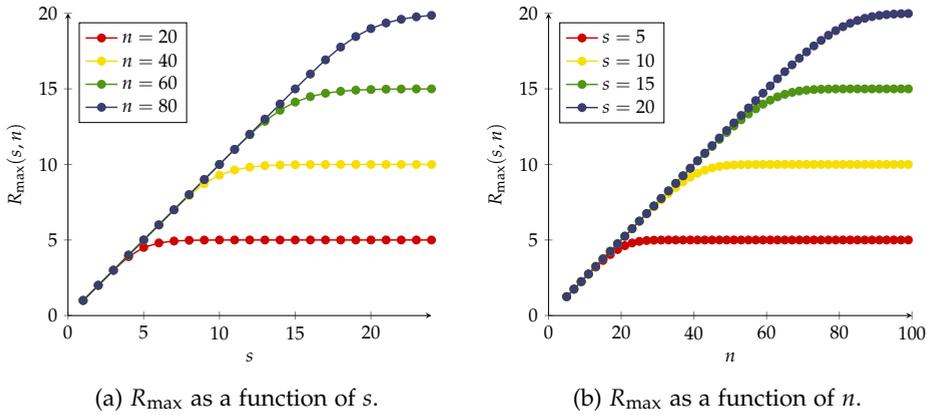
\begin{figure}
\centering
\begin{subfigure}{0.48\textwidth}
\begin{tikzpicture}[scale = 0.7]
\begin{axis}[
	xmin = 0,
	xmax = 24,
	ymin = 0,
	ymax = 20,	
	grid = none, 
	axis line style={->},
	axis lines = left,
	xlabel = $s$,
	ylabel = {$R_{\rm max}(s,n)$},
	legend cell align=left,
	legend pos = north west
]

\addplot[col1,thick,mark=*] table[x=s,y=n20] {Chapter_5/tikz/stability/r025_n_fixed.txt};
\addplot[col3,thick,mark=*] table[x=s,y=n40] {Chapter_5/tikz/stability/r025_n_fixed.txt};
\addplot[col4,thick,mark=*] table[x=s,y=n60] {Chapter_5/tikz/stability/r025_n_fixed.txt};
\addplot[col5,thick,mark=*] table[x=s,y=n80] {Chapter_5/tikz/stability/r025_n_fixed.txt};

\legend{$n=20$,$n=40$,$n=60$,$n=80$}
\end{axis}
\end{tikzpicture}
\caption{$R_{\rm max}$ as a function of $s$.}
\label{fig:Rmax_a}
\end{subfigure}
\begin{subfigure}{0.48\textwidth}
\begin{tikzpicture}[scale = 0.7]
\begin{axis}[
	xmin = 0,
	xmax = 100,
	ymin = 0,
	ymax = 20,	
	grid = none, 
	axis line style={->},
	axis lines = left,
	xlabel = $n$,
	ylabel = {$R_{\rm max}(s,n)$},
	legend cell align=left,
	legend pos = north west
]

\addplot[col1,thick,mark=*,mark repeat = 2] table[x=n,y=s5] {Chapter_5/tikz/stability/r025_s_fixed.txt};
\addplot[col3,thick,mark=*,mark repeat = 2] table[x=n,y=s10] {Chapter_5/tikz/stability/r025_s_fixed.txt};
\addplot[col4,thick,mark=*,mark repeat = 2] table[x=n,y=s15] {Chapter_5/tikz/stability/r025_s_fixed.txt};
\addplot[col5,thick,mark=*,mark repeat = 2] table[x=n,y=s20] {Chapter_5/tikz/stability/r025_s_fixed.txt};

\legend{$s=5$,$s=10$,$s=15$,$s=20$}
\end{axis}
\end{tikzpicture}
\caption{$R_{\rm max}$ as a function of $n$.}
\label{fig:Rmax_b}
\end{subfigure}
\caption{The maximum achievable workload in the restricted Erlang-R model with holding for $r=0.25$.}
\label{fig:Rmax}
\end{figure}

Observe that $\rho_{\max}(s,n)$ poses an upper bound on the occupancy level of the servers in the holding model, which is clearly smaller than 1 for all $s$ and $n$. 
In addition, this implies that the maximum workload $R_{\max}(s,n) := s\cdot\rho_{\max}(s,n)$ the system is able to handle is strictly less than $s$.
If we compare this to the open Erlang-R model, in which the maximal attainable workload equals $s$, we observe the effect of finite-size constraints on operational performance. 
Figure \ref{fig:Rmax} shows the influence of both $s$ and $n$ on the maximum feasible workload in case $r=0.25$. 
From these graphs, note that if $s\ll rn$, $R_{\max}$ grows almost linearly with $s$. 
Furthermore, $R_{\rm max}(s,n)$ is increasing in $n$ for $s$ fixed.
A logical practical consequence is that a larger number of beds allows for a larger patient volume to enter the ED with the same number of nurses.
Moreover, $R_{\rm max}(s,n)$ is increasing in $s$, but as in Figure \ref{fig:Rmax_a}, adding an extra nurse does not increase the stability region in case $n$ is too tight. 
 Conversely, adding extra beds does not increase $R_{\rm max}(s,n)$ if the number of nurses does not allow for an increase in offered load, see Figure \ref{fig:Rmax_b}. 
Additionally, it is easily verified that $R_{\rm max}(s,n)$ is upper bounded by both $s$ and $R_{\rm max}(n,n) = rn$. Therefore, a careful balance is called for between servers (nurses) and beds, so that resources will be efficiently utilized. We observe that when the ratio $s/n\approx r$, the system is better balanced. 
We will propose an appropriate balance between resources by defining a synchronized QED capacity recommendation for both servers and beds in Section  \ref{sec:QED_scaling}.

Provided that the system is stable, the stationary distribution of the QBD process $X$ can be obtained numerically by the matrix geometric method \citep{Neuts1981}.
Subsequently, we can derive the stationary distribution of the original $Q(t)$, denoted by $\pi_h(\cdot,\cdot,\cdot)$.

\subsection{Performance measures}
\label{sec:performance_metrics}
In this work, we concentrate on five performance measures that are central to our analysis. 
In the definitions that follow, we present expressions for these measures in terms of a general three-dimensional measure $\pi$, which one can replace by either $\pi_b$ or $\pi_h$, depending on the scenario considered. 
In the remainder of this work, we will augment the measures related to the Erlang-R model with blocking and holding by the superscript $b$ and $h$, respectively\footnote{In line with $H(t)=0$, we use $\pi_b(i,j,k) = \pi_b(j,k)$ if $i=0$, with $\pi_b(j,k)$ as in \eqref{eq:pih(i,j)}, and $\pi_b(i,j,k) = 0$ otherwise, when considering the model with blocking.}.

As relevant performance measures, we consider the probability of holding (cq.\\ \noindent blocking) at entering the system, the probability of delay at the needy queue, expected waiting time for a nurse, utilization of nurses and utilization of beds:
\begin{equation}
\P({\rm hold}) = \sum_{i=0}^\iy \sum_{j=0}^n \pi(i,j,n-j), \qquad
\P({\rm delay}) \approx \sum_{i=0}^\infty\sum_{j=s}^{n}\sum_{k=0}^{n-j} \pi(i,j,k), 
\label{eq:delay_probability}
\end{equation}
\begin{equation}
\label{eq:EW_exact}
\E [W] \approx \sum_{i=0}^\infty\sum_{j=s}^{n}\sum_{k=0}^{n-j} \frac{\max\{0,j-s+1\}}{\mu}\,\pi(i,j,k), 
\end{equation}
\begin{equation}
\label{eq:utilization}
\rho_s = \frac{1}{s}\,\sum_{i=0}^\infty \sum_{j=0}^n \sum_{k=0}^{n-j} \min\{j,s\} \pi(i,j,k), \qquad
\rho_n = \frac{1}{n}\,\sum_{i=0}^\infty \sum_{j=0}^n \sum_{k=0}^{n-j} \min\{i,n\} \pi(i,j,k).
\end{equation}

It should be stressed that the above expression for the delay probability and the expected waiting time for a nurse are not exact. For the blocking model one can use the Arrival Theorem, see e.g.~\cite{Chen2001a}, whereby the exact expression sums up to $n-1$ instead of $n$. 
Since we consider the system as $n\to\infty$, this discrepancy becomes negligible. 
For the holding model, a similar argument holds. 
We will therefore use the expressions in \eqref{eq:delay_probability}-\eqref{eq:utilization} as definitions for the performance measures. 

\subsection{Stochastic bounds}
\label{sec:bounds}


Although the two variants of the Erlang-R model differ with respect to the admission policy, and require different mathematical treatment, we would like to be able to capture their relative performance.  
We substantiate the intuition that the holding room leads to more patients in the ED, in the following result. 

\begin{proposition}\label{thm:stochasticordering}
Let $Q_1^b$, $Q_2^b$, $Q_1^h$, $Q_2^h$ denote the nurse and content queue length processes in the Erlang-R model with blocking and holding, respectively.
Let $H(0) = 0$,  $Q_1^b(0)=Q_1^h(0)$ and $Q_2^b(0)=Q_2^h(0)$. For all $t\geq 0$, 
\begin{align}
Q_1^b(t) + Q_2^b(t) &\preceq_{\rm st} Q_1^h(t) + Q_2^h(t) \preceq_{\rm st} n ,\\
Q_2^b(t) &\preceq_{\rm st} Q_2^h(t),\\
Q_1^b(t) &\preceq_{\rm st} Q_1^h(t) + H(t),
\end{align}
where $X\preceq_{\rm st} Y$ implies $\P(X\geq k) \leq \P(Y\geq k)$ for all $k\geq 0$.
\end{proposition}
\noindent
The proof of Proposition \ref{thm:stochasticordering} uses sample path coupling and can be found in Appendix \ref{app:stochastic_ordering}.
Note that as an immediate consequence, we have
\[ \P^b( {\rm block}) = \lim_{t\to\iy} \P\big( Q_1^b(t)+Q_2^b(t) \geq n \big) \leq \lim_{t\to\iy} \P\big( Q_1^h(t) + Q_2^h(t) \geq n \big) = \P^h( {\rm hold }) \]
and by similar reasoning $\rho^b_n \leq  \rho_n^h$.
In other words, under similar offered load and capacity constraints, utilization levels for the nurses in the Erlang-R model with blocking are lower than in the Erlang-R model with holding. 
Moreover, the total number of waiting patients in the setting with holding is stochastically larger than in the setting with blocking, and in the open Erlang-R model. 
We further discuss the differences between both models in Section \ref{sec:dimensioning} and Section \ref{sec:analysis_chapter5}.

\section{Two-fold QED regime}
\label{sec:QED_scaling}

We do not want to waste capacity of either servers or beds without getting significant advantage in terms of performance. 
We therefore take an asymptotic approach that lets the external arrival rate $\l$ grow to infinity, while scaling $s$ and $n$ accordingly.
In doing so, we intend to establish QED-type system behavior, i.e.\ high occupancy levels of both nurses and beds and good quality-of-service.

\subsection{Two-fold scaling rule}

In order to identify the scaling of $s$ and $n$ as $\l\to\infty$, we draw inspiration from the two-fold scaling rule used by Jennings \& de V\'ericourt \cite{Jennings2008} and Khudyakov et al.~\cite{Khudyakov2010}, which follows the celebrated square-root staffing principle.
This principle suggests that, in the most general setting, capacity should be equal to the expected offered load entering the system, let us say $R$, plus an additional variability hedge that is proportional to $\sqrt{R}$.
In the restricted Erlang-R model, we have two capacity sources, namely $s$ and $n$, which experience different relevant amounts of work.

The offered load the servers in the needy queue experience is given by $R_{\rm nurse} = R_1$, as in the regular Erlang-R model; 
it does not change due to the finite-size effects, since all patients are served eventually. Hence, we only need to account for the interrupted services. It follows that the appropriate staffing rule for the nurses in the QED regime remains $s=R_1+\beta \sqrt{R_1}$ for some constant $\beta >0$.

To establish the bed capacity level, we need to reflect on the load offered to the beds. Observe that beds remain occupied both in needy and content states. This suggests that $R_{\rm bed} :=R_1+R_2=R_1/r$, with $R_1$ and $R_2$ as in \eqref{eq:R1_R2} and $r$ is the expected fraction of time a patient spends at the nurse station defined in \eqref{eq:delta}.
%
As a result, the appropriate staffing rule is $n=R_{\rm bed}+\gamma \sqrt{R_{\rm bed}}$ for some constant $\gamma>0$. In conclusion, the two-fold QED scaling rule is given by
\begin{equation}\label{eq:twofoldscaling}
\begin{array}{ll}
s &= R_1 + \beta \sqrt{R_1} + o(\sqrt{R_1}) \\
n &= \frac{R_1}{r}+\gamma \sqrt{\frac{R_1}{r}} + o(\sqrt{R_1})
\end{array}
\end{equation}
with $\beta,\gamma>0$ constants and $R_1:=\lambda/((1-p)\mu)$.

Recall that we saw in Figure \ref{fig:Rmax} that resources seem efficiently utilized if $s/n\approx r$. 
Scaling \eqref{eq:twofoldscaling} is in line with this reasoning since 
\[
\frac{s}{n} = r\left(1+ \frac{\beta - \gamma\sqrt{r}}{\sqrt{R_1}}+ O(1/R_1) \right) . 
\]

\begin{remark}
In \cite{Jennings2008}, a similar scaling regime is considered, which only relates $s$ and $n$ through a square-root scaling, namely the regime $s = r n + \hat\g\sqrt{n}$,
which is equivalent to the second relation in \eqref{eq:twofoldscaling} if $\hat\g = \b\sqrt{r} - \g r$.
Due to the absence of external arrivals in this closed system, they let the number of beds $n$ approach infinity as opposed to $\l$ in our settings. 
Nevertheless, this results in the same asymptotic regime.
\end{remark}

Before turning to asymptotic expressions for the performance measures concerning the Erlang-R model with blocking or holding, we conduct a few numerical experiments to confirm that the scaling in \eqref{eq:twofoldscaling} indeed leads to desired QED behavior. 

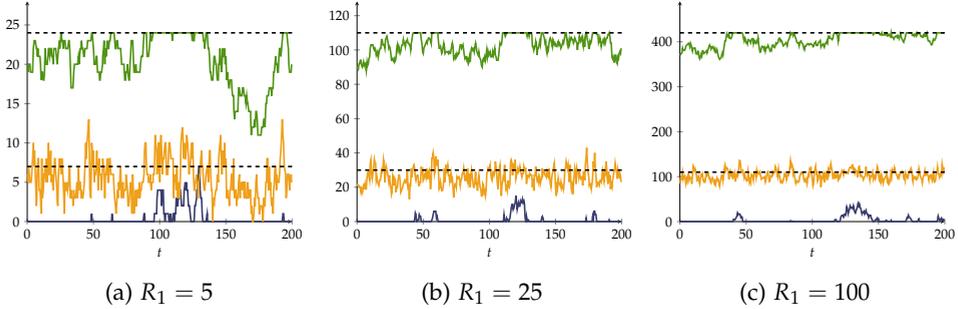
\begin{figure} 
\centering
\begin{subfigure}{0.32\textwidth}
\begin{tikzpicture}[scale=0.51]

\begin{axis}[
	xmin = 0,
	xmax = 200,
	ymin = 0.0,
	ymax = 28,
	ytick = {0,5,10,15,20,25},	
	grid = none, 
	axis line style={->},
	axis lines = left,
	xlabel = $t$,
	legend cell align=left,
	legend pos = north east
]

\definecolor{col1}{rgb}{0.368417, 0.506779, 0.709798}
 	
\addplot[very thick,col5] file {Chapter_5/tikz/sample_paths/R5_holding.txt};
\addplot[very thick,col2] file {Chapter_5/tikz/sample_paths/R5_service.txt};
\addplot[very thick,col4] file {Chapter_5/tikz/sample_paths/R5_total.txt};
\addplot[very thick,dashed] coordinates {
	(0,7)
	(200,7)
	};
\addplot[very thick,dashed] coordinates {
	(0,24)
	(200,24)
	};

\end{axis}

\end{tikzpicture}
\caption{$R_1=5$}
\end{subfigure}
\begin{subfigure}{0.32\textwidth}
\begin{tikzpicture}[scale=0.51]

\begin{axis}[
	xmin = 0,
	xmax = 200,
	ymin = 0.0,
	ymax = 128.333,
	ytick = {0,20,40,60,80,100,120},	
	grid = none, 
	axis line style={->},
	axis lines = left,
	xlabel = $t$,
	legend cell align=left,
	legend pos = north east
]

\addplot[very thick,col5] file {Chapter_5//tikz/sample_paths/R25_holding.txt};
\addplot[very thick,col2] file {Chapter_5//tikz/sample_paths/R25_service.txt};
\addplot[very thick,col4] file {Chapter_5/tikz/sample_paths/R25_total.txt};
\addplot[very thick,dashed] coordinates {
	(0,30)
	(200,30)
	};
\addplot[very thick,dashed] coordinates {
	(0,110)
	(200,110)
	};
\end{axis}

\end{tikzpicture}
\caption{$R_1=25$}
\end{subfigure}
\begin{subfigure}{0.32\textwidth}
\begin{tikzpicture}[scale=0.51]

\begin{axis}[
	xmin = 0,
	xmax = 200,
	ymin = 0.0,
	ymax = 490,
	ytick = {0,100,200,300,400},	
	grid = none, 
	axis line style={->},
	axis lines = left,
	xlabel = $t$,
	legend cell align=left,
	legend pos = north east
]

\addplot[very thick,col5] file {Chapter_5/tikz/sample_paths/R100_holding.txt};
\addplot[very thick,col2] file {Chapter_5/tikz/sample_paths/R100_service.txt};
\addplot[very thick,col4] file {Chapter_5/tikz/sample_paths/R100_total.txt};
\addplot[very thick,dashed] coordinates {
	(0,110)
	(200,110)
	};
\addplot[very thick,dashed] coordinates {
	(0,420)
	(200,420)
	};
	
\end{axis}

\end{tikzpicture}
\caption{$R_1=100$}
\end{subfigure}
\caption{Sample paths of $H(t)$ (blue), $Q_1(t)$ (orange) and $Q_1(t)+Q_2(t)$ (green) of the Erlang-R model with holding with parameters $\mu = 1$, $\d=0.25$, $p=0.75$ and $\b=\g=1$. The staffing levels $s$ and $n$ are depicted by the dashed lines.}
\label{fig:sample_paths}
\end{figure}

In Figure \ref{fig:sample_paths}, we plotted the sample paths of the three-dimensional queue length process of the holding model in which $\b$ and $\g$ are fixed, and $R_1$ is increased. 
Observe that the needy queue length $Q_1(t)$, plotted in orange in Figure \ref{fig:sample_paths}, fluctuates around the values $s$, and stabilizes for larger values of $R_1$. 
This naturally implies that the server (nurses) utilization approaches 100\%, while the number of patients waiting is $O(\sqrt{R_1})$.
Furthermore, we see that the percentage of occupied beds also tends to 100\%, while the holding queue length remains small. 
The holding queue is of much smaller order than $R_1$, which implies that the holding time of a patient becomes negligible as $R_1\to\iy$. 
From these empirical findings we deduce that  under scaling \eqref{eq:twofoldscaling} the restricted Erlang-R model exhibits QED behavior on two levels: Outside the facility while waiting for an available bed, and inside the facility while waiting for attention of a nurse.

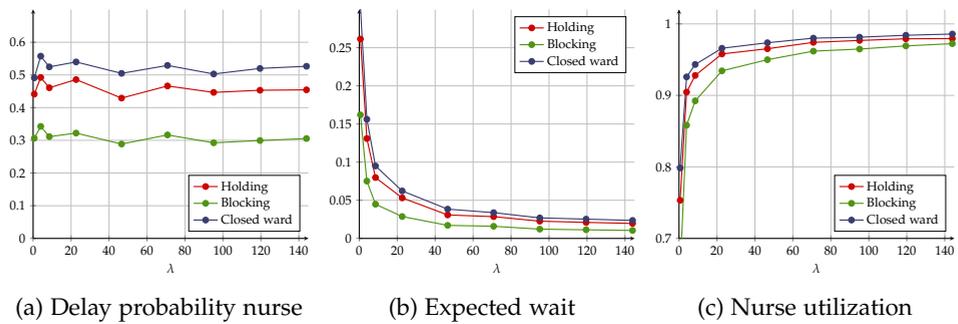
\begin{figure}
\centering
\begin{subfigure}{0.32\textwidth}
\centering
\begin{tikzpicture}[scale=0.53]
\begin{axis}[
	xmin = 0,
	xmax = 145,
	ymin = 0.0,
	ymax = 0.7,
	ytick = {0,0.1,...,0.7},
	xlabel = $\l$,
	grid = both, 
	axis line style={->},
	axis lines = left,
	legend cell align=left,
	legend pos = south east
]

\addplot[thick,col1,mark=*] file {Chapter_5/tikz/empirical/delayProbErlangH.txt};
\addplot[thick,col4,mark=*] file {Chapter_5/tikz/empirical/delayProbYomTov.txt};
\addplot[thick,col5,mark=*] file {Chapter_5/tikz/empirical/delayProbJennings.txt};
\small
\legend{Holding,Blocking,Closed ward};
\end{axis}

\end{tikzpicture}
\caption{Delay probability nurse}
\end{subfigure}
\begin{subfigure}{0.32\textwidth}
\centering
\begin{tikzpicture}[scale=0.53]

\begin{axis}[
	xmin = 0,
	xmax = 145,
	ymin = 0.0,
	ymax = 0.3,
	ytick = {0,{0.05},0.1,0.15,0.2,0.25,3},	
	grid = both, 
	axis line style={->},
	tick label style={/pgf/number format/fixed},
	axis lines = left,
	xlabel = $\l$,
	legend cell align=left,
	legend pos = north east
]

\addplot[thick,col1,mark=*] file {Chapter_5/tikz/empirical/EWErlangH.txt};
\addplot[thick,col4,mark=*] file {Chapter_5/tikz/empirical/EWYomTov.txt};
\addplot[thick,col5,mark=*] file {Chapter_5/tikz/empirical/EWJennings.txt};
\small
\legend{Holding,Blocking,Closed ward};
\end{axis}

\end{tikzpicture}
\caption{Expected wait}
\end{subfigure}
\begin{subfigure}{0.32\textwidth}
\centering
\begin{tikzpicture}[scale=0.53]

\begin{axis}[
	xmin = 0,
	xmax = 145,
	ymin = 0.7,
	ymax = 1.02,
	grid = both, 
	axis line style={->},
	axis lines = left,
	xlabel = $\l$,
	legend cell align=left,
	legend pos = south east
]
 	
\addplot[thick,col1,mark=*] file {Chapter_5/tikz/empirical/rhoErlangH.txt};
\addplot[thick,col4,mark=*] file {Chapter_5/tikz/empirical/rhoYomTov.txt};
\addplot[thick,col5,mark=*] file {Chapter_5/tikz/empirical/rhoJennings.txt};
\small
\legend{Holding,Blocking,Closed ward};
\end{axis}

\end{tikzpicture}
\caption{Nurse utilization}
\end{subfigure}
\caption{Asymptotic behavior of the restricted Erlang-R model with holding and blocking, and the closed ward model for $\m=1$, $\d = 0.2$, $p=0.8$ and $\b=\g=0.5$.}
\label{fig:empiricalAsymptotics}
\end{figure}

We also check how the Erlang-R model with blocking or holding  and the closed ward model of \cite{Jennings2008} relate under scaling \eqref{eq:twofoldscaling}. 
In Figure~\ref{fig:empiricalAsymptotics}, we plot the performance measures, obtained through simulation, for the three models in which we fix $\b=\g=0.5$ and vary the arrival rate $\l$. 
First, we see that $\P({\rm delay})$ stabilizes as $\l\to\iy$ in all three models under scaling \eqref{eq:twofoldscaling}, and the delay probability of the model with holding lies in between the other two. 
Second, note that the expected waiting time for a nurse in all models converges to 0 as $\l$ increases. In fact, the rate of decay is similar in all three models.  
We observe that $\rho_s$ approaches unity in all models, and the rate of convergence seems again comparable. 
Finally, and most importantly, we notice an ordering between the three models.
Namely, in all performance measures considered in Figure \ref{fig:empiricalAsymptotics}, Erlang-R with holding appears to be upper bounded by the closed ward and lower bounded by the Erlang-R with blocking. 
In a multitude of parameter settings of $(\b,\g)$, we have seen the same ordering, leading to the following conjecture: 
\begin{conjecture}\label{conj:stochorder}
Let $Q^b_1(\iy)$, $Q_1^h(\iy)$ and $Q_1^J(\iy)$ denote the stationary number of needy patients in the Erlang-R model with blocking, holding and the closed ward, respectively. Then,
\begin{equation}
Q_1^b(\iy) \preceq_{\rm st} Q_1^h(\iy) \preceq_{\rm st} Q_1^J(\iy).
\end{equation}
\end{conjecture}
Observe that Conjecture \ref{conj:stochorder} poses a stronger statement than the third assertion in Proposition \ref{thm:stochasticordering}. 
The latter does give an upper bound to $Q_1^h(\iy)$ in terms of $Q_1^b(\iy)$, albeit supplemented with the stationary holding queue length. 

\subsection{QED limits for Erlang-R with blocking}
\label{sec:QED_limit_block}

We now continue our analysis by examining its limiting behavior under scaling \eqref{eq:twofoldscaling}, and obtain QED limits for some performance measures of the Erlang-R model with blocking.
Using the explicit expressions for the blocking model in \eqref{eq:pih(i,j)}, we derive the limiting values of the relevant performance measures defined in Section \ref{sec:performance_metrics} in terms of $\b$ and $\g$. 

\begin{theorem}\label{thm:limits_YT}
Let $s$ and $n$ scale as in \eqref{eq:twofoldscaling} with ${-}\infty<\b<\infty,\,\g>0$ as $\l\to\infty$. Then, if $\b \neq 0$,
\begin{align}
g^b(\b,\g)
&:= \lim_{\l\to\iy} \P^b({\rm delay})\nonumber \\
\label{eq:yt_limit_delay}
&= 
\left(1 + 
\frac{ \b \, \int_{-\iy}^\b \Phi\left(\frac{\g-t\sqrt{r}}{\sqrt{1-r}}\right)\, \dd\Phi(t) }
{\f(\b)\Phi(\eta) -  \f(\sqrt{\b^2+\eta^2}){\rm e}^{\tfrac{1}{2} \omega^2} \Phi(\omega)}
\right)^{-1},\\
f^b(\b,\g) 
&:= \lim_{\l\to\iy} \sqrt{R_1}\cdot\P^b({\rm block}) \nonumber\\
\label{eq:yt_limit_block}
&= 
\frac{
\sqrt{r}\f(\g)\Phi(-\omega\sqrt{r}) + \f(\sqrt{\b^2+\eta^2})\,{\rm e}^{\frac{1}{2} \omega^2} \Phi(\omega)
}{
\int_{-\iy}^\b \Phi\left(\frac{\g-t\sqrt{r}}{\sqrt{1-r}}\right)\, \dd\Phi(t) +
\frac{\f(\b)\Phi(\eta)}{\b} -  \frac{\f(\sqrt{\b^2+\eta^2})}{\b}{\rm e}^{\tfrac{1}{2} \omega^2} \Phi(\omega)
 },\\
h^b(\b,\g) &:= \lim_{\l\to\iy} \sqrt{R_1}\cdot\E[W] \nonumber\\
\label{eq:yt_limit_Edelay}
&=
\frac{
\frac{\f(\b)\Phi(\eta)}{\b^2} +
 \left(\frac{\b}{r}-\frac{\g}{\sqrt{r}}-\frac{1}{\b}\right)\,\frac{\f(\sqrt{\eta^2+\b^2})}{\b}\, {\rm e}^{\tfrac{1}{2}\omega^2}\, \Phi(\omega) 
 - \sqrt{\frac{1-r}{r}}\,\frac{\f(\b)\f(\eta)}{\b}
}{
\int_{-\iy}^\b \Phi\left(\frac{\g-t\sqrt{r}}{\sqrt{1-r}}\right)\, \dd\Phi(t) +
\frac{\f(\b)\Phi(\eta)}{\b} -  \frac{\f(\sqrt{\b^2+\eta^2})}{\b}{\rm e}^{\tfrac{1}{2} \omega^2} \Phi(\omega)
 },
\end{align}
and if $\b=0$,
\begin{align}
g^b_0(\g) 
&:= \lim_{\l\to\iy} \P^b({\rm delay}) \nonumber\\
\label{eq:yt_limit_delay_beta0}
&=
\left(1+
\frac{
\int_{-\iy}^0 \Phi\left(\frac{\g-t\sqrt{r}}{\sqrt{1-r}}\right)\, \dd\Phi(t)
}{
\sqrt{\frac{1-r}{r}} \frac{1}{\sqrt{2\pi}}\,\left(\eta \,\Phi(\eta) + \f(\eta) \right)
}
\right)^{-1},\\
f^b_0(\g) 
&:= \lim_{\l\to\iy} \sqrt{R_1}\cdot\P^b({\rm block}) \nonumber \\
\label{eq:yt_limit_block_beta0}
&= 
\frac{
\sqrt{r}\,\f(\g)\Phi(-\omega\sqrt{r}) + \frac{1}{\sqrt{2\pi}} \Phi(\eta)
}{
\int_{-\iy}^\b \Phi\left(\frac{\g-t\sqrt{r}}{\sqrt{1-r}}\right)\, \dd\Phi(t) +
\sqrt{\frac{1-r}{r}} \frac{1}{\sqrt{2\pi}}\,\left(\eta \,\Phi(\eta) + \f(\eta) \right)
 },\\
h_0^b(\g) &:= \lim_{\l\to\iy} \sqrt{R_1}\cdot\E[W] \nonumber\\
\label{eq:yt_limit_Edelay_beta0}
&= \frac{1}{2\mu}\, \frac{ \left( \gamma^2/r+1\right) \Phi(\eta) + \eta \f(\eta) }
{ \frac{r}{1-r} \sqrt{2\pi} \int_{-\infty}^0 \Phi\left(\frac{\g-t\sqrt{r}}{\sqrt{1-r}}\right)\, \dd\Phi(t) + \sqrt{\frac{r}{1-r}} \left(\eta \Phi(\eta)+\f(\eta)\right)},
\end{align}
where $\eta = \frac{\g - \b\sqrt{r}}{\sqrt{1-r}}$ and $\omega := \frac{\g - \b/\sqrt{r}}{\sqrt{1-r}}$.
\end{theorem}
The proof of Theorem \ref{thm:limits_YT} is given in Appendix C of \cite{YomTov2010} under a parameter transformation.

\begin{figure}
\centering
\begin{subfigure}{0.48\textwidth}
\centering
\begin{tikzpicture}[scale = 0.8]
\small
\begin{axis}[
	xmin = -2,
	xmax = 2,
	ymin = 0,
	ymax = 1,
	grid = both, 
	axis line style={->},
	axis lines = left,
	tick label style={/pgf/number format/fixed},
	xlabel = $\beta$,
	ylabel = {$g(\beta,\gamma)$},
	y label style = {at = {(axis cs: -2.4,0.5)}},
	xscale = 0.9,
	yscale = 0.75,
	legend cell align=left,
	legend style = {at = {(axis cs: 1.95,0.99)},anchor = north east}
]

\addplot[thick,col1] table[x=beta,y=g_min1] {Chapter_5/tikz/limit_probabilities_delay.txt};
\addplot[thick,col3] table[x=beta,y=g_0] {Chapter_5/tikz/limit_probabilities_delay.txt};
\addplot[thick,col4] table[x=beta,y=g_1] {Chapter_5/tikz/limit_probabilities_delay.txt};
\addplot[thick,col5] table[x=beta,y=g_2] {Chapter_5/tikz/limit_probabilities_delay.txt};

\legend{{$\gamma = -1$},{$\gamma = 0$},{$\gamma = 1$},{$\gamma = 2$}};
\end{axis}
\end{tikzpicture}
\caption{Delay probability}
\label{fig:pdelay_pblock_a}
\end{subfigure}
\begin{subfigure}{0.48\textwidth}
\centering
\begin{tikzpicture}[scale = 0.8]
\small
\begin{axis}[
	xmin = -2,
	xmax = 2,
	ymin = 0,
	ymax = 2,
	grid = both, 
	axis line style={->},
	axis lines = left,
	xlabel = $\beta$,
	tick label style={/pgf/number format/fixed},
	ylabel = {$f(\beta,\gamma)$},
	y label style = {at = {(axis cs: -2.4,1)}},
	xscale = 0.9,
	yscale = 0.75,
	legend cell align=left,
	legend style = {at = {(axis cs: 1.95,1.98)},anchor = north east}
]

\addplot[thick,col1] table[x=beta,y=g_min1] {Chapter_5/tikz/limit_probabilities_block.txt};
\addplot[thick,col3] table[x=beta,y=g_0] {Chapter_5/tikz/limit_probabilities_block.txt};
\addplot[thick,col4] table[x=beta,y=g_1] {Chapter_5/tikz/limit_probabilities_block.txt};
\addplot[thick,col5] table[x=beta,y=g_2] {Chapter_5/tikz/limit_probabilities_block.txt};


\legend{{$\gamma = -1$},{$\gamma = 0$},{$\gamma = 1$},{$\gamma = 2$}};
\end{axis}
\end{tikzpicture}
\caption{Scaled blocking probability}
\label{fig:pdelay_pblock_b}
\end{subfigure}
\caption{Asymptotic delay and scaled blocking probability for $r=0.5$ based on Theorem \ref{thm:limits_YT}. }
\label{fig:pdelay_pblock}
\end{figure}
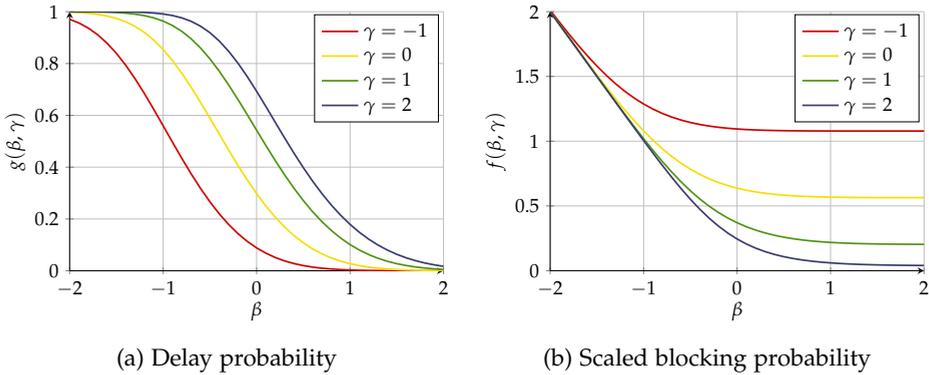

Theorem \ref{thm:limits_YT} proves that the scaling \eqref{eq:twofoldscaling} results in QED behavior: the probability of waiting in Equations \eqref{eq:yt_limit_delay} and \eqref{eq:yt_limit_delay_beta0} converges to a limit that is strictly between 0 and 1. 
Notice that all limits in Theorem \ref{thm:limits_YT} are functions of three parameters: $\beta$ and $\gamma$, which are decision variables, and the fraction of needy time $r$, which is dictated by the physics of the system. Furthermore, the theorem also shows that the probability of blocking (Equations \eqref{eq:yt_limit_block} and \eqref{eq:yt_limit_block_beta0}) is of order $1/\sqrt{R_1}$. 
For example, assume that the fraction of needy time $r$ is $0.5$ and the system is large (100 servers). 
Using Figure \ref{fig:pdelay_pblock}, we observe that, by choosing the pair $\gamma = 1$ and $\beta = 0.245$, we actually aim at a probability of getting served immediately to be 40\%. At the same time, the probability of getting immediately a bed is 97\%. 
Thus, waiting inside the ED occurs at a reasonable level, while wait outside the facility becomes negligible.

Theorem \ref{thm:limits_YT} further shows that the expected waiting (Equations \eqref{eq:yt_limit_Edelay} and \eqref{eq:yt_limit_Edelay_beta0})
is of order $1/\sqrt{R_1}$ too and hence vanishes in the large-system limit.

We see from Theorem \ref{thm:limits_YT} that achieving target service levels is always an interplay between $\beta$ and $\gamma$.
Figure \ref{fig:pdelay_pblock_a} shows for instance that in order to keep $\P({\rm delay})\in (0.25,0.75)$,  choosing $\gamma=-1$ requires $\beta$ to stay within the range $[-1.4,-0.5]$, while $\gamma=1$ corresponds to values of $\beta$ in $[-0.4,0.5]$.

While the two-fold scaling rule in \eqref{eq:twofoldscaling} automatically captures the right dimensioning ratio as the system scales up, Theorem \ref{thm:limits_YT} shows that the parameters $\beta$ and $\gamma$ provide a means to fine-tune the performance.
Figure \ref{fig:pdelay_pblock_b} confirms how adding nurses, i.e.~increasing $\beta$, does not improve the blocking probability if the number of beds, i.e.~$\gamma$, is too tight.
This is in accordance with our previous observations in Figure \ref{fig:Rmax} for the exact steady-state distribution.

\begin{figure}
\centering
\begin{subfigure}{0.48\textwidth}
\begin{tikzpicture}[scale = 0.87]
\small
\begin{axis}[
	xmin = 0,
	xmax = 16,
	ymin = 0,
	ymax = 1,
	grid = both, 
	axis line style={->},
	axis lines = left,
	xlabel = $s$,
	xscale = 0.9,
	yscale = 0.75,
	legend cell align=left,
	legend style = {at = {(axis cs: 0.1,0.01)},anchor = south west}
]

\addplot[thick,col1,mark = *] table[x=s,y=n24] {Chapter_5/tikz/accuracy/accuracy_pdelay_case2.txt};
\addplot[thick,col2,mark = *] table[x=s,y=n28] {Chapter_5/tikz/accuracy/accuracy_pdelay_case2.txt};
\addplot[thick,col3,mark = *] table[x=s,y=n32] {Chapter_5/tikz/accuracy/accuracy_pdelay_case2.txt};
\addplot[thick,col4,mark = *] table[x=s,y=n36] {Chapter_5/tikz/accuracy/accuracy_pdelay_case2.txt};
\addplot[thick,col5,mark = *] table[x=s,y=n40] {Chapter_5/tikz/accuracy/accuracy_pdelay_case2.txt};

\addplot[thick,col1,dashed] table[x=s,y=approx_n24] {Chapter_5/tikz/accuracy/accuracy_pdelay_case2.txt};
\addplot[thick,col2,dashed] table[x=s,y=approx_n28] {Chapter_5/tikz/accuracy/accuracy_pdelay_case2.txt};
\addplot[thick,col3,dashed] table[x=s,y=approx_n32] {Chapter_5/tikz/accuracy/accuracy_pdelay_case2.txt};
\addplot[thick,col4,dashed] table[x=s,y=approx_n36] {Chapter_5/tikz/accuracy/accuracy_pdelay_case2.txt};
\addplot[thick,col5,dashed] table[x=s,y=approx_n40] {Chapter_5/tikz/accuracy/accuracy_pdelay_case2.txt};

\legend{{$n=24$},{$n=28$},{$n=32$},{$n=36$},{$n=40$}};
\end{axis}
\end{tikzpicture}
\caption{Delay probability}
\end{subfigure}
\begin{subfigure}{0.48\textwidth}
\begin{tikzpicture}[scale = 0.87]
\small
\begin{axis}[
	xmin = 0,
	xmax = 16,
	ymin = 0,
	ymax = 2,
	grid = both, 
	axis line style={->},
	axis lines = left,
	xlabel = $s$,
	xscale = 0.9,
	yscale = 0.75,
	legend cell align=left,
	legend style = {at = {(axis cs: 0.1,0.01)},anchor = south west}
]

\addplot[thick,col1,mark = *] table[x=s,y=n24] {Chapter_5/tikz/accuracy/accuracy_pblock_case2.txt};
\addplot[thick,col2,mark = *] table[x=s,y=n28] {Chapter_5/tikz/accuracy/accuracy_pblock_case2.txt};
\addplot[thick,col3,mark = *] table[x=s,y=n32] {Chapter_5/tikz/accuracy/accuracy_pblock_case2.txt};
\addplot[thick,col4,mark = *] table[x=s,y=n36] {Chapter_5/tikz/accuracy/accuracy_pblock_case2.txt};
\addplot[thick,col5,mark = *] table[x=s,y=n40] {Chapter_5/tikz/accuracy/accuracy_pblock_case2.txt};

\addplot[thick,col1,dashed] table[x=s,y=approx_n24] {Chapter_5/tikz/accuracy/accuracy_pblock_case2.txt};
\addplot[thick,col2,dashed] table[x=s,y=approx_n28] {Chapter_5/tikz/accuracy/accuracy_pblock_case2.txt};
\addplot[thick,col3,dashed] table[x=s,y=approx_n32] {Chapter_5/tikz/accuracy/accuracy_pblock_case2.txt};
\addplot[thick,col4,dashed] table[x=s,y=approx_n36] {Chapter_5/tikz/accuracy/accuracy_pblock_case2.txt};
\addplot[thick,col5,dashed] table[x=s,y=approx_n40] {Chapter_5/tikz/accuracy/accuracy_pblock_case2.txt};

\legend{{$n=24$},{$n=28$},{$n=32$},{$n=36$},{$n=40$}};
\end{axis}
\end{tikzpicture}
\caption{Scaled blocking probability}
\end{subfigure}
%
\caption{Comparison of exact performance measures (solid) against asymptotic approximations (dashed) with $\beta=(s-R_1)/\sqrt{R_1}$ and $\gamma=(n-R_1/r)/\sqrt{R_1/r}$ for $\lambda = 2$, $\mu=1$, $\delta=0.25$ and $p=0.75$.}
\label{fig:accuracy_blocking}
\end{figure}
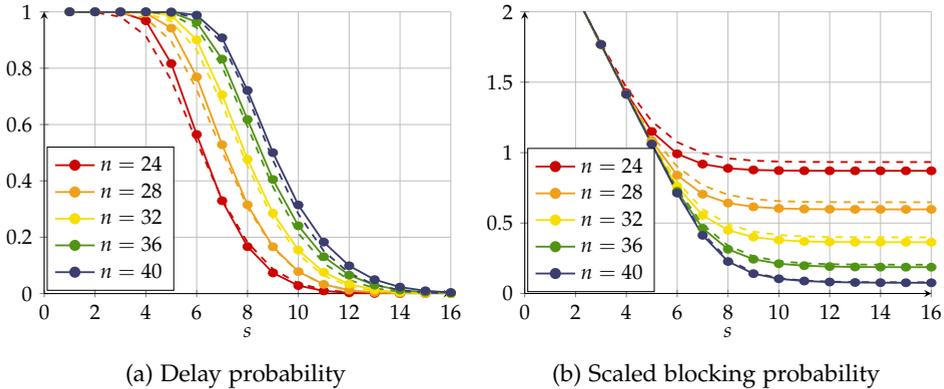

To test the accuracy of the asymptotic results in Theorem \ref{thm:limits_YT} as approximations in a realistic setting, we plot in Figure \ref{fig:accuracy_blocking} the exact probability of delay and blocking for an Erlang-R model with $R=8$ and $r=0.25$, as a function of $s$. The exact probabilities are given by Equation 
\eqref{eq:delay_probability}, and their respective asymptotic approximations are based on Theorem \ref{thm:limits_YT}. 
Despite the realistic moderate size of the system ($R=8$), we see that the QED approximations are remarkably accurate for many settings $(s,n)$. 
This fast relaxation is in line with observations made earlier in the QED literature \cite{Borst2004,Janssen2011}.
\begin{table}[htb]
\centering
\begin{tabular}{|r|rrrr|}
\hline
		&	$\mu$	&	$\d$	&	$p$	&	$r$ \\
\hline
Case 1  &   1   & 0.10 & 0.90 & 0.10 \\
Case 2 	& 	1	& 0.25 &  0.75 &  0.25\\
Case 3  &   1	& 0.50  &  0.50  &  0.50 \\
\hline
\end{tabular}
\caption{Parameter settings for numerical experiments.}
\label{tab:parameter_settings}
\end{table}

\begin{table}[h] \centering
\begin{tabular}{|r|rrr|rrr|}
\cline{2-7}\multicolumn{1}{r|}{} & \multicolumn{3}{c|}{$\beta = 1,\ \g = 1$} & \multicolumn{3}{c|}{$\beta = 1,\ \g = 2$} \bigstrut\\
\hline
$R_1$   & $\P({\rm d})$ & $\sqrt{R_1}\P({\rm b})$ & $\sqrt{R_1}\E[W]$ & $\P({\rm d})$ & $\sqrt{R_1}\P({\rm b})$ & $\sqrt{R_1}\E[W]$ \bigstrut\\
\hline
5     & 0.1270 & 0.0900 & 0.2283 & 0.1553 & 0.0212 & 0.1085 \bigstrut[t]\\
10    & 0.1340 & 0.0910 & 0.1919 & 0.1628 & 0.0206 & 0.1205 \\
25    & 0.1981 & 0.0945 & 0.1614 & 0.2356 & 0.0216 & 0.2145 \\
50    & 0.1513 & 0.0963 & 0.1588 & 0.1830 & 0.0205 & 0.1496 \\
100   & 0.1880 & 0.0956 & 0.1532 & 0.2231 & 0.0224 & 0.2055 \\
250   & 0.1797 & 0.0971 & 0.1399 & 0.2143 & 0.0219 & 0.2057 \\
\hline
\multicolumn{1}{r|}{} & \textit{0.1767} & \textit{0.0981} & \textit{0.1437} & \textit{0.2108} & \textit{0.0217} & \textit{0.1947} \bigstrut\\
\cline{2-7}\end{tabular}%
\vspace{5mm}
\begin{tabular}{|r|rrr|rrr|}
\cline{2-7}\multicolumn{1}{r|}{} & \multicolumn{3}{c|}{$\beta = 2,\ \g = 1$} & \multicolumn{3}{c|}{$\beta = 2,\ \g = 2$} \bigstrut\\
\hline
$R_1$   & $\P({\rm d})$ & $\sqrt{R_1}\P({\rm b})$ & $\sqrt{R_1}\E[W]$ & $\P({\rm d})$ & $\sqrt{R_1}\P({\rm b})$ & $\sqrt{R_1}\E[W]$ \bigstrut\\
\hline
5     & 0.0237 & 0.0868 & 0.0282 & 0.0322 & 0.0192 & 0.0391 \bigstrut[t]\\
10    & 0.0206 & 0.0872 & 0.0188 & 0.0278 & 0.0183 & 0.0264 \\
25    & 0.0277 & 0.0876 & 0.0123 & 0.0363 & 0.0174 & 0.0174 \\
50    & 0.0185 & 0.0913 & 0.0116 & 0.0249 & 0.0175 & 0.0166 \\
100   & 0.0232 & 0.0888 & 0.0103 & 0.0303 & 0.0183 & 0.0145 \\
250   & 0.0203 & 0.0905 & 0.0079 & 0.0267 & 0.0179 & 0.0109 \bigstrut[b]\\
\hline
\multicolumn{1}{r|}{} & \textit{0.0188} & \textit{0.0914} & \textit{0.0084} & \textit{0.0247} & \textit{0.0177} & \textit{0.0118} \bigstrut\\
\cline{2-7}\end{tabular}%

\caption{Exact numerical results for Erlang-R model with blocking for Case 1. The last row presents the asymptotic approximations.}
\label{tab:numerics_case1}
\end{table}
\begin{table}[h] \centering
\begin{tabular}{|r|rrr|rrr|}
\cline{2-7}\multicolumn{1}{r|}{} & \multicolumn{3}{c|}{$\beta = 1,\ \g = 1$} & \multicolumn{3}{c|}{$\beta = 1,\ \g = 2$} \bigstrut\\
\hline
$R_1$   & $\P({\rm d})$ & $\sqrt{R_1}\P({\rm b})$ & $\sqrt{R_1}\E[W]$ & $\P({\rm d})$ & $\sqrt{R_1}\P({\rm b})$ & $\sqrt{R_1}\E[W]$ \bigstrut\\
\hline
5     & 0.0911 & 0.1538 & 0.0479 & 0.1431 & 0.0345 & 0.0909 \bigstrut[t]\\
10    & 0.1010 & 0.1498 & 0.0560 & 0.1520 & 0.0326 & 0.1025 \\
25    & 0.1594 & 0.1509 & 0.1058 & 0.2192 & 0.0405 & 0.1785 \\
50    & 0.1201 & 0.1506 & 0.0726 & 0.1697 & 0.0381 & 0.1248 \\
100   & 0.1514 & 0.1539 & 0.1001 & 0.2088 & 0.0398 & 0.1704 \\
250   & 0.1459 & 0.1524 & 0.0957 & 0.2003 & 0.0397 & 0.1618 \\
\hline
\multicolumn{1}{r|}{} & \textit{0.1429} & \textit{0.1569} & \textit{0.0940} & \textit{0.1976} & \textit{0.0391} & \textit{0.1617} \bigstrut\\
\cline{2-7}\end{tabular}%

\vspace{5mm}

\begin{tabular}{|r|rrr|rrr|}
\cline{2-7}\multicolumn{1}{r|}{} & \multicolumn{3}{c|}{$\beta = 2,\ \g = 1$} & \multicolumn{3}{c|}{$\beta = 2,\ \g = 2$} \bigstrut\\
\hline
$R_1$   & $\P({\rm d})$ & $\sqrt{R_1}\P({\rm b})$ & $\sqrt{R_1}\E[W]$ & $\P({\rm d})$ & $\sqrt{R_1}\P({\rm b})$ & $\sqrt{R_1}\E[W]$ \bigstrut\\
\hline
5     & 0.0130 & 0.1484 & 0.0044 & 0.0277 & 0.0294 & 0.0109 \bigstrut[t]\\
10    & 0.0121 & 0.1432 & 0.0042 & 0.0244 & 0.0267 & 0.0098 \\
25    & 0.0182 & 0.1383 & 0.0070 & 0.0319 & 0.0295 & 0.0141 \\
50    & 0.0119 & 0.1415 & 0.0043 & 0.0216 & 0.0301 & 0.0090 \\
100   & 0.0154 & 0.1413 & 0.0059 & 0.0270 & 0.0290 & 0.0119 \\
250   & 0.0136 & 0.1403 & 0.0051 & 0.0236 & 0.0291 & 0.0103 \bigstrut[b]\\
\hline
\multicolumn{1}{r|}{\textit{}} & \textit{0.0126} & \textit{0.1445} & \textit{0.0048} & \textit{0.0220} & \textit{0.0284} & \textit{0.0097} \bigstrut\\
\cline{2-7}\end{tabular}%
\caption{Exact numerical results for Erlang-R model with blocking for Case 2. The last row presents the asymptotic approximations.}
\label{tab:numerics_case2}
\end{table}

\begin{table}[h] \centering
\begin{tabular}{|r|rrr|rrr|}
\cline{2-7}\multicolumn{1}{r|}{} & \multicolumn{3}{c|}{$\beta = 1,\ \g = 1$} & \multicolumn{3}{c|}{$\beta = 1,\ \g = 2$} \bigstrut\\
\hline
$R_1$   & $\P({\rm d})$ & $\sqrt{R_1}\P({\rm b})$ & $\sqrt{R_1}\E[W]$ & $\P({\rm d})$ & $\sqrt{R_1}\P({\rm b})$ & $\sqrt{R_1}\E[W]$ \bigstrut\\
\hline
5     & 0.0547 & 0.1945 & 0.0221 & 0.1181 & 0.0604 & 0.0617 \bigstrut[t]\\
10    & 0.0579 & 0.2158 & 0.0237 & 0.1325 & 0.0526 & 0.0746 \\
25    & 0.1113 & 0.2086 & 0.0544 & 0.1959 & 0.0641 & 0.1311 \\
50    & 0.0813 & 0.2050 & 0.0363 & 0.1523 & 0.0562 & 0.0933 \\
100   & 0.1060 & 0.2146 & 0.0509 & 0.1873 & 0.0632 & 0.1250 \\
250   & 0.1006 & 0.2179 & 0.0475 & 0.1820 & 0.0596 & 0.1214 \\
\hline
\multicolumn{1}{r|}{} & \textit{0.1011} & \textit{0.2185} & \textit{0.0478} & \textit{0.1792}& \textit{0.0605} & \textit{0.1199} \bigstrut\\
\cline{2-7}\end{tabular}%

\vspace{5mm}

\begin{tabular}{|r|rrr|rrr|}
\cline{2-7}\multicolumn{1}{r|}{} & \multicolumn{3}{c|}{$\beta = 2,\ \g = 1$} & \multicolumn{3}{c|}{$\beta = 2,\ \g = 2$} \bigstrut\\
\hline
$R_1$   & $\P({\rm d})$ & $\sqrt{R_1}\P({\rm b})$ & $\sqrt{R_1}\E[W]$ & $\P({\rm d})$ & $\sqrt{R_1}\P({\rm b})$ & $\sqrt{R_1}\E[W]$ \bigstrut\\
\hline
5     & 0.0034 & 0.1888 & 0.0009 & 0.0175 & 0.0510 & 0.0057 \bigstrut[t]\\
10    & 0.0030 & 0.2093 & 0.0008 & 0.0172 & 0.0416 & 0.0058 \\
25    & 0.0070 & 0.1937 & 0.0020 & 0.0243 & 0.0440 & 0.0089 \\
50    & 0.0043 & 0.1946 & 0.0011 & 0.0163 & 0.0414 & 0.0056 \\
100   & 0.0061 & 0.1999 & 0.0017 & 0.0207 & 0.0431 & 0.0076 \\
250   & 0.0052 & 0.2037 & 0.0014 & 0.0185 & 0.0401 & 0.0067 \bigstrut[b]\\
\hline
\multicolumn{1}{r|}{} & \textit{0.0052} & \textit{0.2039} & \textit{0.0014} & \textit{0.0173} & \textit{0.0404} & \textit{0.0063} \bigstrut\\
\cline{2-7}\end{tabular}%

\caption{Exact numerical results for Erlang-R model with blocking for Case 3. . The last row presents the asymptotic approximations.}
\label{tab:numerics_case3}
\end{table}

We furthermore compare the asymptotic delay and blocking probability in the three scenarios given in Table \ref{tab:parameter_settings}. 
In Tables \ref{tab:numerics_case1}--\ref{tab:numerics_case3} we compute the exact probabilities of delay and blocking through the explicit forms in \eqref{eq:delay_probability} for increasing values of the offered load, $R_1$.

The numerical results show that $g^b(\b,\g)$, $f^b(\b,\g)$ and $h^b(\b,\g)$ provide accurate approximations to $\P({\rm delay})$, $\sqrt{R_1}\P({\rm block})$ and $\sqrt{R_1}\,\E[W]$ in pre-limit systems.
The quality of the approximations increases with $R_1$.  
Naturally, fluctuations occur for relatively small values of $R_1$, because $s$ and $n$ need to be rounded to an integer. 

\subsection{QED limits for Erlang-R with holding}
\label{sec:QED_limit_holding}

As explained in Section \ref{sec:QED_scaling}, the model with holding has no product-form steady-state distribution, which makes it hard (if not impossible) to obtain QED limits.
Instead, we derive QED approximations by exploiting a connection with the blocking model. 

We first prove that under scaling \eqref{eq:twofoldscaling}, the upper bound on the utilization level of the nurses needed to achieve stability in the model with holding, as given in Proposition \ref{prop:StabilityCondition}, converges to unity as $R\to\infty$. 
This facilitates high utilization levels of both nurses and beds, a key characteristic of the QED regime.

\begin{proposition}\label{prop:stability_convergence}
Let $s$ and $n$ scale with $R_1\to\infty$ as in \eqref{eq:twofoldscaling}. Then for $\l\to\infty$,
\[
\rho_{\max}(s,n) \to 1.
\]
\end{proposition}

The proof can be found in Appendix \ref{app:proof_stability_convergence}. 
Combining Proposition \ref{prop:stability_convergence} with Proposition \ref{prop:StabilityCondition} shows that indeed the scaling we use results in a highly utilized system.

As observed before, the nature of the two variants of the model is similar up to the fact that a fraction of the patients is deferred on arrival in the setting with blocking, whereas all the arriving patients are eventually admitted into the system in the holding model. 
This implies that, given $s$ and $n$, the nurses face an increased workload in case of a holding room. 
In fact, Theorem \ref{thm:limits_YT} shows that the blocking probability is of order $1/\sqrt{R_1}$, yielding a volume of blocked patients of order $\sqrt{R_1}$ in setting with blocking. 
Accordingly, if $R^b = R_1$ and $R^h$ denote the nominal load arriving to the nurses in the model with blocking and holding, respectively, we can argue that 
\[R^h = R^b + \a \sqrt{R^b} + o(\sqrt{R^b}),\] 
for some $\a>0$. 
Notice that this additional load is of the same order as the safety staffing in the blocking model staffing rule \eqref{eq:twofoldscaling}. 
As $s$ and $n$ remain unchanged, we rewrite \eqref{eq:twofoldscaling} in terms of $R^h$,
\begin{align}
s &= R^h + (\b-\a)\sqrt{R^h} + o(\sqrt{R^h}), \nonumber \\
n &= \frac{R^h}{r} + \left(\g-\a/\sqrt{r}\right)\sqrt{\frac{R^h}{r}} + o(\sqrt{R^h}),
\label{eq:fixed_point_scaling}
\end{align} 
where we have used $R^b = O(R^h)$.
Observe that the square-root principle prevails also after this substitution, albeit with different hedging parameters.
We therefore heuristically argue that the holding model under scaling \eqref{eq:twofoldscaling} with parameters $\b$ and $\g$ mimics the blocking model with parameters $\b-\a$ and $\g-\a/\sqrt{r}$, respectively. 

Observe, however, that we have not yet specified the value of $\a$. 
By definition, $\a\sqrt{R^b}$ is the expected volume of patients that would be rejected in the model with blocking, that is, $R^h$ times the probability of not being admitted to the ED directly. 
By the construction in \eqref{eq:fixed_point_scaling}, this volume asymptotically equals $R^h \cdot \P^b({\rm block})$, with parameters $\b-\a$ and $\g-\a/\sqrt{r}$, which by Theorem \ref{thm:limits_YT} is approximated by
\[f^b\left(\b-\a,\gamma-\a/\sqrt{r}\right) / \sqrt{R^h}\]
as $R^h$ grows large.
In conclusion, $\a$ is characterized as the solution of the fixed-point equation
\begin{equation}
 \label{eq:fixedpoint}
 \a = f^h\left(\b-\a,\gamma-\a/\sqrt{r}\right),
 \end{equation} 
 and as a result, we are able to approximate the nurse delay probability in the Erlang-R model with holding as
 \begin{equation}
 \P^h({\rm delay}) \approx g^b(\b-\a,\g-\a/\sqrt{r}) =: g^h(\b,\g).
 \label{eq:fixed_point_Pwait}
 \end{equation}
Likewise, the scaled mean waiting time for a nurse can be approximated by
 \begin{equation}
 \sqrt{R_1} \cdot \E[W] \approx h^b(\b-\a,\g-\a/\sqrt{r}) =: h^h(\b,\g).
 \label{eq:fixed_point_Ewait}
 \end{equation}

This also implies that the holding queue is $O(\sqrt{R_1})$.
Subsequently, we argue that the expected holding time (pre-entering wait) under the QED policy is $O(1/\sqrt{R_1})$ and hence asymptotically negligible. 
We justify this claim numerically in Section \ref{sec:analysis_chapter5}. 

\begin{remark}
\label{rem:holding_limit}
Notice that in the reasoning leading to \eqref{eq:fixedpoint}, we implicitly assumed that the additional volume $\a\sqrt{R^b}$ is an independent Poisson process, which is obviously not the case. Therefore, \eqref{eq:fixed_point_Pwait}-\eqref{eq:fixed_point_Ewait} are approximations for pre-limit systems that are not asymptotically correct as $R_1\to\iy$.
Nevertheless, the heuristic approach seems to performs well as we confirm numerically next. 
\end{remark} 
 
\begin{figure}[h]
\centering
\begin{subfigure}{0.48\textwidth}
\begin{tikzpicture}[scale = 0.87]
\small
\begin{axis}[
	xmin = 0,
	xmax = 16,
	ymin = 0,
	ymax = 1,
	grid = both, 
	axis line style={->},
	axis lines = left,
	xlabel = $s$,
	xscale = 0.9,
	yscale = 0.75,
	legend cell align=left,
	legend style = {at = {(axis cs: 0.1,0.01)},anchor = south west}
]
 
\addplot[thick,col3,mark = *] table[x=s,y=sim32] {Chapter_5/tikz/accuracy/accuracy_pdelay_holding_case2.txt};
\addplot[thick,col4,mark = *] table[x=s,y=sim36] {Chapter_5/tikz/accuracy/accuracy_pdelay_holding_case2.txt};
\addplot[thick,col5,mark = *] table[x=s,y=sim40] {Chapter_5/tikz/accuracy/accuracy_pdelay_holding_case2.txt};

\addplot[thick,col3,dashed] table[x=s,y=approx32] {Chapter_5/tikz/accuracy/accuracy_pdelay_holding_case2.txt};
\addplot[thick,col4,dashed] table[x=s,y=approx36] {Chapter_5/tikz/accuracy/accuracy_pdelay_holding_case2.txt};
\addplot[thick,col5,dashed] table[x=s,y=approx40] {Chapter_5/tikz/accuracy/accuracy_pdelay_holding_case2.txt};

\legend{{$n=32$},{$n=36$},{$n=40$}};
\end{axis}
\end{tikzpicture}
\caption{Delay probability}
\label{fig:accuracy_holding_a}
\end{subfigure}
\begin{subfigure}{0.48\textwidth}
\begin{tikzpicture}[scale = 0.87]
\small
\begin{axis}[
	xmin = 0,
	xmax = 16,
	ymin = 0,
	ymax = 1,
	grid = both, 
	axis line style={->},
	axis lines = left,
	xlabel = $s$,
	xscale = 0.9,
	yscale = 0.75,
	legend cell align=left,
	legend style = {at = {(axis cs: 0.1,0.01)},anchor = south west}
]
 
\addplot[thick,col3,mark = *] table[x=s,y=sim32] {Chapter_5/tikz/accuracy/accuracy_holding_probability.txt};
\addplot[thick,col4,mark = *] table[x=s,y=sim36] {Chapter_5/tikz/accuracy/accuracy_holding_probability.txt};
\addplot[thick,col5,mark = *] table[x=s,y=sim40] {Chapter_5/tikz/accuracy/accuracy_holding_probability.txt};

\legend{{$n=32$},{$n=36$},{$n=40$}};
\end{axis}
\end{tikzpicture}
\caption{Holding probability}
\label{fig:accuracy_holding_b}
\end{subfigure}
\caption{Comparison of simulated delay probability (solid) against asymptotic approximations (dashed) with $\beta = (s-R_1)/\sqrt{R_1}$ and $\gamma = (n-R_1/r)/\sqrt{R_1/r}$ for $\lambda = 2$, $\mu=1$, $\delta=0.25$ and $p=0.75$.}
\label{fig:accuracy_holding}
\end{figure}
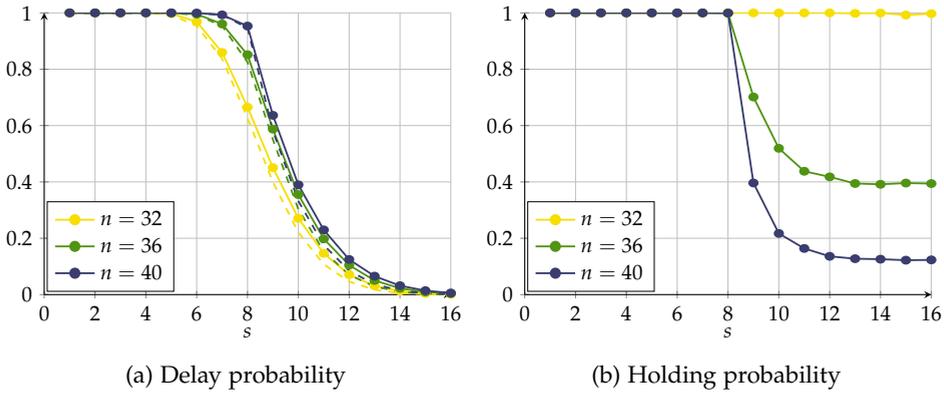

In Figure \ref{fig:accuracy_holding}, we repeat the numerical experiments of Figure \ref{fig:accuracy_blocking} for the model with holding. 
Since the heuristic does not provide an approximation for the holding probability, Figure \ref{fig:accuracy_holding_b} only plots the simulated holding probabilities. 
Those are provided to better understand the implication of operational decisions.
Recall that the holding system is only stable (i.e. $\P({\rm hold})<1$) if both $s>R_1=8$ and $n > R_1/r  = 32$.
We nevertheless included the boundary case $n=32$ for illustrative purposes. 
The graphs in Figure \ref{fig:accuracy_holding} show that the heuristic captures the congestion levels well, even for this moderate-size system.

To see how this heuristic approach performs under different settings, and particularly if $R_1\to \infty$, we again compare the approximated delay probability in the Erlang-R model with holding as solution of the fixed-point procedure to the outcomes of simulation experiments for the three scenarios in Table \ref{tab:parameter_settings}. 
We performed 100 runs of length $10^4$ for each parameter setting and all values of $R$, yielding the results presented in Tables \ref{tab:heuristic_case1}--\ref{tab:heuristic_case3}, which are accurate up to a 95\% confidence interval of width $10^{-3}$.

\begin{table}[h] \centering
\begin{tabular}{|r|rr|rr|}
\cline{2-5}\multicolumn{1}{r|}{} & \multicolumn{2}{c|}{$\b=1,\ \g=1$} & \multicolumn{2}{c|}{$\b=1,\ \g=2$} \bigstrut\\
\hline
$R_1$     & $\P({\rm d})$ & $\sqrt{R_1}\E[W]$ & $\P({\rm d})$ & $\sqrt{R_1}\E[W]$ \bigstrut\\
\hline
5     & 0.1532 & 0.1031 & 0.1628 & 0.1216 \bigstrut[t]\\
10    & 0.1622 & 0.1272 & 0.1697 & 0.1331 \\
25    & 0.2340 & 0.2116 & 0.2413 & 0.2342 \\
50    & 0.1817 & 0.1468 & 0.1890 & 0.1678 \\
100   & 0.2199 & 0.1931 & 0.2304 & 0.2269 \\
250   & 0.2110 & 0.1852 & 0.2176 & 0.2230 \bigstrut[b]\\
\hline
\multicolumn{1}{r|}{} & \textit{0.2076} & \textit{0.1777} & \textit{0.2187} & \textit{0.2050} \bigstrut\\
\cline{2-5}\end{tabular}%

\vspace{5mm}

\begin{tabular}{|r|rr|rr|}
\cline{2-5}\multicolumn{1}{c|}{} & \multicolumn{2}{c|}{$\b=2,\ \g=1$} & \multicolumn{2}{c|}{$\b=2,\ \g=1$} \bigstrut\\
\hline
$R_1$     & $\P({\rm d})$ & $\sqrt{R_1}\E[W]$ & $\P({\rm d})$ & $\sqrt{R_1}\E[W]$ \bigstrut\\
\hline
5     & 0.0310 & 0.0121 & 0.0344 & 0.0148 \bigstrut[t]\\
10    & 0.0267 & 0.0123 & 0.0292 & 0.0128 \\
25    & 0.0348 & 0.0171 & 0.0373 & 0.0184 \\
50    & 0.0240 & 0.0108 & 0.0258 & 0.0125 \\
100   & 0.0293 & 0.0143 & 0.0317 & 0.0163 \\
250   & 0.0256 & 0.0120 & 0.0276 & 0.0145 \\
\hline
\multicolumn{1}{r|}{\textit{}} & \textit{0.0229} & \textit{0.0104} & \textit{0.0257} & \textit{0.0124} \bigstrut\\
\cline{2-5}\end{tabular}%

\caption{Simulated probability of delay and scaled expected waiting time in Erlang-R model with holding for Case 1. The last row gives the asymptotic approximations.}
\label{tab:heuristic_case1}
\end{table}
\begin{table}[h]\centering

\begin{tabular}{|r|rr|rr|}
\cline{2-5}\multicolumn{1}{r|}{} & \multicolumn{2}{c|}{$\b=1,\ \g=1$} & \multicolumn{2}{c|}{$\b=1,\ \g=2$} \bigstrut\\
\hline
$R_1$     & $\P({\rm d})$ & $\sqrt{R_1}\E[W]$ & $\P({\rm d})$ & $\sqrt{R_1}\E[W]$ \bigstrut\\
\hline
5     & 0.1327 & 0.0740 & 0.1620 & 0.1096 \bigstrut[t]\\
10    & 0.1446 & 0.0894 & 0.1683 & 0.1207 \\
25    & 0.2204 & 0.1631 & 0.2442 & 0.2203 \\
50    & 0.1694 & 0.1122 & 0.1888 & 0.1507 \\
100   & 0.2098 & 0.1524 & 0.2322 & 0.2111 \\
250   & 0.2033 & 0.1534 & 0.2190 & 0.1979 \bigstrut[b]\\
\hline
\multicolumn{1}{r|}{} & \textit{0.1840} & \textit{0.1277} & \textit{0.2109} & \textit{0.1759} \bigstrut\\
\cline{2-5}\end{tabular}%
\vspace{5mm}
\begin{tabular}{|r|rr|rr|}
\cline{2-5}\multicolumn{1}{c|}{} & \multicolumn{2}{c|}{$\b=2,\ \g=1$} & \multicolumn{2}{c|}{$\b=2,\ \g=1$} \bigstrut\\
\hline
$R_1$     & $\P({\rm d})$ & $\sqrt{R_1}\E[W]$ & $\P({\rm d})$ & $\sqrt{R_1}\E[W]$ \bigstrut\\
\hline
5     & 0.0219 & 0.0079 & 0.0322 & 0.0137 \bigstrut[t]\\
10    & 0.0199 & 0.0073 & 0.0284 & 0.0115 \\
25    & 0.0283 & 0.0128 & 0.0375 & 0.0163 \\
50    & 0.0190 & 0.0078 & 0.0255 & 0.0107 \\
100   & 0.0244 & 0.0097 & 0.0314 & 0.0151 \\
250   & 0.0214 & 0.0083 & 0.0272 & 0.0134 \bigstrut[b]\\
\hline
\multicolumn{1}{r|}{\textit{}} & \textit{0.0169} & \textit{0.0066} & \textit{0.0234} & \textit{0.0104} \bigstrut\\
\cline{2-5}\end{tabular}%

\caption{Simulated probability of delay and scaled expected waiting time in Erlang-R model with holding for Case 2. The last row gives the asymptotic approximations.}
\label{tab:heuristic_case2}
\end{table}

\begin{table}
\centering
\begin{tabular}{|r|rr|rr|}
\cline{2-5}\multicolumn{1}{r|}{} & \multicolumn{2}{c|}{$\b=1,\ \g=1$} & \multicolumn{2}{c|}{$\b=1,\ \g=2$} \bigstrut\\
\hline
$R_1$     & $\P({\rm d})$ & $\sqrt{R_1}\E[W]$ & $\P({\rm d})$ & $\sqrt{R_1}\E[W]$ \bigstrut\\
\hline
5     & 0.0977 & 0.0413 & 0.1521 & 0.0851 \bigstrut[t]\\
10    & 0.1070 & 0.0469 & 0.1648 & 0.1028 \\
25    & 0.1926 & 0.1076 & 0.2421 & 0.1874 \\
50    & 0.1431 & 0.0727 & 0.1876 & 0.1342 \\
100   & 0.1855 & 0.1012 & 0.2282 & 0.1714 \\
250   & 0.1775 & 0.0963 & 0.2217 & 0.1765 \bigstrut[b]\\
\hline
\multicolumn{1}{r|}{} & \textit{0.1442} & \textit{0.0711} & \textit{0.1981} & \textit{0.1354} \bigstrut\\
\cline{2-5}\end{tabular}%

\vspace{5mm}
\begin{tabular}{|r|rr|rr|}
\cline{2-5}\multicolumn{1}{c|}{} & \multicolumn{2}{c|}{$\b=2,\ \g=1$} & \multicolumn{2}{c|}{$\b=2,\ \g=2$} \bigstrut\\
\hline
$R_1$     & $\P({\rm d})$ & $\sqrt{R_1}\E[W]$ & $\P({\rm d})$ & $\sqrt{R_1}\E[W]$ \bigstrut\\
\hline
5     & 0.0072 & 0.0019 & 0.0250 & 0.0081 \bigstrut[t]\\
10    & 0.0067 & 0.0018 & 0.0235 & 0.0082 \\
25    & 0.0148 & 0.0043 & 0.0325 & 0.0133 \\
50    & 0.0092 & 0.0025 & 0.0217 & 0.0081 \\
100   & 0.0132 & 0.0038 & 0.0277 & 0.0105 \\
250   & 0.0114 & 0.0033 & 0.0246 & 0.0099 \bigstrut[b]\\
\hline
\multicolumn{1}{r|}{\textit{}} & \textit{0.0078} & \textit{0.0022} & \textit{0.0188} & \textit{0.0069} \bigstrut\\
\cline{2-5}\end{tabular}%

\caption{Simulated probability of delay and scaled expected waiting time in Erlang-R model with holding for Case 3. The last row gives the asymptotic approximations.}
\label{tab:heuristic_case3}
\end{table}

We conclude from these tables that the approximation is good. As $R$ increases, the simulated values move closer to the heuristic approximation. Extensive numerical experiments suggest that load is slightly underestimated in the limit. 
The best results in terms of accuracy are attained for small $r$. 
This suggests that the quality of the heuristic method improves as $r$ gets smaller. 
These are exactly the parameter settings for which this model is relevant. 

\begin{remark}
The approximation technique that evolves around the fixed-point\\ \noindent method can be adapted to accommodate balking behavior of external arrivals. If we assume that an arriving patient finding all beds occupied leaves the system instantly with probability $1-q$, for some $q\in(0,1)$, independently of the rest of the arrivals, with the same argumentation, the volume of arrivals blocked is still $\a\sqrt{R_1}$, while the volume that will enter the ED eventually is $q\cdot\a\sqrt{R_1}$. Therefore, we may argue that in the QED regime, the system with holding and balking behaves as the system with blocking but with corrected parameters $(\b-q\a,\g-q\a/\sqrt{r})$, where $\a$ satisfies
\begin{equation}
\a = f^b(\b-q\a,\g-q\a/\sqrt{r}).
\end{equation}
Note that the choice of $q$ interpolates between the two system variants with holding ($q=0$) and blocking ($q=1$).
\end{remark}

\section{Dimensioning}
\label{sec:dimensioning}

We will now use the accurate asymptotic approximations of the previous section to define a procedure that determines resource capacity in the restricted Erlang-R models. 
That is, we aim to set the number of nurses $s$ and the number of beds $n$, such that a preset performance level is achieved. 
We take the probability of delay at the needy queue and the probability of blocking/holding at the pre-entrant queue as the target performance objectives.

\subsection{Capacity setting for Erlang-R with blocking}
\label{sec:dimensioning_block}
In the setting with blocking, we can readily use the asymptotic results of Theorem \ref{thm:limits_YT} to (numerically) find a pair of parameters $(\b^*,\g^*)$ to meet the performance requirements.
For instance, given that we want the delay probability to be at most $\varepsilon$, we first solve the equation $g^b(\b^*,\g^*)=\varepsilon$ and then assign $s = \lceil R_1 + \b^*\sqrt{R_1}\rceil$ and $n = \lceil R_1/r+\g^*\sqrt{R_1/r}\rceil$.  Note that there could be multiple solutions to that problem, i.e.\ there could be multiple combinations of number of beds and number of nurses that can result in the same value of a single performance level.
The ED manager can ultimately decide which of these feasible solutions fits the environment best, for instance taking into account space and cost constraints.

We illustrate the resource allocation decisions in an MU setting, using data originated from two articles: \cite{LS2001} and \cite{GY2011}. Green \& Yankovic describe an MU that has 42 beds, with average occupancy level 78\%, and Average Length of Stay (ALOS) 4.3 days. Lundgren \& Segesten studied nurses' service times in a medical-surgical ward. They found that the average service time in their unit was 15.3 minutes per service, and that the average demand rate for each patient is 0.38 requests per hour. Therefore, we take an average service time of 15 minutes and assume that there are 0.4 requests per hour from each patient. Fitting this data to our model results in the following parameters values: $\lambda =0.32, \mu =4, \delta =0.4$, $p=0.975$ and the fraction of needy time is then approximately $r=0.09$. 
This yields nominal offered load $R_1 = 3.2$ and $R_1/r = 34.4$.

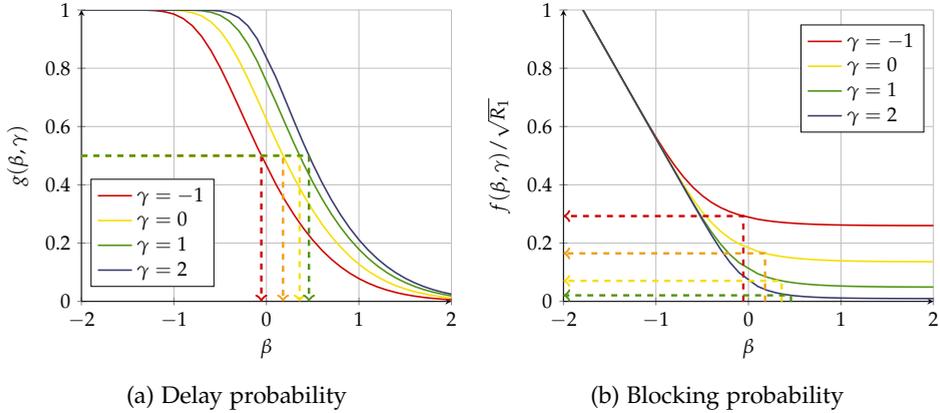
\begin{figure}
\centering
\begin{subfigure}{0.48\textwidth}
\begin{tikzpicture}[scale = 0.75]
\begin{axis}[
	xmin = -2,
	xmax = 2,
	ymin = 0,
	ymax = 1,
	grid = both, 
	axis line style={->},
	axis lines = left,
	xlabel = $\beta$,
	ylabel = {$g(\beta,\gamma)$},
	y label style = {at = {(axis cs: -2.5,0.5)}},
	legend cell align=left,
	legend style = {at = {(axis cs: -1.9,0.05)},anchor = south west},
	xscale = 0.95,
	yscale = 0.9
]

\addplot[thick,col1] table[x=beta,y=delay_gmin1] {Chapter_5/tikz/staffing_example/staffing_example_with_blocking1.txt};
\addplot[thick,col3] table[x=beta,y=delay_g0] {Chapter_5/tikz/staffing_example/staffing_example_with_blocking1.txt};
\addplot[thick,col4] table[x=beta,y=delay_g1] {Chapter_5/tikz/staffing_example/staffing_example_with_blocking1.txt};
\addplot[thick,col5] table[x=beta,y=delay_g2] {Chapter_5/tikz/staffing_example/staffing_example_with_blocking1.txt};

\draw[->,col1,very thick,dashed] (axis cs: -2,0.5) -- (axis cs: -0.0552366,0.5) -- (axis cs: -0.0552366,0);
\draw[->,col2,very thick,dashed] (axis cs: -2,0.5) -- (axis cs: 0.179728,0.5) -- (axis cs: 0.179728,0);
\draw[->,col3,very thick,dashed] (axis cs: -2,0.5) -- (axis cs: 0.359034,0.5) -- (axis cs: 0.359034,0);
\draw[->,col4,very thick,dashed] (axis cs: -2,0.5) -- (axis cs: 0.459825,0.5) -- (axis cs: 0.459825,0);

\legend{$\gamma = -1$, $\gamma =0$,$\gamma=1$, $\gamma=2$}
\end{axis}

\end{tikzpicture}
\caption{Delay probability}
\label{fig:ratio01_delay}
\end{subfigure}
\begin{subfigure}{0.48\textwidth}
\begin{tikzpicture}[scale = 0.75]
\begin{axis}[
	xmin = -2,
	xmax = 2,
	ymin = 0,
	ymax = 1,
	grid = both, 
	axis line style={->},
	axis lines = left,
	xlabel = $\beta$,
	ylabel = {$f(\beta,\gamma)/\sqrt{R_1}$},
	y label style = {at = {(axis cs: -2.5,0.5)}},
	legend cell align=left,
	legend style = {at = {(axis cs: 1.9,0.95)},anchor = north east},
	xscale = 0.95,
	yscale = 0.9
]

\addplot[thick,col1] table[x=beta,y=delay_gmin1] {Chapter_5/tikz/staffing_example/staffing_example_with_blocking2.txt};
\addplot[thick,col3] table[x=beta,y=delay_g0] {Chapter_5/tikz/staffing_example/staffing_example_with_blocking2.txt};
\addplot[thick,col4] table[x=beta,y=delay_g1] {Chapter_5/tikz/staffing_example/staffing_example_with_blocking2.txt};
\addplot[thick,col5] table[x=beta,y=delay_g2] {Chapter_5/tikz/staffing_example/staffing_example_with_blocking2.txt};

\draw[very thick, col1,dashed,->] (axis cs: -0.0552,0) -- (axis cs: -0.0552,0.292798) -- (axis cs: -2,0.292798);
\draw[very thick,col2,dashed,->] (axis cs: 0.179728,0)  -- (axis cs: 0.179728,0.164903)   -- (axis cs: -2,0.164903);
\draw[very thick,col3,dashed,->] (axis cs:  0.359034,0)  -- (axis cs:  0.359034,0.0705547)   -- (axis cs: -2,0.0705547);
\draw[very thick,col4,dashed,->] (axis cs: 0.459825,0)  -- (axis cs: 0.459825,0.0207909)   -- (axis cs: -2,0.0207909);

\legend{$\gamma = -1$, $\gamma =0$,$\gamma= 1$, $\gamma=2$}
\end{axis}

\end{tikzpicture}
\caption{Blocking probability}
\label{fig:ratio01_block}
\end{subfigure}
\caption{Approximate performance of restricted Erlang-R with blocking for $r \approx 0.09$ and $R_1 = 3.2$, as functions of $\beta$.}
\label{fig:ratio01}
\end{figure}

Figure \ref{fig:ratio01} visualizes the dimensioning procedure for this particular MU. 
The hospital management can find a pair of $n$ and $s$ to meet certain criteria, for example to achieve target delay probability $\varepsilon = 0.5$ with reasonable blocking probability.
Figure \ref{fig:ratio01}a indicates that this target can be achieved by a variety of pairs, for instance $(\beta_1,\gamma_1) = (-0.06,-1)$, $(\beta_2,\gamma_2) = (0.16,0)$, $(\beta_3,\gamma_3) = (0.36,1)$ or $(\beta_4,\gamma_4) = (0.46,2)$, among infinitely many others.
According to Figure \ref{fig:ratio01}b, the pairs named above lead to blocking probabilities 0.293, 0.165, 0.071 and 0.021, respectively. 
If the manager decides that probability of blocking of more than 10 percent is not acceptable, this leaves the choices $(\beta_3,\gamma_3) = (0.36,1)$ or $(\beta_4,\gamma_4) = (0.46,2)$ as candidate parameter pairs.
Using the two-fold square-root staffing rule $s_i = \lceil R_1 + \beta_i \sqrt{R_1}\rceil$ and $n_i = [R_1/r + \gamma_i\sqrt{R_1/r}]$, this yields feasible staffing levels $(s_3,n_3) = (4,40)$ and $(s_4,n_4)=(5,46)$. 
The ultimate decision to apply any of these solutions can be based on external factors, such as operational costs or space limitations on the number of beds.

\subsection{Capacity setting for Erlang-R with holding}
For the holding model, we need a more sophisticated approach, exploiting the asymptotic approximation with the fixed-point equation in \eqref{eq:fixedpoint}.  We propose the following dimensioning procedure to achieve a preset target delay probability at the needy queue.

\begin{algorithm}
\hspace{1cm}\rule{10cm}{1pt}\\
\hspace{1.1cm}\KwIn{Target delay probability $\varepsilon$. Parameters $\l,\m,\d$ and $p$.}
\hspace{1.1cm}\KwOut{Staffing levels $s$ and $n$.}
\vspace{-1mm}
\hspace{1cm}\rule{10cm}{0.5pt}\\
\vspace{-1mm}
\begin{enumerate}
\item[] \hspace{0.5cm} 1. Set $R_1:= \frac{\l}{(1-p)\m}$ and $r = \frac{\d}{\d+p\m}$. 
\item[] \hspace{0.5cm} 2. Determine parameters $(\b^*,\g^*)$ such that $g^b(\b^*,\g^*) = \varepsilon$. 
\item[] \hspace{0.5cm} 3. Set $\b = \b^* + f^b(\b^*,\g^*)$ and $\g = \g^* + f^b(\b^*,\g^*)/\sqrt{r}$.
\item[] \hspace{0.5cm} 4. Return $s = \left\lceil R_1 + \b\sqrt{R_1}\right\rceil$ and $n = \left\lfloor R_1/r + \g \sqrt{R_1/r}\right\rfloor$.
\end{enumerate}
\vspace{-3 mm}
\hspace{1cm}\rule{10cm}{1pt}\\
\vspace{2 mm}
\caption{Stationary dimensioning algorithm for ED with holding.}
\label{alg:stationarydimensioning}
\end{algorithm}

\begin{remark}\label{rem:upperboundHW}
In Step 2 of Algorithm \ref{alg:stationarydimensioning} infinitely many pairs $(\b^*,\g^*)$ satisfy the delay probability equation. 
For practical purposes, it is convenient to fix either $\b^*$ or $\g^*$ beforehand, and then solve $g^b(\b^*,\g^*) = \varepsilon$ for the remaining unknown. 
The preset value should however be chosen with care, since $g^b(\b^*,\g^*)$ is upper bounded by the Halfin-Whitt delay probability formula 
\[g_{\rm HW}(\b^*) = \left( 1 + \frac{\b^* \Phi(\b^*)}{\f(\b^*)}\right)^{-1}.\]
Hence, if $\varepsilon > g_{\rm HW}(\b^*)$, then no feasible solution to $g^b(\b^*,\g^*)=\varepsilon$ exists.
This should be considered when choosing $\beta^*$. 
Furthermore, it is evident from Step 3 that the final values $(\b,\g)$ are always larger than $(\b^*,\g^*)$. 
\end{remark}


We now use the same example as in Section \ref{sec:dimensioning_block} to demonstrate capacity allocation decisions for the model with holding. This can be viewed as the additional capacity the medical unit needs in terms of nurses and beds, in order to account for the fact that patients are waiting in the ED to be admitted instead of being blocked and transferred to a less preferred unit. 
Observe that the holding model leaves less flexibility for management in choosing system parameters due to stability constraints. For example, the policy with $n=30$ ($\gamma=-0.75$) is infeasible in the holding model.
For similar reasons, only nurse staffing levels with $\beta>0$, or $s > R_1=3.2$ are feasible.

\begin{figure}[h]
\centering
\begin{tikzpicture}[scale = 0.8]
\begin{axis}[
	xmin = 0,
	xmax = 2,
	ymin = 0,
	ymax = 1,
	grid = both, 
	axis line style={->},
	axis lines = left,
	xlabel = $\beta$,
	ylabel = {$g^h(\beta,\gamma)$},
	legend cell align=left,
	legend style = {at = {(axis cs: 1.9,0.95)},anchor = north east},
	xscale = 0.95,
	yscale = 0.9
]

\addplot[thick,col5] table[x=beta,y=delay_n35] {Chapter_5/tikz/staffing_example/staffing_example_with_holding1.txt};
\addplot[thick,col4] table[x=beta,y=delay_n40] {Chapter_5/tikz/staffing_example/staffing_example_with_holding1.txt};
\addplot[thick,col1] table[x=beta,y=delay_n45] {Chapter_5/tikz/staffing_example/staffing_example_with_holding1.txt};

\draw[very thick,col4,dashed,->] (axis cs: -2,0.5) -- (axis cs: 0.475,0.5) -- (axis cs: 0.475,0);

\legend{$\gamma = -0.75$, $\gamma =0.102$,$\gamma= 0.955$, $\gamma=1.807$};

\end{axis}
\end{tikzpicture}
\caption{Approximate delay probability of restricted Erlang-R system with holding for $r\approx 0.09$ and $R_1=3.2$ }
\label{fig:ratio01_hold}
\end{figure}
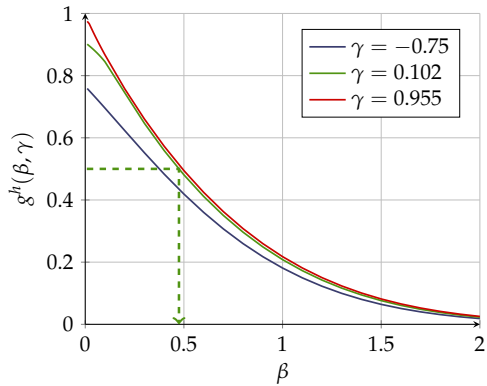

Targeting a delay probability of $0.5$ with $n=40$, Figure \ref{fig:ratio01_hold} shows that operating a MU with holding room requires $\beta = 0.475$ or $s=5$. 
Recall that under the blocking policy, only $s=4$ nurses were needed to achieve a delay probability of $0.5$. 
This example hence shows how the managerial decision to have a holding room, rather than deferring patients to less preferred medical units, requires additional workforce in that unit (as well as the ED). 
This example also shows that the facility with holding room is able to treat fewer patients simultaneously than under blocking constraints, in line with the bounds in Section \ref{sec:bounds} and Conjecture \ref{conj:stochorder}.

\section{Model analysis and managerial implications}
\label{sec:analysis_chapter5}
In this section, we use the analysis and algorithms developed in earlier sections to gain insights into the importance of the capacity restrictions and patient returns in a restricted Erlang-R system by drawing a comparison to related models studied in the literature. 
\subsection{The influence of patient returns or the role of $r$}
Here we study how the parameter $r$ affects the service level in the restricted Erlang-R model with blocking, on the basis of the asymptotic expressions in Theorem \ref{thm:limits_YT}.


To better understand the connection with the single-station model and the importance of returns we examine the role of $r$. 
Recall the interpretation of $r$ as the fraction of time a patient is needy during his stay within the system in the idealized scenario with infinite capacity, i.e. for $r\in(0,1)$. 
The case $r=1$ corresponds to the setting in which patients are needy all the time, in this case patients get service in one time.
When $r=1$ the infinite-server queue, describing the number of content patients, disappears from the queueing system and we end up with a standard loss model---$M/M/s/n$ queue---in which capacity is scaled as
\[ s = R_1+\beta\sqrt{R_1}, \qquad n = R_1+\gamma\sqrt{R_1}. \]
This staffing rule only makes sense in case $\beta<\gamma$, since no delay is experienced if $n\leq s$. 
If indeed $\gamma>\beta$, then the asymptotic delay probability and scaled blocking probability are given by \cite{masseywallace},
\begin{align*}
g_B(\beta,\gamma) &= \frac{1-{\rm e}^{-\beta(\gamma-\beta)}}{1-{\rm e}^{-\beta(\gamma-\beta)}+\beta\Phi(\beta)/\f(\beta)}, \\
f_B(\beta,\gamma) &= \frac{\beta{\rm e}^{-\beta(\gamma-\beta)}}{1-{\rm e}^{-\beta(\gamma-\beta)}+\beta\Phi(\beta)/\f(\beta)}. 
\end{align*}

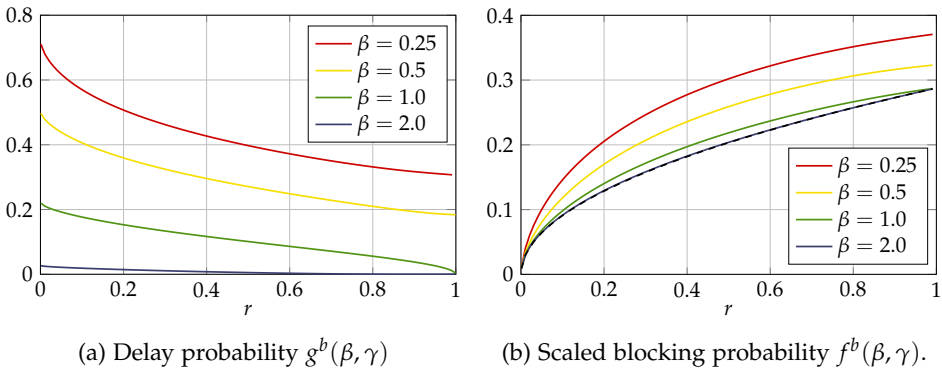
\begin{figure}
\centering
\begin{subfigure}{0.48\textwidth}
\centering
\begin{tikzpicture}[scale=0.8]
\begin{axis}[
	xmin = 0,
	xmax = 1,
	ymin = 0,
	ymax = 0.8,
	grid = both,
	xlabel = $r$,
	axis line style={->},
	legend cell align=left,
	legend style={at={(0.98,1.28)},anchor= north east},
	yscale = 0.75
]

\addplot[thick,col1] file {Chapter_5/tikz/influence_r/PdelayB_g1_b025.txt};
\addplot[thick,col3] file {Chapter_5/tikz/influence_r/PdelayB_g1_b05.txt};
\addplot[thick,col4] file {Chapter_5/tikz/influence_r/PdelayB_g1_b1.txt};
\addplot[thick,col5] file {Chapter_5/tikz/influence_r/PdelayB_g1_b2.txt};

\legend{$\beta = 0.25$,$\beta=0.5$,$\beta=1.0$,$\beta=2.0$}
\end{axis}
\end{tikzpicture}
\caption{Delay probability $g^b(\beta,\gamma)$}
\end{subfigure}
\begin{subfigure}{0.48\textwidth}
\centering
\begin{tikzpicture}[scale=0.8]
\begin{axis}[
	xmin = 0,
	xmax = 1,
	ymin = 0,
	ymax = 0.4,
	grid = both,
	xlabel = $r$,
	axis line style={->},
	legend cell align=left,
	legend style={at={(0.98,0.05)},anchor= south east},
	yscale = 0.75
]

\addplot[thick,col1] file {Chapter_5/tikz/influence_r/PblockB_g1_b025.txt};
\addplot[thick,col3] file {Chapter_5/tikz/influence_r/PblockB_g1_b05.txt};
\addplot[thick,col4] file {Chapter_5/tikz/influence_r/PblockB_g1_b1.txt};
\addplot[thick,col5] file {Chapter_5/tikz/influence_r/PblockB_g1_b2.txt};
\addplot[thick,dashed] file {Chapter_5/tikz/influence_r/PblockB_g1_inf.txt};
\legend{$\beta = 0.25$,$\beta=0.5$,$\beta=1.0$,$\beta=2.0$}
\end{axis}
\end{tikzpicture}
\caption{Scaled blocking probability $f^b(\beta,\gamma)$.}
\label{fig:influence_of_r_b}
\end{subfigure}
\caption{Asymptotic performance measures as a function of $r$ in the restricted Erlang-R model with blocking for $\gamma=1$.}
\label{fig:influence_of_r}
\end{figure}

We can see that $f^b(\beta,\gamma)$ for increasing $\b$ approaches a lower bound that is a function of $r$. 
To understand this, observe that as $\beta$ grows, delays at the nurse queue vanish. 
Then the sojourn time of an admitted patient only consists of a geometric number of needy and content periods with mean $(1/\mu+p/\delta)/(1-p) = 1/(r\mu(1-p))$. 
The blocking model can in this case be modeled as an $M/G/n/n$ queue, with offered load $\l/(r\mu(1-p)) =R_1/r$, in which the scaled blocking probability is known to be, see \cite{Avram2013},
\[\sqrt{R_1} \, \P({\rm block}) = \sqrt{R_1} \, \frac{(R_1/r)^n/n!}{\sum_{k=0}^n (R_1/r)^k / k!} \to \sqrt{r} \, \frac{\f(\gamma)}{\Phi(\gamma)},\]
as $R_1\to\infty$. 
This function of $r$ is plotted in Figure \ref{fig:influence_of_r_b} as the dashed line.

We observe that in general the probability of blocking increases with $r$, regardless of the capacity constraints on the needy station. 
We can explain this by observing that $r$ influences only $n$ in the QED staffing rule. When $n$ reduces, more patients are blocked. Therefore, if patients spend relatively more time in needy state, which usually indicates services that are less interrupted, blocking will increase. Delays, on the other hand, will decrease in such situations---the minimal delay possible can be achieved if service is given in one time ($r=1$). Returns or interruptions increase delays significantly under QED staffing.

\subsection{Comparing restricted and unrestricted Erlang-R models}

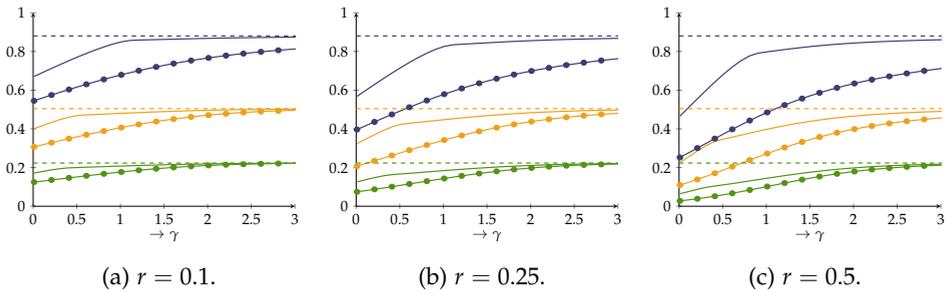
\begin{figure}
\centering
\begin{subfigure}{0.32\textwidth}
\begin{tikzpicture}[scale = 0.56]
\begin{axis}[
	xmin = 0,
	xmax = 3,
	ymin = 0,
	ymax = 1,
	axis line style={->},
	axis lines = left,
	xlabel = $\to \gamma$,
	xscale = 0.9,
	yscale = 0.8
]

\addplot[thick,col5,dashed] file {Chapter_5/tikz/comparison/PdelayR_b01.txt};
\addplot[thick,col5,mark=*, mark repeat = 2] file {Chapter_5/tikz/comparison/PdelayB_r01_b01.txt};
\addplot[thick,col5] file {Chapter_5/tikz/comparison/PdelayH_r01_b01.txt};

\addplot[thick,col2,dashed] file {Chapter_5/tikz/comparison/PdelayR_b05.txt};
\addplot[thick,col2,mark=*, mark repeat = 2] file {Chapter_5/tikz/comparison/PdelayB_r01_b05.txt};
\addplot[thick,col2] file {Chapter_5/tikz/comparison/PdelayH_r01_b05.txt};

\addplot[thick,col4,dashed] file {Chapter_5/tikz/comparison/PdelayR_b1.txt};
\addplot[thick,col4,mark=*, mark repeat = 2] file {Chapter_5/tikz/comparison/PdelayB_r01_b1.txt};
\addplot[thick,col4] file {Chapter_5/tikz/comparison/PdelayH_r01_b1.txt};

\end{axis}
\end{tikzpicture}
\caption{$r=0.1$.}
\end{subfigure}
\begin{subfigure}{0.32\textwidth}
\begin{tikzpicture}[scale = 0.56]
\begin{axis}[
	xmin = 0,
	xmax = 3,
	ymin = 0,
	ymax = 1,
	axis line style={->},
	axis lines = left,
	xlabel = $\to \gamma$,
	xscale = 0.9,
	yscale = 0.8
]

\addplot[thick,col5,dashed] file {Chapter_5/tikz/comparison/PdelayR_b01.txt};
\addplot[thick,col5,mark=*, mark repeat = 2] file {Chapter_5/tikz/comparison/PdelayB_r025_b01.txt};
\addplot[thick,col5] file {Chapter_5/tikz/comparison/PdelayH_r025_b01.txt};

\addplot[thick,col2,dashed] file {Chapter_5/tikz/comparison/PdelayR_b05.txt};
\addplot[thick,col2,mark=*, mark repeat = 2] file {Chapter_5/tikz/comparison/PdelayB_r025_b05.txt};
\addplot[thick,col2] file {Chapter_5/tikz/comparison/PdelayH_r025_b05.txt};

\addplot[thick,col4,dashed] file {Chapter_5/tikz/comparison/PdelayR_b1.txt};
\addplot[thick,col4,mark=*, mark repeat = 2] file {Chapter_5/tikz/comparison/PdelayB_r025_b1.txt};
\addplot[thick,col4] file {Chapter_5/tikz/comparison/PdelayH_r025_b1.txt};

\end{axis}
\end{tikzpicture}
\caption{$r=0.25$.}
\end{subfigure}
\begin{subfigure}{0.32\textwidth}
\begin{tikzpicture}[scale = 0.56]
\begin{axis}[
	xmin = 0,
	xmax = 3,
	ymin = 0,
	ymax = 1,
	axis line style={->},
	axis lines = left,
	xlabel = $\to \gamma$,
	xscale = 0.9,
	yscale = 0.8
]

\addplot[thick,col5,dashed] file {Chapter_5/tikz/comparison/PdelayR_b01.txt};
\addplot[thick,col5,mark=*, mark repeat = 2] file {Chapter_5/tikz/comparison/PdelayB_r05_b01.txt};
\addplot[thick,col5] file {Chapter_5/tikz/comparison/PdelayH_r05_b01.txt};

\addplot[thick,col2,dashed] file {Chapter_5/tikz/comparison/PdelayR_b05.txt};
\addplot[thick,col2,mark=*, mark repeat = 2] file {Chapter_5/tikz/comparison/PdelayB_r05_b05.txt};
\addplot[thick,col2] file {Chapter_5/tikz/comparison/PdelayH_r05_b05.txt};

\addplot[thick,col4,dashed] file {Chapter_5/tikz/comparison/PdelayR_b1.txt};
\addplot[thick,col4,mark=*, mark repeat = 2] file {Chapter_5/tikz/comparison/PdelayB_r05_b1.txt};
\addplot[thick,col4] file {Chapter_5/tikz/comparison/PdelayH_r05_b1.txt};

\end{axis}
\end{tikzpicture}
\caption{$r=0.5$.}
\end{subfigure}
\caption{Asymptotic delay probability in open Erlang-R (dashed), restricted Erlang-R with blocking (marked) and restricted Erlang-R with holding (solid), as function of $\gamma$ with $\b=0.1$ (blue), $\b=0.5$ (orange) and $\b=1$ (green) fixed.}
\label{fig:comparison_delay}
\vspace{3mm}
\end{figure}

Given the expressions for the asymptotic delay probability in the open Erlang-R model, and its restricted versions with blocking and holding, we compare the three policies for various values of $\beta$, $\gamma$ and $r$. 
Figure \ref{fig:comparison_delay} plots the delay probability for blocking ($g^b(\b,\g)$), holding ($g^h(\b,\g)$) and Erlang-R ($g_{\rm HW}(\beta)$) models,  as functions of $\gamma$, while keeping $\beta$ fixed, for three values of $r$. 
We make a couple of observations.
Notice that
\[ g^b(\beta,\gamma) \leq g^h(\beta,\gamma) \leq g_{\rm HW}(\beta) \]
for all $\beta,\gamma>0$ and $r$. 
In that sense, the holding model is an interpolation between the blocking and the open model. 
As expected, the delay probabilities in the restricted models converge to those of the open Erlang-R model, because increasing $\gamma$ is tantamount to lifting the stringent constraints on the system size. Note that the rate of conversion is fast---one can provide probability of waiting close to that of the open model with small values of $\gamma$. Indeed, the fact that when using QED staffing not much of excessive delay results from the beds restriction is important by itself.
Also, we observe that the difference between delay probabilities increases with $r$.

\subsection{The impact of visit number}
\label{subsec:num_visit}
We next reflect on the impact of operational capacity decisions on different patient populations. We measure patient's complexity by the number of times she needs to see the nurse or the physician during her stay. In the ED context, simple-to-treat patients will need to see the physician once, while complex ones will need multiple visits. Hence, we divide the patients into complexity groups by the number of visits in the Needy station. Since the number of visits is geometrically distributed, we have a higher proportion of simple patients than complex ones; that fits well the health care environment.

\begin{figure}
\centering
\begin{subfigure}{0.38\textwidth}
\centering
\begin{tikzpicture}[scale=0.66]
\begin{axis}[
	xmin = 33,
	xmax = 60,
	ymin = 0,
	ymax = 10,
	grid = both, 
	axis line style={->},
	axis lines = left,
	xlabel = $n$,
	yscale = 0.8
]

\addplot[mark=*,black,opacity=0.1] table[x=n,y=in1] {Chapter_5/tikz/inner_vs_outer_wait.txt} node[right,pos=1] {$N=1$};	
\addplot[mark=*,black,opacity=0.2] table[x=n,y=in2] {Chapter_5/tikz/inner_vs_outer_wait.txt};
\addplot[mark=*,black,opacity=0.3] table[x=n,y=in3] {Chapter_5/tikz/inner_vs_outer_wait.txt};	
\addplot[mark=*,black,opacity=0.4] table[x=n,y=in4] {Chapter_5/tikz/inner_vs_outer_wait.txt};	
\addplot[mark=*,black,opacity=0.5] table[x=n,y=in5] {Chapter_5/tikz/inner_vs_outer_wait.txt};	
\addplot[mark=*,black,opacity=0.6] table[x=n,y=in6] {Chapter_5/tikz/inner_vs_outer_wait.txt};	
\addplot[mark=*,black,opacity=0.7] table[x=n,y=in7] {Chapter_5/tikz/inner_vs_outer_wait.txt};
\addplot[mark=*,black,opacity=0.8] table[x=n,y=in8] {Chapter_5/tikz/inner_vs_outer_wait.txt};	
\addplot[mark=*,black,opacity=0.9] table[x=n,y=in9] {Chapter_5/tikz/inner_vs_outer_wait.txt};	
\addplot[mark=*,black,opacity=1] table[x=n,y=in10] {Chapter_5/tikz/inner_vs_outer_wait.txt};	
 	
\addplot[mark=*,col1, thick] table[x=n,y=hold] {Chapter_5/tikz/inner_vs_outer_wait.txt};	

\end{axis}
\end{tikzpicture}
\caption{Expected pre-entering waiting (red) and needy waiting times (black)}
\end{subfigure}
\begin{subfigure}{0.6\textwidth}
\centering
\begin{tikzpicture}[scale=0.66]
\begin{axis}[
	xmin = 33,
	xmax = 60,
	ymin = 0,
	ymax = 10,
	grid = both, 
	axis line style={->},
	axis lines = left,
	xlabel = $n$,
	yscale = 0.8,
	legend cell align=left,
	legend style = {at = {(1.05,0.58)}, anchor = west}
]

\addplot[mark=*,black,opacity=0.1] table[x=n,y=in1] {Chapter_5/tikz/total_wait.txt} node[right,pos=1] {$N=1$};	
\addplot[mark=*,black,opacity=0.2] table[x=n,y=in2] {Chapter_5/tikz/total_wait.txt};
\addplot[mark=*,black,opacity=0.3] table[x=n,y=in3] {Chapter_5/tikz/total_wait.txt};	
\addplot[mark=*,black,opacity=0.4] table[x=n,y=in4] {Chapter_5/tikz/total_wait.txt};	
\addplot[mark=*,black,opacity=0.5] table[x=n,y=in5] {Chapter_5/tikz/total_wait.txt};	
\addplot[mark=*,black,opacity=0.6] table[x=n,y=in6] {Chapter_5/tikz/total_wait.txt};	
\addplot[mark=*,black,opacity=0.7] table[x=n,y=in7] {Chapter_5/tikz/total_wait.txt};
\addplot[mark=*,black,opacity=0.8] table[x=n,y=in8] {Chapter_5/tikz/total_wait.txt};	
\addplot[mark=*,black,opacity=0.9] table[x=n,y=in9] {Chapter_5/tikz/total_wait.txt};	
\addplot[mark=*,black,opacity=1] table[x=n,y=in10] {Chapter_5/tikz/total_wait.txt};	

\legend{{\small $N=1$},{\small $N=2$},{\small $N=3$},{\small $N=4$},{\small $N=5$},{\small $N=6$},{\small $N=7$},{\small $N=8$},{\small $N=9$},{\small $N=10$}};
\end{axis}
\end{tikzpicture}
\caption{Total expected waiting times\\
\quad \\
\quad }
\label{fig:wait_by_visit_b}
\end{subfigure}
\caption{Expected waiting times as a function of $n$ given the number of visits $N$ in the Erlang-R model with holding with $\lambda=2$ $\mu=1$, $\delta=0.25$, $p=0.75$ and $s=9$.}
\label{fig:wait_by_visit}
\end{figure}
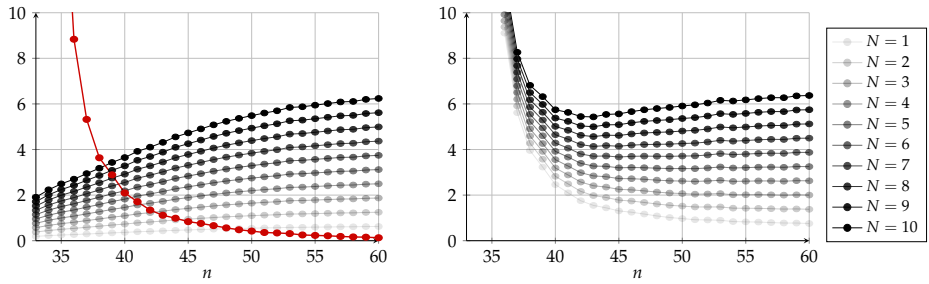

Figure \ref{fig:wait_by_visit} shows the waiting time in the needy and pre-entering queues, and the total waiting time, as a function of $n$ (number of beds), for each complexity group. 
Obviously, the expected waiting time in the pre-entering queue decreases with $n$, while the needy waiting time increases.
For patients who require a relative large number of visits of the physician, in this case more than 6, the total needy wait is the dominant part of the total waiting time. Therefore, as $n$ grows, the total waiting time first decreases and then increases.
In fact, Figure \ref{fig:wait_by_visit_b} suggests that there is an optimal number of beds $n$ that minimizes the total wait for each complexity type. 
Thus, size restrictions reduce the length-of-stay of patients with complex health conditions (given that the constraint is not too tight).
On the other hand, this figure also shows that no such $n$ exists for patients who only require little assistance.
Hence, there is no $n$ that improves the sojourn time of all patients in the ED simultaneously. 
This leaves the decision to the hospital manager to weigh the importance of patients of different complexity levels.

\begin{remark}
From a different perspective, note that in queueing systems such as communication systems, the partitioning of a job to sizable quantities and scheduling those jobs in a similar dynamic to the Erlang-R model became a popular way for increasing throughput. This is because this effectively schedules jobs by their size even though the total job requirements are uncertain. This in fact creates a shortest-job-first policy without prior knowledge of job size \citep{Comte2016}. Considering that perspective we note that the Erlang-R model actually prioritizes simple jobs over complex ones. But without restrictions, when load is too high, such procedures may lead to very long LOS of long jobs. The capacity restriction we analyze in this chapter, in both of its versions, limits such delays. Hence, even in cases in which the returns themselves are created by a managerial decision, imposing the additional managerial restriction on entering the system has benefits. 
\end{remark}

\subsection{Case study: comparison of operational decisions}
\label{sec:case_study}
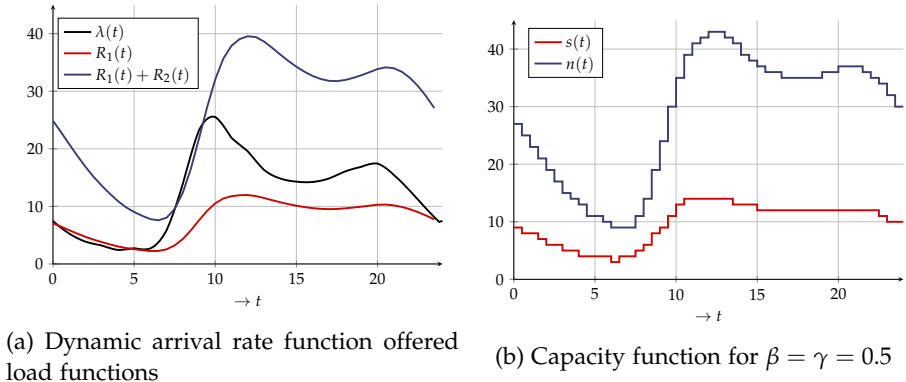
\begin{figure} 
\centering
\begin{subfigure}{0.46\textwidth}
\centering
\begin{tikzpicture}[scale=0.6]

\begin{axis}[
	xmin = 0,
	xmax = 24,
	ymin = 0.0,
	ymax = 45,
	ytick = {0,10,20,30,40},	
	grid = both, 
	axis line style={->},
	axis lines = left,
	xlabel = $\to t$,
	xscale=1.25,
	yscale=1,
	legend cell align=left,
	legend style = {at = {(0.01,0.95)},anchor = north west}
]

\addplot[very thick,black] file {Chapter_5/tikz/lambdaFunction.txt};
\addplot[very thick,col1] file {Chapter_5/tikz/R1.txt};
\addplot[very thick,col5] file {Chapter_5/tikz/R1R2.txt};

\legend{ $\lambda(t)$, $R_1(t)$, $R_1(t)+R_2(t)$};
\end{axis}
\end{tikzpicture}
\caption{Dynamic arrival rate function offered load functions}
\label{fig:Case_study_arrival_pattern_a}
\end{subfigure}
\begin{subfigure}{0.46\textwidth} 
\centering
\begin{tikzpicture}[scale=0.6]

\begin{axis}[
	xmin = 0,
	xmax = 24,
	ymin = 0.0,
	ymax = 45,
	ytick = {0,10,20,30,40},	
	grid = both, 
	axis line style={->},
	axis lines = left,
	xlabel = $\to t$,
	xscale=1.25,
	yscale=1,
	legend cell align=left,
	legend pos = north west
]

\addplot[very thick,col1] file {Chapter_5/tikz/casestudy_new/sFunction.txt};
\addplot[very thick,col5] file {Chapter_5/tikz/casestudy_new/nFunction.txt};

\legend{ $s(t)$, $n(t)$};
\end{axis}
\end{tikzpicture}
\caption{Capacity function for $\beta=\gamma=0.5$}
\label{fig:Case_study_arrival_pattern_b}
\end{subfigure}
\caption{Empirical arrival rate and offered load functions $R_1(t)$ and $R_1(t)+R_2(t)$ in Israeli ED and corresponding capacity functions.}
\label{fig:Case_study_arrival_pattern}
\end{figure}

We now illustrate how the managerial decision to operate under a specific operational regime affects ED performance in terms of efficiency and quality-of-care, through a case study. 
The practical environment we investigate is the ED of a moderately-sized hospital, which faces the arrival pattern $\lambda(t)$ plotted in Figure \ref{fig:Case_study_arrival_pattern_a} on a typical workday. 
Other parameters of the model are estimated to be $\mu = 6.67,\ \delta = 2.18$ and $p = 0.76$, so that $r = 0.301$. These parameters were taken from \cite{YomTov2014}. In order to set time-varying staffing levels $s(t)$ and $n(t)$, we adopt the \textit{mean-offered load} (MOL) approximation of the demand process of~\cite{Jennings1996}.
This approach initially presumes infinite capacity to obtain the number of patients $R(t)$ in the queueing system as a function of time. 
This offered load function then replaces the (constant) value of $R$ in the stationary dimensioning scheme under consideration, to determine the adequate number of servers at each point in time. 
Following this idea in our two-dimensional queueing system, we find the offered load function for the nurses $R_1(t)$ and the offered load function for the beds $R_1(t)+R_2(t)$ as the solution of the system of ODEs,
\begin{align} \label{eq:offeredloadODE}
\frac{d}{dt} R_1(t) &= \l(t) + \d R_2(t) - \m R_1(t),\\
\frac{d}{dt} R_2(t) &= p\m R_1(t) - \d R_2(t),
\end{align}
see \cite[Thm.~2]{YomTov2014} for details.
For this case study's parameters, these offered load functions are also plotted in Figure \ref{fig:Case_study_arrival_pattern_a}.
While the number of nurses can be adjusted in a relatively flexible manner, the value of $n$, which echoes a hard restriction on the ED capacity, is naturally less amenable to fluctuations. The reason is that the maximum ED capacity is to a large extent determined by its hardware, such as beds and medical equipment.
However, the ED manager might deliberately consider reducing $n$ during more quiet periods of the day, e.g.\ during the night, by imposing bed-to-physician constraints. This is done, for example, when setting a case management constraint \citep{EDexperiment,Campello2016}. 
Therefore, we consider the scenario in which both $s$ and $n$ are time-dependent but we do not force a constant case management quantity, rather let our new methodology recommend an appropriate one. 

Extrapolating Algorithm \ref{alg:stationarydimensioning} to the time-varying case, Step 4 is replaced by
\begin{align*}
s(t) &= R_1(t) + \b\sqrt{R_1(t)},\\
n(t) &= R_1(t)+R_2(t) + \g\sqrt{R_1(t)+R_2(t)},
\end{align*}
for some $\b,\g>0$. 
Since $R_1(t)$ and $R_2(t)$ are given, the QED staffing problem again reduces to finding the pair $(\beta,\gamma)$.

Figure \ref{fig:Case_study_arrival_pattern_b} plots the capacity functions for $\beta = 0.5$ and $\gamma=0.5$, assuming capacity can only be adjusted every 30 minutes. 
In this case study, we consider three pairs of parameters $(\beta,\gamma)$.
For each of these we investigate, using simulation, the differences in the time-varying performance indicators between the policy with blocking and holding.

\begin{figure}
\centering

\begin{subfigure}{0.32\textwidth}
\begin{tikzpicture}[scale=0.47]

\begin{axis}[
	xmin = 0,
	xmax = 24,
	ymin = 0.0,
	ymax = 1.0,
	ytick = {0,0.2,0.4,0.6,0.8,1.0},	
	xtick = {0,3,6,9,12,15,18,21,24},
	grid = both, 
	axis line style={->},
	axis lines = left,
	xlabel = $\to t$,
	xscale=1.1,
	yscale=1
]

\addplot[very thick,col1] table[x=t,y=delay_b01g2] {Chapter_5/tikz/casestudy_new/case_study_block.txt};
\addplot[very thick,col1,dashed] table[x=t,y=delay_b01g2] {Chapter_5/tikz/casestudy_new/case_study_hold.txt};

\addplot[very thick,col4] table[x=t,y=delay_b1g15] {Chapter_5/tikz/casestudy_new/case_study_block.txt};
\addplot[very thick,col4,dashed] table[x=t,y=delay_b1g15] {Chapter_5/tikz/casestudy_new/case_study_hold.txt};

\addplot[very thick,col5] table[x=t,y=delay_b2g1] {Chapter_5/tikz/casestudy_new/case_study_block.txt};
\addplot[very thick,col5,dashed] table[x=t,y=delay_b2g1] {Chapter_5/tikz/casestudy_new/case_study_hold.txt};
\end{axis}
\end{tikzpicture}
\caption{$\P({\rm delay})$}
\label{fig:simulation_results_a}
\end{subfigure}
\begin{subfigure}{0.32\textwidth}
\begin{tikzpicture}[scale=0.47]

\begin{axis}[
	xmin = 0,
	xmax = 24,
	ymin = 0.0,
	ymax = 0.52,
	ytick = {0,0.1,0.2,0.3,0.4,0.5},
	xtick = {0,3,6,9,12,15,18,21,24},	
	grid = both, 
	axis line style={->},
	axis lines = left,
	xlabel = $\to t$,
	xscale=1.1,
	yscale=1,
	legend cell align=left,
	legend style = {at = {(0.9,0.95)}, anchor = north east}
]

\addplot[very thick,col1] table[x=t,y=block_b01g2] {Chapter_5/tikz/casestudy_new/case_study_block.txt};
\addplot[very thick,col4] table[x=t,y=block_b1g15] {Chapter_5/tikz/casestudy_new/case_study_block.txt};
\addplot[very thick,col5] table[x=t,y=block_b2g1] {Chapter_5/tikz/casestudy_new/case_study_block.txt};

\addplot[very thick,col1,dashed] table[x=t,y=block_b01g2] {Chapter_5/tikz/casestudy_new/case_study_hold.txt};
\addplot[very thick,col4,dashed] table[x=t,y=block_b1g15] {Chapter_5/tikz/casestudy_new/case_study_hold.txt};
\addplot[very thick,col5,dashed] table[x=t,y=block_b2g1] {Chapter_5/tikz/casestudy_new/case_study_hold.txt};

\legend{{$(\beta,\gamma)=(0.1,2)$},{$(\beta,\gamma)=(1,1.5)$},{$(\beta,\gamma)=(2,1)$}};
\end{axis}
\end{tikzpicture}
\caption{$\P({\rm block})$ or $\P({\rm hold})$}
\label{fig:simulation_results_b}
\end{subfigure}
\begin{subfigure}{0.32\textwidth}
\begin{tikzpicture}[scale=0.47]

\begin{axis}[
	xmin = 0,
	xmax = 24,
	ymin = 0.0,
	ymax = 5.5,
	xtick = {0,3,6,9,12,15,18,21,24},	
	grid = both, 
	axis line style={->},
	axis lines = left,
	xlabel = $\to t$,
	xscale=1.1,
	yscale=1
]

\addplot[very thick,col1] table[x=t,y=ratio_b01g2] {Chapter_5/tikz/casestudy_new/case_study_block.txt};
\addplot[very thick,col1,dashed] table[x=t,y=ratio_b01g2] {Chapter_5/tikz/casestudy_new/case_study_hold.txt};

\addplot[very thick,col4] table[x=t,y=ratio_b1g15] {Chapter_5/tikz/casestudy_new/case_study_block.txt};
\addplot[very thick,col4,dashed] table[x=t,y=ratio_b1g15] {Chapter_5/tikz/casestudy_new/case_study_hold.txt};

\addplot[very thick,col5] table[x=t,y=ratio_b2g1] {Chapter_5/tikz/casestudy_new/case_study_block.txt};
\addplot[very thick,col5,dashed] table[x=t,y=ratio_b2g1] {Chapter_5/tikz/casestudy_new/case_study_hold.txt};
\end{axis}
\end{tikzpicture}
\caption{Nurse-to-patient ratio.}
\label{fig:simulation_results_c}
\end{subfigure}

\caption{Simulation results for case study. Solid and dashed lines represent time-varying performance in the blocking and holding model, respectively.}
\label{fig:simulation_results}
\end{figure}
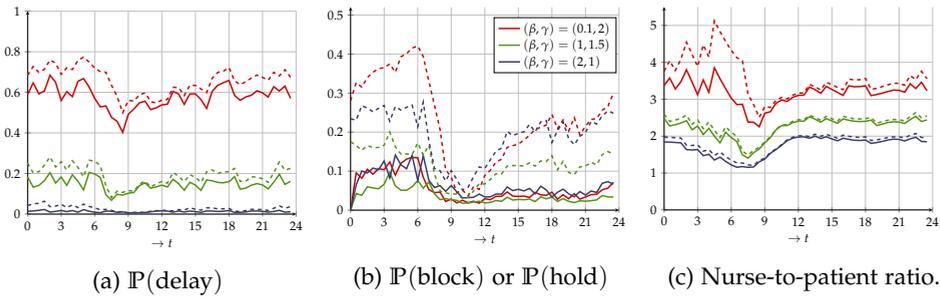

The simulation results are presented in Figure \ref{fig:simulation_results}. 
Figure \ref{fig:simulation_results_a} shows that the MOL approach for capacity allocation roughly stabilizes the delay probability. 
Figure \ref{fig:simulation_results_b} shows that the fraction of patients not entering the ED on arrival in the blocking model is reasonable for all parameter pairs considered and the graphs are ordered according to $\gamma$.
We also see a significant difference with holding.
Observe also that the holding probability drops in the period 8--13, which is exactly the period when the system is experiencing peak offered load. 
Hence, this temporary reduction is in line with our asymptotic findings that the probability of blocking/holding is $O(1/\sqrt{R_1})$.

Finally note that the three parameter settings lead to different nurse-to-patient ratios. 
Clearly, larger $\beta$ leads to small nurse-to-patient ratios (due do larger staffing). 
Figure \ref{fig:simulation_results_c} demonstrates that for $(\beta,\gamma) = (1,1.5)$ and $(\beta,\gamma) = (2,1)$ the difference between the holding policy and the blocking policy is small. However, for $(\beta,\gamma) = (0.1,2)$ we see a significant increase in the ratio during night hours. 
This may be due to the tight nurse schedule, that causes the holding queue to build up just before midnight. 
This queue then empties latter on, causing an increase in workload per nurse in the period 24--7.

To see the direct effect of the size restriction on the queue lengths, we plotted the mean holding and service queue lengths in the holding model as a function of the parameter $\gamma$ in Figure \ref{fig:simulation_queuelengths}.
\begin{figure}
\centering
\begin{subfigure}{0.48\textwidth}
\centering
\begin{tikzpicture}[scale=0.6]
\begin{axis}[
	xmin = 0,
	xmax = 24,
	ymin = 0.0,
	ymax = 5.2,
	xtick = {0,3,6,9,12,15,18,21,24},	
	grid = both, 
	axis line style={->},
	axis lines = left,
	xlabel = $\to t$,
	xscale=1.25,
	yscale=1,
	legend cell align=left,
	legend pos = north west
]

\addplot[very thick,col1] file {Chapter_5/tikz/casestudy_new/holdingQueue_g01.txt};
\addplot[very thick,col3] file {Chapter_5/tikz/casestudy_new/holdingQueue_g025.txt};
\addplot[very thick,col4] file {Chapter_5/tikz/casestudy_new/holdingQueue_g05.txt};
\addplot[very thick,col5] file {Chapter_5/tikz/casestudy_new/holdingQueue_g1.txt};

\legend{$\gamma=0.1$, $\gamma=0.25$,$\gamma=0.5$,$\gamma=1$}
\end{axis}
\end{tikzpicture}
\caption{Mean holding queue length}
\end{subfigure}
\begin{subfigure}{0.48\textwidth}
\centering
\begin{tikzpicture}[scale=0.6]
\begin{axis}[
	xmin = 0,
	xmax = 24,
	ymin = 0.0,
	ymax = 15,
	xtick = {0,3,6,9,12,15,18,21,24},	
	grid = both, 
	axis line style={->},
	axis lines = left,
	xlabel = $\to t$,
	xscale=1.25,
	yscale=1,
	legend cell align=left,
	legend style = {at = {(0.82,0.05)},anchor = south east}
]

\addplot[very thick,col1] file {Chapter_5/tikz/casestudy_new/serviceQueue_g01.txt};
\addplot[very thick,col3] file {Chapter_5/tikz/casestudy_new/serviceQueue_g025.txt};
\addplot[very thick,col4] file {Chapter_5/tikz/casestudy_new/serviceQueue_g05.txt};
\addplot[very thick,col5] file {Chapter_5/tikz/casestudy_new/serviceQueue_g1.txt};
\addplot[very thick,dashed] file {Chapter_5/tikz/casestudy_new/serviceQueue_R.txt};
\legend{$\gamma=0.1$, $\gamma=0.25$,$\gamma=0.5$,$\gamma=1$,Erlang-R}
\end{axis}
\end{tikzpicture}
\caption{Mean service queue length}
\end{subfigure}
\caption{Simulated queue length of holding model with different values of $\gamma$.}
\label{fig:simulation_queuelengths}

\vspace{-6mm}

\end{figure}
Note that for all $\gamma$ considered, the holding queue lengtsh are, as expected, of a smaller order than the number of patients admitted. 
Also, the holding queue length decreases as we increase $\gamma$. 
The service queue lengths naturally approach the expected queue lengths in the Erlang-R model as $\gamma$ is increased.

\section{Conclusion \& future research}
\label{sec:conclusion}

In this chapter we developed and analyzed a queueing network tailored to a health care environment with finite-size restrictions. 
Using the asymptotic approximations, numerical analysis and simulation, we gained insight into staffing problems that arise in EDs, and proposed an efficient, flexible, and easy to implement methodology to dimension medical facilities through a two-fold staffing rule.

The dimensioning scheme we developed provides a powerful and elegant way of finding adequate staffing levels in emergency departments. 
Nonetheless, we see some avenues for further research. 

The asymptotic approximations we developed enabled us to take the first step towards characterizing the pre-entering queue behavior in the 
QED regime.  
We observed how the holding queue length vanishes at rate $1/\sqrt{R_1}$ as $R_1\to\infty$.
Yet, our analysis did not yield explicit characteristics on the holding queue and holding times.
These performance indicators are naturally important to study if one wants to consider the trade-off between waiting time inside the ED and waiting time outside the ED time (pre-entering time). 

Furthermore, it is worthwhile to study the robustness of our approximations against the service and content time distributions. Since the content phase of a patient is modeled after an infinite-server queue, we expect our approximations to be useful for content time distributions beyond the exponential distribution as well, due to distributional insensitivity of the service time in infinite-server queues. For the needy phase, modeled after a multi-server queue, this insensitivity result does not hold and hence this needs further research. 

Finally, the restricted Erlang-R model obviously gives a highly simplified view of the complex reality of the ED. 
In practice, distinctive features such as a triage system (with patient priorities), patient boarding time and availability of medical equipment may play a decisive role on ED dynamics.
However, we think the analysis and dimensioning algorithms presented in this chapter can serve as a building block for staffing procedures that do account for these case-specific factors.

\section*{Appendix}

\addcontentsline{toc}{section}{\hspace{7.1mm} Appendix}

\begin{subappendices}

\settocdepth{chapter}

\section{Description of the QBD process}
\label{app:QBDdescription}

\subsection{The QBD-process}
\label{app:theQBDprocess}
We consider the QBD-process $X(t)=(N(t),Q_1(t))$ in stationarity. Let $\nu(i)=\min\{i,s\}\mu$. To determine the (outgoing) transition rates of the process $X$ we distinguish between the following cases:
\begin{itemize}
\item \emph{Transitions from $(0,0)$:} There are no patients in the Emergency Department and thus the only possible occurrence is when a new patient arrives. This results in a transition to $(1,1)$ and occurs with rate $\lambda$.
\item \emph{Transitions from $(i,0), 1 \leq i < n$:} There are exactly $i$ patients assigned to a bed of which none are seen by a nurse. Then either one of those patients becomes needy, or a new patient arrives at the Emergency Department that can immediately be seen by a nurse. The first results in a transition to $(i,1)$ and occurs at rate $i \delta$, and the second results in a transition to $(i+1,1)$ and occurs with rate $\lambda$.
\item \emph{Transitions from $(i,0), i \geq n$:} Again, the only possible transitions arise from either a newly arrived patient or a patient assigned to a bed becoming needy. However, a newly arrived patient finds all beds occupied and needs to wait. Thus, with rate $\lambda$ we have a transition to $(i+1,0)$ and with rate $n \delta$ a transition to $(i,1)$.
\item \emph{Transitions from $(i,i), i < n$:} In this case all patients assigned to a bed are in need of service. With rate $\lambda$ a new patient arrives at the Emergency Department. She joins the (possible) queue to be seen by a nurse immediately, so this results in a transition to $(i+1,i+1)$. Moreover, since there are only $s < n$ nurses, a service completion occurs with rate $\nu(i)$. With probability $p$ the patient turns to the holding phase, so in total we still have $i$ patients with one patient less in queue for a nurse. With probability $1-p$ the patient leaves the Emergency Department, decreasing both $N$ and $Q_1$ by one. In other words, with rate $p \nu(i)$ we have a transition to $(i,i-1)$ and with rate $(1-p)\nu(i)$ we have a transition to $(i-1,i-1)$.
\item \emph{Transitions from $(n,n)$:} Similar to the previous case, we have a transition to $(n,n-1)$ with rate $p s \mu$ and with rate $(1-p)s \mu$ we have a transition to $(n-1,n-1)$. In this case however, a newly arrived patient finds all beds occupied, resulting in a transition to $(n+1,n)$ with rate $\lambda$.
\item \emph{Transitions from $(i,n), i > n$:} We have a transition to $(i+1,n)$ with rate $\lambda$ and a transition to $(i,n-1)$ with rate $p s \mu$. In case that a patient leaves the Emergency Department there are $i-n>0$ patients in the holding room waiting for an available bed. Thus, one of the $i-n$ patients in the holding room is assigned to the available bed in need of service. That is, with rate $(1-p) s \mu$ we have a transition to $(i-1,n)$.
\item \emph{Transitions from $(i,j), 1 \leq j < i  < n$:} There are four possible transitions. First, with rate $\lambda$ there is a new arrival which results in a transition to $(i+1,j+1)$. Second, with rate $(i-j) \delta$ a patient in one of the beds becomes needy, which results in a transition to $(i,j+1)$. Third, with rate $p \nu(j)$ a patient turns to the content state after service completion, which results in a transition to $(i,j-1)$. Last, with rate $(1-p) \nu(j)$ a patient leaves the Emergency Department after service completion, which results in a transition to $(i-1,j-1)$.
\item \emph{Transitions from $(n,j), 1 \leq j  < n$:} This case is similar to the previous one. The only difference arises when a new patient arrives, since all $n$ beds are already occupied. Thus, with rate $\lambda$ we have a transition to $(n+1,j)$.
\item \emph{Transitions from $(i,j), i > n, 1 \leq j  \leq n$:} This case is the previous one, except when a patient leaves the Emergency Department after service completion. Then one of the $(i-n)$ patients in the holding room will be assigned to a bed in need of service. This results in a transition to $(i-1,j)$ with rate $(1-p) \nu(j)$.
\end{itemize}

\begin{figure}
 \centering
\begin{tikzpicture}[scale=0.66]
\draw[step=1cm,gray!50!,very thin] (0,0) grid (15.5,8.5);
\draw[thick,->] (0,0) -- (15.5,0);
\draw[thick,->] (0,0) -- (0,8.5);
\draw[thick] (0,0) -- (8,8);
\draw[thick,dashed,black!50!] (8,0) -- (8,8);
\draw[thick] (8,8) -- (15.5,8);
\foreach \x in {0,1,2,3,4,5,6,7,8,9,10,11,12,13,14,15}
	\foreach \y in {0,1,2,3,4,5,6,7,8}
    \draw[fill] (\x,\y) circle [radius=0.025];

\node [below left] at (0,0) {$0$};
\node [left] at (0,2) {$j$};
\node [left] at (0,5) {$k$};
\node [left] at (0,8) {$n$};
\node [below] at (6,0) {$i$};
\node [below] at (8,0) {$n$};
\node [above left] at (0,8.5) {$Q_1$}; 
\node [below right] at (15.5,0) {$N$}; 

\path [->,thick,-latex] (0,0) edge [bend right] (1,1);

\path [->,thick,-latex] (3,0) edge (3,1);
\path [->,thick,-latex] (3,0) edge (4,1);

\path [->,thick,-latex] (4,4) edge [bend right] (5,5);
\path [->,thick,-latex] (4,4) edge (4,3);
\path [->,thick,-latex] (4,4) edge [bend right] (3,3);

\path [->,thick,-latex] (6,2) edge (6,3);
\path [->,thick,-latex] (6,2) edge (7,3);
\path [->,thick,-latex] (6,2) edge (6,1);
\path [->,thick,-latex] (6,2) edge (5,1);

\path [->,thick,-latex] (8,5) edge (8,6);
\path [->,thick,-latex] (8,5) edge (8,4);
\path [->,thick,-latex] (8,5) edge (9,5);
\path [->,thick,-latex] (8,5) edge (7,4);

\path [->,thick,-latex] (8,8) edge (8,7);
\path [->,thick,-latex] (8,8) edge [bend right] (7,7);
\path [->,thick,-latex] (8,8) edge [bend left] (9,8);

\path [->,thick,-latex] (11,8) edge (11,7);
\path [->,thick,-latex] (11,8) edge [bend left] (12,8);
\path [->,thick,-latex] (11,8) edge [bend left] (10,8);

\path [->,thick,-latex] (11,0) edge (11,1);
\path [->,thick,-latex] (11,0) edge [bend left] (12,0);

\path [->,thick,-latex] (12,5) edge (12,6);
\path [->,thick,-latex] (12,5) edge (12,4);
\path [->,thick,-latex] (12,5) edge (13,5);
\path [->,thick,-latex] (12,5) edge (11,5);

\node [above] at (12.75,5) {\scriptsize $\lambda$};
\node [above] at (11.25,5) {\scriptsize $(1-p)\nu(k)$};
\node [right] at (12,5.75) {\scriptsize $(n-k)\delta$};
\node [right] at (12,4.25) {\scriptsize $p \nu(k)$};

\node [above] at (6.75,2.25) {\scriptsize $\lambda$};
\node [below] at (5,2) {\scriptsize $(1-p)\nu(j)$};
\node [above] at (6,3) {\scriptsize $(i-j)\delta$};
\node [right] at (6,1.25) {\scriptsize $p \nu(j)$};

\end{tikzpicture}
\caption{Transition diagram for the Erlang-R model with holding.}
\label{fig:QBDIllustration}
\vspace{-3mm}
\end{figure}

\noindent
The state space and transition rates of the Erlang-R model with holding are illustrated in Figure~\ref{fig:QBDIllustration}.
The state space can be partitioned according to its levels, where level $i$ corresponds to a total queue length $N=i$ patients. This results in an infinite-sized matrix consisting of blocks, where each block corresponds to the transition flow from one level to another. Since the only transitions allowed are within the same level or between two adjacent levels in a QBD-process, we obtain a tridiagonal block structure. Each block consists of elements representing the transition rate of one state to another, and therefore each block is a matrix of size at most $(n+1) \times (n+1)$.

For the Erlang-R model with holding this gives the following result. Let $P$ denote the transition matrix of the process $(N(t),Q_1(t))$. We have the boundary levels $\{1,2,...,n\}$ and $P$ is of the form

\[
P = \left( \begin{array}{cccccccccc}
B_{00} & B_{01} & & & & & & & & \\
B_{10} & B_{11} & B_{12} & & & & & & & \\
 & B_{21} & B_{22} & B_{23} & & & & & & \\
 & & \ddots & \ddots &\ddots & & & & & \\
 & & & & & B_{n \, n-1} & & & & \\
 & & & & B_{n-1 \, n} & B_{nn} & A_0 & & & \\
 & & & & & A_2 & A_1 & A_0 & &  \\
 & & & & & & A_2 & A_1 & A_0 & \\
 & & & & & & & \ddots & \ddots & \ddots \\
\end{array} \right),
\]
where $B_{ii} \in \mathbb{R_1}^{(i+1) \times (i+1)}$, $B_{i \, i-1} \in \mathbb{R_1}^{(i+1) \times i}$, $B_{i-1 \, i} \in \mathbb{R_1}^{i \times (i+1)}$, and $A_0,A_1,A_2 \in \mathbb{R_1}^{(n+1)\times(n+1)}$. The matrices of transition rates for the boundary states are given by
\[
B_{00}=(-\lambda), 
\qquad
B_{i-1 \, i} = \left( \begin{array}{ccccc}
0 & \lambda & & &  \\
& \ddots & \lambda & & \\
 & & \ddots & \ddots &\\
 & & & 0 & \lambda  \\
\end{array} \right), 
\]
\[
B_{i \, i-1} = \left( \begin{array}{cccc}
0 & & &   \\
 (1-p)\mu & 0 & &  \\
 & (1-p)\nu(2)& \ddots &  \\
 & & \ddots & 0 \\
 & & & (1-p)\nu(i)\\
\end{array} \right),
\]
and
\[
\scriptsize
B_{ii} = \left(
\begin{array}{ccccccccc}
-(\lambda+i \delta) & i \delta & & & &\\
p \mu & -(\lambda+\mu+(i-1)\delta) & (i-1)\delta & & &\\
& \ddots & \ddots & \ddots & & \\
& & p \nu(i-1) & -(\lambda+\nu(i-1)+\delta) & \delta \\
& & & & p \nu(i) & -(\lambda+\nu(i)) \\
\end{array} \right).
\]
Moreover, the transition rates are given by
\[
A_0 = \left( \begin{array}{ccccc}
\lambda & & & & \\
& \lambda & & & \\
 & & & \ddots & \\
 & & & & \lambda  \\
\end{array} \right)
\]
\[ A_2 = \left( \begin{array}{ccccccc}
0 & & & & & & \\
& (1-p)\mu & & & & & \\
& & 2(1-p)\mu & & & & \\
& & & \ddots & & & \\
& & & & s(1-p)\mu & & \\
& & & & & \ddots & \\
& & & & & & s(1-p)\mu \\
\end{array} \right),
\]
and
\[
\scriptsize
A_1 = \left( \arraycolsep=0.55pt
\begin{array}{cccccccc}
-(\lambda+n \delta) & n \delta & & & & & & \\
p \mu & -(\lambda+\mu+(n-1)\delta) & (n-1)\delta & & & & & \\
& \ddots & \ddots & \ddots & & & & \\
& & s p\mu & -(\lambda+s\mu+(n-s)\delta) & (n-s)\delta & & & \\
& & & \ddots & \ddots & \ddots & & \\
& & & & s p\mu & -(\lambda+s\mu+\delta) & \delta & \\
& & & & & s p\mu & -(\lambda+s\mu)\\
\end{array} \right).
\]

\subsection{Stability condition}
\label{app:stability}
From the general theory of QBD processes \citep{Neuts1981} follows that the Markov process $(N(t),Q_1(t))$ is ergodic (stable) if and only if
\begin{equation}
\pi A_0 e < \pi A_2 e,
\label{eq:QBDstableCondition}
\end{equation}
where $e$ is the all one column vector and $\pi=(\pi_0,...,\pi_n)$ is the equilibrium distribution of the Markov process with generator $A_0+A_1+A_2$. In other words, $\pi$ is such that
\begin{equation}
\begin{array}{ll}
\pi(A_0+A_1+A_2) =0, & \pi e =1,
\end{array}
\label{eq:QBDstableProbabilityVector}
\end{equation}
and
\[
A_0+A_1+A_2 = \qquad\qquad\qquad
\]
\begin{align*} 
{\scriptsize
\left(
\begin{array}{cccccccc}
-n \delta & n \delta & & & & & & \\
p \mu & -(p\mu+(n-1)\delta) & (n-1)\delta & & & & & \\
& \ddots & \ddots & \ddots & & & & \\
& & s p\mu & -(ps\mu+(n-s)\delta) & (n-s)\delta & & & \\
& & & \ddots & \ddots & \ddots & & \\
& & & & p s \mu & -(ps\mu+\delta) & \delta & \\
& & & & & p s \mu & -ps\mu\\
\end{array} \right).
}
\end{align*}

Then $\pi$ must satisfy the balance equations
\begin{align*}
- n \delta \pi_0 + p \mu \pi_1 &= 0, \\
(n-j+1)\delta \pi_{j-1} - (p\nu(j) +(n-j)\delta) \pi_j + p \nu(j+1) \pi_{j+1} &= 0, \\
\delta \pi_{n-1} - p s \mu \pi_n &= 0,
\end{align*}
with $\nu(j)=\min\{j,s\}\mu$, and the normalization condition
\[
\sum_{i=0}^n \pi_i=1.
\]
It is readily verified that
\begin{equation}
\pi_i =
\left\{\begin{array}{ll}
\pi_0 \binom{n}{i} \left(\frac{\delta}{p \mu}\right)^i & \textrm{\normalfont for } 0 \leq i \leq s, \\
\pi_0 \binom{n}{i} \frac{i!}{s!} s^{s-i} \left(\frac{\delta}{p \mu}\right)^i & \textrm{\normalfont for } s+1 \leq i \leq n \\
\end{array} \right.
\label{eq:eqdistr}
\end{equation}
with
\begin{align*}
\pi_0= \left(\sum_{i=0}^{s} \binom{n}{i} \left(\frac{\delta}{p \mu}\right)^i + \sum_{i=s+1}^{n} \binom{n}{i} \frac{i!}{s!} s^{s-i} \left(\frac{\delta}{p \mu}\right)^i \right)^{-1},
\end{align*}
satisfies the balance equations and the normalization condition.

\begin{proposition}
The distribution of the closed two-node Jackson network illustrated in Figure~\ref{fig:Jennings} is given by
\begin{equation}
\hat{\pi_i} =
\left\{\begin{array}{ll}
\hat{\pi}_0 \binom{n}{i} \left(\frac{\delta}{p \mu}\right)^i & \textrm{\normalfont for } 0 \leq i \leq s, \\
\hat{\pi_0} \binom{n}{i} \frac{i!}{s!} s^{s-i} \left(\frac{\delta}{p \mu}\right)^i & \textrm{\normalfont for } s+1 \leq i \leq n \\
\end{array} \right.
\label{eq:eqdistrTwistedJennings}
\end{equation}
with
\begin{align*}
\hat{\pi}_0= \left[\sum_{i=0}^{s} \binom{n}{i} \left(\frac{\delta}{p \mu}\right)^i + \sum_{i=s+1}^{n} \binom{n}{i} \frac{i!}{s!} s^{s-i} \left(\frac{\delta}{p \mu}\right)^i \right]^{-1}.
\end{align*}
\label{prop:CriticalTilburgdistr}
\end{proposition}

\begin{proof}
We have a two-node closed Jackson network, with probability transition matrix
\[
P = \left(
\begin{array}{cc}
1-p & p \\
1 & 0
\end{array} \right).
\]
Let $r_i(m)$ denote the rate of service when there are $m$ patients at queue $i$, so $r_1(m)=\min\{m,s\}$ and $r_2(m)=m$. The throughput vector $\gamma = (\gamma_1,\gamma_2) \in \mathbb{R_1}^2$ must satisfy $\gamma = \gamma P$ and we find that $\gamma=(p,1)$ suffices. From the general theory of Jackson networks, see \cite{Jackson1963}, it follows that the stationary distribution is given by
\begin{align*}
\pi_i = G^{-1} g_1(i) g_2(n-i)
\end{align*}
with
\begin{align*}
\begin{array}{ll}
g_1(i)= \frac{(\gamma_1/\mu)^i}{\prod_{m=1}^i r_1(m)}, & g_2(n-i)= \frac{(\gamma_2/\delta)^{n-i}}{\prod_{m=1}^{n-i} r_2(m)},
\end{array}
\end{align*}
and normalization constant $G= \sum_{i=0}^n g_1(i) g_2(n-i)$. Then,
\begin{align*}
g_1(i) &=  \left\{\begin{array}{ll}
\frac{1}{i! \mu^i} & \textrm{\normalfont for } 0 \leq i \leq s, \\
\frac{1}{s! s^{i-s} \mu^i} & \textrm{\normalfont for } s+1 \leq i \leq n, \\
\end{array} \right.\\
g_2(n-i) &=\frac{1}{(n-i)!} \left(\frac{p}{\delta}\right)^n \left(\frac{\delta}{p}\right)^i,
\end{align*}
and rewriting the expressions yields~\eqref{eq:eqdistrTwistedJennings}.
\end{proof}

\subsection{Stationary distribution}
\label{app:StationaryDistributrion}
Assuming that the stability condition is satisfied, we can determine the unique stationary distribution of the Markov process $(N(t),Q_1(t))$. The vector $\pi_i$ can be written as $\pi_{n+i}= \pi_n G^{i}$ for $i=0,1,...$, where $G$ is the minimal nonnegative solution of the non-linear matrix equation
\begin{equation}
A_0+G A_1 + G^2 A_2=0.
\label{eq:MG-G}
\end{equation}
The balance equations can be written as
\[
\begin{array}{ll}
\pi_{i-1} A_0+ \pi_i A_1 + \pi_{i+1} A_2=0, & i=n+1,n+2,...
\end{array}
\]
and using $\pi_{n+i}= \pi_n G^{i-n}$ for $i=0,1,...$, this find
\[
\begin{array}{ll}
\pi_n G^{i-n-1} \left(A_0+ G A_1 + G A_2\right)=0, & i=n+1,n+2,....
\end{array}
\]

\noindent
Moreover, we have the boundary equations
\begin{align*}
\pi_0 B_{00} + \pi_1 B_{10} &= 0 \\
\pi_0 B_{01} + \pi_1 B_{11} + \pi_2 B_{21} &= 0 \\
\pi_1 B_{12} + \pi_1 B_{22} + \pi_2 B_{32} &= 0 \\
 &\vdots&  \\
\pi_{n-2} B_{n-2 \, n-1} + \pi_{n-1} B_{n-1 \, n-1} + \pi_{n} B_{n \, n-1} &= 0 \\
\pi_{n-1} B_{n-1 \, n} + \pi_{n} B_{nn} + \pi_{n+1} A_2 &= 0, 
\end{align*}
along with the normalization equation
\[
1 = \sum_{i=0}^{\infty} \pi_i e = \sum_{i=0}^{n-1} \pi_i e + \pi_n(I-G)^{-1}e,
\]
where we slightly abuse notation by using $e$ as the all ones vector of appropriate size. We note that the matrix $G$ has a spectral radius less than one and therefore $(I-G)$ is invertible.

These equations provide the tools for finding the equilibrium probabilities. Although it is hard to solve $G$ analytically from Equation~\eqref{eq:MG-G}, it is easy to solve numerically by using the following algorithm (matrix-geometric method). Rewriting~\eqref{eq:MG-G} gives
\[
G=-(A_0+G^2 A_2) A_1^{-1},
\]
where $A_1$ is invertible, since it is a transient generator matrix. Let
\[
G_{k+1}=-(A_0+G_k^2 A_2) A_1^{-1},
\]
starting with $G_0=0$. We note that $G_k \uparrow G$ as $k$ grows to infinity \citep{Neuts1981}. Once $||G_{k+1}-G_{k}||_2$ is below a certain preset threshold, we approximate $G$ by $G_{k+1}$.

\section{Proof of Proposition \ref{thm:stochasticordering}}\label{app:stochastic_ordering}

First, note that by definition of the Erlang-R model with holding, in which no more than $n$ patients can be admitted in the ED simultaneously, that $Q_1^h(t)+Q_2^h(t) \leq n = Q_1^J(t) + Q_2^J(t)$ follows directly. 
Therefore, we only consider the relation between the states in the blocking and holding variants Erlang-R model.

As noted Section \ref{sec:Markov_process}, the model with holding can be characterized as a three-dimensional Markov chain $X^h(t) = (H(t),Q^h_1(t),Q^h_2(t))$ in which the components denote the number of holding, needy and content patients respectively. The Erlang-R model with blocking similarly admits a Markov process description, but with two dimensions, namely $X^b(t) = (Q^b_1(t),Q^b_2(t))$.

We prove the result by constructing a coupling between the Markov processes $X^h$ and $X^b$. Let $Z(t) := \big(\hat{X}^h(t),\hat X^b(t)\big) = \big(\hat{H}(t),\hat{Q}_1^h(t),\hat{Q}_2^h(t),\hat{Q}^b_1(t),\hat{Q}^b_2(t)\big)$.

We first define the transition rates of this five-dimensional Markov process, which naturally only depend on the current state of the system. 
After that we show that the transition rates relevant to $\hat{X}^h(t)$ and $\hat{X}^{b}(t)$ coincide with those of $X^h(t)$ and $X^b(t)$, respectively. The latter implies that the marginal transitions of $\hat{X}^h(t)$ and $X^h(t)$ (and $\hat{X}^b(t)$ and $X^b(t)$) are equal, and hence so are their probability distribution of the Markov processes. 

Let $Z(t) = (h,\qh,\ch,\qy,\cy)$. While defining the reachable states from this state and associated transition rates, we distinguish four transition types, and further differentiate the transition rates depending on the current state.\\
\\*
\textbf{Arrival.}
Arrivals  occur in both models simultaneously, but are handled differently according to the current queue lengths.
\begin{enumerate}
\item If $\qh+\ch < n$ and $\qy+\cy < n$,
\begin{equation}
\label{eq:arr1}
(h,\qh+1,\ch,\qy+1,\cy) \qquad \text{with rate }\l,
\end{equation}
\item if $\qh+\ch = n$ and $\qy+\cy < n$,
\begin{equation}
\label{eq:arr2}
(h+1,\qh,\ch,\qy+1,\cy) \qquad \text{with rate }\l,
\end{equation}
\item if $\qh+\ch < n$ and $\qy+\cy = n$,
\begin{equation}
\label{eq:arr3}
(h,\qh+1,\ch,\qy,\cy) \qquad \text{with rate }\l,
\end{equation}
\item if $\qh+\ch = n$ and $\qy+\cy = n$,
\begin{equation}
\label{eq:arr4}
(h+1,\qh+1,\ch,\qy,\cy) \qquad \text{with rate }\l,
\end{equation}
\end{enumerate}
\noindent \textbf{Departure.}
Basically, we align service completions in the two models, but allow a completion occurring solely in either of one of the two models, only if the queue length in this model is strictly larger than in the other one.
\begin{enumerate}
\item If $\qh \geq \qy$ and $h > 0$
\begin{equation}
\label{eq:dep1}
\left\{
\begin{array}{ll}
(h-1,\qh,\ch,\qy-1,\cy) & \text{with rate }(\qy \wedge s)(1-p)\mu,\\
(h-1,\qh,\ch,\qy,\cy) & \text{with rate }[(\qh \wedge s)-(\qy \wedge s)](1-p)\mu.
\end{array}
\right.
\end{equation}

\item If $\qh < \qy$ and $h > 0$
\begin{equation}
\label{eq:dep2}
\left\{
\begin{array}{ll}
(h-1,\qh,\ch,\qy-1,\cy) & \text{with rate }(\qh \wedge s)(1-p)\mu,\\
(h,\qh,\ch,\qy-1,\cy) & \text{with rate }[(\qy \wedge s)-(\qh \wedge s)](1-p)\mu.
\end{array}
\right.\end{equation}

\item If $\qh \geq \qy$ and $h = 0$
\begin{equation}
\label{eq:dep3}
\left\{
\begin{array}{ll}
(0,\qh-1,\ch,\qy-1,\cy) & \text{with rate }(\qy \wedge s)(1-p)\mu,\\
(0,\qh-1,\ch,\qy,\cy) & \text{with rate }[(\qh \wedge s)-(\qy \wedge s)](1-p)\mu.
\end{array}
\right.\end{equation}

\item If $\qh < \qy$ and $h = 0$
\begin{equation}
\label{eq:dep4}
\left\{
\begin{array}{ll}
(0,\qh-1,\ch,\qy-1,\cy) & \text{with rate }(\qh \wedge s)(1-p)\mu,\\
(0,\qh,\ch,\qy-1,\cy) & \text{with rate }[(\qy \wedge s)-(\qh \wedge s)](1-p)\mu.
\end{array}
\right.\end{equation}
\end{enumerate}

\noindent\textbf{Become content.} 
The differentiation between transitions is similar to those in the \textit{departure} transition type.
\begin{enumerate}
\item If $\qh \geq \qy$,
\begin{equation}
\label{eq:con1}
\left\{
\begin{array}{ll}
(h,\qh-1,\ch+1,\qy-1,\cy+1) & \text{with rate }(\qy \wedge s)p\mu,\\
(h,\qh-1,\ch+1,\qy,\cy) & \text{with rate }[(\qh \wedge s)-(\qy \wedge s)]p\mu.
\end{array}
\right.\end{equation}

\item If $\qh < \qy$,
\begin{equation}
\label{eq:con2}
\left\{
\begin{array}{ll}
(h,\qh-1,\ch+1,\qy-1,\cy+1) & \text{with rate }(\qh \wedge s)p\mu,\\
(h,\qh,\ch,\qy-1,\cy+1) & \text{with rate }[(\qy \wedge s)-(\qh \wedge s)]p\mu.
\end{array}
\right.\end{equation}
\end{enumerate}
\noindent
\textbf{Become needy.}
 \begin{enumerate}
\item If $\ch \geq \cy$,
\begin{equation}
\label{eq:ne1}
\left\{
\begin{array}{ll}
(h,\qh+1,\ch-1,\qy+1,\cy-1) & \text{with rate } \cy\delta,\\
(h,\qh+1,\ch-1,\qy,\cy) & \text{with rate }(\ch-\cy)\delta,\\
\end{array}
\right.\end{equation}

\item If $\ch < \cy$,
\begin{equation}
\label{eq:ne2}
\left\{
\begin{array}{ll}
(h,\qh+1,\ch-1,\qy+1,\cy-1) & \text{with rate } \ch\delta,\\
(h,\qh,\ch,\qy+1,\cy-1) & \text{with rate }(\cy-\ch)\delta,\\
\end{array}
\right.\end{equation}
\end{enumerate}

This set of transitions defines the dynamics of the Markov process $Z(t) = (\hat{X}^h(t),\hat{X}^b(t))$.
Let us now restrict our attention to the transitions in which (at least one of) the first three coordinates of $Z(t)$ changes, that is, the marginal transitions of the process $\hat{X}^h$. 
Let $\hat{X}^h(t) = (h,\qh,\ch)$, then according to the transition scheme above, $\hat{X}^h$ moves to state

\begin{enumerate}
\item If $\qh+\ch < n$ (and hence necessarily $h=0$),
\[
\left\{
\begin{array}{ll}
(0,\qh+1,\ch) & \text{with rate } \l,\\
(0,\qh-1,\ch) & \text{with rate }(\qh\wedge s)(1-p)\mu,\\
(0,\qh-1,\ch+1) & \text{with rate }(\qh\wedge s)p\mu,\\
(0,\qh+1,\ch-1) & \text{with rate }\ch \delta.
\end{array}
\right.\]

\item if $\qh+\ch = n$ and $h=0$,
\[
\left\{
\begin{array}{ll}
(1,\qh,\ch) & \text{with rate } \l,\\
(0,\qh,\ch) & \text{with rate }(\qh\wedge s)(1-p)\mu,\\
(0,\qh-1,\ch+1) & \text{with rate }(\qh\wedge s)p\mu,\\
(0,\qh+1,\ch-1) & \text{with rate }\ch \delta.
\end{array}
\right.\]

\item if $h>0$ (and hence necessarily $\qh+\ch = n$),
\[
\left\{
\begin{array}{ll}
(h+1,\qh,\ch) & \text{with rate } \l,\\
(h-1,\qh,\ch) & \text{with rate }(\qh\wedge s)(1-p)\mu,\\
(h,\qh-1,\ch+1) & \text{with rate }(\qh\wedge s)p\mu,\\
(h,\qh+1,\ch-1) & \text{with rate }\ch \delta.
\end{array}
\right.\]
\end{enumerate}
One can check that these transitions indeed coincide with the transitions in the original holding model, hence $\hat{X}^h(t) \equalD X^h(t)$.

Similarly, when focusing on transitions of $Z(t)$ that are relevant for $\hat{X}^b(t)$, we deduce the following transition scheme. If $\hat{X}^b(t) = (\qy,\cy)$, then the next move according to the transitions of $Z(t)$ is
\[
\left\{
\begin{array}{ll}
(\qy+\mathbbm{1}_{\{\qy + \cy < n\}},\cy) & \text{with rate } \l,\\
(\qy-1,\cy) & \text{with rate }(\qy\wedge s)(1-p)\mu,\\
(\qy-1,\cy+1) & \text{with rate }(\qy\wedge s)p\mu,\\
(\qy+1,\cy-1) & \text{with rate }\cy \delta.
\end{array}
\right.\]
These transition rates clearly coincide with the original Erlang-R model with blocking, and also hence $\hat{X}^b(t) \equalD X^h(t)$.

Next, we show that under this coupling scheme we have that if $\hat{H}(0) = 0$,  $\hat{Q}_1^h(0)=\hat{Q}_1^b(0)$ and $\hat{Q}_2^h(0)=\hat{Q}_2^b(0)$ then for all $t\geq 0$, $Z(t)$ satisfies the hypothesis:
\begin{itemize}
\item[(i)] $\hat{Q}_1^b(t) + \hat{Q}_2^b(t) \leq \hat{Q}_1^h(t) + \hat{Q}_2^h(t)$,
\item[(ii)] $\hat{Q}_2^b(t) \leq \hat{Q}_2^h(t)$,
\item[(iii)] $\hat{Q}_1^b(t) \leq \hat{Q}_1^h(t) + H(t)$.
\end{itemize}
We do so by induction on the next state reached after a transition of the joint Markov process $Z=(\hat{X}^h,\hat{X}^b)$.
First of all, $Z(0)$ clearly satisfies (i)-(iii). 
Next, assume $Z(t^-) = (h,\qh,\ch,\qy,\cy)$ satisfies the hypothesis and a transition occurs at $t$. 
We show that under the specified coupling scheme, the state reached after the next transition, $Z(t)$ must satisfy (i)-(iii) as well. To do so, we differentiate between the four types of transitions that could occur: arrival, departure, become content and become needy.\\
\\*
\noindent\textbf{Arrival.}
Recall that under our coupling scheme an arrival always occurs in both the holding and blocking model simultaneously, see \eqref{eq:arr1}--\eqref{eq:arr4}. Furthermore, $\ch$ and $\cy$ are unchanged during this transition, rendering (ii) trivial.

By hypothesis $\qy + \cy \leq \qh+\cy$, hence the event $\qh+\ch < n$ and $\qh+\cy =n$, with resulting state $(0,\qh+1,\ch,\qy,\cy)$,  can be excluded from our analysis.
We check the conditions for the remaining three cases.

\begin{enumerate}[noitemsep]
\item If $Z(t)= (0,\qh+1,\ch,\qy+1,\cy)$, then $\qy + \cy<n$ and $\qh + \ch <n$.
\begin{itemize}[noitemsep]
\item[(i)] $\QY +\CY = \qy+\cy+1 \less[i] \qh+\ch+1 =\QH+\CH$.
\item[(iii)] $\QY  = \qy+1 \less[iii] \qh+1 = \QH = \QH+\hat{H}(t)$.
\end{itemize}

\item If $Z(t)= (h+1,\qh,\ch,\qy+1,\cy)$, then $\qy + \cy<n$ and $\qh + \ch =n$.
\begin{itemize}[noitemsep]
\item[(i)] $\QY +\CY = \qy+\cy+1 \leq n = \qh+\ch =\QH+\CH$.
\item[(iii)] $\QY = \qy+1 \less[iii] \qh +1= \QH +\hat{H}(t)$.
\end{itemize}

\item If $Z(t)= (h+1,\qh,\ch,\qy,\cy)$, then $\qy + \cy = \qh+\ch=n$.
\begin{itemize}[noitemsep]
\item[(i)] $\QY +\CY = \qy+\cy \less[i] \qh+\ch =\QH+\CH$.
\item[(iii)] $\QY = \qy \less[iii] \qh+h < \qh+h+1 = \hat{H}(t)$.
\end{itemize}
\end{enumerate}

\noindent\textbf{Departure.}
By carefully examining the possible state transitions of $Z(t)$ following a departure, we list six reachable states. However, by (iii), we have that if $h=0$, then $\qy \leq \qh$, which excludes the state $(0,\qh,\ch,\qy,\cy)$ in \eqref{eq:dep4} from the reachability graph. 
We check the remaining states for conditions (i)--(iii). Again, during a departure, $\cy$ and $\ch$ are unchanged, so (ii) is automatically satisfied by the induction hypothesis.
\begin{enumerate}[noitemsep]
\item If $Z(t) = (h-1,\qh,\ch,\qy-1,\cy)$, then $h>0$. 
\begin{itemize}
\item[(i)] $\QY+\CY = \qy+\cy-1 \less[i] \qh+\ch-1 < \qh+\ch = \QH+\CH$. 
\item[(iii)] $\QY = \qy-1 \less[iii] \qh + h-1 = \QH +\hat{H}(t)$.
\end{itemize}

\item If $Z(t) = (h-1,\qh,\ch,\qy,\cy)$, then $h>0$ and $\qh \geq \qy$ (*). 
\begin{itemize}
\item[(i)] $\QY+\CY = \qy+\cy \less[i] \qh+\ch = \QH+\CH$.
\item[(iii)] $\QY = \qy \less[*] \qh-1 \leq \qh+h-1 = \QH+\hat{H}(t)$.
\end{itemize}

\item If $Z(t) = (h,\qh,\ch,\qy-1,\cy)$, then $h>0$ and $\qh < \qy$. 
\begin{itemize}
\item[(i)] $\QY+\CY = \qy+\cy-1 < \qy+\cy \less[i] \qh+\ch = \QH+\CH$.
\item[(iii)] $\QY = \qy-1 < \qy \less[iii] \qh + h = \QH+\hat{H}(t)$.
\end{itemize}

\item If $Z(t) = (h,\qh-1,\ch,\qy-1,\cy)$, then $h=0$. 
\begin{itemize}
\item[(i)] $\QY+\CY = (\qy-1)+\cy-1 <  \less[i] (\qh-1)+\ch = \QH+\CH$.
\item[(iii)] $\QY  = \qy-1 \less[iii] \qh-1 = \QH + \hat{H}(t)$.
\end{itemize}

\item If $Z(t) = (0,\qh-1,\ch,\qy,\cy)$, then $h=0$ and $\qh>\qy$ (*). 
\begin{itemize}
\item[(i)] $\QY+\CY = \qy+\cy  \less[i] (\qh-1)+\cy \less[ii] (\qh-1)+\ch = \QH+\CH$.
\item[(iii)] $\QY = \qy \less[*] \qh-1 =\QH+ \hat{H}(t)$.
\end{itemize}

\end{enumerate}

\noindent\textbf{Content start.}
On the event of a patient becoming content, it is clear that the sums $\QH+\CH$ and $\QY+\CY$ and $H(t)$ are unaffected. This means that (i) is directly satisfied by the induction hypothesis.
According to \eqref{eq:con1}--\eqref{eq:con2}, three states can be reached.  
\begin{enumerate}[noitemsep]
\item If $Z(t) = (h,\qh-1,\ch+1,\qy-1,\cy+1)$,
\begin{itemize}[noitemsep]
\item[(ii)] $\CY = \cy+1 \less[ii] \ch+1 = \CH$. 
\item[(iii)] $\QY = \qy-1 \less[iii] \qh+h-1 = \QH+\hat{H}(t)$.
\end{itemize}
\item If $Z(t) = (h,\qh-1,\ch+1,\qy,\cy)$, then $\qh > \qy$,
\begin{itemize}[noitemsep]
\item[(ii)] $\CY = \cy \less[ii] \ch < \ch+1 = \CH$. 
\item[(iii)] $\QY = \qy \less[iii] \qh+h < \qh+1+h = \QH+ \hat{H}(t)$.
\end{itemize}
\item If $Z(t) = (h,\qh,\ch,\qy-1,\cy+1)$, then $\qy > \qh$ (*) and hence by (iii) $h > 0$. The latter is only possible if $\qh+\ch=n$,
\begin{itemize}[noitemsep]
\item[(ii)] $\CY = \cy+1 \leq n-\qy+1 = (\qh+\ch)-\qy+1 \less[*] \ch = \CH$. 
\item[(iii)] $\QY = \qy-1 < \qh+h-1 \less[*] \qh+h = \QH+\hat{H}(t)$.
\end{itemize}
\end{enumerate}

\noindent \textbf{Become needy.}
Just as in the event of content start, the sums $\QH+\CH$ and $\QY+\CY$ and $H(t)$ are unaffected, whereby  (i) is directly satisfied by the induction hypothesis.
By (ii), we have $\ch \geq \cy$. This excludes the state $(h,\qh,\ch,\qy+1,\cy-1)$ from being reached, see \eqref{eq:ne2}. 
We check the remaining two possibilities.

\begin{enumerate}
\item If $Z(t) = (h,\qh+1,\ch-1,\qy+1,\cy-1)$. 
\begin{itemize}[noitemsep]
\item[(ii)] $\CY = \cy-1 \less[ii] \ch-1 = \CH$. 
\item[(iii)] $\QY = \qy+1 \less[iii] \qh+h+1 = \QH+\hat H(t)$.
\end{itemize}
\item If $Z(t) = (h,\qh+1,\ch-1,\qy,\cy)$, then $\ch > \cy$ (*). 
\begin{itemize}
\item[(ii)] $\CY = \cy \less[*] \ch-1 = \CH$.
\item[(iii)] $\QY = \qy \less[iii] \qh+h < \qh+1+h =\QH + \hat H(t)$.
\end{itemize}
\end{enumerate}

Hence, the state reached after any feasible transition under the coupling scheme satisfies the conditions (i)--(iii).
Thus we conclude that the joint process\\  $(\hat{H}(t),\QH,\CH,\QY,\CY)$ adheres to (i)--(iii) for all $t$. Consequently, we have that (i) implies
\begin{align*}
\P\left(Q_1^b(t) + Q_2^b(t) \geq k\right) &= \P\left(\hat{Q}_1^b(t) + \hat{Q}_2^b(t) \geq k\right)\\
&=\sum_{j=0}^n \P\left( \QY + \CY \geq k , \QH+\CH = j \right) \\
&=\sum_{j=k}^n \P\left( \QY + \CY \geq k , \QH+\CH = j \right) \\
&\leq \sum_{j=h}^n \P\left( \QH+\CH = j \right)\\
&= \P\left( Q_1^h(t) + Q_2^h(t) \geq k\right) = \P\left(Q_1^h(t) + Q_2^h(t) \geq k\right).
\end{align*}
The other two orderings follow similarly.

\begin{remark}
Note that under this coupling scheme we cannot get the ordering $\QH(t) \geq \QY(t)$ for all $t\geq 0$. A minimal counter example occurs for $s=n=1$. Let $Z(0) = ((0,0,0),(0,0))$. First, two arrivals occur, such that state $((1,1,0),(1,0))$ is reached, followed by a departure transition, yielding $((0,1,0),(0,0))$. Next, the one patient left in the model with holding system becomes content, so that we obtain $((0,0,1),(0,0))$. 
At this stage, if an arrival occurs, the arriving patient will be put in the holding queue in the model with holding, and admitted to nurse queue in the model with blocking. Hence we end up in state $((1,0,1),(1,0))$, in which $\QH < \QY$. 
\end{remark}

\section{Proof of Proposition \ref{prop:stability_convergence}}\label{app:proof_stability_convergence}
 
Define 
\[
A(s,n) = \sum_{k=0}^s \frac{k}{s} \, \binom{n}{k} b^k ,\quad 
B(s,n) = \sum_{k=s+1}^n \frac{k!}{s!} \, \binom{n}{k} s^{s-k} b^k, \quad
C(s,n) = \sum_{k=0}^s \binom{n}{k} \, b^k,
\]
\[
\]
where $b = \delta/p\mu = r/(1-r)$. Then
\[
\rho_{\rm max}(s,n) = \frac{A(s,n)+B(s,n)}{C(s,n)+B(s,n)}.
\]
Proving that $\rho_{\rm max}(s,n) \to 1$ as $R_1\to\infty$ with $s$ and $n$ as in \eqref{eq:twofoldscaling} is equivalent to showing that
\begin{equation}\label{eq:proof_stab_1}
1-\rho_{\rm max}(s,n) = \frac{C(s,n)-A(s,n)}{C(s,n)+B(s,n)} = \frac{(1+b)^{-n}[C(s,n)-A(s,n)]}{(1+b)^{-n}[C(s,n)+B(s,n)]} \to 0.
\end{equation}
First, we rewrite
\begin{align*}
(1+b)^{-n} A(s,n) 
&= (1+b)^{-n} \sum_{k=1}^s \frac{n}{s} \binom{n-1}{k-1} b^k  \\
&= \frac{n}{s}\left(\frac{b}{1+b}\right)\sum_{k=0}^{s-1} \binom{n-1}{k} \left(\frac{b}{1+b}\right)^k \left(\frac{1}{1+b}\right)^{n-1-k}\\
&= \frac{r n}{s}\sum_{k=0}^{s-1} \binom{n-1}{k} r^k (1-r)^{n-1-k}\\
&= \frac{r n}{s} \P( {\rm Bin}(n-1,r) \leq s-1 ) \\
&= \frac{rn}{s} \P\left( \frac{{\rm Bin}(n-1,r) - (n-1)r}{\sqrt{nr(1-r)}} \leq \frac{s-1 - (n-1)r}{\sqrt{nr(1-r)}} \right)\\
&\to \Phi\left(\frac{\beta-\gamma\sqrt{r}}{\sqrt{1-r}}\right),
\end{align*}
since $nr/s = 1 + O(1/\sqrt{R_1})$.
Also, 
\begin{align*}
(1+b)^{-n} C(s,n) 
&= \sum_{k=0}^s \binom{n}{k} \left(\frac{b}{1+b}\right)^k \left(\frac{1}{1+b}\right)^{n-k}\\
&= \sum_{k=0}^s \binom{n}{k} r^k (1-r)^{n-k}\\
&= \P( {\rm Bin}(n,r) \leq s) \to \Phi\left(\frac{\beta-\gamma\sqrt{r}}{\sqrt{1-r}}\right).
\end{align*}
Therefore, we have $(1+b)^{-n}[C(s,n)-A(s,n)] \to 0$ as $\l\to\infty$.
For the remaining term,

\begin{align*}
(1+b)^{-n} B(s,n) 
&= (1+b)^{-n}\sum_{k=s+1}^n \binom{n}{k}\,\frac{k!}{s!} s^{s-k} b^k  \\
&= (1+b)^{-n}\frac{n!}{s!}\, s^s\sum_{k=s+1}^n \frac{1}{(n-k)!} \left(\frac{s}{b}\right)^{-k}\\
&= (1+b)^{-n} \frac{n!}{s!}\, s^s\, \left(\frac{b}{s}\right)^n \sum_{k=s+1}^n \frac{1}{(n-k)!} \left(\frac{s}{b}\right)^{n-k}\\
&=  r^n\, \frac{n!}{s!} s^{s-n} \sum_{m=0}^{n-s-1} \frac{1}{m!} \left(\frac{s}{b}\right)^m\\
&= \left(\frac{r}{s}\right)^n \frac{n!}{s!} s^s \,{\rm e}^{s/b} \, \P({\rm Pois}(s/b)\leq n-s-1),
\end{align*}
in which 
\begin{align*}
\P({\rm Pois}(s/b)\leq n-s-1) 
&= \P\left(\frac{{\rm Pois}(s/b)-s/b}{\sqrt{s/b}} \leq \frac{n-s-1-s/b}{\sqrt{s/b}}\right) \\
&\to \Phi\left(\frac{\gamma-\beta/\sqrt{r}}{\sqrt{1-r}}\right),
\end{align*}
as $\l\to\infty$.
By Stirling's approximation,
\begin{align*}
\left(\frac{r}{s}\right)^n \frac{n!}{s!} s^s \,{\rm e}^{s/b} 
&\sim \left(\frac{r}{s}\right)^n \sqrt{\frac{n}{s}} \,\frac{n^n {\rm e}^{-n}}{s^s {\rm e}^{-s}}\, s^s \,{\rm e}^{s/b} \\
&= \left(\frac{rn}{s}\right)^n \sqrt{\frac{n}{s}} {\rm e}^{-n+s+s/b} = \left(\frac{rn}{s}\right)^n \sqrt{\frac{n}{s}} {\rm e}^{-n+s/r}.
\end{align*}
Since, 
\[
\frac{rn}{s} = 1 + \frac{\gamma\sqrt{r}-\beta}{\sqrt{R_1}} + O(1/R_1),
\]
we find $\sqrt{n/s} = 1/\sqrt{r} + O(1/\sqrt{R_1})$ and 
\begin{align*}
\log\left[ \left(\frac{rn}{s}\right)^n \sqrt{\frac{n}{s}} {\rm e}^{-n+\tfrac{s}{r}} \right]
&= n \log\left[ \frac{rn}{s}\right] - n+\frac{s}{r}\\
&= -n \left[ \left(1-\frac{rn}{s}\right) + \frac{1}{2}\left(1-\frac{rn}{s}\right)^2 + O(R^{-\tfrac{3}{2}}) \right] + \frac{s}{r}\left(1-\frac{rn}{s}\right)\\
&= \frac{s}{r}\left(1-\frac{rn}{s}\right)^2 - \frac{n}{2}\left(1-\frac{rn}{s}\right)^2 + O(1/\sqrt{R_1})\\
&= \frac{(\gamma\sqrt{r} - \beta)^2}{2r} + O(1/\sqrt{R_1}),
\end{align*}
as $\l\to\infty$ and hence, 
\[
(1+b)^{-n} B(s,n) \to \f\left(\frac{\gamma\sqrt{r}-\beta}{\sqrt{r}}\right)\Phi\left(\frac{\gamma-\beta/\sqrt{r}}{\sqrt{1-r}}\right).
\]
Hence, we conclude that the denominator of \eqref{eq:proof_stab_1} converges to a constant value as $R_1$ grows, and hence the $1-\rho_{\rm max}(s,n)\to 0$ as $\l\to\infty$.

\resettocdepth

\end{subappendices}

\chapter{Transient error approximation in a L\'evy queue}

\begin{chapterstart}
Motivated by a capacity allocation problem within a finite planning period, we conduct a transient analysis of a single-server queue with L\'evy input. From a cost minimization perspective, we investigate the error induced by using stationary congestion measures as opposed to time-dependent measures. Invoking recent results from fluctuation theory of L\'evy processes, we derive a refined cost function, that accounts for transient effects. This leads to a corrected capacity allocation rule for the transient single-server queue. Extensive numerical experiments indicate that the cost reductions achieved by this correction can be significant.
\end{chapterstart}

\begin{flushright}
Based on\\
\textbf{Transient error approximation in a L\'evy queue}\\
\textit{Britt Mathijsen \& Bert Zwart}\\
\textit{Queueing Systems, 85(3), 269-304 (2017)}
\end{flushright}
\newpage

\section{Introduction}

The issue of matching a service system's capacity to stochastic demand induced by its clients arises in many practical settings. Typically, the resources available to satisfy demand are scarce and hence expensive. This forces the manager to consider a trade-off between the system efficiency and the quality of service perceived by its clients. In this chapter, we focus on this trade-off in the context of the $M/G/1$ queue, in which the variable amenable for optimization is the server speed $\mu$.

In general, optimizing the server speed $\m$ in a single-server queue in a time-homogeneous environment, while trading off congestion levels against capacity allocation costs, does not pose any technical challenges. Typically, the objective function to be minimized, the total cost function, has the shape
\begin{equation}\label{eq:intro}
\Pi_\iy(\mu) = \E[Q_\mu(\infty)] + \aaa\mu = \frac{\la\E[B^2] }{2(\mu-\la\E[B])} + \aaa\mu,
\end{equation}
where $\E[Q_\mu(\infty)]$ denotes the expected steady-state amount of work given server speed $\m$, and $B$ describes the service requirement per arrival. The parameter $\aaa>0$ represents the relative capacity allocation costs incurred by deploying service rate $\mu$. This one-dimensional optimization problem yields the optimizer
\begin{equation*}
\mui = \lambda \E[B] + \sqrt{\frac{\la\E[B^2]}{2\aaa}}.
\end{equation*}
Despite the simplicity and tractability of the problem described above, the presence of the \emph{steady-state} measure in the cost function in \eqref{eq:intro} should be handled carefully. By employing this particular cost structure, one automatically agrees with the underlying assumption of the system being sufficiently close to its steady state. 
However, referring to practical applications of the single-server model, system parameters rarely remain constant over time. Moreover, planning periods for the optimization problem are naturally finite. Hence, the \emph{true} expected costs incurred, which we denote by $\Pi_T(\mu)$, in addition depend on the length of the planning period $T$. Consequently, the usage of steady-state models for decision making needs to be justified by a more elaborate time-dependent or \emph{transient} analysis for these type of settings.\\
\\*
\noindent
\textbf{Related literature}. 
The time-dependent behavior of the single-server queue received much attention in queueing theory. First efforts to analyze the time-dependent properties of the $M/G/1$ queue date back to the 1950s and 1960s, e.g. \cite{Benes1957,Gaver1959,Kendall1951,Takacs1955,Takacs1962}. The analyses in these papers mostly yield implicit expressions for performance characteristics through Laplace transforms, integro-differential equations and infinite convolutions.
More specifically, there is vast literature on the transient analysis of the $M/M/1$ queue, with the goal to derive explicit expressions for queue length characteristics, see e.g. \cite{Abate1987,Cohen1982,Pegden1982,Prabhu1964}. 
These works provide a variety of explicit expressions for the transient dynamics, although the complexity of the resulting expressions, typically involving Bessel functions, exposes the intricate intractability of the matter. Consequently, approximation methods for insightful quantification of the dynamics based on numerical \cite{Neuts1966} or asymptotic methods, have become prevalent in more recent literature.
The asymptotic methods either exploit knowledge on the evolution of the queueing process as time $t$ grows large \cite{Abate1987,Newell1982,Odoni1983}, or as the arrival rate $\la$ is increased to infinity \cite{Abate1987a,Abate1987b,Gaver1968}. 
It is noteworthy that a substantial contribution to the transient literature is made by Abate and Whitt \cite{Abate1987a,Abate1987b,Abate1987,Abate1994}, who exploit the existence of a decomposition of the mean transient queue length and obtain expressions for the moments of the queue length and virtual waiting through probabilistic arguments in several queueing models. 
More recently, asymptotic methods have been used to justify the application of stationary performance measures in Markovian environments or to refine them, see e.g. \cite{Green1991,Whitt1991}.
Other approximative methods known as uniform acceleration expansions \cite{Massey1998} have been developed to reveal the asymptotic behavior of the single-server queue as a function of $t$, which are moreover able to capture time-varying arrival rates.
The majority of the works mentioned above do reflect on the error imposed by usage of steady-state performance metrics instead of the correct time-dependent counterpart. However, no light has been shed on the accumulation of this error over a finite period of time. To the best of our knowledge, the only work that addresses this issue is the paper by Steckley and Henderson \cite{Steckley2007}, who compute an approximation for the error accumulated between the steady-state and transient delay probability. Our analysis on the other hand is centered around the mean workload, which requires a different approach. In addition, the focus in \cite{Steckley2007} is on performance measures only, while the main goal of our work is to investigate the quality of staffing rules. \\
\\*
\noindent\textbf{L\'evy input}.
Although the $M/G/1$ queue serves as the leading example in our analysis, we choose to use a more general framework for the arrival process of the queue. Namely, we let the server face a L\'evy process.
This gives the advantage that once we have obtained the results, we can apply them to broader queue input classes, such as Brownian motion and the Gamma process.
To shed light on the influence of the transience of the queueing process on traditional staffing questions, we will study the capacity allocation problem in the context of cost minimization in which the objective function is $\Pi_T(\mu)$, i.e. a function of both $\mu$ and $T$. We investigate how the invalidity of the stationary assumption is echoed through the operational cost accounting for congestion-related penalties. 
Furthermore, we establish a result on the strict convexity of the function $\Pi_T(\mu)$, for almost all values of $T$ (with a few minor exceptions for certain deterministic initial states), which is an essential property for convergence of both cost function and corresponding minimizer to their stationary counterparts.
\\
\\*
\noindent\textbf{Corrected staffing rule}. 
As it will appear that an exact analysis of this disparity is intractable, we will present an explicit approximate correction to the conventional stationary objective function given by $\Psi(\mu)/T$ and prove that
\begin{equation*}
\Pi_T(\mu) = \Pi_\iy(\mu) + \frac{\Psi(\mu)}{T} + O(1/T^2),
\end{equation*}
with the help of recent results from the fluctuation theory of L\'evy processes. 
 Based on this refinement we ultimately examine how incorporating transient effects\\ \noindent changes the optimal capacity level and propose a refinement to the steady-state capacity allocation rule, 
\begin{equation*}
\muT = \mui + \frac{\mu_\bullet}{T} + o(1/T).
\end{equation*}
We moreover deduce an explicit expression for $\mu_\bullet$ in terms of the initial state and the first three moments of the service requirement per arrival.
It is noteworthy that similar refined square-root staffing rules have been proposed for multi-server queues in the Halfin-Whitt regime, see e.g. \cite{Janssen2015,Janssen2008,Janssen2011,Randhawa2014,Zhang2012}. In those cases, the relevant decision value is the number of servers  and refinements are derived for $\la\to\iy$, whereas we consider the regime $T\to\infty$. 

  Building upon the insights gained through the analysis of this optimality gap, we reflect on the parameter settings of the underlying queueing process in which our refined capacity sizing rule yields significant improvement and in which cases it has little effect. Special emphasis is put on the relationship between the accuracy of the standard procedure and the length of the planning period.
\\
\\*
\noindent\textbf{Structure of the chapter}.
The remainder of this chapter is structured as follows. Section \ref{sec:model_description} is devoted to the model description and presents some preliminary results. The main result will be given in Section \ref{sec:analysis} and results regarding the optimization problem will be discussed in Section \ref{sec:optimization}, followed by the validation of our novel techniques through numerical experiments in Section \ref{sec:numerics}. We will give some concluding remarks and topics for further research in Section \ref{sec:conclusion_chapter6}. We have deferred all proofs to the appendix.

\section{Model description}
\label{sec:model_description}

\subsection{A queueing model with L\'evy input \label{sec:levymodel}}
The model that inspired our study is the standard $M/G/1$ queue starting out of equilibrium. Customers arrive to the queue according to a Poisson process with rate $\la$. 
All arrivals have i.i.d. service requirement $B_i$, stemming from a common random variable $B$. 
Without loss of generality we will assume $\E[B] = 1$ throughout. The server is able to remove $\mu$ amounts of work from the system per time unit; a variable we will refer to as the \emph{server speed}. 
E.g. if $\mu = 3$ and two customers are in the system with remaining service times $4$ and $2$, then  the queue will be empty 2 time units later, provided that no new arrivals occur in the meantime.
Let $N_\la(t)$ denote the number of arrivals until time $t$.
Accordingly, the total work generated by the customers is given by
\begin{equation*}
Z_\la(t) = \sum_{i=1}^{N_\la(t)} B_i.
 \end{equation*} 
Furthermore, define $X_{\la,\mu}(t) := Z_\la(t) - \mu t$. We call $X_{\la,\mu}$ the \emph{net-input process}. 
More generally, we assume throughout the chapter that $X_{\la,\mu}$ is a L\'evy process.
Specifically, we let $Z_\la$ be of the form $Z_\la(t) = U(\la t)$, where $U$ is a spectrally positive L\'evy process generated by the triplet $(a,\s,\nu)$ and $\E[U(1)] = 1$. 
This restriction to spectrally positive processes is equivalent to stating $\nu(-\infty,0)=0$ and is a vital assumption to our analysis. 
Subsequently, we assume the net-input process $X_{\la,\mu}$ to be
\begin{equation}
\label{eq:Xlmprocess}
X_{\la,\mu}(t) = U(\la t) - \mu t, \qquad t \geq 0.
\end{equation}
Note that by setting $a=\s=0$ and $\nu = \la\, F_B$, where $F_B$ is the cumulative distribution function of $B$, we retrieve the original $M/G/1$ queue. 
The stochastic process central to our analysis is the \emph{workload process} $Q_{\la,\mu}(t)$, $t\geq 0$, which describes the amount of work the server is facing at time $t$. 
The net-input process $X_{\la,\mu}$ completely determines the trajectory of $Q_{\la,\mu}$, namely
\begin{equation}\label{eq:Qlm}
Q_{\la,\mu}(t) = \max\left\{ Q(0) + X_{\la,\mu}(t), \sup_{s\in[0,t]} [X_{\la,\mu}(t)-X_{\la,\mu}(s)]\right\}, \qquad t\geq 0,
\end{equation}
where $Q(0)$ is the initial workload in the system. 
In fact, $Q_{\la,\mu}$ is the reflected version of $X_{\la,\mu}$ with reflection barrier at zero.
Careful inspection of the structure also reveals that $X_{\la,\m}(t) \equiv X_{\la/\mu,1}(\mu t) \equiv X_{1,\mu/\la}(\la t)$, so that
\begin{equation}
\label{eq:Qidentity}
Q_{\la,\m}(t) \equalD Q_{\la/\mu,1}(\mu t) \equalD Q_{1,\mu/\la}(\la t)
\end{equation}
for all $\la,\m,t>0$.
This identity will prove to be convenient for the numerical analysis in Section \ref{sec:numerics}. For reasons of clarity, we omit the subscript $\la$ in our expressions if no ambiguity is possible.

The process $Q_{\m}$ is a natural indicator of the level of congestion in the system and therefore a good choice for quantifying the Quality of Service (QoS) received by a client.
We remark that alternative processes characterizing congestion in the system can be deduced directly from $Q_{\m}(t)$. For example, consider the virtual waiting time process $V_{\mu}(t)$, which is the waiting time a customer would experience if he arrives at time $t$. This, under the first-come-first-served policy, satisfies the relation $\E[V_{\mu}(t)] = \E[Q_{\m}(t)]/\mu$ for all $t\geq 0$.
Likewise, the expected number of the customers in the system $L_{\m}(t)$ at time $t\geq 0$ is given by Little's law
\begin{equation*}
\E[L_{\m}(t)] = \la\, \E[V_{\m}(t)] = \frac{\la}{\mu}\, \E[Q_{\mu}(t)].
\end{equation*}
To facilitate our investigation of the queueing model, we end this subsection by introducing some notation regarding the net-input and workload process and by stating a useful preliminary result concerning the stationary process $\Qlm(\iy)$. 
Throughout the chapter we assume $\mu>\la$ to ensure ergodicity of the queue and convergence in distribution to the limit
\begin{equation*}
\Qlm(\iy) := \lim_{t\to\iy} \Qlm(t),
\end{equation*}
for any initial state $Q(0)<\iy$.
The distribution of $Q_{\mu}(\infty)$ coincides with the stationary distribution of $\Qlm(t)$. 
By $\ka_U(\cdot)$ and $\ka_{\mu}(\cdot)$ we denote the L\'evy exponents of the processes $U$ and $\Xlm$, respectively:
\begin{equation*}
\ka_{\m}(\thh) = \log \E[\ee^{\thh \Xlm(1)}] = \log \E[\ee^{\thh(U(\la) - \mu)}] = \la \ka_U(\thh) - \mu \thh.
\end{equation*}
Furthermore, define $u_k = \E[\{U(1) - \E U(1)\}^k]$ for $k=2,3,...$. 
Using this representation we obtain the following preliminary result.
\begin{lemma}\label{lemma:workloadmoments}
Let $\E|U(1)|<\infty$, $u_2, u_3 < \iy$ and $\mu > \la$. If $Q_{\mu}(\infty)$ represents the steady-state distribution of the workload process, then
\begin{equation*}
\E[\Qlm(\infty)] = \frac{\la u_2}{2(\mu-\la)},\qquad \E[Q_{\mu}^2(\iy)]=\frac{\la^2u_2^2}{2(\mu-\la)^2} + \frac{\la u_3}{3(\mu-\la)}.
\end{equation*}
\end{lemma}
The proof of Lemma \ref{lemma:workloadmoments} follows directly by differentiation of the Laplace transform of $Q_\mu(\iy)$ and is given in Appendix \ref{app:proof_lemma_workload_moments}.
\subsection{Finite horizon}

\begin{figure}
\centering
\begin{tikzpicture}[xscale=0.15,yscale=0.225]
\draw (0,0) -- coordinate (x axis mid) (50,0);
    \draw (0,0) -- coordinate (y axis mid) (0,21);
	\node[right] at (51,0) {$t$};
	\node[rotate=90, above=0.7 cm] at (y axis mid) {$\E[\Qlm(t)]$};
	
	\draw[dashed, thick, gray] (0,10) -- coordinate (eq) (51,10);

	\draw[->] (24,6.4) --  coordinate (a1) (21.65,8.49574);
	\node[right=0.6cm,below=0.3cm] at (a1) {$Q(0)\equiv 0$};
	\draw[->] (14,6.2) --  coordinate (a2) (13.,8.49434);
	\node[right=0.2cm,below=0.3cm] at (a2) {$Q(0)\equiv10$};
	\draw[->] (9,15.5) --  coordinate (a3) (7.5,13.7648);
	\node[right=0.6cm,above=0.1cm] at (a3) {$Q(0)\equiv20$};
	\draw[->] (40,12.5) -- coordinate (a4) (38,10.9712);
	\node[right,above=0.3cm] at (a4) { $Q(0)\sim \exp\left(\tfrac{1}{15}\right)$ };

    	\foreach \x in {0,10,...,50}
     		\draw (\x,1pt) -- (\x,-10pt)
			node[anchor=north] {\x};
    	\foreach \y in {5,10,15,20}
     		\draw (1pt,\y) -- (-20pt,\y)
     			node[anchor=east] {\y};

	\draw[thick,color = col1] plot
			file {Chapter_6/tikz/means0.txt};
	\draw[thick,color = col3] plot
			file {Chapter_6/tikz/means10.txt};
    \draw[thick,color = col4] plot
			file {Chapter_6/tikz/means20.txt};
	\draw[thick,color = col5] plot
			file {Chapter_6/tikz/meansExp.txt};		
\end{tikzpicture}
\caption{Time-dependent mean workload in a $M/M/1$ queue with $\la = 10$ and $\mu=11$ for different initial states $Q(0)$. The dashed line depicts $\E\Qlm(\iy)$.}
\label{fig:transientmeans}
\end{figure}
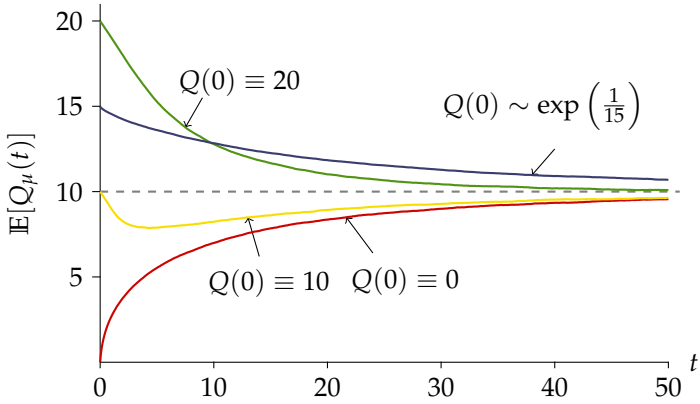

For the purpose of our research, we are interested in the dynamics of the workload process within a fixed time frame of length $T>0$. 
For all $0\leq t \leq T$, we assume that the parameters of the queue, $\la,\m,u_2,u_3$, remain unchanged. 
If at $t=0$ the queue is not in steady-state corresponding to the specified parameters of the starting period, the process $\{\Qlm(t)\}_{t\in[0,T]}$ differs from its stationary counterpart $\Qlm(\infty)$. 
To illustrate this, Figure \ref{fig:transientmeans} depicts the expected value $\Qlm$ in a $M/M/1$ queue as a function of time for several initial workloads $Q(0)$ for a particular setting of $\la$ and $\m$. 
Clearly, transient behavior of $\E[\Qlm(t)]$, for $Q(0) \neq \Qlm(\iy)$, differs significantly from the steady-state mean with the same system parameters. 
Note that even if $Q(0) \equiv \E[\Qlm(\iy)]$, the time-dependent mean does not coincide with the steady-state mean. Moreover, $\E[\Qlm(t)]$ is not even a strictly increasing nor decreasing function of time. This phenomenon is a consequence of the decomposition of the transient mean into one strictly increasing, and a strictly decreasing term for $Q(0)>0$, as discussed in \cite{Abate1987}.
Nonetheless, $\Qlm(t)$ converges in distribution to $\Qlm(\infty)$ as $t\to\iy$, if $\mu>\la$. 
 
Since the time horizon of our analysis is limited to $t\leq T$, the process may not approach the steady-state distribution sufficiently close to appropriately use its steady-state properties for capacity allocation. 
To overcome this disparity, we propose a way to include the influence of this transient phase in the capacity allocation problem.
 
\subsection{Cost structure}

As mentioned before, we are interested in balancing the QoS and efficiency of the queue by choosing the optimal server speed $\mu$. 
The adjective \emph{optimal} indicates that we intend to choose the speed according to some objective function. 
In our case, we conduct our analysis based on a cost function, which consists of a part accounting for the penalty for congestion in the system and a part for staffing cost. The cost value of both parts is governed by the variable $\mu$. 
The instantaneous cost incurred at time $t$ equals
\begin{equation*}
\E[\Qlm(t)] + \aaa \mu,
\end{equation*}
where $\aaa$ is a positive constant defining the \emph{relative staffing cost}. 
Hence, the cost structure we apply is a combination of the transient mean of the workload process and a linear staffing cost. 
Accumulated and normalized over the period $[0,T]$, the cost function on which the rest of this chapter will be based equals
\begin{equation}\label{eq:PiT}
\Pi_{T}(\mu) := \frac{1}{T}\int_0^T\left(  \E[\Qlm(t)] + \aaa\mu \right) \dd t 
 = \frac{1}{T} \int_0^T \E[\Qlm(t)] \dd t + \aaa\mu.
\end{equation}
We use shorthand notation for the normalized congestion costs:
\begin{equation}\label{eq:CTmu}
C_{T}(\mu) := \frac{1}{T}\int_0^T \E[Q_{\mu}(t)] \dd t,
\end{equation}
and $C_{\iy}(\mu) := \E[\Qlm(\iy)]$. 
In order to compare the actual costs incurred over the interval $[0,T]$ to the cost function of the queue in stationary conditions, we define
\begin{equation}\label{eq:PiInf}
\Pi_{\iy}(\mu) := C_{\iy}(\mu) + \aaa \mu = \E[Q_\mu(\iy)] + \aaa\mu,
\end{equation}
which allows an explicit expression by Lemma \ref{lemma:workloadmoments}. 
Under mild conditions on the net-input process and the distribution of the initial state, the cost functions coincide for $T\to\iy$.
\begin{proposition}\label{prop:cost_convergence}
Let $\mu>\la$ and assume $\E[U(1)],\, \E[Q(0)] < \iy$. Then
\begin{equation*}
\lim_{T\to\iy} \Pi_{T}(\mu) = \Pi_{\iy}(\mu).
\end{equation*} 
\end{proposition}
\noindent
The proof of Proposition \ref{prop:cost_convergence} can be found in Appendix \ref{app:proof_prop1}.
\noindent 
Define 
\[
\Omega_T := \frac{1}{T}\int_0^{T} \left( \E[\Qlm(t)] - \E[\Qlm(\iy)] \right) \dd t \]
We can then rewriting \eqref{eq:PiT} as
\begin{align}
\Pi_{T}(\mu) &= \frac{1}{T}\int_0^{T} \left( \E[\Qlm(t)] - \E[\Qlm(\iy)] \right) \dd t + \E[\Qlm(\iy)] + \aaa\mu =  \Omega_{T}(\mu) + \Pi_{\infty}(\mu).
\label{eq:decomp}
\end{align}
Section \ref{sec:analysis} is concerned with the analysis of the correction factor $\Omega_{T}(\mu)$.

Ultimately, we are concerned with the additional costs incurred by choosing the server speed through minimization of $\Pi_{\iy}(\mu)$ instead of $\Pi_{T}(\mu)$. 
Therefore, we formulate the exact and approximate optimization problems as follows
\begin{equation}\label{eq:muStar}
\mu_T^\star := \arg\min_{\mu\geq 0} \Pi_{T}(\mu), \qquad \qquad \mu_\infty^\star := \arg\min_{\mu\geq 0} \Pi_{\iy}(\mu),
\end{equation}
\begin{equation}\label{eq:piStar}
\Pi_{T}^\star := \Pi_{T}(\mu_T^\star), \qquad \qquad \Pi_{\iy}^\star := \Pi_{T}(\mu_\iy^\star).
\end{equation}
In Section \ref{sec:optimization} we turn to the comparison of $\mu_T^{\star}$ and $\mu_\iy^\star$ as well as the \emph{optimality gap} $\Pi_{\iy}^\star - \Pi_{T}^\star$.

\section{Analysis of the objective function}
\label{sec:analysis}
From \eqref{eq:decomp} it is evident that, for finding an explicit characterization of $\Pi_{T}(\mu)$, it suffices to study the term $\Omega_T(\mu)$ in more detail. We start by stating the main result of this section, which describes the leading order behavior of $\Omega_T(\mu)$ as $T$ increases.
\begin{theorem}\label{thm:mainresult}
Let $X_\mu(t)$ be of the form \eqref{eq:Xlmprocess}. If $\E[\max(Q(0),Q_\mu(\infty))^3] < \iy$ and $u_2,u_3 < \iy$, then
\begin{align*}
\Omega_T(\mu) &= \frac{\E[Q(0)^2] - \E[Q_\mu(\iy)^2]}{2T(\mu-\la)} + O\left(\frac{1}{T^2}\right) \nonumber\\
&= \frac{1}{2T(\mu-\la)}\left( \E[Q(0)^2] - \frac{\la^2 u_2^2}{2(\mu-\la)^2} - \frac{\la u_3}{3(\mu-\la)}\right) + O\left(\frac{1}{T^2}\right),
\end{align*} 
for $\mu>\la$.
\end{theorem}
Note that this expression provides an \emph{approximation} of the actual cost function
$\Pi_T(\mu)$. We elaborate on the implications of this additional information on the optimization problem in Section \ref{sec:optimization}. 

In the remainder of this section we provide a detailed description of the steps taken to obtain this outcome.
We assume a fixed service rate $\mu$ throughout the analysis in this section and therefore omit the subscript $\mu$. Proofs of the intermediate results can be found in Appendix \ref{app:proofs_analysis}.
\subsection{Constructing a coupling}

Before starting our analysis of the correction term $\Omega_{T}(\mu)$ we introduce some auxiliary notation. 
By $Q^A(t)$ we denote the workload process as described in Subsection \ref{sec:levymodel} with initial state $A$ and $\E_A$ the expectation with respect to any non-negative random variable $A$, which is independent of the net-input process $X$.
To be able to compare $\E[Q(t)]$ and $\E[Q(\iy)]$ as in $\Omega_T(\mu)$, we will use a coupling technique.
Observe that by definition of the stationary distribution $Q(\iy) \equalD Q^{Q(\iy)}(t)$ for all $t \geq 0$ and therefore $\E[Q(\iy)] = \E_{Q(\iy)}[Q^{Q(\iy)}(t)]$. Furthermore, $\E[Q(t)] = \E_{Q(0)}[Q^{Q(0)}(t)]$.
Hence, quantifying the difference between the transient and stationary mean is equivalent to comparing the workload processes of two queues starting in two different (random) states at $t=0$. 

We begin our analysis for two queues starting in two \textit{deterministic} states $x,y\geq 0$, respectively. At the end of our analysis we will obtain the original form by replacing $x$ with $Q(0)$ and $y$ with $Q(\iy)$.

Equation \eqref{eq:Qlm} shows that all randomness in the workload process originates from the process $X(t)$. 
With this in mind, we couple the processes $Q^x(t)$ and $Q^{y}(t)$ on a sample path level by feeding both queues the same net-input process $X(t)$ for $t\geq 0$. 
This allows us to compare the processes in the same probability space, so that $\E[Q^x(t)] - \E[Q^y(t)] = \E[Q^x(t) - Q^y(t)]$ for all $t\geq 0$.
Define
\begin{equation*}
Y^{x,y}(t) := Q^x(t) - Q^y(t)
\end{equation*}
and 
\begin{equation*}
\Omega_{T}^{x,y} := \frac{1}{T}\,\int_0^T \E\left[Y^{x,y}(t)\right] \, \dd t.
\end{equation*}

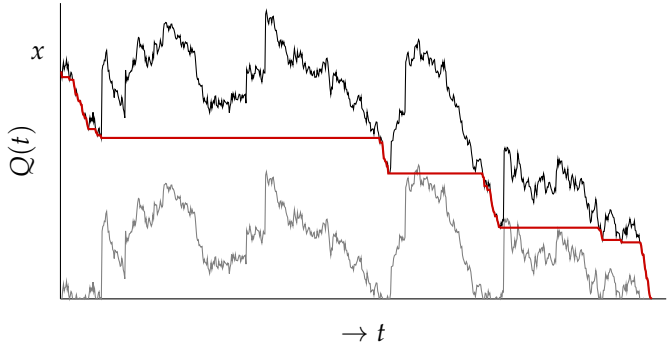
\begin{figure}
\centering
\begin{tikzpicture}[y=0.6cm, x=0.01cm]
	\draw (0,0) -- coordinate (x axis mid) (800,0);
    \draw (0,0) -- coordinate (y axis mid) (0,6.5);
	\node[below=0.2cm] at (x axis mid) {$\to t$};
	\node[rotate=90, above=0.2cm] at (y axis mid) {$Q(t)$};
	\node[above=1.3cm,left =0.08 cm] at (y axis mid) {$x$};

	\draw plot
			file {Chapter_6/tikz/samplePathLevy.txt};
	\draw[color = gray] plot
			file {Chapter_6/tikz/samplePathLevy2.txt};
	\draw[thick,color=col1] plot
			file {Chapter_6/tikz/runningMinimumLevy.txt};
\end{tikzpicture}
\caption{Sample path visualization of the processes $Q^x(t)$ (solid), $Q^0(t)$ (gray) and $Y^{x,0}(t)$ (red).}
\label{fig:samplePaths}
\end{figure}

A possible sample path triple for $Q^x(t)$, $Q^0(t)$ and $Y^{x,0}(t)$ is depicted in Figure \ref{fig:samplePaths}.
 As we see from this figure, $Y^{x,0}(t)$ has nice structural properties which we will exploit in the next subsection. 

\subsection{Difference process and leading order behavior of the correction term}
\label{sec:difference}
We further examine the \emph{difference process} $Y^{x,y}(t)$ with $x>y$. Recall from \eqref{eq:Qlm},
\begin{equation}\label{eq:Wz}
Q^z(t) = \max\{ z + X(t),\, \sup_{0<s\leq t} [X(t)-X(s)]\} = X(t) + \max\{ z, -\inf_{0\leq s\leq t} X(s)\},
\end{equation}
for any initial state $z\geq 0$, where $X(t)$ is a L\'evy process with no negative jumps.
Let $\tau^z(w)$, $0\leq w<z$ denote the first passage time of level $w$ by the process starting in $z$, i.e.
\begin{equation*}
\tau^z(w) := \inf \left\{ t \geq 0\, |\, Q^z(t) \leq w \,\right\}.
\end{equation*}
Then it is easily seen that for all $z\geq 0$,
\begin{equation*}
Q^z(t) = \left\{
\begin{array}{ll}
z + X(t), & {\rm if }\ t <\tau^z(0), \\
\sup_{0<s\leq t} [X(t)-X(s)], & {\rm if }\ t \geq \tau^z(0).
\end{array}\right.
\end{equation*}
Consequently,
\begin{equation}\label{eq:Yxy}
Y^{x,y}(t) = \left\{
\begin{array}{ll}
x - y, & \text{if }t < \tau^y(0),\\
\inf_{0<s\leq t} \{ x+X(s)\}, & \text{if }\tau^y(0) \leq t < \tau^x(0),\\
0, & \text{if }t \geq \tau^x(0).
\end{array}\right.
\end{equation}
Using this representation we can identify
\begin{equation*}
\Omega^{x,y}_T = \frac{1}{T}\,\E\left[\int_0^{\tau^x(0)\wedge T} Y^{x,y}(t) \dd t\right],
\end{equation*}
where $\wedge$ denotes the minimum operator, due to the fact $Y^{x,y}(t) = 0$ for $t\geq \tau^x(0)$. 
Subsequently, we decompose $\Omega_T^{x,y}$ into two terms 
\begin{equation}
\Psi^{x,y}_T := \frac{1}{T} \int_0^\infty \E[Y^{x,y}(t)]\, \dd t  \qquad
\text{and}
\qquad
\Delta_T^{x,y} := \Omega_T^{x,y} - \Psi_T^{x,y}.
\label{eq:Deltaxy}
\end{equation}
Note that $\Psi_T^{x,y}$ is obtained by replacing $T$ by $\infty$ only in the integration bound.
It is customary in the literature, particularly in the area of stochastic simulation, to compare the truncated integral to its natural expansion of the integration range to a  semi-infinite interval, see e.g. \cite[Prop.~2.1]{Awad2007}.The truncated integral connects to the long-run average estimator of a certain performance metric, whereas the infinite integral reflects its exact expectation.
The decomposition in \eqref{eq:Deltaxy} is insightful, because $\Psi_T^{x,y}$ prescribes the leading order behavior of $\Omega_T^{x,y}$, while $\Delta_T^{x,y}$ captures the smaller order error term.
In this section, we only consider $\Psi_T^{x,y}$. Subsection \ref{sec:trunc} investigates the magnitude of $\Delta_T^{x,y}$.
The next preliminary result presents a useful property of $\Psi_T^{x,y}$.
\begin{lemma}\label{lemma:psixy}
Let $x>y$. If $\E[\tau^x(0)]<\iy$, then
\begin{equation}\label{eq:H(x,y)}
\Psi^{x,y}_T = \frac{1}{T}\,\E[\tau^{y}(0)](x-y) + \Psi^{x-y,0}_T.
\end{equation}
\end{lemma}
The proof can be found in Appendix \ref{app:psixy}.
This leaves us with two unknowns $\E[\tau^y(0)]$ and $\Psi_T^{x-y,0}$.
The next lemma gives an equivalent form for the latter.
\begin{lemma}\label{lemma:psiz0}
If $\E[\tau^z(0)] < \iy$, then for all $z\geq 0$
\begin{equation}\label{eq:H(x,0)}
\Psi^{z,0}_T = \int_0^z \E[\tau^w(0)]\, \dd w.
\end{equation}
\end{lemma}
The proof can be found in Appendix \ref{app:psiz0}.
Since the term $\E[\tau^z(0)]$, for several values of $z$, appears in many of the preliminary results, we devote our attention to this in the next subsection.\\
\\*
\noindent
\textbf{First passage time}.
When studying the first passage time of level $0\leq w < z$, $\tau^z(w)$, of the workload process starting in $z$, we first observe that $\{\tau^z(z-w)\}_{w=0}^z$ is a spectrally positive L\'evy process itself, also visible through Figure \ref{fig:samplePaths}.
More precisely, it is a subordinator, i.e. a L\'evy process whose paths are almost surely non-decreasing \cite{Kyprianou2006}. 
In order to calculate $\E[\tau^z(z-w)]$ we use theory presented in \cite[Section 46]{Sato1999}, although results presented there are valid for spectrally \emph{negative} L\'evy processes, as opposed to the absence of negative jumps in our case. 
Nonetheless, our setting is easily transformed into this framework by observing that $\hat{X} \equiv -X$, that is $\hat{X}(t) = -X(t)$ for all $t\geq 0$, is spectrally negative. 
Furthermore, let
\begin{equation}
\label{eq:transformedTau}
\hat{\tau}^0(w) := \inf\{ t \geq 0\,:\, \hat{X}(t) \geq w\} =  \inf\{ t \geq 0\,:\, z+X(t) \leq z-w\} = \tau^z(z-w).
\end{equation}
For completeness, we cite \cite[Thm.~46.3]{Sato1999}.
\begin{theorem}
Let $\hat{X}(t)$ be a spectrally negative L\'evy process with generating triplet $(-a,\s,\hat{\nu})$ and $\hat{\tau}^0(y)$ its corresponding hitting time process. Define $\Upsilon(\thh)$ for $\thh\geq 0$ as
\begin{equation}\label{eq:thmCharExp}
\Upsilon(\thh) = -a\thh + \tfrac{1}{2}\s^2\thh^2 + \int_{-\infty}^0 (\ee^{\thh x}-1-\thh x{\bf 1}_{[-1,0)}(x))\, \hat{\nu}(\dd x).
\end{equation}
Then $\Upsilon(\thh)$ is strictly increasing and continuous, $\Upsilon(0)=0$, and $\Upsilon(\thh)\to\infty$ as $\thh\to\infty$. For $w\geq 0$ and $0\leq u < \infty$ we have
\begin{equation}\label{eq:invCharExp}
\E[\exp(-u\hat{\tau}^0(w))] = \exp(-w\,\Upsilon^{-1}(u)),
\end{equation}
where $\thh=\Upsilon^{-1}(u)$ is the inverse function of $u=\Upsilon(\thh)$. 
\end{theorem}
\noindent This immediately induces an expression for $\E[\tau^w(0)]$ and henceforth $\Psi^{z,0}$.
\begin{corollary}\label{cor:Psixy}
Let $X(t)$ be a spectrally positive L\'evy process defined as in \eqref{eq:Xlmprocess} with $\mu > \la$. Let $\Psi^{z,0}_T$ as in \eqref{eq:H(x,0)}. Then 
\begin{equation*}
\Psi^{z,0}_T = \frac{z^2}{2T(\mu-\la)}.
\end{equation*}
Furthermore, if $x,y\geq 0$, then
\begin{equation}\label{eq:mainResult}
\Psi^{x,y}_T = \frac{x^2-y^2}{2T(\mu-\la)}.
\end{equation}
\end{corollary}
The proof of Corollary \ref{cor:Psixy} can be found in Appendix \ref{app:Psixy}.
\noindent\textbf{Randomization}.
As we stated before, we easily obtain the original $\Omega_T$ from $\Omega_T^{x,y}$ through substitution of $x$ and $y$ by $Q(0)$ and $Q(\iy)$, respectively, and taking the expectation. 
In the previous paragraph, we deduced an explicit expression for $\Psi_T^{x,y}$, the leading order term for $\Omega_T^{x,y}$. 
Therefore we equivalently get an approximation for $\Omega_T$, given by
\begin{equation*}
\Psi_T := \frac{1}{T} \int_0^\iy \left( \E[Q(t)]-\E[Q(\iy)] \right)\, \dd t,
\end{equation*}
through randomization of $x$ and $y$ in $\Psi_T^{x,y}$.
By combining the results in Corollary \ref{cor:Psixy}, Lemma \ref{lemma:workloadmoments} and Proposition \ref{prop:truncation_error}, which is given at the end of this section, we directly prove the result in Theorem \ref{thm:mainresult}.

\subsection{Truncation error}\label{sec:trunc}

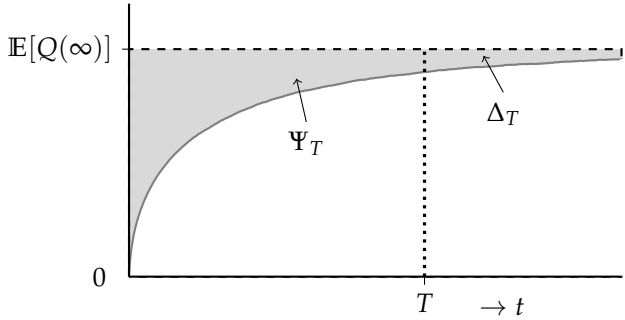
\begin{figure}
\centering
\begin{tikzpicture}[xscale=0.13,yscale=0.3]
	\node[below=0.4cm,right=0.5cm] at (x axis mid) {$\to t$};
	
	\draw[dashed, thick, fill =gray!30] (0,0) rectangle coordinate (eq) (50,10);
	\node[] at (-7,10) {$\E[Q(\infty)]$};
	\node[] at (-3,0) {$0$};

	\draw[->] (18,6.4) coordinate (a1) --   (21.65,8.49574);
	\node[below] at (a1) {$x=0$}; 

    	\foreach \x in {30}
     		\draw (\x,1pt) -- (\x,-10pt)
			node[anchor=north] {$T$};
    	\foreach \y in {10}
     		\draw (1pt,\y) -- (-20pt,\y);

	\draw[thick,color = gray,fill=white] plot
			file {Chapter_6/tikz/means0_2.txt};
		
	\draw[thick] (0,0) -- coordinate (x axis mid) (50,0);
    \draw[thick] (0,0) -- coordinate (y axis mid) (0,12);		
    \draw[color=white,very thick] (50,0.05) -- (50,9.56);
    
    \draw[very thick, dotted] (30,0) -- (30,10);
    
    \draw[->] (18,6.8) coordinate (delta) -- (17,8.7);
    \node[below] at (delta) {$\Psi_{T}$};
    
    \draw[->] (38,8.1) coordinate (delta) -- (36,9.7);
    \node[below] at (delta) {$\Delta_{T}$};
\end{tikzpicture}
\caption{Visualization of $\Omega_T$ and $\Psi_T$ as the area between the curves $\E[Q(t)]$, $\E[Q(\iy)]$ for $Q(0) = 0$.}
\label{fig:PsiVisualization}
\end{figure}

In order to get a better comprehension of the properties of $\Psi_T$, we depict the value in terms of the (infinite) region between the curves $\E[Q(t)]$, $\E[Q(\iy)]$ and the vertical axis for the case $Q(0)\equiv 0$ in Figure \ref{fig:PsiVisualization}.
In this figure, $\Omega_T$ is given by the area enclosed by the two curves, the vertical axis and the line $t=T$. 
One can see that the main contribution to the correction term $\Omega_T$ is given for small $t$. 
As $t$ increases, the difference between transient and stationary mean decreases.
Hence for moderate values of $T$, the contribution to the integral in \eqref{eq:Deltaxy} is only minor compared to the contribution over the interval $[0,T]$. 

Recall the definition of $\Delta^{x,y}_T$ as in \eqref{eq:Deltaxy}. As we alluded to in Subsection \ref{sec:difference} we claim the contribution of $\Delta^{x,y}_T$ to $\Omega_T^{x,y}$ is negligible compared to $\Psi^{x,y}_T$. Also note that
\begin{equation}
\label{eq:Delta}
\Delta_T := \Omega_T - \Psi_T = {-}\frac{1}{T} \int_T^\iy \big( \E[Q(t)] - \E[Q(\iy)] \big)\,\dd t.
\end{equation}
can be derived through $\Delta^{x,y}_T$ in a similar manner as we did for $\Psi^{x,y}_T$ to obtain $\Psi_T$.
To substantiate our claim, we compute an upper bound for $\Delta^{x,y}_T$ of order $1/T^2$. The existence of such an upper bound poses a limit on the error this tail integral contributed to the cost structure as a whole. 
\begin{proposition} \label{prop:truncation_error}
Let $x,y\geq 0$ and $\E[\max(Q(0),Q_\m(\iy))^3] < \iy$. Then
\begin{equation*}
|\Delta^{x,y}_T| \leq \frac{1}{T^2}\left(\frac{\max(y,x)^3}{3(\mu-\la)^2}+\frac{u_2 \max(y,x)^2}{2(\m-\la)^3}\right) 
\end{equation*}
and
\begin{equation*}
|\Delta_T| \leq \frac{1}{T^2}\left(\frac{\E[\max(Q(0),Q_\m(\iy))^3]}{3(\mu-\la)^2}+\frac{u_2 \E[\max(Q(0),Q_\m(\iy))^2]}{2(\m-\la)^3}\right).
\end{equation*}
\end{proposition}
The proof of Proposition \ref{prop:truncation_error} is given in Appendix \ref{sec:proof_truncation}.

\begin{remark}
In case the net-input process $X$ is light-tailed, that is there exists $u>0$ such that $\E[{\rm e}^{u X(1)}] < \iy$, it can be shown that the truncation error is of order ${\rm e}^{-\beta T}/T$ for some $\beta>0$. See Appendix \ref{sec:proof_truncation} for details.
\end{remark}
\section{Optimization}
\label{sec:optimization}
The result in Theorem \ref{thm:mainresult}, characterizing the leading order behavior of $\Omega_T(\mu)$, also reveals the behavior of $\Pi_T(\mu)$ in leading order. Namely, 
\begin{equation*}
\Pi_T(\mu) = \Pi_\iy(\mu) + \Psi_T(\mu) + O(1/T^2).
\end{equation*}
In fact, this representation naturally gives rise to an \emph{approximation} of the actual cost function:
\begin{align}\label{eq:decomposition}
\hat{\Pi}_{T}(\mu) := \Pi_{\iy}(\mu) + \Psi_T(\mu) 
\end{align}
Denote the corresponding minimizer of $\Pih$ by
\begin{equation}\label{eq:muhat}
\hat{\mu}_T^\star := \arg\min_{\m\geq 0} \Pih(\mu), \qquad \Pih^\star := \Pih(\hat{\mu}_T^\star)
\end{equation}
in addition to the definitions in \eqref{eq:muStar} and \eqref{eq:piStar}.
This section is devoted to the analysis of the minimizers $\muT$, $\muh$ and $\mui$, and the optimality gap for the two approximations. 

Throughout this section, we assume that $u_2, u_3 <\iy$ and $\E[Q(0)^2] <\iy$.

By its definition in \eqref{eq:PiInf} and Lemma \ref{lemma:workloadmoments}, we have an exact expression for the steady-state cost function:
\begin{equation*}
\Pi_{\iy}(\mu) = \frac{\la u_2}{2(\mu-\la)} + \aaa\mu.
\end{equation*}
It is easily verified that $\Pi_{\iy}$ is strictly convex in $\mu$, for instance by observing that $\Pi_{\iy}''(\mu) > 0$ for all $\mu > \la$. Therefore $\Pi_{\iy}$ has a unique global minimizer and 
\begin{equation}
\label{eq:muInf}
\mui = \la + \sqrt{\frac{\la u_2}{2\aaa}}, \qquad  \Pi_{\iy}^\star = \aaa\la + \sqrt{2\aaa\la u_2}.
\end{equation}
We are interested in the relation between $\mui$ and $\muT$, and between $\muh$ and $\muT$. 
Since $\Pi_{T}(\mu) = \Pi_{\iy}(\mu) + O(1/T)$ for all $\mu > \la$, we have pointwise convergence of the sequence $\Pi_{T}$, as well as $\hat{\Pi}_{T}$, to $\Pi_{\iy}$ for $T\to\iy$, we also expect $\muT \to \mui$ and $\muh\to\mui$ for $T\to\iy$.
Before proving that this convergence indeed holds, we present a result on the strict convexity of the function $\Pi_{T}$. 

\begin{lemma}\label{lemma:strict_convexity}
Let $\mu\geq 0$. The function $\Pi_{T}(\mu)$ is 
\begin{itemize}
\item convex in $\m$, if $Q(0)\equiv x$, $T<x/\mu$ and $\sigma=0$,
\item strictly convex in $\mu$, otherwise. 
\end{itemize}
\end{lemma}
Building upon strict convexity of both $\Pi_T(\mu)$ and $\Pi_\iy(\mu)$ for $\mu>\la$, we derive the following convergence result. 
\begin{proposition}\label{prop:min_convergence_mu}
Let $\muT$, $\muh$ and $\mui$ be as defined in \eqref{eq:muStar} and \eqref{eq:muhat}. Then 
\begin{equation*}
\muT \to \mui\, \qquad \text{\rm and } \qquad \muh \to \mui,
\end{equation*}
for $T\to\infty$.
\end{proposition}
The next result describes a refinement of $\muT$ in terms of $\mui$.
\begin{proposition}\label{prop:muBullet}
For $T$ sufficiently large,
\begin{equation*}
\muT = \mui + \frac{\mu_\bullet}{T} + o(1/T),
\end{equation*}
where
\begin{equation}\label{eq:muBullet}
\mu_\bullet = \frac{\E[Q(0)^2]}{\sqrt{8\la u_2\aaa}} - \frac{u_3}{3 u_2} - 3\sqrt{\frac{\aaa\la u_2}{8}}.
\end{equation}
\end{proposition}
\noindent
The proofs of the three results above can be found in Appendix \ref{app:proofs_optimization}.
Based on Proposition \ref{prop:muBullet} we propose a \emph{corrected staffing rule}, accounting for the finite horizon
\begin{equation}
\label{eq:correctedMu}
\tilde{\mu}_T^\star = \left[\mui + \frac{\mu_\bullet}{T}\right]^+, 
\end{equation}
with $\mu_\bullet$ as in \eqref{eq:muBullet}.
Here $[x]^+ := \max\{x,0\}$, which ensures the value of $\tilde{\mu}_T^\star$ is non-negative and thus is a feasible solution of the optimization problem.
This refined capacity allocation rule is expected to reduce the costs incurred in transient settings.
 However, the value of particular interest to us is the cost penalty for using either one of the approximations rather than the actual minimum $\muT$, that is, the \emph{optimality gap}.
 As it happens, we deduce the order of the optimality gap for $\mui$ with the help of the explicit form of $\mu_\bullet$ given in \eqref{eq:muBullet}, which is stated in the next proposition. The proof is given in Appendix \ref{sec:proofProp4}.

\begin{proposition}\label{prop:optimalitygap_mui}
Let $\mui$ be as in \eqref{eq:muInf}. Then,
\begin{equation*}
\Pi_\iy^\star- \Pi_T^\star = O(1/T^2).
\end{equation*}
\end{proposition}

\section{Numerical experiments}
\label{sec:numerics}

\subsection{Influence of $\Omega_{T}(\mu)$}
\label{sec:influence_omega}
We first assess the contribution of the correction to the cost function provided by Theorem 1. In other words, we investigate whether $\hat{\Pi}_{T}(\mu)$ as in \eqref{eq:PiT} yields a significantly better fit to $\Pi_{T}(\mu)$, than $\Pi_{\iy}(\mu)$ does. 
Note that these three functions only differ in the costs describing the congestion. 
Therefore, we limit our study in this subsection to the evaluation of $C_T(\mu)$ as in \eqref{eq:CTmu} with stationary equivalent $C_{\iy}(\mu) = \E[Q_{\m}(\iy)]$. 
Our novel approximation hence reads 
\begin{equation*}
\hat{C}_{T}(\mu) := C_{\infty}(\mu) + \Omega_{T}(\mu),
\end{equation*}
with $\Omega_{T}(\mu)$ given in Theorem \ref{thm:mainresult}.
We conduct our numerical experiments based on three models, namely: 
\begin{enumerate}
\item \underline{$M/M/1$ queue}: $U(t)$ is a unit rate compound Poisson process with exponentially distributed increments. We have $u_2 = 2$, $u_3=3$, so that
\begin{equation}\label{eq:MM1cor}
\hat{C}_{T}(\mu) = \frac{\la} {\mu-\la}  + \frac{1}{T(\mu-\la)} \left(\frac{x^2}{2} - \frac{\la^2}{(\mu-\la)^2} - \frac{\la} {\mu-\la} \right).
\end{equation}
\item \underline{$M/{\rm Pareto}/1$ queue}: $U(t)$ is a unit rate compound Poisson process with Pareto increments. The Pareto distribution deserves special attention due to its heavy-tailed nature, having tail probability $\bar{F}(x) = (x/k)^{-\g}$, if $x\geq k$ and 1 otherwise. 
It is well-known that heavy-tailed service times lead to long relaxation time. For our purposes, we fix shape parameter $\g = 16/5$ and scale parameter $k=11/16$, so that $\beta = 1$, $u_2 = 121/96$, $u_3 = 1331/256$ and $u_k=\iy$ for all $k>3$. Hence, 
\begin{equation}
\label{eq:MP1cor}
\hat{C}_{T}(\mu) = \frac{121\la} {192(\mu-\la)} + \frac{1}{2T(\mu-\la)}
\left( x^2 - \frac{(121\la/96)^2}{2(\mu-\la)^2} - \frac{ 1331\la/256 }{2(\mu-\la)}\right)
\end{equation}
\item \underline{Reflected Brownian motion}: $U(t)$ is Brownian motion with drift 1 and infinitesimal variance $\s^2$. We have $u_2 = \sigma^2$, $u_3=0$, so that
\begin{equation}\label{eq:RBMcor}
\hat{C}_{T}(\mu) = \frac{\la\sigma^2}{2(\mu-\la)} + \frac{1}{2T(\mu-\la)} \left( x^2 - \frac{\la^2\sigma^4}{2(\mu-\la)^2}\right).
\end{equation}
\end{enumerate} 
In light of the equivalence relations in \eqref{eq:Qidentity} we only consider the case $\la=1$. The cost values for general values of $\la$ follow by appropriate rescaling of $\mu$ and $T$.\\
\\*
\noindent
For the $M/M/1$ and $M/{\rm Pareto}/1$ queue, we obtained the function $C_{T}(\mu)$ through simulation and the results are accurate up until a 95\% confidence interval of width $10^{-3}$. For reflected Brownian motion, we used the explicit distribution function given in \cite{Harrison1985} for double numerical integration. The results for several values of $T$ and two different starting states are depicted in Figures 4-6. These plots also include the approximated functions $\hat{C}_{T}(\mu)$. 

\begin{figure}%
\centering
\begin{subfigure}{0.48\textwidth}
\centering
\begin{tikzpicture}[scale=0.75]
\begin{axis}[
	xmin = 1,
	xmax = 5,
	ymin = 0,
	ymax = 2,
	xlabel = {$\to \mu$},
	axis line style={->},
	axis lines = left,
	legend cell align=left,
	legend style = {at = {(axis cs: 5,2)},anchor = north east},
	legend columns=2,
	yscale = 0.8,
	xscale = 1
]

\addplot[col1, thick] table[x = mu,y=T2] {Chapter_6/tikz/mm1_0.txt};
\addplot[col1, dashed, thick] table[x = mu,y=ap2] {Chapter_6/tikz/mm1_0.txt};
\addplot[col4, thick] table[x = mu,y=T5] {Chapter_6/tikz/mm1_0.txt};
\addplot[col4, dashed, thick] table[x = mu,y=ap5] {Chapter_6/tikz/mm1_0.txt};
\addplot[col5, thick] table[x = mu,y=T10] {Chapter_6/tikz/mm1_0.txt};
\addplot[col5, dashed, thick] table[x = mu,y=ap10] {Chapter_6/tikz/mm1_0.txt};
\addplot[thick] table[x = mu,y=PSA] {Chapter_6/tikz/mm1_0.txt};

\legend{{$C_2(\mu)$},{$\hat C_2(\mu)$},{$C_5(\mu)$},{$\hat C_5(\mu)$},{$C_{10}(\mu)$},{$\hat C_{10}(\mu)$},{$C_\infty(\mu)$}}
\end{axis}

\end{tikzpicture}
\caption{$x=0$}
\end{subfigure}
\begin{subfigure}{0.48\textwidth}
\centering
\begin{tikzpicture}[scale=0.75]
\begin{axis}[
	xmin = 1,
	xmax = 5,
	ymin = 0,
	ymax = 3,
	xlabel = {$\to \mu$},
	axis line style={->},
	axis lines = left,
	legend cell align=left,
	legend style = {at = {(axis cs: 5,3)},anchor = north east},
	legend columns=2,
	yscale = 0.8,
	xscale = 1
]

\addplot[col1, thick] table[x = mu,y=T2] {Chapter_6/tikz/mm1_25.txt};
\addplot[col1, dashed, thick] table[x = mu,y=ap2] {Chapter_6/tikz/mm1_25.txt};
\addplot[col4, thick] table[x = mu,y=T5] {Chapter_6/tikz/mm1_25.txt};
\addplot[col4, dashed, thick] table[x = mu,y=ap5] {Chapter_6/tikz/mm1_25.txt};
\addplot[col5, thick] table[x = mu,y=T10] {Chapter_6/tikz/mm1_25.txt};
\addplot[col5, dashed, thick] table[x = mu,y=ap10] {Chapter_6/tikz/mm1_25.txt};
\addplot[thick] table[x = mu,y=PSA] {Chapter_6/tikz/mm1_25.txt};

\legend{{$C_2(\mu)$},{$\hat C_2(\mu)$},{$C_5(\mu)$},{$\hat C_5(\mu)$},{$C_{10}(\mu)$},{$\hat C_{10}(\mu)$},{$C_\infty(\mu)$}}
\end{axis}

\end{tikzpicture}
\caption{$x=2.5$}
\end{subfigure}
\caption{Comparison of exact waiting cost function $C_T(\mu)$ against corrected cost function $\hat{C}_T(\mu)$ and PSA cost function $C_\infty(\mu)$ for $T=2,5$ and 10 for the $M/M/1$ queue with $\la=1$.}
\label{fig:cont}%
\end{figure}
\begin{figure}[h!]
\centering
\begin{subfigure}{0.48\textwidth}
\centering
\begin{tikzpicture}[scale=0.75]
\begin{axis}[
	xmin = 1,
	xmax = 5,
	ymin = 0,
	ymax = 2,
	xlabel = {$\to \mu$},
	axis line style={->},
	axis lines = left,
	legend cell align=left,
	legend style = {at = {(axis cs: 5,2)},anchor = north east},
	legend columns=2,
	yscale = 0.8,
	xscale = 1
]

\addplot[col1, thick] table[x = mu,y=T2] {Chapter_6/tikz/mp1_0.txt};
\addplot[col1, dashed, thick] table[x = mu,y=ap2] {Chapter_6/tikz/mp1_0.txt};
\addplot[col4, thick] table[x = mu,y=T5] {Chapter_6/tikz/mp1_0.txt};
\addplot[col4, dashed, thick] table[x = mu,y=ap5] {Chapter_6/tikz/mp1_0.txt};
\addplot[col5, thick] table[x = mu,y=T10] {Chapter_6/tikz/mp1_0.txt};
\addplot[col5, dashed, thick] table[x = mu,y=ap10] {Chapter_6/tikz/mp1_0.txt};
\addplot[thick] table[x = mu,y=PSA] {Chapter_6/tikz/mp1_0.txt};

\legend{{$C_2(\mu)$},{$\hat C_2(\mu)$},{$C_5(\mu)$},{$\hat C_5(\mu)$},{$C_{10}(\mu)$},{$\hat C_{10}(\mu)$},{$C_\infty(\mu)$}}
\end{axis}
\end{tikzpicture}
\caption{$x=0$}
\end{subfigure}
\begin{subfigure}{0.48\textwidth}
\centering
\begin{tikzpicture}[scale=0.75]
\begin{axis}[
	xmin = 1,
	xmax = 5,
	ymin = 0,
	ymax = 3,
	xlabel = {$\to \mu$},
	axis line style={->},
	axis lines = left,
	legend cell align=left,
	legend style = {at = {(axis cs: 5,3)},anchor = north east},
	legend columns=2,
	yscale = 0.8,
	xscale = 1
]

\addplot[col1, thick] table[x = mu,y=T2] {Chapter_6/tikz/mp1_25.txt};
\addplot[col1, dashed, thick] table[x = mu,y=ap2] {Chapter_6/tikz/mp1_25.txt};
\addplot[col4, thick] table[x = mu,y=T5] {Chapter_6/tikz/mp1_25.txt};
\addplot[col4, dashed, thick] table[x = mu,y=ap5] {Chapter_6/tikz/mp1_25.txt};
\addplot[col5, thick] table[x = mu,y=T10] {Chapter_6/tikz/mp1_25.txt};
\addplot[col5, dashed, thick] table[x = mu,y=ap10] {Chapter_6/tikz/mp1_25.txt};
\addplot[thick] table[x = mu,y=PSA] {Chapter_6/tikz/mp1_25.txt};

\legend{{$C_2(\mu)$},{$\hat C_2(\mu)$},{$C_5(\mu)$},{$\hat C_5(\mu)$},{$C_{10}(\mu)$},{$\hat C_{10}(\mu)$},{$C_\infty(\mu)$}}
\end{axis}
\end{tikzpicture}
\caption{$x=2.5$}
\end{subfigure}
\caption{Comparison of exact waiting cost function $C_T(\mu)$ against corrected cost function $\hat{C}_T(\mu)$ and PSA cost function $C_\infty(\mu)$ for $T=2,5$ and 10 for the $M/$Pareto$/1$ queue with $\la=1$.}
\end{figure}

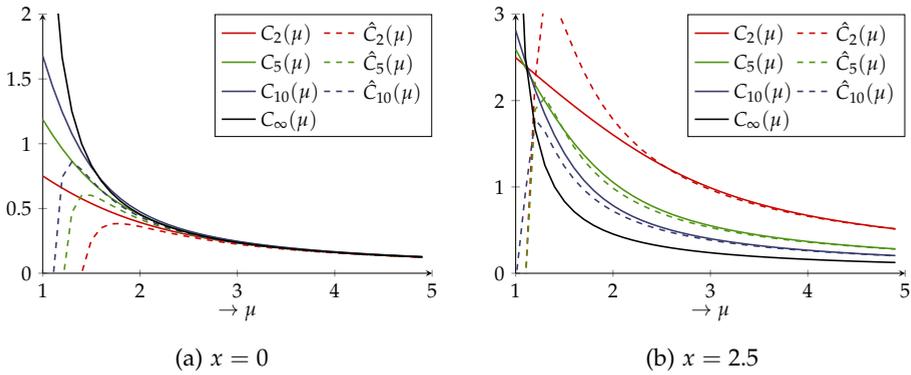
\begin{figure}[h!]
\centering
\begin{subfigure}{0.48\textwidth}
\centering
\begin{tikzpicture}[scale=0.75]
\begin{axis}[
	xmin = 1,
	xmax = 5,
	ymin = 0,
	ymax = 2,
	xlabel = {$\to \mu$},
	axis line style={->},
	axis lines = left,
	legend cell align=left,
	legend style = {at = {(axis cs: 5,2)},anchor = north east},
	legend columns=2,
	yscale = 0.8,
	xscale = 1
]

\addplot[col1, thick] table[x = mu,y=T2] {Chapter_6/tikz/rbm1_0.txt};
\addplot[col1, dashed, thick] table[x = mu,y=ap2] {Chapter_6/tikz/rbm1_0.txt};
\addplot[col4, thick] table[x = mu,y=T5] {Chapter_6/tikz/rbm1_0.txt};
\addplot[col4, dashed, thick] table[x = mu,y=ap5] {Chapter_6/tikz/rbm1_0.txt};
\addplot[col5, thick] table[x = mu,y=T10] {Chapter_6/tikz/rbm1_0.txt};
\addplot[col5, dashed, thick] table[x = mu,y=ap10] {Chapter_6/tikz/rbm1_0.txt};
\addplot[thick] table[x = mu,y=PSA] {Chapter_6/tikz/rbm1_0.txt};

\legend{{$C_2(\mu)$},{$\hat C_2(\mu)$},{$C_5(\mu)$},{$\hat C_5(\mu)$},{$C_{10}(\mu)$},{$\hat C_{10}(\mu)$},{$C_\infty(\mu)$}}
\end{axis}
\end{tikzpicture}
\caption{$x=0$}
\end{subfigure}
\begin{subfigure}{0.48\textwidth}
\centering
\begin{tikzpicture}[scale=0.75]
\begin{axis}[
	xmin = 1,
	xmax = 5,
	ymin = 0,
	ymax = 3,
	xlabel = {$\to \mu$},
	axis line style={->},
	axis lines = left,
	legend cell align=left,
	legend style = {at = {(axis cs: 5,3)},anchor = north east},
	legend columns=2,
	yscale = 0.8,
	xscale = 1
]

\addplot[col1, thick] table[x = mu,y=T2] {Chapter_6/tikz/rbm1_25.txt};
\addplot[col1, dashed, thick] table[x = mu,y=ap2] {Chapter_6/tikz/rbm1_25.txt};
\addplot[col4, thick] table[x = mu,y=T5] {Chapter_6/tikz/rbm1_25.txt};
\addplot[col4, dashed, thick] table[x = mu,y=ap5] {Chapter_6/tikz/rbm1_25.txt};
\addplot[col5, thick] table[x = mu,y=T10] {Chapter_6/tikz/rbm1_25.txt};
\addplot[col5, dashed, thick] table[x = mu,y=ap10] {Chapter_6/tikz/rbm1_25.txt};
\addplot[thick] table[x = mu,y=PSA] {Chapter_6/tikz/rbm1_25.txt};

\legend{{$C_2(\mu)$},{$\hat C_2(\mu)$},{$C_5(\mu)$},{$\hat C_5(\mu)$},{$C_{10}(\mu)$},{$\hat C_{10}(\mu)$},{$C_\infty(\mu)$}}
\end{axis}
\end{tikzpicture}
\caption{$x=2.5$}
\end{subfigure}
\caption{Comparison of exact waiting cost function $C_T(\mu)$ against corrected cost function $\hat{C}_T(\mu)$ and PSA cost function $C_\infty(\mu)$ for $T=2,5$ and 10 for reflected Brownian motion with $\sigma=1$.}
\end{figure}

We name a few observations based on these figures. 
First, we indeed note the pointwise convergence of $\hat{C}_{T}(\mu)$ to $\hat{C}_{\iy}(\mu)$ as $T$ grows, for all $\mu$ in all three cases. However, the difference between the stationary costs and those for small values of $T$ can be significant. This is most clear in the plots with $x=2.5$ and when $\mu$ is close to $\la$, i.e. it is in heavy-traffic. In these scenarios, it is evident that refinements to the stationary cost function are needed. $\hat{C}_{T}(\mu)$ does a fairly good job at providing such correction, especially for moderate values of $\mu$. 

Furthermore, we note that $C_{T}(\mu)$ approaches $C_{\iy}(\mu)$ from below for $x=0$ for any value of $\mu$, while this is not strictly the case for $x>0$. 
$\hat{C}_{T}(\mu)$ correctly captures the sign of this correction. 

Finally, observe that $\hat{C}_{T}(\mu)\to -\iy$ as $\mu$ approaches $\la$ from above. This divergence is clear from the expressions in \eqref{eq:MM1cor}-\eqref{eq:RBMcor}. 
Our correction term relies on the premise that under the coupling scheme, the sample paths of the two queues starting from different states have hit with high probability.
This is equivalent to stating that the `largest' of the two queues has emptied at least once before time $T$. However, as $\mu$ approaches $\la$, the system enters heavy traffic, and hence the hitting time of the zero barrier is set to run off to infinity. 
Consequently, this causes our approximation to be inaccurate for small values of $\mu$.

\subsection{Validation of corrected staffing rule}
\label{sec:num_opt}
In this section, we examine whether the corrected staffing rule $\tilde{\mu}_T^\star$ as in \eqref{eq:correctedMu} indeed yields a significant cost reduction over the choice of $\mui$ by comparing their true costs $\Pi_{T}(\tilde{\mu}_T^\star)$ and $\Pi_{T}(\mui)$.
We conduct this comparison for different values of the parameters, $\aaa$, $T$ and starting state $x$ through numerical experiments.  
The three models on which we do our calculations are the $M/M/1$ queue, the $M/$Pareto$/1$ queue and the reflected Brownian motion, as introduced in the previous subsection. 
We again focus on $\la=1$ only.

For each of the three models, we adhere to the following set-up. The quality of both staffing rules is assessed for $\aaa = 0.1, 1$ and 2, resembling three modes of valuation of the QoS in the system. 
As a benchmark, observe that the expected workload in steady-state conditions with staffing level $\mui$ equals
\begin{equation*}
C_\iy(\mui) = \sqrt{\frac{\aaa\la u_2}{2}}.
\end{equation*}
For each value of $\aaa$, we consider two scenarios: one in which the system starts empty, i.e. $x=0$, and one in which the initial state is double this benchmark value, thus $x=\sqrt{2\aaa\la u_2}$.
The numerics are presented for each model separately. We discuss general conclusions drawn from these results afterwards.\\
\\*
\noindent\textbf{$M/M/1$ queue.}
As we discussed before, if $U$ is a unit rate compound Poisson process with exponentially distributed increments, then $\Qlm$ describes the workload process in an $M/M/1$ queue. 
For this setting we get
\begin{equation*}
\mui = \la + \sqrt{\frac{\la} {\aaa}},\qquad \tilde{\mu}_T^\star = \left[\la + \sqrt{\frac{\la} {\aaa}} + \frac{1}{T}\left( \frac{x^2}{4\sqrt{\la\aaa}} - 1 - \frac{3}{2} \sqrt{\la\aaa}\right)\right]^+.
\end{equation*} 
Table \ref{tab:mm1} presents the actual costs corresponding to these two staffing levels for different value of $x$ and $\aaa$. 
\begin{table}[h!]
\centering
\resizebox{13cm}{!} {
\begin{tabular}{|c|r|r @{}r|r@{}r|r||r@{}r|r@{}r|r|}
\cline{3-12}
\multicolumn{1}{c}{} & \multicolumn{1}{r|}{} & \multicolumn{5}{c||}{$x = 0$} & \multicolumn{5}{c|}{$x = 2\sqrt{\aaa}$} \\
\hline
$\aaa$ & $T$     & $\mui$  & \hspace{1pt} $\Pi_T(\mui)$ & $\tilde{\mu}_T^\star$  & \hspace{1pt} $\Pi_T(\tilde{\mu}_T^\star)$ & r.c.i. & $\mui$  & \hspace{1pt}$\Pi_T(\mui)$ & $\tilde{\mu}_T^\star$  & \hspace{1pt}$\Pi_T(\tilde{\mu}_T^\star)$ & r.c.i. \\
\hline
\multicolumn{1}{|c|}{\multirow{4}[2]{*}{0.1}} & 1     & 4.162 & 0.620 & 2.688 & 0.536 & 0.136 & 4.162 & 0.682 & 2.688 & 0.536 & 0.214 \\
\multicolumn{1}{|c|}{} & 2     & 4.162 & 0.669 & 3.425 & 0.641 & 0.041 & 4.162 & 0.700 & 3.425 & 0.641 & 0.085 \\
\multicolumn{1}{|c|}{} & 5     & 4.162 & 0.706 & 3.867 & 0.703 & 0.005 & 4.162 & 0.719 & 3.867 & 0.703 & 0.022 \\
\multicolumn{1}{|c|}{} & 10    & 4.162 & 0.719 & 4.015 & 0.719 & 0.001 & 4.162 & 0.726 & 4.015 & 0.719 & 0.010 \\
\hline
\hline
\multicolumn{1}{|c|}{\multirow{4}[2]{*}{1}} & 1     & 2.000 & 2.309 & 0.000 & 0.500 & 0.783 & 2.000 & 3.500 & 0.500 & 2.750 & 0.214 \\
\multicolumn{1}{|c|}{} & 2     & 2.000 & 2.461 & 0.750 & 1.480 & 0.398 & 2.000 & 3.218 & 1.250 & 3.125 & 0.029 \\
\multicolumn{1}{|c|}{} & 5     & 2.000 & 2.675 & 1.500 & 2.400 & 0.103 & 2.000 & 3.043 & 1.700 & 2.968 & 0.025 \\
\multicolumn{1}{|c|}{} & 10    & 2.000 & 2.810 & 1.750 & 2.726 & 0.030 & 2.000 & 3.007 & 1.850 & 2.980 & 0.009 \\
\hline
\hline
\multicolumn{1}{|c|}{\multirow{4}[2]{*}{2}} &1     & 1.707 & 3.744 & 0.000 & 0.500 & 0.866 & 1.707 & 5.889 & 0.000 & 3.328 & 0.435 \\
\multicolumn{1}{|c|}{} &2     & 1.707 & 3.924 & 0.146 & 1.232 & 0.686 & 1.707 & 5.547 & 0.854 & 4.682 & 0.156 \\
\multicolumn{1}{|c|}{} &5     & 1.707 & 4.209 & 1.083 & 3.343 & 0.206 & 1.707 & 5.114 & 1.366 & 4.910 & 0.040 \\
\multicolumn{1}{|c|}{} &10    & 1.707 & 4.424 & 1.395 & 4.108 & 0.071 & 1.707 & 4.945 & 1.536 & 4.868 & 0.016 \\
\hline
\end{tabular}
}
\caption{Comparison of costs for the $M/M/1$ queue for steady-state and corrected staffing rules and relative cost improvement (r.c.i.).}
\label{tab:mm1}
\end{table}

\noindent
\textbf{$M$/Pareto/1 queue.}
In case the service requirements follow a Pareto distribution with shape parameter $\gamma = 16/5$, the staffing rule becomes
\begin{equation*}
 \mui = \la + \frac{11}{8}\sqrt{\frac{ \la }{3 \aaa}}, \  \tilde{\mu}_T^\star = \left[\la + \frac{11}{8}\sqrt{\frac{ \la }{3 \aaa}} + \frac{1}{T}\left( \frac{2 x^2}{11\sqrt{\la\aaa/3}} - \frac{11}{8} - \frac{11\sqrt{3\la\aaa}}{16}\right)\right]^+.
 \end{equation*} 
The numerical results are given in Table \ref{tab:mp1}. 
\begin{table}[h!]
\centering
\resizebox{13cm}{!} {
\begin{tabular}{|c|r|r @{}r|r@{}r|r||r@{}r|r@{}r|r|}
\cline{3-12}
\multicolumn{1}{c}{} & \multicolumn{1}{r|}{} & \multicolumn{5}{c||}{$x = 0$} & \multicolumn{5}{c|}{$x = 11/4\cdot \sqrt{\aaa/3}$} \\
\hline
$\aaa$ & $T$     & $\mui$  & \hspace{1pt} $\Pi_T(\mui)$ & $\tilde{\mu}_T^\star$  & \hspace{1pt} $\Pi_T(\tilde{\mu}_T^\star)$ & r.c.i. & $\mui$  & \hspace{1pt}$\Pi_T(\mui)$ & $\tilde{\mu}_T^\star$  & \hspace{1pt}$\Pi_T(\tilde{\mu}_T^\star)$ & r.c.i. \\
\hline
\multicolumn{1}{|c|}{\multirow{4}[2]{*}{0.1}} & 1     & 3.510 & 0.524 & 1.759 & 0.461 & 0.120 & 3.510 & 0.573 & 2.010 & 0.562 & 0.019 \\
\multicolumn{1}{|c|}{} & 2     & 3.510 & 0.555 & 2.635 & 0.539 & 0.029 & 3.510 & 0.580 & 2.760 & 0.574 & 0.010 \\
\multicolumn{1}{|c|}{} & 5     & 3.510 & 0.580 & 3.160 & 0.578 & 0.003 & 3.510 & 0.591 & 3.210 & 0.589 & 0.002 \\
\multicolumn{1}{|c|}{} & 10    & 3.510 & 0.590 & 3.335 & 0.590 & 0.000 & 3.510 & 0.596 & 3.360 & 0.595 & 0.001 \\
\hline
\hline
\multicolumn{1}{|c|}{\multirow{4}[2]{*}{1}} & 1     & 1.794 & 2.076 & 0.000 & 0.500 & 0.759 & 1.794 & 2.989 & 0.000 & 2.088 & 0.302 \\
\multicolumn{1}{|c|}{} & 2     & 1.794 & 2.190 & 0.511 & 1.291 & 0.411 & 1.794 & 2.790 & 0.610 &  2.588 & 0.072 \\
\multicolumn{1}{|c|}{} & 5     & 1.794 & 2.345 & 1.281 & 2.108 & 0.101 & 1.794 & 2.638 & 1.320 & 2.607 & 0.012 \\
\multicolumn{1}{|c|}{} & 10    & 1.794 & 2.441 & 1.537 & 2.371 & 0.029 & 1.794 & 2.597 & 1.557 & 2.585 & 0.005 \\
\hline
\hline
\multicolumn{1}{|c|}{\multirow{4}[2]{*}{2}} & 1     & 1.561 & 3.427 & 0.000 & 0.500 & 0.854 & 1.561 & 5.087 & 0.000 & 2.745 & 0.460 \\
\multicolumn{1}{|c|}{} & 2     & 1.561 & 3.567 & 0.032 & 1.050 & 0.706 & 1.561 & 4.832 & 0.172 & 3.417 & 0.293 \\
\multicolumn{1}{|c|}{} & 5     & 1.561 & 3.779 & 0.950 & 3.012 & 0.203 & 1.561 & 4.499 & 1.006 & 4.313 & 0.041 \\
\multicolumn{1}{|c|}{} & 10    & 1.561 & 3.935 & 1.255 & 3.356 & 0.147 & 1.561 & 4.351 & 1.284 & 4.304 & 0.011 \\
\hline
\end{tabular}
}
\caption{Comparison of costs for the $M/{\rm Pareto}/1$ queue for steady-state and corrected staffing rules and relative cost improvement (r.c.i.).}
\label{tab:mp1}
\end{table}
Just as in the results for the $M/M/1$ queue, we observe a higher reduction for larger value of $\aaa$ and $T$. Also, again $\tilde{\mu}_T < \mui$. Hence, the conclusions for the $M/{\rm Pareto}/1$ queue are similar to those of the $M/M/1$ queue. \\
\\*
\noindent\textbf{Reflected Brownian motion}.
In case the input process $U$ is Brownian motion with drift 1 and infinitesimal variance $\s^2$, the steady-state staffing rule and its corrected version reduce to 

\begin{equation*}
\mui = \la + \sqrt{\frac{\la\s^2}{2\aaa}}, \qquad  
\tilde{\mu}_T^\star = \left[\la + \sqrt{\frac{\la\s^2}{2\aaa}} + \frac{1}{2\sqrt{2}\,T}\left(\frac{x^2}{\sqrt{\la \aaa}\s} - 3\s\sqrt{\aaa\la}  \right)\right]^+.
\end{equation*}

In Tables \ref{tab:rbm1} and \ref{tab:rbm2}, the costs obtained through numerical evaluation are presented for several values of $x$, $T$. We also vary $\s$ to examine the influence of the volatility of arrival process on the quality of the staffing rules.

The observations on the influence of $\aaa, x$ and $T$ are similar to those of the $M/M/1$ queue and the $M/{\rm Pareto}/1$ queue. 
However, here we see little improvement by the corrected staffing rule for small values of $\aaa$ for both values of $x$. 
The results in Tables \ref{tab:rbm1}-\ref{tab:rbm2} also suggest that the reduction is smaller for larger values of $\s$. 

\begin{table}
\centering
\resizebox{13cm}{!} {
\begin{tabular}{|c|r|r @{}r|r@{}r|r||r@{}r|r@{}r|r|}
\cline{3-12}
\multicolumn{1}{c}{} & \multicolumn{1}{r|}{} & \multicolumn{5}{c||}{$x = 0$} & \multicolumn{5}{c|}{$x = \sqrt{2\aaa} $} \\
\hline
$\aaa$ & $T$     & $\mui$  & \hspace{1pt} $\Pi_T(\mui)$ & $\tilde{\mu}_T^\star$  & \hspace{1pt} $\Pi_T(\tilde{\mu}_T^\star)$ & r.c.i. & $\mui$  & \hspace{1pt}$\Pi_T(\mui)$ & $\tilde{\mu}_T^\star$  & \hspace{1pt}$\Pi_T(\tilde{\mu}_T^\star)$ & r.c.i. \\
\hline
\multicolumn{1}{|c|}{\multirow{4}[2]{*}{0.1}} & 1     & 3.236 & 0.525 & 2.901 & 0.518 & 0.013 & 3.236 & 0.565 & 3.124 & 0.564 & 0.001 \\
\multicolumn{1}{|c|}{} & 2     & 3.236 & 0.536 & 3.068 & 0.534 & 0.003 & 3.236 & 0.556 & 3.180 & 0.556 & 0.000 \\
\multicolumn{1}{|c|}{} & 5     & 3.236 & 0.543 & 3.169 & 0.542 & 0.000 & 3.236 & 0.551 & 3.214 & 0.551 & 0.000 \\
\multicolumn{1}{|c|}{} & 10    & 3.236 & 0.545 & 3.203 & 0.545 & 0.000 & 3.236 & 0.549 & 3.225 & 0.549 & 0.000 \\
\hline
\hline
\multicolumn{1}{|c|}{\multirow{4}[2]{*}{1}} & 1     & 1.500 & 3.420 & 0.000 & 0.833 & 0.756 & 1.500 & 4.741 & 1.000 & 3.984 & 0.160 \\
\multicolumn{1}{|c|}{} & 2     & 1.500 & 3.539 & 0.750 & 2.386 & 0.326 & 1.500 & 4.579 & 1.250 & 4.293 & 0.063 \\
\multicolumn{1}{|c|}{} & 5     & 1.500 & 3.707 & 1.200 & 3.363 & 0.093 & 1.500 & 4.335 & 1.400 & 4.274 & 0.014 \\
\multicolumn{1}{|c|}{} & 10    & 1.500 & 3.820 & 1.350 & 3.705 & 0.030 & 1.500 & 4.190 & 1.450 & 4.175 & 0.004 \\
\hline
\hline
\multicolumn{1}{|c|}{\multirow{4}[2]{*}{2}} & 1     & 1.500 & 3.420 & 0.000 & 0.833 & 0.756 & 1.500 & 4.741 & 1.000 & 3.984 & 0.160 \\
\multicolumn{1}{|c|}{} & 2     & 1.500 & 3.539 & 0.750 & 2.386 & 0.326 & 1.500 & 4.579 & 1.250 & 4.293 & 0.063 \\
\multicolumn{1}{|c|}{} & 5     & 1.500 & 3.707 & 1.200 & 3.363 & 0.093 & 1.500 & 4.335 & 1.400 & 4.274 & 0.014 \\
\multicolumn{1}{|c|}{} & 10    & 1.500 & 3.820 & 1.350 & 3.705 & 0.030 & 1.500 & 4.190 & 1.450 & 4.175 & 0.004 \\
\hline
\end{tabular}
}
\caption{Comparison of costs for RBM with $\sigma = 1$ for steady-state and corrected staffing rules and relative cost improvement (r.c.i.).}
\label{tab:rbm1}
\end{table}


\begin{table}
\centering
\resizebox{13cm}{!} {
\begin{tabular}{|c|r|r @{}r|r@{}r|r||r@{}r|r@{}r|r|}
\cline{3-12}
\multicolumn{1}{c}{} & \multicolumn{1}{r|}{} & \multicolumn{5}{c||}{$x = 0$} & \multicolumn{5}{c|}{$x = 2\sqrt{2\aaa} $} \\
\hline
$\aaa$ & $T$     & $\mui$  & \hspace{1pt} $\Pi_T(\mui)$ & $\tilde{\mu}_T^\star$  & \hspace{1pt} $\Pi_T(\tilde{\mu}_T^\star)$ & r.c.i. & $\mui$  & \hspace{1pt}$\Pi_T(\mui)$ & $\tilde{\mu}_T^\star$  & \hspace{1pt}$\Pi_T(\tilde{\mu}_T^\star)$ & r.c.i. \\
\hline
\multicolumn{1}{|c|}{\multirow{4}[2]{*}{0.1}} & 1     & 5.472 & 0.950 & 4.801 & 0.936 & 0.015 & 5.472 & 1.030 & 5.249 & 1.029 & 0.001 \\
\multicolumn{1}{|c|}{} & 2     & 5.472 & 0.972 & 5.137 & 0.968 & 0.003 & 5.472 & 1.012 & 5.360 & 1.012 & 0.000 \\
\multicolumn{1}{|c|}{} & 5     & 5.472 & 0.985 & 5.338 & 0.985 & 0.000 & 5.472 & 1.002 & 5.427 & 1.002 & 0.000 \\
\multicolumn{1}{|c|}{} & 10    & 5.472 & 0.990 & 5.405 & 0.990 & 0.000 & 5.472 & 0.998 & 5.450 & 0.998 & 0.000 \\
\hline
\hline
\multicolumn{1}{|c|}{\multirow{4}[2]{*}{1}} & 1     & 2.414 & 3.176 & 0.293 & 1.546 & 0.513 & 2.414 & 4.633 & 1.707 & 4.228 & 0.087 \\
\multicolumn{1}{|c|}{} & 2     & 2.414 & 3.356 & 1.354 & 2.690 & 0.199 & 2.414 & 4.375 & 2.061 & 4.247 & 0.029 \\
\multicolumn{1}{|c|}{} & 5     & 2.414 & 3.573 & 1.990 & 3.411 & 0.045 & 2.414 & 4.094 & 2.273 & 4.073 & 0.005 \\
\multicolumn{1}{|c|}{} & 10    & 2.414 & 3.689 & 2.202 & 3.646 & 0.012 & 2.414 & 3.966 & 2.344 & 3.962 & 0.001 \\
\hline
\hline
\multicolumn{1}{|c|}{\multirow{4}[2]{*}{2}} & 1     & 2.000 & 4.839 & 0.000 & 1.339 & 0.723 & 2.000 & 7.481 & 1.000 & 5.967 & 0.202 \\
\multicolumn{1}{|c|}{} & 2     & 2.000 & 5.078 & 0.500 & 2.773 & 0.454 & 2.000 & 7.158 & 1.500 & 6.585 & 0.080 \\
\multicolumn{1}{|c|}{} & 5     & 2.000 & 5.414 & 1.400 & 4.726 & 0.127 & 2.000 & 6.670 & 1.800 & 6.549 & 0.018 \\
\multicolumn{1}{|c|}{} & 10    & 2.000 & 5.639 & 1.700 & 5.409 & 0.041 & 2.000 & 6.380 & 1.900 & 6.349 & 0.005 \\
\hline
\end{tabular}
}
\caption{Comparison of costs for RBM with $\sigma = 2$ for steady-state and corrected staffing rules and relative cost improvement (r.c.i.).}
\label{tab:rbm2}
\end{table}

\subsection{Discussion}

Based upon these numerical results in Tables \ref{tab:mm1}-\ref{tab:rbm2}, we make a few remarks. The three models roughly exhibit similar behavior as $T$, $x$ and $\aaa$ are varied.

Non-surprisingly, we note that $\tilde{\mu}_T$ approaches $\mui$ with increasing $T$, which also implies that the cost reduction achieved by the corrected staffing rule vanishes as $T\to\iy$. 
Also, we observe that in all scenarios examined, the cost reduction increases with $\aaa$. This can be explained through investigation of the objective function $\Pi_T$ as function of $\m$. Namely, for $\aaa$ small, the curve is relatively flat around the true optimum $\muT$. Hence, in this case a moderate deviation from $\muT$ will likely not lead to a significant cost increase. However, as $\aaa$ becomes larger, i.e. server efficiency is valued more than minimization of congestion, the curve becomes more sharp around $\muT$, and hence more accurate approximations of $\muT$ are required to achieve an acceptable cost level. Hence, the corrected staffing rule \eqref{eq:correctedMu} proves particularly useful in these cases.

Another point we highlight is that the relative improvement is higher for $x=0$ than for $x=\sqrt{2\aaa\la u_2}$. Moreover, even though the initial state of the system is above the optimal equilibrium, $\tilde{\mu}_T$ is smaller than $\mui$. This is somewhat counter-intuitive. In fact, from \eqref{eq:muBullet} it follows that $\mu_\bullet$ positively contributes to the corrected staffing function if
\begin{equation*}
\E[Q^2(0)] >  3\aaa\la u_2 + \frac{2 u_2}{3 u_3}\,\sqrt{2\aaa\la u_2}.
\end{equation*}
%
\section{Conclusion \& further research}
\label{sec:conclusion_chapter6}
Motivated by the time-varying nature of queues in practical applications, we studied the impact that the transient phase has on traditional capacity allocation questions.
By defining a cost minimization problem, in which the objective function contains a correction accounting for the transient period, we identified the leading and second-order behavior of the cost function as a function of the interval length $T$. 
As a by-product, this result yields an approximation for the actual cost function, which is a refinement to its stationary counterpart.
Our numerical experiments in Section \ref{sec:influence_omega} demonstrate the improved accuracy achieved by this approximation in a number of settings. 
By perturbation analysis of the optimization problem, this furthermore gives rise to a correction to the steady-state optimal capacity allocation of order $1/T$. 
The necessity of the refined capacity allocation level is substantiated by the numerics in Section \ref{sec:num_opt}, which show the cost reduction that can be achieved in a number of settings, compared to settings in which stationary metrics are used.  
Especially for small values of $T$ and large values of $\alpha$ this reduction is significant.
Additionally, these results also indicate that it is relatively safe to use the stationary cost when $T$ is moderate, or $\alpha$ is small. 
The latter reflects the scenario in which QoS is much more valued than service efficiency. 
This observation links to the flat nature of the cost function around its optimal value for $\alpha$ small, a statement on the optimality gap that we formally proved in Proposition \ref{prop:optimalitygap_mui}.

Besides the validation of our theoretical results of Sections \ref{sec:analysis} and \ref{sec:optimization}, the numerical results also reveal some phenomena that require more investigation. 
As noted, our corrected capacity allocation level $\tilde{\mu}_T^\star$ is in most studied cases less than the steady-state optimal value $\mu_{\iy}^\star$. This implies that congestion levels tends to be higher under our staffing scheme then under stationary staffing. 
A possible explanation for this may be the fact that the planning period under consideration is finite. 
Clearly, in the setting we analyzed, anything that happens after time $T$ is neglected. 
Therefore, it might be beneficial from the cost perspective to end the period with a higher expected congestion level, as it does not need to be canceled out in the future. 
Related to this observation, it would be interesting to look at the setting in which staffing decisions need to be made in consecutive periods of equal length, in which the arrival rate changes at the start of each period. 
This case requires careful consideration of the correlation among the staffing decisions within the separate periods.

Another question that arises concerns the translation of our (qualitative) findings to more general queues, in particular the $M/G/s$ queue. 
Whereas in our analysis, the central decision variable is the server speed $\mu$, the variable of interest in multi-server queues is typically the number of servers. 
It may well be that similar explicit corrections to staffing levels can be deduced to account for transience.
Since our analysis heavily relies on the comparibility of the sample paths of two single-server queues, which is due to the equal negative drift for the two processes, another approach must be taken to tackle this extension.

The analysis and findings for the single-server queue with L\'evy input presented in this chapter may serve a stepping stone for investigation of these more elaborate problems.

\section*{Appendix}

\addcontentsline{toc}{section}{\hspace{7.1mm} Appendix}

\begin{subappendices}

\settocdepth{chapter}

\section{Proofs of Section \ref{sec:model_description}}
\label{app:proofs_model}

\subsection{Proof of Lemma \ref{lemma:workloadmoments}}
\label{app:proof_lemma_workload_moments}
\begin{proof}
The conditions of \cite[Cor.IX3.4]{Asmussen2003} are satisfied and therefore $Q_{\mu}(t)\Rightarrow \Qlm(\infty)$ in distribution for $t\rightarrow\infty$. Furthermore, its Laplace transform is for ${\rm Re}(s) < 0$
\[\tilde{Q}_{\mu}(s) = \E\Big[\ee^{s \Qlm(\infty)}\Big] = \frac{s \ka_{\mu}'(0)}{\ka_{\mu}(s)} = \frac{s(\la\ka_U'(0) - \mu)}{\la\ka_U(s) - \mu s} = \frac{s(\mu-\la)}{\mu s-\la \ka_U(s)}.\]
It can be checked that $\ka_U'(0) = \E[U(1)] = 1$, $\ka_U''(0) = u_2$ and $\ka_U'''(0) = u_3$, and $\klm'(0) = \la-\mu$, $\klm''(0) = \la u_2$ and $\klm'''(0) = \la u_3$.
Using l'H\^opital's rule we obtain the first moment of $\Qlm(\infty)$:
\begin{align*}
\E[\Qlm(\iy)] &= \lim_{s\to 0} \frac{\dd}{\dd s} \tilde{Q}_{\mu}(s) 
= \lim_{s\to 0} \klm'(0)\, \frac{\klm(s)-s\klm'(s)}{\klm(s)^2}\\
&= \klm'(0)\,\lim_{s\to 0} \frac{{-}s\klm''(s)}{2\klm(s)\klm'(s)} 
= \klm'(0)\,\lim_{s\to 0} \frac{ -s\klm'''(s)-\klm''(s)}{2\klm'(s)^2 + 2\klm(s)\klm''(s)} \\
&= {-}\frac{\klm''(0)}{2\klm'(0)} = \frac{\la u_2}{2(\mu-\la)}.
\end{align*}
Similarly, we derive the second moment:
\begin{align*}
\E[\Qlm^2(\iy)] &= \lim_{s\to 0} \frac{\dd^2}{\dd s^2} \tilde{Q}_{\mu}(s) 
= \lim_{s\to 0} \klm'(0)\, \frac{2 s \klm'(s)^2-2\klm'(s)\klm(s) - s\klm''(s)\klm(s)}{\klm(s)^3},
\end{align*}
We apply l'H\^opital's rule twice, to find
\begin{align*}
\E[\Qlm^2(\iy)] &= 
\klm'(0) \lim_{s\to 0} \frac{ 3s\klm''(s)\klm'(s) - 3 \klm''(s)\klm(s) - s \klm'''(s)\klm(s)}{ 3\klm'(s)\klm(s)^2 }\\
&= \klm'(0) \lim_{s\to 0} 
\frac{ 2 s \klm'''(s)\klm'(s) + 3 s\klm''(s)^2 - 4\klm'''(s)\klm(s) - s \klm^{(4)}(s)\klm(s)}
{6\klm'(s)^2\klm(s) + 3\klm''(s)\klm(s)^2}\\
&= \klm'(0) \lim_{s\to 0} 
\frac{ s\big[ 2\klm'''(s)\klm'(s) + 3 \klm''(s)^2 - \klm^{(4)}(s)\klm(s)\big] - 4\klm'''(s)\klm(s) }
{\klm(s)\big[6\klm'(s)^2 + 3\klm''(s)\klm(s)\big]}\\
&= \klm'(0) \lim_{s\to 0} \frac{s}{\klm(s)}\,
\frac{  2\klm'''(s)\klm'(s) + 3 \klm''(s)^2 - \klm^{(4)}(s)\klm(s)}{6\klm'(s)^2 + 3\klm''(s)\klm(s)}\\
&\qquad\qquad -\klm'(0) \lim_{s\to 0} \frac{4\klm'''(s)}{6\klm'(s)^2 + 3\klm''(s)\klm(s)}.
\end{align*}
Since $\klm(0) = 0$ and $\lim_{s\to 0} s/\klm(s) = 1/\klm'(0)$, we have
\begin{align}\label{eq:lemma_eq1}
& \klm'(0) \lim_{s\to 0} \frac{s}{\klm(s)}\,
\frac{  2\klm'''(s)\klm'(s) + 3 \klm''(s)^2 - \klm^{(4)}(s)\klm(s)}{6\klm'(s)^2 + 3\klm''(s)\klm(s)}
\nonumber\\
&\qquad \qquad \qquad = \frac{ 2\klm'''(0)\klm'(0) + 3 \klm''(0)^2}{6\klm'(0)^2} 
= \frac{\klm'''(0)}{3\klm'(0)} + \frac{\klm''(0)^2}{2\klm'(0)^2}
\end{align}
and 
\begin{equation}\label{eq:lemma_eq2}
\klm'(0) \lim_{s\to 0} \frac{4\klm'''(s)}{6\klm'(s)^2 + 3\klm''(s)\klm(s)} = \frac{ 2\klm'''(0)}{3\klm'(0)}.
\end{equation}
Combining \eqref{eq:lemma_eq1} and \eqref{eq:lemma_eq2} yields
\[
\E[\Qlm^2(\iy)] = \frac{\klm''(0)^2}{2\klm'(0)^2} - \frac{\klm'''(0)}{3\klm'(0)} 
= \frac{ \l^2u_2^2}{ 2(\mu-\l)^2} + \frac{\l u_3}{3(\mu-\l)}.
\]
\end{proof}

\subsection{Proof of Proposition \ref{prop:cost_convergence}}
\label{app:proof_prop1}
\begin{proof}
We prove the limit by showing that the difference
\[
\Pi_T(\mu) - \Pi_\iy(\mu) = \frac{1}{T} \int_0^T \Big(\E[Q_\mu(t)] - \E[Q_\mu(\iy)]\Big) \, \dd t
\]
converges to zero as $T\to\iy$ for $\mu>\la$ fixed. The assumption $\E[U(1)], \E[Q(0)] < \iy$ implies by \cite[Prop.~1]{Abate1994} that $\E[Q_\mu(t)]<\iy$ for all $t\geq 0$.
Following \cite{Abate1994}, we use the decomposition
\[
\E[Q_\mu(t)] = \E[Q^0_\mu(t)] + \left\{ \E[Q_\mu(t)] - \E[Q_\mu^0(t)]\right\},
\]
where $Q_\mu^0(t)$ represents the workload process if the system starts empty. 
From this decomposition it is revealed that $\E[Q^0_\mu(t)]$ and $\left\{ \E[Q_\mu(t)] - \E[Q_\mu^0(t)]\right\}$ are non-negative monotonically increasing and decreasing functions of $t$, respectively, see \cite[Prop.~2,Thm.~11]{Abate1994}. 
 Recall $\E[Q_\mu(t)] \to \E[Q_\mu(\iy)]$ for $t\to\iy$ by ergodicity of the workload process for any initial state $\E[Q(0)]< \iy$, if $\mu>\la$.
 Henceforth,
\begin{align*}
\E[Q_\mu(t)] &\leq \sup_t \E[Q_\mu^0(t)] + \sup_t  \left\{ \E[Q_\mu(t)] - \E[Q_\mu^0(t)]\right\} \\
&= \E[Q_\mu(\iy)] + \left\{\E[Q_\mu(0)] - \E[Q_\mu^0(0)]\right\} = \E[Q_\mu(\iy)] + \E[Q(0)],
\end{align*}
for all $t\geq 0$, which proves that the expected workload is bounded.
Fix $\varepsilon>0$. By convergence of $\E[Q_\mu(t)]$ for $t\to\iy$, there exists a value $t^* := t^*(\varepsilon)$ such that for all $t\geq t^*$
\begin{equation}
\left| \E[Q_\mu(t)] - \E[Q_\mu(\iy)] \right| < \varepsilon/2.
\end{equation}
Next, set
\[
T^* := T^*(\varepsilon) = \frac{2\,t^*(\varepsilon)}{\varepsilon}\, ( 2 \E[Q_\mu(\iy)] + \E[Q(0)]).
\]
Then for $T\geq \hat{T}:= \max\{ t^*,T^* \}$, we have
\begin{align*}
\left| \frac{1}{T} \int_0^T \big(\E[Q_\mu(t)] - \E[Q_\mu(\iy)]\big)  \dd t \right| 
&\leq \frac{1}{T} \int_0^{t^*} \left| \E[Q_\mu(t)] - \E[Q_\mu(\iy)] \right|  \dd t \\
& \qquad+ \frac{1}{T} \int_{t^*}^T \left| \E[Q_\mu(t)] - \E[Q_\mu(\iy)] \right| \, \dd t \\
&\leq \frac{1}{T} \int_0^{t^*} \big(\E[Q_\mu(t)] + \E[Q_\mu(\iy)]\big) \dd t + \frac{1}{T} \int_{t^*}^T \frac{\varepsilon}{2}\, \dd t \\
&< \frac{t^*}{T}\,( 2 \E[Q_\mu(\iy)] + \E[Q(0)]) + \frac{T-t^*}{T} \,\frac{\varepsilon}{2}\\
&< \frac{t^*}{T^*} \,( 2 \E[Q_\mu(\iy)] + \E[Q(0)]) + \frac{\varepsilon}{2} = \varepsilon.
\end{align*}
Hence, for any choice of $\varepsilon>0$ we can find a value $\hat{T}$ such that $\Pi_{\hat{T}}(\mu)$ approaches $\Pi_\iy(\mu)$ within distance $\varepsilon$, which proves the limit.
\end{proof}

\section{Proofs of Section \ref{sec:analysis}}
\label{app:proofs_analysis}
\subsection{Proof of Lemma \ref{lemma:psixy}}\label{app:psixy}
\begin{proof}
Using the representation in \eqref{eq:Yxy} we write
\begin{align*}
\Psi^{x,y}_T &= \frac{1}{T}\int_0^{\infty} \E[Y^{x,y}(t)]\dd t \\
&= \frac{1}{T}\,\E\left[\int_0^{\tau^y(0)}Y^{x,y}(t)\right] \dd t + \frac{1}{T}\,\E\left[\int_{\tau^y(0)}^{\tau^x(0)} Y^{x,y}(t) \dd t\right]
+ \frac{1}{T}\,\E\left[\int_{\tau^y(0)}^\infty Y^{x,y}(t) \, \dd t\right]\\
&= \frac{1}{T}\,\E\left[\int_0^{\tau^y(0)}(x-y) \dd t\right] + \frac{1}{T}\,\E\left[\int_{\tau^y(0)}^{\tau^x(0)} Y^{x,y}(t) \dd t\right] \\
&= \frac{1}{T}\,\E[\tau^y(0)](x-y) + \frac{1}{T}\,\E\left[\int_{\tau^y(0)}^{\tau^x(0)} Y^{x,y}(t) \dd t\right].
\end{align*}
By \eqref{eq:Yxy} and the Strong Markov property holding for L\'evy processes \cite{Asmussen2003}, observe that \\* $Y^{x-y,0}(t) \equalD Y^{x,y}(\tau^y(0)+t)$, whereby
\begin{equation*}
\frac{1}{T}\,\E\left[\int_{\tau^y(0)}^{\tau^x(0)} Y^{x,y}(t)\,\dd t\right] = \frac{1}{T}\,\E\left[\int_{0}^{\tau^{x-y}(0)} Y^{x-y,0}(t) \dd t\right] = \Psi^{x-y,0}_T,
\end{equation*}
which completes the proof.
\end{proof}

\subsection{Proof of Lemma \ref{lemma:psiz0}}\label{app:psiz0}
\begin{proof}
Note that $Y^{z,0}(t)$ and $\tau^z(w)$ are intimately related. Namely, due to the fact that $X$ has no negative jumps 
\begin{equation*}
\{ \tau^z(w) \leq t\} = \{Y^{z,0}(t) \leq w \}.
\end{equation*}
In fact, $Y^{z,0}(\tau^z(w)) = w$, which implies that $\tau^z$ is a right inverse for $Y^{z,0}(t)$. Therefore, the following equality holds
\begin{equation*}
\int_0^{\tau^z(0)} Y^{z,0}(t)\, \dd t = \int_0^z \tau^z(w)\, \dd w,
\end{equation*}
which implies with the help of Fubini's theorem
\begin{equation*}
\Psi^{z,0}_T = \frac{1}{T}\,\int_0^z \E[\tau^z(w)]\, \dd w = \frac{1}{T}\,\int_0^z \E[\tau^{z-w}(0)]\,\dd w =\frac{1}{T}\, \int_0^z \E[\tau^{w}(0)] \,\dd w.
\end{equation*}
\end{proof}

\subsection{Proof of Corollary \ref{cor:Psixy}}\label{app:Psixy}
\begin{proof}
From \eqref{eq:invCharExp},
\begin{equation}\label{eq:corEq1}
\E[\hat{\tau}^0(w)] = -\tfrac{\dd}{\dd u} \left. \E[\exp(-u\,\hat{\tau}^0(w))]\right|_{u=0} = w\left.\frac{\dd}{\dd u} \Upsilon^{-1}(u)\right|_{u=0}.
\end{equation}
Since $\Upsilon(\theta)$ is strictly increasing and $\Upsilon(0)=0$, we get $\Upsilon^{-1}(0)=0$ and
\begin{equation*}
\left.\tfrac{\dd}{\dd u}\Upsilon^{-1}(u)\right|_{u=0} = \frac{1}{\Upsilon'(\Upsilon^{-1}(0))} = \{ \Upsilon'(0) \}^{-1}.
\end{equation*}
Furthermore,
\begin{align*}
\Upsilon'(\thh) &= -a+ \s^2\thh + \int_{-\infty}^0 (x\, \ee^{\thh x} - x{\bf 1}_{[-1,0)}(x)) \hat{\nu}(\dd x) \\
&=  -a + \s^2\thh - \int_0^\infty (y\, \ee^{-\thh y} - y{\rm 1}_{(0,1]}(y)) \nu(\dd y).
\end{align*}
Thus, $\Upsilon'(0) = -\E[X(1)] = \mu-\la$ and $\E[\hat{\tau}^0(w)] = w/(\mu-\la)$. By \eqref{eq:H(x,0)} and \eqref{eq:transformedTau}, we deduce that
\begin{equation*}
\Psi^{z,0}_T = \frac{1}{T}\, \int_0^z \E[\tau^w(0)] \,\dd w = \frac{1}{T}\, \int_0^z \E[\hat{\tau}^0(w)] \dd w = \frac{z^2}{2T(\mu-\la)}.
\end{equation*}
For $x>y$, we use Lemma \ref{lemma:psixy} to conclude
\begin{equation*}
\Psi^{x,y}_T =  \frac{y(x-y)}{T(\mu-\la)} + \frac{(x-y)^2}{2 T(\mu-\la)} = \frac{x^2-y^2}{2T(\mu-\la)}.
\end{equation*}
The result for $x<y$ follows directly by the observation $\Psi^{y,x}_T = -\Psi_T^{x,y}$.
\end{proof}

\subsection{Proof of Proposition \ref{prop:truncation_error}}\label{sec:proof_truncation}

\begin{proof}
To derive the upper bound for $\Delta^{x,y}_T$, we apply the same coupling argument as described in Section \ref{sec:analysis}. Let us assume without loss of generality $x>y$. 
In this case, 
\begin{equation*}
|\Delta^{x,y}_T| = \frac{1}{T} \int_T^\iy \E[Q^x(t)-Q^y(t)] \dd t \leq \frac{1}{T}\int_T^\iy \E[Q^x(t)-Q^0(t)]\dd t.
\end{equation*}
By the decomposition in \eqref{eq:Yxy},
\begin{align}
 \int_T^\infty \E[Q^x(t) - Q^0(t)] \dd t
	&= \int_T^\infty \E[(x+\inf_{s\leq t} X(s))\mathbbm{1}_{\{\tau^x(0)>t\}}] \dd t \nonumber\\
	&= \int_T^\infty \int_0^x P( x-u + \inf_{s\leq t}X(s) > 0) \dd u \dd t \nonumber\\
	&= \int_T^\infty \int_0^x P( \tau^{x-u}(0) > t ) \dd u\dd t \\
	&\leq \int_T^\iy \int_0^x \frac{\E[\tau^{x-u}(0)^2]}{t^2}\dd u \dd t \nonumber\\
	&=  \int_0^x \int_T^\iy \frac{\E[\tau^{x-u}(0)^2]}{t^2}\dd t\dd u
	= \int_0^x \frac{\E[\tau^{w}(0)^2]}{T}\,\dd w. \nonumber
	\label{eq:tailprobIntegral}
\end{align}
We obtain $\E[\tau^w(0)^2]$ with the help of its Laplace transform in \eqref{eq:invCharExp}. Namely,
\begin{align*}
\E[\tau^w(0)^2] &= \left.\tfrac{\dd^2}{\dd u^2}\E[\exp(-u \tau^w(0))]\right|_{u=0} \\
&= w^2\,\left(\left.\tfrac{\dd}{\dd u}\Upsilon^{-1}(u)\right|_{u=0}\right)^2 - w\left. \tfrac{\dd^2}{\dd u^2}\Upsilon^{-1}(u)\right|_{u=0}.
\end{align*}
As in the previous subsection we have $\left.\tfrac{\dd}{\dd u}\Upsilon^{-1}(u)\right|_{u=0} = (\mu-\la)^{-1}$, and 
\begin{equation*}
\left.\tfrac{\dd^2}{\dd u^2}\Upsilon^{-1}(u)\right|_{u=0} = {-}\frac{\Upsilon''(\Upsilon^{-1}(0))}{\Upsilon'(\Upsilon^{-1}(0))^3} = {-}\frac{\Upsilon''(0)}{\Upsilon'(0)^3}.
\end{equation*}
Since $\Upsilon'(0) = \mu-\la$ and 
\begin{equation*}
\Upsilon''(0) = \s^2 + \int_0^\infty x^2\,\nu(\dd x) = u_2,
\end{equation*}
we conclude
\begin{equation*}
\E[\tau^w(0)^2] = \frac{w^2}{(\m-\la)^2} + \frac{u_2w}{(\m-\la)^3},
\end{equation*}
so that
\begin{equation}\label{eq:delta_upper}
|\Delta^{x,y}_T| \leq \frac{1}{T^2}\int_0^x \left(\frac{w^2}{(\m-\la)^2} + \frac{u_2w}{(\m-\la)^3} \right)\dd w = \frac{1}{T^2}\left(\frac{x^3}{3(\mu-\la)^2}+\frac{u_2 x^2}{2(\m-\la)^3}\right).
\end{equation}
For general $x,y\geq 0$, 
\begin{equation*}
|\Delta^{x,y}_T| \leq \frac{1}{T^2}\left(\frac{\max(y,x)^3}{3(\mu-\la)^2}+\frac{u_2 \max(y,x)^2}{2(\m-\la)^3}\right).
\end{equation*}
As a direct consequence,
\begin{equation*}
|\Delta_T| \leq \frac{1}{T^2}\left(\frac{\E[\max(Q(0),Q_\m(\iy))^3]}{3(\mu-\la)^2}+\frac{u_2 \E[\max(Q(0),Q_\m(\iy))^2]}{2(\m-\la)^3}\right).
\end{equation*}
\end{proof}
\begin{remark}
Observe that if $X$ is light-tailed, that is $\E[\exp\{ -\theta X(1) \}]$ $= \E[\exp\{\kappa(\theta)\}] < \iy$ for some $\theta<0$, then $\Upsilon(\theta)$ as in \eqref{eq:invCharExp} has an analytic continuation in the negative half-plane, and in this region $\Upsilon(\theta)<0$. Consequently, we can replace the upper bound on the tail probability of $\tau^{x-u}(0)$ by 
\begin{equation*}
\mathbb{P}\left( \tau^{x-u}(0) > t\right) = \mathbb{P}\left( \ee^{\beta \tau^{x-u}(0)} > \ee^{\beta t} \right) \leq \ee^{-\beta t} \, \ee^{ (x-u)\Upsilon^{-1}(-\beta)},
\end{equation*}
for some $\beta > 0$, so that
\[ \int_T^\infty \E[Q^x(t) - Q^0(t)]\, \dd t \leq \ee^{-\beta T}\, \frac{\ee^{x\Upsilon^{-1}(-\beta)}-1}{\beta\,\Upsilon^{-1}(-\beta)}. \]
Along similar lines we deduce
\[ |\Delta^{x,y}_T| \leq \frac{ \ee^{-\beta T}}{T}\, \frac{\ee^{x\Upsilon^{-1}(-\beta)} + \ee^{y\Upsilon^{-1}(-\beta)} -2}{\beta\,\Upsilon^{-1}(-\beta)}
\]
and 
\[ |\Delta_T| \leq \frac{\ee^{-\beta T}}{T}\, \frac{\E[\ee^{Q(0)\Upsilon^{-1}(-\beta)}] + \E[\ee^{Q_\mu(\iy)\Upsilon^{-1}(-\beta)}] -2}{\beta\,\Upsilon^{-1}(-\beta)},\]
assuming that $\E[{\rm e}^{-y Q(0)}] < \iy$ for all $y>0$. The condition $\E[{\rm e}^{Q_\mu(\iy)\Upsilon^{-1}(-\beta)}]<\iy$ follows from Lemma \ref{lemma:workloadmoments}. 
Hence, the error decays exponentially fast for light-tailed input processes.
\end{remark}

\section{Proofs of Section \ref{sec:optimization}}
\label{app:proofs_optimization}

\subsection{Proof of Lemma \ref{lemma:strict_convexity}}
\begin{proof}

Since the term $\aaa\mu$ is convex, the strictness should come from  the term $C_T(\mu)$. 
Furthermore, observe that if a function $f_\mu(t)$ is convex for all $t\geq 0$, and strictly convex for all $t\geq\e$ for some $\e\in[0,T)$, i.e. for any $\mu_1,\mu_2>0$ and $a\in (0,1)$
\begin{equation*}
a\, f_{\mu_1}(t) + (1-a) f_{\mu_2}(t) > f_{a\mu_1+(1-a)\mu_2}(t),
\end{equation*}
then,
\begin{equation*}
a\int_0^T \, f_{\mu_1}(t)\, \dd t + (1-a)\int_0^T f_{\mu_2}(t) \dd t =
\int_0^T \,\big( a f_{\mu_1}(t) + (1-a)f_{\mu_2}(t) \big) \dd t 
\end{equation*}
\begin{align*}
&=  \int_0^\e \left( a f_{\mu_1}(t) + (1-a)f_{\mu_2}(t)\right) \dd t  + \int_\e^T \, \left(a f_{\mu_1}(t) + (1-a)f_{\mu_2}(t)\right) \dd t \\
&> \int_0^\e  f_{a\mu_1+(1-a)\mu_2}(t) \dd t + \int_\e^T  f_{a\mu_1+(1-a)\mu_2}(t)  \dd t. \\
&= \int_0^T f_{a\mu_1+(1-a)\mu_2}(t)  \dd t.
\end{align*}
Hence, it suffices to prove the convexity of $\E[Q_\mu(t)]$ as a function of $\mu$ for all $t\geq 0$, and strict convexity for $t\geq \e$ for some $\e\in[0,T)$.
Let $\tau^x_{\mu}(0)$ denote the first passage time of level 0 in the process $Q_\mu$ with $Q(0)=x$. Then,
\begin{align}
Q_\mu(t) &= U(t)-\mu t + \max\left\{ x , -\inf_{s\leq t} [ U(s)-\mu s] \right\}\\
&=
\left\{
\begin{array}{ll}
x+U(t)-\mu t, & \text{if } t<\tau^x_{\mu}(0),\\
U(t)-\mu t  -\inf_{s\leq t} [ U(s)-\mu s], & \text{if } t\geq \tau^x_{\mu}(0) ,
\end{array}\right.\label{eq:Qrep}
\end{align}
where
\begin{equation*}
 \tau^x_{\mu}(0) := \inf\{ t \geq 0\,:\, x+U(t)-\mu t \leq 0\}
 \end{equation*} 
and $U(t)$ is a spectrally positive L\'evy process.
Fix $\mu_1, \mu_2>0$ and $a\in(0,1)$. Define $\mu_3 := a\mu_1+(1-a)\mu_2$, and
\begin{equation*}
D(t) := a Q_{\mu_1}(t) + (1-a) Q_{\mu_2}(t) - Q_{\mu_3}(t).
\end{equation*}
In order to prove strict convexity we have to show that $D(t) \geq 0$ for all $t\geq 0$, thereby implying $\E [D(t)] \geq 0$, i.e. convexity, for all $t\geq 0$, and $D(t)>0$ with positive probability for $t\in[\e,T]$, for some $\e \in[0,T)$. 
We distinguish two cases: $x>0$ and $x=0$. \\
\\*
\textbf{The case $x>0$.}
We start by noticing that if $Q_{\mu_1}$, $Q_{\mu_2}$ and $Q_{\mu_3}$ experience the same input process $U(t)$, then by absence of negative jumps in $U(t)$, it holds that
\begin{equation}\label{eq:stochDom}
\tau^x_{\mu_2}(0) < \tau^x_{\mu_3}(0) < \tau^x_{\mu_1}(0).
\end{equation}
We use shorthand notation 
\begin{equation*}
I_k(t) := \inf_{0\leq s\leq t}[U(s)-\mu_k s],
\end{equation*}
for $k=1,2,3$.
Using representation \eqref{eq:Qrep} of the workload process, we obtain 
\begin{equation*}
D(t) = \left\{
\begin{array}{ll}
0,  
& \text{if } t < \tau^x_{\mu_2}(0),\\
-(1-a)\left(x+I_2(t) \right), 
& \text{if } \tau^x_{\mu_2}(0) \leq t < \tau^x_{\mu_3}(0),\\
a x - (1-a)I_2(t) + I_3(t),
& \text{if } \tau^x_{\mu_3}(0) \leq t < \tau^x_{\mu_1}(0),\\
- a I_1(t) - (1-a) I_2(t)
+ I_3(t),
& \text{if } t \geq \tau^x_{\mu_1}(0).
\end{array}
\right.
\end{equation*}
This partition allows us to spot when strict convexity can occur.
Note that by definition $t \geq \tau^x_{\mu_2}(0)$, $I_2(t) = \inf_{0\leq s\leq t}[U(s)-\mu_2s]\leq -x$, so that $D(t)\geq 0$ if $\tau^x_{\mu_2}(0) \leq t < \tau^x_{\mu_3}(0)$. 
Moreover, by subadditivity of the infimum, 
\begin{align*}
I_3(t) &= \inf_{0\leq s\leq t}[U(s)-\mu_3s] = \inf_{0\leq s\leq t}[a(U(s)-\mu_1s)+(1-a)(U(s)-\mu_2s)] \\
&\geq a \inf_{0\leq s\leq t}[U(s)-\mu_1s] + (1-a) \inf_{0\leq s\leq t}[U(s)-\mu_2s] = a I_1(t) + (1-a) I_2(t),
\end{align*}
and hence $D(t)\geq 0$ for $t \geq \tau^x_{\mu_1}(0)$.
Using the same argument, we deduce
\begin{equation*}
ax - (1-a)I_2(t) + I_3(t) \geq a x - (1-a) I_2(t) + a I_1(t) + (1-a) I_2(t) = a(x + I_1(t)).
\end{equation*}
In particular for $t < \tau^x_{\mu_1}(0)$, this value is strictly positive. 
As a result, $D(t)\geq 0$ for all $t\geq 0$.
On top of that $D(t) > 0$ for $t\in[\tau^x_{\mu_3}(0),\tau^x_{\mu_1}(0))$.
Accordingly, the latter implies strict positivity of $\E D(t)$, and therefore strict convexity of $\E Q_\mu(t)$, if the event $\{\tau^x_{\mu_3}(0)\leq t< \tau^x_{\mu_1}(0)\}$ occurs with positive probability. 
That is, 
\begin{align}
P(D(t)>0) &\geq P\left( a(x+I_1(t))\mathbbm{1}_{\{\tau^x_{\mu_3}(0)\leq t < \tau^x_{\mu_1}(0)\}}  > 0 \right)\nonumber\\
&= P\left( x+ I_1(t) > 0 , \tau^x_{\mu_3}(0)\leq t < \tau^x_{\mu_1}(0)\right)\nonumber\\
&= P\left( x+ I_1(t) > 0 | \tau^x_{\mu_3}(0)\leq t < \tau^x_{\mu_1}(0)\right)P\left(\tau^x_{\mu_3}(0)\leq t < \tau^x_{\mu_1}(0)\right)\nonumber\\
&= P\left(\tau^x_{\mu_3}(0)\leq t < \tau^x_{\mu_1}(0)\right) = P(\tau^x_{\mu_3}(0)\leq t) - P(\tau^x_{\mu_1}(0) \leq t) > 0, \label{eq:strictConv}
\end{align}
by the stochastic dominance in \eqref{eq:stochDom}. To ensure the strict inequality in \eqref{eq:strictConv} we have to enforce the condition
\begin{equation}\label{eq:condition}
P(\tau^x_{\mu_1}(0)<T) > 0.
\end{equation}
\begin{remark}
An example illustrating the need for this condition is the case in which $U(t)$ is a compound Poisson process and $T < x/\mu_2 < x/\mu_1$. Then 
\[Q_{\mu_k}(t) = x + U(t) - \mu_k t,\]
for all $t\in[0,T]$, since $U(t)\geq 0$ and therefore $\tau^x_{\mu_1}(0) > T$. Consequently, for all $a\in(0,1)$, 
\[ a\,Q_{\mu_1} + (1-a)\,Q_{\mu_2}(t) = Q_{\mu_3}(t),\]
proving only convexity of $\E Q_{\mu}(t)$ and subsequently convexity of $\int_0^T \E[Q_\mu(t)]\,\dd t$. In case $\sigma>0$, the probability in \eqref{eq:condition} is necessarily positive.
\end{remark}

\noindent \textbf{The case $x=0$.}
By the fact that $\tau_{\mu}(0) = 0$ for all $\mu>0$, proving that $D(t)>0$ in the case $x=0$ reduces to showing that the probability of 
\begin{equation*}
D(t) = a I_1(t) + (1-a) I_2(t) - I_3(t)>0
\end{equation*}
 happening is positive for all $t>0$. Define 
\begin{equation*}
 t_0 := \inf\{ t > 0\, :\, U(t) > 0 \},
 \end{equation*} 
and
\begin{equation*}
\tilde{\tau}_\mu := \inf\{ t > t_0\,: U(t) - \mu t \leq 0\}.
\end{equation*}
We note that $t_0$ as defined above, also defines the epoch of the start of a new excursion of the reflection $Q_\mu$ for all $\mu>0$. Namely,
\[U(s) \leq 0  \quad \Rightarrow\quad  U(s) - \mu s \leq -\mu s  \qquad \text{for all }0\leq s< t_0\]
\[\Rightarrow \inf_{0\leq s < t_0} [U(s)-\mu s] \leq -\mu t_0 \quad 
\Rightarrow U(t_0) - \mu t_0 -  \inf_{0\leq s < t_0} [U(s)-\mu s] \geq U(t_0) > 0.\]
Then $Q_\mu(t_0-) = 0$ for all $\mu>0$. 
By virtue of the Strong Markov Property, note that $Q_\mu(t_0+t) \equalD Q_\mu(t)$. 
Hence we assume without loss of generality $t_0=0$. 
Again, we have a stochastic dominance relation similar to \eqref{eq:stochDom}:
\begin{equation*}
\tilde{\tau}_{\mu_2} < \tilde{\tau}_{\mu_3} < \tilde{\tau}_{\mu_1},
\end{equation*}
for all $\mu_1<\mu_3<\mu_2$. 
Then 

\begin{equation*}
D(t) \equalD \left\{
\begin{array}{ll}
0,  
& \text{if } t < \tilde{\tau}_{\mu_2},\\
-(1-a)I_2(t), 
& \text{if } \tilde{\tau}_{\mu_2} \leq t < \tilde{\tau}_{\mu_3},\\
(1-a)I_2(t) + I_3(t),
& \text{if } \tilde{\tau}_{\mu_3} \leq t < \tilde{\tau}_{\mu_1},\\
- a I_1(t) - (1-a) I_2(t)
+ I_3(t),
& \text{if } t \geq \tilde{\tau}_{\mu_1}.
\end{array}
\right.
\end{equation*}
Clearly, $D(t)\geq 0$ for all $t\geq 0$ and 
\begin{equation*}
-(1-a)I_2(t) + I_3(t) \geq a I_1(t) > 0,
\end{equation*}
for $\tilde{\tau}_{\mu_3} \leq t < \tilde{\tau}_{\mu_1}$. 
Hence, in a similar manner to \eqref{eq:strictConv},
\begin{align}
P(D(t)>0) &\geq P\left( aI_1(t)\mathbbm{1}_{\{\tilde{\tau}_{\mu_3} \leq t < \tilde{\tau}_{\mu_1}\}}  > 0 \right)\nonumber\\
&= P\left( I_1(t) > 0 , \tilde{\tau}_{\mu_3} \leq t < \tilde{\tau}_{\mu_1}\right)\nonumber\\
&= P\left( I_1(t) > 0 | \tilde{\tau}_{\mu_3} \leq t < \tilde{\tau}_{\mu_1}\right)P\left(\tilde{\tau}_{\mu_3} \leq t < \tilde{\tau}_{\mu_1}\right)\nonumber\\
&= P\left(\tilde{\tau}_{\mu_3} \leq t < \tilde{\tau}_{\mu_1}\right) = P(\tilde{\tau}_{\mu_3}\leq t) - P(\tilde{\tau}_{\mu_1} \leq t) > 0. \label{eq:strictConv2}
\end{align}
The last inequality is satisfied it $P(\tilde{\tau}_{\mu_1} < T) >0$, which is  equivalent to $P( U(T) - \mu T \leq 0 ) >0$, a condition that is clearly true for all our choices of $U$. 
In conclusion, for $x=0$, $\E[D(t)] >0$ and therefore $\E[Q_\mu(t)]$ is a strictly convex function of $\mu$. 
\end{proof}

\subsection{Proof of Proposition \ref{prop:min_convergence_mu}}
The proof of the proposition relies on the following auxiliary lemma, of which we include the proof for completeness.
\begin{lemma}\label{lemma:minimizerConvergence}
Consider the sequence of functions $f_n:\, [x_0,\infty) \to \mathbb{R}$ and let $f: [x_0,\infty) \to\mathbb{R}$ be the pointwise limit for some $x_0\in \mathbb{R}$. 
Assume $f$ and $f_n$ are strictly convex for all $n$. 
Furthermore, let $f(y) \to \infty$ for both $y\to x_0^+$ and $y\to \infty$.  
If $x_n$ and $x$ are the minimizers for $f_n$ and $f$, respectively, then $x_n\to x$ for $n\to\infty$.  
\end{lemma}
\begin{proof}
We start by showing that the sequence $x_n$ is bounded. Fix $u_l, u_r$ such that $x_0<u_l < x < u_r$. We claim that there exists a $N\in\mathbb{N}$ such that $x_n\in[u_l,u_r]$ for all $n \geq N$. First, we prove the upper bound on $x_n$. For any strictly convex function $h$ with minimizer $x_h$, the following statement holds true:
\begin{equation}\label{eq100}
x_h < u_r \quad \Leftrightarrow \quad h \text{ is strictly increasing at } u_r.
\end{equation}
The first implication follows from observing that $h(x_h) < h(y)$ for all $y> x^*$ and definition of convexity:
\[ 0<\frac{h(u_r)-h(x_h)}{u_r-x_h} \leq \frac{h(u_r+\de)-h(u_r)}{\de}, \]
for all $\de>0$. Hence $h(u_r)<h(u_r+\de)$, i.e. $h$ is increasing at $u_r$. The converse follows immediately by observing that $h(u_r) < h(u_r+\de)$ for all $\de>0$, so that $x_h < u_r$. 
Next, we show that $f_n$ must be increasing at $u_r$ for $n$ sufficiently large. By pointwise convergence of $f_n$ we have
\[ \lim_{n\to\infty} [f_n(u_r+\de) - f_n(u_r)] = f(u_r+\de) - f(u_r).\]
Let $w_r:= f(u_r+\de) - f(u_r)>0$. Then 
\[ \exists N_r \in \mathbb{N}:\, \forall n\geq N_r:\, |[f_n(u_r+\de) - f_n(u_r)] - [f(u_r+\de)-f(u_r)] | < w_r/2.\]
Hence for $n\geq N_r$,
\[f(u_r+\de)-f(u_r) - w_r/2 <  f_n(u_r+\de) - f_n(u_r) < f(u_r+\de)-f(u_r) + w_r/2\]
\[\Rightarrow 0 < w_r/2 < f_n(u_r+\de) - f_n(u_r).
\]
Hence by \eqref{eq100}, $x_n < u_r$ for sufficiently large $n$. Similarly, we argue
\begin{equation*}
x_h > u_l \quad \Leftrightarrow \quad h \text{ is strictly decreasing at } u_l,
\end{equation*}
for any strictly convex function $h$ with minimizer $x_h$. Note that $x_h > u_l$ implies $h(x_h) - h(u_l) < 0$ and for all $\de>0$ we get by strict convexity
\[\frac{h(u_l)-h(u_l-\de)}{\de} < \frac{h(x_h)-h(u_l)}{x_h-u_l} < 0,\]
by which $h(u_l-\de)>h(u_l)$, i.e. $h$ is decreasing in $u_l$. Moreover, if $h$ is decreasing at $u_l$, then it is decreasing for all $y < u_l$, by arguments similar to the above. Therefore, $h(u_l-\de)> h(u_l)$ for all $\de>0$ and it must hold that $x_h>u_l$. Define $f(u_l) - f(u_l-\de) :=w_l < 0$, then again by pointwise convergence, we have that
\[ \exists N_l \in \mathbb{N}:\, \forall n\geq N_l:\, |[f_n(u_l) - f_n(u_l-\de)] - [f(u_l)-f(u_l-\de)] | < w_l,\]
whereupon
\[ f_n(u_l) - f_n(u_l-\de) < f(u_l) - f(u_l-\de) + w_l = 2w_l < 0.\]
Hence, for sufficiently large $n$, we also have $x_n > u_l$. Fix $N = \max\{N_l,N_r\}$, then for $n\geq N$, $x_n\in( u_l,u_r)$. That is, the sequence $x_n$ is bounded. Therefore, by the theorem of Bolzano-Weierstrass, $x_n$ has to have a convergent subsequence. That is, there exists a sequence $n_k$ such that $n_k \to\infty$ and $x_{n_k}\to a$ as $k\to \infty$ for some $a \in [u_l,u_r]$. 
We prove that every subsequence must converge to $x$ by contradiction. Suppose there exists a subsequence $n_k$ such that $x_{n_k}\to a\neq x$. Since, $x_n\in [u_l,u_r]$ for $n\geq N$, we may restrict our attention to the sequence of functions $\hat{f}_n:[u_l,u_r] \to \mathbb{R^+}$, consisting of the original function $f_n$ restricted to the domain $[u_l,u_r]$. To be precise $x_n = \arg\min_y f_n(y) = \arg\min_y \hat{f}_n(y)$ for $n\geq N$. Because $\hat{f}_n$ and $\hat{f}$ are bounded, we furthermore have that $\hat{f}_n \to \hat{f}$ uniformly.

Fix $\e>0$. By uniform convergence, there exists an $K_0 \in\mathbb{N}$ such that
\[  | \hat{f}_{n_k}( y ) - \hat{f}( y)| < \e /2,\quad \forall k\geq K_0,\ y \in[u_l,u_r].\]
Also, because $\hat{f}$ is convex, it is continuous, so that there exists a $\de := \de(\e)$ so that
\[ |z-y| < \de \quad \Rightarrow \quad |\hat{f}(z) - \hat{f}(y)| < \e/2.\]
Let $K_1$ be such that $|x_{n_k}-a| < \de$ for all $k\geq K_1$. Then for $k \geq K= \max\{K_0,K_1\}$ this implies 
\begin{align*}
|f_{n_k}(x_{n_k}) - f(a)| &= |\hat{f}_{n_k}(x_{n_k}) - \hat{f}(a)| \\
&\leq  |\hat{f}_{n_k}(x_{n_k}) - \hat{f}(x_{n_k})| + | \hat{f}(x_{n_k}) - f(a)| < \e/2 + \e/2 = \e.
\end{align*} 
Hence we conclude $\lim_{k\to\infty} \hat{f}_{n_k}(x_{n_k}) = f(a)$. 
Therefore, 
\[ \limsup_{n\to \infty} f_n(x_n) \geq f(a) > f(x),\]
by minimality of $x$. However, $f_n(x_n) \leq f_n(x)$, which implies $\limsup_{n\to\infty} f_n(x_n) \leq \lim_{n\to\infty} f_n(x) = f(x)$, contradicting the strict inequality above. Hence we deduce $x=a$. Consequently, every subsequence of $x_n$ converges to $x$ and therefore $x_n\to x$ as $n\to \infty$.
Applying Lemma \ref{lemma:minimizerConvergence} to the functions $\Pi_T$ and $\Pi_\iy$ with $x_0=\la$, together with Lemma \ref{lemma:strict_convexity}, we obtain the result immediately.

\end{proof}

\subsection{Proof of Proposition \ref{prop:muBullet}}
\begin{proof}
Note that $\Pi_\infty$ is a smooth function. 
By the first optimality condition $\Pi_\infty'(\mui)$ $= 0$. 
We first prove that also $\Pi_T(\mu)$ is differentiable with respect to $\mu$ for all $\mu\geq 0$. 
Recall \eqref{eq:PiT}, which defines the cost function as a combination of the accumulated expected transient queue length, and linear staffing costs.
The latter term is clearly differentiable, hence it remains to be proved that 
\begin{equation*}
C_T(\mu) = \frac{1}{T}\int_0^\iy \E[Q_\mu(t)] \, \dd t,
\end{equation*}
admits a derivative for all $\mu\geq 0$ with $T$ fixed.
This holds if and only if $\E[Q_\mu(t)]$ is differentiable for all $t\geq 0$. 
Let $Q(0)= x\geq 0$. 
Following \eqref{eq:Qlm},
\begin{align*}
\E[Q_\mu(t)] &= \E[X_\mu(t)] + \E\Big[ \max\{ x, \sup_{s\in[0,t]}\{- X_\mu(s)\} \}\Big]\\
&= (\la-\mu)t+ \E\Big[ \max\{ x, \sup_{s\in[0,t]}\{- X_\mu(s)\} \}\Big],
\end{align*}
where the first term is differentiable. 
Furthermore, 
\begin{align*}
 \E[ \max\{ x, \sup_{s\in[0,t]} \{ - X_\mu(s) \} \} ]
  &= x + \int_x^\iy P(\sup_{s\in[0,t]} \{ - X_\mu(s) \} > u )\dd u \\
 &= x+\int_x^\infty P(\hat\tau^0(u) \leq t ) \dd u,
\end{align*}
with $\hat{\tau}^0(u)$ as defined in \eqref{eq:transformedTau}.

Since $-X_\mu$ is a process with no positive jumps, we may apply \cite[Cor.~VII.3]{Bertoin1996}, which states that the following equivalence between measures holds:
\begin{equation}
s\,P( \hat\tau^0(u) \in ds ) du = u\,P( -X_\mu(s) \in du ) ds,
\end{equation}
so that 
\begin{align}
 \int_{u=x}^\infty P(\hat{\tau}^0(u) \leq t )\, \dd u 
&= 
\int_{u=x}^\iy \int_{s=0}^t P( \hat\tau^0(u) \in ds ) \dd u \nonumber \\
&= \int_{u=x}^\iy \int_{s=0}^t\,s^{-1} u\,\,P( -X_\mu(s) \in \dd u ) \dd s \nonumber  \\
&= \int_{u=x}^\iy \int_{s=0}^t\,s^{-1}  u\,P( X_\mu(s) \in \dd u ) \dd s \nonumber  \\
&= \int_{s=0}^t s^{-1} \E[ \max\{x, X_\mu(s)\} ]  \dd s \nonumber \\
&= \int_{s=0}^t \int_{v=x/s}^\iy P( X_\mu(s)/s > v ) \dd v \dd s \nonumber \\
&= \int_{s=0}^t \int_{v=x/s}^\iy P( U(\la s)/s > v + \mu ) \dd v \dd s \nonumber \\
&= \int_{s=0}^t \int_{w=x/s+\mu}^\iy P( U(\la s)/s > w) \dd w \dd s,
\end{align}
where the interchange of integrals is justified by Fubini's theorem and this last form is differentiable with respect to $\mu$. 
Substituting $Q(0)$ for $x$ straightforwardly yields differentiability of the complete cost function $\Pi_T$ for all $T$. 

Consequently we invoke the first optimality condition for $\muT$ to find
\begin{align*}
0=\Pi_T'(\muT)  
&= \Pi_\infty'(\muT) + \Psi_T'(\muT) + O(1/T^2)\\
&= \Pi_\infty'(\mui) + \Psi_T'(\mui) + (\muT-\mui)\left[ \Pi_\infty''(\mui) + \Psi_T''(\mui) \right] \\
&\qquad + \frac{1}{2}(\mu_T-\mui)^2\left[\Pi_T'''(\xi)+\Psi_T'''(\xi) \right] + O(1/T^2)\\
&= \Psi_T'(\mui) + (\muT-\mui)\left[ \Pi_\infty''(\mui) + \Psi_T''(\mui) \right] \\
&\qquad + \frac{1}{2}(\mu_T-\mui)^2\left[\Pi'''(\xi)+\Psi_T'''(\xi)\right] + O(1/T^2),
\end{align*}
for some $\xi \in [\muT,\mui]$. Rearranging this gives
\begin{align*}
\muT-\mui &= \frac{-\Psi_T'(\mui)}{\Pi_\infty''(\mui)+\Psi_T''(\mui) + \frac{1}{2}(\muT-\mui)(\Pi_\infty'''(\muT)+\Psi_T'''(\xi))} + O(1/T)\\
&= {-}\frac{\Psi_T'(\mui)}{\Pi_\iy''(\mui)} \left[1 - \frac{\Psi_T''(\mui)}{\Pi_\infty''(\mui)} - \frac{\muT-\mu_\infty}{2}\frac{\Pi_\infty'''(\mui)+\Psi_T'''(\mui)}{\Pi_\infty''(\mui)}\right] + O(1/T)\\
&= {-}\frac{\Psi_T'(\mui)}{\Pi_\iy''(\mui)} [1 + o(1)]
\end{align*}
for $T\to\infty$, since both $\mu_T - \mu_\infty$ and $\Psi_T''(\mui)$ are $o(1)$. 
Let
\begin{equation*}
\mu_\bullet := \lim_{T\to\iy} \frac{T \Psi_T'(\mui)}{\Pi_\iy''(\mui)}.
\end{equation*}
By \eqref{eq:mainResult} we have
\begin{equation*}
T \Psi'_T(\mu) = {-} \frac{\E[Q(0)^2]}{2(\mu-\la)^2} + \frac{\la u_3}{3(\mu-\la)^3} + \frac{3\la^2u_2^2}{4(\mu-\la)^4}.
\end{equation*}
Together with 
\begin{equation*}
\Pi_\iy''(\mu) = \frac{\la u_2}{(\mu-\la)^3}
\end{equation*}
and \eqref{eq:muInf} we obtain the expression for $\mu_\bullet$ in \eqref{eq:muBullet}.
\end{proof}
\subsection{Proof of Proposition \ref{prop:optimalitygap_mui}}\label{sec:proofProp4}
\begin{proof}
We upper bound the optimality gap by using the decomposition in \eqref{eq:decomposition}.
\begin{align}
|\Pi_\iy^\star - \Pi_T^\star| &= \left|\hat{\Pi}_T(\mu_\infty) + \Delta_T(\mui) - \hat{\Pi}_T(\muT) - \Delta_T(\muT)\right|\nonumber\\
&\leq |\hat{\Pi}_T(\mui) - \hat{\Pi}_T(\muT)| + |\Delta_T(\mui)| + |\Delta_T(\muT)|\nonumber\\
&= |\hat{\Pi}_T(\mui) - \hat{\Pi}_T(\muT)| + O(1/T^2),
\end{align}
since $\Delta_T(\mu) = O(1/T^2)$ by Proposition \ref{prop:truncation_error}. Next, we find an upper bound for $|\hat{\Pi}_T(\gamma) - \hat{\Pi}_T(\beta)|$, with $\hat{\Pi}_T(\cdot)$ as in \eqref{eq:decomposition}, in terms of the difference between $\gamma$ and $\beta$. 
For simplicity, denote $\hat{\gamma} = \gamma - \la$ and $\hat{\beta} = \beta-\la$, implying $\hat{\gamma}-\hat{\beta}=\gamma-\beta$. Then, using \eqref{eq:mainResult}, we get
\begin{align*}
|\hat{\Pi}_T(\mui) - \hat{\Pi}_T(\muT)| &= 
\left| \aaa(\hat{\gamma}-\hat{\beta})
+\left(\frac{\la u_2}{2} + \frac{\E[Q(0)^2]}{2T}\right)\left(\frac{1}{\hx}-\frac{1}{\hy}\right) \right. \\
& \qquad \left. -\frac{\la^2 u_2^2}{4T}\left(\frac{1}{\hx^3}-\frac{1}{\hy^3}\right) 
-\frac{\la u_3}{6T}\left(\frac{1}{\hx^2} - \frac{1}{\hy^2}\right)
\right|.
\end{align*}
Furthermore, we have 
\begin{align*}
\frac 1 \hx - \frac 1 \hy &= -\frac{\hx-\hy}{\hy^2} + \frac{(\hx-\hy)^2}{\hy^3} + O\left((\gamma-\beta)^3\right),\\
\frac 1 {\hx^2} - \frac 1 {\hy^2} &= -\frac{2(\hx-\hy)}{\hy^3} + \frac{3(\hx-\hy)^2}{\hy^4} + O\left((\gamma-\beta)^3\right),\\
\frac 1 {\hx^3} - \frac  1 {\hy^3} &= -\frac{3(\hx-\hy)}{\hy^4} + \frac{6(\hx-\hy)^2}{\hy^5} + O\left((\gamma-\beta)^3\right).\\
\end{align*}
Substituting these yields 
\begin{align*}
|\hat{\Pi}_T(\gamma) - \hat{\Pi}_T(\beta)| &= \left|(\gamma-\beta)\left[ \aaa - \frac{\la u_2}{2 \hy^2} + \frac{1}{2T\hy^2}\left(\E[Q(0)^2] + \frac{3\la^2 u_2^2}{2\hy^2} + \frac{2\la u_3}{3 \hy}\right)\right]\right. \\
&\qquad \left. - (\gamma-\beta)^2\left[ \frac{\la u_2}{2 \hy^3} + \frac{1}{2T\hy^3}\left(\E[Q(0)^2] - \frac{3\la^2 u_2^2}{\hy^2} - \frac{\la u_3}{\hy}\right)\right]\right| \\
& \qquad \qquad + O\left((\gamma-\beta)^3\right).
\end{align*}
Given that $\muT = \mui + \mu_\bullet/T + o(1/T)$, we find
\begin{align*}
|\hat{\Pi}_T(\mui) - \hat{\Pi}_T(\muT)| &= \frac{|\mu_\bullet|}{T}\left(\aaa - \frac{\la u_2}{2(\mui-\la )^2}\right) + O(1/T^2)\\
&= \frac{|\mu_\bullet|}{T}\left(\aaa - \frac{\la u_2}{2(\sqrt{\la u_2/2\aaa})^2}\right) + O(1/T^2) = O(1/T^2),
\end{align*}
which concludes the proof. 
\end{proof}

\resettocdepth 

\end{subappendices}

\chapter{A blood bank model}

\begin{chapterstart}
We consider a stochastic model for a blood bank, in which amounts of blood are offered and demanded according to
independent compound Poisson processes.
Blood is perishable, that is, blood can only be kept in storage for a limited amount of time.
Furthermore, demand for blood is impatient, that is, a demand for blood may be canceled if it cannot be satisfied soon enough.
For a range of perishability functions and demand impatience functions,
we derive the steady-state distribution of the blood inventory level.
Moreover, we deduce fluid and diffusion limits for the inventory process as the arrival rates of of the compound Poisson processes grow indefinitely. 
These scaling limits in turn provide normal approximations for the performance of large-scale systems.
\end{chapterstart}

\begin{flushright}
Based on\\
\textbf{A blood bank model with perishable blood and demand impatience}\\
\textit{Shaul Bar-Lev, Onno Boxma, Britt Mathijsen \& David Perry}\\
Submitted to \textit{Stochastic Systems}
\end{flushright}

\newpage

\section{Introduction}

This chapter is devoted to the study of a stochastic blood bank model in which amounts of blood are offered and demanded according to stochastic processes, and in which blood is perishable (that is, blood can only be kept for a limited amount of time) and demand for blood is impatient (that is, a demand request for blood may be canceled if it cannot be satisfied soon enough).
Let us first provide some background, and subsequently sketch the blood bank model in some more detail.\\
\\*
\textbf{Practical background.}
One of the major issues in securing blood supply to patients worldwide is to
provide blood of the best achievable quality, in the needed quantities. 
In most countries, blood, which is collected as whole blood units from human
donors, is separated into different components which are subsequently stored under different storage conditions according to their biological
characteristics, functions and respective expiration dates. Blood units and components are ordered by local hospital blood banks (LBB) from the Central Blood Bank (CBB) according to their operational needs. The CBB
has to run its inventory and supply according to these requests and to the
need to keep sufficient stock for immediate release in emergency situations.
It also has to perform tests to determine the unit's blood type
and to detect the presence of various pathogens which are able to cause
transfusion-transmitted diseases, such as Hepatitis B, Hepatitis C, Human Immunodeficiency Virus (HIV) and Syphilis, see e.g.~Steiner et al.~\cite{Steiner2010}.

Blood consists of several components: red blood cells, plasma and
plate-lets. 
In addition, there are $8$ blood groups (types):
$O^{+},O^{-},A^{+},A^{-},B^{+}$ ,$B^{-},AB^{+}$, $AB^{-}$ ($-$ means Rh
negative) where the interrelationship between the transfusion issuing
policies among the $8$ types is quite intricate. 
It turns out that each of
the negative types can satisfy the corresponding $+$ type, but not vice
versa.
Blood components are perishable as red blood cells can be used for only $35$ to $42$
days and platelets for only $5$ days (plasma, however, can be frozen and
kept for one year). 
Accordingly, if red blood cells and particularly platelets are not
used for blood transfusion within their expiration dates, then they perish. 

In most developed countries demand requirements of about $50.000$
blood donations are needed per one million persons per year. About 95\% of
these donations are aggregated by CBBs and the remaining 5\% by LBBs.
Blood units stored at the CBB are usually ordered by LBBs for planned elective
surgeries. However, as it happens rather frequently, elective surgeries
turn out to become emergency ones due to various conditions of the
patient involved. In such cases, hospitals use their own local blood banks
to supply the demand, and they cancel the required demand
from the CBB; this is what we refer to as demand impatience.
A good review on supply chain management in blood products appears in Beli\"{e}n \& Forc\'{e} \cite{Belieen2012} and the
references cited therein. Other relevant studies are
Ghandforoush and Sen \cite{Ghandforoush2010} \&
Stanger et al.\ \cite{Stanger2012}.
\newpage

\noindent
\textbf{Inventory model.}
In this chapter we consider the analysis of blood perishability and demand impatience, concentrating on only one 
blood type. We do this by considering the stochastic inventory processes $\{X_b(t)\}_{t \geq 0}$, with $X_b(t)$ the amount of blood
kept in storage at time $t$, and $\{X_d(t)\}_{t \geq 0}$, with $X_d(t)$ the amount of demand for blood (the shortage) at time $t$.
If $X_b(t)>0$ then $X_d(t)=0$, and if $X_d(t)>0$ then $X_b(t)=0$.
We assume that amounts of blood arrive according to a Poisson process,
and that requests for blood arrive according to another, independent, Poisson process.
The delivered and requested amounts of blood are assumed to be random variables.
We represent the perishability of blood by letting the amount of blood, when positive, decrease in a state-dependent way:
if the amount is $v$, then the decrement rate is $\xi_b v + \alpha_b$.
The $\xi_b$ factor is motivated by the fact that a large amount of blood 
suggests that some of the blood has been present for quite a while -- and hence there is a relatively high perishability rate when much blood is in inventory.
The $\alpha_b$ factor provides additional modeling flexibility. 
One can in this way represent the blood perishability more accurately;
but the $\alpha_b$ term could also, e.g., represent a fluid demand rate of individuals or organizations,
which contact the CBB directly, and that is only satisfied when there is inventory.
Similarly, we represent the demand impatience by a decrement rate $\xi_d v + \alpha_d$.
The $\xi_d$ factor is motivated by the following fact. When there is a large shortage (demand) of blood, 
there are probably many patients waiting for blood, so many patients that might become impatient
(that is, they could recover, or die, or become in need of emergency surgery)
leading to a cancellation of the required demand from the CBB.
Again, the $\alpha_d$ factor provides additional modeling flexibility;
it not only allows us to represent demand impatience more accurately, but it could also, e.g.,
represent additional donations of individuals in times of blood shortage.

The inclusion of both the perishability factor $\xi_b v + \alpha_b $ and the demand impatience factor
$\xi_d v + \alpha_d$ makes the analysis of the ensuing model mathematically
quite challenging, but 
leads to a very general model that contains many well-known models as special cases.
Our two-sided stochastic process, with both upward and downward jumps,
and with the rather general slope factors $\xi v + \alpha$, could represent a quite large class of stochastic phenomena.
It should for example be noted that this model is a two-sided generalization of the well-known shot-noise model that describes certain physical phenomena, see \cite{Keilson1959}).
In some of our calculations we remove either the $\xi$ factors or the
$\alpha$ factors, and this results in easier calculations and more explicit results.

Our main results are: (i) Determination of the steady-state distributions
of the amounts of blood and of demand in inventory; in particular, we present a detailed analysis of the case in which
the delivered and requested amounts of blood are both exponentially distributed. (ii) 
Expressions for mean amounts of blood and demand in storage, and for the probability of not being able to satisfy demand.
(iii) We obtain the fluid and diffusion limits of the blood inventory process, providing in particular sufficient conditions for the limit process to be an Ornstein-Uhlenbeck process.
\\*
\\
\noindent
\textbf{Structure of the chapter.}
The chapter is organized as follows:
Section~\ref{modeldesc} presents a detailed model description.
A steady-state analysis of the densities of demand and of blood amount in storage is 
contained in Section~\ref{analysis}, including the special case of exponentially distributed delivered and requested blood amounts when $\alpha_b=\alpha_d=0$ (i.e., pure proportionality).
The fluid and diffusion scalings are discussed in Section~\ref{sectionscaling},
and in Section~\ref{numericals} we present numerical results for certain performance measures like mean net amount of blood and the probability that there is a shortage of blood.
These results indicate, among other things, that the probability that there is a shortage of blood can be accurately approximated via a normal approximation, based on the Ornstein-Uhlenbeck process appearing in the diffusion scaling.
Section~\ref{conclus} contains some conclusions and suggestions for further research.

\section{Model description}
\label{modeldesc}
We consider the following highly simplified model of a blood bank, restricting ourselves to only one type of blood.

Blood amounts arrive according to a Poisson process with rate $\lambda_b$. The amounts which successively arrive are independent, identically
distributed random variables $B_1,B_2,\dots$ with distribution $F_b(\cdot)$; $\bar{F}_b(x) = 1 - F_b(x)$.
Demands for blood arrive according to a Poisson process with rate $\lambda_d$. The successive demand amounts are independent, identically
distributed random variables $D_1,D_2,\dots$ with distribution $F_d(\cdot)$; $\bar{F}_d(x) = 1 - F_d(x)$.
We view these amounts as continuous quantities, measured in, for instance liters.
If there is enough blood for a demand, then that demand is immediately satisfied.
If there is some blood, but not enough to fully satisfy a demand, then that demand is partially satisfied, using all the available blood.
The remainder of the demand may be satisfied later.

Blood has a finite expiration date. We make the assumption
that if the total amount of blood present is $v>0$,
then blood is discarded -- because of its finite expiration date -- at a rate $\xi_b v + \alpha_b$, so linear in $v$.
Blood demands have a finite patience.
We make the assumption
that if the total amount of demand present is $v>0$,
then demand disappears -- because of its finite patience -- at a rate $\xi_d v + \alpha_d$, so linear in $v$.

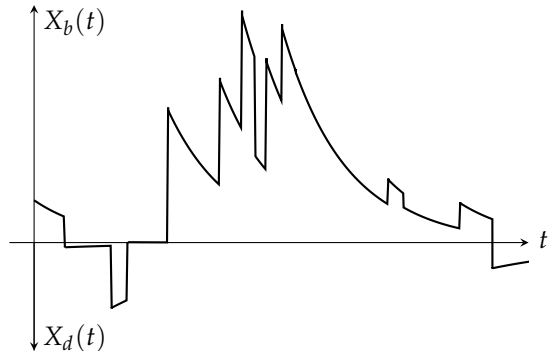
\begin{figure}
\centering
\begin{tikzpicture}
\begin{axis}[
	xmin = -0.2,
	xmax = 4,
	ymin = -5,
	ymax = 11,
	ticks = none,
	axis line style={->},
	axis lines = middle,
	yscale = 0.8,
	xscale = 1
	]

\addplot[black,thick] file {Chapter_7/tikz/sample_path2.txt};

\node at (axis cs: 0.33,-4.4) {$X_d(t)$};
\node at (axis cs: 0.33,10.2) {$X_b(t)$};
\draw[-stealth] (axis cs: 0,0) -- (axis cs:0,-5);

\end{axis}

\node at (7.05,1.47) {$t$};
\end{tikzpicture}
\caption{Sample path of net amount of blood available as a function of time.}
\label{fig:samplePath}
\end{figure}

Notice that {\em either} the total amount of blood present, {\em or} the total amount of demands,
is zero, {\em or} both are zero; they cannot be both positive. Hence we can easily in one figure depict
the two-sided process $\{X(t)\}_{t \geq 0}$ $=\{(X_b(t),X_d(t))\}_{t \geq 0}$
of total blood and total demand amounts present at any time $t$, as we have done in Figure \ref{fig:samplePath}.
For our purposes, we are mainly interested in the characteristics of the process described above in stationarity. 
Let us denote by $X_d$ the steady-state total amount of demand and by $X_b$ the steady-state total amount of blood present, with corresponding density functions $f(\cdot)$ and $g(\cdot)$, respectively. Notice that these are defective densities;
we have $\int_{0^+}^{\infty} f(v) \dd v = \pi_d = \P({\rm demand} ~ > ~ 0)$
and
$\int_{0^+}^{\infty} g(v) \dd v = \pi_b = \P({\rm blood} ~ > ~ 0)$.
If $\alpha_b=\alpha_d=0$, then neither $X_b$ nor $X_d$ has probability mass at zero, and $\pi_b+\pi_d=1$
(when there is only a very small amount $v$ present, the "decay" rate $\xi_b v$ or $\xi_d v$ is very small).
However, if $\alpha_b$ and/or $\alpha_d$ is positive, then there is a positive probability $\pi_0$ of being in $0$.

When $\xi_d$ and $\xi_b$ are positive, existence of these steady-state densities is obvious;
otherwise, the conditions for the existence of the steady-state distributions require some discussion, see Section~\ref{sectionvariant}.

\section{Steady-state analysis}
\label{analysis}
In this section we present a global approach towards determining $f(\cdot)$ and $g(\cdot)$ in the most general form of our model.
Using the Level Crossing Technique (LCT),
we derive two integral equations in $f(\cdot)$ and $g(\cdot)$.
Before attempting to solve these equations,
we consider a few important performance measures which can be expressed in $f(\cdot)$
and $g(\cdot)$, $\pi_0$
and the mean length of time during which, uninterruptedly, there is a positive amount of blood (respectively demand).
The latter could be viewed as the busy period of the $X_b$ process (respectively of the $X_d$ process).

First we consider the density $g(\cdot)$ of the amount of blood.
We equate the rate at which some positive blood level $v$ is upcrossed and downcrossed, respectively.
LCT leads to the following integral equation: for $v>0$,
\begin{align}
&\lambda_b  \int_0^v g(y) \bar{F}_b(v-y) {\rm d}y 
+
\lambda_b  \int_0^{\infty} f(y) \bar{F}_b(v+y) {\rm d}y 
+ \pi_0 \lambda_b \bar{F}_b(v)
\nonumber
\\
&\qquad =
\lambda_d  \int_v^{\infty} g(y) \bar{F}_d(y-v) {\rm d}y 
+
(\xi_b v + \alpha_b) g(v).
\label{eq:blood}
\end{align}
Here the three terms in the left-hand side represent the rate of crossing level $v$ from below;
the first term corresponds to a jump from a blood inventory level between $0$ and $v$,
whereas the second term corresponds to a jump from a shortage level, and the third term
corresponds to a jump from level $0$.
The two terms in the right-hand side represent the rate of crossing level $v$ from above;
the first term corresponds to a jump from above $v$, and the second term to a smooth crossing.

Next, we consider the density $f(\cdot)$ of the amount of demand (shortage).
We equate the rate at which some positive demand level $v$ is upcrossed and downcrossed, respectively.
LCT leads to the following integral equation: for $v>0$,
\begin{align}
& \lambda_d  \int_0^v f(y) \bar{F}_d(v-y) {\rm d}y 
+
\lambda_d  \int_0^{\infty} g(y) \bar{F}_d(v+y) {\rm d}y 
+ \pi_0 \lambda_d \bar{F}_d(v)
\nonumber
\\
&\qquad =
\lambda_b  \int_v^{\infty} f(y) \bar{F}_b(y-v) {\rm d}y 
+
(\xi_d v + \alpha_d) f(v).
\label{eq:demand}
\end{align}
It should be noted that these two, coupled, equations are symmetric (swap $f$ and $g$, and the $b$ and $d$ parameters).

In general, it appears to be very difficult to solve these integral equations.
In Section~\ref{sectionexp} we assume that
both $F_b(\cdot)$ and $F_d(\cdot)$ are exponential.
In that case, we are able to obtain explicit expressions of $f(\cdot)$ and $g(\cdot)$, in terms of hypergeometric functions.
In Section~\ref{gener} we consider the case that $F_b(\cdot)$ and $F_d(\cdot)$ are Coxian  distributions,
a class of distributions that lies dense in the class of all distributions of non-negative random variables,
and that is suitable for handling the above coupled integral equations via Laplace transforms (LT).
We are able to transform (\ref{eq:blood}) and (\ref{eq:demand})
into inhomogeneous first-order differential equations in the LTs of $f(\cdot)$ and $g(\cdot)$,
and thus to obtain those LTs.\\
\\*
\noindent
\textbf{A few simple performance measures.}
Without solving \eqref{eq:blood}-\eqref{eq:demand} explicitly, we are able to deduce some characteristics of the steady-state inventory level.

First, we can relate $\pi_0$ to the densities $f(\cdot)$ and $g(\cdot)$; see Proposition~\ref{prop. emptiness} below.
Subsequently we express the mean length of time during which there is, uninterruptedly, a positive amount of blood present (we call this the non-emptiness period of the inventory system), into $f(\cdot)$, $g(\cdot)$ and $\pi_0$.
We do the same for the mean length of time during which there is, uninterruptedly, a positive demand, i.e., the non-emptiness period of the demand process, see Proposition~\ref{prop: empty}.

\begin{proposition}
\label{prop. emptiness}Let $\pi _{0}\ $be\ the steady-state atom probability
of the zero period. Then%
\[
\pi _{0}=\frac{\alpha _{d}f(0)+\alpha _{b}g(0)}{\lambda
_{d}+\lambda _{b}}.
\]
\end{proposition}

\begin{proof}
Substitute $v=0$\ in \eqref{eq:blood} and \eqref{eq:demand} and take the sum. The result is
obtained after several steps of elementary algebra.   
\end{proof}

The result introduced in the proposition above is very intuitive. By LCT, $%
\alpha _{d}f(0)+\alpha _{b}g(0)$\ is the rate at which
level $0$ is reached (i.e., the process will now really stay at $0$ for a while),
so that $[\alpha _{d}f(0)+\alpha
_{b}g(0)]^{-1}$\ is the expected length of time between two successive times level $0$ is reached by the fluid. 
More precisely, the \textit{zero periods} and \textit{non-zero periods}
generate an alternating renewal process whose expected cycle length is $%
[\alpha _{d}f(0)+\alpha _{b}g(0)]^{-1}$. The expected length of the zero
period is $[\lambda _{d}+\lambda _{b}]^{-1}$, since the end of the zero
period is terminated at the moment of the next jump. But the jump process is
a Poisson process with rate $\lambda _{d}+\lambda _{b}$. Now the renewal
reward theorem simply says that 
\[
\pi _{0}=\frac{\E[\text{zero period}]}{\E[\text{cycle}]}.
\]%

In preparation of the next proposition, 
for the  process $\{X(t)\}_{t\geq 0}$ we
define a modified process $\{X_{m}(t)\}_{t\geq 0}$, where $X_{m}$ is constructed by deleting the zero-periods (only the zero periods, not the emptiness periods) from $X$ and gluing together the \textit{non-zero periods}. The modified
process is $X_m$ such that $X_{m}(t)=X_{d}(t)\mathbbm{1}_{\{X_{d}(t)>0\}}+X_{b}(t)\mathbbm{1}_{\{X_{b}(t)>0\}}$ where by definition
of the model $\{X_{d}(t)>0\}\Rightarrow \{X_{b}(t)=0\}$ and $%
\{X_{b}(t)>0\}\Rightarrow \{X_{d}(t)=0\}$.

\begin{proposition}
\label{prop: empty}Let $B_{b}$\ and $I_{b}$ be the generic non-emptiness
period and the emptiness period, respectively, of the inventory system.
Similarly, let $B_{d}$\ and $I_{d}$ be the generic non-emptiness period and the
emptiness period, respectively, of the demand process. Then%
\[
{\rm (i)}\left\{ 
\begin{array}{l}
\ \E [B_{b}]=\frac{1-\pi _{0}}{\alpha _{b}g(0)+\lambda _{d}\int_{0}^{\infty }%
\bar{F}_{d}(y)g(y)\dd y},  \\ 
\ \E [B_{d}]=\frac{1-\pi _{0}}{\alpha _{d}f(0)+\lambda _{b}\int_{0}^{\infty }%
\bar{F}_{b}(y)f(y)\dd y} 
\end{array}%
\right. 
\]%
and
\\
\[
{\rm (ii)}\left\{ 
\begin{array}{l}
\ \E[I_{b}]=\frac{1}{\lambda _{b}\int_{0}^{\infty }\bar{F}_{b}(y)f(y)dy+\lambda
_{b}\pi _{0}}
- \E[B_b] ,\\ 
\ \E[I_{d}]=\frac{1}{\lambda _{d}\int_{0}^{\infty }\bar{F}_{d}(y)g(y)dy+\lambda
_{d}\pi _{0}}
-\E[B_d] .
\end{array}%
\right. 
\]
\end{proposition}

\begin{proof}

(i) Consider the non-emptiness period of the inventory system. 
The steady-state densities of the inventory system and the demand process of $X_m$\ are given by%
\[
g_{m}(x)=\frac{g(x)}{1-\pi _{0}},\ \ \ \ \ f_{m}(x)=\frac{f(x)}{1-\pi _{0}},
\]%
respectively. At the end of the non-emptiness period of the inventory system there are two disjoint
ways (disjoint events) to downcross level $0+$. Either level $0$ is
downcrossed by a negative jump or level $0+$ is reached by the fluid
reduction (both in $X_m$). The rate of the first
event\ is $\lambda _{d}\int_{0}^{\infty }\bar{F}_{d}(y)g_{m}(y)\dd y$\ and by
LCT the rate of the second event is $\alpha _{b}g_{m}(0)$. Since the
events are disjoint, the rate of downcrossings of level $0+$ is $\lambda
_{d}\int_{0}^{\infty }\bar{F}_{d}(y)g_{m}(y)\dd y+\alpha _{b}g_{m}(0)$. That
means that the expected length of the non-emptiness period is given by $%
[\lambda _{d}\int_{0}^{\infty }\bar{F}_{d}(y)g_{m}(y)\dd y+\alpha
_{b}g_{m}(0)]^{-1}$. Thus%
\[
\E[B_{b}]=\frac{1-\pi _{0}}{\alpha _{b}g(0)+\lambda _{d}\int_{0}^{\infty }\bar{%
F}_{d}(y)g(y)\dd y}.
\]%
The expression for $\E[B_{d}]$ is obtained by symmetry.
\\
(ii) Define\ a \textit{cycle} in the real process $X$ (not the
modified process $X_m$) as the time between two upcrossings of
level $0+$. By definition, the emptiness period plus the non-emptiness
period is a cycle in $X$. That means that the expected length of
the emptiness period is the expected length of the cycle\ minus the expected
length of the non-emptiness period. The non-emptiness period in $X$
and in $X_m$ are identical and the length of the expected cycle
is $[\lambda _{b}\int_{0}^{\infty }\bar{F}_{b}(y)f(y)\dd y+\lambda _{b}\pi
_{0}]^{-1}$, since $\lambda _{b}\int_{0}^{\infty }\bar{F}%
_{b}(y)f(y)\dd y+\lambda _{b}\pi _{0}$\ \ is the rate of the upcrossings of
level $0+$.\ \ We obtain%
\[
\E[I_{b}] + \E[B_b] = \frac{1}{\lambda _{b}\int_{0}^{\infty }\bar{F}_{b}(y)f(y)dy+\lambda
_{b}\pi _{0}},
\]
yielding $\E[I_b]$.
$\E[I_{d}]$\ is obtained by symmetry.
\end{proof}

For the special case in which $\xi_b=\xi_d=\xi$ and $\alpha_b=\alpha_d = 0$, we are able to deduce that the expected steady-state inventory level $\E[X]$ has a simple form. 
\begin{proposition}\label{prop:mean_inventory}
If $\xi_b=\xi_d=\xi$ and $\alpha_b=\alpha_d=0$, then 
\begin{equation}
\E[X] = m/\xi,
\end{equation} 
where $m=\l_b\E[B] - \l_d\E[D]$.
\end{proposition}
\begin{proof}
We study the discrete-time embedding of the blood inventory process \\ \noindent $\{X_k\}_{k\geq 1}$, where $X_k$ denotes the blood inventory level \textit{just before} the $k^{th}$ arrival (either blood or demand).
Suppose the process is in steady state. 
By the PASTA property, we have that $X_k \equalD X$ for all $k\geq 1$.
Also, the process $\{X_k\}_{k\geq 1}$ constitutes a Markov chain, of which the evolution is characterized by the recursion
\begin{equation}
X_{k+1} = \left( X_k + \mathbbm{1}_{k,b}B_k - \mathbbm{1}_{k,d} D_k \right)\cdot \ee^{-\xi A_k},
\label{eq:X_recursion}
\end{equation}
where $\mathbbm{1}_{k,b}$ and $\mathbbm{1}_{k,d}$ denote the indicator function of the event that the $k^{th}$ arrival is a blood or demand arrival, respectively. 
Remark that the relation holds for both $X_k  \geq 0$ and $X_k <0$. 
Furthermore, $B_k$ and $D_k$ denote the amount of blood or demand in the $k^{th}$ jump, respectively, and $A_k$ denotes the interarrival time between the $k^{th}$ and $(k+1)^{th}$ arrival. 
Note that $A_k$ is the minimum of two exponentially distributed random variables with rate $\l_b$ and $\l_d$, so that $A_k$ is exponentially distributed with rate $\l_b+\l_d$. 
Next, we take the expectation on both sides of \eqref{eq:X_recursion}, which gives
\begin{equation}
\E[X_{k+1}] = \big( \E[X_k] + p_{k,b}\E[B] - p_{k,d}\E[D]\big)\,\E\big[\ee^{-\xi A_k}\big].
\label{eq:X_recursion_mean}
\end{equation}
Here, we used independence between Poisson processes and their jump sizes, and their memoriless property, and $p_{k,b} = \l_b/(\l_b+\l_d)$ and $p_{k,d} = \l_d/(\l_b+\l_d)$ denote probability of the $k^{th}$ jump being either a blood delivery or demand, respectively. 
Since $X_{k} \equalD X$, we have $\E[X_{k+1}] = \E[X_k] = \E[X]$, and thus we may rewrite \eqref{eq:X_recursion_mean} as
\begin{equation}
\E[X] = \left( \E[X] + \frac{\l_b\E[B] - \l_d\E[D]}{\l_b+\l_d} \right)\cdot\frac{\l_b+\l_d}{\l_b+\l_d+\xi}, 
\end{equation}
from which we easily deduce $\E[X] = (\l_b\E[B] - \l_d \E[D])/\xi = m/\xi$. 
\end{proof}

\subsection{The exponential case}
\label{sectionexp}

\textbf{Density functions.}
We assume in this section that $\bar{F}_b(x) = {\rm e}^{-\mu_b x}$
and
$\bar{F}_d(x) = {\rm e}^{-\mu_d x}$.
Let $\rho_d:= \lambda_d/\mu_d$ and $\rho_b:= \lambda_b/\mu_b$ denote the expected amount of demand requested, and amount of blood delivered into the system, per time unit. 
Moreover, we take $\alpha_b = \alpha_d = 0$.
Under these assumptions, we can
solve (\ref{eq:demand}) and (\ref{eq:blood}) explicitly.

Equations (\ref{eq:blood}) and (\ref{eq:demand}) reduce to:
\begin{align}
&\lambda_d  \int_0^v f(y) {\rm e}^{-\mu_d(v-y)} {\rm d}y 
+
\lambda_d  {\rm e}^{-\mu_d v} \int_0^{\infty} g(y) {\rm e}^{-\mu_d y} {\rm d}y \nonumber \\
&\qquad =
\lambda_b  \int_v^{\infty} f(y) {\rm e}^{-\mu_b(y-v)} {\rm d}y 
+
\xi_d v f(v),
\label{eq:demand2}
\end{align}
\begin{align}
&\lambda_b  \int_0^v g(y) {\rm e}^{-\mu_b(v-y)} {\rm d}y 
+
\lambda_b  {\rm e}^{-\mu_b v} \int_0^{\infty} f(y) {\rm e}^{-\mu_b y} {\rm d}y \nonumber \\
& \qquad =
\lambda_d  \int_v^{\infty} g(y) {\rm e}^{-\mu_d(y-v)} {\rm d}y 
+
\xi_b v g(v),
\label{eq:blood2}
\end{align}
for $v>0$.
%
%
In our analysis, we concentrate on the derivation of $f(v)$. Notice that, once $f(\cdot)$ has been determined, $g(\cdot)$ follows by swapping parameters (symmetry).

In Appendix \ref{app:transformation_int} we show how the integral equations \eqref{eq:demand2}-\eqref{eq:blood2} can be translated into the following decoupled second order differential equations:
\begin{align}
\xi_d v f''(v) &+ \left(2\xi_d -\lambda_d -\lambda_b + \mu_d\xi_dv -\mu_b \xi_d v\right)f'(v)  \nonumber \\
& \qquad + \left(\mu_d\xi_d -\mu_b\xi_d -\mu_d\lambda_b + \mu_b \lambda_d -\mu_b \mu_d \xi_d v\right)f(v) =0
\label{eq:demand6}
\end{align}
and 
\begin{align}
\xi_b v g''(v) &+ \left(2\xi_b -\lambda_d -\lambda_b + \mu_b\xi_dv -\mu_d \xi_b v\right)g'(v)  \nonumber \\
& \qquad + \left(\mu_b\xi_b -\mu_d\xi_b -\mu_b\lambda_d + \mu_d \lambda_b -\mu_d \mu_b \xi_b v\right)g(v) =0,
\label{eq:blood6_1}
\end{align}

with the additional constraint (obtained by applying the level crossing identity for level $v=0$ in either (\ref{eq:demand2}) or (\ref{eq:blood2})):
\begin{equation}
\lambda_b  \int_0^{\infty} f(y) {\rm e}^{-\mu_b y} {\rm d}y 
=
\lambda_d  \int_0^{\infty} g(y) {\rm e}^{-\mu_dy} {\rm d}y .
\label{eq:blood2a}
\end{equation}
Equation (\ref{eq:demand6}) describes a known type of second order differential equation, namely the \textit{extended confluent hypergeometric equation} \cite{Slater1960}, which allows an explicit solution. 
A detailed deduction of the solution to \eqref{eq:demand6} is given in Appendix \ref{app:proof_prop_density}, and yields the following result.
\\
\begin{proposition}\label{densityProp}
The probability density functions of the amount of demand $X_d$ and the amount of blood present $X_b$ are given by
\begin{align}
f(v) &= \pi_d\, \frac{\Gamma\left(1+\frac{\lambda_b}{\xi_d}\right)}{\Gamma\left(\frac{\lambda_b+\lambda_d}{\xi_d}\right)}\,\frac{\ee^{-\mu_d v}U\left( 1-\tfrac{\lambda_d}{\xi_d}, 2-\tfrac{\lambda_b+\lambda_d}{\xi_d},(\mu_b+\mu_d)v\right)}{ _2F_1\left(1-\tfrac{\lambda_d}{\xi_d},1,1+\tfrac{\lambda_b}{\xi_d},-\tfrac{\mu_b}{\mu_d}\right)}\label{eq:fullf},\\
g(v) &= \pi_b\, \frac{\Gamma\left(1+\frac{\lambda_d}{\xi_b}\right)}{\Gamma\left(\frac{\lambda_b+\lambda_d}{\xi_b}\right)}\,\frac{\ee^{-\mu_b v}U\left( 1-\tfrac{\lambda_b}{\xi_b}, 2-\tfrac{\lambda_b+\lambda_d}{\xi_b},(\mu_b+\mu_d)v\right)}{ _2F_1\left(1-\tfrac{\lambda_b}{\xi_b},1,1+\tfrac{\lambda_d}{\xi_b},-\tfrac{\mu_d}{\mu_b}\right)}\label{eq:fullg},
\end{align}
for $v>0$, respectively.
\end{proposition}
\noindent 
Here, $\Gamma(\cdot)$ denotes the gamma function, $ _2F_1(a,b,c,z)$ is the Gaussian hypergeometric function, defined as
\begin{equation}
_2F_1(a,b,c,z) = \sum_{n=0}^\infty \frac{(a)_n(b)_n}{(c)_n\, n!}\, z^n
\end{equation}
and $U(a,b,z)$ is Tricomi's confluent hypergeometric function, see \cite{Slater1960},
\begin{align}
U(a,b,x) &= \frac{\Gamma(b-1)}{\Gamma(1+a-b)}\,\sum_{n=0}^\infty \frac{(a)_n}{(b)_n n!} x^n + \frac{\Gamma(b-1)}{\Gamma(a)}\,x^{1-b}\,\sum_{n=0}^\infty \frac{(1+a-b)_n}{(2-b)_n n!} x^n ,
\end{align}
in which $(a)_n$ is the Pochhammer symbol, defined as $(a)_n = a\cdot(a+1)\cdots(a+n-1)$.
As a direct consequence of Proposition \ref{densityProp}, we obtain expressions for the LTs $\phi(s) = \int_0^{\infty} {\rm e}^{-sv} f(v) {\rm d}v$ and $\gamma(s) = \int_0^{\infty} {\rm e}^{-sv} g(v) {\rm d}v$ for ${\rm Re}\,s \geq 0$ through \cite[Eq.~(3.2.51)]{Slater1960}, which we state here for future use.
\begin{corollary}\label{cor:lsts}
The Laplace transforms for $X_d$ and $X_b$, for ${\rm Re}\,s \geq 0$, are given by
\begin{align}
\label{eq:diffLST}
\phi(s) &= \pi_d\, \frac{\mu_d}{\mu_d+s}\frac{ _2F_1\left(1-\frac{\lambda_d}{\xi_d},1,1+\frac{\lambda_b}{\xi_d},\frac{s-\mu_b}{s+\mu_d}\right)}{ _2F_1\left(1-\frac{\lambda_d}{\xi_d},1,1+\frac{\lambda_b}{\xi_d},-\frac{\mu_b}{\mu_d}\right)},\\
\gamma(s) &= \pi_b\, \frac{\mu_b}{\mu_b+s}\frac{ _2F_1\left(1-\frac{\lambda_b}{\xi_b},1,1+\frac{\lambda_d}{\xi_b},\frac{s-\mu_d}{s+\mu_b}\right)}{ _2F_1\left(1-\frac{\lambda_b}{\xi_b},1,1+\frac{\lambda_d}{\xi_b},-\frac{\mu_d}{\mu_b}\right)},
\label{eq:diffLST_2}
\end{align}
respectively.
\end{corollary}
 Last, we obtain expressions for $\pi_d$ and $\pi_b$. These follow immediately by
using the normalization equation $\pi_b + \pi_d=1$
and \eqref{eq:blood2a}, or equivalently,
$\lambda_b\phi(\mu_b) = \lambda_d\gamma(\mu_d)$. By filling in $s=\mu_b$ in \eqref{eq:diffLST},
\begin{align}\label{eq:LSTequal2}
&\pi_d\,\,\frac{\lambda_b\mu_d}{\mu_b+\mu_d}\, _2F_1\left(1-\frac{\lambda_d}{\xi_d},1,1+\frac{\lambda_b}{\xi_d},-\frac{\mu_b}{\mu_d}\right)^{-1} \nonumber\\
&\qquad = \pi_b\,\,\frac{\lambda_d\mu_b}{\mu_b+\mu_d}\,_2F_1\left(1-\frac{\lambda_b}{\xi_b},1,1+\frac{\lambda_d}{\xi_b},-\frac{\mu_d}{\mu_b}\right)^{-1},
\end{align}
where we used that $_2F_1(a,b,c,0) = 1$. Using the normalization equation, we obtain
\begin{equation}
\label{eq:piD}
\pi_d = \frac{\rho_b,_2F_1\left(1-\frac{\lambda_d}{\xi_d},1,1+\frac{\lambda_b}{\xi_d},-\frac{\mu_b}{\mu_d}\right)}
{\rho_d,_2F_1\left(1-\frac{\lambda_d}{\xi_d},1,1+\frac{\lambda_b}{\xi_d},-\frac{\mu_b}{\mu_d}\right) + 
\rho_b\,_2F_1\left(1-\frac{\lambda_b}{\xi_b},1,1+\frac{\lambda_d}{\xi_b},-\frac{\mu_d}{\mu_b}\right)}.
\end{equation}
By substituting this result into both \eqref{eq:fullf} and \eqref{eq:diffLST}, we obtain the full pdf for the blood inventory process in steady-state.

\begin{theorem}\label{thm:full_pdf}
The steady-state pdf of the net inventory level $X$ is given by
\begin{equation}
h(v) = 
\left\{
\begin{array}{ll}
f(-v), & \text{if }v<0,\\
g(v), & \text{if }v\geq 0,
\end{array}
\right.
\end{equation}
where
\begin{align}
\label{eq:ftotal}
f(v) &= \bar{C}^{-1}\,\frac{\,\Gamma\left(1+\frac{\lambda_b}{\xi_d}\right)}{\Gamma\left(\frac{\lambda_b+\lambda_d}{\xi_d}\right)}\, \rho_d \, \ee^{-\mu_d v}\, U\left( 1-\tfrac{\lambda_d}{\xi_d}, 2-\tfrac{\lambda_b+\lambda_d}{\xi_d},(\mu_b+\mu_d)v\right),\\
g(v) &= \bar{C}^{-1}\,\frac{\,\Gamma\left(1+\frac{\lambda_d}{\xi_b}\right)}{\Gamma\left(\frac{\lambda_b+\lambda_d}{\xi_b}\right)}\, \rho_b \, \ee^{-\mu_b v}\, U\left( 1-\tfrac{\lambda_b}{\xi_b}, 2-\tfrac{\lambda_b+\lambda_d}{\xi_b},(\mu_b+\mu_d)v\right),
\end{align}
with
\begin{equation}
\bar{C} = \rho_d \,_2F_1\left(1-\frac{\lambda_d}{\xi_d},1,1+\frac{\lambda_b}{\xi_d},-\frac{\mu_b}{\mu_d}\right) + 
\rho_b\,_2F_1\left(1-\tfrac{\lambda_b}{\xi_b},1,1+\tfrac{\lambda_d}{\xi_b},-\frac{\mu_d}{\mu_b}\right).
\end{equation}
\end{theorem}

\begin{remark}
By applying the Pfaff transformation $_2F_1(a,b,c,z)=$ \\
$(1-z)^{-b}\,_2F_1\left(c-a,b,c,\frac{z}{1-z}\right)$, we may reformulate
\begin{equation}
\label{eq:pfaffTransform}
_2F_1\left(1-\frac{\lambda_d}{\xi_d},1,1+\frac{\lambda_b}{\xi_d},-\frac{\mu_b}{\mu_d}\right) = \frac{\mu_d}{\mu_b+\mu_d}\, _2F_1\left( \frac{\lambda_b+\lambda_d}{\xi_d},1,\frac{\lambda_b}{\xi_d},\frac{\mu_b}{\mu_b+\mu_d}\right),
\end{equation}
so that
\begin{equation}
\label{eq:piDalternative}
\pi_d = \frac{\lambda_d\,_2F_1\left( \tfrac{\lambda_b+\lambda_d}{\xi_d},1,\tfrac{\lambda_b}{\xi_d},\tfrac{\mu_b}{\mu_b+\mu_d}\right)}
{\lambda_d\,_2F_1\left( \tfrac{\lambda_b+\lambda_d}{\xi_d},1,\tfrac{\lambda_b}{\xi_d},\tfrac{\mu_b}{\mu_b+\mu_d}\right) + 
\lambda_b\,_2F_1\left( \tfrac{\lambda_b+\lambda_d}{\xi_b},1,\tfrac{\lambda_d}{\xi_b},\tfrac{\mu_d}{\mu_b+\mu_d}\right)}.
\end{equation}
By also transforming the hypergeometric term in the numerator of \eqref{eq:fullf}, we get an equivalent form of \eqref{eq:ftotal}, namely
\begin{equation}
\label{eq:ftotalAlternative}
f(v) = \bar{C}^{-1}_{\rm alt}\frac{\Gamma\left(1+\frac{\lambda_b}{\xi_d}\right)}{\Gamma\left(\frac{\lambda_b+\lambda_d}{\xi_d}\right)}\, \rho_b\mu_b(\mu_b+\mu_d)\, \ee^{-\mu_d v}\,U\left( 1-\frac{\lambda_d}{\xi_d}, 2-\frac{\lambda_b+\lambda_d}{\xi_d},(\mu_b+\mu_d)v\right),
\end{equation}
with 
\begin{equation}
\bar{C}_{\rm alt} = \lambda_d\,_2F_1\left( \frac{\lambda_b+\lambda_d}{\xi_d},1,\frac{\lambda_b}{\xi_d},\frac{\mu_b}{\mu_b+\mu_d}\right) + 
\lambda_b\,_2F_1\left( \frac{\lambda_b+\lambda_d}{\xi_b},1,\frac{\lambda_d}{\xi_b},\frac{\mu_d}{\mu_b+\mu_d}\right).
\end{equation}
As a consequence, \eqref{eq:diffLST} is given by
\begin{equation}
\phi(s) = \pi_d \,\frac{ _2F_1\left( \tfrac{\lambda_b+\lambda_d}{\xi_d},1,\tfrac{\lambda_b}{\xi_d},\tfrac{\mu_b-s}{\mu_b+\mu_d}\right)}{ _2F_1\left( \tfrac{\lambda_b+\lambda_d}{\xi_d},1,\tfrac{\lambda_b}{\xi_d},\tfrac{\mu_b}{\mu_b+\mu_d}\right) } = \bar{C}^{-1}_{\rm alt}\,\lambda_d \,_2F_1\left( \frac{\lambda_b+\lambda_d}{\xi_d},1,\frac{\lambda_b}{\xi_d},\frac{\mu_b-s}{\mu_b+\mu_d}\right).
\end{equation}
\end{remark}

Based on the density functions in Theorem \ref{thm:full_pdf}, we make some comments on its properties, and discuss parameter settings that leads to special cases.

By close inspection of these derived density functions, we can observe the following on the distribution shape around $z=0$.
The confluent hypergeometric function $U(a,b,z)$ has limiting form as $z\rightarrow 0$, 
\begin{equation}\label{eq:limit0}
U(a,b,z) = \frac{\Gamma(1-b)}{\Gamma(a-b+1)} +  \frac{\Gamma(b-1)}{\Gamma(a)}\,z^{1-b} + O(z^{2-b}), \qquad b\leq 2, 
\end{equation}
see \cite[Sub.~13.2]{NIST}. 
Note that in our model, $b = 2-(\lambda_b+\lambda_d)/\xi_d<2$ for all parameter settings. 
Equation \eqref{eq:limit0} shows that $U(a,b,z)$  has a singularity at $z=0$ if Re$(b)>1$, which in our case translates to $f(v)$ and $g(v)$ being analytic at $v=0$ if $\lambda_b+\lambda_d > \xi_d$ and $\lambda_b+\lambda_d > \xi_b$, respectively. Assuming $\lambda_b+\lambda_d > \max\{\xi_b,\xi_d\}$, \eqref{eq:limit0} also implies that
\begin{align}\label{eq:contLimit1}
\lim_{v\rightarrow 0} f(v) &=  \bar{C}^{-1}\,\frac{\,\Gamma\left(1+\frac{\lambda_b}{\xi_d}\right)}{\Gamma\left(\frac{\lambda_b+\lambda_d}{\xi_d}\right)}\, \lambda_d\mu_b\cdot \frac{\Gamma\left(\frac{\lambda_b+\lambda_d}{\xi_d}-1\right)}{\Gamma\left(\frac{\lambda_b}{\xi_d}\right)}\\ 
&= \bar{C}^{-1}\,\frac{\frac{\lambda_b}{\xi_d}}{\frac{\lambda_b+\lambda_d}{\xi_d}-1}\, \lambda_d\mu_b= \bar{C}^{-1}\,\frac{\lambda_b\lambda_d\mu_b\mu_d}{\lambda_b+\lambda_d-\xi_d}.
\nonumber
\end{align}
Similarly,
\begin{equation}
\lim_{v\rightarrow 0} g(v) = \bar{C}^{-1}\,\frac{\lambda_b\lambda_d\mu_b\mu_d}{\lambda_b+\lambda_d-\xi_b}.
\label{eq:contLimit2}
\end{equation}
By equating these two expressions, we conclude that $\lim_{v\rightarrow 0} f(v) = \lim_{v\rightarrow 0} g(v)< \infty$, i.e. the overall density function $h(v)$ is continuous at $v=0$, if and only if $\xi_b = \xi_d$.  \\
\noindent
The asymptotic behavior of $U$ as $z\to\infty$ is given by \cite[p.~60]{Slater1960},
\begin{equation}
U(a,b,z) \sim z^{-a}, ~~~~~~~~~~~~~~~~ z \rightarrow \infty,
\end{equation}
which implies that the density function tail decays as
\begin{equation}
\label{eq:asympt}
f(v) \sim C^*\, e^{-\mu_d v}\, v^{\lambda_d/\xi_d-1}, ~~~~~~~~~~~~~~~~~~~
v\rightarrow\infty ,
\end{equation}
for some constant $C^*$.

\noindent\textbf{Special cases.}
Equation \eqref{eq:asympt} suggests that the case $\lambda_d = \xi_d$ is special.
Indeed, then \eqref{eq:diffLST} reduces to
\begin{equation}
\label{eq:lambdaisxi}
\phi(s) = \bar{C}^{-1}\,\lambda_d\mu_b\, \frac{\mu_d}{\mu_d+s} = \pi_d\, \frac{\mu_d}{\mu_d+s},
\end{equation}
where we used that $_2F_1(0,a,b,z) = 1$ for all $a,b,z$. Hence, conditioned on being positive, the amount of demand present is exponentially distributed with parameter $\mu_d$, regardless of the values of $\lambda_d = \xi_d$, as well as $\lambda_b,\, \xi_b,$ and $\mu_b$. 
If we moreover let $\lambda_b = \xi_b$, then 
\[
\pi_d = \frac{\lambda_d/\mu_d}{\lambda_b/\mu_b + \lambda_d/\mu_d} = \frac{\rho_d}{\rho_b+\rho_d},
\]
and $X$ has exponential distribution both above and below 0, with parameters $\mu_b$ and $\mu_d$, respectively. \\

A second special case arises when the process is symmetric, that is, $\lambda_b=\lambda_d=\lambda$, $\mu_b=\mu_d=\mu$ and $\xi_b=\xi_d=\xi$. Obviously, we get $\pi_b= \pi_d = \tfrac{1}{2}$ due to the symmetry. If we define $\eta := \lambda/\xi$, 
\begin{align}
f(v) &= \frac{\Gamma(1+\eta)\, \mu e^{-\mu v}\, U\left(1-\eta,2(1-\eta),2\mu v\right)}{2\,\Gamma(2\eta) _2F_1\left(2\eta,1,1+\eta,\tfrac{1}{2}\right)}\\
&= \frac{\Gamma(1+\eta)}{2\,\Gamma(2\eta) _2F_1\left(2\eta,1,1+\eta,\tfrac{1}{2}\right)}\, \frac{\mu}{2\sqrt{\pi}}\, \left(2\mu v\right)^{\eta-\tfrac{1}{2}}\, K_{\tfrac{1}{2}-\eta}\left(\mu v\right),
\nonumber
\end{align}
where $K_\alpha(\cdot)$ is the modified Bessel function of the second kind, see \cite[Eq.~(13.6.10)]{NIST}.\\ 
\\*
\noindent
\textbf{Performance measures.}
Based on Theorem \ref{thm:full_pdf}, we can directly derive a couple of characteristics of the process.
First, we consider the mean inventory level

\begin{figure}
\centering
\begin{subfigure}[b]{0.45\textwidth}
\centering
\begin{tikzpicture}[scale=0.78]
\begin{axis}[
	xmin = -0.02,
	xmax = 1,
	ymin = -0.02,
	ymax = 5,
	axis line style={->},
	axis lines = middle,
	yscale = 0.8,
	xscale = 1,
	xlabel = {$\xi$},
	xlabel near ticks,
	legend cell align = left,
	legend style = {at = {(axis cs: 1,5)},anchor = north east}
	]

\addplot[black,dashed,thick] table[x=x,y=Xd] {Chapter_7/tikz/means1.txt};
\addplot[black,dotted,thick] table[x=x,y=Xb] {Chapter_7/tikz/means1.txt};
\addplot[black,thick] table[x=x,y=Q] {Chapter_7/tikz/means1.txt};

\legend{{$\E[X_d|X_d>0]$},{$\E[X_b|X_b>0]$},{$\E[X]$}};

\end{axis}
\end{tikzpicture}
\caption{As a function of $\xi$}
\end{subfigure}
\hspace{5mm}
\begin{subfigure}[b]{0.45\textwidth}
\centering
\begin{tikzpicture}[scale=0.78]
\begin{axis}[
	xmin = -0.02,
	xmax = 2,
	ymin = -0.02,
	ymax = 5,
	axis line style={->},
	axis lines = middle,
	yscale = 0.8,
	xscale = 1,
	xlabel = {$\theta$},
	xlabel near ticks,
	legend cell align = left,
	legend style = {at = {(axis cs: 2,5)},anchor = north east}
	]

\addplot[black,dashed,thick] table[x=x,y=Xd] {Chapter_7/tikz/means2.txt};
\addplot[black,dotted,thick] table[x=x,y=Xb] {Chapter_7/tikz/means2.txt};
\addplot[black,thick] table[x=x,y=Q] {Chapter_7/tikz/means2.txt};

\legend{{$\E[X_d|X_d>0]$},{$\E[X_b|X_b>0]$},{$\E[X]$}};

\end{axis}
\end{tikzpicture}
\caption{As a function of $\theta$}
\end{subfigure}
\caption{Expected mean amount of blood, demand, and net blood present.}
\label{fig:means}
\end{figure}
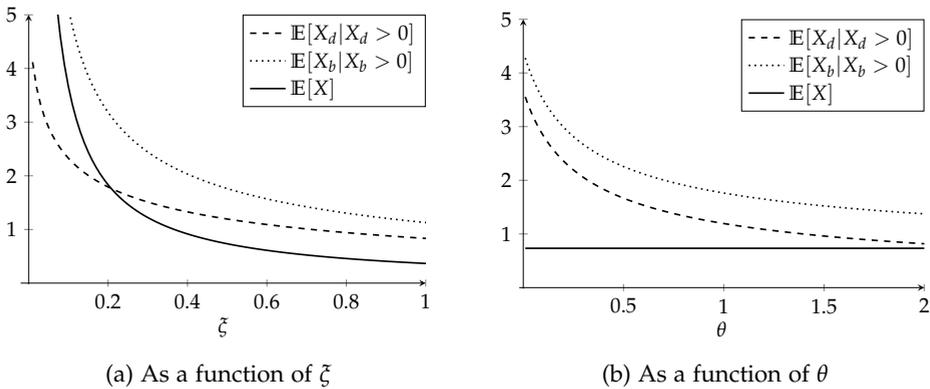

\begin{corollary}\label{cor:means}
The expected amount of demand (blood) present, given that it is positive equals
\begin{align}
\E[X_d|X_d>0] &= \frac{1}{\xi_d}\left[ \rho_d - \rho_b + \rho_b\, _2F_1\left(1-\frac{\lambda_d}{\xi_d},1,1+\frac{\lambda_b}{\xi_d},{-}\frac{\mu_b}{\mu_d}\right)^{-1}\right],\label{EXd>0}\\
\E[X_b|X_b>0] &= \frac{1}{\xi_b}\left[ \rho_b - \rho_d + \rho_d \, _2F_1\left(1-\frac{\lambda_b}{\xi_b},1,1+\frac{\lambda_d}{\xi_b},{-}\frac{\mu_d}{\mu_b}\right)^{-1}\right]\label{EXb>0}.
\end{align}
Accordingly, the expected net amount of blood present equals
\begin{equation}
\E[X] = \left(\rho_b-\rho_d \right)\left(\frac{\pi_b}{\xi_b}+\frac{\pi_d}{\xi_d}\right)+\frac{\lambda_b\lambda_d}{\bar{C}}\left(\frac{1}{\xi_b}-\frac{1}{\xi_d}\right).
\end{equation}
\end{corollary}
\begin{proof}
Let us use shorthand notation 
\[
F(s) = \Big(1-\frac{\lambda_d}{\xi_d},1,1+\frac{\lambda_b}{\xi_d},\frac{s-\mu_b}{s+\mu_d}\Big)
,\]
so that 
\[\phi(s)=\pi_d\, \frac{\mu_b}{\mu_b+s}\, \frac{F(s)}{F(0)}.\]
Through \cite[Eq.~(15.5.20)]{NIST},
\begin{equation}\label{eq:proof_a}
\frac{\dd}{\dd z}\, _2F_1(a,1,c,z) = \frac{c-1}{z(1-z)} + \frac{1-c+az}{z(1-z)} \, _2F_1(a,1,c,z),
\end{equation}
where we also used that $ _2F_1(a,1,c,z) = 1$. 
Then, 
\begin{align*}
\frac{\phi'(0)}{\pi_d} &= \left[ \frac{-\mu_d}{(\mu_d+s)^2} \, \frac{F(s)}{F(0)} 
+ \frac{\mu_d}{\mu_d+s}\, \frac{F'(s)}{F(0)} \right]_{s=0} = {-}\frac{1}{\mu_d} + \frac{F'(0)}{F(0)}.
\end{align*}
By \eqref{eq:proof_a}, we find
\begin{align*}
F'(s) &= \Big( \frac{\lambda_b/\xi_d}{\frac{s-\mu_b}{s+\mu_d}\cdot \frac{\mu_b+\mu_d}{s+\mu_d}} + \frac{-\lambda_b/\xi_d + (1-\lambda_d/\xi_d)\frac{s-\mu_b}{s+\mu_d}}{\frac{s-\mu_b}{s+\mu_d}\cdot \frac{\mu_b+\mu_d}{s+\mu_d}}\, F(s)\Big)
 \,\frac{\dd}{\dd s} \Big[ \frac{s-\mu_b}{s+\mu_d} \Big]
 \\
&= \Big( \frac{ \lambda_b}{\xi_d} + \left[\frac{{-}\lambda_b}{\xi_d} + \Big(1-\frac{\lambda_d}{\xi_d}\Big)\frac{s-\mu_b}{s+\mu_d}\right] F(s) \Big) \frac{ (s+\mu_d)^2} {(s-\mu_b)(\mu_b+\mu_d)}\cdot \frac{\mu_b+\mu_d}{(s+\mu_d)^2}\\
&= \Big( \frac{ \lambda_b}{\xi_d} + \left[{-}\frac{\lambda_b}{\xi_d} + \Big(1-\frac{\lambda_d}{\xi_d}\Big)\frac{s-\mu_b}{s+\mu_d}\right] F(s) \Big)
\frac{1}{s-\mu_b},
\end{align*}
so that
\begin{align*}
F'(0) &= {-} \frac{\lambda_d/\mu_b}{\xi_d} + \left( \frac{\lambda_d/\mu_b}{\xi_d}
+ \frac{1}{\mu_d} - \frac{\lambda_d/\mu_d}{\xi_d}\right)F(0)\\
&= {-}\frac{\rho_b}{\xi_d} + \left( \frac{\rho_b-\rho_d}{\xi_d}
+ \frac{1}{\mu_d}\right)F(0).
\end{align*}
Hence, we find
\begin{align*}
\E[X_d|X_d>0] &= {-}\frac{\phi'(0)}{\pi_d} = \frac{1}{\mu_d}- \frac{1}{F(0)}\left[{-}\frac{\rho_b}{\xi_d} + \left( \frac{\rho_b-\rho_d}{\xi_b}
+ \frac{1}{\mu_d}\right)F(0)\right]\\
&= \frac{1}{\xi_d}\left( \rho_d-\rho_b + \rho_b/F(0)\right) = \frac{1}{\xi_d}\left( -m + \rho_b/F(0)\right),
\end{align*}
which equals \eqref{EXd>0}. 
The expression for \eqref{EXb>0} follows by symmetry.
Furthermore,
 \begin{align*}
 \E[X] &= \pi_b \E[X_b|X_b>0] + \pi_d \E[-X_d|X_d>0]\\
 &= m\left[\frac{\pi_b}{\xi_b}+\frac{\pi_d}{\xi_d}\right] + \frac{\lambda_d}{\mu_d\,\xi_b}\frac{\pi_b}{_2F_1\left(1-\tfrac{\lambda_b}{\xi_b},1,1+\tfrac{\lambda_d}{\xi_b},{-}\tfrac{\mu_d}{\mu_b}\right)}\\
 &\qquad - \frac{\lambda_b}{\mu_b\,\xi_d}\,\frac{\pi_d}{_2F_1\left(1-\tfrac{\lambda_d}{\xi_d},1,1+\tfrac{\lambda_b}{\xi_d},{-}\tfrac{\mu_b}{\mu_d}\right)}.
 \end{align*}
 Note that $\pi_d\,_2F_1\left(1-\tfrac{\lambda_d}{\xi_d},1,1+\tfrac{\lambda_b}{\xi_d},{-}\tfrac{\mu_b}{\mu_d}\right)^{-1} = \lambda_d\mu_b\bar{C}^{-1}$. 
 Hence,
 \begin{align*}
 \E[X] &= m\left[\frac{\pi_b}{\xi_b}+\frac{\pi_d}{\xi_d}\right] + \frac{\lambda_b\lambda_d}{\bar{C}}\left(\frac{1}{\xi_b}-\frac{1}{\xi_d}\right),
 \end{align*}
 which completes the proof.
\end{proof}
\begin{remark}
Note that if $\xi_b = \xi_d = \xi$, we get $\E[X] = m(\pi_b+\pi_d)/\xi = m/\xi$, which is consistent with Proposition \ref{prop:mean_inventory}. 
The expression in (\ref{EXd>0}) contains no $\xi_b$. Indeed, while the value of $\xi_b$ influences the probability that $X_d>0$,
it does not influence the mean of $X_d$ given that $X_d >0$.
\end{remark}
In Figure \ref{fig:means}, we plot the behavior of the three performance metrics in Corollary \ref{cor:means} while keeping $m$ fixed. In Figure \ref{fig:means}(a) we set $\lambda_b = 1.2$, $\lambda_d = 1$, $\mu_b=1$, $\mu_d=1.2$, so that $m = 11/30$ and vary $\xi_b=\xi_d=\xi$ between 0 and 1. In Figure \ref{fig:means}b, we fix $\xi_b=\xi_d=0.5$ and take $\lambda_b = 1.2\theta$, $\lambda_d = \theta$, $\mu_b=\theta$, $\mu_d=1.2\theta$, so that still $m=11/30$, and vary $\theta$.
Observe that in Figure \ref{fig:means}b, $\E[X]$ is constant, since the value of $m/\xi$ if unaffected by the parameter $\theta$. 

Secondly, we present the probability of positive (cq. negative) inventory.
\begin{corollary}\label{cor:pid}
The probability of positive (cq.~negative) inventory is given by,
\begin{align}
\pi_b &= \frac{\rho_b\,_2F_1\left(1-\frac{\lambda_b}{\xi_b},1,1+\frac{\lambda_d}{\xi_b},-\frac{\mu_d}{\mu_b}\right)}
{\rho_d\,_2F_1\left(1-\frac{\lambda_d}{\xi_d},1,1+\frac{\lambda_b}{\xi_d},-\frac{\mu_b}{\mu_d}\right) + 
\rho_b\,_2F_1\left(1-\frac{\lambda_b}{\xi_b},1,1+\frac{\lambda_d}{\xi_b},-\frac{\mu_d}{\mu_b}\right)},\\
\pi_d &= \frac{\rho_d\,_2F_1\left(1-\frac{\lambda_d}{\xi_d},1,1+\frac{\lambda_b}{\xi_d},-\frac{\mu_b}{\mu_d}\right)}
{\rho_d\,_2F_1\left(1-\frac{\lambda_d}{\xi_d},1,1+\frac{\lambda_b}{\xi_d},-\frac{\mu_b}{\mu_d}\right) + 
\rho_b\,_2F_1\left(1-\frac{\lambda_b}{\xi_b},1,1+\frac{\lambda_d}{\xi_b},-\frac{\mu_d}{\mu_b}\right)},
\end{align}
respectively.
\end{corollary}
\begin{proof}
The expressions follow directly from \eqref{eq:piD} and $\pi_b = 1-\pi_d$. 
\end{proof}
The last relevant performance indicator we consider is the fraction of demand that is immediately satisfied from stock. 

\begin{corollary}
The probability that a demand request can be fully satisfied from stock is given by
\begin{align}
\P({\rm demand\  satisfied}) = \bar{C}^{-1}\rho_b \left( _2F_1\left(1-\tfrac{\lambda_d}{\xi_d},1,1+\tfrac{\lambda_b}{\xi_d},{-}\tfrac{\mu_b}{\mu_d}\right) - \frac{\mu_b}{\mu_b+\mu_d}\right).
\end{align}
\end{corollary}
\begin{proof}
Using the PASTA property of the Poisson process, we get
\begin{align*}
\P({\rm demand\ satisfied}) &= \P( X > D ) = \P(X_b > D) \\
&= \int_0^\infty g(u) (1-{\rm e}^{-\mu_d u})\, \dd u  = \pi_b - \gamma(\mu_d).
\end{align*}
Substituting the expressions for $\pi_b$ as in Corollary \ref{cor:pid} and $\gamma(\mu_b)$ as in \eqref{eq:diffLST_2} yields the result. 
\end{proof}

\subsection{The general case}
\label{gener}
In this section we outline how the integral equations (\ref{eq:blood}) and (\ref{eq:demand})
can be solved using Laplace transforms, when we make the restriction that $F_b(\cdot)$ and $F_d(\cdot)$ 
are Coxian distributions.
This is not a major restriction, because the class of Coxian distributions lies dense in the class of all distributions
of non-negative random variables, see e.g.~\cite[Sec.~III.4]{Asmussen2003}.
Hence, one can approximate $F_b(\cdot)$ arbitrarily closely by a Coxian distribution.

If $X_i$, $i=1,2,\dots,K$ are independent, exponentially distributed random variables,
and $\E[X_i] = \frac{1}{\b_i}$, $i=1,2,\dots,K$, then a Coxian amount of blood $B$ can be represented
as:
\begin{equation}
B = \sum_{j=1}^i X_j\quad {\rm  with ~ probability } \quad p_i \prod_{j=1}^{i-1} (1-p_j), \quad i=1,2,\dots,K.
\end{equation}
In the above case, it is easily verified that one can represent $\bar{F}_b(x)$ as follows:
\begin{equation}
\bar{F}_b(x) = \P(B>x) = \sum_{i=1}^K p_i \prod_{h=1}^{i-1} (1-p_h) \sum_{j=1}^i
{\rm e}^{-\b_j x}
\prod_{l=1; l \neq j}^i \frac{\b_l}{\b_l - \b_j}
, 
\label{Fbarb}
\end{equation}
if all $\b_j$ are different. If two $\b_j$ coincide, then a term with $x {\rm e}^{-\b_j x}$ (Erlang-$2$) must be added.
We leave this to the reader, but in Remark~\ref{RmErlang} below we outline how Erlang terms can be handled in solving the integral equations
(\ref{eq:demand}) and (\ref{eq:blood}).
The counterpart of (\ref{Fbarb}) for the case that $F_d(\cdot)$ is Coxian, is
\begin{equation}
\bar{F}_d(x) = \P(D>x) = \sum_{i=1}^L q_i \prod_{h=1}^{i-1} (1-q_h) \sum_{j=1}^i
{\rm e}^{-\d_j x}
\prod_{l=1; l \neq j}^i \frac{\d_l}{\d_l - \d_j}
. 
\label{Fdarb}
\end{equation}

Taking Laplace transforms $\phi(s) = \int_0^{\infty} {\rm e}^{-sy} f(y) {\rm d}y$ and
$\gamma(s) = \int_0^{\infty} {\rm e}^{-sy} g(y) {\rm d}y$ 
in (\ref{eq:blood}) and (\ref{eq:demand})
results in first-order inhomogeneous differential equations in $\phi(s)$ and $\gamma(s)$, respectively, which can be solved in a straightforward way.

\begin{equation}
\phi'(s) = A_H(s) \phi(s) + A_I(s),
\label{diffeq}
\end{equation}
with the homogeneous term $A_H(s)$ being given by
\begin{align}
A_H(s) &:= -\frac{1}{\xi_d} 
\left[
\lambda_d  
\sum_{i=1}^L q_i \prod_{h=1}^{i-1} (1-q_h) \sum_{j=1}^i
\frac{1}{\d_j+s}
\prod_{l=1; l \neq j}^i \frac{\d_l}{\d_l - \d_j}
\right. 
\nonumber
\\
&\qquad - \lambda_b 
\left.\sum_{i=1}^K p_i \prod_{h=1}^{i-1} (1-p_h) \sum_{j=1}^i
\frac{1}{\b_j-s}
\prod_{l=1; l \neq j}^i \frac{\b_l}{\b_l - \b_j}
- \alpha_d \right],
\end{align}
and the inhomogeneous term $A_I(s)$ being given by
\begin{align}
A_I(s) &:=
- \frac{1}{\xi_d}
\left[
\lambda_d
\sum_{i=1}^L q_i \prod_{h=1}^{i-1} (1-q_h) \sum_{j=1}^i
\frac{1}{\d_j+s} [\gamma(\d_j) + \pi_0]
\prod_{l=1; l \neq j}^i \frac{\d_l}{\d_l - \d_j}
\right.
\nonumber
\\
&\qquad + \left.\lambda_b 
\sum_{i=1}^K p_i \prod_{h=1}^{i-1} (1-p_h) \sum_{j=1}^i
\frac{1}{\b_j-s} \phi(\b_j)
\prod_{l=1; l \neq j}^i \frac{\b_l}{\b_l - \b_j}
\right].
\end{align}
\noindent
The solution of (\ref{diffeq}) is given by the following expression:
\begin{equation}
\phi(s) =  \phi(0) {\rm e}^{\int_0^s A_H(z) {\rm d}z} + \int_0^s A_I(u) 
{\rm e}^{\int_u^s A_H(z) {\rm d}z} {\rm d}u , ~~~~  s \geq 0.
\label{diffeqsoln}
\end{equation}
$\gamma(s)$ is given by a mirror expression, where $\phi(0)$ is replaced by $\gamma(0)$
and where $A_H(s)$ and $A_I(s)$ are replaced by expressions in which $K$ and $L$ are swapped, and $p$ and $q$, and $\beta_i$ and $\delta_i$.

It should be noticed that $\phi(0)$, $\gamma(0)$ and $\pi_0$ still have to be determined.
Furthermore, it should be noticed that $A_H(s)$ and $A_I(s)$ have singularities at $s=\beta_1,\dots,\beta_K$.
These singularities are removable, but handling Equation \eqref{diffeqsoln} clearly requires some care.
Instead of working out the details, we shall below return to the case
of exponentially distributed amounts of blood and demand  -- so $K=L=1$.
For that case, we shall not only work out the solution of the differential equation for $\phi(s)$ in detail,
including the determination of the missing constants, but
we also relate the results to those obtained in Section~\ref{sectionexp} 
without resorting to Laplace transforms.
%
Taking $K=1, p_1=1, \delta_1 = \mu_d$, and $L=1, q_1=1, \beta_1 = \mu_b$, we obtain
the following two inhomogeneous first order differential equations in the LTs $\phi(s)$ and $\gamma(s)$:
\begin{equation}\label{eq:firstPhis}
\phi'(s) = \phi(s)\left[\frac{\lambda_b}{\xi_d} \frac{1}{\mu_b-s} - \frac{\lambda_d}{\xi_d} \frac{1}{\mu_d+s}\right]
-\frac{\lambda_b}{\xi_d} \frac{\phi(\mu_b)}{\mu_b-s} -\frac{\lambda_d}{\xi_d} \frac{\gamma(\mu_d)}{\mu_d+s} ,
\end{equation}
\begin{equation}
\gamma'(s) = \gamma(s)\left[\frac{\lambda_d}{\xi_b} \frac{1}{\mu_d-s} - \frac{\lambda_b}{\xi_b} \frac{1}{\mu_b+s}\right]
-\frac{\lambda_d}{\xi_b} \frac{\gamma(\mu_d)}{\mu_d-s} -\frac{\lambda_b}{\xi_b} \frac{\phi(\mu_b)}{\mu_b+s} .
\end{equation}
They are routinely solved:
\begin{align}
\phi(s) &= \left(\frac{\mu_b}{\mu_b-s}\right)^{\frac{\lambda_b}{\xi_d}}
\left(\frac{\mu_d}{\mu_d+s}\right)^{\frac{\lambda_d}{\xi_d}}
\left[\phi(0) \frac{}{}\right. 
\nonumber \\
&\qquad -
\frac{\lambda_d}{\xi_d} \gamma(\mu_d) \int_0^s \left(\frac{\mu_b-z}{\mu_b}\right)^{\frac{\lambda_b}{\xi_d}}
\left(\frac{\mu_d+z}{\mu_d}\right)^{\frac{\lambda_d}{\xi_d}-1} \frac{{\rm d}z}{\mu_d}
\nonumber
\\
&\qquad \qquad -
\frac{\lambda_b}{\xi_d} \phi(\mu_b) \int_0^s \left(\frac{\mu_b-z}{\mu_b}\right)^{\frac{\lambda_b}{\xi_d}-1}
\left(\frac{\mu_d+z}{\mu_d}\right)^{\frac{\lambda_d}{\xi_d}} \frac{{\rm d}z}{\mu_b}\left.\frac{}{}\right].
\label{phis1} 
\end{align}
Similarly,
\begin{align}
\gamma(s) &= \left(\frac{\mu_d}{\mu_d-s}\right)^{\frac{\lambda_d}{\xi_b}}
\left(\frac{\mu_b}{\mu_b+s}\right)^{\frac{\lambda_b}{\xi_b}}
\left[\gamma(0) \frac{}{}\right.
\nonumber
\\
&\qquad -
\frac{\lambda_b}{\xi_b} \phi(\mu_b) \int_0^s \left(\frac{\mu_d-z}{\mu_d}\right)^{\frac{\lambda_d}{\xi_b}}
\left(\frac{\mu_b+z}{\mu_b}\right)^{\frac{\lambda_b}{\xi_b}-1} \frac{{\rm d}z}{\mu_b}
\nonumber
\\
&\qquad \qquad -
\frac{\lambda_d}{\xi_b} \gamma(\mu_d) \int_0^s \left(\frac{\mu_d-z}{\mu_d}\right)^{\frac{\lambda_d}{\xi_b}-1}
\left(\frac{\mu_b+z}{\mu_b}\right)^{\frac{\lambda_b}{\xi_b}} \frac{{\rm d}z}{\mu_d}\left.\frac{}{}\right] .
\label{gammas1}
\end{align}
Notice that the exponents in the above integrals have powers which are larger than $-1$ (e.g., $\frac{\lambda_d}{\xi_d}-1$),
so that these integrals do not lead to singularities.
We still need to determine the two constants $\phi(0)=\pi_d$ and $\gamma(0)=\pi_b$.
Together with $\phi(\mu_b)$ and $\gamma(\mu_d)$, we have four unknowns.
We determine these unknowns using the following four equations:
(i) From (\ref{eq:blood2a}), we get
$\lambda_b \phi(\mu_b) = \lambda_d \gamma(\mu_d)$, while (ii) $\pi_d + \pi_b =1$.
Finally, we take (iii) $s=\mu_b$ in (\ref{phis1}) and (iv) $s=\mu_d$ in (\ref{gammas1}).

Notice that the identity
$\lambda_b \phi(\mu_b) = \lambda_d \gamma(\mu_d)$ allows us to
reduce the two integrals in (\ref{phis1}) to one integral (and similarly in (\ref{gammas1})):
\begin{align}
\phi(s) &= \left(\frac{\mu_b}{\mu_b-s}\right)^{\frac{\lambda_b}{\xi_d}}
\left(\frac{\mu_d}{\mu_d+s}\right)^{\frac{\lambda_d}{\xi_d}}
\left[\phi(0) \frac{}{}\right.
\nonumber
\\
&\quad -
\frac{\lambda_d}{\xi_d} \gamma(\mu_d)\, \frac{\mu_b+\mu_d}{\mu_b\mu_d} \int_0^s \left(\frac{\mu_b-z}{\mu_b}\right)^{\frac{\lambda_b}{\xi_d}-1} \left(\frac{\mu_d+z}{\mu_d}\right)^{\frac{\lambda_d}{\xi_d}-1} {\rm d}z \left.\frac{}{}\right].
\label{phis11}
\end{align}

\begin{remark}
We have numerically verified that
the expressions in (\ref{phis1}) and (\ref{eq:diffLST}) coincide. 

\end{remark}
\begin{remark}
If $\lambda_b=0$ then we have a known queueing model or shot-noise model
with state-dependent service rate, see Keilson \& Mermin \cite{Keilson1959}
and Bekker et al.~\cite{Bekker2004} for the so-called shot noise model.
\end{remark}
\begin{remark}
\label{R7}
The case $\lambda_d = \xi_d$ is special. Equation \eqref{phis1} now reduces to
\begin{align}
\phi(s) &= \left(\frac{\mu_b}{\mu_b-s}\right)^{\frac{\lambda_b}{\lambda_d}}
\frac{\mu_d}{\mu_d+s}
\left[\phi(0) \frac{}{}\right.
\label{phis1A}
-
\gamma(\mu_d) \int_0^s \left(\frac{\mu_b-z}{\mu_b}\right)^{\frac{\lambda_b}{\lambda_d}}
\frac{{\rm d}z}{\mu_d}
\nonumber
\\
&\qquad -
\frac{\lambda_b}{\lambda_d} \phi(\mu_b) \int_0^s \left(\frac{\mu_b-z}{\mu_b}\right)^{\frac{\lambda_b}{\lambda_d}-1}
\frac{\mu_d+z}{\mu_d} \frac{{\rm d}z}{\mu_b}\left.\frac{}{}\right].
\nonumber
\end{align}
Both integrals are easily evaluated (rewrite,
in the last integral, $\mu_d + z = \mu_d + \mu_b -(\mu_b - z)$).
We find
\begin{align}
\phi(s) &=
\left(\frac{\mu_b}{\mu_b-s}\right)^{\frac{\lambda_b}{\lambda_d}}
\frac{\mu_d}{\mu_d+s}\nonumber\\
&\ \cdot 
\left[\phi(0)
+ \frac{\gamma(\mu_d)}{\mu_d} \frac{\lambda_d}{\lambda_b + \lambda_d} \mu_b - \phi(\mu_b) \frac{\mu_d+\mu_b}{\mu_d} - \frac{\phi(\mu_b)}{\mu_d} \frac{\lambda_b}{\lambda_b + \lambda_d} \mu_b \right]
\nonumber
\\
&\ + \frac{\mu_d}{\mu_d+s} \left[
\frac{\gamma(\mu_d)}{\mu_d} \frac{\lambda_d}{\lambda_b + \lambda_d} (\mu_b -s) + \phi(\mu_b) \frac{\mu_d+\mu_b}{\mu_d} - \frac{\phi(\mu_b)}{\mu_d} \frac{\lambda_b}{\lambda_b + \lambda_d} (\mu_b -s)\right] .
\end{align} 
Now observe through \eqref{eq:blood2a}, that  $\lambda_b \phi(\mu_b) = \lambda_d \gamma(\mu_d)$.
Hence, in both lines of the above formula, two terms cancel.
Moreover, $\phi(s)$ should be analytic for $s=\mu_b$, yielding
\begin{equation}
\phi(0) = \phi(\mu_b) \frac{\mu_d + \mu_b}{\mu_d}.
\end{equation}
Finally we obtain, see also \eqref{eq:lambdaisxi},
\begin{equation}
\phi(s) = \frac{\mu_d}{\mu_d+s} \phi(\mu_b) \frac{\mu_d+\mu_b}{\mu_d} = \phi(0) \frac{\mu_d}{\mu_d+s} = \pi_d \frac{\mu_d}{\mu_d + s},
\end{equation}
and hence
\begin{equation}\label{eq:exp}
f(x) = \pi_d  \mu_d {\rm e}^{-\mu_d x}, ~~~ x > 0;
\end{equation}
the shortage (amount of demand present) is exponentially distributed
when $\lambda_d = \xi_d$.
\\
It should be noticed that, if $\lambda_d = \xi_d$, then the first and last term of (\ref{eq:demand2})
are equal when (\ref{eq:exp}) holds; and using (\ref{eq:blood2a})
it is also readily verified that the second and third term of (\ref{eq:demand2}) are equal.
The constant $\pi_d$ will in general still depend on the parameters
$\lambda_d = \xi_d$, $\lambda_b$, $\mu_b$ and $\xi_b$.
\\
We end this remark with the observation that in the one-sided shot-noise process
(so $\lambda_b=0$), Bekker et al.\ \cite{Bekker2004} also observe that $\lambda_d = \xi_d$
results in an exponential density. 
\end{remark}

\subsection{A variant}
\label{sectionvariant}
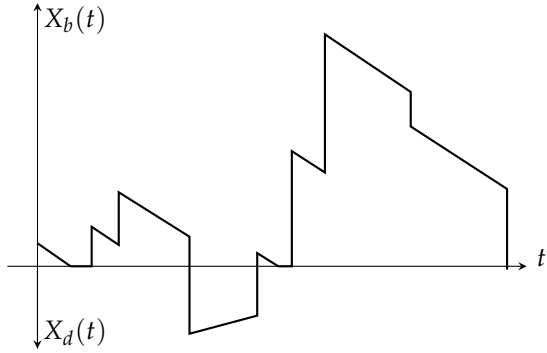
\begin{figure}
\centering
\begin{tikzpicture}
\begin{axis}[
	xmin = -0.2,
	xmax = 3.25,
	ymin = -0.5,
	ymax = 1.6,
	ticks = none,
	axis line style={->},
	axis lines = middle,
	yscale = 0.8,
	xscale = 1
	]

\addplot[black,thick] file {Chapter_7/tikz/sample_path_xi0.txt};

\node at (axis cs: 0.25,-0.4) {$X_d(t)$};
\node at (axis cs: 0.25,1.5) {$X_b(t)$};
\draw[-stealth] (axis cs: 0,0) -- (axis cs:0,-0.5);

\end{axis}

\node at (7.05,1.2) {$t$};
\end{tikzpicture}
\caption{Sample path of the net amount of blood present if $\xi_b = \xi_d = 0$.}
\label{FIG2}
\end{figure}

In this section, we assume that the expiration rate of blood and the patience rate of demand are constant. 
So, we take $\xi_b = \xi_d = 0$. 
A visualization of a possible sample path is depicted in Figure \ref{FIG2}.

We again restrict ourselves to the case of exponentially distributed amounts of demand and of blood deliveries.
We now need to impose stability conditions. 
In the case of positive demand, the drift is towards zero if $\lambda_d \E[D] < \alpha_d + \lambda_b \E [B]$,
while in the case of a positive amount of blood, the drift is towards zero if
$\lambda_b \E [B] < \alpha_b + \lambda_d \E [D]$.
If these two conditions are violated, either the amount of demand or the amount of blood present increases without bound
(see also
below).
In this case, \eqref{eq:demand6} reduces to
\begin{equation}
\alpha_d f''(v) + (-\lambda_d -\lambda_b + \mu_d \alpha_d -\mu_b \alpha_d)f'(v)
+
(-\mu_d\lambda_b + \mu_b \lambda_d -\mu_b \mu_d \alpha_d)f(v) =0.
\label{eq:demand6b}
\end{equation}
Hence $f(\cdot)$ is a mixture of two exponential terms: $f(v) = R_+ {\rm e}^{- x_+ v} + R_- {\rm e}^{- x_- v}$,
where $x_+$ and $x_-$ are the positive and negative root of the equation
\begin{equation}
\alpha_d x^2 -(\mu_d \alpha_d - \mu_b \alpha_d -\lambda_d -\lambda_b)x +(-\mu_d \lambda_b + \mu_b \lambda_d -\mu_b \mu_d \alpha_d) = 0.
\label{zeros}
\end{equation}
Notice that the last term in the left-hand side of (\ref{zeros}) is negative if the stability condition
$\lambda_d \E [D] < \alpha_d + \lambda_b \E[B]$ holds,
that is, if $\mu_b \lambda_d < \mu_d \lambda_b + \mu_b \mu_d \alpha_d$,
thus guaranteeing that the product of the two roots $x_+$ and $x_-$ is negative,
and hence that there is a positive and a negative root.
One should subsequently observe that
$R_-$ must be zero to have a probability density.
Hence $f(v)$ is simply (a constant times) an exponential;
similarly for $g(v)$.
In addition, the steady-state amounts of demand and of blood have an atom at $0$ (since $\xi_d$ and  $\xi_b$ are no longer zero, the demand and blood processes
can reach $0$).

Interestingly, the model of this section is closely related to the model with workload removal
that is considered in \cite{Boucherie1996}. There an $M/G/1$ queue is studied with the extra feature
that, at Poisson epochs, a stochastic amount of work is removed.
In the $M/M/1$ case with removal of exponential amounts of work, see \cite[Sec.~5.1]{Boucherie1996}, one has the model of the present section
when we concentrate on the amount of demand present. 
One difference with the model in \cite{Boucherie1996}
is that, when the workload in that model has become zero, the work becomes positive at rate $\lambda_d$,
whereas in the present model the amount of blood can become positive (so zero demand is present)
and the amount of demand does not have to become positive when demands arrive (because they are immediately satisfied, see Figure~\ref{FIG2}). So the atom at zero is in the present model larger than in the model of \cite{Boucherie1996}.
In our model a positive demand level may be reached from below zero (by a jump, i.e., a demand arriving at an epoch
that there is some, but not enough, blood present). The memoryless property of the exponential
demand requirement distribution implies that this jump results in a demand level that is exp($\mu_d$),
just as if the initial demand level had been zero.
In the case of non-exponential demand requirements, our model becomes equivalent with an $M/G/1$ queue with
exponential amounts of work removed, and with the special feature that the first service requirement of a busy period has a different distribution.
Lemmas 4.1 and 4.2 of \cite{Boucherie1996} present the stability condition
of that $M/G/1$ queue with work removal; it amounts to 
$\lambda_d \E [D] < \alpha_d + \lambda_b \E [B]$, which indeed is one of the two stability conditions
of the present demand/blood model.

Finally we observe that Equation (5.1) of \cite{Boucherie1996} coincides with \eqref{zeros}
(take $\alpha_d=1$, $\lambda_d = \lambda_+$, $\lambda_b = \lambda_-$, $\mu_d=1/\beta$ and $\mu_b = 1/\gamma$).

\section{Asymptotic analysis}
\label{sectionscaling}
We finally study the model with $\alpha_b = \alpha_d = 0$ from an asymptotic perspective, by obtaining the fluid and diffusion limits of the blood inventory process. That is, we will create a sequence of processes, indexed by $n=1,2,...$, in which we let the rates of blood and demand arrivals grow large. If we then scale the process in a proper manner, we are able to deduce a non-degenerate limiting process, that provides insight in the overall behavior of the arrival volume when the system grows large, which only relies on the first two moments of the blood and demand distributions.

\subsection{Identification of the limiting process}
First, we introduce some additional notation. Let $X_b(t)$ and $X_d(t)$ denote the amount of blood and demand, respectively, at time $t>0$. Let 
\begin{equation}
X(t) := X_b(t) - X_d(t),
\end{equation}
be the net amount of blood available at time $t$. Remember that $X_b(t), X_d(t)\geq 0$, and $X_b(t)>0$ \emph{or} $X_d(t)>0$ for all $t$,
since $\alpha_d=\alpha_b=0$. Let $N_b(t)$, $N_d(t)$ be the two independent Poisson processes counting the number of arrivals of blood and demand, respectively. Then the following integral representation holds for $X(t)$,
\begin{equation}\label{eq:integralRep}
X(t) = X(0) - \xi_b\int_0^t X_b(s)\, ds + \xi_d \int_0^t X_d(s)\, ds + \sum_{i=1}^{N_b(t)} B_i - \sum_{i=1}^{N_d(t)} D_i.
\end{equation}
For the sake of exhibition, we will concentrate on the case $\xi_b=\xi_d =: \xi$. 
Our analysis may be extended to the general case.
A sketch of this generalization is given at the end of this section without going into the technical difficulties that arise when rigorously proving these limits. 

Define 
\begin{equation}\label{eq:defX}
Z(t) = \sum_{i=1}^{N_b(t)} B_i - \sum_{i=1}^{N_d(t)} D_i,
\end{equation}
that is, the difference between two compound Poisson processes, so that \eqref{eq:integralRep} reduces to 
\begin{equation}\label{eq:simpleRep}
X(t) = X(0) - \xi\int_0^t X(s)\, \dd s + Z(t).
\end{equation}
The first step in the definition of the sequence of processes under investigation is defining the asymptotic scheme we are interested in. As mentioned above, we intend to let the arrival rates grow to infinity. Therefore, in the $n^{th}$ process $X_n(t)$, we replace the rates of the arrival processes by $n\lambda_b$ and $n\lambda_d$. This induces Poisson processes $N^{(n)}_b(t)$ and $N^{(n)}_d(t)$ with arrival rates $n\lambda_b$ and $n\lambda_d$, respectively. 
However, we have 
\begin{equation}
N_b^{(n)}(t) \equalD N_b(n t)\qquad \text{and}\qquad N_d^{(n)}(t) \equalD N_d(n t),
\end{equation}
so that the term $Z(t)$ in \eqref{eq:simpleRep} in this asymptotic scheme can be replaced by
\begin{equation}
Z_n(t) = \sum_{i=1}^{N_b(nt)} B_i - \sum_{i=1}^{N_d(nt)} D_i.
\end{equation}
The first step in our analysis is obtaining the fluid limit of the process. Bearing in mind application of the Functional Law of Large Numbers (FLLN), we scale the process as $\bar{X}_n(t) = X_n(t)/n$,
so that with \eqref{eq:simpleRep}
\begin{equation}\label{eq:fluidRep}
\bar{X}_n(t) = \bar{X}_n(0) -\xi\int_0^t \bar{X}_n(s)\, ds + \bar{Z}_n(t),
\end{equation}
where $\bar{Z}_n(t) = Z_n(t)/n$.
%

The essential step in establishing a result on the convergence of $\bar{X}_n$ is the application of \cite[Thm.~4.1]{Pang2007}, which we cite here for completeness, slightly rewritten to fit our setting.

\begin{theorem}[{\cite[Thm.~4.1]{Pang2007}}]
\label{thm:pang}Let $D[0,\infty)$ be the space of all one-dimensional real-valued c\`adl\`ag functions defined on $[0,\infty)$, endowed with the usual $J_1$-Skorohod topology.
Consider the integral representation 
\begin{equation}\label{eq:pangInt}
x(t) = y(t) + \int_0^t u(x(s))\, ds, \qquad t \geq 0,
\end{equation}
where $u:\mathbb{R}\to\mathbb{R}$ satisfies $u(0)=0$ and is Lipschitz continuous. 
The integral representation in \eqref{eq:pangInt} has a unique solution $x$, so that the integral representation constitutes a function $H_u: D[0,\infty) \to D[0,\infty)$ mapping $y$ into $x\equiv H_u(y)$. 
In addition, the function $H_u$ is continuous, and if $y$ is continuous, then so is $x$. 
\end{theorem}

In our case, we set $u(x) = -\xi x$, to be able to write $\bar{X}_n = H_u\left(\bar{X}_n(0) +  \bar{Z}_n\right)$. Since $u$ is clearly Lipschitz continuous, the mapping $H_u$ is indeed continuous. 
Let us rewrite \eqref{eq:fluidRep}, by observing 
\begin{equation}
\E \bar{Z}_n(t) = \frac{1}{n}\Big(\E [N_b(nt)]\E [B] - \E [N_d(n t)]\E[D]\Big) = \lambda_b\E [B] t - \lambda_d \E [D] t, 
\end{equation}
where the expectation is taken with respect to the compound Poisson processes. 
Since $m=\lambda_b \E[B] - \lambda_d \E[D]$,
\begin{equation}\label{eq:fluidRep2}
\bar{X}_n(t) = \bar{X}_n(0)  -\xi\int_0^t\left(\bar{X}_n(s)-\frac{m}{\xi}\right)\, \dd s + \bar{Y}_n(t),
\end{equation}
where $\bar{Y}_n(t) := \bar{Z}_n(t) - mt$ is now a centered process. 
This allows us to state the next result.

\begin{proposition}[Fluid limit]\label{fluidProp}
Let $\E[B],\, \E[D] <\infty$ and $\bar{X}_n(0) = X_n(0)/n \rightarrow q_0 \in \mathbb{R}$, as $n\rightarrow \infty$. Then for $n\rightarrow\infty$,
\begin{equation}\label{eq:fluidLimit}
\bar{X}_n \Rightarrowd q,
\end{equation}
where
\begin{equation}
q(t) = \frac{m}{\xi} + \left(q_0 - \frac{m}{\xi}\right) \ee^{-\xi t}.
\label{eq:fluid_function}
\end{equation}
\end{proposition}
\begin{proof}
First, we concentrate on the process $\bar{Y}_n$. Observe that, by the FLLN for renewal-reward processes, which follows from \cite[Thm.~7.4.1]{Whitt2002}, we have
\begin{equation}
\frac{1}{nt}\sum_{i=1}^{N_b(nt)} B_i \Rightarrowd \lambda_b\E[B],\qquad \frac{1}{nt}\sum_{i=1}^{N_d(nt)} D_i \Rightarrowd  \lambda_d\E[D],
\end{equation}
for $n\rightarrow\infty$ and for all $t>0$. Hence, $\bar{Z}_n(t) \Rightarrowd \lambda_b\E[B]t-\lambda_d\E[D]t = mt$. By definition of $\bar{Y}_n$ and the assumption of convergence of $\bar{X}_n(0)$, this implies 
\begin{equation}
\bar{Y}_n + \bar{X}_n \Rightarrowd q_0
\end{equation}
as $n\rightarrow\infty$.
Next, note $\bar{X}_n = H_u\left(\bar{X}_n(0) + \bar{Z}_n \right) = H_u\left(\bar{X}_n(0) + \bar{Y}_n +I t \right)$, where $I$ denotes the identity map, i.e.~ $I(t) = t$ for all $t \geq 0$. Due to Lipschitz continuity of $u$, $H_u$ constitutes a continuous mapping, and hence we can apply the Continuous Mapping Theorem (CMT), to find 
\begin{equation}
\bar{X}_n = H_u\left(\bar{X}_n(0) + \bar{Y}_n +m I\right) \Rightarrow H_u\left(q_0+m I\right)\equiv q,
\end{equation}
for all $t\geq 0$, where $q(\cdot)$ is the solution of 
\begin{align*}
q(t) &= q_0 + \int_0^t u(q(s))\, \dd s = q_0 + mt - \xi \int_0^t q(s) \, \dd s\\
&= q_0  - \xi \int_0^t \left(q(s) - \frac{m}{\xi}\right) \, \dd s.
\end{align*}
The unique solution of this integral equation is given in \eqref{eq:fluid_function}.
\end{proof}

According to Proposition~\ref{fluidProp}, the fluid limit approaches $\E[X] = m/\xi$ exponentially fast.
To obtain an expression for the {\em diffusion limit} of the process, we analyze the fluctuations of the process around the fluid limit in (\ref{eq:fluidLimit}), again by scaling the process in a proper manner. First, we subtract $q(t)$ on both sides of \eqref{eq:fluidRep2}, and multiply by $\sqrt{n}$:
\begin{equation}\label{eq:hatEquation}
\sqrt{n}\left( \bar{X}_n(t) - q(t) \right) = \sqrt{n}\left( \bar{X}_n(0) - q_0 \right) -\xi \int_ 0^t \sqrt{n}\left( \bar{X}_n(s) - q(s) \right)\, ds + \sqrt{n}\,\bar{Y}_n(t).
\end{equation}
Let $\hat{X}_n \equiv \sqrt{n}\left( \bar{X}_n - q \right)$ and $\hat{Y}_n \equiv \sqrt{n}\,\bar{Y}_n$, then this reduces to 
\begin{equation}
\hat{X}_n(t) = \hat{X}_n(0) - \xi \int_0^t \hat{X}_n(s) \,ds + \hat{Y}_n(t).
\end{equation}
Again the process $\hat{Y}_n$ needs special attention. 

\begin{lemma}\label{diffLemma}
Let $\E[B^2], \E[D^2] < \infty$. Then $\hat{Y}_n \Rightarrowd \sigma W$ as $n\rightarrow \infty$, where $\sigma^2 := \lambda_b \E[B^2]+\lambda_d\E[D^2]$ and $W$ is a standard Brownian motion.
\end{lemma}
\begin{proof}
Recall that 
\begin{equation}
\hat{Y}_n(t) \equalD \sqrt{n}\left[ \Big(\frac{1}{n}\sum_{i=1}^{N_b(n t)} B_i -\lambda_b\E [B]t \Big) - 
 \Big(\frac{1}{n}\sum_{i=1}^{N_d(nt)} D_i - \lambda_d \E[D] t \Big)
 \right] .
\end{equation}
By the Functional Central Limit Theorem (FCLT) for renewal-reward processes given in \cite[Thm.~7.4.1]{Whitt2002}, the process
\begin{equation}
\hat{Y}_n^b(t) = \sqrt{n} \Big(\frac{1}{n}\sum_{i=1}^{N_b(n t)} B_i -\lambda_b\E[ B] t\Big),
\end{equation}
converges weakly to $\sigma_b W_b$, where $W_b$ is a standard Brownian motion, and 
\begin{equation}
\sigma_b^2 = \lambda_b\,\text{Var}\,B + \lambda_b(E[B])^2 = \lambda_b\E[B^2].
\end{equation}
Similarly, $\hat{Y}_n^d \Rightarrow \sigma_d W_d$ as $n\to\infty$, with the obvious parameter switches and $W_d$ is standard Brownian motion. Since the processes $\hat{Y}_n^b$ and $\hat{Y}_n^d$ are independent, so are their limits, and 
\begin{equation}
\hat{Y}_n \Rightarrow \sqrt{\lambda_b \E[B^2]}\, W_b + \sqrt{\lambda_d \E[D^2]}\, W_d \equalD \sqrt{\lambda_b\E[B^2]+\lambda_d\E[D^2]}\,W,
\end{equation}
for $n\rightarrow\infty$ and $W$ a standard Brownian motion.
\end{proof}
Now, we are ready to prove the diffusion counterpart of Proposition \ref{fluidProp}. 

\begin{proposition}[Diffusion limit]\label{diffProp}
Let $\E[B^2], \E[D^2] < \infty$. If $\hat{X}_n(0) \rightarrow \hat{X}(0)$, then $\hat{X}_n \Rightarrow \hat{X}$ as $n\rightarrow \infty$, where $\hat{X}$ satisfies the integral equation 
\begin{equation}\label{diffLimit}
\hat{X}(t) = \hat{X}(0) - \xi \int_0^t \hat{X}(s) \, \dd s + \sigma W(t).
\end{equation}
In other words, $\hat{X}$ is an Ornstein-Uhlenbeck diffusion process with infinitesimal mean $\xi$ and infinitesimal variance $\sigma^2 := \lambda_b\E[B^2] + \lambda_d\E[D^2]$.
\end{proposition}
\begin{proof}
We again rely on the result that the mapping $H_u$ as in the proof of Proposition \ref{fluidProp} is continuous if $u$ is Lipschitz continuous. Here, we set $u(x) = -\xi x$ which again clearly satisfies this condition. We have $\hat{X}_n \equiv H_u(\hat{X}_n(0)+ \hat{Y}_n)$. From Lemma \ref{diffLemma}, we know
\begin{equation}
\hat{X}_n(0) + \hat{Y}_n \Rightarrow \hat{X}(0) + \sigma W,
\end{equation}
for $n\rightarrow\infty$. As a consequence of the CMT, we conclude 
\begin{equation}
\hat{X}_n = H_u\left( \hat{X}_n(0) + \hat{Y}_n\right) \Rightarrow H_u\left(\hat{X}(0) + \sigma W\right) \equiv \hat{X},
\end{equation}
where $\hat{X}$ solves \eqref{diffLimit}.
\end{proof}

\subsection{Generalization for $\xi_b\neq \xi_d$}
We now sketch the scaling approach towards fluid and diffusion limits for the general case in which
$\xi_b$ may differ from $\xi_d$.
In case $\xi_b \neq \xi_d$, the integral equation for $\bar{X}_n$ as in \eqref{eq:fluidRep} becomes
\begin{align}
\label{eq:fluidRepNEQ}
\bar{X}_n(t) &= \bar{X}_n(0) + \int_0^t ( -\xi_b \bar{X}_n^+(s) + \xi_d \bar{X}_n^-(s) - m)\, \dd s + \bar{Y}_n(t)\\
&= \bar{X}_n(0) - \int_0^t  (\left[\xi_b \mathbbm{1}_{\{ \bar{X}_n(s)\geq 0\}}+\xi_d \mathbbm{1}_{\{ \bar{X}_n(s)<0\}}\right] \bar{X}_n(s) + m) \, \dd s + \bar{Y}_n(t),
\nonumber
\end{align}
where $\bar{Y}_n(t)$ is defined as before. Note that $\hat{X}_n \equiv H_u(\bar{X}_n(0)+\bar{Y}_n)$, where we now have
\begin{equation}
u(x) = - \left[\xi_b \mathbbm{1}_{\{ x\geq 0\}}+\xi_d \mathbbm{1}_{\{x<0\}}\right]x + m,
\end{equation}
which is still Lipschitz continuous. Following the same reasoning of the proof of Proposition \ref{fluidProp}, we obtain the fluid limit $\bar{X}_n \Rightarrowd q$, where $q$ is the solution of
\begin{equation}
q(t) = q_0 - \int_0^t (\left[\xi_b \mathbbm{1}_{\{ q(s)\geq 0\}}+\xi_d \mathbbm{1}_{\{q(s)<0\}}\right]q(s) - m)\, \dd s .
\end{equation}
The solution to this integral equation is more elaborate than \eqref{eq:fluidLimit} and depends on the sign of $m$ and $q_0$. Assuming $m\geq 0$, one can check that, 
\begin{align}
q(t) &= \frac{m}{\xi_b} + \left(q_0- \frac{m}{\xi_b}\right) e^{-\xi_b t}, & \text{if }q_0\geq 0,\\
q(t) &= \left\{
\begin{array}{ll}
\frac{m}{\xi_d} + \left(q_0- \frac{m}{\xi_d}\right) e^{-\xi_d t}, & \text{if } 0\leq t < t_d^*,\\
\frac{m}{\xi_b}\left(1-e^{-\xi_b (t-t^*_d)}\right), & \text{if } t \geq t_d^*,
\end{array}\right. & \text{if } q_0 <0,
\label{eq:fluid_1}
\end{align}
where \begin{equation}
t_d^* = - \frac{1}{\xi_d}\,\log\left(\frac{m/\xi_d}{m/\xi_d-q_0}\right).
\end{equation}
If $m < 0$,
\begin{align}
q(t) &= \frac{m}{\xi_d} + \left(q_0- \frac{m}{\xi_d}\right) e^{-\xi_d t}, & \text{if }q_0\leq 0,\\
q(t) &= \left\{
\begin{array}{ll}
\frac{m}{\xi_b} + \left(q_0- \frac{m}{\xi_b}\right) e^{-\xi_b t}, & \text{if } 0\leq t < t_b^*,\\
\frac{m}{\xi_d}\left(1-e^{-\xi_d (t-t^*_b)}\right), & \text{if } t \geq t_b^*,
\end{array}\right. & \text{if } q_0>0,
\label{eq:fluid_2}
\end{align}
where \begin{equation}
t_b^* = - \frac{1}{\xi_b}\,\log\left(\frac{m/\xi_b}{m/\xi_b-q_0}\right).
\end{equation}
Note that the equilibrium of the fluid limit also depends on the sign of $m$:
\begin{equation}
\lim_{t\rightarrow\infty} q(t) = \left\{\begin{array}{ll}
m/\xi_b, & \text{if }m\geq 0,\\
m/\xi_d, & \text{if }m <0.
\end{array}\right.
\end{equation}
In the remainder, without loss of generality $m\geq 0$. Furthermore, set $q_0= m/\xi_b$ so that $q \equiv m/\xi_b$. Subtracting $q(t)$ on both sides of \eqref{eq:fluidRepNEQ} yields,
\begin{align}
\left(\bar{X}_n(t)-q(t)\right) &= \left(\bar{X}_n(0)-q_0\right) - 
\int_0^t  \Big\{ \left[\xi_b \mathbbm{1}_{\{\bar{X}_n(s)\geq 0\}}+\xi_d \mathbbm{1}_{\{\bar{X}_n(s)<0\}}\right] \bar{X}_n(s) \nonumber \\
& \qquad \qquad - \xi_b\,q(s)\Big\} \, \dd s + \bar{Y}_n(t)\\
&= \left(\bar{X}_n(0)-q_0\right) - \int_0^t \xi_b \left(\bar{X}_n(s)-q(s)\right)\, \dd s \nonumber \\
&\qquad \qquad + 
\int_0^t \mathbbm{1}_{\{\bar{X}_n(s) < 0\}}(\xi_b-\xi_d)\bar{X}_n(s)\, \dd s.
\end{align}
Let $\hat{X}_n(t) = \sqrt{n}\left(\bar{X}_n(t)-q(t)\right)$. Then 
\begin{equation}
\hat{X}_n(t) = \hat{X}_n(0) - \xi_b \int_0^t \hat{X}_n(s)\,\dd s +\int_0^t \mathbbm{1}_{\{\bar{X}_n(s) < 0\}}(\xi_b-\xi_d)\bar{X}_n(s)\, \dd s + \hat{Y}_n(t)
\end{equation}
Now, we argue non-rigorously that the one-but-last term vanishes as $n\rightarrow\infty$. Namely, by defining the function $G: D[0,\infty)\rightarrow D[0,\infty)$ by the integration operator:
\begin{equation}
G(u) = \int_0^t \mathbbm{1}_{\{u(s)<0\}} (\xi_b-\xi_d) u(s)\, \dd s,
\end{equation}
this term can be expressed as $G(\bar{X}_n)$. Hence by the fact that $\hat{X}_n\Rightarrowd m/\xi_b$ and the CMT we see $G(\hat{X}_n)\Rightarrow 0$.

Under this claim, we deduce by the approach of Proposition \ref{diffProp}, that if $\hat{X}_n \Rightarrow \hat{X}$ for $n\rightarrow \infty$, then $\hat{X}$ satisfies the stochastic integral equation 
\begin{equation}
\hat{X}(t) = \hat{X}(0) - \xi_b \int_0^t \hat{X}(s)\, \dd s + \sigma W(t),
\end{equation}
which implies that $\hat{X}$ is an Ornstein-Uhlenbeck process with infinitesimal mean $\xi_b$ and variance $\sigma^2 := \lambda_b\E[B^2] + \lambda_d\E[D^2]$.
\\
The result that the scaled process converges to an Ornstein-Uhlenbeck process can be intuitively justified by the so-called \textit{mean-reverting} behavior of the original process. That is, the further the process is away from its mean, the greater the drift towards that equilibrium. This is the defining feature of the OU diffusion process. The decay rates $\xi_b$ and $\xi_d$ are responsible for the original process being `forced' towards 0 and therefore the similarities should not be surprising. However, note that in the diffusion limit $X_n$ has drift $\xi_b$ (cq. $\xi_d$) towards $nm/\xi_b$ (cq. $nm/\xi_d$), if $m>0$ (cq. $<0$) at \emph{any} position of the process. This implies that if $X_n \in(0,nm/\xi_b)$, it has an upward drift equal to $\xi_b$, which is at first sight counter-intuitive.
However, we can argue that in case $X_n(t) = v \in (0, nm/\xi_b)$, the mean upward drift of the process $X_n$ equals $n\lambda_b\E[B]$, and the mean downward drift equals $n\lambda_d\E[D] + \xi_b v$, since $v>0$.
Rewrite $v = nm/\xi_b - w\sqrt{n}$ for some $w \in (0,\sqrt{n} m/\xi_b)$. 
Then, the mean net drift equals 
\[
n\lambda_b\E[B] - n\lambda_d\E[D] - \xi_b \left( \frac{n m}{\xi_b} - w\sqrt{n} \right) = \xi_b w \sqrt{n} >0,
\]
which explains both the sign and magnitude of the drift factor in the scaled process.

\subsection{Related literature}
The Ornstein-Uhlenbeck process is a diffusion process that often arises as the limit of a sequence of stochastic systems, in which the system size tends to infinity. 
Particularly in queueing settings with mean reverting behavior, the OU process appears in so-called heavy traffic, i.e.~the arrival rate grows without bound. 
We mention a couple of models that exhibit limiting behavior that is similar to ours.

First, it is well-known that the properly normalized $M/M/\infty$ queue length process converges weakly to a OU process as the arrival rate tends to infinity, see e.g.~\cite[Sec.~10.3]{Whitt2002}. 
This limiting behavior continues to hold in case the queueing process is modulated by a Markovian background process, see \cite{Anderson2016}.

Another well-known queueing model in which a (piecewise) OU process appears in the limit is the multi-server queue with abandonments.
For the $M/M/s+M$ queue, where $+M$ denotes the exponentially distributed patience of customers, Garnett et al.~\cite{Garnett2002} showed that in the Halfin-Whitt regime, the queue length process, centered and scaled around the number of servers $s$, approaches a hybrid OU process, of which the drift parameter depends on the current state: If the queue length is larger (cq.~smaller) than zero, then the drift is governed by the abandonment rate (cq.~service rate). 
Dai et al.~\cite{Dai2010} find a similar piecewise diffusion process under more general assumptions on the model primitives. 

For the single-server queue with abandoning customers, Ward \& Glynn \cite{Ward2003,Ward2005} showed that in conventional heavy traffic, the queue length process converges to a OU process with reflecting barrier 0.  
 
Since we in our setting assumed both demand impatience and perishability of inventory (which can be seen as a kind of impatience as well), it should not come as a surprise that we also find our limiting process to be a OU process. 
Observe however that in our model, unless $m=0$, we find a OU process with constant, rather than piecewise, parameters, and no reflection barrier, since our (scaled) inventory process can go both positive and negative. 

Last, we mention that there is a connection between our blood inventory process and the work of Reed \& Zwart \cite{Reed2011}.
Rather than looking at the OU process as the limit of a sequence of stochastic processes, Reed and Zwart \cite{Reed2011} study a stochastic differential equation that is closely related to Equation \eqref{eq:integralRep}, in the sense that the process has a different (constant) drift term in the upper and lower half plane.
Under the assumption that the input process is a L\'evy process with only one-sided jumps, they develop a methodology to derive the invariant distribution of the solution of the SDE. 
Unfortunately, the input in our scenario exhibits both positive and negative jumps, which prevents us from applying their results directly to \eqref{eq:integralRep}. 

\section{Numerical evaluation}
\label{numericals}

\subsection{Approximation scheme}
The asymptotic results of the previous section regarding the fluid and diffusion limits allude to the fact that for large arrival rates, the normalized inventory process $\{\hat X_n(t)\, |\, t\geq 0\}$, resembles that of the Ornstein-Uhlenbeck process. 
Indeed, the sample paths of the scaled process $\bar X_n$ for increasing values of $n$ in Figures \ref{fig:sample_paths_fluid1} and \ref{fig:sample_paths_fluid2} show that the mean-reverting behavior around $m/\xi^*$, that is typical of OU processes, kicks in rather quickly. 
Moreover, the fluid limits $q(t)$ as presented by Proposition \ref{fluidProp} and \eqref{eq:fluid_1}-\eqref{eq:fluid_2} predict the mean well for both $\xi_b=\xi_d$ and $\xi_b \neq \xi_d$. 
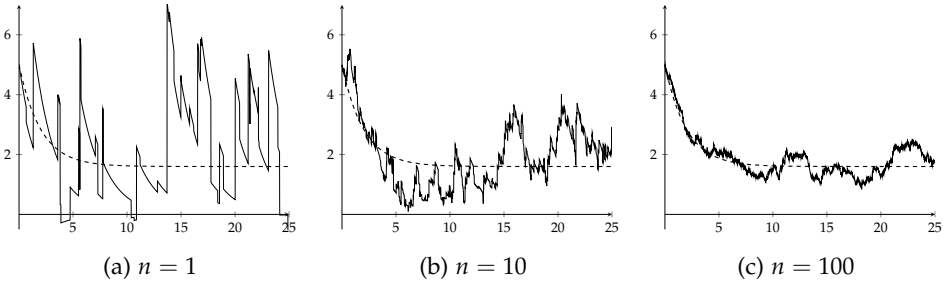
\begin{figure}
\begin{subfigure}{0.32\textwidth}
\centering
\begin{tikzpicture}[scale=0.52]
\begin{axis}[
	xmin = 0,
	xmax = 25,
	ymin = -0.5,
	ymax = 7,
	axis line style={->},
	axis lines = middle]
	
\addplot[] table[x=t,y=sp1] {Chapter_7/tikz/sample_paths_fluid.txt};
\addplot[dashed,thick] table[x=t,y=fluid] {Chapter_7/tikz/sample_paths_fluid.txt};
\end{axis}
\end{tikzpicture}
\caption{$n=1$}
\end{subfigure}
\begin{subfigure}{0.32\textwidth}
\centering
\begin{tikzpicture}[scale=0.52]
\begin{axis}[
	xmin = 0,
	xmax = 25,
	ymin = -0.5,
	ymax = 7,
	axis line style={->},
	axis lines = middle]
	
\addplot[] table[x=t,y=sp10] {Chapter_7/tikz/sample_paths_fluid.txt};
\addplot[dashed,thick] table[x=t,y=fluid] {Chapter_7/tikz/sample_paths_fluid.txt};
\end{axis}
\end{tikzpicture}
\caption{$n=10$}
\end{subfigure}
\begin{subfigure}{0.32\textwidth}
\centering
\begin{tikzpicture}[scale=0.52]
\begin{axis}[
	xmin = 0,
	xmax = 25,
	ymin = -0.5,
	ymax = 7,
	axis line style={->},
	axis lines = middle]
	
\addplot[] table[x=t,y=sp100] {Chapter_7/tikz/sample_paths_fluid.txt};
\addplot[dashed,thick] table[x=t,y=fluid] {Chapter_7/tikz/sample_paths_fluid.txt};
\end{axis}
\end{tikzpicture}
\caption{$n=100$}
\end{subfigure}
\caption{Sample paths of the process $\bar{X}_n(t) = X_n(t)/n$ with $\bar{X}_n(0) = 5$, $\lambda_b = 1.2$, $\lambda_d = 1$, $\xi_b = \xi_d = 0.5$ and $\mu_b=0.5$ and $\mu_d=1$. The fluid limit is depicted by the dashed line.} 
\label{fig:sample_paths_fluid1}
\end{figure}
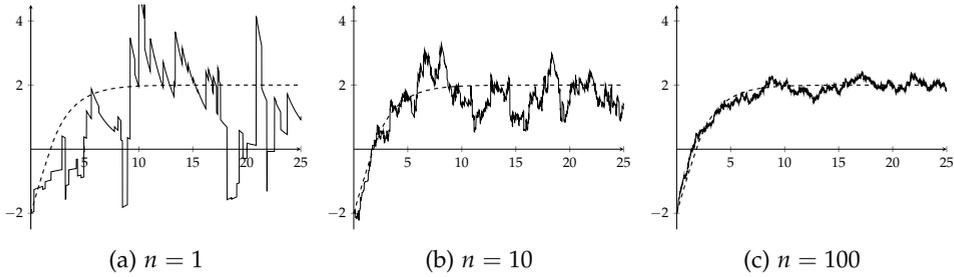
\begin{figure}
\begin{subfigure}{0.32\textwidth}
\centering
\begin{tikzpicture}[scale=0.52]
\begin{axis}[
	xmin = 0,
	xmax = 25,
	ymin = -2.5,
	ymax = 4.5,
	axis line style={->},
	axis lines = middle]
	
\addplot[] table[x=t,y=sp1] {Chapter_7/tikz/sample_paths_fluid_min.txt};
\addplot[dashed,thick] table[x=t,y=fluid] {Chapter_7/tikz/sample_paths_fluid_min.txt};
\end{axis}
\end{tikzpicture}
\caption{$n=1$}
\end{subfigure}
\begin{subfigure}{0.32\textwidth}
\centering
\begin{tikzpicture}[scale=0.52]
\begin{axis}[
	xmin = 0,
	xmax = 25,
	ymin = -2.5,
	ymax = 4.5,
	axis line style={->},
	axis lines = middle]
	
\addplot[] table[x=t,y=sp10] {Chapter_7/tikz/sample_paths_fluid_min.txt};
\addplot[dashed,thick] table[x=t,y=fluid] {Chapter_7/tikz/sample_paths_fluid_min.txt};
\end{axis}
\end{tikzpicture}
\caption{$n=10$}
\end{subfigure}
\begin{subfigure}{0.32\textwidth}
\centering
\begin{tikzpicture}[scale=0.52]
\begin{axis}[
	xmin = 0,
	xmax = 25,
	ymin = -2.5,
	ymax = 4.5,
	axis line style={->},
	axis lines = middle]
	
\addplot[] table[x=t,y=sp100] {Chapter_7/tikz/sample_paths_fluid_min.txt};
\addplot[dashed,thick] table[x=t,y=fluid] {Chapter_7/tikz/sample_paths_fluid_min.txt};
\end{axis}
\end{tikzpicture}
\caption{$n=100$}
\end{subfigure}
\caption{Sample paths of the process $\bar{X}_n(t) = X_n(t)/n$ with $\bar{X}_n(0) = -2$, $\lambda_b = 2$, $\lambda_d = 1$, $\xi_b = 0.5, \xi_d = 0.1$ and $\mu_b=1$ and $\mu_d=1$. The fluid limit is depicted by the dashed line.} 
\label{fig:sample_paths_fluid2}
\end{figure}
Furthermore, we observe that steady state is attained fairly quickly.
This is suggestive of the claim that the steady-state distribution of the normalized process $\hat{X}_n$ is well-described by the steady-state distribution of the OU process $\hat{X}$. 
Since the OU process with mean 0, infinitesimal variance $\sigma^2$ and drift $\xi^*$ is known to be normally distributed with mean 0 and variance $\sigma^2/2\xi^*$ in steady-state, this leads to a simpler approximation scheme based on the first two moments of $B$ and $D$ only. 
In non-rigorous mathematical terms, we use the approximation that
\begin{equation}
\hat{X}_n = \frac{ X_n - nm/\xi^*}{\sqrt{n}} {\;\buildrel{d}\over\approx \;} Z^*,
\label{eq:normal_approximation}
\end{equation}
where $Z^*$ is a normally distributed random variable with mean 0 and variance $\sigma^2/2\xi^*$. 

Note that justification of the conjecture that the normal approximation is indeed an asymptotically correct approximation for systems with large arrival rates requires proof that the interchange-of-limits between $t\to\infty$ and $n\to\infty$ is indeed valid. 
Rather than going into the technical details, we provide in the remainder of this section numerical evidence that this interchange indeed holds, and that the normal approximation is able to capture characteristics of processes with exponential jumps as well as generally distributed jumps. 

\subsection{Distribution functions}

Since we obtained an explicit expression for the steady-state density function of the net inventory process $X$ in case $B$ and $D$ are exponential, see Theorem \ref{thm:full_pdf}, we will exploit this formula for numerical comparison to the normal approximation arising from the OU process.

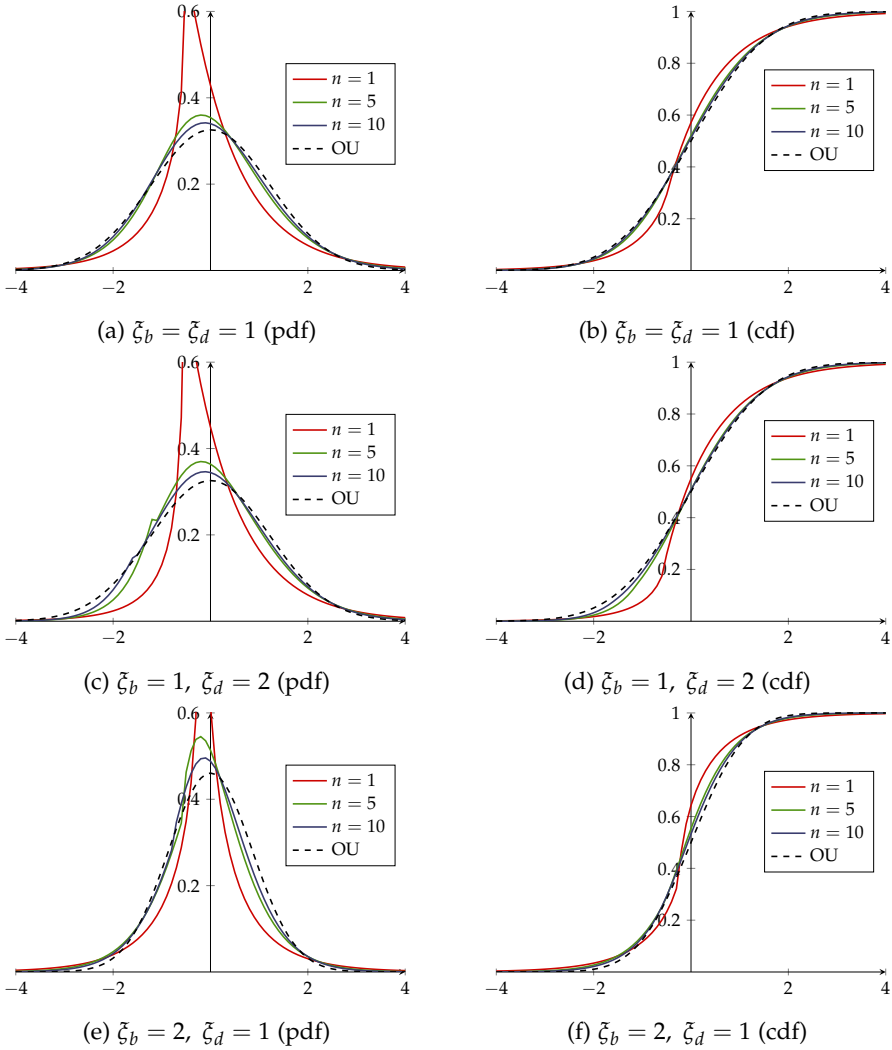
\begin{figure}
\centering
\begin{subfigure}{0.48\textwidth}
\centering
\begin{tikzpicture}[scale=0.75]
\small
\begin{axis}[
	xmin = -4,
	xmax = 4, 
	ymin = 0,
	ymax = 0.6,
	axis line style={->},
	axis lines = middle,
	legend style = {at = {(axis cs: 3.8,0.48)},anchor=north east},
	legend cell align = left,
	yscale = 0.8
	]
	
\addplot[thick,col1] table[x=v,y=h1] {Chapter_7/tikz/convergence_1_1.txt};
\addplot[thick,col4] table[x=v,y=h5] {Chapter_7/tikz/convergence_1_1.txt};
\addplot[thick,col5] table[x=v,y=h10] {Chapter_7/tikz/convergence_1_1.txt};
\addplot[thick,dashed] table[x=v,y=ou] {Chapter_7/tikz/convergence_1_1.txt};

\legend{$n=1$,$n=5$,$n=10$,OU};
\end{axis}
\end{tikzpicture}
\caption{$\xi_b=\xi_d=1$ (pdf)}
\end{subfigure}
\begin{subfigure}{0.48\textwidth}
\centering
\begin{tikzpicture}[scale=0.75]
\begin{axis}[
	xmin = -4,
	xmax = 4, 
	ymin = 0,
	ymax = 1,
	axis line style={->},
	axis lines = middle,
	legend style = {at = {(axis cs: 3.8,0.04)},anchor=south east},
	legend cell align = left,
	legend pos = north east,
	yscale = 0.8
	]
\small	
\addplot[thick,col1] table[x=v,y=H1] {Chapter_7/tikz/convergence_1_1.txt};
\addplot[thick,col4] table[x=v,y=H5] {Chapter_7/tikz/convergence_1_1.txt};
\addplot[thick,col5] table[x=v,y=H10] {Chapter_7/tikz/convergence_1_1.txt};
\addplot[thick,dashed] table[x=v,y=OU] {Chapter_7/tikz/convergence_1_1.txt};

\legend{$n=1$,$n=5$,$n=10$,OU};

\end{axis}
\end{tikzpicture}
\caption{$\xi_b=\xi_d=1$ (cdf)}
\end{subfigure}

\begin{subfigure}{0.48\textwidth}
\centering
\begin{tikzpicture}[scale=0.75]
\small
\begin{axis}[
	xmin = -4,
	xmax = 4, 
	ymin = 0,
	ymax = 0.6,
	axis line style={->},
	axis lines = middle,
	legend style = {at = {(axis cs: 3.8,0.48)},anchor=north east},
	legend cell align = left,
	yscale = 0.8
	]
	
\addplot[thick,col1] table[x=v,y=h1] {Chapter_7/tikz/convergence_1_2.txt};
\addplot[thick,col4] table[x=v,y=h5] {Chapter_7/tikz/convergence_1_2.txt};
\addplot[thick,col5] table[x=v,y=h10] {Chapter_7/tikz/convergence_1_2.txt};
\addplot[thick,dashed] table[x=v,y=ou] {Chapter_7/tikz/convergence_1_2.txt};

\legend{$n=1$,$n=5$,$n=10$,OU};
\end{axis}
\end{tikzpicture}
\caption{$\xi_b=1,\ \xi_d=2$ (pdf)}
\end{subfigure}
\begin{subfigure}{0.48\textwidth}
\centering
\begin{tikzpicture}[scale=0.75]
\begin{axis}[
	xmin = -4,
	xmax = 4, 
	ymin = 0,
	ymax = 1,
	axis line style={->},
	axis lines = middle,
	legend style = {at = {(axis cs: 3.8,0.04)},anchor=south east},
	legend cell align = left,
	legend pos = north east,
	yscale = 0.8
	]
\small	
\addplot[thick,col1] table[x=v,y=H1] {Chapter_7/tikz/convergence_1_2.txt};
\addplot[thick,col4] table[x=v,y=H5] {Chapter_7/tikz/convergence_1_2.txt};
\addplot[thick,col5] table[x=v,y=H10] {Chapter_7/tikz/convergence_1_2.txt};
\addplot[thick,dashed] table[x=v,y=OU] {Chapter_7/tikz/convergence_1_2.txt};

\legend{$n=1$,$n=5$,$n=10$,OU};

\end{axis}
\end{tikzpicture}
\caption{$\xi_b=1,\ \xi_d=2$ (cdf)}
\end{subfigure}

\begin{subfigure}{0.48\textwidth}
\centering
\begin{tikzpicture}[scale=0.75]
\small
\begin{axis}[
	xmin = -4,
	xmax = 4, 
	ymin = 0,
	ymax = 0.6,
	axis line style={->},
	axis lines = middle,
	legend style = {at = {(axis cs: 3.8,0.48)},anchor=north east},
	legend cell align = left,
	yscale = 0.8
	]
	
\addplot[thick,col1] table[x=v,y=h1] {Chapter_7/tikz/convergence_2_1.txt};
\addplot[thick,col4] table[x=v,y=h5] {Chapter_7/tikz/convergence_2_1.txt};
\addplot[thick,col5] table[x=v,y=h10] {Chapter_7/tikz/convergence_2_1.txt};
\addplot[thick,dashed] table[x=v,y=ou] {Chapter_7/tikz/convergence_2_1.txt};

\legend{$n=1$,$n=5$,$n=10$,OU};
\end{axis}
\end{tikzpicture}
\caption{$\xi_b=2,\ \xi_d=1$ (pdf)}
\end{subfigure}
\begin{subfigure}{0.48\textwidth}
\centering
\begin{tikzpicture}[scale=0.75]
\begin{axis}[
	xmin = -4,
	xmax = 4, 
	ymin = 0,
	ymax = 1,
	axis line style={->},
	axis lines = middle,
	legend style = {at = {(axis cs: 3.8,0.04)},anchor=south east},
	legend cell align = left,
	legend pos = north east,
	yscale = 0.8
	]
\small	
\addplot[thick,col1] table[x=v,y=H1] {Chapter_7/tikz/convergence_2_1.txt};
\addplot[thick,col4] table[x=v,y=H5] {Chapter_7/tikz/convergence_2_1.txt};
\addplot[thick,col5] table[x=v,y=H10] {Chapter_7/tikz/convergence_2_1.txt};
\addplot[thick,dashed] table[x=v,y=OU] {Chapter_7/tikz/convergence_2_1.txt};

\legend{$n=1$,$n=5$,$n=10$,OU};

\end{axis}
\end{tikzpicture}
\caption{$\xi_b=2,\ \xi_d=1$ (cdf)}
\end{subfigure}
\caption{Probability functions of $\hat{X}_n$ for $n=1,5$ and $10$ with $\lambda_b=1,\lambda_d=0.5,\mu_b=\mu_d=1$, and the probability function of the OU process.}
\label{fig:distributions}
\end{figure}

Let $h(\cdot)$ as in Theorem \ref{thm:full_pdf} be the pdf of $X$ with parameters $\lambda_b$, $\lambda_d$, $\mu_b$, $\mu_d$, $\xi_b$ and $\xi_d$, and the corresponding cdf $H$, defined as $H(v) = \int_{-\infty}^v h(x) \dd x$. 
We denote by $h_n$ and $H_n$ the pdf and cdf, respectively, of the inventory process $X_n$ with arrival rates $n\lambda_b$ and $n\lambda_d$, and the remaining parameters unchanged.
Then, the pdf and cdf of the normalized process are given by $\hat{h}_n(v) = \sqrt{n}\,h_n(v_n)$ and $\hat{H}_n(v) = H_n(v_n)$, respectively, with $v_n = nm/\xi^* + v\sqrt{n}$ for all $v\in\mathbb{R}$.
By the normal approximation scheme, we expect
\begin{equation}
\hat{h}_n(v) \approx \frac{\sqrt{2\xi^*}}{\sigma}\,\f\left(\frac{\sqrt{2\xi^*}}{\sigma}v\right),\quad \text{ and } \quad
\hat{H}_n(v) \approx \F\left(\frac{\sqrt{2\xi^*}}{\sigma}v\right).
\end{equation}
We perform this numerical comparison of probability functions in Figure \ref{fig:distributions} for three cases: $\xi_b=\xi_d$, $\xi_b>\xi_d$ and $\xi_b<\xi_d$.

From Figure \ref{fig:distributions}, in which $m = 1$, so that $\xi^* = \xi_b$, the convergence of the pdf and cdf is evident. 
For $n=10$, the distribution functions of the scaled processes are almost aligned with the normal distribution already. 
For $\xi_b=\xi_d$, the convergence is fastest. 
This can be explained by observing that in cases where $\xi_b\neq \xi_d$, the parameter $\xi_d$ still plays a role in pre-limit systems, whereas it does not appear in the normal limit. 
In the cases where $\xi_b\neq \xi_d$ we furthermore see that the functions are not smooth around $v_n=0$ or $v^* = -\sqrt{n}m/\xi^*$, which is the zero-inventory level in the original (unscaled) process.
As $n$ increases, this point of irregularity goes to $-\infty$ and therefore disappears.

\subsection{Approximations to performance metrics}

The plots in the previous section indicate that the normal approximation gives simple yet accurate approximations to the stationary distribution of the inventory process. 
We now assess if this also translates to the performance measures.
Again, we choose to fix the parameters $\lambda_b$ and $\lambda_d$, and evaluate the system with arrival rates $n\lambda_b$ and $n\lambda_d$ for increasing $n$. 
First, the normal approximation in \eqref{eq:normal_approximation} yields the following approximation for the expected inventory level:
\begin{equation}
\E[X_n] \approx \frac{ n m }{\xi^*} = \frac{n (\lambda_b \E[B] - \lambda_d \E[D]) }{\xi^*}.
\label{eq:approx1}
\end{equation}
For the probability of negative inventory, we have
\begin{equation}
\pi_d = \P(X_n < 0) \approx \P\left( Z^* < -\sqrt{n}\,m/\xi^* \right) = \F\left(-\sqrt{n/2\xi^*}\,m/\sigma\right).
\label{eq:approx2}
\end{equation}
Last, the probability of demand being satisfied immediately is approximately
\begin{equation}
\P({\rm demand\ satisfied }) = \P( X_n > D ) \approx 1 - \int_0^\infty \F\left({-}\frac{\sqrt{2\xi^*}}{\sigma}\, \frac{x-nm/\xi^*}{\sqrt{n}}\right) \,\dd F_d(x).
\label{eq:approx3}
\end{equation}

\begin{remark}
Note that if $\lambda_b$ and $\lambda_d$ are large themselves, the parameter $n$ can be eliminated from \eqref{eq:approx1}-\eqref{eq:approx3}, so that
\begin{equation*}
\E[X] \approx \frac{m}{\xi^*}, \qquad
\pi_d \approx \F\left({-}m/(\sigma\sqrt{2\xi^*})\right),
\end{equation*}
\begin{equation*}
\P({\rm demand\ satisfied}) \approx 1-\int_0^\infty \F\left({-} \sqrt{2\xi^*}\,\frac{x-m/\xi^*}{\sigma}\right)\, \dd F_d(x),
\end{equation*}
where $m = \lambda_b \E[B] - \lambda_d \E[D]$ and $\sigma^2 = \lambda_b \E[B^2] + \lambda_d \E[D^2]$.
\end{remark}

We will now test these approximations under various assumptions on the distribution of $B$ and $D$.
In Tables \ref{tab:accuracy_deterministic}-\ref{tab:accuracy_gamma} we compare the values obtained through the normal approximation against the true values obtained through numerical evaluation (for exponential jump sizes only) and simulation. 
All simulation results are accurate up to a 95\% confidence interval of width $10^{-4}$.
We set $\lambda_b=1$ and $\lambda_d=0.5$ and let the mean jump sizes be equal to 1, i.e.~ $\E[B] = 1$ and $\E[B] = 1$ in all numerical experiments. 
In Table \ref{tab:accuracy_deterministic}, we let the jump sizes be deterministic, so that $\Var B = \Var D = 0$. 
Table \ref{tab:accuracy_exponential} shows the results in case of exponential jump sizes, so that $\Var B = \Var D = 1$. 
Last, in Table \ref{tab:accuracy_gamma} we investigate the quality of the approximation for jump sizes that follow a Gamma$(0.25,0.25)$ distribution, yielding $\Var B = \Var D = 4$. 
With this set-up we cover jump distributions of increasing variance, so that we are able to study the impact of increased variability on the accuracy of the approximations.
Moreover, we investigate the influence of the decay parameters $\xi_b$ and $\xi_d$ by considering the scenarios $\xi_b=\xi_d$, $\xi_b<\xi_d$ and $\xi_b>\xi_d$.

We make a couple of observations based on the numbers in Tables \ref{tab:accuracy_deterministic}-\ref{tab:accuracy_gamma}. 
First, we see that the approximation for the mean blood inventory level $\E[X_n]$ is exact if $\xi_b=\xi_d$, see Proposition \ref{prop:mean_inventory}.
This obviously does not extend to $\pi_d$ and $\P({\rm demand\ satisfied})$, since these performance measures are based on the entire distribution of $X_n$ rather than the mean. 
Nonetheless, the normal approximation appears to be very accurate in the case $\xi_b=\xi_d$. 
We may explain this by observing that in the approximations \eqref{eq:approx1}-\eqref{eq:approx3}, only $\xi^*$ appears. 
In our setting, we have $m = \l_b - \l_d = 0.5$, so that $\xi^* = \xi_b$.
If $\xi_b\neq \xi_d$, then the value of $\xi_d$ plays a role in pre-limit systems, which induces inaccuracies in the approximation of performance measures. 
In case $\xi_b = \xi_d$, we have $\xi^* = \xi_b = \xi_d$, so that this discrepancy is overcome. 

Moreover, since $m>0$, we see that $\pi_d\to 0$ and $\P({\rm demand\ satisfied}) \to 1$ as $n$ increases.
This is due to the observation that as $n$ grows large, the inventory process concentrates around the level $nm$ with fluctuations of order $\sqrt{n}$, so that the process stays away from level zero, see Figure \ref{fig:sample_paths_fluid1}.
The approximations \eqref{eq:approx2}-\eqref{eq:approx3} adequately capture this convergence.  

As expected, the accuracy of the approximations increases with $n$. 
Moreover, increased variability in the jump distributions appears to cause a decrease in accuracy. 
However, for all cases considered in Tables \ref{tab:accuracy_deterministic}-\ref{tab:accuracy_gamma}, the normal approximations \eqref{eq:approx1}-\eqref{eq:approx3} seem to yield relatively sharp estimates for the relevant performance measures under various assumptions on the distributions of the jump sizes.

\begin{table}
\centering
\begin{subtable}{0.99\textwidth}\centering
\begin{tabular}{|r|rr|rr|rr|}
\cline{2-7}\multicolumn{1}{r|}{} & \multicolumn{2}{c|}{$\E[X_n]$} & \multicolumn{2}{c|}{$\pi_d$} & \multicolumn{2}{c|}{$\P({\rm dem.sat.})$} \bigstrut\\
\hline
$n$     & Sim.   & \eqref{eq:approx1} & Sim.   & \eqref{eq:approx2} & Sim.   & \eqref{eq:approx3} \bigstrut\\
\hline
1     & 0.500 & 0.500 & 0.2702 & 0.2819 & 0.2598 & 0.2819 \bigstrut[t]\\
2     & 1.000 & 1.000 & 0.2014 & 0.2071 & 0.4859 & 0.5000 \\
5     & 2.500 & 2.500 & 0.0943 & 0.0984 & 0.7814 & 0.7807 \\
10    & 5.000 & 5.000 & 0.0316 & 0.0339 & 0.9306 & 0.9279 \\
20    & 10.000 & 10.000 & 0.0043 & 0.0049 & 0.9908 & 0.9899 \\
50    & 25.000 & 25.000 & 0.0000 & 0.0000 & 1.0000 & 1.0000 \bigstrut[b]\\
\hline
\end{tabular}%
\caption{$\xi_b = 1$, $\xi_d=1$.}
\end{subtable}

\vspace{5mm}

\begin{subtable}{0.99\textwidth}\centering
\begin{tabular}{|r|rr|rr|rr|}
\cline{2-7}\multicolumn{1}{r|}{} & \multicolumn{2}{c|}{$\E[X_n]$} & \multicolumn{2}{c|}{$\pi_d$} & \multicolumn{2}{c|}{$\P({\rm dem.sat.})$} \bigstrut\\
\hline
$n$     & Sim.   & \eqref{eq:approx1} & Sim.   & \eqref{eq:approx2} & Sim.   & \eqref{eq:approx3} \bigstrut\\
\hline
1     & 0.584 & 0.500 & 0.2522 & 0.2819 & 0.2712 & 0.2819 \bigstrut[t]\\
2     & 1.086 & 1.000 & 0.1809 & 0.2071 & 0.5020 & 0.5000 \\
5     & 2.558 & 2.500 & 0.0837 & 0.0984 & 0.7911 & 0.7807 \\
10    & 5.024 & 5.000 & 0.0286 & 0.0339 & 0.9335 & 0.9279 \\
20    & 10.006 & 10.000 & 0.0040 & 0.0049 & 0.9912 & 0.9899 \\
50    & 25.000 & 25.000 & 0.0000 & 0.0000 & 1.0000 & 1.0000 \bigstrut[b]\\
\hline
\end{tabular}%
\caption{$\xi_b =1$, $\xi_d=2$.}
\end{subtable}

\vspace{5mm}

\begin{subtable}{0.99\textwidth}\centering
\begin{tabular}{|r|rr|rr|rr|}
\cline{2-7}\multicolumn{1}{r|}{} & \multicolumn{2}{c|}{$\E[X_n]$} & \multicolumn{2}{c|}{$\pi_d$} & \multicolumn{2}{c|}{$\P({\rm dem.sat.})$} \bigstrut\\
\hline
$n$     & Sim.   & \eqref{eq:approx1} & Sim.   & \eqref{eq:approx2} & Sim.   & \eqref{eq:approx3} \bigstrut\\
\hline
1     & 0.158 & 0.250 & 0.3308 & 0.3415 & 0.1006 & 0.1103 \bigstrut[t]\\
2     & 0.397 & 0.500 & 0.2973 & 0.2819 & 0.2465 & 0.2819 \\
5     & 1.164 & 1.250 & 0.1952 & 0.1807 & 0.5482 & 0.5724 \\
10    & 2.447 & 2.500 & 0.1036 & 0.0984 & 0.7729 & 0.7807 \\
20    & 4.980 & 5.000 & 0.0340 & 0.0339 & 0.9283 & 0.9279 \\
50    & 12.497 & 12.500 & 0.0017 & 0.0019 & 0.9964 & 0.9960 \bigstrut[b]\\
\hline
\end{tabular}%
\caption{$\xi_b = 2$, $\xi_d=1$.}
\end{subtable}
\caption{Accuracy of diffusion approximation for the blood inventory process $\E[X_n]$, the probability of negative inventory $\pi_d$ and the probability of demand being fully satisfied $\P(dem.sat)$, with arrival rates $n\l_b = n$ and $n\l_d = 0.5n$ and deterministic jump sizes, $B\equiv 1$ and $D\equiv 1$.}
\label{tab:accuracy_deterministic}
\end{table}

\begin{table}
\centering
\begin{subtable}{0.99\textwidth}\centering
\begin{tabular}{|r|rr|rr|rr|}
\cline{2-7}\multicolumn{1}{r|}{} & \multicolumn{2}{c|}{$\E[X_n]$} & \multicolumn{2}{c|}{$\pi_d$} & \multicolumn{2}{c|}{$\P({\rm dem.sat.})$} \bigstrut\\
\hline
$n$     & Exact   & \eqref{eq:approx1} & Exact   & \eqref{eq:approx2} & Exact   & \eqref{eq:approx3} \bigstrut\\
\hline
1     & 0.500 & 0.500 & 0.2929 & 0.3415 & 0.3536 & 0.3925 \\
2     & 1.000 & 1.000 & 0.2500 & 0.2819 & 0.5000 & 0.5135 \\
5     & 2.500 & 2.500 & 0.1642 & 0.1807 & 0.7062 & 0.7009 \\
10    & 5.000 & 5.000 & 0.0898 & 0.0984 & 0.8491 & 0.8418 \\
20    & 10.000 & 10.000 & 0.0307 & 0.0339 & 0.9506 & 0.9467 \\
50    & 25.000 & 25.000 & 0.0017 & 0.0019 & 0.9974 & 0.9970 \bigstrut[b]\\
\hline
\end{tabular}%
\caption{$\xi_b = 1$, $\xi_d=1$.}
\end{subtable}

\vspace{5mm}

\begin{subtable}{0.99\textwidth}\centering
\begin{tabular}{|r|rr|rr|rr|}
\cline{2-7}\multicolumn{1}{r|}{} & \multicolumn{2}{c|}{$\E[X_n]$} & \multicolumn{2}{c|}{$\pi_d$} & \multicolumn{2}{c|}{$\P({\rm dem.sat.})$} \bigstrut\\
\hline
$n$     & Exact   & \eqref{eq:approx1} & Exact   & \eqref{eq:approx2} & Exact   & \eqref{eq:approx3} \bigstrut\\
\hline
1     & 0.621 & 0.500 & 0.2589 & 0.3415 & 0.3705 & 0.3925 \bigstrut[t]\\
2     & 1.153 & 1.000 & 0.2164 & 0.2819 & 0.5224 & 0.5135 \\
5     & 2.656 & 2.500 & 0.1414 & 0.1807 & 0.7254 & 0.7009 \\
10    & 5.113 & 5.000 & 0.0784 & 0.0984 & 0.8598 & 0.8418 \\
20    & 10.050 & 10.000 & 0.0275 & 0.0339 & 0.9538 & 0.9467 \\
50    & 25.004 & 25.000 & 0.0016 & 0.0019 & 0.9975 & 0.9970 \bigstrut[b]\\
\hline
\end{tabular}%
\caption{$\xi_b =1$, $\xi_d=2$.}
\end{subtable}

\vspace{5mm}

\begin{subtable}{0.99\textwidth}\centering
\begin{tabular}{|r|rr|rr|rr|}
\cline{2-7}\multicolumn{1}{r|}{} & \multicolumn{2}{c|}{$\E[X_n]$} & \multicolumn{2}{c|}{$\pi_d$} & \multicolumn{2}{c|}{$\P({\rm dem.sat.})$} \bigstrut\\
\hline
$n$     & Exact   & \eqref{eq:approx1} & Exact   & \eqref{eq:approx2} & Exact   & \eqref{eq:approx3} \bigstrut\\
\hline
1     & 0.125 & 0.250 & 0.3548 & 0.3864 & 0.2168 & 0.2942 \bigstrut[t]\\
2     & 0.333 & 0.500 & 0.3333 & 0.3415 & 0.3333 & 0.3925 \\
5     & 1.059 & 1.250 & 0.2647 & 0.2593 & 0.5264 & 0.5570 \\
10    & 2.333 & 2.500 & 0.1856 & 0.1807 & 0.6881 & 0.7009 \\
20    & 4.893 & 5.000 & 0.0995 & 0.0984 & 0.8400 & 0.8418 \\
50    & 12.475 & 12.500 & 0.0198 & 0.0206 & 0.9692 & 0.9678 \bigstrut[b]\\
\hline
\end{tabular}%
\caption{$\xi_b = 2$, $\xi_d=1$.}
\end{subtable}
\caption{Accuracy of diffusion approximation for the blood inventory process $\E[X_n]$, the probability of negative inventory $\pi_d$ and the probability of demand being fully satisfied $\P(dem.sat)$, with arrival rates $n\l_b = n$ and $n\l_d = 0.5n$ and exponentially distributed jump sizes, $B\sim \exp(1)$ and $D\sim\exp(1)$.}
\label{tab:accuracy_exponential}
\end{table}

\begin{table}
\centering
\begin{subtable}{0.99\textwidth}\centering
\begin{tabular}{|r|rr|rr|rr|}
\cline{2-7}\multicolumn{1}{r|}{} & \multicolumn{2}{c|}{$\E[X_n]$} & \multicolumn{2}{c|}{$\pi_d$} & \multicolumn{2}{c|}{$\P({\rm dem.sat.})$} \bigstrut\\
\hline
$n$     & Sim.   & \eqref{eq:approx1} & Sim.   & \eqref{eq:approx2} & Sim.   & \eqref{eq:approx3} \bigstrut\\
\hline
1     & 0.500 & 0.500 & 0.3118 & 0.3981 & 0.4412 & 0.4636 \\
2     & 1.000 & 1.000 & 0.2894 & 0.3575 & 0.5343 & 0.5288 \\
5     & 2.500 & 2.500 & 0.2375 & 0.2819 & 0.6590 & 0.6381 \\
10    & 5.000 & 5.000 & 0.1785 & 0.2071 & 0.7592 & 0.7385 \\
20    & 10.000 & 10.000 & 0.1090 & 0.1241 & 0.8593 & 0.8454 \\
50    & 25.000 & 25.000 & 0.0303 & 0.0339 & 0.9624 & 0.9583 \bigstrut[b]\\
\hline
\end{tabular}%
\caption{$\xi_b = 1$, $\xi_d=1$.}
\end{subtable}

\vspace{5mm}

\begin{subtable}{0.99\textwidth}\centering
\begin{tabular}{|r|rr|rr|rr|}
\cline{2-7}\multicolumn{1}{r|}{} & \multicolumn{2}{c|}{$\E[X_n]$} & \multicolumn{2}{c|}{$\pi_d$} & \multicolumn{2}{c|}{$\P({\rm dem.sat.})$} \bigstrut\\
\hline
$n$     & Sim.   & \eqref{eq:approx1} & Sim.   & \eqref{eq:approx2} & Sim.   & \eqref{eq:approx3} \bigstrut\\
\hline
1     & 0.667 & 0.500 & 0.2695 & 0.3981 & 0.4636 & 0.4636 \bigstrut[t]\\
2     & 1.253 & 1.000 & 0.2469 & 0.3575 & 0.5632 & 0.5288 \\
5     & 2.863 & 2.500 & 0.2009 & 0.2819 & 0.6895 & 0.6381 \\
10    & 5.385 & 5.000 & 0.1518 & 0.2071 & 0.7834 & 0.7385 \\
20    & 10.328 & 10.000 & 0.0938 & 0.1241 & 0.8739 & 0.8454 \\
50    & 25.124 & 25.000 & 0.0269 & 0.0339 & 0.9658 & 0.9583 \bigstrut[b]\\
\hline
\end{tabular}%
\caption{$\xi_b =1$, $\xi_d=2$.}
\end{subtable}

\vspace{5mm}

\begin{subtable}{0.99\textwidth}\centering
\begin{tabular}{|r|rr|rr|rr|}
\cline{2-7}\multicolumn{1}{r|}{} & \multicolumn{2}{c|}{$\E[X_n]$} & \multicolumn{2}{c|}{$\pi_d$} & \multicolumn{2}{c|}{$\P({\rm dem.sat.})$} \bigstrut\\
\hline
$n$     & Sim.   & \eqref{eq:approx1} & Sim.   & \eqref{eq:approx2} & Sim.   & \eqref{eq:approx3} \bigstrut\\
\hline
1     & 0.081 & 0.250 & 0.3694 & 0.4276 & 0.3270 & 0.4104 \bigstrut[t]\\
2     & 0.238 & 0.500 & 0.3593 & 0.3981 & 0.4137 & 0.4636 \\
5     & 0.857 & 1.250 & 0.3237 & 0.3415 & 0.5311 & 0.5528 \\
10    & 2.045 & 2.500 & 0.2739 & 0.2819 & 0.6282 & 0.6381 \\
20    & 4.568 & 5.000 & 0.2039 & 0.2071 & 0.7361 & 0.7385 \\
50    & 12.231 & 12.500 & 0.0966 & 0.0984 & 0.8797 & 0.8779 \bigstrut[b]\\
\hline
\end{tabular}%
\caption{$\xi_b = 2$, $\xi_d=1$.}
\end{subtable}
\caption{Accuracy of diffusion approximation for the blood inventory process $\E[X_n]$, the probability of negative inventory $\pi_d$ and the probability of demand being fully satisfied $\P(dem.sat)$, with arrival rates $n\l_b = n$ and $n\l_d = 0.5n$ and Gamma distributed jump sizes, $B\sim\text{Gamma}(0.25,0.25)$ and $D \sim\text{Gamma}(0.25,0.25)$.}
\label{tab:accuracy_gamma}
\end{table}

\section{Conclusions \& suggestions for further research}
\label{conclus}
In this chapter, we studied a stochastic model for a blood bank.
We have presented a global approach to the model in its full generality,
and obtained very detailed exact expressions for the densities of amount of inventory and amount of demand (shortage)
in special cases (exponential amounts of donated and requested blood; and either $\xi_b=\xi_d=0$ or $\alpha_b=\alpha_d=0$).
Moreover, we have shown how an appropriate scaling, for the model in full generality, leads to an Ornstein-Uhlenbeck diffusion process,
which can be used as a tool to obtain simple yet accurate approximations for some key performance measures.

Our model is a two-sided model, in the sense that we simultaneously consider
the amount of blood in inventory and the amount of demand (shortage), one of the two at any time being zero.
Such two-sided processes arise in many different settings, and thus are of considerable interest.
The present setting is reminiscent of an organ transplantation problem, where there is
either a queue of persons waiting to receive an organ,
or a queue of donor organs. The perishability/impatience aspect features there too \cite{Boxma2011}.
A quite different setting is that of insurance risk. We refer to Albrecher \& Lautscham \cite{Albrecher2013}
who extend the classical Cram\'er-Lundberg insurance risk model by allowing the capital of an insurance company
to become negative -- a situation that is usually indicated by ``ruin" in the insurance literature.
Their process thus becomes two-sided. The capital might become positive again; however, 
at a rate $\omega(x)$ when the capital has a negative value $-x$, bankruptcy is declared and the process ends.
Interestingly, similar special functions (like hypergeometric functions) play a role in \cite{Albrecher2013} and in the present study.

The analyses performed in this chapter, which evolved around a simplified version of the inventory process of a blood bank, revealed some interesting avenues for further research. 
We name a couple of them.

First, we remark that our results are restricted to one type of blood. 
Naturally, it would be very interesting to extend the analysis to multiple types of blood.

Another important extension would be to use our results to facilitate the decision process that is faced by the CBB on a daily basis:
Which amounts of blood, and of which types, should today be sent to the local blood banks (hospitals)?
Knowing that, e.g., blood types $O^-,A^-,B^-,AB^-$ can satisfy the corresponding $+$ type (but not vice versa),
one may try to optimize the blood allocation process on the basis of actual amounts of blood present.

Finally, we mention a significant open research question regarding the process limits that we derived in Section \ref{sectionscaling}, of which the steady-state distributions were used to approximate steady-state performance measures in pre-limit systems. 
As we pointed out earlier, the justification that the steady-state distribution of the scaled inventory process indeed converges to the steady-state distribution of the fluid (cq.~diffusion) limit requires a rigorous argument why the order of limits $n\to\infty$ and $t\to \infty$ may be interchanged. 
Proving interchange-of-limits statements typically raises many technical challenges, see e.g.~\cite{Dai2014a,Gamarnik2013a,Gurvich2013,Gamarnik2006} for works tackling this issue in the context of queues in heavy traffic. 
The usual approach is to prove tightness of the sequence of steady-state distributions of pre-limit, followed by applying Prokhorov's theorem, see e.g.~\cite[Sec.~1.5]{Billingsley1995}.
For our model, such an approach seems to be straightforward for the fluid scaling, since our inventory process can be upper (cq.~lower) bounded by a shot-noise process with only positive (cq.~negative) jumps. 
Of the latter, the steady-state behavior is known. 
This allows us to derive a uniform bound on the absolute mean of the stationary fluid-scaled process, which gives tightness.
The final step uses the deterministic nature of the differential equation governing the dynamics of the fluid limit, by which the steady-state distribution must be unique.
For the diffusion-scaled process, the steps towards proving the interchange-of-limits are not obvious and hence this needs further investigation. 
Our numerical results for various jump size distributions, however, support the conjecture that  this interchange is indeed valid. 

\section*{Appendix}

\addcontentsline{toc}{section}{\hspace{7.1mm} Appendix}

\begin{subappendices}

\settocdepth{chapter}

\section{Transformation integral equation}
\label{app:transformation_int}
In this appendix we show how integral equation (\ref{eq:demand2}) can be transformed into a second-order differential equation,
in the case of exponential $F_b(\cdot)$ and $F_d(\cdot)$.
Differentiate (\ref{eq:demand2}) w.r.t.\ $v$:
\begin{align}
& \lambda_d f(v) - \mu_d \left[\lambda_d \int_0^v f(y) {\rm e}^{-\mu_d(v-y)} {\rm d}y +
\lambda_d \int_0^{\infty} g(y) {\rm e}^{-\mu_d (v+y)} {\rm d}y\right]
\nonumber
\\
&\qquad = 
-\lambda_b f(v) + \lambda_b \mu_b \int_v^{\infty} f(y) {\rm e}^{-\mu_b(y-v)} {\rm d}y + \xi_d f(v) +\xi_d
v f'(v) .
\label{eq:demand3}
\end{align}
Using (\ref{eq:demand2}) once more, now to replace the term between square brackets in (\ref{eq:demand3}),
we get:
\begin{align}
\xi_d v f'(v) &= (\lambda_d +\lambda_b - \xi_d) f(v)
\nonumber
\\
&\qquad - \mu_d \left(\lambda_b \int_v^{\infty} f(y) {\rm e}^{-\mu_b(y-v)} {\rm d}y + \xi_d v f(v)\right)
\nonumber
\\
&\qquad \qquad - \mu_b \lambda_b \int_v^{\infty} f(y) {\rm e}^{-\mu_b(y-v)} {\rm d}y,
\label{eq:demand4}
\end{align}
and once more differentiating w.r.t.\ $v$ then gives:
\begin{align}
&\xi_d v f''(v) + \xi_d f'(v) - (\lambda_d +\lambda_b -\xi_d -\mu_d \xi_d v) f'(v)
\nonumber
\\
& \qquad = -\mu_d \xi_d f(v) +(\mu_b+\mu_d) \lambda_b f(v) -\mu_b(\mu_b+\mu_d) \lambda_b \int_v^{\infty} f(y) {\rm e}^{-\mu_b(y-v)} {\rm d}y .
\label{eq:demand5}
\end{align}
The integral that appears in \eqref{eq:demand4} can be eliminated by using (\ref{eq:demand5}),
and we thus finally obtain the following second order homogeneous differential equation:
\begin{align}
&\xi_d v f''(v) + \left(2\xi_d -\lambda_d -\lambda_b + \mu_d\xi_dv -\mu_b \xi_d v\right)f'(v) \nonumber\\
& +\left(\mu_d\xi_d -\mu_b\xi_d -\mu_d\lambda_b + \mu_b \lambda_d -\mu_b \mu_d \xi_d v\right)f(v) =0.
\end{align}

\section{Proof of Proposition \ref{densityProp}}
\label{app:proof_prop_density}
In the proof, we concentrate on the derivation of $f(v)$, which is the solution to 
\begin{align}
\xi_d v f''(v) &+ \left(2\xi_d -\lambda_d -\lambda_b + \mu_d\xi_dv -\mu_b \xi_d v\right)f'(v)  \nonumber \\
& \qquad + \left(\mu_d\xi_d -\mu_b\xi_d -\mu_d\lambda_b + \mu_b \lambda_d -\mu_b \mu_d \xi_d v\right)f(v) =0
\label{eq:demand6_1}
\end{align}
The expression for $g(v)$ follows directly from exchanging $\lambda_b$ with $\lambda_d$, $\mu_b$ with $\mu_d$, $\xi_b$ with $\xi_d$, and $\pi_b$ with $\pi_d$ in $f(v)$.
We rewrite \eqref{eq:demand6_1} as follows:
\begin{equation}\label{diffvgl}
vf''(v) + (A+Bv)f'(v)+(C+Dv)f(v)=0,
\end{equation}
where 
\begin{equation*}
A = 2-\frac{\lambda_b+\lambda_d}{\xi_d}, ~~
B = \mu_d-\mu_b, ~~
C = \mu_d-\mu_b + \frac{\lambda_d\mu_b-\lambda_b\mu_d}{\xi_d}, ~~
D = -\mu_b\mu_d.
\end{equation*}
Note that we divided both sides of equation \eqref{eq:demand6_1} by $\xi_d$ here.
We will try to transform the differential equation into one for which the solution is easily derived. In order to do so, we first guess $f$ to be of the form $f(v) = \ee^{\beta v}h(v)$, where $\beta$ is a constant and $h$ another real-valued function. Substituting this into \eqref{diffvgl} gives
\begin{equation}\label{eq:diff2}
v h''(v) + \left[ (2\beta +B)v + A\right]h'(v) + \left[(\beta^2+B\beta+D)v+ A\beta+C\right] h(v)=0.
\end{equation}
Next, we would like to choose $\beta$ such that $\beta^2+B\beta+D=0$, that is 
\begin{equation*}
\beta = \frac{-B \pm \sqrt{B^2-4D}}{2},
\end{equation*}
which equals either $-\mu_d$ or $\mu_b$. Since the solution of \eqref{diffvgl} we are looking for is a density, and necessarily $f(v) = \ee^{\beta v}h(v) \rightarrow 0 $ as $v\rightarrow \infty$, we set $\beta$ equal to the negative root $-\mu_d$. Lastly, we apply a change of variable, $x = \delta v$, and $h(v) = w(x)$, so that \eqref{eq:diff2} is transformed into 
\begin{equation*}
x w''(x) + \left[ (2\beta+B)\delta^{-1} x+A\right] w'(x) + \delta^{-1}\left[A\beta+C\right] w(x) = 0.
\end{equation*}
By choosing $(2\beta+B)\delta^{-1} = -1$, i.e.
\begin{equation*}
\delta = {-}(2\beta+B) = \mu_b+\mu_d,
\end{equation*}
we obtain
\begin{equation*}
x w''(x) +[ A - x ] w'(x) + \delta^{-1}\left[A\beta+C\right] w(x) = 0,
\end{equation*}
which is known as Kummer's equation, $x w''(x) + (b-x) w'(x) - aw(x) = 0$, see \cite{Slater1960}, with parameters
\begin{align*}
a &= -\delta^{-1}\left[A\beta+C\right] = 1-\frac{\lambda_d}{\xi_d},\\
b &= A = 2-\frac{\lambda_b+\lambda_d}{\xi_d}.
\end{align*}

Kummer's equation has two linearly independent solutions, namely $w(x) =$\\ \noindent $M(a,b,x)$, where $M$ is Kummer's hypergeometric function, also denoted by \\ \noindent $ _1F_1(a,b,x)$, and $U(a,b,x)$, Tricomi's hypergeometric function. These are defined as, see \cite[Eq.~(1.3.1)]{Slater1960},
\begin{align*}
M(a,b,x) &= \sum_{n=0}^\infty \frac{(a)_n}{(b)_n n!} x^n,\\
U(a,b,x) &= \frac{\Gamma(b-1)}{\Gamma(1+a-b)}\,M(a,b,x) + \frac{\Gamma(b-1)}{\Gamma(a)}\,x^{1-b}\,M(1+a-b,2-b,x),
\end{align*}
where $(.)_n$ is the Pochhammer symbol, which is used to represent $(y)_n = y\cdot(y+1)\cdot...\cdot (y+n-1)$.
We can therefore deduce that $f(v)$ is of the form
\begin{equation*}
\ee^{\beta v}\left[ c_1\, M(a,b,\delta v) + c_2\, U(a,b,\delta v)\right],
\end{equation*}
or
\begin{equation*}
\ee^{-\mu_d v}\left[ c_1 M\left( 1-\tfrac{\lambda_d}{\xi_d}, 2-\tfrac{\lambda_b+\lambda_d}{\xi_d},(\mu_b+\mu_d) v\right) + 
c_2 U\left( 1-\tfrac{\lambda_d}{\xi_d}, 2-\tfrac{\lambda_b+\lambda_d}{\xi_d},(\mu_b+\mu_d) v\right)\right],
\end{equation*}
where $c_1$ and $c_2$ are constants. From \cite[p.~60]{Slater1960}, we have
\begin{equation*}
M(a,b,x) \sim \frac{\Gamma(b)}{\Gamma(a)}\ee^x x^{a-b}, \qquad \text{as } x\rightarrow \infty.
\end{equation*}
Hence, 
\begin{align*}
&\ee^{-\mu_d v} M\left( 1-\frac{\lambda_d}{\xi_d}, 2-\frac{\lambda_b+\lambda_d}{\xi_d},(\mu_b+\mu_d)v\right) \\
&\qquad \qquad \sim 
\frac{\Gamma(2-\frac{\lambda_b+\lambda_d}{\xi_d})}{\Gamma(1-\frac{\lambda_d}{\xi_d})}\ee^{\mu_b v}\left((\mu_b+\mu_d)v\right)^{\lambda_b/\xi_d-1}
\to \infty
\end{align*}
for all $\mu_b>0$, which leads us to conclude $c_1 = 0$. We deduce $c_2$ by exploiting the restriction that 
\begin{equation*}
\int_0^\infty f(v)\, \dd v = \pi_d,
\end{equation*}
where $\pi_d$ is the probability of positive demand. Hence
\begin{equation*}
\pi_d c_2^{-1} = \int_0^\infty \ee^{-\mu_d v} U\left( 1-\tfrac{\lambda_d}{\xi_d}, 2-\tfrac{\lambda_b+\lambda_d}{\xi_d},(\mu_b+\mu_d)v\right)\, dv.
\end{equation*}
By slightly transforming \cite[Eq.~(3.2.51)]{Slater1960}, we find
\begin{equation*}
c_2^{-1} = \frac{1}{\pi_d}\,\frac{\Gamma\left(\frac{\lambda_b+\lambda_d}{\xi_d}\right)}{\Gamma\left(1+\tfrac{\lambda_b}{\xi_d}\right) }\, _2F_1\left(1-\tfrac{\lambda_d}{\xi_d},1,1+\tfrac{\lambda_b}{\xi_d},-\tfrac{\mu_b}{\mu_d}\right)
,
\end{equation*}
where $_2F_1(a_1,a_2,a_3,x) := \sum_{n=0}^{\infty} \frac{(a_1)_n (a_2)_n}{(a_3)_n n!} x^n$
is the hypergeometric function of Gauss.

\section{Laplace Transforms for Coxian jumps}
\label{app:coxian}

We outline how the differential equation (\ref{diffeq}) is obtained.
We take Laplace transforms in (\ref{eq:demand}), considering its five terms and calling them $T_1, T_2, T_3, T_4$ and $T_5$, successively.
Equation (\ref{eq:demand})  then translates into
\begin{equation*}
T_1 + T_2 +T_3 = T_4 + T_5,
\label{Tequ}
\end{equation*}
where
\begin{align}
T_1 
&= \lambda_d \int_{v=0}^{\infty} {\rm e}^{-sv} \int_{y=0}^v f(y) \bar{F}_d(v-y) {\rm d}y {\rm d}v \nonumber \\
&=
\lambda_d \phi(s) \frac{1 - \E[{\rm e}^{-sD}]}{s},
\label{T-1}\\
T_2 
&= \lambda_d \int_{v=0}^{\infty} {\rm e}^{-sv} \int_{y=0}^{\infty} g(y) \bar{F}_d(v+y) {\rm d}y {\rm d}v \nonumber\\
&= \lambda_d \int_{y=0}^{\infty} {\rm e}^{sy} g(y) \int_{z=y}^{\infty} {\rm e}^{-sz} \bar{F}_d(z) {\rm d}z {\rm d}y ,
\label{T-2}\\
T_3 
&= \pi_0 \lambda_d \int_0^{\infty} {\rm e}^{-sy} \bar{F}_d(y) {\rm d}y ,\\
T_4 
&= \lambda_b \int_{v=0}^{\infty} {\rm e}^{-sv} \int_{y=v}^{\infty} f(y) \bar{F}_b(y-v) {\rm d}y {\rm d}v \nonumber\\
&= \lambda_b \int_{y=0}^{\infty} {\rm e}^{-sy} f(y) \int_{z=0}^{y} {\rm e}^{sz} \bar{F}_b(z) {\rm d}z {\rm d}y ,
\label{T-4}\\
T_5 
&= \xi_d \int_{v=0}^{\infty} v {\rm e}^{-sv} f(v) {\rm d}v + \alpha_d \phi(s) \nonumber\\
&= - \xi_d \phi'(s) + \alpha_d \phi(s) .
\label{T-5}
\end{align}

We now evaluate the terms appearing in the righthand sides of (\ref{T-1})-(\ref{T-4}) for the Coxian case
of (\ref{Fbarb}) and (\ref{Fdarb}):
\begin{align}
\int_{z=0}^y {\rm e}^{sz} \bar{F}_b(z) {\rm d}z &= 
\sum_{i=1}^K p_i \prod_{h=1}^{i-1} (1-p_h) \sum_{j=1}^i \prod_{l=1; l \neq j}^i \frac{\b_l}{\b_l - \b_j}
\frac{1}{\b_j-s} (1 - {\rm e}^{(s-\b_j)y}),
\label{hulp1}\\
\int_{z=y}^{\infty} {\rm e}^{-sz} \bar{F}_b(z) {\rm d}z &= 
\sum_{i=1}^K p_i \prod_{h=1}^{i-1} (1-p_h) \sum_{j=1}^i \prod_{l=1; l \neq j}^i \frac{\b_l}{\b_l - \b_j}
\frac{1}{\b_j+s} {\rm e}^{-(s+\b_j)y} ,
\label{hulp2}\\
\E[{\rm e}^{-sB}] &= 
\sum_{i=1}^K p_i \prod_{h=1}^{i-1} (1-p_h) \sum_{j=1}^i \prod_{l=1; l \neq j}^i \frac{\b_l}{\b_l - \b_j}
\frac{\b_j}{\b_j+s} ,
\label{hulp3}
\end{align}
and hence
\begin{equation}
\frac{1-\E[{\rm e}^{-sB}]}{s} = 
\sum_{i=1}^K p_i \prod_{h=1}^{i-1} (1-p_h) \sum_{j=1}^i \prod_{l=1; l \neq j}^i \frac{\b_l}{\b_l - \b_j}
\frac{1}{\b_j+s} .
\label{hulp4}
\end{equation}

Combining \eqref{Tequ} with \eqref{T-1}-\eqref{T-5}, and using \eqref{hulp1} and the counterparts of \eqref{hulp2} and \eqref{hulp4}
for $\bar{F}_d(\cdot)$, we find:
\begin{align}
& \lambda_d \phi(s) 
\sum_{i=1}^K q_i \prod_{h=1}^{i-1} (1-q_h) \sum_{j=1}^i \prod_{l=1; l \neq j}^i \frac{\d_l}{\d_l - \d_j}
\frac{1}{\d_j+s} 
\nonumber
\\
&\qquad + \lambda_d
\sum_{i=1}^L q_i \prod_{h=1}^{i-1} (1-q_h) \sum_{j=1}^i \prod_{l=1; l \neq j}^i \frac{\d_l}{\d_l - \d_j}
\frac{1}{\d_j+s} [\gamma(\d_j)  + \pi_0]
\nonumber
\\
&= \lambda_b 
\sum_{i=1}^K p_i \prod_{h=1}^{i-1} (1-p_h) \sum_{j=1}^i \prod_{l=1; l \neq j}^i \frac{\b_l}{\b_l - \b_j}
\frac{1}{\b_j-s} (\phi(s) - \phi(\b_j))
\nonumber
\\
&\qquad - \xi_d  \phi'(s) + \alpha_d \phi(s),
\label{eq:3star}
\end{align}
which is readily rewritten into \eqref{diffeq}.
\begin{remark}
If $\xi_d=0$, then $\phi(s)$ is obtained from \eqref{eq:3star}
in a standard manner, see also Section~\ref{sectionvariant}. 
\end{remark}
\begin{remark}
We now outline how \eqref{hulp2} and \eqref{hulp3} change when the $B_i$ have an Erlang-($l+1,\b$) distribution,
and when the $D_i$ have an Erlang-($k+1,\d$) distribution (see also \eqref{Fbarb} and the line below it);
\eqref{hulp1} and \eqref{hulp4} do not change (but of course $\E[{\rm e}^{-sD}]$ changes).
Firstly,
\begin{equation*}
\int_{z=0}^y {\rm e}^{sz} \bar{F}_b(z) {\rm d}z = 
\sum_{j=0}^l \frac{\b^j}{(\b-s)^{j+1}}
\left[1 - \sum_{i=0}^j {\rm e}^{-(\b-s)y} \frac{((\b-s)y)^i}{i!}\right] .
\end{equation*}
Term $T_4$ now becomes:
\begin{align*}
T_4 &= \lambda_b \int_{v=0}^{\infty} {\rm e}^{-sv} \int_{y=v}^{\infty} f(y) \bar{F}_b(y-v) {\rm d}y {\rm d}v\nonumber\\
 &= \lambda_b \sum_{j=0}^l \frac{\b^j}{(\b-s)^{j+1}}
\left[\phi(s) - \sum_{i=0}^j \frac{(\b-s)^i}{i!} \int_{y=0}^{\infty} y^i {\rm e}^{-\b y} f(y) {\rm d}y \right].
\end{align*}
It should be noted that $s = \b$ is a removable singularity. E.g., for $l=0$ one has
$T_4 = \lambda_b \frac{\phi(s) - \phi(\b)}{\b - s}$.
\\
Secondly,
\begin{equation*}
\int_{z=y}^{\infty} {\rm e}^{-sz} \bar{F}_b(z) {\rm d}z = 
\sum_{j=0}^k \frac{\d^j}{(s+\d)^{j+1}} \sum_{i=0}^j  {\rm e}^{-(s+\d)y} \frac{((s+\d)y)^i}{i!} .
\end{equation*}
Term $T_2$ now becomes:
\begin{align*}
T_2 &= \lambda_d \int_{v=0}^{\infty} {\rm e}^{-sv} \int_{y=0}^{\infty} g(y) \bar{F}_d(v+y) {\rm d}y\, {\rm d}v \nonumber\\
&=
\lambda_d  
\sum_{j=0}^k \frac{\d^j}{(s+\d)^{j+1}} \sum_{i=0}^j \frac{(s+\d)^i}{i!} 
\int_{y=0}^{\infty} y^i {\rm e}^{-\d y} g(y)\, {\rm d}y .
\end{align*}
It is readily seen that the resulting counterpart of \eqref{eq:3star} can again be written in the form \eqref{diffeq},
and hence the solution is formally still given by \eqref{diffeqsoln}.
\label{RmErlang}
\end{remark}

\resettocdepth

\end{subappendices}

\addcontentsline{toc}{chapter}{Bibliography}
\bibliography{bibliography}

\newpage

\fancyhead[LO]{\itshape Summary}
\fancyhead[RE]{\itshape Summary}

\addcontentsline{toc}{chapter}{Summary}

\chapter*{Summary}

\vspace{-1cm}
\noindent\rule[0.5ex]{\linewidth}{2pt}


\section*{Asymptotic dimensioning of stochastic service systems}

Stochastic service systems describe situations in which customers
compete for service from scarce resources. Think of check-in
lines at airports, waiting rooms in hospitals or queues in supermarkets,
where the scarce resource is human manpower. Next
to these traditional settings, resource sharing is also important
in large-scale service systems such as the internet, wireless networks
and cloud computing facilities. In these virtual environments,
geographical conditions do not restrict the system size, paving the way for the emergence of large-scale
resource sharing networks. 
This thesis investigates how to design large-scale systems in order to achieve the dual goal of operational efficiency and quality-of-service, by which we mean that the system is highly occupied and hence efficiently utilizes the expensive resources, while at the same time, the level of service, experienced by customers, remains high. 

The intrinsic stochastic variability of arrival and service processes is the predominant cause of delays experienced by customers. 
Queueing theory and stochastics provide the tools to describe and evaluate congestion in these systems.
An important insight obtained through queueing analysis is the effect of resource pooling for systems with many servers and corresponding economies-of-scale that can be achieved by increasing the scale of the system.
Although classical queueing theory allows for exact evaluation of the performance of queueing systems of moderate size,
exact analysis becomes intractable as demand $R$ and capacity $s$ become large. In those cases, one typically resorts to
asymptotic approximation techniques, such as heavy-traffic diffusion approximations: the analysis of a
sequence of queueing processes, scaled in space, in which the server utilization level approaches 100\%. 
The resulting probabilistic limiting processes are easier to analyze. Moreover, the diffusion approximations have direct interpretations in
terms of the original systems and lead to tractable characterizations of their performance.

The heavy-traffic regime that plays a central role in this thesis is the Halfin-Whitt regime, also known as the Quality-and-Efficiency Driven (QED) regime, which dictates that capacity should be equal to the nominal demand plus an additional
variability hedge which is proportional to the square-root of the nominal load, i.e. $s = R + \beta\sqrt{R}$ for some $\beta>0$. The driving force behind this scaling regime is the central limit theorem (CLT). 
The rule $s=R+\beta\sqrt{R}$, commonly known as the square-root staffing
principle, has been proved to secure both efficiency (utilization approaches 100\%) and quality-of-service, since the mean waiting time is negligible under this scaling as the system grows large. 
Since the QED regime allows coexistence of the
two seemingly conflicting objectives in large-scale service systems, the paradigm has been implemented in a wide variety of
operational settings.
However, the standard QED regime fails to acknowledge features that play a dominant role in practice. This thesis
contributes to the existing literature by identifying these distinctive traits and showing how to account for them in a modified QED framework.

In Chapters 2 \& 3, we study how the limiting behavior of many-server queues is affected when one deviates from the standard square-root staffing principle.
In Chapter 2 we investigate a novel family of scaling regimes, in which the amount of overcapacity $s-R$ is not necessarily of the order $\sqrt{R}$, which gives rise to a novel family of heavy-traffic regimes and corresponding scaling limits.
Continuing our study of alternative scaling regimes, we investigate in Chapter 3 how to adapt the square-root staffing
paradigm in case the system faces demand patterns that are stochastically more volatile than anticipated. 
This phenomenon is known as overdispersion and can be caused by e.g.~the existence of correlation between the sources generating demand, or uncertainty about the arrival volume. 

In Chapters 4 \& 5, we review a family of queueing models in the QED regime in which the total number of customers that can reside in the system simultaneously is limited. 
As a result, customers may be denied access in case they find a full
system on arrival. 
This fraction of arrivals may either reattempt later or leave the system directly.
The impact of retrials on scaling rules in the QED regime is the focus of Chapter 4.
Since the volume of initially blocked customers is proportional to $\sqrt{R}$, that is, the same order as the variability hedge in the
staffing rule, retrials are prone to have a non-negligible effect on performance. 
We propose a heuristic method for the performance analysis of these types
of queueing models with finite-size restrictions, which is based on a fixed-point equation. As a by-product this yields a two-fold square-root staffing principle, which prescribes a synchronous scaling for both the system capacity and waiting space. 
Chapter 5 describes how these ideas can be applied in the context of an emergency department.

Chapter 6 studies a cost minimization problem in a single-server queue with non-stationary input. 
The bulk of the queueing literature concerns performance analysis assuming that steady state is reached. However, the validity of this assumption in practice is questionable, for the simple fact that no service system runs infinitely long. Moreover,
system parameters, such as the arrival volume, are likely to change over time. 
In this chapter, we characterize the error in performance metrics that follows from this transient nature of queues, and present a correction to the original staffing rule to account for the finite time horizon.

Finally, we analyze in Chapter 7 a specific stochastic service system: an inventory model of a blood bank with backlogs, perishable goods and consumer impatience. 
We obtain the stationary distribution of the inventory level, and deduce under appropriate scaling the stochastic process limit in terms of a diffusion process. 
This process limit allows for a more tractable approximate analysis of the model in case the number of blood deliveries and demand is large.



\end{document}